\documentclass[11pt,leqno]{article}

\usepackage{amsmath} 
\usepackage{amssymb}

 \makeatletter
 \@addtoreset{equation}{subsection}
 \makeatother

\newif\ifpdf 
\ifx\pdfoutput\undefined \pdffalse 
\else \pdfoutput=1 
\pdftrue \fi

\ifpdf
\usepackage[pdftex]{graphicx}
\usepackage[pdftex]{hyperref}
\else
\usepackage{graphicx}
\fi

\newtheorem{theorem}{Theorem}[subsection]
\newtheorem{lemma}[theorem]{Lemma}
\newtheorem{proposition}[theorem]{Proposition}
\newtheorem{definition}[theorem]{Definition}
\newtheorem{remark}[theorem]{Remark} 
\newtheorem{corollary}[theorem]{Corollary}
\newtheorem{example}[theorem]{Example} 
\newtheorem{notation}[theorem]{Notation}
\newtheorem{consequence}[theorem]{Consequence}

\def\C{{\mathbb C}}
 \def\N{{\mathbb N}} 
 \def\R{{\mathbb R}}

\def\T{{\mathbb T}}

\def\BB{\mathcal {B}} 
\def\CC{\mathcal {C}} 
\def\CD{\mathcal {D}}
\def\CE{\mathcal {E}}
\def\CF{\mathcal {F}}
\def\CG{\mathcal {G}}
\def\CH{\mathcal {H}}

\def\CJ{\mathcal {J}}

\def\CL{\mathcal {L}}

\def\CO{\mathcal {O}}
\def\CP{\mathcal {P}}
\def\CR{\mathcal {R}}
\def\CS{\mathcal {S}}

 \def\Im{\mathop{\rm Im}\nolimits}
\def\Re{\mathop{\rm Re}\nolimits}

\def\Op{\mathop{\rm Op}\nolimits}

\def\supp{\mathop {\rm supp}\nolimits} 

\def\sgn{\mathop{\rm sgn}\nolimits}

\def \som{\mathop{\sum}\limits}
\def \prodi{\mathop{\prod}\limits}
\def\loc{{\mathop{\rm loc}\limits}}

\def\bra{\langle}
 \def\ket{\rangle}

\def\v{\vert}

\def\Id{\mathop{\rm Id}\nolimits}

\def\cotg{\mathop{\rm cotg}\nolimits}
\def\Arc{\mathop{\rm Arc}\nolimits}

\def\square{\hbox{\vrule\vbox{\hrule\phantom{o}\hrule}\vrule}}
\def\cqfd{\hfill\square}

\def\ds{\displaystyle}

\title{Brief Article}
\author{The Author}

\title{Strichartz estimates for Schr\"odinger equations with 
variable coefficients}
\author{Luc Robbiano\thanks{Universit\'e de Versailles, UMR 8100, B\^at.
Fermat, 45, Avenue des Etats-Unis, 78035 Versailles}
\  and  Claude Zuily\thanks{Universit\'e Paris Sud, UMR 8628, D\'epartement de
Math\'ematiques, B\^at. 425, 91406 Orsay Cedex} }
\date{}

\begin{document}

\tableofcontents

\maketitle

\begin{abstract}
We prove the (local in time) Strichartz estimates (for the full range of parameters given by the
scaling unless the end point) for asymptotically flat and non trapping perturbations of the flat
Laplacian in $\R^n$, $n\geq 2$. The main point of the proof, namely the
dispersion estimate, is obtained in constructing a parametrix. The main tool for
this construction is the use of the FBI transform.
\end{abstract}

\section{Introduction and statement of the result}\label{sI}
The purpose of this work is to provide a proof of the full (local in time)
Strichartz estimates for the Schr\"odinger operator related to a non trapping
asymptotically flat perturbation of the usual Laplacian in $\R^n$.

To be more precise let us first introduce the following space. Let $\sigma_0$ be
in $]0,1[$. We set
\begin{equation}\label{eqI.0.1}
\BB_{\sigma_0}=\Big\{a\in C^\infty (\R^n) :\forall \alpha\in\N^n\,,
\enskip\exists  C_\alpha>0:  \vert \partial ^\alpha a(x)\vert \leq
\frac{C_\alpha}{\langle x\rangle^{1+\vert \alpha\vert +\sigma_0}}\,, \forall
x\in\R^n\Big\}
\end{equation}
Here we have set $\langle x\rangle=(1+\vert x\vert ^2)^{1/2}$.

Let $P$ be a second order differential operator,
\begin{equation}\label{eqI.0.2}
P=\sum^n_{j,k=1} D_j\big (g^{jk}(x)\, D_k\big)+\sum^n_{j=1}
(D_j\,b_j(x)+b_j(x)\,
D_j)+V(x)\,, D_j=\frac{1}{i} \frac{\partial }{\partial
x_j}\,,
\end{equation}
with principal symbol $p(x,\xi)=\som^n_{j,k=1}g^{jk}(x)\,\xi_j\,\xi_k$. (Here
$g^{jk}=g^{kj}$).

We shall make the following assumptions.
\begin{equation}\label{eqI.0.3}\left\{
\begin{array}{ll}
\textrm{(i)}\enskip &\textrm{The coefficients}\enskip  g^{jk},\enskip
b_j,\enskip V\enskip \textrm{are real valued,}\enskip 1\leq j\leq k\leq n\,.\\
\textrm{(ii)}\enskip &\textrm{There exists}\enskip \sigma_0>0\enskip
\textrm{such that}\enskip g^{jk}-\delta_{jk}\in\BB_{\sigma_0},\enskip
b_j\in\BB_{\sigma_0}.\\
& \textrm{Here}\enskip \delta_{jk}\enskip \textrm{is the Kronecker symbol}\,.\\
\textrm{(iii)}\enskip &V\in L^\infty (\R^n)\,.
\end{array}\right.
\end{equation}
\begin{equation}\label{eqI.0.4}
\textrm{There exists }\nu>0 \textrm{ such that for every } (x,\xi)\textrm{ in
}\R^n\times\R^n\,,\enskip p(x,\xi)\geq \nu\,\v\xi\v^2\,.
\end{equation}
Then $P$ has a self-adjoint extension with domain $H^2(\R^n)$.

Now we associate to the symbol $p$ the bicharacteristic flow given by the
following equations for $j=1,\ldots ,n$,
\begin{equation}\label{eqI.0.5}\left\{
\begin{array}{ll}
\dot x_j(t)=\frac{\partial p}{\partial \xi_j}\, (x(t),\,\xi(t))\,,\quad
&x_j(0)=x_j\,,\\
\dot \xi_j(t)=-\frac{\partial p}{\partial x_j}\, (x(t),\,\xi(t))\,,\quad
&\xi_j(0)=\xi_j\,.
 \end{array}\right.
\end{equation}
We shall denote by $(x(t,x,\xi),\xi(t,x,\xi))$ the solution, whenever it exists,
of the system (\ref{eqI.0.5}).
 In fact it is an easy consequence of (\ref{eqI.0.3}) and (\ref{eqI.0.4}) that
this flow exists for all $t$ in $\R$. Indeed by (\ref{eqI.0.4}) we have
$$
\nu\,\vert \xi(t)\vert ^2\leq p(x(t),\,\xi(t))=p(x,\xi)\,,
$$ 
and it follows from (\ref{eqI.0.4}) that
$$
\vert \dot x_j(t)\vert \leq 2\sum^n_{k=1} \vert g^{jk}(x)\,\xi_k(t)\vert \leq
C\,\vert \xi(t)\vert \leq C\,\nu^{-1/2}\, p(x,\xi)^{1/2}\,.
$$
Our last assumption will be the following.
\begin{equation}\label{eqI.0.6}
\textrm{For all }(x,\xi)\enskip \textrm{in}\enskip
T^*\R^n\setminus\{0\}\textrm{ we have }\lim_{t\rightarrow \pm\infty } \vert
x(t,x,\xi)\vert =+\infty \,.
\end{equation}
This means that the flow is not trapped backward nor forward. Now let us denote
by $e^{-itP}$ the solution of the following initial value problem
\begin{equation}\label{eqI.0.7}\left\{
\begin{array}{l}
i\, \frac{\partial u}{\partial t}-Pu=0\\
u(0,\cdot )=u_0\,.
\end{array}\right.
\end{equation}
Then the main result of this work is the following.

\begin{theorem}\sl  \label{tI.0.1} Assume that the operator $P$ satisfies the
conditions (\ref{eqI.0.3}), (\ref{eqI.0.4}), (\ref{eqI.0.6}). Let $T>0$ and
$(q,r)$ be a couple of real numbers such that $q>2$ and
$\frac{2}{q}=\frac{n}{2}-\frac{n}{r}$. Then there exists a positive
constant $C$ such that
\begin{equation}\label{eqI.0.8}
\Vert e^{-itP}u_0\Vert _{L^q([-T,T],L^r(\R^n))}\leq C\,\Vert u_0\Vert
_{L^2(\R^n)}\,,
\end{equation}
for all $u_0$ in $L^2(\R^n)$.
\end{theorem}

Such estimates are known in the litterature under the name of Strichartz
estimates. They have been proved for the flat Laplacian by Strichartz 
\cite{St} when $p=q=\frac{2n+4}{n}$ and extended to the full range of
$(p,q)$ given by the scaling by Ginibre-Velo \cite {GV}  and Yajima 
\cite{Y}. The limit case $q=2$  (the end point) when $n\geq 3$ is due to
Keel-Tao \cite {KT} . These estimates have been a key tool in the study of
non linear equations. Very recently several works appeared showing a new
interest for such estimates in the case of variable coefficients.
Steffilani-Tataru \cite {ST} proved Theorem \ref{tI.0.1} under conditions
(\ref{eqI.0.4}) and (\ref{eqI.0.6}) for compactly supported perturbations of
the flat Laplacian. In  \cite {B} Burq gave an alternative proof of this
result using the work of Burq-G\'erard-Tzvetkov \cite {BGT}. In the same
work Burq announced without proof that if you accept to replace in the right
hand side of (\ref{eqI.0.8}) the
$L^2$ norm by an $H^\varepsilon$ norm, for any small $\varepsilon>0$, then you
can weaken the decay hypotheses on the coefficients of $P$ in the sense that you
may replace in the definition (\ref{eqI.0.1}) of $\BB_{\sigma_0}$ the power
$\vert
\alpha\vert +1+\sigma_0$ by $\vert \alpha\vert +\sigma_0$. We have also to
mention a recent work of Hassel-Tao-Wunsch  \cite {HTW1} who proved in
dimension $n=3$ a weaker form of our result corresponding to the case where
$q=4$, $r=3$, under conditions similar to ours. Still more recently these three
authors announced the same result as ours under hypotheses on the coefficients
similar to ours (see \cite{HTW2}).

It is also worthwhile to mention the work of Burq-G\'erard-Tzvetkov who
investigate the Strichartz estimates on compact Riemannian manifolds. In that
case they show that such estimates hold with the $L^2$ norm replaced by the
$H^{1/q}$ norm. In the same paper these authors show that the same result holds
on $\R^n$ when the coefficients of their Laplacian (and its derivatives) are
merely bounded. Let us note also that these estimates concern also the wave
equation and many works have been devoted to this case. However we would like
to emphasize that, due to the finite speed of propagation, the extension to the
variable coefficients case appear to be much less technical (see \cite{SS}).

Let us now give some ideas on the proof. It is by now well known that a proof of
the Strichartz estimates can be done using a dispersion result, duality
arguments and the Hardy-Littlewood-Sobolev lemma. This has been formulated as an
abstract result in the paper  \cite {KT} as follows. Assume that for every
$t\in\R$ we have an operator $U(t)$ which maps $L^2(\R^n)$ to
$L^2(\R^n)$ and satisfies,
\begin{equation*}\left\{
\begin{array}{ll}
\textrm{(i)}\quad &\Vert U(t)\,f\Vert _{L^2(\R^n)}\leq C\,\Vert f\Vert
_{L^2(\R^n)}\,,\enskip
\forall t\in\R\,,\enskip C\textrm{ independent of }t\,.\\
\textrm{(ii)}\quad &\Vert U(s)(U(t))^*\,g\Vert _{L^\infty (\R^n)}\leq C\,\vert
t-s\vert ^{-\frac{n}{2}}\, \Vert g\Vert _{L^1(\R^n)}\,,\enskip t\not =s\,,
\end{array}\right.
\end{equation*}
then the Strichartz estimates (\ref{eqI.0.5}) hold for $U(t)$. 
It is not difficult to see that the serious estimate to be proved
is (ii). In the case when $U(t)=e^{it\Delta_0}$ (the flat Laplacian) this
estimate is obtained by the explicit formula giving the solution in term of the
data $u_0$. In the variable coefficients case such a formula is of course out of
hope and the better we can have is a parametrix. However due to strong technical
difficulties (which we try to explain below) which seem to be serious we are not
able to write a parametrix for $e^{-itP}$ so we have to explain what we do
instead. First of all let $\varphi _0\in C^\infty _0(\R^n)$ be such that
$\varphi _0(x)=1$ if $\vert x\vert \leq \frac{3}{2}$ and $\supp \varphi
_0\subset [-1,1]$. With a large $R>0$ we write
$$
e^{-itP}\, u_0(x)=\varphi _0\Big(\frac{x}{R}\Big )\,
e^{-itP}\,u_0(x)+\Big(1-\varphi _0\Big
(\frac{x}{R}\Big)\Big)\,e^{-itP}\,u_0(x)=v+w\,.
$$
It is not difficult to see that the Strichartz estimates for $v$ will be ensured
by the result of Staffilani-Tataru  \cite{ST} while the same estimate for
$w$ leads to consider an operator which is a small perturbation of the Laplacian
(see Section \ref{sII}).

Now it is not a surprise that microlocal analysis is strongly needed in our
proof. So let $\xi_0\in\R^n$, $\vert \xi_0\vert =1$ be a fixed direction. Let
$\chi_0\in C^\infty (\R)$, $\chi_0(s)=1$ if $s\leq \frac{3}{4}$, $\chi_0(s)=0$
if $s\geq 1$, $0\leq \chi_0\leq 1$ and let us set
$\chi_+(x)=\chi_0\big(\frac{-x\cdot \xi_0}{\delta_1}\big )$,
$\chi_-(x)=\chi_0\big(\frac{x\cdot \xi_0}{\delta_1}\big )$, $\delta_1>0$. We set
$U_+(t)=\chi_+\,e^{-itP}$, $U_-(t)=\chi_-\,e^{-itP}$. Now since
$\chi_+(x)+\chi_-(x)\geq 1$ for all $x$ in $\R^n$ then Strichartz estimates
separately for $U_+(t)$ and $U_-(t)$ will give the result. It is therefore
sufficient to prove the estimate (ii) above for $U_+(s)$
$(U_+(t))^*=\chi_+\,e^{i(s-t)P}\,\chi_+$ (and for $U_-(s)(U_-(t))^*$). In our
proof we shall construct a parametrix for these operators.

Our construction relies heavily on the theory of FBI transform (see Sj\"ostrand
 \cite{Sj} and Melin-Sj\"ostrand \cite{MS}) viewed as a Fourier integral
operator with complex phase. One of the reason of our choice is that in our
former works on the analytic smoothing effect  \cite{RZ2} we have already done
similar constructions (but only near the outgoing points: see below). Let us
explain very roughly the main ideas. The standard FBI transform is given by
\begin{equation}\label{eqI.0.9}
Tv(\alpha,\lambda)=c_n\, \lambda^{3n/4} \int_{\R^n} e^{i\lambda(y-\alpha_x)\cdot
\alpha_\xi-\frac{\lambda}{2} \vert y-\alpha_x\vert ^2+\frac{\lambda}{2}
\vert \alpha_\xi\vert ^2}\,v(y)\,dy
\end{equation}
where $\alpha=(\alpha_x,\alpha_\xi)\in\R^n\times\R^n$ and $c_n$ is a
positive constant.

Let us note that the phase can be written $i\lambda\varphi _0$ where $\varphi
_0(y,\alpha)=\frac{i}{2}$ $(y-(\alpha_x+i\alpha_\xi))^2$. Then $T$ maps
$L^2(\R^n)$ into the space $L^2(\R^{2n},e^{-\lambda\vert \alpha_\xi\vert
^2}\,d\alpha)$. The adjoint $T^*$ of $T$ is given by a similar formula (see
(\ref{eqVI.1.2}))   and we have,
\begin{equation}\label{eqI.0.10}
T^*T\textrm{ is the identity operator on }L^2(\R^n)\,.
\end{equation}
We embed the transform $T$ into a continuous family of FBI transform
\begin{equation}\label{eqI.0.11}\left\{
\begin{array}{l}
T_\theta v(\alpha,\lambda)=\lambda^{3n/4} \int_{\R^n} e^{i\lambda\varphi
(\theta,y,\alpha)} a(\theta,y,\alpha)\,v(y)\,dy\textrm{ with}\\
\varphi (0,y,\alpha)=\frac{1}{2}
(y-(\alpha_x+i\alpha_\xi))^2\,, a(0,y,\alpha)=c_n\,.
\end{array}\right.
\end{equation}
Let us set $U(\theta,t,\alpha,\lambda)=T_\theta[K_\pm(t)\,u_0](\alpha,\lambda)$,
where $K_\pm(t)=\chi_\pm\, e^{-itP}\,\chi_\pm$. Then it is shown that if
$\varphi $ satisfies the eikonal equation,
\begin{equation}\label{eqI.0.12}
\Big[\frac{\partial \varphi }{\partial \theta}+p\Big(x,\frac{\partial \varphi
}{\partial x}\Big)\Big](\theta,x,\alpha)=0\,,
\end{equation}
and if the symbol $a$ satisfies appropriate transport equations then $U$ is a
solution of the following equation
$$\Big(\frac{\partial U}{\partial t}+\lambda\,\frac{\partial U}{\partial
\theta}\Big)(\theta,t,\alpha,\lambda)\sim 0\,.
$$
It follows that essentially we have,
$U(\theta,t,\alpha,\lambda)=V(\theta-\lambda t,\alpha,\lambda)$. In particular
this shows that $U(0,t,\alpha,\lambda)=U(-\lambda t,0,\alpha,\lambda)$. Written
in terms of the transformations $T_\theta$ this reads
$$
T[K_\pm(t)u_0](\alpha,\lambda)=T_{-\lambda t}[\chi^2_\pm\,u_0](\alpha,\lambda)\,.
$$
Applying $T^*$ to both members and using (\ref{eqI.0.10}) we obtain
$$
K_{\pm}(t)\,u_0(x)=T^*\{T_{-\lambda t}[\chi^2_\pm\,u_0](\cdot ,\lambda)\}(t,x)\,.
$$
Thus we have expressed the solution in terms of the data through a Fourier
integral operator with complex phase.

This short discussion shows that as usual the main point of the proof is to
solve the eikonal and transport equations. Let us point out the main
difficulties which occur in solving these equations. They are of three types:
the bad behavior of the flow from incoming points and for large time, the global
(in $\theta,x)$ character of all our constructions and the mixing of $C^\infty $
coefficients and complex variables (coming from the non real character of our
phase). Let us discuss each of them. First of all whatever the method you use to
solve an eikonal equation (symplectic geometry or another one) a precise
description of the flow of the symbol $p$ is needed. Let us recall (see
(\ref{eqI.0.5})) that our flow $(x(t,x,\xi),\xi(t,x,\xi))$, issued from the
point
$(x,\xi)\in T^*\R^n\setminus 0$, is defined for all $t\in\R$. In the case of the
flat Laplacian we have $\xi(t,x,\xi)=\xi$ and $x(t,x,\xi)=x+2t\xi$. Let now
$(x,\xi)\in T^*\R^n\setminus\{0\}$ and assume that $x\cdot \xi\geq 0$. Then it
is easy to see that $\vert x(t,x,\xi)\vert ^2\geq \vert x\vert ^2+4t^2\,\vert
\xi\vert ^2$ for $t\geq 0$ so that $\vert x(t,x,\xi)\vert $ becomes larger and
larger while $x(t,x,\xi)$ may vanish for a large $t<0$. Such a point is called
``outgoing for $t\geq 0$'' and ``incoming for $t<0$''. In the case of a
perturbed Laplacian this distinction between the directions is very important.
Indeed although the flow from outgoing points for $t\geq 0$ is very well
described for $t\geq 0$ and has very similar properties to the flat case, it has
a bad behavior for $t<0$ in what concerns its derivatives with respect to
$(x,\xi)$. For instance $\frac{\partial x_j}{\partial \xi_k} (t,x,\xi)$ does not
behave at all as $2t\,\delta_{jk}$. This is of great importance and causes some
trouble in the proof. However still when $t<0$, the flow behaves correctly as
long as the point $(x(t,x,\xi),\xi(t,x,\xi))$ is outgoing for $t\geq 0$. Roughly
speaking that is the reason why we are not able to construct a parametrix for
$e^{-itP}$ while it is possible for the operator
$\chi_\pm\,e^{-itP}\chi_\pm$. The Section \ref{sIII} is entirely devoted to
a careful study of the flow. Let us now describe our method of resolution of the
eikonal equation. The classical method uses the ideas of symplectic geometry.
Roughly speaking the manifold constructed from the flow is a Lagrangian manifold
on which the symbol $\tau+p(x,\xi)$ is constant. If it projects (globally) and
clearly on the basis then it is a graph of some function $\varphi $ which is the
desired phase. However this general method leads immediately to a difficulty in
our case. Indeed since we want that for $\theta=0$
 the phase $\varphi $ coincides with the phase $\varphi _0$ of the FBI transform
(see (\ref{eqI.0.9})) which is non real, we should take, in solving the flow,
data which are non real, so the flow itself would be non real; but our symbol
has merely $C^\infty $ coefficients. To circumvent this difficulty a method has
been proposed by Melin-Sj\"ostrand \cite{MS} which uses the almost analytic
 machinery. Another method, different in spirit, that the one described
above and known under the name of ``Lagrangian ideals'', has been introduced by
H\"ormander  \cite{H}. Here the initial data in the flow are kept real. Let
us set $u_j(x,\xi)=\xi_j-\frac{\partial \varphi _0}{\partial x_j}
(x,\xi)=\xi_j-\alpha^j_\varepsilon-i(x_j-\alpha^j_x)$. Then obviously we have
$\{u_j,u_k\}=0$ if $j\not =k$ (where $\{,\}$ denotes the Poisson bracket). Now
let us set $v_j(\theta,x,\xi)=u_j(x(-\theta,x,\xi),\xi(-\theta,x,\xi))$,
$j=1,\ldots ,n$. Then for every $\theta$ in $\R$ the Poisson bracket of $v_j$
and $v_k$ still vanishes if $j\not =k$. Thus the ideal generated by the $v_j$'s
is closed under the Poisson bracket. The main step in H\"ormander's method is to
show that this ideal is generated by functions of the form
$\xi_j-\Phi_j(\theta,x,\alpha)$. This will imply that one can find a function
$\varphi =\varphi (\theta,x,\alpha)$ such that $\frac{\partial \varphi
}{\partial x_j} (\theta,x,\alpha)=\Phi_j(\theta,x,\alpha)$ and it turns out
that $\varphi $ is the desired phase. To achieve its main step, H\"ormander uses
a precise version of the Malgrange preparation theorem which is discussed in
 \cite{H}, tome~1. This is the way we chosed to use in our case. It occupies
all Section \ref{sIV} of the paper. The proof is made separately for
outgoing and incoming points. Since the $v_j$'s are defined by mean of the
backward flow, in both cases we encounter the difficulty caused by the bad
behavior of the flow from incoming points. As it can be seen many technical
difficulties arise in the procedure.

The next step in the proof is the resolution of the transport equations. Here
also the cases of outgoing and incoming points have to be separated. We have
also to be careful since these are first order equations with non real $C^\infty
$ coefficients. The first case is easier. Indeed due to the good behavior of the
flow and the decay of the perturbation one can cut the Taylor expansion of the
coefficients of the vector field to some order and thus reduce ourselves to the
case of polynomial coefficients. Then by classical holomorphic methods  one can
solve the equations modulo flat terms which will be enough for our purpose. In
the second case there is no more such an asymptotic and the situation is much
more intricate. So we use the classical idea which consists in straightening
the vector field. This forces us to enter in the almost analytic machinery of
Melin-Sj\"ostrand  \cite {MS} (see Section \ref{sV}). Of course all the
constructions made above are done microlocally and in a neighborhood of the
bicharacteristic. Therefore to define the general FBI transform $T_\theta$ (see
(\ref{eqI.0.11})) as well as to pass from the standard $T$ to $T_{-\lambda t}$
we have to insert many microlocal cut-off. Of course we have to check at each
microlocalization that the remainder leads to an acceptable error. This is the
goal of Section
\ref{sVI}. At this stage of the proof the operator
$K_\pm(t)=\chi_\pm\,e^{-itP}\,\chi_\pm$ is writen as
$$
K_\pm(t)\,u_0(x)=\int k_\pm (t,x,y)\,u_0(y)\,dy
$$
where
$$
k_\pm (t,x,y)=\int e^{i\lambda F(-\lambda t,x,y,\alpha)}\, a(\lambda
t,x,y,\alpha)\,d\alpha\,.
$$
Thus the dispersion estimate would follow from the bound
$$
\vert k_\pm(t,x,y)\vert \leq \frac{C}{\vert t\vert ^{n/2}}
$$
for $0<\vert t\vert \leq T$.

Here we have two regimes according to the fact that $\vert \lambda t\vert \geq
1$ or $\vert \lambda t\vert \leq 1$. In the first case on the support of
$a(\lambda t,x,y,\alpha)$ we could be very far from the critical point of $F$.
Fortunately the phase $F$ has enough convexity to produce the desired bound of
$k_\pm$. In the second regime we are close to the critical point of $F$ so we
expect a stationnary phase method to work. However since the phase $F$ is non
real and since the determinant of its Hessien in $\alpha$ degenerates in some
direction when $\vert \lambda t\vert \rightarrow 0$ we cannot apply the standard
results as they appear in \cite{H}. Instead, after a careful study of the
phase $F$ we use merely an integration by part method with an appropriate vector
field to conclude. This is done in Section \ref{sVII}. The rest of this part
is devoted, using the Littlewood-Paley theory, to the end of the proof of our
main Theorem.

Finally an Appendix gathers the proofs of some technical results used in the
paper.

{\sl Aknowledgments.} We would like to thank Nicolas Burq for useful
discussions at an earlier stage of the work.

\section{Preliminaries and reduction to the case of a small
perturbation of the Laplacian}\label{sII}

\subsection{Preliminaries}\label{ssII.1}

We begin by recalling several earlier results which will be used in the sequel.

The first result concerns the case of compactly supported perturbations of the
Laplacian.
\begin{theorem}\sl (Staffilani-Tataru  \cite{ST})\label{tII.1.1} 
Let $P$ be defined by (\ref{eqI.0.2}). Assume that $P$ satisfies
(\ref{eqI.0.4}), (\ref{eqI.0.6}) and
\begin{equation}\label{eqA.1'}
\textrm{for } j,k=1,\ldots ,n\,, \enskip g^{jk}-\delta_{jk}\,,\enskip b_j,V
\textrm{ are compactly supported.}
\end{equation}
Then the Strichartz estimates (\ref{eqI.0.8}) hold. 
\end{theorem}

The second result which we recall is the extension to the variable coefficients
case by Do\"i \cite{D} of the Kato smoothing effect. Let us introduce the
following space. We set for $s,\mu$ in $\R$
$$
H^s_{\mu}(\R^n)=\{u\in \CS' : \bra x\ket^\mu\,(I-\Delta)^{\frac s2}\,u\in
L^2(\R^n)\}
$$
with its standard norm.
\begin{theorem}\sl (Do\"i  \cite{D})\label{tII.1.2}
Let $P$ be defined by (\ref{eqI.0.2}) and assume it satisfies the conditions
(\ref{eqI.0.3}), (\ref{eqI.0.4}), (\ref{eqI.0.6}). Then for all $T>0$  and all
$\sigma>\frac 12$ one can find a constant $C\geq 0$ such that,
\begin{equation}\label{eqII.1.1}
\Vert e^{-itP}\,u_0\Vert _{L^2([-T,T],H^{\frac 12}_{-\sigma}(\R^n))}\leq
C\,\Vert u_0\Vert _{L^2(\R^n)}\,,
\end{equation} 
for all $u_0$ in $L^2(\R^n)$.
\end{theorem}
We shall also use the following result.
\begin{lemma}\sl  (Keel-Tao \cite{KT}) \label{lII.1.3}
Let $(X,dx)$ be a measure space, $H$ a Hilbert space and $T>0$. Suppose that for
each time $t\in [-T,T]$ we have an operator $U(t):H\rightarrow L^2(X)$ which
satisfies the following estimates.
\begin{itemize}
\item[(i)] There exists $C_1\geq 0$ such that for all $t\in [-T,T]$ and all
$f\in H$,
$$
\Vert U(t)\,f\Vert _{L^2(X)}\leq C_1\, \Vert f\Vert _H\,.
$$
\item[(ii)] There exists $C_2\geq 0$ such that for all $t,s\in[-T,T]$, $t\not
=s$ and all $g\in L^1(X)$,
$$
\Vert U(t)(U(s))^*\,g\Vert _{L^\infty (X)}\leq C_2\,\vert t-s\vert
^{-n/2}\,\Vert g\Vert _{L^1(X)}\,.
$$
\end{itemize}
Let $(q,r)$ be a couple of real numbers such that $q\geq 2$, $r<+\infty $ and
$\frac{2}{q}=\frac{n}{2}-\frac{n}{r}$. Then  there exists $C\geq 0$ such that
for all $f$ in $H$
$$
\Vert U(t)\,f\Vert _{L^q([-T,T],L^r(X))}\leq C\,\Vert f\Vert _H\,.
$$
\end{lemma}
This result will be used in the sequel with $H=L^2(\R^n)$, $X=\R^n$.

Finally let's recall the following technical lemma.

\begin{lemma}\sl (Christ-Kiselev \cite{CK}) \label{lII.1.4}
Let $X,Y$ be two Banach spaces and $K(t,s)$ be a continuous function taking its
values in $B(X,Y)$, the space of bounded linear mappings from $X$ to $Y$. Let
$-\infty \leq a<b\leq +\infty $ and set
\begin{equation*}
\begin{aligned}
Sf(t)&=\int^b_a K(t,s)\,f(s)\,ds\\
Wf(t)&=\int^t_a K(t,s)\,f(s)\,ds\,.
\end{aligned}
\end{equation*}
Let $1\leq p<q\leq +\infty $. Then if we can find a constant $C>0$ such that
$$
\Vert Sf\Vert _{L^q((a,b),Y)}\leq C\, \Vert f\Vert _{L^p((a,b),X)}
$$
it follows that
$$
\Vert Wf\Vert _{L^q((a,b),Y)}\leq \frac{2^{-2(\frac{1}{p}-\frac{1}{q})}\cdot
2C}{1-2^{-(\frac{1}{p}-\frac{1}{q})}}\, \Vert f\Vert _{L^p((a,b),X)}\,.
$$
\end{lemma}
Using these results we shall see that Theorem \ref{tI.0.1} will be a
consequence of the following Theorem.
\begin{theorem}\sl \label{tII.1.5}
Let us set $\Delta_g=\som^n_{j,k=1} \frac{\partial }{\partial x_j}\big
(g^{jk}\,\frac\partial {\partial x_k}\big)$ and assume that the conditions
(\ref{eqI.0.3}), (\ref{eqI.0.4}), (\ref{eqI.0.6}) are satisfied by $\Delta_g$.
Let $T>0$ and $(q,r)$ be a couple of real numbers such that $q>2$ and $\frac
2q=\frac n2-\frac nr$. Then there exists a positive constant $C$ such that
$$
\big\Vert e^{it\Delta_g}\,u_0\big\Vert _{L^q([-T,T],L^r(\R^n))}\leq C\,
\Vert u_0\Vert _{L^2(\R^n)}
$$
for all $u_0$ in $L^2(\R^n)$.
\end{theorem}

Let us  show how Theorem \ref{tII.1.5} implies Theorem \ref{tI.0.1}.

Let us set $I=[0,T]$. (The case $I=[-T,0]$ is symmetric). Using (\ref{eqI.0.2})
we can write
\begin{equation}\label{eqII.1.3}
i\,\partial _t\,u+\Delta_gu=-\Big (\sum^n_{j=1}
(D_j\,b_j)+V\Big)\,u-2\,\sum^n_{j=1} b_j\,D_ju=:F=F_1+F_2\,.
\end{equation}
It follows from Duhamel formula that
\begin{equation}\label{eqII.1.4}
e^{-itP}\,u_0=e^{it\Delta_g}\,u_0+i \int^t_0 e^{i(t-s)\Delta_g}\,[F(s,\cdot
)]\,ds\,.
\end{equation}
Using Theorem \ref{tII.1.5} we obtain
\begin{equation}\label{eqII.1.5}
\big\Vert e^{it\Delta_g}\,u_0\big\Vert _{L^q(I,L^r(\R^n))}\leq C\,\Vert
u_0\Vert _{L^2(\R^n)}\,.
\end{equation}
Let us set now
\begin{equation}\label{eqII.1.6}
Sf(t)=\int^T_0 e^{i(t-s)\Delta_g}\,[f(s,\cdot )]\,ds\,.
\end{equation}
Since $Sf(t)=e^{it\Delta_g}\int^T_0 e^{-is\Delta_g}\,[f(s,\cdot )]\,ds$
we can use Theorem \ref{tII.1.5} to write 
\begin{equation*}\begin{aligned}
\Vert Sf(t)\Vert _{L^q(I,L^r(\R^n))}&\leq C \Big\Vert \int^T_0
e^{-is\Delta_g}\,[f(s,\cdot )]\,ds\Big\Vert _{L^2(\R^n)}\\
&\leq C \int^T_0 \big\Vert e^{-is\Delta_g}\,[f(s,\cdot )]\big\Vert
_{L^2(\R^n)}\,ds\\
&\leq C \int^T_0 \Vert f(s,\cdot )\Vert _{L^2(\R^n)}\,ds=C\,\Vert f\Vert
_{L^1(I,L^2(\R^n))}\,.
\end{aligned}\end{equation*}
Using Lemma \ref{lII.1.4} with $p=1$, $q>2$, $Y=L^r(\R^n)$, $X=L^2(\R^n)$ we
deduce that
$$
\Big\Vert \int^t_0 e^{i(t-s)\,\Delta_g}[F_1(s,\cdot )]\,ds\Big\Vert
_{L^q(I,L^r(\R^n))}\leq C\,\Vert F_1\Vert _{L^1(I,L^2(\R^n))}
$$
where $F_1=-\big(\som^n_{j=1} (D_j\,b_j)+V)\,u$.
Since $\som^n_{j=1}\v D_j\,b_j\v+\v V\v$ is bounded (by condition
(\ref{eqI.0.3})) we have
$$
\Vert F_1\Vert _{L^1(I,L^2(\R^n))}\leq C\, \int^T_0 \Vert u(s,\cdot )\Vert
_{L^2(\R^n)}\,ds\leq C'\,T\,\Vert u_0\Vert _{L^2}\,.
$$
Therefore we have
\begin{equation}\label{eqII.1.7} 
\Big\Vert \int^t_0 e^{i(t-s)\,\Delta_g}\,[F_1(s,\cdot )]\,ds\Big\Vert \leq
C(T)\,\Vert u_0\Vert _{L^2(\R^n)}\,.
\end{equation}
Let us look to the term corresponding to $F_2$ in (\ref{eqII.1.3}),
(\ref{eqII.1.4}). Let us fix $\sigma=\frac 12+\frac 12\,\sigma_0$. Then by
Theorem \ref{tII.1.2} the operator $e^{it\,\Delta_g}$ is continuous from
$L^2(\R^n)$ to $L^2(I,H^{\frac 12}_{-\sigma}(\R^n))$. Its adjoint is defined by
$$
((e^{it\,\Delta_g}\,u_0,f))=(u_0,U^*f)_{L^2(\R^n)}
$$
where $((\kern 3pt ,\kern 1pt))$ denotes the duality between $L^2(I,H^{\frac
12}_{-\sigma})$ and $L^2(I,H^{-\frac 12}_{\sigma})$. It satisfies the estimate
$$
\Vert U^*f\Vert _{L^2(\R^n)}\leq C\,\Vert f\Vert _{L^2(I,H^{-\frac
12}_\sigma (\R^n))}\,.
$$
A straightforward computation shows that
$$
U^*f(x)=\int^T_0 e^{-is\,\Delta_g}\,[f(s,\cdot )]\,ds\,.
$$
Using Theorem \ref{tII.1.5} for $\Delta_g$ we see that the operator $S$
introduced in (\ref{eqII.1.6}) satisfies the estimate
$$
\Vert Sf(t)\Vert _{L^q(I,L^r(\R^n))}\leq C\,\Vert f\Vert _{L^2(I,H^{-\frac
12}_\sigma(\R^n))}\,.
$$
Using Lemma \ref{lII.1.4} with $p=2$, $q>2$, $Y=L^r(\R^n)$, $X=H_\sigma^{-\frac
12}(\R^n)$ we see that
\begin{equation}\label{eqII.1.8}
\Big\Vert \int^t_0 e^{i(t-s)\,\Delta_g}\,[F_2(s,\cdot )]\,ds\Big\Vert
_{L^q(I,L^r(\R^n))}\leq C\,\Vert F_2\Vert _{L^2(I,H^{-\frac 12}_\sigma(\R^n))}
\end{equation}
where $F_2=-2\som^n_{j=1} b_j\,D_j\,u$.
If we set, with $\Delta=\som^n_{j=1}\,\frac{\partial ^2}{\partial x^2_j}$,
\begin{equation}\label{eqII.1.9}
A=\bra x\ket^\sigma\,(I-\Delta)^{-\frac 14} \sum^n_{j=1}
b_j\,D_j(I-\Delta)^{-\frac 14}\,\bra x\ket^\sigma
\end{equation}
then we can write
\begin{equation}\label{eqII.1.10}
\Vert F_2\Vert ^2_{L^2(I,H^{-\frac 12}_\sigma(\R^n))}=4\,\int^T_0 \big\Vert
A\bra x\ket^{-\sigma}\,(I-\Delta)^{\frac 14}\,u(s,\cdot )\big\Vert
^2_{L^2(\R^n)}\,ds\,.
\end{equation}
Let us consider the metric on the cotangent space
$$
G=\frac{dx^2}{\bra x\ket^2}+\frac{d\xi^2}{\bra\xi\ket^2}\,.
$$
It is a H\"ormander's metric and we have
$\bra x\ket^\sigma\in\Op S(\bra x\ket^\sigma,G)$, $(I-\Delta)^{-\frac
14}\in\Op S\big (\bra\xi\ket^{-\frac 12},G\big)$, $b_j\in\Op S(\bra
x\ket^{-2\sigma},G)$, $D_j\in\Op S(\bra\xi\ket,G)$. It follows that the
operator $A$ introduced in (\ref{eqII.1.9}) belongs to $\Op S(1,G)$ and
therefore is $L^2$ continuous. It follows then from (\ref{eqII.1.8}),
(\ref{eqII.1.10}) that
$$
\Big\Vert \int^t_0 e^{i(t-s)\,\Delta_g}\,[F_2(s,\cdot )]\,ds\Big\Vert
_{L^q(I,L^r(\R^n))}\leq C\,\Big[\int^T_0 \Vert u(s,\cdot )\Vert ^2_{H^{\frac
12}_{-\sigma}(\R^n)}\, ds\Big]^{1/2}\,.
$$
Using Theorem \ref{tII.1.2} for $P$ we deduce that
\begin{equation}\label{eqII.1.11}
\Big\Vert \int^t_0 e^{i(t-s)\,\Delta_g}\,[F_2(s,\cdot )]\,ds\Big\Vert
_{L^q(I,L^r(\R^n))}\leq C'\,\Vert u_0\Vert _{L^2(\R^n)}\,.
\end{equation}
Gathering the informations given by (\ref{eqII.1.4}), (\ref{eqII.1.5}), 
(\ref{eqII.1.7}) and (\ref{eqII.1.11}) we obtain the conclusion of Theorem
\ref{tI.0.1}. So we are left with the proof of Theorem \ref{tII.1.5}.


\subsection{Reduction to a small perturbation}\label{ssII.2}

The purpose of this Section is to show that, using the result of \ref{ssII.1}
one can reduce the proof of Theorem \ref{tII.1.5} to the case of a small
perturbation of the flat Laplacian.

Let $\varphi $ be in $C^\infty _0(\R^n)$. We write
$e^{it\,\Delta_g}u_0=u$ and
\begin{equation}\label{eqII.2.1}
u=\varphi  u+(1-\varphi )\,u=v+w\,.
\end{equation}

\noindent {\bf (i) Estimate of $\boldsymbol{v}$}

Since $v=\varphi u$ it follows from (\ref{eqI.0.7}) that $(i\partial
_t+\Delta_g)\,v=[\Delta_g,\varphi ]\,u$. Let $\varphi _1\in C^\infty _0(\R^n)$
be such
$\varphi _1=1$ on the support of $\varphi $ then setting $\Delta=\som^n_{j=1}
\frac{\partial ^2}{\partial x^2_j}$ one can write
\begin{equation}\label{eqII.2.2}
(i\,\partial _t+\Delta_g)\,v=(i\,\partial _t+\Delta+\varphi
_1(\Delta_g-\Delta)\varphi _1)\,v=[\Delta_g,\varphi ]\,u
\end{equation}
and $\varphi _1(\Delta_g-\Delta)$ is a compactly supported perturbation of the
flat Laplacian. Let us set $\tilde P=-\Delta-\varphi
_1(\Delta_g-\Delta)\,\varphi _1$. We have, from (\ref{eqII.2.2})
\begin{equation}\label{eqII.2.3}
v=e^{-it\tilde P}\,\varphi u_0+\int^t_0 e^{-i(t-s)\tilde P}\,[f(s,\cdot )]\,ds
\end{equation}
where $f=[\Delta_g,\varphi ]\,u$.

It follows from Theorem \ref{tII.1.1} that
\begin{equation} \label{eqII.2.4}
\big\Vert e^{-it\tilde P}\varphi \,u_0\big\Vert_{L^q([-T,T],L^r(\R^n))}\leq C\,
\Vert u_0\Vert _{L^2(\R^n)}\,. 
\end{equation}
To estimate the second term in the right-hand side of (\ref{eqII.2.3}) we shall
use Lemma \ref{lII.1.4} with $a=-T$, $b=T$, $Y=L^r(\R^n)$, $p=2$,
$X=H^{-1/2}(\R^n)$. For this one first remark that if $U=e^{-it\tilde P}$ then
Theorem \ref{tII.1.2} shows that $U$ is continuous from $L^2(\R^n)$ to
$L^2([-T,T],H^{1/2}_{\loc}(\R^n))$. Then it is easy to see that
$U^*:L^2([-T,T],H^{-1/2}_c(\R^n))\rightarrow L^2(\R^n)$ is continuous and is
given by $U^*f(x)=\int^T_0 e^{-is\tilde P}[f(s,\cdot )]\,ds$. It follows that
$$
\Big\Vert \int^T_0 e^{-i(t-s)\tilde P}[f(s,\cdot )]\,ds\Big\Vert
_{L^q([-T,T],L^r(\R^n))}=\Vert U\,U^*f\Vert _{L^q([-T,T],L^r(\R^n))}\,.
$$
Then, using again Theorem \ref{tII.1.1} and the above continuity of $U^*$ we get
$$
\Vert U\,U^*f\Vert _{L^q([-T,T],L^r(\R^n))}\leq C\,\Vert U^*f\Vert
_{L^2(\R^n)}\leq C'\,\Vert f\Vert _{L^2([-T,T],H^{-1/2}(\R^n))}
$$
since $f=[\Delta_g,\varphi ]\,u$ has compact support in $x$.

Now we use Lemma \ref{lII.1.4} to deduce that
$$
\Big\Vert \int^t_0 e^{-i(t-s)\tilde P}[f(s,\cdot )]\,ds\Big\Vert
_{L^q([-T,T],L^r(\R^n))}\leq C''\,\Vert f\Vert _{L^2([-T,T],H^{-1/2}(\R^n))}
$$
since $f(s,\cdot )$ has compact support in $x$ and $q>2$.

Moreover since $[\Delta_g,\varphi ]$ is first order we have, using again Theorem
\ref{tII.1.2},
$$
\Vert f\Vert _{L^2([-T,T],H^{-1/2}(\R^n))}\leq C\,\Vert \psi u\Vert
_{L^2([-T,T],H^{1/2}(\R^n))}\leq C'\, \Vert u_0\Vert _{L^2(\R^n)}
$$
where $\psi\in C^\infty _0(\R^n)$, $\psi=1$ on the support of $\varphi $. This
gives the estimate of the second term in the right hand side of (\ref{eqII.2.3})
which, together with (\ref{eqII.2.4}) shows that $v$ satisfies the Strichartz
estimate.

\bigskip

\noindent {\bf (ii) Estimate of $\boldsymbol{w}$}

We shall take the function $\varphi $, introduced above, of the following form.
Let $R>0$ (which will be chosen large enough)
 and $\varphi _0\in C^\infty _0(\R^n)$   such that $\varphi _0(x)=1$ if $\vert
x\vert \leq \frac{3}{2}$, $\supp \varphi _0\subset [-2,2]$. We shall take
$\varphi (x)=\varphi _R(x)=\varphi _0\big (\frac{x}{R}\big)$.

Let $\tilde \varphi _0\in C^\infty _0(\R^n)$ be such that $\tilde \varphi
_0(x)=1$ if $\vert x\vert \leq \frac{1}{2}$, $\supp \tilde \varphi _0\subset
[-1,1]$ and let us set $\tilde \varphi _R(x)=\tilde \varphi _0\big
(\frac{x}{R}\big )$.

Let $w=(1-\varphi _R)\,u$ be the second term in the right hand side of
(\ref{eqII.2.1}). Since $1-\tilde \varphi_R=1$ on the support of $1-\varphi _R$
we have according to (\ref{eqI.0.2})
\begin{equation}\label{eqII.2.5}
(i\,\partial _t+\Delta_g)\,w=\Big (i\partial _t+\Delta+\sum^n_{j,k=1}
\frac\partial {\partial x_j}\Big[(1-\tilde
\varphi _R)\,b_{jk}\,\frac{\partial }{\partial x_k}\Big]\Big)\,
w=-[\Delta_g,\varphi _R]\,u
\end{equation}
where $b_{jk}=g^{jk}-\delta_{jk}$.

Now if we denote by $f$ one of the coefficients $b_{jk}$  we claim
that we have
\begin{equation}\label{eqII.2.6}
\vert \partial ^\alpha_x[(1-\tilde \varphi _R)\,f](x)\vert \leq
\frac{1}{R^{\frac{\sigma_0}2}}\, \frac{C_\alpha}{\langle x\rangle^{\vert
\alpha\vert +1+\frac{\sigma_0}{2}}}\,,\quad \forall x\in\R^n\,.
\end{equation}
Using (\ref{eqI.0.1}) and denoting by $A$ the left hand side of (\ref{eqII.2.6})
we see that
\begin{equation*}
\begin{aligned}
A&\leq \Big(1-\tilde \varphi \Big(\frac{x}{R}\Big)\Big)\vert \partial
^\alpha_x\,f(x)\vert +\sum_{0<\beta\leq \alpha} \begin{pmatrix}\alpha\\
\beta\end{pmatrix}
\frac{1}{R^{\vert \beta\vert }}\,\Big\vert (\partial ^\beta_x\,\tilde \varphi
)\Big(\frac{x}{R}\Big)\Big\vert \vert \partial ^{\alpha-\beta}_x\,f(x)\vert \\
A&\leq \Big (1-\tilde \varphi \Big (\frac xR\Big )\Big )\frac{C_\alpha}{\langle
x\rangle^{\vert \alpha\vert +1+\sigma_0}}+\sum_{0<\beta\leq \alpha}
\frac{C'_{\alpha\beta}}{R^{\vert \beta\vert }} \,\Big\vert (\partial
^\beta_x\,\tilde \varphi )\Big (\frac xR\Big )\Big\vert \frac{1}{\langle
x\rangle^{\vert \alpha\vert -\vert \beta\vert +1+\sigma_0}}\,.
\end{aligned}
\end{equation*}
Now, on the support of $1-\tilde \varphi (\frac xR)$ we have $\langle
x\rangle>\vert x\vert \geq \frac 32\,R$ so
$$
(1)\leq \frac{C'_\alpha}{R^{\frac{\sigma_0}{2}}}\, \frac{1}{\langle
x\rangle^{\vert \alpha\vert +1+\frac{\sigma_0}{2}}}\,.
$$
On the support of $\partial ^\beta\tilde \varphi (\frac xR)$, with $\beta\not
=0$, we have $\frac 12\,R\leq \vert x\vert \leq R$ so $\langle x\rangle\leq
\sqrt 2\,R$ if $R>1$. Therefore
$$
(2)\leq \frac{1}{R^{\frac{\sigma_0}{2}}}\sum_{0<\beta\leq \alpha}
C''_{\alpha\beta}\frac{1}{\langle x\rangle^{\vert \beta\vert
-\frac{\sigma_0}{2}+\vert \alpha\vert -\vert \beta\vert +1+\sigma_0}}\leq
\frac{1}{R^{\frac{\sigma_0}{2}}}\, \frac{C''_\alpha}{\langle x\rangle^{\vert
\alpha\vert +1+\frac{\sigma_0}{2}}}\,.
$$
It follows from (\ref{eqII.2.6}) that we can work in the rest of the paper with
a non negative self adjoint operator $P$ such that
\begin{equation}\label{eqII.2.7}\left\{
\begin{array}{l}
P=-\Delta+\varepsilon\,Q\,,\textrm{ where } Q=\som_{\vert \beta\vert \leq 2}
a^\varepsilon_\beta\, D^\beta\,, \\
\varepsilon\textrm{ is a small constant and } \vert D^\alpha_x\,
a^\varepsilon_\beta(x)\vert \leq \frac{C_\alpha}{\langle x\rangle^{\vert
\alpha\vert +1+\frac{\sigma_0}{2}}}\,,\enskip \forall \alpha\in\N^n\,,\\
\textrm{uniformly for }x\in\R^n\textrm{ with } C_\alpha\textrm{ independent of }
\varepsilon\,.
\end{array}\right.
\end{equation}
Since the estimates on the coefficients are uniform in $\varepsilon$ we shall  write
$a_\beta$ instead of $a^\varepsilon_\beta$. The principal symbol $p$ of $P$ will be
written as
$$
p(x,\xi)=\vert \xi\vert ^2+\varepsilon\,q(x,\xi)\,, \quad q(x,\xi)=\sum^n_{j,k=1}
b_{jk}(x)\,\xi_j\,\xi_k
$$
and we shall take $\varepsilon$ so small that
$$
\frac{9}{10}\, \vert \xi\vert ^2\leq p(x,\xi)\leq \frac{11}{10}\, \vert \xi\vert ^2\,.
$$
Finally without loss of generality we shall take $\sigma_0$ instead of $\frac{\sigma_0}2$
in (\ref{eqII.2.7}).

We assume that $P$ satisfies the condition (\ref{eqI.0.3}) and (\ref{eqI.0.6}).
Let
$T>0$.

\begin{theorem}\sl \label{tII.2.1}
Let $(q,r)$ be  such  $q\geq 2$ and $\frac{2} {q}=\frac {n} {2}-\frac {n} {r}$.
If $\varepsilon$ is small enough then there exists $C>0$ such that
$$
\Vert e^{-itP}\,v_0\Vert _{L^q([-T,T],L^r(\R^n))}\leq C\,\Vert v_0\Vert _{L^2(\R^n)}
$$
\end{theorem}
for all $v_0\in L^2(\R^n)$.

Let us assume that we have proved this result. Then we can applied it to the operator
occuring in (\ref{eqII.2.5}) with $R$ large enough. We have,
$$
w=e^{-itP}(1-\varphi _R)\,u_0+\int^t_0 e^{-i(t-s)P}[f_R(s,\cdot )]\,ds\,.
$$
It follows from Theorem \ref{tII.2.1} that
\begin{equation}\label{eqII.2.8}
\Vert e^{-itP}\,(1-\varphi _R)\,u_0\Vert _{L^q([-T,T],L^r(\R^n))}\leq C\,\Vert
u_0\Vert _{L^2(\R^n)}
\end{equation}
and the same argument as used in the estimate of $v$, namely the use of Theorem
\ref{tII.2.1},and Lemma \ref{lII.1.4} shows that
\begin{equation}\label{eqII.2.9}
\Vert \int^t_0 e^{-i(t-s)P}\,[f_R(s,\cdot )]\,ds\Vert _{L^q([-T,T],L^r(\R^n))}\leq
C(R)\Vert u_0\Vert _{L^2(\R^n)}\,.
\end{equation}
Then using (\ref{eqII.2.8}) we see that the second term $w$ in the right hand
side of (\ref{eqII.2.1}) satisfies the Strichartz estimate which completes the
proof of Theorem \ref{tII.1.5}.

Our goal now is to prove Theorem \ref{tII.2.1}. The first step is to make a carefull
study of the flow.
 
\section{Study of the flow}\label{sIII}
\subsection{Preliminaries}\label{ssIII.1}

Let $p(x,\xi)=\vert \xi\vert ^2+\varepsilon\,q(x,\xi)$, $q(x,\xi)=\som^n_{j,k=1}
b_{jk}(x)\,\xi_j\,\xi_k$ where,
\begin{equation}\label{eqIII.1.1}\left\{
\begin{array}{l}
\textrm{there exists }\sigma_0>0 \textrm { such that for every } \ell \in\N 
\textrm{ one
can find } A_\ell >0\textrm{ such that}\\
\som_{\vert \alpha\vert =\ell }\, \som^n_{j,k=1} \vert \partial ^\alpha_x\,b_{jk}(x)\vert
\leq \frac{A_\ell}{\bra x\ket^{1+\ell +\sigma_0}}\textrm{ for all } x\textrm{ in
}\R^n\,.
\end{array}\right.
\end{equation}

We introduce the equations of the bicharacteristic flow issued from a point $(x,\xi)$ in
$T^*\R^n\setminus\{0\}$. They are given for $j=1,\ldots ,n$, by
\begin{equation}\label{eqIII.1.2}\left\{
\begin{array}{ll}
\dot x_j(t)=\frac{\partial p}{\partial \xi_j}\, (x(t),\xi(t))\,, &x_j(0)=x_j\,,\\
\dot \xi_j(t)=-\frac{\partial p}{\partial x_j}\, (x(t),\xi(t))\,,&\xi_j(0)=\xi_j\,,
\end{array}\right.
\end{equation}
and we denote by $(x(t,x,\xi),\xi(t,x,\xi))$ the solution of (\ref{eqIII.1.2}) whenever
it exists (or $(x(t),\xi(t))$ for short if no confusion is possible).

Assuming $\varepsilon$ so small that $\varepsilon\,A_0\leq \frac 1{10}$ we see that
$\frac 9{10}\, \vert \xi\vert ^2\leq p(x,\xi)\leq \frac{11}{10}\,\vert \xi\vert ^2$. It
follows that
$$
\frac 9{10} \,\vert \xi(t,x,\xi)\vert ^2\leq p(x(t),\xi(t))=p(x,\xi)\leq \frac{11}{10}\,
\vert \xi\vert ^2\,,
$$
so that
\begin{equation}\label{eqIII.1.3}
\vert \xi(t,x,\xi)\vert \leq 2\,\vert \xi\vert \,.
\end{equation}
Using the first equation of (\ref{eqIII.1.2}) we see then, that the solution of
(\ref{eqIII.1.2}) exists for all $t$ in $\R$ and is a $C^\infty $ function with respect to
$(x,\xi)$. Moreover we have the following lemma.
\begin{lemma}\sl \label{lIII.1.1}
For all $t$ in $\R$ we have
$$
x(t,x,\xi)\cdot \xi(t,x,\xi)=x\cdot \xi+2t\,p(x,\xi)+f(t,x,\xi)
$$
where
$$
\vert f(t,x,\xi)\vert \leq 4\varepsilon\,A_1\, \vert \xi\vert ^2 \Big\vert
\int^t_0
\frac{ds}{\bra x(s)\ket^{1+\sigma_0}}\Big\vert\leq 4\varepsilon\,A_1\,\vert t\vert \,\vert
\xi\vert ^2\,.
$$
\end{lemma}

\noindent {\bf Proof }
We have by (\ref{eqIII.1.2}) 
$$
\frac d{dt} [x(t)\cdot \xi(t)]=\xi(t)\cdot \frac{\partial p}{\partial \xi}\,
(x(t),\xi(t))-\varepsilon\,x(t)\cdot \frac{\partial q}{\partial x}\,
(x(t),\xi(t))\,.
$$
Using Euler's identity we obtain
$$
\xi(t)\cdot \frac{\partial p}{\partial \xi}\,
(x(t),\xi(t))=2p(x(t),\xi(t))=2p(x,\xi)\,.
$$
We set $f(t,x,\xi)=-\varepsilon\int^t_0 x(s)\cdot \frac{\partial
q}{\partial x}\,(x(s),\xi(s))\,ds$. Now since
$$
\Big\vert \frac{\partial q}{\partial x}\,(x(s),\xi(s))\Big\vert \leq
\frac{A_1}{\bra x(s)\ket^{2+\sigma_0}}\,\vert \xi(s)\vert ^2\,,
$$
it follows from (\ref{eqIII.1.3}) that
$$
\vert f(t,x,\xi)\vert \leq 4\varepsilon\,A_1\,\vert \xi\vert ^2 \Big\vert
\int^t_0 \frac{ds}{\bra x(s)\ket^{1+\sigma_0}}\Big\vert \leq
4\varepsilon\,A_1\,\vert t\vert \,\vert \xi\vert ^2\,.
$$
\cqfd

We shall use later on the result given by the following lemma.

For $t\in\R$ and $(x,\xi)\in T^*\R^n\setminus\{0\}$ let us set,
\begin{equation}\label{eqIII.1.4}
\rho(t,x,\xi)=(x(t,x,\xi),\xi(t,x,\xi))\,.
\end{equation}
\begin{lemma}\sl \label{lIII.1.2}
We have the following identities for $j,k=1,\ldots n,$.
\begin{equation*}
\begin{aligned}
\frac{\partial x_j}{\partial x_k}\, (t,x,\xi)&=\frac{\partial \xi_k}{\partial
\xi_j}\, (-t,\rho(t,x,\xi))\\
\frac{\partial x_j}{\partial \xi_k}\, (t,x,\xi)&=-\frac{\partial x_k}{\partial
\xi_j}\, (-t,\rho(t,x,\xi))\\
\frac{\partial \xi_j}{\partial x_k}\, (t,x,\xi)&=-\frac{\partial \xi_k}{\partial
x_j}\, (-t,\rho(t,x,\xi))\\
\frac{\partial \xi_j}{\partial \xi_k}\, (t,x,\xi)&=\frac{\partial
x_k}{\partial x_j}\, (-t,\rho(t,x,\xi))\,.
\end{aligned}
\end{equation*}
\end{lemma}

\noindent {\bf Proof } For $j=1,\ldots ,n$ we have
\begin{equation*}\left\{
\begin{array}{ll}
x_j(-t;\rho(t;x,\xi))&=x_j\\
\xi_j(-t;\rho(t,x,\xi))&=\xi_j
\end{array}\right.
\end{equation*}
Differentiating both sides with respect to $x_k$ and $\xi_k$ we obtain
\begin{equation*}
\begin{aligned}
\sum^n_{\ell =1} \frac{\partial x_j}{\partial x_\ell }\,
(-t;\rho(t;x,\xi))\,\frac{\partial x_\ell }{\partial
x_k}\,(t;x,\xi)&+\sum^n_{\ell =1} \frac{\partial x_j}{\partial \xi_\ell
}\,(-t;\rho(t;x,\xi))\,\frac{\partial \xi_\ell }{\partial
x_k}\,(t;x,\xi)=\delta_{jk}\\
\sum^n_{\ell =1} \frac{\partial x_j}{\partial x_\ell }\,
(-t;\rho(t;x,\xi))\,\frac{\partial x_\ell }{\partial
\xi_k}\,(t;x,\xi)&+\sum^n_{\ell =1} \frac{\partial x_j}{\partial \xi_\ell
}\,(-t;\rho(t;x,\xi))\,\frac{\partial \xi_\ell }{\partial
\xi_k}\,(t;x,\xi)=0\\
\sum^n_{\ell =1} \frac{\partial \xi_j}{\partial x_\ell }\,
(-t;\rho(t;x,\xi))\, \frac{\partial x_\ell }{\partial
x_k}\,(t;x,\xi)&+\sum^n_{\ell =1} \frac{\partial \xi_j}{\partial \xi_\ell
}\,(-t;\rho(t;x,\xi))\,\frac{\partial \xi_\ell }{\partial
x_k}\,(t;x,\xi)=0\\
\sum^n_{\ell =1} \frac{\partial \xi_j}{\partial x_\ell }\,
(-t;\rho(t;x,\xi))\, \frac{\partial x_\ell }{\partial
\xi_k}\,(t;x,\xi)&+\sum^n_{\ell =1} \frac{\partial \xi_j}{\partial \xi_\ell
}\,(-t;\rho(t;x,\xi))\,\frac{\partial \xi_\ell }{\partial
\xi_k}\,(t;x,\xi)=\delta_{jk}
\end{aligned}
\end{equation*}
where $\delta_{jk}$ is the Kronecker symbol.

If we set
$$
M(t;\rho)=\begin{pmatrix} \big(\frac{\partial x_j}{\partial
x_k}\big)(t;\rho)&\big(\frac{\partial x_j}{\partial \xi_k}\big)(t;\rho)\\
\big(\frac{\partial \xi_j}{\partial x_k}\big)(t;\rho)&\big(\frac{\partial
\xi_j}{\partial\xi_k}\big)(t;\rho)
\end{pmatrix}
$$
then the above relations can be written
\begin{equation}\label{eqIII.1.5}
M(-t;\rho(t;x,\xi))\cdot M(t;x,\xi)=\begin{pmatrix}I_n&0\\0&I_n\end{pmatrix}
\end{equation}
where $I_n$ denotes the $n\times n$ identity matrix.

Let us introduce for $s\in\R$ the following matrix.
\begin{equation}\label{eqIII.1.6}
A(s;\rho)=\begin{pmatrix}{}^t\big(\frac{\partial \xi_j}{\partial
\xi_k}\big)\,(s;\rho) & -{}^t\big(\frac{\partial x_j}{\partial
\xi_k}\big)\,(s;\rho)\\-{}^t\big(\frac{\partial \xi_j}{\partial
x_k}\big)\,(s;\rho) & {}^t\big(\frac{\partial x_j}{\partial
x_k}\big)\,(s;\rho)
\end{pmatrix}.
\end{equation}
We claim that for $s\in \R$ and $\rho\in T^*\R^n$ we have
\begin{equation}\label{eqIII.1.7}
A(s;\rho)\,M(s;\rho)=I_{2n}
\end{equation}
where $I_{2n}$ is the $2n\times 2n$ identity matrix.

Indeed let us set $A(s;\rho)\cdot M(s;\rho)=(C_{\alpha\beta})_{1\leq
\alpha,\beta\leq 2n}$. We have for $j,k=1\ldots n$,
\begin{equation}\label{eqIII.1.8}\left\{
\begin{array}{ll}
C_{j,k}&=\som^n_{\ell =1} \big (\frac{\partial \xi_\ell }{\partial \xi_j}
\,\frac{\partial x_\ell }{\partial x_k}-\frac{\partial x_\ell }{\partial
\xi_j}\, \frac{\partial \xi_\ell }{\partial x_k}\big )(s;\rho)\\
C_{j,k+n}&=\som^n_{\ell =1} \big (\frac{\partial \xi_\ell }{\partial \xi_j}
\,\frac{\partial x_\ell }{\partial \xi_k}-\frac{\partial x_\ell }{\partial
\xi_j}\, \frac{\partial \xi_\ell }{\partial \xi_k}\big )(s;\rho)\\
C_{j+n,k}&=\som^n_{\ell =1} \big (\frac{\partial x_\ell }{\partial x_j}
\,\frac{\partial \xi_\ell }{\partial x_k}-\frac{\partial x_\ell }{\partial
x_k}\, \frac{\partial \xi_\ell }{\partial x_j}\big )(s;\rho)\\
C_{j+n,k+n}&=\som^n_{\ell =1} \big (\frac{\partial x_\ell }{\partial x_j}
\,\frac{\partial \xi_\ell }{\partial \xi_k}-\frac{\partial \xi_\ell }{\partial
x_j}\, \frac{\partial x_\ell }{\partial \xi_k}\big )(s;\rho)\,.
\end{array}\right.
\end{equation}
Let us remark that $C_{j+n,k+n}=C_{k,j}$.

Now we recall that for every $s\in\R$ the map $(x,\xi)\mapsto \rho(s;x,\xi)$ is
symplectic which means that
\begin{equation}\label{eqIII.1.9}
\sum^n_{\ell =1} d(\xi_\ell (s;x,\xi))\wedge d(x_\ell (s;x,\xi))=\sum^n_{j=1}
d\xi_j\wedge dx_j\,.
\end{equation}
Writting $u(s)=u(s;x,\xi)$ for short we have
\begin{equation*}
\begin{split}
(1)=\sum^n_{\ell =1} d(\xi_\ell (s))\wedge d(x_\ell (s))=\sum^n_{\ell =1}
\bigg(\sum^n_{j=1} \Big (\frac{\partial \xi_\ell }{\partial x_j}\,
(s)\,dx_j+\frac{\partial \xi_\ell }{\partial \xi_j}\,(s)\,d\xi_j\Big)\\
\wedge\sum^n_{k=1}\Big (\frac{\partial x_\ell }{\partial
x_k}\,(s)\,dx_k+\frac{\partial x_\ell }{\partial \xi_k}\,(s)\,d\xi_k\Big )\bigg
)
\end{split}
\end{equation*}
It follows that
\begin{equation*}
\begin{aligned}
(1)=&\sum_{j\leq k} \bigg(\sum^n_{\ell =1} \Big (\frac{\partial \xi_\ell
}{\partial x_j}\,(s)\,\frac{\partial x_\ell }{\partial x_k}\,(s)-\frac{\partial
\xi_\ell }{\partial x_k}\,(s)\,\frac{\partial x_\ell }{\partial x_j}\,(s)\Big
)\bigg)\,dx_j\wedge dx_k+\\
&\sum^n_{j,k=1} \bigg(\sum^n_{\ell =1} \Big (\frac{\partial \xi_\ell
}{\partial \xi_j}\,(s)\,\frac{\partial x_\ell }{\partial
x_k}\,(s)-\frac{\partial
\xi_\ell }{\partial x_k}\,(s)\,\frac{\partial x_\ell }{\partial \xi_j}\,(s)\Big
)\bigg)\,d\xi_j\wedge dx_k+\\
&\sum_{j\leq k} \bigg(\sum^n_{\ell =1} \Big (\frac{\partial \xi_\ell
}{\partial \xi_j}\,(s)\,\frac{\partial x_\ell }{\partial
\xi_k}\,(s)-\frac{\partial
\xi_\ell }{\partial \xi_k}\,(s)\,\frac{\partial x_\ell }{\partial
\xi_j}\,(s)\Big )\bigg)\,d\xi_j\wedge d\xi_k\,.
\end{aligned}
\end{equation*}
Using (\ref{eqIII.1.8}) and (\ref{eqIII.1.9}) we see easily that
$$
C_{j,k}=\delta_{jk}\,,\quad  C_{j,k+n}=C_{j+n,k}=0\,,\quad 
C_{j+n,k+n}=C_{k,j}=\delta_{jk}\,.
$$
This proves (\ref{eqIII.1.7}).

It follows from (\ref{eqIII.1.5}) and (\ref{eqIII.1.7}) that
\begin{equation}\label{eqIII.1.10}
M(t;x,\xi)=A(-t;\rho(t;x,\xi))
\end{equation}
which by (\ref{eqIII.1.6}) proves the Lemma \ref{lIII.1.2}.\cqfd

\subsection{The flow for short time}\label{ssIII.2}

Here is a description of the flow for short time.
\begin{proposition} \sl \label{pIII.2.1}
Let us set
\begin{equation*}\left\{
\begin{array}{l}
r(t,x,\xi)=x(t,x,\xi)-(x+2t\xi)\\
\zeta(t,x,\xi)=\xi(t,x,\xi)-\xi\,.
\end{array}\right.
\end{equation*}
Let $T>0$. Then for all $A,B$ in $\N^n$ one can find $C_{A,B}>0$ such that
\begin{equation*}\left\{
\begin{array}{ll}
(i)\quad \vert \partial ^A_x\,\partial
^B_\xi\,Z(t,x,\xi)\vert &\leq C_{A,B}\, \varepsilon\,\vert t\vert \\
(ii)\quad \vert \partial_t\,\partial
^A_x\,\partial ^B_\xi\,Z(t,x,\xi)\vert &\leq C_{A,B}\,
\varepsilon\, 
\end{array}\right.
\end{equation*}
if $Z=r$ or $\zeta$, for all $\vert t\vert \leq T$ and all $(x,\xi)\in T^*\R^n$
with $\vert \xi\vert \leq 3$.
\end{proposition}

\noindent {\bf Proof } See Appendix, Paragraph \ref{ssVIII.2}.

We introduce now the following definition which distinguish microlocally the
points in the cotangent bundle.
\begin{definition}\sl \label{dIII.2.2}
Let
\begin{equation*}
\begin{array}{l}
\CS_+=\big\{(x,\xi)\in T^*\R^n\setminus \{0\} : x\cdot \xi\geq -\frac 1 4\, \bra
x\ket \,\vert \xi\vert \big\}\\
\CS_-=\big\{(x,\xi)\in T^*\R^n\setminus \{0\} : x\cdot \xi\leq \frac 1 4\, \bra
x\ket \,\vert \xi\vert \big\}\,.
\end{array}
\end{equation*}
Then $\CS_+$ (resp. $\CS_-$) is called the set of outgoing points for $t\geq 0$
(resp. $t\leq 0$).
\end{definition}

Of course the constant $\frac 14$ in the above definition is unimportant and
could be replaced by any fixed small constant. The reason for this definition
is the following. If $(x,\xi)\in\CS_+$ then, for $t\geq 0$
$$
1+\vert x+2t\xi\vert ^2\geq \frac12\, (\bra x\ket^2+t^2\,\vert \xi\vert ^2)\,.
$$
Since $x+2t\xi$ will be an approximation of $x(t;x,\xi)$, then $\CS_+$, will be
the set of points $(x,\xi)$ for which the projection of the bicharacteristic
goes to $+\infty $  when $t\rightarrow +\infty $ in staying away from the
origin.

\subsection{The forward flow from points in $\CS_+$ and
backward from $\CS_-$}\label{ssIII.3}

Our goal is to obtain for these points a nice global representation of the flow
together with precise estimates of its derivatives with respect to $x$ and
$\xi$.
\begin{proposition} \sl \label{pIII.3.1}
There exists $\varepsilon_0>0$ depending on the constants $A_0,A_1$ in
(\ref{eqIII.1.1}) such that for $\varepsilon$ in $]0,\varepsilon_0[$ the
solution of (\ref{eqIII.1.2}) with $(x,\xi)$ in $\CS_+$ (resp. $\CS_-$) and
$\frac 12\leq \vert \xi\vert \leq 2$ can be written for all $t\geq 0$ (resp.
$t\leq 0$)
\begin{equation}\label{III.3.1}\left\{
\begin{array}{l}
x(t;x,\xi)=x+2t\,\xi(t;x,\xi)+z(t;x,\xi)\\
\xi(t;x,\xi)=\xi+\zeta(t;x,\xi)
\end{array}\right.
\end{equation}
with
\begin{equation}\label{eqIII.3.2}
\vert z_j(t;x,\xi)\vert \leq \frac{2\cdot 10^2}{\sigma_0}\,\varepsilon\,\max
(A_0,A_1)\,,\quad  
\vert \zeta_j(t;x,\xi)\vert \leq \frac{2\cdot 10^2}{\sigma_0}\,\varepsilon\,\max
(A_0,A_1)\,,
\end{equation}
where $A_0A_1$ are the constants arising in (\ref{eqIII.1.1}) and $j=1,\ldots
,n$.

Moreover for all $t\geq 0$ (resp. $t\leq 0$) we have
\begin{equation}\label{eqIII.3.3}
\frac13\leq \frac{1+\vert x(t;x,\xi)\vert ^2}{1+\vert x\vert ^2+t^2}\leq 40\,.
\end{equation}
\end{proposition}

\noindent {\bf Proof } Let $I=\big\{T>0:\vert z_j(t)\vert \leq \frac{2\cdot
10^2}{\sigma_0}\,\varepsilon\,\max (A_0,A_1),\enskip \v\zeta_j(t)\vert \leq
\frac{2\cdot 10^2}{\sigma_0}\,\varepsilon\,\max (A_0,A_1)$ for $j=1,\ldots ,n$
and all $t\in[0,T]\big\}$. Then $I$ is an interval which is non empty by the
local Cauchy-Lipschitz Theorem. Let $T^*=\sup I$. If $T^*=+\infty $ we are
done. Otherwise let $T<T^*$. Since $\frac12\leq \vert \xi\vert \leq 2$ we have
for
$t\in[0,T]$, $\frac13\leq \vert \xi(t)\vert \leq 3$ if $\varepsilon\,\max
(A_0,A_1)$ is small enough. Indeed we have
\begin{equation*}
\begin{aligned}
(1-\varepsilon\,A_0)\vert \xi\vert ^2&\leq p(x,\xi)=p(x(t),\xi(t))\leq
(1+\varepsilon\,A_0)\vert \xi(t)\vert ^2\\
(1-\varepsilon\,A_0)\vert \xi(t)\vert ^2&\leq p(x(t),\xi(t))=p(x,\xi)\leq
(1+\varepsilon\,A_0)\vert \xi\vert ^2\,.
\end{aligned}
\end{equation*}
Now, for $t$ in $[0,T]$ we have
\begin{equation*}
\begin{split}
1+\vert x(t)\vert ^2=\bra x\ket^2+4t^2\,\vert \xi\vert ^2+4t^2\,\vert
\zeta(t)\vert ^2+\vert z(t)\vert ^2 
+\underbrace{4t\, x\cdot \xi}_{(1)}+\underbrace{4t\,x\cdot \zeta(t)}_{(2)}\\
+\underbrace{2x\cdot z(t)}_{(3)}+\underbrace{8t^2\,\xi\cdot \zeta(t)}_{(4)}
+\underbrace{4t\,\xi\cdot z(t)}_{(5)}+\underbrace{4t\,\zeta(t)\cdot
z(t)}_{(6)}\,.
\end{split}
\end{equation*}
Since $(x,\xi)\in\CS_+$ we have for $t\geq 0$, $(1)\geq -\frac12\,(\bra
x\ket^2+t^2\,\vert \xi\vert ^2)$. Now, by the definition of $I$ we have on
$[0,T]$ if $\varepsilon\, \max (A_0,A_1)$ is small enough. 
\begin{equation*}
\begin{aligned}
\vert (2)\vert &\leq C_1(n)\,t\,\vert x\vert \,\varepsilon\,\max (A_0,A_1)\leq
10^{-2}\,(\vert x\vert ^2+t^2)\\
\vert (3)\vert &\leq C_2(n)\,\vert x\vert \,\varepsilon\,\max (A_0,A_1)\leq
10^{-2}\,\bra x\ket^2\\
\vert (4)\vert &\leq C_3(n)\,t^2 \,\varepsilon\,\max (A_0,A_1)\leq
10^{-2}\,t^2\\
\vert (5)\vert &\leq C_4(n)\,t \,\varepsilon\,\max (A_0,A_1)\leq
10^{-2}\,(1+t^2)\\
\vert (6)\vert &\leq C_5(n)\,t (\varepsilon\,\max (A_0,A_1))^2\leq
10^{-2}\,(1+t^2)\,.
\end{aligned}
\end{equation*}
It follows that
$$
\bra x(t)\ket^2\geq \frac 12 (\bra x\ket^2+t^2)-4\cdot 10^{-2} (\bra
x\ket^2+t^2)\geq \frac 13 \,(\bra x\ket^2+t^2)\,.
$$
The same computation shows that $\bra x(t)\ket^2\leq 40 (\bra x\ket^2+t^2)$. It
follows that on $[0,T]$ we have
\begin{equation}\label{eqIII.3.4}
\frac 1{\sqrt 6}\,(1+t)\leq \frac 1{\sqrt 3}\, (\bra x\ket^2+t^2)^{1/2}\leq
\bra x(t)\ket\leq 7 (\bra x\ket^2+t^2)^{1/2}\,.
\end{equation}
Now it follows from (\ref{eqIII.1.2}) that $(z(t),\zeta(t))$ satisfy the
equations
\begin{equation}\label{III.3.5}\left\{
\begin{array}{l}
\dot z_j(t)=-\varepsilon\,\frac{\partial q}{\partial
\xi_j}\,(x(t),\xi(t))+2t\,\varepsilon\,\frac{\partial q}{\partial
x_j}\,(x(t),\xi(t))\\
\dot\zeta_j(t)=-\varepsilon\,\frac{\partial q}{\partial x_j}\,(x(t),\xi(t))
\end{array}\right.
\end{equation}
with $z_j(0)=\zeta_j(0)=0$.

We deduce from (\ref{eqIII.1.1}), (\ref{eqIII.3.4}) and the bounds $\frac
13\leq \vert \xi(t)\vert \leq 3$ that
\begin{equation*}
\begin{array}{ll}
\Big\vert \frac{\partial q}{\partial \xi_j}\,(x(t),\xi(t))\Big\vert \leq
\frac{3\,A_0}{\bra x(t)\ket^{1+\sigma_0}}\leq \frac{3(\sqrt
6)^{1+\sigma_0}A_0}{(1+t)^{1+\sigma_0}}&\leq
\frac{12\,A_0}{(1+t)^{1+\sigma_0}}\,,\\
\Big\vert \frac{\partial q}{\partial x_j}\,(x(t),\xi(t))\Big\vert \leq
\frac{9\,A_1}{\bra x(t)\ket^{2+\sigma_0}}\leq \frac{9(\sqrt
6)^{1+\sigma_0}\,\sqrt 3\,A_1}{(1+t)^{1+\sigma_0}\bra t\ket}&\leq
\frac{60\,A_1}{(1+t)^{1+\sigma_0}\bra t\ket}\,,
\end{array}
\end{equation*}
(since we may assume that $1+\sigma_0<\frac 32$ and  $(\sqrt
6)^{1+\sigma_0}<4$).  It follows from (\ref{III.3.5}) that
$$
\vert \dot z_j(t)\vert \leq \frac{132\,\varepsilon}{(1+t)^{1+\sigma_0}}\, \max
(A_0,A_1)\,,\quad \vert \dot \zeta_j(t)\vert \leq
\frac{60\,\varepsilon}{(1+t)^{1+\sigma_0}}\, \max (A_0,A_1)\,.
$$
Therefore we have on $[0,T]$
$$
\vert z_j(t)\vert \leq \frac{132}{\sigma_0}\,\varepsilon\, \max
(A_0,A_1)\,,\quad \vert \zeta_j(t)\vert \leq
\frac{60}{\sigma_0}\,\varepsilon\, \max (A_0,A_1)\,.
$$
Since $z(t)$ and $\zeta(t)$ exist for all $t\geq 0$ and are smooth we still
have the above estimates on $[0,T^*]$. By continuity it will exist $\eta>0$
such that $\vert z_j(t)\vert \leq \frac{2\cdot
10^2}{\sigma_0}\,\varepsilon\,\max (A_0,A_1)$ and $\vert \zeta_j(t)\vert \leq
\frac{2\cdot 10^2}{\sigma_0}\,\varepsilon\,\max (A_0,A_1)$ on $[0,T^*+\eta]$.
This contradicts the maximality of $T^*$ and proves that $T^*=+\infty $.

\cqfd

Now we estimate the derivatives of the flow with respect to $(x,\xi)$. 
\begin{proposition} \sl \label{pIII.3.2}
With the notations of Proposition \ref{pIII.3.1}, for every integer $k$ one can
find a positive constant $M_k$ such that for all $(A,B)\in\N^n\times\N^n$ such
that $\vert A\vert +\vert B\vert \leq k$, all $t\geq 0$ (resp. $t\leq 0$) and
$(x,\xi)$ in $\CS_+$ (resp. $\CS_-$) we have,
\begin{equation*}\left\{
\begin{array}{l}
\big\vert \partial ^A_x\,\partial ^B_\xi\,z(t,x,\xi)\big\vert \leq
\frac{\varepsilon\,M_k}{\bra x\ket^{\vert A\vert +\sigma_0}}\,,\\
\big\vert \partial ^A_x\,\partial ^B_\xi\,\zeta(t,x,\xi)\big\vert \leq
\frac{\varepsilon\,M_k}{\bra x\ket^{1+\vert A\vert +\sigma_0}}\,.
\end{array}\right.
\end{equation*}
\end{proposition}
 
\noindent {\bf Proof } See Appendix VIII.3.
\begin{corollary} \sl \label{cIII.3.3}
Keeping the notations of Propositon \ref{pIII.3.1} we have, for all $t\geq 0$
(resp. $t\leq 0$) and all $(x,\xi)\in\CS_+$ (resp. $\CS_-$)
$$
\frac{\partial x_j}{\partial
\xi_k}\,(t,x,\xi)=2t\,\delta_{jk}+\CO(\varepsilon\bra t\ket)\,,\quad
\frac{\partial x_j}{\partial x_k}\,(t,x,\xi)=\delta_{jk}+\CO(\varepsilon\bra
t\ket )
$$
$$
\frac{\partial \xi_j}{\partial
\xi_k}\,(t,x,\xi)=\delta_{jk}+\CO(\varepsilon)\,,\enskip  \frac{\partial 
\xi_j}{\partial x_k}\,(t,x,\xi)=\CO(\varepsilon)\,,\quad j,k=1,\ldots ,n\,,
$$
where $\delta_{jk}$ is the Kronecker symbol and $\CO(\varepsilon)$ means
``bounded by $C\,\varepsilon$ where $C$ is independent of $(x,\xi)$''. In
particular we have
\begin{equation}\label{eqIII.3.6}
\frac{\partial \xi_j}{\partial \xi_k}\,(t,x,\xi)-i\,\frac{\partial
x_j}{\partial \xi_k}\,(t,x,\xi)=(1-2it)\,\delta_{jk}+\CO(\varepsilon\bra
t\ket)\,,\quad j,k=1,\ldots ,n~.
\end{equation}
\end{corollary}

\subsection{Precisions on the flow in the general case}\label{ssIII.4}

The results obtained above allow us to give a rough form of the flow through
any point in $T^*\R^n\setminus\{0\}$  for $t\in\R$.
\begin{proposition} \sl \label{pIII.4.1}
Let $(x,\xi)\in T^*\R^n\setminus\{0\}$ with $\vert \xi\vert \leq 2$. Then,
\begin{itemize}
\item[(i)] the function $s\mapsto \bra x(s,x,\xi)\ket^{-(1+\sigma_0)}$ belongs
to
$L^1(\R)$,
\item[(ii)] for $t\in\R$ we have,
\begin{equation*}\left\{
\begin{array}{l}
x(t,x,\xi)=x+2t\,\xi+r(t,x,\xi)\,,\\
\xi(t,x,\xi)=\xi+\zeta(t,x,\xi)\,,
\end{array}\right.
\end{equation*}
where $\vert r(t,x,\xi)\vert \leq C\,\varepsilon\,\bra t\ket$, $\vert
\zeta(t,x,\xi)\vert \leq C\,\varepsilon$ with $C$ independent of $(x,\xi)$.
\end{itemize}
\end{proposition}

Before going into the proof let us note that in general we do not have good
estimates on the derivatives of $r$ with respect to $(x,\xi)$ (in the spirit
of those given in Proposition \ref{III.3.1} for instance). In particular we do
not have a good control of $\frac{\partial x_j}{\partial \xi_k}\,(t,x,\xi)$.
This occurs for instance for points $(x,\xi)$ such that $\vert x\vert $ is very
large and the bicharacteristic crosses back a neighborhood of the origin.
That's why we used the term rough for this description.

\noindent {\bf Proof of Proposition \ref{pIII.4.1} } If $\vert x\cdot \xi\vert
\leq -\frac 14
\bra x\ket\, \vert \xi\vert $  then Proposition
\ref{pIII.3.1} gives the claimed description of the flow for $t\geq 0$ and
$t\leq 0$. If $x\cdot \xi\leq -\frac 14\bra x\ket\, \vert \xi\vert$ the same
Proposition applies for $t\leq 0$ so we are left with the case $t\geq 0$. (The
case $x\cdot \xi\geq -\frac 14\bra x\ket\, \vert \xi\vert$ is symmetric).
It follows from Lemma \ref{lIII.1.1} that, if $\varepsilon\,A_1$ is small
enough, we have ${\displaystyle\lim_{t\rightarrow +\infty }} x(t)\cdot
\xi(t)=+\infty
$. Since
$x\cdot \xi\leq 0$ one can find $t^*>0$ such that $x(t^*,x,\xi)\cdot
\xi(t^*,x,\xi)=0$. If we set $x^*=x(t^*,x,\xi)$
$\xi^*=\xi(t^*,x,\xi)$ then, according
to Definition \ref{dIII.2.2}, we have $(x^*,\xi^*)\in\CS_+\cap\CS_-$ so we can
use Proposition \ref{pIII.3.1} for $t\in\R$. Now we have by the flow property
for $t\geq 0$,
$$
x(t,x,\xi)=x(t-t^*,x^*,\xi^*)\,.
$$
Using Proposition \ref{pIII.3.1} we deduce the following lower bound
\begin{equation}\label{eqIII.4.1}
\bra x(t,x,\xi)\ket=\bra x(t-t^*,x^*,\xi^*)\ket\geq \frac{1}{\sqrt 3}\, \bra
t-t^*\ket\,.
\end{equation}
This proves the part (i) in Proposition \ref{pIII.4.1}. To prove part (ii) we
use the formulas (\ref{eqIII.1.2}) for the flow. Then we see that for $t\geq 0$,
$$
\xi_\ell (t,x,\xi)=\xi_\ell +\zeta_\ell (t,x,\xi)\,,\enskip \zeta_\ell
(t,x,\xi)=-\varepsilon \int^t_0 \sum^n_{jk=1} \frac{\partial b_{jk}}{\partial
x_\ell }\, (x(s))\,\xi_j(s)\,\xi_k(s)\,ds\,.
$$
Then using (\ref{eqIII.1.1}), (\ref{eqIII.4.1}), (\ref{eqIII.1.3}) and the fact
that $\v\xi\v\leq 2$ we see  that $\vert \zeta_\ell (t,x,\xi)\vert
\leq C\,\varepsilon$, where $C$ depends only on $A_1$.

On the other hand we have
$$
\dot x_j(t,x,\xi)=2\xi_j+2\zeta_j(t,x,\xi)+2\varepsilon \sum^n_{k=1}
b_{jk}(x(t,x,\xi))\,\xi_k(t,x,\xi)\,.
$$
Integrating between $0$ and $t$ and using the above estimates we obtain the
claimed description of $x(t,x,\xi)$.\cqfd

\subsection{The flow from points in $(\CS_+\cap\CS_-)^c$}\label{ssIII.5}

We study now, more carefully the flow from points $(x,\xi)\in
T^*\R^n\setminus\{0\}$, such that, 
\begin{equation}\label{eqIII.5.1}
\vert x\cdot \xi\vert >c_0\,\bra x\ket\,\vert \xi\vert \textrm{ and } \frac
12\leq \vert \xi\vert \leq 2\,.
\end{equation}
Even if, as we said before, we do not have a nice representation of the flow
for all $t$ in $\R$ we shall see that such a representation is available for
limited values of $t$.

Since the description is symmetric, we shall assume that
\begin{equation}\label{eqIII.5.2}
x\cdot \xi\leq -c_0\,\langle x\rangle \,\vert \xi\vert \,.
\end{equation}
Then $(x,\xi)\in\CS_-$ and Proposition \ref{pIII.3.1} give a good description
of the flow for $t\leq 0$.
\begin{definition}\sl \label{dIII.5.1}
Let $(x,\xi)$ satisfying \ref{eqIII.5.2}. We set
$$
I_+=\Big\{t\geq 0:x(t,x,\xi)\cdot \xi(t,x,\xi)\leq \frac 14 \bra x
(t,x,\xi)\ket\,\vert \xi(t,x,\xi)\vert \Big\}~.
$$
In other word $I_+$ is the set of $t\geq 0$ such that
$(x(t,x,\xi),\xi(t,x,\xi))$ belongs to $\CS_-$.
\end{definition}

The main result of this Section is the following description of the flow
on
$I_+$.
\begin{proposition} \sl \label{pIII.5.2}
Let $(x,\xi)$ satisfying \ref{eqIII.5.2}. Then for $t$ in $I_+$ we have
\begin{equation*}
\begin{aligned}
x(t,x,\xi)&=x+2t\xi-z(-t,x(t,x,\xi),\xi(t,x,\xi))\,,\\
\xi(t,x,\xi)&=\xi-\zeta(-t,x(t,x,\xi),\xi(t,x,\xi))\,,
\end{aligned}
\end{equation*}
where $z$ and $\zeta$ have been defined in Proposition \ref{pIII.3.1}. Moreover
for $j,k=1,\ldots ,n$ we have,
\begin{equation*}
\begin{array}{ll}
\frac{\partial x_j}{\partial
x_k}\,(t,x,\xi)=\delta_{jk}+\CO(\varepsilon)\,,\enskip &\frac{\partial
x_j}{\partial \xi_k}\,(t,x,\xi)=2t\,\delta_{jk}+\CO(\varepsilon\bra t\ket)\\
\frac{\partial \xi_j}{\partial
\xi_k}\,(t,x,\xi)=\delta_{jk}+\CO(\varepsilon\bra t\ket)\,,\enskip
&\frac{\partial \xi_j}{\partial
x_k}\,(t,x,\xi)=\CO(\varepsilon)
\end{array}
\end{equation*}
where $\delta_{jk}$ is the Kronecker symbol and $\CO(A)$ means ``bounded by
$CA$'' with $C$ independent of $(x,\xi)$. In particular we have
$$
\frac{\partial \xi_j}{\partial \xi_k}\, (t,x,\xi)-i\,\frac{\partial
x_j}{\partial \xi_k}\,(t,x,\xi)=(1-2it)\,\delta_{jk}+\CO(\varepsilon\bra
t\ket)\,.
$$
\end{proposition}

\noindent {\bf Proof } As said before, for $t\in I_+$ the point
$\rho(t,x,\xi)=(x(t,x,\xi),\xi(t,x,\xi))$ belongs to $\CS_-$. Therefore we can
apply Proposition \ref{pIII.3.1} for $\theta\leq 0$. We get 
\begin{equation*}
\begin{array}{l}
x(\theta,\rho(t,x,\xi))=x(t,x,\xi)+2\theta\,\xi(\theta,\rho(t,x,\xi))+z(\theta,\rho(t,x,\xi))\\
\xi(\theta,\rho(t,x,\xi))=\xi(t,x,\xi)+\zeta(\theta,\rho(t,x,\xi))\,.
\end{array}
\end{equation*}
Taking $\theta=-t$ with $t\geq 0$ we obtain
\begin{equation*}
\begin{array}{l}
x=x(t,x,\xi)-2t\xi+z(-t,\rho(t,x,\xi))\,,\\
\xi=\xi(t,x,\xi)+\zeta(-t,\rho(t,x,\xi))~.
\end{array}
\end{equation*}
This proves the first part of Proposition \ref{pIII.5.2}. To prove the claim on
the derivatives we use Lemma \ref{lIII.1.2} and Corollary \ref{cIII.3.3}.
\cqfd
\begin{remark} \sl \label{rIII.5.3}
Since the points $(x,\xi)$ satisfying \ref{eqIII.5.2} belong to $\CS_-$,
Propositions \ref{pIII.3.1} and \ref{pIII.5.2} provide a description of the
flow on $(-\infty ,0)\cup I_+$.
\end{remark}

\section{The phase equation}\label{sIV}

\subsection{Statement of the result}\label{ssIV.1}

Let $p(x,\xi)=\v\xi\v^2+\varepsilon\,q(x,\xi)$, $q(x,\xi)=\som^n_{j,k=1}
b_{jk}(x)\,\xi_j\,\xi_k$ where the coefficients $b_{jk}$ satisfy the condition
(\ref{eqIII.1.1}).

In this Section $\alpha=(\alpha_x,\alpha_\xi)$ will be a fixed point in
$T^*\R^n$ such that $\frac 12\leq \v\alpha_\xi\v\leq 2$. Let us recall that
$(x(t,\alpha),\xi(t,\alpha))$ denotes the flow of $p$ starting for $t=0$ at the
point $\alpha$.

We introduce now several sets.
\begin{definition}\sl \label{dIV.1.1}
Let $\delta>0$, $c_0>0$, $c_1>0$ be small constants (chosen later on).
\begin{itemize}
\item[(i)] If $\v\alpha_x\cdot \alpha_\xi\v\leq c_0\,\bra
\alpha_x\ket\,\v\alpha_\xi\v$ we set,
\begin{equation}\label{eqIV.1.1}
\Omega_\delta=\big\{(\theta,x)\in\R\times\R^n : \v
x-x(\theta,\alpha)\v<\delta\,\bra\theta\ket\big\}\,.
\end{equation}
\item[(ii)] If $\alpha_x\cdot \alpha_\xi>
c_0\,\bra\alpha_x\ket\,\v\alpha_\xi\v$ we set,
\begin{equation}\label{eqIV.1.2}
\Omega_\delta=\big\{(\theta,x)\in\R\times\R^n : \v x-x(\theta,\alpha)\v\leq
\delta\,\bra\theta\ket\,,\enskip x\cdot \alpha_\xi\geq -c_1\,\bra
x\ket\,\v\alpha_\xi\v\big\}\,.
\end{equation}
\item[(iii)] If $\alpha_x\cdot \alpha_\xi<
-c_0\,\bra\alpha_x\ket\,\v\alpha_\xi\v$ we set,
\begin{equation}\label{eqIV.1.3}
\Omega_\delta=\big\{(\theta,x)\in\R\times\R^n : \v x-x(\theta,\alpha)\v\leq
\delta\,\bra\theta\ket\,,\enskip x\cdot \alpha_\xi\leq c_1\,\bra
x\ket\,\v\alpha_\xi\v\big\}\,.
\end{equation}
\end{itemize}
\end{definition}

Let us give some explanations on this Definition.

Taking $c_0$ and $c_1$ small with respect to $\frac 14$ we see from Definition
\ref{dIII.2.2} that the case (i) corresponds to points $(\alpha_x,\alpha_\xi)$
which are outgoing for $\theta\geq 0$ and $\theta\leq 0$. Then $\Omega_\delta$
is simply a conic neighborhood of the projection of the bicharacteristic. In
the case (ii) the point $(\alpha_x\,\alpha_\xi)$ is outgoing for $\theta\geq
0$ and $\Omega_\delta$ can be written as follows.
\begin{equation}\label{eqIV.1.4}
\begin{split}
\Omega_\delta=\big\{(\theta,x)\in (0,+\infty )\times\R^n : \v
x-x(\theta,\alpha)\v\leq \delta\,\bra\theta\ket\big\}\cup\big\{(\theta,x)\in
(-\infty ,0)\times\R^n :\\
\v x-x(\theta,\alpha)\v\leq \delta\,\bra\theta\ket\enskip \textrm{and}\enskip
x\cdot \alpha_\xi\geq -c_1\,\bra x\ket\,\v\alpha_\xi\v\big\}\,.
\end{split}
\end{equation}
Indeed if $\v x-x(\theta,\alpha)\v\leq \delta\,\bra\theta\ket$ and $\theta\geq
0$ we have by Proposition \ref{pIII.4.1}, $x\cdot
\alpha_\xi=(x-x(\theta,\alpha))\cdot \alpha_\xi+\alpha_x\cdot
\alpha_\xi+2\theta\,\v\alpha_\xi\v^2+\CO(\varepsilon\,\bra\theta\ket)$.  Since
$\v\alpha_\xi\v\geq \frac 12$ and we are in case (ii) we deduce that
$x\cdot \alpha_\xi\geq c_0\,\bra\alpha_x\ket\,\v\alpha_\xi\v+\frac
12\,\theta-C(\varepsilon+\delta)\bra\theta\ket\geq 0$ if $\varepsilon+\delta$
is small enough. Therefore when $\theta\geq 0$ the condition $x\cdot
\alpha_\xi\geq -c_1\,\bra x\ket\,\v\alpha_\xi\v$ is automatically satisfied.

In the case (iii) we have the same discussion changing $\theta\geq 0$ to
$\theta\leq 0$.

The purpose of this Section is to prove the following result.
\begin{theorem}\sl \label{tIV.1.2}
There exist $\delta>0$, $c_0>0$, $c_1>0$ such that for any $\alpha\in T^*\R^n$
with $\frac 12\leq \v\alpha_\xi\v\leq 2$ one can find a function $\varphi
=\varphi (\theta,x,\alpha)$ on $\Omega_\delta$ which is $C^\infty $ and
satisfies the following.
\begin{itemize}
\item[(i)] $\displaystyle{\varphi (0,x,\alpha)=(x-\alpha_x)\cdot
\alpha_\xi+\frac i2\,
\v x-\alpha_x\v^2+\frac 1{2i}\,\v\alpha_\xi\v^2+g(x,\alpha)}$

 where $\v
g(x,\alpha)\v\leq C_N\,\v x-\alpha_x\v^N$ for all $N\in\N$.

\item[(ii)] For any $N\in\N$ there exists $C_N\geq 0$ such that
$$
\Big\v\frac{\partial \varphi }{\partial \theta}\,(\theta,x,\alpha)+p\Big
(x,\,\frac{\partial \varphi }{\partial x}\,(\theta,x,\alpha)\Big)\Big\v\leq
C_N\Big (\frac{\v x-x(\theta,\alpha)\v}{\bra\theta\ket}\Big)^N
$$
for all $(\theta,x)$ in $\Omega_\delta$.

Moreover for $(\theta,x)$ in $\Omega_\delta$ we have
\item[(iii)] 
$\displaystyle{
\Big\v\frac{\partial \varphi }{\partial
x}\,(\theta,x,\alpha)-\alpha_\xi\Big\v\leq
C(\varepsilon+\sqrt\delta)\,.}$
\item[(iv)] $\displaystyle{\Im \varphi (\theta,x,\alpha)=\frac 12\,\frac{\v
x-x(\theta,\alpha)\v^2}{1+4\theta^2}+o(1)\,\frac{\v
x-x(\theta,\alpha)\v^2}{\bra\theta\ket^2}-\frac 12 \, \v\alpha_\xi\v^2}$.
\item[(v)] $\displaystyle{\v \partial ^A_x\,\varphi (\theta,x,\alpha)\v\leq
C_A\,,\enskip \textrm{for every}\enskip A\enskip \textrm{in}\enskip
\N^n\setminus \{0\}}$

where $o(1)$ and $C_A$ are independent of $(\theta,x,\alpha)$.
\end{itemize}
\end{theorem}
The proof of this result is based on the theory of Lagrangian ideals of L.
H\"ormander (\cite{H}, vol~4, chap.~XXV). It will require several steps.
The first one is a slight extension of Theorem~7.5.4 in  \cite{H}, vol.~1 to
the case of higher dimensions.

\subsection{The preparation theorem}\label{ssIV.2}

The aim of this Section is to prove the following result.
\begin{lemma}\sl \label{lIV.2.1}
Let $g\in\CS(\R^n_\xi)$ and $z\in\C^n$. Then there exist functions
$q_j(\xi,z,g)$, $j=1,\ldots ,n$, $r(z,g)$ which are $C^\infty $ with respect to
$\xi$ and $z$, which depend linearly on $g$ such that
\begin{equation}\label{eqIV.2.1}
g(\xi)=\sum^n_{j=1} q_j(\xi,z,g)(\xi_j+z_j)+r(z,g)
\end{equation}
\begin{equation}\label{eqIV.2.2}\left\{
\begin{array}{l}
\v\partial ^\alpha_\xi\,\partial ^\beta_z\,q_j(\xi,z,g)\v\leq C_{\alpha\beta}
\som_{\v\gamma\v\leq \v\alpha\v+\v\beta\v+4n} \int \v\partial
^\gamma_\eta\,g(\eta)\v\,d\eta\\
\v\partial ^\beta_z\,r(z,g)\v\leq C_\beta \som_{\v\gamma\v\leq
\v\beta\v+3n}\int 
\v\partial ^\gamma\,g(\eta)\v\,d\eta\,.
\end{array}\right.
\end{equation}
\end{lemma}

{\bf Proof } We proceed by induction on the dimension $n$. If $n=1$ this
follows from Theorem 7.5.4 of \cite{H}. Let $n\geq 2$ and let us set
$\xi'=(\xi_1,\ldots ,\xi_{n-1})$. For fixed $\xi'\in \R^{n-1}$ we apply Theorem
7.5.4 of \cite{H} to the function $\xi_n\mapsto g(\xi',\xi_n)$. We get
\begin{equation}\label{eqIV.2.3}
g(\xi',\xi_n)=q(\xi_n,z_n,g(\xi',\cdot ))(\xi_n+z_n)+r(z_n,g(\xi',\cdot ))\,.
\end{equation}
Let us set $Q_n(\xi,z_n,g)=q(\xi_n,z_n,g(\xi',\cdot ))$ and $\tilde
r(z_n,\xi',g)=r(z_n,g(\xi',\cdot ))$. 
Since $r$ is linear in $g$ we have $\partial ^\alpha_{\xi'}\,\tilde
r(z_n,\xi',g)=r(z_n,\partial ^\alpha_{\xi'}\,g(\xi',\cdot ))$ and the estimates
(\ref{eqIV.2.2}) for $n=1$ show that $\xi'\mapsto \tilde r(z_n,\xi',g)$ is in
$\CS(\R^{n-1})$. Therefore we can apply, by the induction, the Lemma to the
function $\xi'\mapsto\tilde r(z_n,\xi',g)$ and to $z'=(z_1,\ldots ,z_{n-1})$.
We obtain the existence of $q_j$, $j=1,\ldots ,n-1$ and $R$ satisfying the
estimates (\ref{eqIV.2.2}) such that
$$
\tilde r(z_n,\xi',g)=\sum^{n-1}_{j=1}q_j(\xi',z',\tilde r(z_n,\cdot
,g)(\xi_j+z_j)+R(z',\tilde r(z_n,\cdot ,g))\,.
$$
Using (\ref{eqIV.2.3}) we obtain therefore
\begin{equation*}
\begin{aligned}
g(\xi)=Q_n(\xi,z_n,g)(\xi_n+z_n)&+\som^n_{j=1} q_j(\xi',z',\tilde r(z_n,\cdot
,g))(\xi_j+z_j)\\
&+R(z',\tilde r(z_n,\cdot ,g))\,.
\end{aligned}
\end{equation*}
If we set
\begin{equation}\label{eqIV.2.4}\left\{
\begin{array}{l}
Q_j(\xi,z,g)=q_j(\xi',z',\tilde r(z_n,\cdot ,g))\,,\quad j=1,\ldots ,n-1\\
r(z,g)=R(z',\tilde r(z_n,\cdot ,g))
\end{array}\right.
\end{equation}
we obtain (\ref{eqIV.2.1}) at the level $n$. Moreover $Q_j$ and $r$ are linear
in $g$ since $q_j,R$ are linear in $\tilde r(z_n,\cdot ,g)$ and $\tilde r$ is
linear in $g$. Let us look to the estimate (\ref{eqIV.2.2}) for $r$. We have
$$
\partial ^{\beta'}_{z'}\,\partial ^{\beta_n}_{z_n}\,r(z,g)=\partial
^{\beta'}_{z'}\,R\big (z',\partial ^{\beta_n}_{z_n}\,\tilde r(z_n,\cdot ,g)\big)
$$
so,
$$
\v\partial ^\beta_z\,r(z,g)\v\leq C\sum_{\v\gamma'\v\leq \v\beta'\v+3(n-1)} \int
\big\v\partial ^{\gamma'}_{\xi'}\,\partial ^{\beta_n}_{z_n}\,\tilde
r(z_n,\xi',g)\big\v\,d\xi'\,.
$$
Now $\v\partial ^{\gamma'}_{\xi'}\,\partial ^{\beta_n}_{z_n}\,\tilde
r(z_n,\xi',g)=\partial ^{\beta_n}_{z_n}\,r(z_n,\partial
^{\gamma'}_{\xi'}\,g(\xi',\cdot ))$ and from the case $n=1$ we have
$$
\big\v\partial ^{\beta_n}_{z_n}\,r(z_n,\partial ^{\gamma'}_{\xi'}\,g(\xi',\cdot
))\big\v\leq C \sum_{\v\gamma_n\v\leq \beta_n+3} \int \big\v\partial
^{\gamma'}_{\xi_n}\,\partial
^{\gamma'}_{\xi'}\,g(z_n,\xi',\xi_n)\big\v\,d\xi_n\,.
$$
It follows that
\begin{equation}\label{eqIV.2.5)}
\v\partial ^\beta_z\,r(z,g)\v\leq \sum_{\v\gamma'\v\leq
\v\beta'\v+3(n-1)}\,\sum_{\v\gamma_n\v\leq \beta_n+3} \int \big\v\partial
^{\gamma'}_{\xi'}\,\partial
^{\gamma_n}_{\xi_n}\,g(\xi',\xi_n)\big\v\,d\xi'\,d\xi_n\,.
\end{equation}
The proof of the estimates for the $q'_js$ is the same. \cqfd

\begin{remark}\sl \label{rIV.2.2}
Let us set  $z=a+ib$ and let us write $r(z,g)=r(a,b,g)$ and
$q_j(\xi,z,g)=q_j(\xi,a,b,g)$. If we take in (\ref{eqIV.2.1}) $b=0$, $\xi=-a$
we obtain
\begin{equation}\label{eqIV.2.6}
r(a,0,g)=g(-a)\,.
\end{equation}
If we differentiate (\ref{eqIV.2.1}) with respect to $b_k$ and then take
$z=a\in\R^n$, $\xi=-a$, we get
\begin{equation}\label{eqIV.2.7}
\frac{\partial r}{\partial b_k}\,(a,0,g)=-i\,q_k(-a,a,0,g)\,,\quad k=1,\ldots
,n\,.
\end{equation}
Finally if we differentiate (\ref{eqIV.2.1}) with respect to $\xi_\ell $ and
then take $z=a\in\R^n$, $\xi=-a$ we obtain
\begin{equation}\label{eqIV.2.8}
\frac{\partial g}{\partial \xi_\ell }\,(-a)=q_\ell (-a,a,0,)\,,\quad \ell
=1,\ldots ,n\,.
\end{equation}
\end{remark}
We introduce now the following notations which will be used in the next
sections.
\begin{notation}\sl \label{nIV.2.3}

Let $\alpha=(\alpha_x,\alpha_\xi)\in T^*\,\R^n\setminus 0$. We introduce
\begin{equation}\label{eqIV.2.9}\left\{
\begin{array}{l}
\varphi _0(x,\alpha)=(x-\alpha_x)\,\alpha_\xi+\frac i2\,(x-\alpha_x)^2+\frac
1{2i}\,\alpha^2_\xi\,,\\
u_j(x,\xi,\alpha)=\xi_j-\frac{\partial \varphi _0}{\partial
x_j}\,(x,\alpha)=\xi_j-\alpha^j_\xi-i(x_j-\alpha^j_x)\,.
\end{array}\right.
\end{equation}
Let $p(x,\xi)=\v\xi\v^2+\varepsilon\,q(x,\xi)$~; we denote by $H_p$ its
hamiltonian and we introduce
\begin{equation}\label{eqIV.2.10}\left\{
\begin{array}{ll}
v_j(\theta;x,\xi,\alpha)&=u_j(\exp (-\theta\,H_p)(x,\xi))\\
&=\xi_j(-\theta;x,\xi)-\alpha^j_\xi-i(x_j(-\theta;x,\xi)-\alpha^j_x)\,,
\end{array}\right.
\end{equation}
where $(x,\xi)$ is close to $(x(\theta;\alpha),\xi(\theta,\alpha))$.

We split the proof of Theorem \ref{tIV.1.2} according to the different values
of $\alpha$ described in Definition \ref{dIV.1.1}.
\end{notation}

\subsection{The case of outgoing points}\label{ssIV.3}
Let us set
\begin{equation}\label{eqIV.3.1}
\CS =\big\{\alpha\in\R^{2n} : \frac 12 \leq \v\alpha_\xi\v\leq 2\,,\enskip 
\v\alpha_x\cdot \alpha_\xi\v\leq c_0\,\bra\alpha_x\ket\,\v\alpha_\xi\v\big\}\,.
\end{equation}
 We shall use the following notations
\begin{equation}\label{eqIV.3.2}\left\{
\begin{array}{l}
\tilde \Omega_\delta=\{(\theta,y)\in\R\times\R^n : \v y\v\leq
\delta\,\bra\theta\ket\}\\
\sgn \theta=1\enskip \textrm{(resp.} -1)\enskip \textrm{if}\enskip
\theta>0\enskip \textrm{(resp.}\enskip \theta<0)\,.
\end{array}\right.
\end{equation}
Let now $\chi_0\in C^\infty _0(\R^n)$, $\chi_1\in C^\infty _0(\R)$ be such that,
$$
\chi_0(t)=1\enskip \textrm{if}\enskip \v t\v\leq 1\,, \enskip \chi_0(t)=0\enskip
\textrm{if}\enskip \v t\v\geq 2\enskip \textrm{and}\enskip 0\leq \chi_0\leq
1\,,
$$
$$
\chi_1(\theta)=1\enskip \textrm{if}\enskip \v \theta\v\leq 1\,,\enskip 
\chi_1(\theta)=0\enskip
\textrm{if}\enskip \v \theta\v\geq 2\enskip \textrm{and}\enskip 0\leq \chi_1\leq
1\,.
$$
Then we can state the following result.
\begin{theorem}\sl \label{tIV.3.1}
There exist small positive constants $\mu_0,\delta$ such that if we set for
$\alpha\in\CS$, $\theta\in\R$, $y\in\R^n$, $\eta\in\R^n$, $j=1,\ldots ,n$,
$$
g_j(\eta)=\chi_0\Big(\frac
1{\mu_0}\,\eta\Big)\,v_j\Big(\theta,y+x(\theta,\alpha),\,\eta\,\chi_1(\theta)+(1-\chi_1(\theta))
\Big[\frac\eta{\bra\theta\ket}+\frac 12\,\frac{\sgn
\theta}{\bra\theta\ket}\,y\Big]+\xi(\theta,\alpha),\alpha\Big)
$$
there exist smooth functions $a_j=a_j(\theta,y,\alpha)$,
$b_j=b_j(\theta,y,\alpha)$ defined on $\tilde \Omega_\delta$ such that, with
$a=(a_j)_{j=1,\ldots ,n}$, $b=(b_j)_{j=1,\ldots ,n}$ we have for $\eta$ in
$\R^n$ and $(\theta,y)\in\tilde \Omega_\delta$,
\begin{itemize}
\item[(i)] $\displaystyle{g_j(\eta)=\som^n_{k=1}
q_k(\eta,a,b,g_j)(\eta_k+a_k(\theta,y,\alpha)+i\,b_k(\theta,y,\alpha)) }$

where the $q'_ks$ have been introduced in Lemma \ref{lIV.2.1}. 

Moreover in the
set $\tilde \Omega_\delta$ we have
\item[(ii)]  $\displaystyle{\v a(\theta,y,\alpha)\v\leq 10\,\frac{\v
y\v}{\bra\theta\ket}}$,\qquad 
$\displaystyle{\big\v
b(\theta,y,\alpha)+\frac{K(\theta)}{1+4\theta^2}\,y\big\v\leq
\sqrt\delta\,\frac {\v y\v}{\bra\theta\ket}}$,

where
$K(\theta)=\frac{\bra\theta\ket}{\bra\theta\ket\,\chi_1(\theta)+1-\chi_1(\theta)}$, 
$\frac 1{\sqrt 5}\,\bra\theta\ket\leq K(\theta)\leq \bra\theta\ket$.

On the other hand we have, uniformly with respect to
$(\theta,y)\in\tilde
\Omega_\delta$ and $\alpha\in\CS$,
\item[(iii)] $\displaystyle{\v\partial ^A_y\,a(\theta,y,\alpha)\v+\v\partial
^\alpha_y\,b(\theta,y,\alpha)\v\leq \frac{C_A}{\bra\theta\ket^{\v
A\v}}\,,\enskip A\in\N^n }$.

Moreover for $j=1,\ldots ,n$, $k=1,\ldots ,n$,
\item[(iv)] $\displaystyle{\big\v
q_k(\eta,a,b,g_j)-(1+2i\theta)\,\frac{k(\theta)}{{\bra\theta\ket}}\,\delta_{jk}\big\v
\leq C\,(\varepsilon+\delta)}$, if $\v\eta\v\leq \delta$

where $k(\theta)={\bra\theta\ket}\,\chi_1(\theta)+1-\chi_1(\theta)$.
\item[(v)] $\v\partial ^A_{(a,b)}\,\partial ^B_\eta\,q_k(\eta,a,b,g_j)\v\leq C
(\mu_0)$, if $\v A\v+\v B\v\geq 1$, $\v\eta\v\leq \mu_0$, $1\leq j,k\leq n$.
\end{itemize}
\end{theorem}

{\bf Proof } According to Definition \ref{dIII.2.2} we have
$\CS\subset\CS_+\cap\CS_-$. Moreover
\begin{equation}\label{eqIV.3.3}
B(\alpha,c_0):=\{\tilde \alpha\in T^*\,\R^n : \v\alpha-\tilde \alpha\v\leq
c_0\}\subset\CS_+\cap\CS_-\,.
\end{equation} 

We start with the following Lemma.
\begin{lemma}\sl\label{lIV.3.2}
There exists a small positive constant $\mu_0$ such that for all $(\theta,y)$
with $\theta\in\R$, $\v y\v\leq \mu_0\,{\bra\theta\ket}$ and all $\eta$ in
$\R^n$ such that $\v\eta\v\leq 2\,\mu_0$ there exists a unique
$\beta=\beta(\theta,y,\alpha,\eta)$ in $B(\alpha,c_0)$ such that
\begin{equation}\label{eqIV.3.4}\left\{
\begin{array}{l}
x(\theta,\beta)=y+x(\theta,\alpha)\\
\xi(\theta,\beta)=\chi_1(\theta)\,\eta+(1-\chi_1(\theta))\big[\frac\eta{\bra\theta\ket}
+\frac 12\,\frac{\sgn \theta}{\bra\theta\ket}\,y\big]+\xi(\theta,\alpha)\,.
\end{array}\right.
\end{equation}
Moreover we have
\begin{equation}\label{eqIV.3.5}\left\{
\begin{array}{l}
\beta_x=\alpha_x+\big[1-\frac{\v\theta\v}{\bra\theta\ket}\,(1-\chi_1(\theta))\big]\,
y-\frac{2\theta}{\bra\theta\ket}\,\eta(\chi_1(\theta)\bra\theta\ket+1-\chi_1(\theta))
+z(\theta,\alpha)-z(\theta,\beta)\\
\beta_\xi=\alpha_\xi+\frac\eta{\bra\theta\ket}\,[\chi_1(\theta)\bra\theta\ket+(1-\chi_1(\theta))]
+\frac 12\, \frac{1-\chi_1(\theta)}{\bra\theta\ket}\,\sgn
\theta\,y+\zeta(\theta,\alpha)-\zeta(\theta,\beta)
\end{array}\right.
\end{equation}
where $z$ and $\zeta$ have been introduced in Proposition \ref{pIII.3.1} and
\begin{equation}\label{eqIV.3.6}\left\{
\begin{array}{ll}
\textrm{(i)} &\v\beta-\alpha\v\leq 10\big (\v\eta\v+\frac{\v
y\v}{\bra\theta\ket}\big )\\
\textrm{(ii)} &\big\v\frac{\partial \beta}{\partial
\eta}\,(\theta,y,\alpha,\eta)\big\v\leq C\\
\textrm{(iii)} &\big\v\big(\frac{\partial \beta^j_\xi}{\partial \eta_k}-i\,
\frac{\partial \beta^j_x}{\partial
\eta_k}\big)(\theta,y,\alpha,\eta)-\frac{1+2i\theta}{\bra\theta\ket}\,(\bra\theta\ket\,
\chi_1(\theta)+(1-\chi_1(\theta)))\,\delta_{jk}\big\v\leq C\,\varepsilon\\
\textrm{(iv)} &\v\partial ^\alpha_\eta\,\partial
^B_y\,\beta(\theta,y,\alpha,\eta)\v\leq C_{AB}\,\varepsilon\,, \enskip
\textrm{if}\enskip \v A\v+\v B\v\geq 2\,.
\end{array}\right.
\end{equation}
\end{lemma}

{\bf Proof }
The system (\ref{eqIV.3.4}) with $\beta\in B(\alpha,c_0)\subset \CS_+\cap\CS_-$
is equivalent by Proposition \ref{pIII.3.1} to the following
\begin{equation}\label{eqIV.3.7}\left\{
\begin{array}{l}
\beta_x+2\theta\xi(\theta,\beta)+z(\theta,\beta)=y+\alpha_x+2\theta\xi
(\theta,\alpha)+z(\theta,\alpha)\\
\beta_\xi+\zeta(\theta,\beta)=\chi_1(\theta)\,\eta+\frac{1-\chi_1(\theta)}{\bra\theta\ket}\,
\big[\eta+\frac 12\,\sgn \theta\,y\big]+\alpha_\xi+\zeta(\theta,\alpha)\,.
\end{array}\right.
\end{equation}
Using again (\ref{eqIV.3.4}) the left hand side of the first line of
(\ref{eqIV.3.7}) can be written
$$
\beta_x+2\theta\,\chi_1(\theta)\,\eta+\frac{2\theta}{\bra\theta\ket}\,(1-\chi_1(\theta))\,\eta
+\frac{\v\theta\v}{\bra\theta\ket}\,(1-\chi_1(\theta))\,y+2\theta\xi(\theta,\alpha)+z(\theta,\beta)\,.
$$
Finally (\ref{eqIV.3.7}) is equivalent to
\begin{equation}\label{eqIV.3.8}\left\{
\begin{array}{l}
\beta_x=\alpha_x+\big[1-\frac{\v\theta\v}{\bra\theta\ket}\,(1-\chi_1(\theta))\big]\,
y-\frac{2\theta}{\bra\theta\ket}\,\eta(\chi_1(\theta)\bra\theta\ket+1-\chi_1(\theta))+z(\theta,\alpha)-z(\theta,\beta)\\
\beta_\xi=\alpha_\xi+\frac\eta{\bra\theta\ket}\,
[\chi_1(\theta)\bra\theta\ket+(1-\chi_1(\theta))]+\frac
12\,\frac{1-\chi_1(\theta)}{\bra\theta\ket}\,\sgn
\theta\,y+\zeta(\theta,\alpha)-\zeta(\theta,\beta)\,.
\end{array}\right.
\end{equation}
Writing this system $\beta_x=\Phi_x(\beta)$\, $\beta_\xi=\Phi_\xi(\beta)$ and
setting $\Phi(\beta)=(\Phi_x(\beta),\Phi_\xi(\beta))$ we are going to solve it
using the fixed point theorem in $B(\alpha,c_0)$.

(i) $\Phi$ maps $B(\alpha,c_0)$ in itself.

We have $0<1-\frac{\v\theta\v}{\bra\theta\ket}\leq \frac 1{\bra\theta\ket^2}$,
$\v\theta\v\,\chi_1(\theta)\leq 2$,
$\chi_1(\theta)\bra\theta\ket+1-\chi_1(\theta)\leq \sqrt 5$, $\v y\v\leq
\mu_0\,\bra\theta\ket$, $\v\eta\v\leq 2\,\mu_0$, $\v z\v+\v\zeta\v\leq
C\,\varepsilon$ by Proposition \ref{pIII.3.2}. It follows that
$\v\Phi(\beta)-\alpha\v\leq 20\,\mu_0+C\,\varepsilon\leq c_0$ if $\mu_0$ and
$\varepsilon$ are small enough.

(ii) Let $\beta,\beta'$ be in $B(\alpha,c_0)$. Then $t\,\beta+(1-t)\,\beta'\in
B(\alpha,c_0)\subset\CS_+\cap \CS_-$ for all $t$ in $(0,1)$.  It follows that
$$
\v\Phi(\beta)-\Phi(\beta')\v\leq \v
z(\theta,\beta)-z(\theta,\beta')\v+\v\zeta(\theta,\beta)-\zeta(\theta,\beta')\v\leq
C\,\varepsilon \,\v\beta-\beta'\v
$$
by Proposition \ref{pIII.3.2}. Here $C$ depends only on the constants $A_0,A_1$
in (\ref{eqIII.1.1}). Taking $\varepsilon$ so small that $C\,\varepsilon<1$ we
see that we can apply the fixed point theorem in $B(\alpha,c_0)$. This proves
the existence of $\beta$ satisfying (\ref{eqIV.3.4}) and (\ref{eqIV.3.5}) by
(\ref{eqIV.3.8}).
Now (i) in (\ref{eqIV.3.6}) follows from (\ref{eqIV.3.5}) taking $\varepsilon$
small enough since $\v Z(\theta,\alpha)-Z(\theta,\beta)\v\leq
C\,\varepsilon\,\v\alpha-\beta\v$ if $Z=z$ or $\zeta$. The claim (ii) is
obtained by differentiating the equations (\ref{eqIV.3.5}) with respect to
$\eta_k$ and using Proposition \ref{pIII.3.2}. Then (iii) follows easily from
(\ref{eqIV.3.5}) and (ii). Finally (iv) is obtained by an induction on $\v
A\v+\v B\v$.
\cqfd

From now on we fix the constant $\mu_0$ occuring in Lemma \ref{lIV.3.2}.

Now for $j=1,\ldots ,n$ let $g_j$ be the function introduced in the statement
of Theorem \ref{tIV.3.1}. Then, according to Lemma \ref{lIV.3.2} and
(\ref{eqIV.2.10}) we have for $\v y\v\leq \mu_0\,\bra\theta\ket$ and
$\eta\in\R$,
\begin{equation}\label{eqIV.3.9}
g_j(\eta)=\chi_0\Big(\frac\eta{\mu_0}\Big)
\big[\beta^j_\xi(\theta,y,\alpha,\eta)-\alpha^j_\xi-i(\beta^j_x(\theta,y,\alpha,\eta)-\alpha^j_x)
\big]\,,
\end{equation}
since $f(-\theta,x(\theta,\beta),\xi(\theta,\beta))=\beta_f$ for $f=x$ and
$\xi$. It follows from Lemma \ref{lIV.2.1} that the existence of $a_j$, $b_j$
in Theorem \ref{tIV.3.1} will be proved if we can solve the equations
\begin{equation}\label{eqIV.3.10}
r(a,b,g_j(\cdot ))=0\,,\quad j=1,\ldots ,n\,.
\end{equation}
Let us now take $(\theta,y)\in\tilde \Omega_\delta$, that is $\theta\in\R$, $\v
y\v\leq \delta\,\bra\theta\ket$ where $0<\delta<\frac 12\,\mu_0$ is to be
chosen. We look for a solution $(a,b)$ of the system (\ref{eqIV.3.10}) in the
set
\begin{equation}\label{eqIV.3.11}\left\{
\begin{array}{l}
E=\big\{(a,b)\in\R^n\times\R^n : \v a\v\leq \frac{10\,\v
y\v}{\bra\theta\ket}\,,\enskip \big\v
b+\frac{K(\theta)}{1+4\,\theta^2}\,y\big\v\leq \sqrt\delta\,\frac{\v
y\v}{\bra\theta\ket}\big\}\\
\textrm{where}\enskip
K(\theta)=\frac{\bra\theta\ket}{\bra\theta\ket\,\chi_1(\theta)+(1-\chi_1(\theta))}\,,
\enskip \frac 1{\sqrt 5}\, \bra\theta\ket\leq K(\theta)\leq \bra\theta\ket\,.
\end{array}\right.
\end{equation}
We shall first give equivalent equations to (\ref{eqIV.3.10}) in the set $E$.
We write,
$$
r(a,b,g_j)=r(a,0,g_j)+\sum^n_{k=1} \frac{\partial r}{\partial
b_k}\,(a,0,g_j)\,b_k+\som^n_{p,q=1} H^j_{p,q}(\theta,y,\alpha,a,b)\,b_p\,b_q
$$
where
\begin{equation}\label{eqIV.3.12}
H^j_{p,q}(\theta,y,\alpha,a,b)=\int^1_0 (1-t)\,\frac{\partial ^2r}{\partial
b_p\,\partial b_q}\,(a,t\,b,g_j(\cdot ))\,dt\,.
\end{equation}
It follows from (\ref{eqIV.2.6}), (\ref{eqIV.2.7}) and (\ref{eqIV.2.8}) that
\begin{equation}\label{eqIV.3.13}
r(a,b,g_j(\cdot ))=g_j(-a)-i \sum^n_{k=1} \frac{\partial g_j}{\partial \eta_k}\,
(-a)\,b_k+\sum^n_{p,q=1} H^j_{p,q}(\theta,y,\alpha,a,b)\,b_p\,b_q\,.
\end{equation}
Now if $(a,b)\in E$ we have $\v a\v\leq \frac{12\,\v y\v}{\bra\theta\ket}\leq
12\, \delta\leq \mu_0$. Therefore $\chi_0\big (-\frac a{\mu_0}\big)=1$,
$\chi'_0\big(-\frac a{\mu_0}\big)=0$. Then by (\ref{eqIV.3.9}) and
(\ref{eqIV.3.5}) we obtain,
\begin{equation}\label{eqIV.3.14}
\begin{split}
g_j(-a)&=\frac 12\,\frac{1-\chi_1(\theta)}{\bra\theta\ket}\,\sgn
\theta\,y_j-\frac{a_j}{\bra\theta\ket}\,
[\chi_1(\theta)\bra\theta\ket+1-\chi_1(\theta)]+\zeta_j(\theta,\alpha)
-\zeta_j(\theta,\beta)\\
& -i\Big
(1-\frac{\v\theta\v}{\bra\theta\ket}\,(1-\chi_1(\theta))\Big)\,y_j-\frac{2i\theta}
{\bra\theta\ket}\,a_j[\chi_1(\theta)\bra\theta\ket+1-\chi_1(\theta)]
-i(z_j(\theta,\alpha)-z_j(\theta,\beta)
\end{split}
\end{equation}
\begin{equation}\label{eqIV.3.15}
\begin{split}
\frac{\partial g_j}{\partial \eta_k}\,(-a)=\frac{1+2i\theta}{\bra\theta\ket}\,
[\chi_1(\theta)\bra\theta\ket+(1-\chi_1(\theta))]-\partial
\zeta_j(\theta,\beta(\theta,y,\alpha,-a))
\frac{\partial \beta}{\partial \eta_k}\,(\theta,y,\alpha,-a)\\
+i\,\partial
z_j(\theta,\beta(\theta,y,\alpha,-a))\cdot \frac{\partial \beta}{\partial
\eta_k}\,(\theta,y,\alpha,-a)
\end{split}
\end{equation}
where $\partial =(\partial _x,\partial _\xi)$ and $\beta=(\beta_x,\beta_\xi)$.

On the other hand we deduce from (\ref{eqIV.3.12}) and (\ref{eqIV.2.2}) that
\begin{equation}\label{eqIV.3.16}
\big\v\partial ^A_{(a,b)}\,\partial
^B_y\,H^j_{pq}(\theta,y,\alpha,a,b)\big\v\leq C_{AB} \sum_{\v\gamma\v\leq \v
A\v+3n+2} \int \v\partial ^\gamma_\eta\,\partial ^B_y\,g_j(\eta)\v\,d\eta\,.
\end{equation}
Using (\ref{eqIV.3.6}) and (\ref{eqIV.3.9}) we obtain
\begin{equation}\label{eqIV.3.17}
\big\v \partial ^A_{(a,b)}\,\partial
^B_y\,H^j_{pq}(\theta,y,\alpha,a,b)\big\v\leq C'_{AB}(\mu_0)\,.
\end{equation}
It follows from (\ref{eqIV.3.13}), (\ref{eqIV.3.14}), (\ref{eqIV.3.15}) that
(\ref{eqIV.3.10}) is equivalent to
\begin{equation*}
\begin{split}
-\frac{a_j}{\bra\theta\ket}\, [\chi_1(\theta)\bra\theta\ket&+1-\chi_1(\theta)]-
\frac{2i\theta}{\bra\theta\ket}\,a_j[\chi_1(\theta)\bra\theta\ket+1-\chi_1(\theta)]
+\frac 12\, \frac{1-\chi_1(\theta)}{\bra\theta\ket}\,\sgn \theta\,y_j\\
&-i\Big (1-\frac{\v\theta\v}{\bra\theta\ket}\,(1-\chi_1(\theta))\Big)\,y_j
+\zeta_j(\theta,\alpha)-\zeta_j(\theta,\beta)-i(z_j(\theta,\alpha)-z_j(\theta,\beta))\\
&-i\,\frac{1+2i\theta}{\bra\theta\ket}\,[\chi_1(\theta)\bra\theta\ket+1-\chi_1(\theta)]\,
b_j+F^j_1(\theta,y,\alpha,a)\cdot b+i\,F^j_2(\theta,y,\alpha,a)\cdot b\\
&+H^j_1(\theta,y,\alpha,a,b)\,b\cdot b+i\,H^j_2(\theta,y,\alpha,a,b)\,b\cdot
b=0\,,
\end{split}
\end{equation*}
where
\begin{equation}\label{eqIV.3.18}\left\{
\begin{array}{l}
\beta=\beta(\theta,y,\alpha,-a)\,,\\
F^j_1(\theta,y,\alpha,a)=\partial \,z_j(\theta,\beta)\cdot \frac{\partial
\beta}{\partial \eta}\,(\theta,y,\alpha,-a)\,,\\
F^j_2(\theta,y,\alpha,a)=\partial \,\zeta_j(\theta,\beta)\cdot \frac{\partial
\beta}{\partial \eta}\,(\theta,y,\alpha,-a)\,,\\
H^j=(H^j_{pq})=H^j_1+i\,H^j_2\,.
\end{array}\right.
\end{equation}
Taking the real and the imaginary parts we are led to the system
\begin{equation*}
\begin{aligned}
a_j-2\theta\,b_j&=K(\theta)\Big[\frac
12\,\frac{1-\chi_1(\theta)}{\bra\theta\ket}\,\sgn
\theta\,y_j+\zeta_j(\theta,\alpha)-\zeta_j(\theta,\beta)+F^j_1\,
b+H^j_1\,b\cdot b\Big ]\\
2\theta\,a_j+b_j&=K(\theta)\Big[-\Big
(1-\frac{\v\theta\v}{\bra\theta\ket}\,
(1-\chi_1(\theta))\Big)\,y_j-(z_j(\theta,\alpha)-z_j(\theta,\beta))+F^j_2\,b+H^j_2\,b\cdot
b\Big]
\end{aligned}
\end{equation*}
where
$K(\theta)=\frac{\bra\theta\ket}{\chi_1(\theta)\bra\theta\ket+1-\chi_1(\theta)}$.

Inverting this system we are led to solve
\begin{equation}\label{eqIV.3.19}\left\{
\begin{array}{l}
a_j=\frac{K(\theta)}{1+4\,\theta^2}\,\Big (\frac
12\,\frac{1-\chi_1(\theta)}{\bra\theta\ket}\,\sgn \theta-2\theta\Big
(1-\frac{\v\theta\v}{\bra\theta\ket}\,(1-\chi_1(\theta))\Big)\Big)\,y_j
+Z^j_1(\theta,\alpha)\,-Z^j_1(\theta,\beta)\\
\hbox to 7,2cm{} +F^j_3\,b+H^j_3\,b\cdot b=:\Phi^j_1(a,b)\\
b_j=-\frac{K(\theta)}{1+4\theta^2}\,y_j+Z^j_2(\theta,\alpha)-Z^j_2(\theta,\beta)
+F^j_4\,b+H^j_4\,b\cdot b=:\Phi^j_2(a,b)
\end{array}\right.
\end{equation}
where
\begin{equation}\label{eqIV.3.20}\left\{
\begin{array}{l}
Z^j_1(\theta,\cdot )=\frac{K(\theta)}{1+4\theta^2}\,(\zeta_j(\theta,\cdot
)-2\theta\,z_j(\theta,\cdot ))\\
Z^j_2(\theta,\cdot
)=-\frac{K(\theta)}{1+4\theta^2}\,(2\theta\,\zeta_j(\theta,\cdot
)+z_j(\theta,\cdot ))\\
F^j_3=\frac{K(\theta)}{1+4\theta^2}\,(F^j_1+2\theta\,F^j_2)\,,\quad \kern
5pt F^j_4=
\frac{K(\theta)}{1+4\theta^2}\, (-2\theta\,F^j_1+F^j_2)\\
H^j_3=\frac{K(\theta)}{1+4\theta^2}\, (H^j_1+2\theta\,H^j_2)\,,\quad
H^j_4=\frac{K(\theta)}{1+4\theta^2}\,(-2\theta\,H^j_1+H^j_2)\,.
\end{array}\right.
\end{equation}
Let us set $\Phi^j=(\Phi^j_1,\Phi^j_2)$ (see (\ref{eqIV.3.19})) and
$\Phi=(\Phi^j)_{j=1,\ldots ,n}$. We have shown that our initial system
(\ref{eqIV.3.10}) is equivalent in $E$ to the equation $\Phi(a,b)=(a,b)$. We
are going to show that this equation has a unique solution in $E$ by using the
fixed point theorem.

(i) $\Phi(E)\subset E$.

We have $2\,\v\theta\v\big(1-\frac{\v\theta\v}{\bra\theta\ket}\big)\leq \frac
1{\bra\theta\ket}$, $2\,\frac{\v\theta\v^2}{\bra\theta\ket}\,\chi_1(\theta)\leq
\frac 8{\bra\theta\ket}$ and by (\ref{eqIV.3.20}), (\ref{eqIV.3.18}),
(\ref{eqIV.3.6}) we see that $\v F^j_3\v+\v F^j_4\v\leq C\,\varepsilon$.
Moreover we deduce from (\ref{eqIV.3.17}) and (\ref{eqIV.3.20}) that $\v
H^j_3\v+\v H^j_4\v\leq C(\mu_0)$. Finally in $E$ we have $\v b\v\leq
\frac{2\,\v y\v}{\bra\theta\ket}\leq 2\delta$. It follows then from
(\ref{eqIV.3.19}) that
$$
\v\Phi_1(a,b)\v\leq \frac{19}2\,\frac{\v
y\v}{1+4\theta^2}+C\,\varepsilon\,\v\alpha-\beta\v+C'\,\varepsilon\,\frac{\v
y\v}{\bra\theta\ket}+C(\mu_0)\,\delta\,\frac{\v y\v}{\bra\theta\ket}\,.
$$
Now using (\ref{eqIV.3.6}) (i) we see that when $(a,b)$ belongs to $E$ we have
$$
\v\beta(\theta,y,\alpha,-a)-\alpha\v\leq 10 \Big (\v a\v+\frac{\v
y\v}{\bra\theta\ket}\Big )\leq 110\, \frac{\v y\v}{\bra\theta\ket}\,.
$$
Therefore
$$
\v\Phi_1(a,b)\v\leq \Big
(\frac{19}2+110\,C\,\varepsilon+C'\,\varepsilon+C(\mu_0)\,\delta\Big)\,\frac{\v
y\v}{\bra\theta\ket}\leq \frac{10\,\v y\v}{\bra\theta\ket}\,,
$$
if $\varepsilon$ and $\delta$ are small enough.

By the same estimates we obtain,
$$
\Big\v\Phi_2(a,b)+\frac{K(\theta)}{1+4\theta^2}\,y\Big\v\leq
C\,\varepsilon\,\v\beta-\alpha\v+C\,\varepsilon\,\frac{\v
y\v}{\bra\theta\ket}+C(\mu_0)\,\delta\,\frac{\v y\v}{\bra\theta\ket}\,;
$$
so if $\varepsilon$ and $\delta$ are small enough we can bound the right hand
side by $\sqrt{\delta}\,\frac{\v y\v}{\bra\theta\ket}$. This shows that $\Phi$
maps $E$ to $E$.

(ii) $\Phi$ is a contraction.

Let now $(a_1,b_1)$, $(a_2,b_2)$ be two points in $E$. Then
\begin{equation*}
\begin{split}
\v\Phi(a_1,b_1)-\Phi(a_2,b_2)\v\leq \som^n_{j=1} \som^2_{\ell
=1}\big\v\underbrace{ Z^j_\ell (\theta,\beta(\theta,y,\alpha,-a_1))-Z^j_\ell
(\theta,\beta(\theta,y,\alpha,-a_2))}_{(1)}\big\v\\
+\som^n_{j=1} \som^4_{\ell =3} \Big\{\underbrace{\v F^j_\ell
(\theta,y,\alpha,a_1)\cdot b_1-F^j_\ell (\theta,y,\alpha,a_2)\cdot
b_2\v}_{(2)}\\
+\v\underbrace{H^j_\ell (\theta,y,\alpha,a_1,b_1)\,b_1\cdot
b_1-H^j_\ell (\theta,y,\alpha,a_2,b_2)\,b_2\cdot b_2\v}_{(3)}\Big\}\,.
\end{split}
\end{equation*}
Using (\ref{eqIV.3.20}) and Proposition \ref{pIII.3.2} we can write
$$
(1)\leq C\,\varepsilon
\,\v\beta(\theta,y,\alpha,-a_1)-\beta(\theta,y,\alpha,-a_2)\v\,.
$$
Then  (\ref{eqIV.3.6}) gives $(1)\leq C'\,\varepsilon\,\v a_1-a_2\v$.

To handle the therm (2) we use (\ref{eqIV.3.18}), (\ref{eqIV.3.20}),
(\ref{eqIV.3.6}) and Proposition \ref{pIII.3.2}. We obtain
$$
(2)\leq C\,\varepsilon(\v a_1-a_2\v+\v b_1-b_2\v)\,.
$$
Finally using (\ref{eqIV.3.17}), (\ref{eqIV.3.18}), (\ref{eqIV.3.20}) and the
fact that in $E$ we have $\v b\v\leq 2\,\frac{\v y\v}{\bra\theta\ket}\leq
2\delta$ we see easily that
$$
(3)\leq C\,(\varepsilon+\delta)(\v a_1-a_2\v+\v b_1-b_2\v)\,.
$$
It follows then that
$$
\v\Phi(a_1,b_1)-\Phi(a_2,b_2)\v\leq C\,(\varepsilon+\delta)(\v a_1-b_1\v+\v
a_2-b_2\v)
$$
where $C$ is an absolute constant depending only on the dimension and a finite
number of $A_\ell $ appearing in (\ref{eqIII.1.1}). Thus we can take
$\varepsilon$ and $\delta$ so small that $C\,(\varepsilon+\delta)<1$.

Therefore we can apply the fixed point theorem to solve (\ref{eqIV.3.19})
which  is equivalent to (\ref{eqIV.3.10}). This proves the claims (i) and (ii)
in Theorem \ref{tIV.3.1}.

Let us now prove the point (iii). We state a Lemma.
\begin{lemma}\sl\label{lIV.3.3}
 There exists $C_0>0$ such that for every $A\in \N^n$ there exist $C_A\geq 0$,
$C'_A\geq 0$ such that with $\beta$ defined in Lemma \ref{lIV.3.2},
\begin{itemize}
\item[a)]\quad  $\displaystyle{\big\v\partial
^A_y[\beta(\theta,y,\alpha,-a(\theta,y,\alpha))]\big\v\leq C_0\,\v\partial
^A_y\,a(\theta,y,\alpha)\v+\frac{C_A}{\bra\theta\ket^{\v A\v}}}$,
\item[b)] \quad $\displaystyle{\v\partial
^A_y\,a(\theta,y,\alpha)\v+\v\partial ^A_y\,b(\theta,y,\alpha)\v\leq
\frac{C'_A}{\bra\theta\ket^{\v A\v}}}$,
for all $(\theta,y)$ in $\tilde \Omega_\delta$.
\end{itemize}
\end{lemma}

{\bf Proof } We shall use an induction on $\v A\v$, starting with the formulas
(\ref{eqIV.3.19}). But before we need some intermediate results. We introduce
the following space of functions.

Let $f=f(\theta,y,\alpha)$ be a function from $\tilde
\Omega_\delta\times\R^{2n}$ to $\C$. We shall say that $f\in\CG_\pm$ if we can
write
\begin{equation}\label{eqIV.3.21}
f(\theta,y,\alpha)=G(\theta,\beta(\theta,y,\alpha,-a(\theta,y,\alpha)))
\end{equation}
where $G : \R^{\pm}_\theta\times\R^{2n}_X\rightarrow \C$ is smooth in $X$ and
satisfies
\begin{equation}\label{eqIV.3.22}
\sup_{\R^\pm\times\R^{2n}}\, \v\partial ^\gamma_X\,G(\theta,X)\v\leq
C_\gamma\,\varepsilon\,,\quad \forall \gamma\in\N^n\,.
\end{equation}
For example Proposition \ref{pIII.3.2} shows that the functions
$(\theta,y,\alpha)\mapsto
z_j(\theta,\beta(\theta,y,\alpha,-a(\theta,y,\alpha))$ (and $\zeta_j$) belong
to $\CG_\pm$ if $\alpha\in\CS_\pm$. (Here we have the sign $+$ if
$\alpha\in\CS_+$ and $-$ if $\alpha\in\CS_-$). 

Then we have the following claim.

{\bf Claim 1 } For all $\nu\in\N^n$, $\v\nu\v\geq 1$, $j=1,\ldots ,n$ we have
\begin{equation}\label{eqIV.3.23}
\partial
^\nu_\eta\,\beta^j_x(\theta,y,\alpha,\eta)=G^j_\nu(\theta,\beta(\theta,y,\alpha,\eta))
\end{equation}
where $G^j_\nu$ has all derivatives uniformly bounded on $\R^\pm\times
\R^{2n}$. The same is true for
$\partial^
\nu_\eta\,\beta^j_\xi$.

{\bf Proof of the claim } We proceed by induction on $\v\nu\v$ beginning with
$\v\nu\v=1$. Let us set
$k(\theta)=\bra\theta\ket\,\chi_1(\theta)+1-\chi_1(\theta)$. Then $1\leq
k(\theta)\leq \sqrt 5$ since $\v\theta\v\leq 2$ on $\supp \chi_1$. It follows
from (\ref{eqIV.3.5}) that for fixed $k$ in $\{1,2,\ldots ,n\}$ we have
\begin{equation*}
\begin{aligned}
\frac{\partial \beta^j_x}{\partial
\eta_k}\,(\theta,y,\alpha,\eta)&=-\frac{2\theta}{\bra\theta\ket}\,k(\theta)\,
\delta_{jk}-\som^n_{\ell =1} \Big(\frac{\partial z_j}{\partial x_\ell
}\,(\theta,\beta(\theta,y,\alpha,\eta))\,\frac{\partial \beta^\ell _x}{\partial
\eta_k}\,(\theta,y,\alpha,\eta)\\
&\enskip \enskip  +\frac{\partial z_j}{\partial \xi_\ell
}\,(\theta,\beta(\theta,y,\alpha,\eta))\,\frac{\partial \beta^j_\xi}{\partial
\eta_k}\,(\theta,y,\alpha,\eta)\Big)\\
\frac{\partial \beta^j_\xi}{\partial
\eta_k}\,(\theta,y,\alpha,\eta)&=\frac{k(\theta)}{\bra\theta\ket}\,
\delta_{jk}-\som^n_{\ell =1} \Big(\frac{\partial \zeta _j}{\partial x_\ell
}\,(\theta,\beta(\theta,y,\alpha,\eta))\,\frac{\partial \beta^\ell _x}{\partial
\eta_k}\,(\theta,y,\alpha,\eta)\\
&\enskip \enskip  +\frac{\partial \zeta_j}{\partial \xi_\ell
}\,(\theta,\beta(\theta,y,\alpha,\eta))\,\frac{\partial \beta^\ell
_\xi}{\partial\eta_k}\,(\theta,y,\alpha,\eta)\Big)\,.
\end{aligned}
\end{equation*}
Let us set $X_j=\frac{\partial \beta^j_x}{\partial \eta_k}$,
$\Xi_j=\frac{\partial \beta^j_\xi}{\partial \eta_k}$, $U=\begin{pmatrix}
X_1\\\Xi_1\\ \vdots\\ X_n\\ \Xi_n\end{pmatrix}$, $M_j=\begin{pmatrix}
\frac{\partial z_j}{\partial x_1} & \frac {\partial z_j}{\partial \xi_1}&\ldots
& \frac{\partial z_j}{\partial x_n} & \frac{\partial z_j}{\partial \xi_n} \\
\frac{\partial \zeta_j}{\partial x_1} & \frac {\partial \zeta_j}{\partial
\xi_1}&  & \frac{\partial \zeta_j}{\partial x_n} & \frac{\partial
\zeta_j}{\partial \xi_n}\end{pmatrix}$, $M=\begin{pmatrix} M_1\\M_2\\\vdots\\M_n
\end{pmatrix}$ and $F(\theta)=\frac
{k(\theta)}{\bra\theta\ket}\,\begin{pmatrix} -2\theta\,\delta_{1k}\\
\delta_{1k}\\ \vdots
\\-2\theta\,\delta_{nk}\\\delta_{nk}\end{pmatrix}$.

Then the above equation imply that
$$
U(\theta,y,\alpha,\eta)=F(\theta)-M(\theta,\beta(\theta,y,\alpha,\eta))\,U(\theta,y,\alpha,\eta)\,.
$$
Now by Proposition \ref{pIII.3.2} the entries of the matrix $M$ are $\CO
(\varepsilon)$. It follows that the matrix $I+M$ is invertible. Therefore we
obtain
$$
U(\theta,y,\alpha,\eta)=(I+M(\theta,\beta(\theta,y,\alpha,\eta)))^{-1}\,F(\theta)\,.
$$
This proves the case $\v\nu\v=1$.

Assume now that our claim is proved for $\v\nu\v\leq N-1$. Then
$$
\partial
^\nu_\eta\,\beta^j_x(\theta,y,\alpha,\eta)=G^j_\nu(\theta,\beta(\theta,y,\alpha,\eta))\,.
$$
Then for $k=1,\ldots ,n$
\begin{equation*}
\begin{aligned}
\frac\partial {\partial \eta_k}\,\partial
^\nu_\eta\,\beta^j_x(\theta,y,\alpha,\eta)&=\som^n_{\ell =1} \,\frac{\partial
G^j_\nu}{\partial X_\ell
}\,(\theta,\beta(\theta,y,\alpha,\eta))\,\frac{\partial \beta^\ell _x}{\partial
\eta_k}\,(\theta,y,\alpha,\eta)\\
&\enskip\enskip  +\som^{2n}_{\ell =n+1}\, \frac{\partial G^j_\nu}{\partial
X_\ell }\,(\theta,\beta(\theta,y,\alpha,\eta))\,\frac{\partial \beta^\ell
_\xi}{\partial \eta_k}\,(\theta,y,\alpha,\eta)\,.
\end{aligned}
\end{equation*}
Using (\ref{eqIV.3.23}) with $\v\nu\v=1$ we obtain the claim up to the order
$\v\nu\v=N$.
\cqfd

\begin{consequence} \sl  \label{cIV.3.4}
With the notations of (\ref{eqIV.3.18}), (\ref{eqIV.3.20}) and
(\ref{eqIV.3.21}) The functions $(\theta,y,\alpha)\mapsto Z^j_\ell
(\theta,\beta(\theta,y,\alpha,-a(\theta,y,\alpha))$, $\ell =1,2$, $F^j_p$,
$p=3,4$ belong to $\CG_+.$
\end{consequence}

Let us now go back to the proof of Lemma \ref{lIV.3.3}. We begin by the case
$\v A\v=1$. For convenience we shall set
\begin{equation}\label{eqIV.3.24}
f(y)=\beta(\theta,y,\alpha,-a(\theta,y,\alpha))\,.
\end{equation}
It follows from (\ref{eqIV.3.5}) that
\begin{equation*}
\begin{aligned}
\frac{\partial f^j_x}{\partial y_k}(y)&=\Big
(1-\frac{\v\theta\v}{\bra\theta\ket}\,(1-\chi_1(\theta))\Big)\,
\delta_{jk}+\frac{2\theta}{\bra\theta\ket}\,k(\theta)\,\frac{\partial
a_j}{\partial y_k}-\som^n_{\ell =1} \Big (\frac{\partial z_j}{\partial x_\ell
}\, (\theta,f(y))\,\frac{\partial f^j_x}{\partial y_k}\,(y)\\
&\hbox to 6cm{} -\frac{\partial z_j}{\partial \xi_\ell
}\,(\theta,f(y))\,\frac{\partial f^j_\xi}{\partial y_k}\,(y)\Big)\\
\frac{\partial f^j_\xi}{\partial y_k}\,(y)&=\frac 12\,(1-\chi_1(\theta))\,
\frac{\sgn \theta}{\bra\theta\ket}\,\delta_{jk}-\frac
1{\bra\theta\ket}\,k(\theta)\,\frac{\partial a_j}{\partial y_k}-\som^n_{
\ell =1} \Big (\frac{\partial \zeta_j}{\partial x_\ell
}\,(\theta,f(y))\,\frac{\partial f^j_x}{\partial y_k}\,(y)\\
&\hbox to 6cm{} -\frac{\partial \zeta_j}{\partial \xi_\ell
}\,(\theta,f(y))\,\frac{\partial f^j_\xi}{\partial y_k}\,(y)\Big)\,.
\end{aligned}
\end{equation*}
It follows from Proposition \ref{pIII.3.2} that,
$$
\Big\v\frac{\partial f}{\partial y_k}\,(y)\Big\v\leq C \Big (\frac
1{\bra\theta\ket}+\Big\v\frac{\partial a}{\partial y_k}\,(\theta,y,\alpha)\Big\v
+\varepsilon\,\Big\v\frac{\partial f}{\partial y_k}\,(y)\Big\v\Big)
$$
where $C$ depends only on the constants $A_0$, $A_1$ appearing in
(\ref{eqIII.1.1}). Therefore taking $\varepsilon$ so small that
$C\,\varepsilon\leq \frac 12$ we obtain the point a) in Lemma \ref{lIV.3.3} for
$\v A\v=1$. Let us prove the point b). First of all we deduce from a) and the
Consequence \ref{cIV.3.4} that
\begin{equation}\label{eqIV.3.25}\left\{
\begin{array}{l}
\frac{\partial }{\partial y_k}\,[Z^j_\ell (\theta,f(y))]\,,\enskip 
\ell =1,2\,,\enskip \frac{\partial }{\partial y_k}\, [F^j_p]\,,\enskip
p=3,4\\
 \textrm{are bounded by}\enskip C\,\varepsilon \big (\frac
1{\bra\theta\ket}+\big\v\frac{\partial a}{\partial
y_k}\,(\theta,y,\alpha)\big\v\big)\,.
\end{array}\right.
\end{equation}
Now we claim that for $p=3,4$,
\begin{equation}\label{eqIV.3.26}
\Big\v\frac\partial {\partial y_k}\,[H^j_p(\theta,y,\alpha,a,b)]\Big\v\leq
C\Big(\Big\v\frac{\partial a}{\partial
y_k}\,(\theta,y,\alpha)\Big\v+\Big\v\frac{\partial b}{\partial y_k}\,(\theta,y,
\alpha)\Big\v+\frac 1{\bra\theta\ket}\Big)\,.
\end{equation}
Indeed the left hand side of (\ref{eqIV.3.26}) can be written
$$
\frac{\partial H^j_p}{\partial a}\cdot \frac{\partial a}{\partial
y_k}+\frac{\partial H^j_p}{\partial b}\cdot \frac{\partial b}{\partial
y_k}+\frac{\partial H^j_p}{\partial y_k}=(1)+(2)+(3)\,.
$$
Now using (\ref{eqIV.3.20}), (\ref{eqIV.3.18}) and (\ref{eqIV.3.17}) we see
that (1) and (2) can be bounded by the right hand side of (\ref{eqIV.3.26}). To
handle the term (3) we use (\ref{eqIV.3.16}) with $A=0$, $\v B\v=1$. We obtain
$$
\v(3)\v\leq C\,\som_{\v\gamma\v\leq 3n+2} \int \Big \v\frac{\partial
}{\partial y_k}\,\partial ^\gamma_\eta\,g_j(\eta)\Big\v\,d\eta\,.
$$
Now we use (\ref{eqIV.3.9}) and (\ref{eqIV.3.23}). We obtain
$$
\partial ^\gamma_\eta\,g_j(\eta)=\partial ^\gamma_\eta\Big[\chi_0\Big
(\frac
1{\mu_0}\,\eta)\Big]\Big[(\beta^j_\xi-i\,\beta^j_x)(\theta,y,\alpha,\gamma)-(\alpha^j_\xi+i\,\alpha^j_x)
\Big]+\som_{\v\gamma'\v<\v\gamma\v}\,\partial
^{\gamma'}_\eta\Big[\chi_0\Big(\frac
1{\mu_0}\,\eta\Big)\Big]\,G^j_{\gamma'}(\theta,\beta(\theta,y,\alpha,\eta))
$$
where $G^j_{\gamma'}$ satisfy (\ref{eqIV.3.22}). Since by (\ref{eqIV.3.5}) we
have
$$
\Big\v\frac{\partial \beta}{\partial y_k}\,(\theta,y,\alpha,\eta)\Big\v\leq
\frac C{\bra\theta\ket}\,.
$$
We obtain
$$
\v(3)\v\leq \frac C{\bra\theta\ket}
$$
which proves (\ref{eqIV.3.26}).

Now we use (\ref{eqIV.3.19}), (\ref{eqIV.3.25}), (\ref{eqIV.3.26}) and the fact
that $\v b\v\leq \frac{2\,\v y\v}{\bra\theta\ket}\leq 2\delta$. We obtain with
$a_j=a_j(\theta,y,\alpha)$, $b_j=b_j(\theta,y,\alpha)$,
$$
\Big\v\frac{\partial a_j}{\partial y_k}\Big\v+\Big\v\frac{\partial
b_j}{\partial y_k}\Big\v\leq
\frac{12}{\bra\theta\ket}+C\,(\varepsilon+\delta)\Big (\frac
1{\bra\theta\ket}+\Big\v\frac{\partial a}{\partial
y_k}\Big\v+\Big\v\frac{\partial b}{\partial y_k}\Big\v\Big)\,.
$$
Taking $\varepsilon$ and $\delta$ small enough we obtain the point b) in Lemma
\ref{lIV.3.3} when $\v A\v=1$.

Assume now that a) and b) in Lemma \ref{lIV.3.3} are true when $\v A\v\leq N$
and let $\v A\v=N+1$. It follows from the induction that
\begin{equation}\label{eqIV.3.27}
\big\v\partial ^B_y[\beta(\theta,y,\alpha,-a(\theta,y,\alpha))]\big\v\leq
\frac{C_B}{\bra\theta\ket^{\v B\v}}\,,\enskip \textrm{if}\enskip \v B\v\leq N\,.
\end{equation}

{\bf Claim 1 } If $\v A\v=N+1\geq 2$,
\begin{equation}\label{eqIV.3.28}
\big\v\partial ^A_y[\beta(\theta,y,\alpha,-a(\theta,y,\alpha))]\big\v\leq
C_0\,\big\v\partial
^A_y\,a(\theta,y,\alpha)\big\v+\frac{C_A}{\bra\theta\ket^{\v A\v}}\,.
\end{equation}
Indeed, according to (\ref{eqIV.3.5}) we have, setting for short
$f(y)=\beta(\theta,y,\alpha,-a(\theta,y,\alpha))$ and
$k(\theta)=\bra\theta\ket\,\chi_1(\theta)+1-\chi_1(\theta)$,
\begin{equation*}\left\{
\begin{array}{l}
f_x(y)=\alpha_x+\big
(1-\frac{\v\theta\v}{\bra\theta\ket}\,(1-\chi_1(\theta))\big)\,y+\frac{2\theta}{\bra
\theta\ket}\,k(\theta)\,a(\theta,y,\alpha)+z(\theta,\alpha)-z(\theta,f(y))\\
f_\xi(y)=\alpha_\xi+ \frac 12\,
(1-\chi_1(\theta))\,\frac{\sgn \theta}{\bra\theta\ket}\,  
y-\frac{k(\theta)}{\bra
\theta\ket}\,a(\theta,y,\alpha)+\zeta(\theta,\alpha)-\zeta(\theta,f(y))\,.
\end{array}\right.
\end{equation*}
Differentiating both size $A$ times with respect to $y$ we obtain since
$\v A\v\geq 2$,
$$
\v\partial ^A_y\,f(y)\v\leq 5\,\v\partial ^A_y\,a(\theta,y,\alpha)\v+\v\partial
^A_y[z(\theta,f(y))]\v+\v\partial ^A_y[\zeta(\theta,f(y))]\v\,.
$$
We use now the Faa di Bruno formula (see Appendix VIII.1). Let $Z$ be $z$ or
$\zeta$ then
$$
\partial ^A_y[Z(\theta,f(y))]=\som^n_{j=1} \Big\{\underbrace{\frac{\partial
Z}{\partial x_j}\,(\theta,f(y))\,\partial ^A_y\,f^j_x(y)+\frac{\partial
Z}{\partial
\xi_j}\,(\theta,f(y))\,\partial ^A_y\,f^j_\xi(y)}_{(1)}\Big\}+(2)
$$
where (2) is a finite sum of terms of the form
$$
(\partial ^\nu_X\,Z)(\theta,f(y)) \prod^s_{j=1} \big (\partial ^{\ell
_j}_y\,f(y)\big)^{k_j}
$$
where $X=(x,\xi)$, $2\leq \v\nu\v\leq \v A\v$, $\v\ell _j\v\geq 1$, $\v
k_j\v\geq 1$ and
$$
\sum^s_{j=1}\,k_j=\nu\,,\quad \sum^s_{j=1}\,\v k_j\v\,\ell _j=A\,.
$$
The term (1) can be bounded by $C_0\,\varepsilon\,\v\partial ^A_y\,f(y)\v$
(where $C_0$ depends on $A_0,A_1$ in (\ref{eqIII.1.1})). Since $\v\nu\v\geq 2$
it is easy to see that $\v\ell _j\v\leq \v A\v-1$. We can therefore use
(\ref{eqIV.3.27}) and Proposition \ref{pIII.3.2} to write
$$
\v(2)\v\leq C_A\,\varepsilon \prod^s_{j=1}\, \frac 1{\bra\theta\ket^{\v\ell
_j\v\,\v k_j\v}}\leq \frac{C_A}{\bra\theta\ket^{\v A\v}}\,.
$$
Thus (\ref{eqIV.3.28}) is proved which implies the part a) of Lemma
\ref{lIV.3.3} when $\v A\v=N+1$.

{\bf Claim 2 } If $F\in\CG_\pm$ (see (\ref{eqIV.3.21})) and $\v A\v=N+1$ we have
\begin{equation}\label{eqIV.3.29}
\v\partial ^A_y\,F(\theta,y,\alpha)\v\leq \varepsilon \Big (C_0\,\v\partial
^A_y\,a(\theta,y,\alpha)\v+\frac{C_A}{\bra\theta\ket^{\v A\v}}\Big)\,.
\end{equation}
Let us set for convenience as in (\ref{eqIV.3.24}),
$$
f(y)=\beta(\theta,y,\alpha,-a(\theta,y,\alpha))\,.
$$
We know by assumption that $F(\theta,y,\alpha)=G(\theta,f(y))$ where $G$
satisfies (\ref{eqIV.3.22}). By the Faa di Bruno formula we have
\begin{equation}\label{eqIV.3.30}
\partial ^A_y\, F(\theta,y,\alpha)=\som^{2n}_{i=1}\underbrace{\frac{\partial
G}{\partial X_i}\,(\theta,f(y))\,\partial ^A_y\,f(y)}_{(1)}+(2)
\end{equation}
where $(2)$ is a finite sum of terms of the form
$$
(\partial ^\nu_X\,G)(\theta,f(y))\,\prod^s_{j=1}\,\big (\partial ^{\ell
_j}_y\,f(y)\big)^{k_j}
$$
where $2\leq \v\nu\v\leq \v A\v$, $\v\ell _j\v\geq 1$, $\v k_j\v\geq 1$, $1\leq
s\leq \v A\v$ and
\begin{equation}\label{eqIV.3.31}
\sum^s_{j=1}\,k_j=\nu\,,\quad \sum^s_{j=1}\,\v k_j\v\,\ell _j=A\,.
\end{equation}
Now by the Claim 1 we have
\begin{equation}\label{eqIV.3.32}
\v(1)\v\leq \varepsilon \Big (C_0\,\v\partial
^A_y\,a(\theta,y,\alpha)\v+\frac{C_A}{\bra\theta\ket^{\v A\v}}\Big)\,.
\end{equation}
On the other hand in the term (2) it is easy to see that $\v\ell _j\v\leq \v
A\v-1$. Indeed if we had a $j_0$ such that $\v\ell _{j_0}\v=\v A\v$ it would
follow from (\ref{eqIV.3.31}) that $j_0=1$, $s=1$ and $\v k_1\v=1$~; but then
$\v k_1\v=1=\v \nu\v$ which is in contradiction with $\v\nu\v\geq 2$. Therefore
we can use (\ref{eqIV.3.27}), (\ref{eqIV.3.31}) to write,
$$
\v(2)\v\leq C_A\,\varepsilon\,\prod^s_{j=1}\, \Big (\frac
1{\bra\theta\ket^{\v\ell _j\v}}\Big)^{\v k_j\v}=C_A\,\varepsilon\, \frac
1{\bra\theta\ket^{\v A\v}}\,.
$$
Then (\ref{eqIV.3.29}) follows from (\ref{eqIV.3.30}) and (\ref{eqIV.3.32}).

{\bf Claim 3 } If $\v A\v=N+1$, $j=1,\ldots ,n$, $\ell =3,4$ we have,
\begin{equation}\label{eqIV.3.33}
\big\v\partial ^A_y(H^j_\ell (\theta,y,\alpha,a,b))\big\v\leq C_0 \big
(\v\partial ^A_y\,a(\theta,y,\alpha)\v+\v\partial
^A_y\,b(\theta,y,\alpha)\v\big)+\frac{C_A}{\bra\theta\ket^{\v A\v}}\,,
\end{equation}
where $H^j_\ell $ is defined in (\ref{eqIV.3.20}), (\ref{eqIV.3.18}),
(\ref{eqIV.3.13}).

The proof is exactly the same as in the Claim 2. We use the Faa di Bruno
formula, the estimates on $a,b$ given by the induction, the estimate
(\ref{eqIV.3.16}) to obtain (\ref{eqIV.3.33}). We are ready now to prove the
part b) of Lemma \ref{lIV.3.3} when $\v A\v=N+1$.

We use the equations (\ref{eqIV.3.19}), (\ref{eqIV.3.20}) which we
differentiate $\v A\v$ times with respect to $y$. Since $Z^j_\ell
(\theta,\beta)$ and $F^j_k$ belong to $\CG_\pm$ we use the Claim 2 to estimate
them~; the term $H^j_\ell \,b\cdot b$ is handled by (\ref{eqIV.3.33}), the
Leibniz formula and the induction hypothesis. Finally we obtain since $\v
b\v\leq
\frac{2\,\v y\v}{\bra\theta\ket}\leq 2\delta$,
$$
\v\partial ^A_y\,a_j\v\leq (\varepsilon+\delta)\,C_0\big (\v\partial
^A_y\,a\v+\v\partial ^A_y\,b\v\big)+\frac{C_A}{\bra\theta\ket^{\v A\v}}\,.
$$
Taking $\varepsilon$ and $\delta$ small enough we obtain the part b) of Lemma
\ref{lIV.3.3}. \cqfd

So far we have proved the points (i), (ii), (iii) in Theorem \ref{tIV.3.1}. 

Let us prove (iv). It follows from (\ref{eqIV.2.1}), (\ref{eqIV.3.6}) and
(\ref{eqIV.3.9}) that
$$
\Big\v\frac{\partial q_k}{\partial \eta_\ell
}\,(\eta,a,b,g_j\Big\v+\Big\v\frac{\partial q_k}{\partial b_\ell
}\,(\eta,a,b,g_j)\Big\v\leq C
$$
where $C$ is an absolute constant.

Since $\v\eta\v\leq \delta$, $\v a\v\leq \frac{10\,\v y\v}{\bra\theta\ket}\leq
10\,\delta$, $\v b\v\leq 2\,\frac{\v y\v}{\bra\theta\ket}\leq 2\delta$ we can
write
\begin{equation}\label{eqIV.3.34}
\big\v q_k(\eta,a,b,g_j)-q_k(-a,a,0,g_j)\big\v\leq C'\,\delta\,.
\end{equation}
Now (\ref{eqIV.2.8}) gives $q_k(-a,a,0,g_j)=\frac{\partial g_j}{\partial
\eta_k}\,(-a)$. It follows then from (\ref{eqIV.3.15}), (\ref{eqIV.3.6}) and
Proposition \ref{pIII.3.2} that
$$
\Big\v
q_k(-a,a,0,g_j)-(1+2i\theta)\,\frac{k(\theta)}{\bra\theta\ket}\,\delta_{jk}
\Big\v\leq C\,\varepsilon
$$
which combined with (\ref{eqIV.3.34}) gives the point (iv).

Finally (v) follows easily from (\ref{eqIV.2.2}), (\ref{eqIV.3.9}) and
(\ref{eqIV.3.6}). This ends the proof of Theorem \ref{tIV.3.1}. \cqfd
\begin{corollary}\sl \label{cIV.3.5}
Let us set $k(\theta)=\bra\theta\ket\,\chi_1(\theta)+(1-\chi_1(\theta))$ and
for $j=1,\ldots ,n$,
$$
\tilde g_j(\eta)=\chi_0 \Big (\frac{\eta-\frac 12 (1-\chi_1(\theta))\,(\sgn
\theta)\,y}{\mu_0\,k(\theta)}\Big)\,v_j\Big(\theta,y+x(\theta,\alpha),\,\frac\eta{\bra\theta\ket}
+\xi(\theta,\alpha),\alpha\Big)\,.
$$
Then we can write
\begin{equation}\label{eqIV.3.35}
\tilde g_j(\eta)=\som^n_{\ell =1}\,\tilde q_\ell (\eta,\tilde a,\tilde b,\tilde
g_j)\Big(\eta_\ell -\frac 12\,(1-\chi_1(\theta))\,(\sgn \theta) y+(\tilde
a_\ell +i\,\tilde b_\ell)(\theta,y,\alpha)\Big) 
\end{equation}
where
$$
\tilde q_\ell (\eta,\tilde a,\tilde b,\tilde g_j)=\frac 1{k(\theta)}\,q_\ell
\Big (\frac{\eta-\frac 12\,(1-\chi_1(\theta))\,(\sgn \theta) y}{k(\theta)}\,,
\frac{\tilde a}{k(\theta)}\,, \frac{\tilde b}{k(\theta)}\,, \tilde g_j\Big)
$$
and $\tilde q_\ell $, $\tilde a_\ell $, $\tilde b_\ell $ satisfy
\begin{itemize}
\item[(i) ] $\displaystyle{\v\tilde a(\theta,y,\alpha)\v\leq 10\,\sqrt
5\,\frac{\v y\v}{\bra\theta\ket} }$

$\displaystyle{\big\v\tilde
b(\theta,y,\alpha)+\frac{\bra\theta\ket}{1+4\theta^2}\,y\big\v}\leq
\sqrt{5\delta}\, \frac{\v y\v}{\bra\theta\ket}$.
\item[(ii) ] $\displaystyle{\v\partial ^A_y\,\tilde
a(\theta,y,\alpha)\v+\v\partial ^A_y\,\tilde b(\theta,y,\alpha)\v\leq
\frac{C_A}{\bra\theta\ket^{\v A\v}}}$,\enskip  $A\in\N^n$.
\item[(iii) ] $\displaystyle{\big\v\tilde q_\ell (\eta,\tilde a,\tilde b,\tilde
g_j)-\frac{(1+2i\theta)}{\bra\theta\ket}\,\delta_{jk}\big\v\leq
C\,(\varepsilon+\delta)}$\enskip  if\enskip  $\v\eta\v\leq \delta$.
\item[(iv) ] $\v\partial ^A_{(a,b)}\,\partial ^B_\eta\,\tilde q_\ell
(\eta,\tilde a,\tilde b,\tilde g_j)\v\leq C(\mu_0)$ \enskip if \enskip $\v
A\v+\v B\v\geq 1$,\enskip 
$\v\eta\v\leq \mu_0$,\enskip  $1\leq j,\ell \leq n$

uniformly with respect to $(\theta,y)\in\tilde \Omega_\delta$ and
$\alpha\in\CS$.
\end{itemize}
\end{corollary}

{\bf Proof } We have  $\chi_1(\theta)\,\eta+(1-\chi_1(\theta))
\big[\frac\eta{\bra\theta\ket} +\frac
12\,\frac{1-\chi_1(\theta)}{\bra\theta\ket}\,(\sgn \theta)\,
y\big]=\frac{k(\theta)}{\bra\theta\ket}\,\eta+\frac 12\,
\frac{1-\chi_1(\theta)}{\bra\theta\ket}\, (\sgn \theta)\,y$. So let us set in
the statement of Theorem \ref{tIV.3.1} $\tilde \eta=k(\theta)\,\eta+\frac
12\,(1-\chi_1(\theta))\,(\sgn \theta)\,y$~; then we obtain the decomposition of
$\tilde g_j$ in Corollary \ref{cIV.3.5} with $\tilde a_\ell =k(\theta)\,a_\ell
$, $\tilde b_\ell =k(\theta)\,b_\ell $ and the estimates on $\tilde q_\ell $,
$\tilde a_\ell $, $\tilde b_\ell $ follow easily from the correspondent one for
$q_\ell $, $a_\ell $, $b_\ell $ stated in Theorem \ref{tIV.3.1}. \cqfd

{\bf Lagrangian ideals and the phase equation}

We pursue here the proof of Theorem \ref{tIV.1.2}. Let us set
\begin{equation}\label{eqIV.3.36}
\CO=\Big\{(\theta,y,\eta)\in\R\times\R^n\times\R^n : \v
y\v<\delta\,\bra\theta\ket\,,\, \frac{\v\eta-\frac
12\,(1-\chi_1(\theta))\,(\sgn \theta)\,y\v}{k(\theta)}<\delta\Big\}\,.
\end{equation}
We introduce now a space of families $(f(\cdot ,\alpha))_{\alpha\in\CS}$ of
functions on $\CO$.
\begin{definition}\sl \label{dIV.3.6}
We say that $(f(\cdot ,\alpha))_{\alpha\in\CS}$ belongs to $\CH$ if
\begin{itemize}
\item[(i) ] For all $\alpha$ in $\CS$, $(\theta,y,\eta)\mapsto
f(\theta,y,\eta,\alpha)$ belongs to $C^\infty (\CO)$.
\item[(ii) ] For every $A,B$ in $\N^n$ there exists $C_{AB}>0$ independent of
$\alpha$ such that 

$\displaystyle{\sup_{(\theta,y,\eta)\in\CO} \v\partial
^A_y\,\partial ^B_\eta\,f(\theta,y,\eta,\alpha)\v\leq C_{AB}}$, for all
$\alpha\in\CS$.
\end{itemize}
\end{definition}
\begin{remark}\sl \label{rIV.3.7}

1)  $\CH$ is closed under multiplication and derivation with respect to
$(y,\eta)$.

2)  If we set, with notation (\ref{eqIV.2.10})
$$
f(\theta,y,\eta,\alpha)=v_j(\theta,y+x(\theta,\alpha)\,,\,\frac\eta{\bra\theta\ket}
+\xi(\theta,\alpha),\alpha)
$$
then $(f(\cdot ,\alpha))_{\alpha\in\CS}$ belongs to $\CH$. This follows from
(\ref{eqIV.3.6}).
\end{remark}
\begin{definition}\sl \label{dIV.3.8}
Let $F=(F(\cdot ,\alpha))_{\alpha\in\CS}$. We say that $F\in \CJ$ if we can
write
$$
F(\theta,y,\eta,\alpha)=\sum^n_{j=1}
f_j(\theta,y,\eta,\alpha)\,v_j\Big(\theta,y+x(\theta,\alpha)\,,
\frac\eta{\bra\theta\ket}+\xi(\theta,\alpha),\alpha\Big )
 $$
for all $(\theta,y,\eta,\alpha)$ in $\CO\times\CS$ with $(f(\cdot
,\alpha))_{\alpha\in\CS}\in\CH$.
\end{definition}
\begin{example}\sl \label{eIV.3.9}
Let us set $F=(F(\cdot ,\alpha))_{\alpha\in\CS}$ where
$$
F(\theta,y,\eta,\alpha)=\eta_k-\frac 12\, (1-\chi_1(\theta)\,(\sgn
\theta)\,y_k+(\tilde a_k+i\,\tilde b_k)(\theta,y,\alpha)\,,\quad
k\in\{1,2,\ldots ,n\}\,.
$$
Then $F\in\CJ$.
\end{example}

This follows from Corollary \ref{cIV.3.5}. Indeed the matrix $(\tilde q_\ell
(\eta,\tilde a,\tilde b,\tilde g_j))_{j,\ell }$ is invertible if
$\varepsilon+\delta$ is small enough and according to the estimate (v) if we
set $(d_{\ell j})=(\tilde q(\cdots ))^{-1}$ then $(d_{jk}(\cdot
,\alpha))_{\alpha\in\CS}$ belongs to $\CH$ so our claim follows from Remark
\ref{rIV.3.7}. 

Now if $F=(F(\cdot ,\alpha))_{\alpha\in\CS}$, $G=(G(\cdot
,\alpha))_{\alpha\in\CS}$ we set
\begin{equation}\label{eqIV.3.37}
\{F,G\}=\Big(\sum^n_{j=1} \Big (\frac{\partial F}{\partial
\eta_j}\,\frac{\partial G}{\partial y_j}-\frac{\partial F}{\partial
y_j}\,\frac{\partial G}{\partial \eta_j}\Big )(\cdot
,\alpha)\Big)_{\alpha\in\CS}\,.
\end{equation}
Then we have the following result.
\begin{lemma}\sl\label{lIV.3.10}
$\CJ$ is closed under the Poisson bracket (\ref{eqIV.3.37}).
\end{lemma}

{\bf Proof } Recall (see (\ref{eqIV.2.10})) that
$v_j(\theta,x,\xi,\alpha)=u_j\circ \chi_{-\theta}(x,\xi)$ where
$u_j(x,\xi,\alpha)=\xi_j-\alpha^j_\xi-i(x_j-\alpha^j_x)$ and
$\chi_{-\theta}(x,\xi)=(x(-\theta;x,\xi),\,\xi(-\theta;x,\xi))$ is the
symplectic map defined by the flow. Since $\{u_j,u_k\}=0$ we have, denoting by
$\{\kern 5pt ,\kern 3pt \}$ the Poisson bracket in $(x,\xi)$,
$\{v_j,v_k\}=\{u_j\circ\chi_{-\theta},u_k\circ\chi_{-\theta}\}=\{u_j,u_k\}
\circ\chi_{-\theta}=0$. It follows that, in the coordinates $(y,\eta)$,
$$
\Big\{v_j(\theta;y+x(\theta;\alpha),\,
\frac\eta{\bra\theta\ket}+\xi(\theta;\alpha)),\,v_k(\theta;y+x(\theta;\alpha),\,
\frac\eta{\bra\theta\ket}+\xi(\theta;\alpha))\Big\}=0\,.
$$
Let $F=\big (\som_j f_j\,v_j(\cdot ,\alpha)\big)_{\alpha\in\CS}$,
$G=\big(\som_k g_k\,v_k(\cdot ,\alpha)\big )_{\alpha\in\CS}$ two elements of
$\CJ$. Then if $\{\kern 5pt ,\kern 3pt \}$ is the Poisson bracket in
$(y,\eta)$, we have,
$$
\Big\{\sum_j f_j\,v_j,\,\sum_k g_k\,v_k\Big\}=\sum_k\Big (\sum_j
f_j\{v_j,g_k\}\Big )\,v_k+\sum_k \Big (\sum_j \{f_j,g_k\}\,v_j\Big )\,
v_k+\sum_j \Big (\sum_k \{f_j,v_k\}\,g_k\Big )\,v_j\,.
$$
Since $f_j, g_k, v_j$ belong to $\CH$ it follows from Remark \ref{rIV.3.7} (i)
and Definition \ref{dIV.3.8} that $\{F,G\}\in\CJ$.

 \cqfd

Here is an important lemma which is a generalization to our context of Lemma
7.5.10 of \cite{H}.

According to Corollary \ref{cIV.3.5}, we shall set
\begin{equation}\label{eqIV.3.38}
\psi_k(\theta,y,\alpha)=\frac 12\,(1-\chi_1(\theta))\,(\sgn
\theta)\,y_k-(\tilde a_k(\theta;y,\alpha)+i\,\tilde
b_k(\theta;y,\alpha))\,,\quad k=1,\ldots ,n\,.
\end{equation}
\begin{lemma}\sl \label{lIV.3.11}
Let $R=(R(\cdot ,\alpha))_{\alpha\in\CS}\in\CJ$ where $R(\theta,y,\alpha)$ is
independent of $\eta$. Then for every $N\in\N$ one can find $C_N>0$ such that
for every $(\theta,y)$ in $\tilde \Omega_\delta$ (see (\ref{eqIV.3.2})) and
$\alpha$ in $\CS$ we have
$$
\v R(\theta,y,\alpha)\v\leq C_N\,\v \Im \psi(\theta,y,\alpha)\v^N\,.
$$
\end{lemma}

{\bf Proof } We are going to show by induction on $N\in\N^*$ that we can write
for $(\theta,y,\alpha)$ in $\tilde \Omega_\delta\times\CS$,
\begin{equation}\label{eqIV.3.39}
R(\theta,y,\alpha)=\sum_{0<\v\gamma\v<N}
h_\gamma(\theta,y,\alpha)(\eta-\psi)^\gamma+\sum_{\v\gamma\v=N}
w_\gamma(\theta,y,\eta,\alpha)(\eta-\psi)^\gamma\,,
\end{equation}
where $(h_\gamma(\cdot ,\alpha))_{\alpha\in\CS}$ and $(w_\gamma(\cdot
,\alpha))_{\alpha\in\CS}$ belong to $\CH$.

For $N=1$ the first sum in (\ref{eqIV.3.39}) is empty and by assumption we have
$$
R(\theta,y,\alpha)=\sum^n_{j=1}
f_j(\theta,y,\eta,\alpha)\,v_j\Big(\theta;y+x(\theta;\alpha),\,\frac\eta{\bra\theta\ket}
+\xi(\theta;\alpha)\Big)\,,
$$
where $(f_j(\cdot ,\alpha))_{\alpha\in\CH}$. Now we use Corollary
\ref{cIV.3.5}. Since $\chi_0\big (\frac{\eta-\frac 1
2\,(1-\chi_1(\theta))\,(\sgn \theta)\,y}{k(\theta)\,\mu_0}\big)=1$ if
$(\theta,y,\eta)\in\CO$ we can write
$$
R(\theta,y,\alpha)=\sum^n_{k=1}
r_k(\theta,y,\eta,\alpha)(\eta_k-\psi_k(\theta,y,\alpha))
$$
where
$$
r_k(\theta,y,\eta,\alpha)=\sum^n_{j=1} f_j(\theta,y,\eta,\alpha)\,\tilde
q_k(\eta,\tilde a,\tilde b,\tilde g_j)\,.
$$
Now it follows from Corollary \ref{cIV.3.5} and Remark \ref{rIV.3.7} that
$(r_k(\cdot ,\alpha))_{\alpha\in\CS}$ belongs to $\CH$. Therefore
(\ref{eqIV.3.39}) holds when $N=1$. Assume now that (\ref{eqIV.3.39}) is true
at the level $N$. Then apply Lemma \ref{lIV.2.1} to the function
$$
g_\gamma(\theta,y,\eta,\alpha)=\chi_0\Big (\frac{\eta-\frac
12\,(1-\chi_1(\theta))\,(\sgn \theta)\,y}{\mu_0\,k(\theta)}\Big
)\,w_\gamma(\theta,y,\eta,\alpha)\,,\quad \v\gamma\v=N\,,
$$
with $z=-\psi(\theta,y,\alpha)$. It follows then that
\begin{equation}\label{eqIV.3.40}
g_\gamma(\theta,y,\eta,\alpha)=\sum^n_{k=1}
q_k(\eta,a,b,g_\gamma)(\eta_k-\psi_k(\theta,y,\alpha))+r(a,b,g_\gamma)\,.
\end{equation}
For the $q_k's$ and $r$ we have the estimates (\ref{eqIV.2.2}). If we set
\begin{equation}\label{eqIV.3.41}\left\{
\begin{array}{l}
h_\gamma(\theta,y,\alpha)=r(a(\theta,y\,\alpha),\,b(\theta,y,\alpha),\,g_\gamma(\theta,y,\cdot 
,\alpha))\\
w_{\gamma
j}(\theta,y,\eta,\alpha)=q_j(\eta,a(\theta,y,\alpha),\,b(\theta,y,\alpha),\,
g_\gamma(\theta,y,\cdot ,\alpha))\,.
\end{array}\right.
\end{equation}
We deduce from (\ref{eqIV.2.2}) and Corollary \ref{cIV.3.5} that
$(h_\gamma(\cdot ,\alpha))_{\alpha\in\CS}$ and $(w_{\gamma j}(\cdot
,\alpha))_{\alpha\in\CS}$ belong to $\CH$. Then using (\ref{eqIV.3.39}) at the
level $N$  and (\ref{eqIV.3.40}), (\ref{eqIV.3.41}) we obtain (\ref{eqIV.3.39})
at the level $N+1$.
Let us take now in (\ref{eqIV.3.39}) $\eta=\Re \psi(\theta,y,\alpha)+s\,\Im
\psi (\theta,y,\alpha)$ where $s\in[0,1]$. We obtain
\begin{equation}\label{eqIV.3.42}
\v R(\theta,y,\alpha)-\sum_{0<\v\gamma\v<N} h_\gamma(\theta,y,\alpha)(\Im
\psi(\theta,y,\alpha))^\gamma\,(s-i)^\gamma\v \leq C_N\,\v \Im
\psi(\theta,y,\alpha)\v^N
\end{equation}
where $C_N$ is independent of $(\theta,y,\alpha)\in \tilde \Omega_\delta\times
\CS$.

Using an interpolation formula we deduce that the coefficients of the
polynomial in $(s-i)$ in the left hand side of (\ref{eqIV.3.42}) satisfy the
same estimate which proves that $R$ has the claimed bound. \cqfd
\begin{corollary}\sl\label{cIV.3.12}
For every $N\in\N$ there exists a constant $C_N>0$ such that
$$
\Big\v\Big (\frac{\partial \psi_j}{\partial y_k}-\frac{\partial
\psi_k}{\partial y_j}\Big)(\theta,y,\alpha)\Big\v\leq C_N\, \v \Im
\psi(\theta,y,\alpha)\v^N
$$
for every $(\theta,y)$ in $\tilde \Omega_\delta$ and $\alpha$ in $\CS$.
\end{corollary}

{\bf Proof } According to (\ref{eqIV.3.38}) and Example \ref{eIV.3.9} we have
$\eta_k-\psi_k(\theta,y,\alpha)\in\CJ$. It follows from Lemma \ref{lIV.3.10}
that
$$
R_{jk}(\theta,y,\alpha):=\Big(\frac{\partial \psi_j}{\partial
y_k}-\frac{\partial \psi_k}{\partial y_j}\Big
)(\theta,y,\alpha)=\{\eta_j-\psi_j,\eta_k-\psi_k\}(\theta,y,\alpha)
$$
defines an element of $\CJ$. Since $R_{jk}$ does not depend on $\eta$ we can
apply Lemma \ref{lIV.3.11} and the conclusion follows. \cqfd

So far we have worked in the coordinates $(y,\eta)$. Let us go back to the
original coordinates $(x,\xi)$ and let us summarize the results already
obtained.

We set $x=y+x(\theta,\alpha)$,
$\xi=\frac\eta{\bra\theta\ket}+\xi(\theta,\alpha)$. Then
$(\theta,x)\in\Omega_\delta$ (see (\ref{eqIV.3.1})). Let us recall that,
$$
v_j(\theta,x,\xi,\alpha)=\xi_j(-\theta,x,\xi)-\alpha^j_\xi-i(x_j(-\theta,x,\xi)-\alpha^j_x)\,,
$$
(see (\ref{eqIV.2.10})).

Now let us introduce
\begin{equation}\label{eqIV.3.43}
\Phi_k(\theta,x,\alpha)=\xi_k(\theta,\alpha)+\frac{\frac 12\,(\sgn \theta)
(1-\chi_1(\theta))(x-x(\theta,\alpha))-(\tilde a_k+i\,\tilde
b_k)(\theta,x-x(\theta,\alpha),\alpha)}{\bra\theta\ket}
\end{equation}
where $\tilde a_k$, $\tilde b_k$ have been introduced in Corollary
\ref{eqIV.3.5}. Then we have,
\begin{theorem}\sl\label{tIV.3.13}
We can write
\begin{itemize}
\item[(i) ] $\displaystyle{\xi_k-\Phi_k(\theta,x,\alpha)=\som^n_{j=1}
d^j_k(\theta,x,\xi,\alpha)\,v_j(\theta,x,\xi,\alpha)}$

where $d^j_k$ are smooth functions defined for $(\theta,x)\in\Omega_\delta$ and
$$\Big\v\xi-\xi(\theta,\alpha)-\frac 12\,(1-\chi_1(\theta))\,(\sgn
\theta)\,\frac{x-x(\theta,\alpha)}{\bra\theta\ket}\Big\v\leq
\frac\delta{\bra\theta\ket}
\,.$$

 Moreover we have in this set,
\item[(ii) ] $\displaystyle{\v\partial ^A_x\,d^j_k(\theta,x,\xi,\alpha)\v\leq
\frac{C_A}{\bra\theta\ket}}$, \enskip $A\in\N^n$,
\item[(iii) ] $\displaystyle{\v\Phi_k(\theta,x,\alpha)-\alpha_\xi\v\leq
C(\varepsilon+\delta)}$,
\item[(iv) ] $\displaystyle{\Big\v\Im
\Phi_k(\theta,x,\alpha)-\frac{x_k-x_k(\theta,\alpha)}{1+4\theta^2}\Big\v\leq
\sqrt{5\delta}\,\enskip \frac{\v x-x(\theta,\alpha)\v}{\bra\theta\ket^2}}$.
\item[(v) ] $\v\partial ^A_x\,\Phi_k(\theta,x,\alpha)\v\leq \begin{cases}
C_A\,\bra\theta\ket^{-\v A\v}&\text{if } \enskip \v A\v\leq 1\\
C_A\,\bra\theta\ket^{-\v A\v-1}&\text{if }  \enskip \v A\v\geq 2\end{cases}$.
\item[(vi) ] $\Phi_k(\theta,x(\theta,\alpha),\alpha)=\xi(\theta,\alpha)$.
\item[(vii) ] $\displaystyle{\Big\v\Big(\frac{\partial \Phi_j}{\partial
x_k}-\frac{\partial \Phi_k}{\partial x_j}\Big)(\theta,x,\alpha)\Big\v\leq
\frac{C_N}{\bra\theta\ket}\,\Big(\frac{\v
x-x(\theta,\alpha)\v}{\bra\theta\ket}\Big)^N }$, \enskip $\forall N\in\N$

where the constants $C_A$, $C_0$, $C_N$ are independent of
$(\theta,x,\xi,\alpha)$.
\end{itemize}
\end{theorem}

{\bf Proof } It follows from Corollary \ref{cIV.3.5} that when
$\big\v\xi-\xi(\theta,\alpha)-\frac 12\,(\sgn
\theta)\,(1-\chi_1(\theta))\,\frac{x-x(\theta,\alpha)}{\bra\theta\ket}\big\v\leq
\frac\delta{\bra\theta\ket}$ then
$$
v_j(\theta,x,\xi,\alpha)=\sum^n_{k=1}
c_{jk}(\theta,x,\xi,\alpha)(\xi_k-\Phi_k(\theta,x,\alpha))
$$
with
$$
c_{jk}(\theta,x,\xi,\alpha)=\bra\theta\ket\,q_k\Big
(\bra\theta\ket(\xi-\xi(\theta,\alpha))-\frac 12\,(\sgn
\theta)\,(1-\chi_1(\theta))(x-x(\theta,\alpha)),\tilde a,\tilde b,\tilde
g_j\Big)
$$
where $\tilde q_k$ is defined in Corollary \ref{cIV.3.5} (i). By (v) of the
same result we have,
\begin{equation}\label{eqIV.3.44}
\big\v c_{jk}(\theta,x,\xi)-(1+2i\,\theta)\,\delta_{jk}\big\v\leq
C\,(\varepsilon+\delta)\,\bra\theta\ket\,.
\end{equation}
It follows then that the matrix $(c_{jk})$ is uniformly invertible.
Let us set $(d^j_k(\theta,x,\xi))=(c_{jk}(\theta,x,\xi))^{-1}$. Then we obtain
(i) in Theorem \ref{tIV.3.13}. The estimate (ii) follows then from
Corollary \ref{cIV.3.5} (v) and (\ref{eqIV.3.44}).

Let us now prove the claimed properties of $\Phi_k$. First of all since
$\xi(\theta,\alpha)=\alpha_\xi+\CO(\varepsilon)$, $\frac{\v
x-x(\theta,\alpha)\v}{\bra\theta\ket}\leq \delta$, $\v\tilde a_k\v+\v\tilde
b_k\v\leq C\,\delta$. We deduce easily from (\ref{eqIV.3.43}) that
$\v\Phi_k(\theta,x,\xi)-\alpha_\xi\v\leq C\,(\varepsilon+\delta)$. On the other
hand it follows from (\ref{eqIV.3.43}) and Corollary \ref{cIV.3.5} (i)
$$
\Im \Phi_k(\theta,x,\xi)=-\frac 1{\bra\theta\ket}\,\tilde
b_k(\theta,x-x(\theta,\alpha),\alpha)=\frac{x_k-x_k(\theta,\alpha)}{1+4\,\theta^2}+R
$$
where $\v R\v\leq \frac{\v x-x(\theta,\alpha)\v}{\bra\theta\ket^2}$.

The point (v) in Theorem \ref{tIV.3.13} follows easily from (\ref{eqIV.3.43})
and Corollary \ref{cIV.3.5} (ii)~; the point (vi) is obvious since $\tilde
a_k(\theta,0,\alpha)=\tilde b_k(\theta,0,\alpha)=0$. Finally for the point
(vii) we remark that according to (\ref{eqIV.3.38}) and (\ref{eqIV.3.43}) we
have
$$
(1)=\Big (\frac{\partial \Phi_j}{\partial x_k}-\frac{\partial \Phi_k}{\partial
x_j}\Big )(\theta,x,\alpha)=\frac 1{\bra\theta\ket}\,\Big (\frac{\partial
\psi_j}{\partial y_k}-\frac{\partial \psi_k}{\partial y_j}\Big
)(\theta,x-x(\theta,\alpha),\alpha)\,.
$$ 
Using Corollary \ref{eqIV.3.12} and the point (iv) we obtain
\begin{equation*}
\begin{aligned}
\v(1)\v&\leq \frac{C_N}{\bra\theta\ket}\,\v\Im
\psi(\theta,x-x(\theta,\alpha),\alpha)\v^N=\frac{C_N}{\bra\theta\ket}\,\big
(\bra\theta\ket\,\v\Im \Phi(\theta,x,\alpha)\v\big)^N\\
\v(1)\v&\leq \frac{C_N}{\bra\theta\ket}\,\Big (\frac{\v
x-x(\theta,\alpha)\v}{\bra\theta\ket}\Big)^N\,,\enskip \textrm{for all}\enskip
N\in\N\,.
\end{aligned}
\end{equation*}\cqfd

We need now to introduce the definition of Lagrangian ideals in the coordinates
$(x,\xi)$.

Let $(F(\cdot ,\alpha))_{\alpha\in\CS}$ a family of functions
$F(\theta,x,\xi,\alpha)$ defined for $(\theta,x)$ in $\Omega_\delta$ and
$\big \v\xi-\xi(\theta,\alpha)-\frac 12\,(\sgn
\theta)\,\frac{(x-x(\theta,\alpha))}{\bra\theta\ket}\big\v
<\frac\delta{\bra\theta\ket}$.
\begin{definition}\sl\label{dIV.3.14}
We shall say that the family $(F(\cdot ,\alpha))_{\alpha\in\CS}$ belongs to
$\CJ_{(x,\xi)}$ if
\begin{equation}\label{eqIV.3.45}
\Big (\bra\theta\ket\, F(\theta,y+x(\theta,\alpha),\,\frac\eta{\bra\theta\ket}+
\xi(\theta,\alpha),\alpha)\Big)_{\alpha\in\CS}\in\CJ
\end{equation}
where $\CJ$ has been introduced in Definition \ref{dIV.3.8}.

For example, $(\xi_j-\Phi_j(\theta,x,\alpha))_{\alpha\in\CS}\in\CJ_{(x,\xi)}$.
\end{definition}

As a consequence of Lemma \ref{lIV.3.11} we have the following result.
\begin{lemma}\sl\label{lIV.3.15}
Let $(R(\cdot ,\alpha)_{\alpha\in\CS}$ be in $\CJ_{(x,\xi)}$ and assume that
$R$ is independent of $\xi$ then for every $N\in\N$ there exists $C_N>0$ such
that,
\begin{equation}\label{eqIV.3.46}
\v R(\theta,x,\alpha)\v\leq \frac{C_N}{\bra\theta\ket}\,\Big (\frac{\v
x-x(\theta,\alpha)\v}{\bra\theta\ket}\Big )^N
\end{equation}
\end{lemma}

We can now pursue the proof of the existence of a phase as described in Theorem
\ref{tIV.1.2}.
\begin{lemma}\sl\label{lIV.3.16}
With $\Phi$ defined in (\ref{eqIV.3.43}) we have
$$
\Big (-\frac{\partial p}{\partial x_k}\,(x,\xi)-\frac{\partial \Phi_k}{\partial
\theta}\,(\theta,x,\alpha)-\sum^n_{j=1}\, \frac{\partial \Phi_k}{\partial
x_j}\,(\theta,x,\alpha)\,\frac{\partial p}{\partial \xi_j}\,(x,\xi)\Big
)_{\alpha\in\CS}\in\CJ_{(x,\xi)}\,.
$$
\end{lemma}

{\bf Proof } We know from (\ref{eqIV.3.43}), (\ref{eqIV.3.38}) and Example
\ref{eIV.3.9}, with $\xi=\frac\eta{\bra\theta\ket}+\xi(\theta,\alpha)$ that,
\begin{equation*}
\begin{aligned}
\xi_k-\Phi_k\,(\theta,x,\alpha)&=\frac 1{\bra\theta\ket}
\big(\eta_k-\psi_k(\theta,x-x(\theta,\alpha),\alpha)\big)\\
&=\frac 1{\bra\theta\ket}\, \som^n_{j=1}
f_{jk}(\theta,x-x(\theta,\alpha),\bra\theta\ket
(\xi-\xi(\theta,\alpha)),\alpha)\,v_j(\theta,x,\xi,\alpha)
\end{aligned}
\end{equation*}
where $(f_{jk}(\cdot ,\alpha))_{\alpha\in\CS}=((q_k(\eta-\frac{1}{2} \,(\sgn
\theta)(1-\chi_1(\theta))\,y,\tilde a,\tilde b,g_j))^{-1})_{\alpha \in\CS}
 \in\CH$ (see Definition
\ref{dIV.3.6}).

Recall now that $v_j(\theta,x,\xi,\alpha)=u_j\circ\chi_{-\theta}(x,\xi)$ where
$u_j(X,\Xi)=\Xi_j-\alpha^j_\xi-i(X_j-\alpha^j_x)$ (see (IV.2.9)). Let
us set
$$
\chi_{-\theta}(x,\xi)=(X,\Xi)\Longleftrightarrow x(\theta,X,\Xi)=x\,,\quad
\xi(\theta,X,\Xi)=\xi\,.
$$
It follows that
\begin{equation}\label{eqIV.3.47}
\begin{split}
\xi_k(\theta,X,\Xi)-\Phi_k(\theta,x(\theta,X,\Xi),\alpha)=\frac
1{\bra\theta\ket}\,\som^n_{j=1}
f_{jk}\big(\theta,x(\theta,X,\Xi)-x(\theta,\alpha),\\
\bra\theta\ket\,(\xi(\theta,X,\Xi)
-\xi(\theta,\alpha)),\alpha\big)\,u_j(X,\Xi)\,.
\end{split}
\end{equation}
Let us set
\begin{equation*}\left\{
\begin{array}{l}
M(\theta,X,\Xi,\alpha)=(\theta,x(\theta,X,\Xi)-x(\theta,\alpha),\,
\bra\theta\ket\,(\xi(\theta,X,\Xi)-\xi(\theta,\alpha)))\,,\\
\rho(\theta,X,\Xi)=(x(\theta,X,\Xi),\,\xi(\theta,X,\Xi))\,,\\
\rho(\theta,\alpha)=(x(\theta,\alpha),\,\xi(\theta,\alpha))\,.
\end{array}\right.
\end{equation*}
Now we differentiate (\ref{eqIV.3.47}) with respect to $\theta$ using the
equation of the flow given in (\ref{eqIII.1.2}). We obtain
\begin{equation*}
\begin{split}
-\frac{\partial p}{\partial x_k}\,(\rho(\theta,X,\Xi))&-\frac{\partial
\Phi_k}{\partial \theta}\,(\theta,x(\theta,X,\Xi),\alpha)-\som^n_{j=1}
\frac{\partial \Phi_k}{\partial x_j}\,(\theta,x(\theta,X,\Xi),\alpha)\cdot 
 \frac{\partial p}{\partial
\xi_j}\,(\rho(\theta,X,\Xi))\\
&\enskip \enskip =\som^n_{j=1}
\bigg[-\frac{\theta}{\bra\theta\ket^3}\,f_{jk}(M(\theta,X,\Xi,\alpha))+\frac
1{\bra\theta\ket}\,\frac{\partial f_{jk}}{\partial
\theta}\,(M(\theta,X,\Xi,\alpha))\\
&+\frac 1{\bra\theta\ket}\som^n_{\ell =1} \Big[\Big (\frac{\partial p}{\partial
\xi_\ell }\,(\rho(\theta,X,\Xi))-\frac{\partial p}{\partial \xi_\ell
}\,(\rho(\theta,\alpha))\Big)\,\frac{\partial f_{jk}}{\partial
y_\ell }\,(M(\theta,X,\Xi,\alpha))\\
&+\Big\{\frac\theta{\bra\theta\ket}\,(\xi_\ell (\theta,X,\Xi)-\xi_\ell
(\theta,\alpha))-\bra\theta\ket\, \Big (\frac{\partial p}{\partial x_\ell }\,
(\rho(\theta,X,\Xi)-\frac{\partial p}{\partial x_\ell
}\,(\rho,\alpha)\Big)\Big\}\,\\
&\hbox to 5cm{}\frac{\partial f_{jk}}{\partial \eta_\ell }\,
(M(\theta,X,\Xi,\alpha))\Big ]\bigg]\,u_j(X,\Xi)
\end{split}
\end{equation*}
We can now go back to the coordinates $(x,\xi)=\rho(\theta,X,\Xi)$ and then to
$(y,\eta)$ where $y=x-x(\theta,\alpha)$,
$\xi-\xi(\theta,\alpha)=\frac\eta{\bra\theta\ket}$. We obtain
\begin{equation*}
\begin{split}
&-\Big[\frac{\partial p}{\partial x_k}-\frac{\partial \Phi_k}{\partial
\theta}-\som^n_{j=1}\, \frac{\partial \Phi_k}{\partial x_j}\,\frac{\partial
p}{\partial
\xi_j}\Big]\Big(\theta,y+x(\theta,\alpha),\,\frac\eta{\bra\theta\ket}
+\xi(\theta,\alpha)\Big)\\
&\qquad=\som^n_{j=1}\bigg[\frac{-\theta}{\bra\theta\ket^3}\,f_{jk}
(\theta,y,\eta,\alpha)
+\frac 1{\bra\theta\ket}\, \frac{\partial f_{jk}}{\partial
\theta}(\theta,y,\eta,\alpha)\\
&+\frac
1{\bra\theta\ket}\,\som^n_{\ell =1}\Big (\frac{\partial
p}{\partial
\xi_\ell}\,(y+x(\theta,\alpha)\,,\frac{\eta}{\bra\theta\ket}+\xi(\theta,\alpha)\Big)
-\frac{\partial p}{\partial \xi_\ell }\,(\rho(\theta,\alpha))\,
\frac{\partial f_{jk}}{\partial y_\ell }\,(\theta,y,\eta,\alpha)\\
& +\som^n_{\ell =1}
\Big(\frac\theta{\bra\theta\ket^3}\,\eta_\ell
-\bigg\{\frac{\partial
p}{\partial x_\ell
}\,(y+x(\theta,\alpha)\,\frac\eta{\bra\theta\ket}+\xi(\theta,\alpha))
-\frac{\partial p}{\partial x_\ell
}\,(\rho(\theta,\alpha))\bigg\}\Big)\,\frac{\partial f_{jk}}{\partial
\eta_\ell
}\,(\theta,y,\eta,\alpha)\bigg]\\
&\qquad v_j\Big(\theta,y+x(\theta,\alpha)\,,
 \frac\eta{\bra\theta\ket}+\xi(\theta,\alpha)\Big)
=:\som^n_{j=1}\,
G_j(\theta,y,\eta,\alpha)\,v_j\Big(\theta,y+x(\theta,\alpha)\,\frac\eta{\bra\theta\ket}
+\xi(\theta,\alpha)\Big)\,.
\end{split}
\end{equation*}
According to Definition \ref{dIV.3.14} the Lemma will be proved if we show that
$(\bra\theta\ket\,G_j(\theta,y,\eta,\alpha))_{\alpha\in \CS}$ belongs to $\CH$
that is all the derivatives with respect to $(y,\eta)$ are uniformly bounded
when $\v y\v\leq \delta\,\bra\theta\ket$ and $\big\v \eta-\frac 12\,\sgn
\theta(1-\chi_1(\theta))\, y\v\leq \delta$.
Recall that $(f_{jk})=((q_k(\eta-\frac 12\,\sgn
\theta(1-\chi_1(\theta))\,y,\tilde a,\tilde b,g_j))^{-1}).$
Using (\ref{eqIV.3.5}) we see that $\frac{\partial \beta}{\partial
\theta}\in\CH$. Then differentiating (\ref{eqIV.3.19}) with respect to
$\theta$ we see that $\frac{\partial \alpha}{\partial \theta}$,
$\frac{\partial b}{\partial \theta}\in\CH$. It follows from (\ref{eqIV.3.9})
that $\frac{\partial g_j}{\partial \theta}\in \CH$. Then we deduce from the
estimates (\ref{eqIV.2.2}) that $\frac\partial {\partial
\theta}\big[q_k(\eta-\frac 12\,\sgn \theta(1-\chi_1(\theta))\,y,\tilde
a,\tilde b,g_j)\big]$ belongs to $\CH$ and we deduce from Corollary
\ref{cIV.3.5} (iii) that $\frac\partial {\partial
\theta}\big[(q_k(\eta-\frac 12\,\sgn \theta(1-\chi_1(\theta))\,y,\tilde
a,\tilde b,g_j))^{-1}\big]\in\CH$. Thus $(\frac{\partial f_{jk}}{\partial
\theta})\in\CH$.
 Since $f_{jk}$, $\frac{\partial
f_{jk}}{\partial \theta}$, $\frac{\partial f_{jk}}{\partial y_\ell }$,
$\frac{\partial f_{jk}}{\partial \eta_\ell }$ belong to $\CH$ and since $\CH$
is closed under multiplication it remains to prove that the functions
$h(\theta,y,\eta,\alpha)=\frac{\partial p}{\partial \xi_\ell
}\,(y+x(\theta,\alpha),\,\frac\eta{\bra\theta\ket}+\xi(\theta,\alpha))$ or
$\bra\theta\ket\,\frac{\partial p}{\partial x_\ell
}\,(y+x(\theta,\alpha),\,\frac\eta{\bra\theta\ket}+\xi(\theta,\alpha))$  belong
to $\CH$. Since $\frac{\eta_\ell }{\bra\theta\ket}+\xi_\ell (\theta,\alpha)$
has all its derivatives in $\eta$ uniformly bounded, we are led to prove that
if $g^{jk}(x)$ are the coefficients of $p(x,\xi)$ then $\partial
^A_x\,g^{jk}(y+x(\theta,\alpha))$ for $A\in\N^n$ and
$\bra\theta\ket\,\partial ^B_x\,g^{jk}(y+x(\theta,\alpha))$ for $B\in\N^n$, $\v
B\v\geq 1$ are uniformly bounded when $\v y\v\leq \delta\,\bra\theta\ket$.
This is obvious if $A=0$  and if $\v B\v\geq 1$ condition (\ref{eqIII.1.1})
shows that
\begin{equation*}
\begin{aligned}
\v\partial ^B_x\,g^{jk}(y+x(\theta,\alpha))\v&\leq \frac{C_B}{\bra
y+x(\theta,\alpha)\ket^{\v B\v+1+\sigma_0}}=\frac{C_B}{\bra
x(\theta,\beta)\ket^{\v B\v+1+\sigma_0}}\\
&\leq \frac{C_B}{\bra\theta\ket^{\v B\v+1+\sigma_0}}
\end{aligned}
\end{equation*}
since $\beta\in\CS\subset\CS_+\cap\CS_-$ (see (\ref{eqIII.3.3})).\cqfd
\begin{corollary}\sl\label{cIV.3.17}
With $\Phi$ defined in (\ref{eqIV.3.43}) we have for $k=1,\ldots ,n$,
$$
\Big(-\frac{\partial p}{\partial x_k}\,(x,\Phi(\theta,x,\alpha))-\frac{\partial
\Phi_k}{\partial \theta}\,(\theta,x,\alpha)-\som^n_{j=1} \frac{\partial
\Phi_k}{\partial x_j}\,(\theta,x,\alpha)\,\frac{\partial p}{\partial
\xi_j}\,(x,\Phi(\theta,x,\alpha))\Big)_{\alpha\in\CS}\in\CJ_{(x,\xi)}\,.
$$
\end{corollary}

{\bf Proof } First of all we show that
$$
(1)=\Big(\frac{\partial p}{\partial
x_k}\,(x,\Phi(\theta,x,\alpha))-\frac{\partial p}{\partial
x_k}\,(x,\xi)\Big)_{\alpha\in \CS}\in \CJ_{(x,\xi)}\,.
$$
Denoting by $g^{ij}$ the coefficients of $p$ we can write
$$
\frac{\partial g^{ij}}{\partial x_k}\,(x)\,\Phi_i\,\Phi_j-\frac{\partial
g^{ij}}{\partial x_k}\,(x)\,\xi_i\,\xi_j=\underbrace{\frac{\partial
g^{ij}}{\partial
x_k}\,(x)\,(\Phi_i-\xi_i)\,\xi_j}_{(a)}-\underbrace{\frac{\partial
g^{ij}}{\partial x_k}\,(x)\,\Phi_i(\xi_j-\Phi_j) }_{(b)}\,.
$$
To see that this belongs to $\CJ_{(x,\xi)}$ we use (\ref{eqIV.3.45}). In the
coordinates $(y,\eta)$ we have
$$
\bra\theta\ket\, (a)=-\frac{\partial g^{ij}}{\partial
x_k}\,(y+x(\theta,\alpha))\,\Big
(\frac{\eta_j}{\bra\theta\ket}+\xi_j(\theta,\alpha)\Big
)(\eta_i-\psi_i(\theta,y,\alpha))\,.
$$
Since $\eta_i-\psi_i(\theta,y,\alpha)\in \CJ$ (see Example \ref{eIV.3.9} and
(\ref{eqIV.3.38})) and
$\big(\frac{\eta_j}{\bra\theta\ket}+\xi_j(\theta,\alpha)\big)\,\frac{\partial
g^{ij}}{\partial x_k}(y+x(\theta,\alpha))$ belongs to $\CH$ we have
$(a)\in\CJ_{(x,\xi)}$.

A similar argument shows that (b) belongs to $\CJ_{(x,\xi)}$.

We show now that
$$
(2)=\Big(\frac{\partial \Phi_k}{\partial
x_j}\,(\theta,x,\alpha)\Big[\frac{\partial p}{\partial
\xi_j}\,(x,\Phi(\theta,x,\alpha))-\frac{\partial p}{\partial
\xi_j}\,(x,\xi)\Big]\Big)_{\alpha\in\CS}\in\CJ_{(x,\xi)}\,.
$$
The coefficients of $p$ can be written $g^{ij}(x)=\delta_{ij}+c_{ij}(x)$. It
follows that
$$
(2)=\Big (2\,\frac{\partial \Phi_k}{\partial x_j}\,(\theta,x,\alpha)\Big
(\Phi_j(\theta,x,\alpha)-\xi_j+\som^n_{\ell =1} c_{j\ell }(x)(\Phi_\ell
(\theta,x,\alpha)-\xi_\ell )\Big)\Big)\,.
$$
Now $\frac{\partial \Phi_k}{\partial x_j}=\frac
1{\bra\theta\ket}\,\frac{\partial \psi_k}{\partial y_\ell }\in\CH$ and
$c_{j\ell }(y+x(\theta,\alpha))\in\CH$. It follows that $(2)\in\CJ_{(x,\xi)}$
and the Corollary follows from Lemma \ref{lIV.3.16}. \cqfd
\begin{corollary}\sl\label{cIV.3.18}
For every $N\in\N$ there exists $C_N>0$ such that
$$
\Big\v-\frac{\partial p}{\partial
x_k}\,(x,\Phi(\theta,x,\alpha))-\frac{\partial \Phi_k}{\partial
\theta}\,(\theta,x,\alpha)-\frac{\partial \Phi_k}{\partial
x}\,(\theta,x,\alpha)\cdot \frac{\partial p}{\partial
\xi}\,(x,\Phi(\theta,x,\alpha))\Big\v
\leq \frac{C_N}{\bra\theta\ket}\,\Big(\frac{\v
x-x(\theta,\alpha)\v}{\bra\theta\ket}\Big)^N\,.
$$
\end{corollary}

{\bf Proof }This follows from Lemma \ref{lIV.3.15} and Corollary \ref{cIV.3.17}
since the left hand side does not depend of $\xi$. \cqfd
\begin{proposition}\sl\label{pIV.3.19}
Let $\alpha\in\CS$ (see (\ref{eqIV.3.1})). Let us set for
$(\theta,x)\in\Omega_\delta$,
\begin{equation}\label{eqIV.3.48}
\varphi (\theta,x,\alpha)=\int^1_0 (x-x(\theta,\alpha))\cdot
\Phi(\theta,s\,x+(1-s)\,x(\theta,\alpha),\alpha)\,ds+\theta\,p(\alpha)+\frac
1{2i}\,\v\alpha_\xi\v^2\,.
\end{equation}
Then we have
\begin{itemize}
\item[(i)] $\displaystyle{\varphi (0,x,\alpha)=(x-\alpha_x)\cdot
\alpha_\xi+\frac i2\,(x-\alpha_x)^2+\frac 1{2i}\,\v\alpha_\xi\v^2+\CO (\v
x-\alpha_x\v^N)}$.

For every $N\in\N$ there exists $C_N>0$ such that
\item[(ii)] $\displaystyle{\Big\v\frac{\partial \varphi }{\partial
x}\,(\theta,x,\alpha)-\Phi(\theta,x,\alpha)\Big\v\leq C_N\Big (\frac{\v
x-x(\theta,\alpha)\v}{\bra\theta\ket}\Big)^N }$.
\item[(iii)] $\displaystyle{\Big\v\frac{\partial \varphi }{\partial
\theta}\,(\theta,x,\alpha)+p\Big (x,\,\frac{\partial \varphi }{\partial
x}\,(\theta,x,\alpha)\Big) \Big\v\leq  C_N\Big (\frac{\v
x-x(\theta,\alpha)\v}{\bra\theta\ket}\Big)^N }$

uniformly with respect to $(\theta,x)\in\Omega_\delta$ and $\alpha\in\CS$.

Moreover, uniformly with respect to
$(\theta,x,\alpha)\in\Omega_\delta\times\CS$, we have
\item[(iv)] $\displaystyle{\Big\v\frac{\partial \varphi }{\partial
x}\,(\theta,x,\alpha)-\alpha_\xi\Big\v\leq C\,(\varepsilon+\sqrt\delta)}$.
\item[(v)] $\v\partial ^A_x\,\varphi (\theta,x,\alpha)\v\leq C_A$,\enskip
$\forall\,A\in\N^n$.
\item[(vi)] $\displaystyle{\Big\v\Im \varphi (\theta,x,\alpha)-\frac
12\,\frac{\v x-x(\theta,\alpha)\v^2}{1+4\,\theta^2}+\frac 12\,\v\alpha_\xi\v^2
\Big\v\leq C\,(\varepsilon+\sqrt\delta)\,\frac{\v
x-x(\theta,\alpha)\v^2}{\bra\theta\ket^2}}$.
\end{itemize}
\end{proposition}

{\bf Proof } If we can prove that for $j=1,\ldots ,n$,
$$
\Phi_j(0,s\,x+(1-s)\,\alpha_x,\alpha)=\alpha^j_\xi+i\,s(x_j-\alpha^j_x)+
\CO(s^N\v x-\alpha_x\v^N)
$$
then (i) will follow according to (\ref{eqIV.3.48}). By (\ref{eqIV.3.43}),
Corollary \ref{cIV.3.5} and Theorem \ref{tIV.3.1} we have
$$
\Phi_j(0,x,\alpha)=\alpha^j_\xi-(\tilde a_j+i\,\,\tilde
b_j)(0,x-\alpha_x,\alpha)=\alpha^j_\xi-(a_k+i\,b_k)(0,x-\alpha_x,\alpha)\,.
$$
Now Example \ref{eIV.3.9} (for $\theta=0$) shows that
$\eta_j+a_j(0,y,\alpha)+i\,b_j(0,y,\alpha)$ belongs to the ideal $\CJ$
introduced in Definition \ref{dIV.3.8}. On the other hand, since by
(\ref{eqIV.2.10}) for $\theta=0$ $g_j(\eta)=\chi_0(\frac
1{\mu_0}\,\eta)(\eta_j-i\,y_j)$ it follows that $\eta_j-i\,y_j$ belongs also
to $\CJ$. Thus the difference $a_j(0,y,\alpha)+i\,b_j(0,y,\alpha)+i\,y_j$
belongs also to $\CJ$ and does not depend on $\eta$. It follows from Lemma
\ref{lIV.3.11} that for all $N\in\N$,
\begin{equation*}
\begin{aligned}
a_j(0,y,\alpha)=\CO(\v \Im \psi(0,y,\alpha)\v^N)\\
b_j(0,y,\alpha)+y_j=\CO(\v\Im \psi(0,y,\alpha)\v^N)\,.
\end{aligned}
\end{equation*}
Now (\ref{eqIV.3.38}) and Theorem \ref{tIV.3.1} (ii) show that $\v \Im
\psi_j(0,y,\alpha)\v=\v b_j(0,y,\alpha)\v\leq C\,\v y\v$. Thus for all
$N\in\N$,
$$
a_j(0,y,\alpha)=\CO(\v y\v^N)\,,\quad b_j(0,y,\alpha)=-y_j+\CO(\v y\v^N)\,.
$$
It follows that for all $N\in\N$,
$$
\Phi_j(0,y,\alpha)=\alpha^j_\xi+i(x-\alpha^j_x)+\CO(\v x-\alpha_x\v^N)
$$
which proves our claim.

 Let us prove (ii). We have, by
(\ref{eqIV.3.48}),
\begin{equation*}
\begin{split}
\frac{\partial \varphi }{\partial x_j}\,(\theta,x,\alpha)=\int^1_0
\Phi_j(\theta,s\,x+(1-s)\,x(\theta,\alpha),\alpha)\,ds+\som^n_{k=1} \int^1_0
s(x_k-x_k(\theta,\alpha))\cdot \\
\cdot \frac{\partial \Phi_k}{\partial
x_j}\,(s\,x+(1-s)\,x(\theta ,\alpha),\alpha)\,ds
\end{split}
\end{equation*}
Now we use Theorem \ref{tIV.3.13} (vii) and the fact that $\v
x-x(\theta,\alpha)\v\leq \delta\,\bra\theta\ket$. We deduce that
\begin{equation*}
\begin{split}
\Big\v\frac{\partial \varphi }{\partial x_j}\,(\theta,x,\alpha)&-\int^1_0
\Big[\Phi_j+s\,\som^n_{k=1} (x_k-x_k(\theta,\alpha))\,\frac{\partial
\Phi_j}{\partial
x_k}\Big](\theta,s\,x+(1-s)\,x(\theta,\alpha),\alpha)\,ds\Big\v\\
&\hbox to 5cm{}\leq C_N\Big (\frac{\v
x-x(\theta,\alpha)\v}{\bra\theta\ket}\Big)^N
\end{split}
\end{equation*}
where $C_N$ is independent of $(\theta,x,\alpha)$.

It follows that
\begin{equation*}
\begin{split}
\Big\v\frac{\partial \varphi }{\partial x_j}\,(\theta,x,\alpha)&-\int^1_0
\Phi_j(\theta,s\,x+(1-s)\,x(\theta,\alpha),\alpha)\,ds\\
&-\int^1_0 s\,\frac
d{ds}\big[\Phi_j(\theta,s\,x +(1-s)\,x(\theta,\alpha),\alpha)\big]\,ds\Big\v
\leq C_N\Big (\frac{\v x-x(\theta,\alpha)\v}{\bra\theta\ket}\Big)^N\,.
\end{split}
\end{equation*}
Integrating by parts in the second integral above, we obtain (ii). As a
consequence of (ii) we have the estimate
\begin{equation}\label{eqIV.3.49}
\Big\v p(x,\Phi(\theta,x,\alpha))-p\Big (x,\,\frac{\partial \varphi }{\partial
x}\,(\theta,x,\alpha)\Big)\Big\v\leq C_N \Big (\frac{\v
x-x(\theta,\alpha)\v}{\bra\theta\ket}\Big)^N\,.
\end{equation}
Let us prove (iii). We deduce from (\ref{eqIV.3.48}) that
\begin{equation}\label{eqIV.3.50}
\frac{\partial \varphi }{\partial \theta}\,(\theta,x,\alpha)=(1)+(2)+(3)+(4)
\end{equation}
where
\begin{equation*}
\begin{aligned}
(1)&=-\som^n_{k=1} \int^1_0 \dot
x_k(\theta,\alpha)\,\Phi_k(\theta,s\,x+(1-s)\,x(\theta,\alpha),\alpha)\,ds\\
(2)&=\som^n_{k=1} \int^1_0 (x_k-x_k(\theta,\alpha))\,\frac{\partial
\Phi_k}{\partial \theta}\,(\theta,s\,x+(1-s)\,x(\theta,\alpha),\alpha)\,ds\\
(3)&=\som^n_{k=1} \int^1_0 (x_k-x_k(\theta,\alpha))\,\som^n_{j=1}
\frac{\partial \Phi_k}{\partial
x_j}\,(\theta,s\,x+(1-s)\,x(\theta,\alpha),\alpha)(1-s)\,\dot
x_j(\theta,\alpha)\,ds\\
(4)&=p(\alpha)\,.
\end{aligned}
\end{equation*}
Let us consider the term (2). We use Corollary \ref{cIV.3.18} to get
\begin{equation*}
\begin{split}
(2)&=\som^n_{k=1} \int^1_0 (x_k-x_k(\theta,\alpha))\Big[-\frac{\partial
p}{\partial
x_k}\,(s\,x+(1-s)\,x(\theta,\alpha),\Phi(\theta,s\,x+(1-s)\,
x(\theta,\alpha),\alpha))\\
&\enskip -\som^n_{j=1} \frac{\partial \Phi_k}{\partial
x_j}\,(\theta,s\,x+(1-s)\,x(\theta,\alpha),\alpha)\,\frac{\partial p}{\partial
\xi_j}\,(s\,x+(1-s)\,x(\theta,\alpha),\Phi(\theta,s\,x+(1-s)\,
x(\theta,\alpha),\alpha))\Big]\\
&\hbox to 6cm{} +\CO \Big (\Big(\frac{\v
x-x(\theta,\alpha)\v}{\bra\theta\ket}\Big)^N\,\Big)\,.
\end{split}
\end{equation*}
Now, by Theorem \ref{tIV.3.13}, (vii),
\begin{equation}\label{eqIV.3.51}
\Big\v\Big(\frac{\partial \Phi_k}{\partial x_j}-\frac{\partial
\Phi_j}{\partial x_k}\Big )(\theta,s\,x+(1-s)\,x(\theta,\alpha),\alpha)\Big\v
\leq C_N\, \frac{s^N}{\bra\theta\ket}\,\frac{\v
x-x(\theta,\alpha)\v^N}{\bra\theta\ket^N}
\end{equation}
and $s^N\leq 1$. Therefore,
$$
(2)=-\int^1_0 \frac
d{ds}\big[p(s\,x+(1-s)\,x(\theta,\alpha),\Phi(\theta,s\,x+(1-s)\,
x(\theta,\alpha),\alpha))\big]\,ds
+\CO\Big (\frac{\v x-x(\theta,\alpha)\v^N}{\bra\theta\ket^N}\Big)\,.
$$
Therefore we obtain
\begin{equation}\label{eqIV.3.52}
\Big\v(2)+p(x,\Phi(\theta,x,\alpha))-p(x(\theta,\alpha),\Phi(\theta,x(\theta,\alpha),\alpha))\Big\v
\leq C_N \Big (\frac{\v x-x(\theta,\alpha)\v}{\bra\theta\ket}\Big)^N\,.
\end{equation}
Let us consider the term (3). Using again (\ref{eqIV.3.51}) we get
$$
(3)=\som^n_{j=1} \int^1_0 (1-s)\,\dot x_j(\theta,\alpha)\,\frac d{ds}
(\Phi_j(\theta,s\,x+(1-s)\,x(\theta,\alpha),\alpha))+\CO\Big (\frac{\v
x-x(\theta,\alpha)\v^N}{\bra\theta\ket^N}\Big )\,.
$$
Integrating by part we obtain
$$
(3)=-\dot x(\theta,\alpha)\cdot
\Phi(\theta,x(\theta,\alpha),\alpha)+\som^n_{j=1} \int^1_0 \dot
x_j(\theta,\alpha)\,\Phi_j(\theta,s\,x+(1-s)\,x(\theta,\alpha),\alpha)\,ds\,.
$$
Comparing with the term (1) we obtain,
$$
\big\v(1)+(3)+\dot
x(\theta,\alpha)\,\Phi(\theta,x(\theta,\alpha),\alpha)\big\v\leq C_N \Big
(\frac{\v x-x(\theta,\alpha)\v}{\bra\theta\ket}\Big)^N\,.
$$
Now by Theorem \ref{tIV.3.13} (vi) and the Euler relation we have,
\begin{equation*}
\begin{aligned}
\dot x(\theta,\alpha)\cdot
\Phi(\theta,x(\theta,\alpha),\alpha)&=\xi(\theta,\alpha)\,\frac{\partial
p}{\partial \xi}\,(x(\theta,\alpha),\xi(\theta,\alpha)\\
&=2p(x(\theta,\alpha),\xi(\theta,\alpha))=2p(\alpha)\,.
\end{aligned}
\end{equation*}
It follows that
\begin{equation}\label{eqIV.3.53}
\big\v(1)+(3)+2p(\alpha)\big\v\leq C_N\Big (\frac{\v
x-x(\theta,\alpha)\v}{\bra\theta\ket}\Big)^N\,.
\end{equation}
Since in (\ref{eqIV.3.52}) we have
$p(x(\theta,\alpha),\Phi(\theta,x(\theta,\alpha),\alpha))
=p(x(\theta,\alpha),\xi(\theta,\alpha))=p(\alpha)$, we deduce from
(\ref{eqIV.3.50}), (\ref{eqIV.3.52}) and (\ref{eqIV.3.53}) that
\begin{equation}\label{eqIV.3.54}
\Big\v\frac{\partial \varphi }{\partial
\theta}\,(\theta,x,\alpha)+p(x,\Phi(\theta,x,\alpha))\Big\v\leq C_N \Big (
\frac{\v x-x(\theta,\alpha)\v}{\bra\theta\ket}\Big)^N\,.
\end{equation}
Therefore (iii) follows from (\ref{eqIV.3.54}) and (\ref{eqIV.3.49}).

Finally (iv), (v), (vi) follow easily from Theorem \ref{tIV.3.13}. \cqfd
\begin{remark}\sl\label{rIV.3.20}
Assume that $\alpha$ is such that
\begin{equation}\label{eqIV.3.55}
\frac 12\leq \v\alpha_\xi\v\leq 2\quad \textrm{and}\quad \alpha_x\cdot
\alpha_\xi\leq c_0\,\bra \alpha_x\ket\,\v\alpha_\xi\v 
\end{equation}
(so $\alpha\in\CS_-$ instead of $\alpha\in \CS$). Then Theorem \ref{tIV.3.1},
Corollary \ref{cIV.3.5} and Proposition \ref{pIV.3.19} are true for
$\theta\leq 0$.

By the same way if $\frac 12\leq \v\alpha_\xi\v\leq 2$ and $\alpha_x\cdot
\alpha_\xi\geq -c_0\,\bra\alpha_x\ket\,\v\alpha_\xi\v$ (which imply that
$\alpha\in\CS_+$) the above results hold for $\theta\geq 0$.
\end{remark}

\subsection{The case of incoming points}\label{ssIV.4}

We are going to prove Theorem \ref{tIV.1.2} when
$$
\v\alpha_x\cdot \alpha_\xi\v>c_0\,\bra\alpha_x\ket\,\v\alpha_\xi\v\quad
\textrm{and}\quad \frac 12\, \leq \v\alpha_\xi\v\leq 2\,.
$$
Since the problem is entirely symetric we can without loss of generality
assume that
\begin{equation}\label{eqIV.4.1}
\frac 12\leq \v\alpha_\xi\v\leq 2\,,\quad \alpha_x\cdot
\alpha_\xi<-c_0\,\bra\alpha_x\ket \v\alpha_\xi\v\,.
\end{equation}
It follows from Definition \ref{dIV.1.1} and  the discussion after, that
\begin{equation*}
\begin{split}
\tilde \Omega_\delta=\big\{(\theta,y)\in\R\times\R^n : \theta\leq 0\,,\enskip
\v y\v\leq \delta\,\bra\theta\ket\big\}\cup\big\{(\theta,y)\in\R\times\R^n :
\theta\geq 0\,,\enskip \v y\v\leq \delta\,\bra\theta\ket\\
\textrm{and}\quad (y+x(\theta,\alpha))\cdot \alpha_\xi\leq c_1\,\bra
y+x(\theta,\alpha)\ket\,\v\alpha_\xi\v\big\}\,.
\end{split}
\end{equation*}

Let now $\chi_0\in C^\infty _0(\R^n)$, $\chi_1\in C^\infty (\R^n)$ be such
that,
\begin{equation*}
\begin{aligned}
\chi_0(t)=1&\quad \textrm{if}\quad \v t\v\leq 1\,,\quad \chi_0(t)=0\quad
\textrm{if}\quad \v t\v\geq  2\quad \textrm{and}\quad 0\leq \chi_0\leq 1\,,\\
\chi_1(\theta)=1&\quad \textrm{if}\quad \theta\geq -1\,,\quad
\chi_1(\theta)=0\quad
\textrm{if}\quad \theta\leq -2\quad \textrm{and}\quad 0\leq \chi_1\leq
1\,.
\end{aligned}
\end{equation*}
For $j=1,\ldots ,n$ we introduce
\begin{equation}\label{eqIV.4.2}
g_j(\eta)=\chi_0\Big (\frac 1{\mu_0}\,\eta\Big
)\,v_j(\theta,y+x(\theta,\alpha),\eta\,\chi_1(\theta)+(1-\chi_1(\theta))\Big[\frac
\eta{\bra\theta\ket}+\frac 12\,\frac{\sgn
\theta}{\bra\theta\ket}\,y\Big]+\xi(\theta,\alpha),\alpha))\,,
\end{equation}
where $\mu_0$ is a small constant to be chosen, $(\theta,y)\in\tilde
\Omega_\delta$,
$\alpha$ satisfies (\ref{eqIV.4.1}) and $v_j$ has been introduced in
(\ref{eqIV.2.9}).
\begin{remark}\sl\label{rIV.4.1}
 (i)
According to Remark \ref{rIV.3.20} since (\ref{eqIV.4.1}) implies
(\ref{eqIV.3.55}) the phase has been already constructed when $(\theta,x) $
belongs to the first part of $\tilde \Omega_\delta$ where $\theta\leq 0$.
Therefore we are left with the case $\theta\geq 0$.

(ii)  If $\alpha$ satisfies (\ref{eqIV.4.1}) and $(\theta,y)\in\tilde
\Omega_\delta$, $\theta\geq 0$, then the point
$(y+x(\theta,\alpha),\eta+\xi(\theta,\alpha))$ belongs to $\CS_-$. Indeed
recall that $\frac 12\leq \v\alpha_\xi\v\leq 2$ implies $\frac 14\leq
\frac12\,\v\alpha_\xi\v\leq \v\eta+\xi(\theta,\alpha)\v\leq 2\,\v\alpha_\xi\v$
if $\v\eta\v\leq 2\,\mu_0$ and $\mu_0,\varepsilon$ are small enough. Therefore
setting $x=x(\theta,\alpha)$ we can write
\begin{equation*}
\begin{aligned}
x\cdot (\eta+\xi(\theta,\alpha))&=x\cdot
(\alpha_\xi+\eta+\zeta(\theta,\alpha))\leq c_1\,\bra x\ket\,\v\alpha_\xi\v+\v
x\v\,(\mu_0+\varepsilon)\\
&\leq 2\,c_1\,\langle x\rangle\,\v
\eta+\xi(\theta,\alpha)\v+4(\mu_0+\varepsilon)\,\langle
x\rangle\,\v\eta+\xi(\theta,\alpha)\v\\
&\leq \frac 14 \,\bra x\ket\,\v\eta+\xi(\theta,\alpha)\v
\end{aligned}
\end{equation*}
if $c_1,\mu_0,\varepsilon$ are small enough.
\end{remark}

Our first step will be the proof of the following result.
\begin{theorem}\sl\label{tIV.4.2}
There exist small positive constants $\mu_0,\delta$ and $C^\infty $ functions
$a=a(\theta,y,\alpha)$, $b_k=b_k(\theta,y,\alpha)$, $k=1,\ldots ,n$, defined
on $\tilde \Omega_\delta$ with $\theta\geq 0$ such that, with $a=(a_k)$,
$b=(b_k)$ we have for $\eta\in\R^n$,
\begin{itemize}
\item[(i)] $\displaystyle{g_j(\eta)=\som^n_{k=1}
q_k(\eta,a,b,g_j)(\eta_k+a_k(\theta,y,\alpha)+i\,b_k(\theta,y,\alpha))}$

where the $q'_k\,s$ have been introduced in Lemma \ref{lIV.2.1}.

Moreover we have for $(\theta,y)\in\tilde \Omega_\delta$, $\theta\geq 0$ and
$k=1,2,\ldots ,n$,
\item[(ii)] $\displaystyle{\Big\v
a_k(\theta,y,\alpha)+\frac{2\theta\,y_k}{1+4\,\theta^2}\Big\v\leq
\sqrt\delta\, \inf (1,\v y\v)}$,

$\displaystyle{\Big\v b(\theta,y,\alpha)+\frac {y_k}{1+4\,\theta^2}\Big\v\leq
\sqrt\delta\, \frac{\v y\v}{\bra\theta\ket^2} }$.
\item[(iii)] If we set $\tilde
a(\theta,y,\alpha)=a(\theta,y,\alpha)+\frac{2\theta\,y}{1+4\,\theta^2}$ and
$\tilde b(\theta,y,\alpha)=\bra\theta\ket \big (b(\theta,y,\alpha)+\frac
y{1+4\,\theta^2}\big )$ then for every $A\in\N^n$, $\v A\v\geq 1$ one can find
$C_A\geq 0$ such that with $x=y+x(\theta,\alpha)$, $\theta\geq 0$,
$\displaystyle{
\v\partial ^A_y\,\tilde a(\theta,y,\alpha)\v+\v\partial ^A_y\,\tilde
b(\theta,y,\alpha)\v\leq C_A \Big[\frac \varepsilon{\bra x\ket^{\sigma_0}}
\Big (\frac 1{\bra x\ket}+\frac 1{\bra\theta\ket}\Big )^{\v A\v+1}+\frac
\delta{\bra\theta\ket^{\v A\v+1}}\Big ]}$
\item[(iv)] $\v q_k(\eta,a,b,g_j)-(1+2i\theta)\,\delta_{jk}\v\leq C
(\varepsilon+\sqrt \delta)\,\bra\theta\ket$ \quad if\quad  $\v\eta\v\leq
\sqrt\delta$.
\item[(v)] $\v\partial ^B_{(a,b)}\,\partial
^\gamma_\eta\,q_k(\eta,a,b,g_j)\v\leq C_{B,\gamma}\,\bra\theta\ket$,\quad
if\quad $B\in\N^n$, \quad $\gamma\in\N^n$, \quad$ 1\leq k\leq n$.
\end{itemize}
\end{theorem}

{\bf Proof } We use the same method as in Theorem \ref{tIV.3.1}. According to
Lemma \ref{lIV.2.1}, the claim (i) is equivalent to solve the system of
equations
\begin{equation}\label{eqIV.4.3}
r(a,,g_j(\theta,y,\alpha;\cdot )=0\,,\quad j=1,\ldots ,n\,.
\end{equation}
We shall solve this system in the set
\begin{equation}\label{eqIV.4.4}
E=\Big\{(a,b)\in\R^n\times\R^n : \Big\v
a+\frac{2\theta y}{1+4\theta^2}\Big\v\leq \sqrt\delta \, \inf (1,\v y\v),\,
\Big \v b+\frac y{1+4\theta^2}\Big\v\leq \sqrt\delta \,\frac {\v
y\v}{\bra\theta\ket^2}\Big\}
\end{equation}
where $0<\delta\ll 1$.

First of all we give equivalent equations to (\ref{eqIV.4.3}) in the set $E$.
We write as in (\ref{eqIV.3.13})
\begin{equation}\label{eqIV.4.5}
r(a,b,g_j)=g_j(-a)-i\som^n_{k=1}\,\frac{\partial g_j}{\partial
\xi_k}\,(-a)\,b_k+\som^n_{p,q=1} \,H^j_{pq}(\theta,y,\alpha,a,b)\,b_p\,b_q
\end{equation}
where
\begin{equation}\label{eqIV.4.6}
H^j_{pq}(\theta,y,\alpha,a,b)=\int^1_0 \frac{\partial ^2r}{\partial
b_p\,\partial b_q}\,(a,t\,b,g_j(\theta,y,\alpha;\cdot ))(1-t)\,dt\,.
\end{equation}
By (\ref{eqIV.2.2}) we have
\begin{equation}\label{eqIV.4.7}
\v\partial ^\alpha_{(a,b)}\,\partial ^B_y\,r(a,t\,b,g_j(\cdots ))\v\leq C_{AB}
\som_{\v\gamma\v\leq \v A\v+3n} \int \v\partial  ^\gamma_\xi\,\partial 
^B_y\,g_j(\theta,y,\alpha,\xi)\v\,d\xi\,.
\end{equation}
Since for $\theta\geq 0$ we have,
$$
g_j(\eta)=\chi\,\Big(\frac\eta{\mu_0}\Big)
\big[(\xi_j-i\,x_j)(-\theta;y+x(\theta;\alpha),\eta+\xi(\theta;\alpha))+
i(\alpha_x+i\,\alpha_\xi)\big]
$$
we deduce from Propositions \ref{pIII.3.1} and \ref{pIII.3.2} that  $\v
\partial ^\gamma_\eta\,\partial ^B_y\,g_j(\theta,y,\alpha,\eta)\v$ is bounded
on the support of
$\chi$, by
\begin{equation*}\left\{
\begin{array}{ll}
C\,\bra\theta\ket&\quad \textrm{if}\quad B=0\,,\\
C\Big (1+\frac{\varepsilon\bra\theta\ket}{\bra x\ket^{2+\sigma_0}}\Big)&\quad
\textrm{if}\quad \v B\v=1\,,\\
\frac{C_B\varepsilon}{\bra x\ket^{\v B\v+\sigma_0}}\Big
(1+\frac{\bra\theta\ket}{\bra x\ket}\big)&\quad \textrm{if}\quad \v B\v\geq
2\,.
\end{array}\right.
\end{equation*}
It follows that
\begin{equation}\label{eqIV.4.8}
\v\partial ^A_{(a,b)}\,\partial ^B_y\,H^j_{p,q}(\theta,y,\alpha,a,b)\v\leq
\begin{cases} C\,\bra\theta\ket &\textrm{if\enskip  $\v B\v=0$}\\
C\big (1+\frac{\varepsilon\,\bra\theta\ket}{\bra
x\ket^{2+\sigma_0}}\big)&\textrm{if\enskip  $\v B\v=1$}\\
\frac{C_B\varepsilon}{\bra x\ket^{\v
B\v+\sigma_0}}\,\big(1+\frac{\bra\theta\ket}{\bra
x\ket}\big)&\textrm{if\enskip  $\v B\v\geq 2$}\\
\end{cases}
\end{equation}
Since $\chi(-a)=1$ and $\chi'(-a)=0$ we see that (\ref{eqIV.4.3}) is
equivalent in $E$ to the vectorial equation,
\begin{equation}\label{eqIV.4.9}
\begin{split}
\Big[\xi-i\,x&-i\,\som^n_{k=1} \Big(\frac{\partial \xi}{\partial
\xi_k}-i\,\frac{\partial x}{\partial \xi_k}\Big
)\,b_k\Big](-\theta;y+x(\theta;\alpha),\xi(\theta;\alpha)-a)\\
&+i(\alpha_x+i\,\alpha_\xi)+\som^n_{p,q=1}
H_{pq}(\theta,y,\alpha,a,b)\,b_p\,b_q=0\,.
\end{split}
\end{equation}
To shorten the notations we shall set
\begin{equation}\label{eqIV.4.10}\left\{
\begin{array}{l}
\rho(\theta;\alpha)=(x(\theta;\alpha),\xi(\theta;\alpha))\\
\rho_y(\theta;\alpha)=(y+x(\theta;\alpha),\xi(\theta;\alpha))\\
\rho_{y,a}(\theta;\alpha)=(y+x(\theta;\alpha),\xi(\theta;\alpha)-a)\,.
\end{array}\right.
\end{equation}
Since, by assumption, the point $\rho_{y,a}(\theta;\alpha)$ belongs to
$\CS_-$ (the outgoing set for $\theta\leq 0$) we can use the Proposition
\ref{pIII.3.1} to write
\begin{equation*}
\begin{aligned}
x(-\theta;\rho_{y,a}(\theta;\alpha))&=y+x(\theta;\alpha)-2\theta\,
\xi(-\theta;\rho_{y,a}(\theta;\alpha))+z(-\theta;\rho_{y,a}(\theta;\alpha))\\
\xi(-\theta;\rho_{y,a}(\theta;\alpha))&=-a+\xi(\theta;\alpha)
+\zeta(-\theta;\rho_{y,a}(\theta;\alpha))\,.
\end{aligned}
\end{equation*}
It follows that (\ref{eqIV.4.9}) is equivalent to
\begin{equation*}
\begin{split}
(1+2i\,\theta)\,\xi(&-\theta;\rho_{y,a}(\theta;\alpha))-i\,y-i\,x(\theta;\alpha)-i\,z
(-\theta;\rho_{y,a}(\theta;\alpha))\\
&-i\Big[(1+2i\,\theta)\som^n_{k=1} \frac{\partial \xi}{\partial \xi_k}\,
(-\theta;\rho_{y,a}(\theta;\alpha))\,b_k-i\,\som^n_{k=1} \frac{\partial
z}{\partial \xi_k}\,(-\theta;\rho_{y,a}(\theta;\alpha))\,b_k\Big ]\\
&+i(\alpha_x+i\,\alpha_\xi)+\som^n_{p,q=1} H_{pq}(\cdots )\,b_p\,b_q=0\,.
\end{split}
\end{equation*}
Taking the real and the imaginary parts, we are led to the $2n$ real equations
\begin{equation*}
\begin{split}
\xi_j(&-\theta;\rho_{y,a}(\theta;\alpha))+2\theta\,b_j+2\theta\,\frac{\partial
\zeta_j}{\partial \xi}\,(-\theta;\rho_{y,a}(\theta;\alpha))\cdot
b-\frac{\partial z_j}{\partial \xi}\,(-\theta;\rho_{y,a}(\theta;\alpha))\cdot
b-\alpha^j_\xi+H^j_1\,b\cdot
b=0\\
&-2\theta\,\xi_j(-\theta;\rho_{y,a}(\theta;\alpha))+b_j+y_j+x_j(\theta;\alpha)
-\alpha^j_x+z_j(-\theta;\rho_{y,a}(\theta;\alpha))\\
&\quad \quad \quad \quad \quad +\frac{\partial \zeta_j}{\partial
\xi}\,(-\theta;\rho_{y,a}(\theta;\alpha))\cdot b
+H^j_2\,b\cdot b=0
\end{split}
\end{equation*}
where
$$
\frac{\partial \zeta_j}{\partial \xi}\cdot b=\som^n_{k=1} \frac{\partial
\zeta_j}{\partial \xi_k}\cdot b_k\,,\quad H^j\,b\cdot b=\som^n_{p,q=1}
H^j_{p,q}\,b_p\,b_q\,,\quad H^j=H^j_1+i\,H^j_2\,.
$$
Setting $X=\begin{pmatrix} \xi_j(-\theta;\rho_{y,a}(\theta;\alpha))\\b_j
\end{pmatrix}$, $A=\begin{pmatrix} 1 & 2\theta\\ -2\theta & 1\end{pmatrix}$,
our system can be written $A\,X=F$. Since $A^{-1}=\frac 1{1+4\theta^2}
\begin{pmatrix} 1 & -2\theta\\ 2\theta & 1\end{pmatrix}$, it is equivalent to
the following system.
\begin{equation*}
\begin{split}
\xi(-\theta;\rho_{y,a}(\theta;\alpha))&=\frac{2 \theta y}{1+4\theta^2}+\frac
1{1+4\theta^2} \,
\Big[2\theta(x(\theta;\alpha)-\alpha_x)+2\theta\,z(-\theta;\rho_{y,a}
(\theta;\alpha))\\
&\quad \quad \quad \quad \quad \quad +\alpha_\xi+\frac{\partial z}{\partial
\xi}\,(-\theta;\rho_{y,a}(\theta;\alpha))\cdot
b-(H_1-2\theta\,H_2)\,b\cdot b\Big]\,.
\end{split}
\end{equation*}
\begin{equation*}
\begin{split}
b=\frac{-y}{1+4\theta^2}-\frac 1{1+4\theta^2}
\big[x(\theta;\alpha)&-\alpha_x-2\theta\,\alpha_\xi+z(-\theta;\rho_{y,a}(\theta;\alpha))
+(2\theta\,H_1+H_2)\,b\cdot b\big]\\
&-\frac{\partial \zeta}{\partial \xi}\,(-\theta;\rho_{y,a}(\theta;\alpha))\cdot
b+\frac{2\theta}{1+4\theta^2}\,\frac{\partial z}{\partial
\xi}\,(-\theta;\rho_{y,a}(\theta;\alpha))\cdot b\,.
\end{split}
\end{equation*}
Finally, since
$\xi(-\theta;\rho_{y,a}(\theta;\alpha))=-a+\xi(\theta;\alpha)+\zeta(-\theta;\rho_{y,a}(\theta;\alpha))$
the system (\ref{eqIV.4.3}) is equivalent to
\begin{equation*}
\begin{aligned}
a&=\frac{-2\theta\,y}{1+4\theta^2}+\xi(\theta;\alpha)-\frac 1{1+4\theta^2}
\Big[\alpha_\xi+2\theta(x(\theta;\alpha)-\alpha_x)+2\theta\,z(-\theta;
\rho_{y,a}(\theta;\alpha))\\
&\hbox to 4cm{} +\frac{\partial
z}{\partial\xi}\,(-\theta;\rho_{y,a}(\theta;\alpha))\cdot
b\Big]+\zeta(-\theta;\rho_{y,a}(\theta;\alpha))+H_3\,b\cdot b\\
b&=\frac{-y}{1+4\theta^2}-\frac 1{1+4\theta^2}\,
\big[x(\theta;\alpha)-\alpha_x-2\theta\,\alpha_\xi+z(-\theta;\rho_{y,a}
(\theta;\alpha))\big]\\
&\hbox to 4cm{} -\frac{\partial \zeta}{\partial
\xi}\,(-\theta;\rho_{y,a}(\theta;\alpha))\cdot b+\frac{2\theta}{1+4\theta^2}\,
\frac{\partial z}{\partial \xi}\,(-\theta;\rho_{y,a}(\theta;\alpha))\cdot b+
H_4\,b\cdot b
\end{aligned}
\end{equation*}
where according to (\ref{eqIV.4.8}) $H_\ell $, $\ell =3,4$, are two matrices
which entries satisfy the following estimates
\begin{equation}\label{eqIV.4.11}
\v\partial ^A_{(a,b)}\,\partial ^B_y\,H^j _{\ell
,p,q}(\theta,y,\alpha,a,b)\v\leq
\begin{cases} C\, , &\textrm{if\enskip  $\v B\v=0$}\\
\frac C{\bra\theta\ket}\big (1+\frac{\varepsilon\,\bra\theta\ket}{\bra
x\ket^{2+\sigma_0}}\big)\,,&\textrm{if\enskip  $\v B\v=1$}\\
\frac{C_B\varepsilon}{\bra \theta\ket\,\bra x\ket^{\v
\alpha\v+1+\sigma_0}}\,\big(1+\frac{\bra\theta\ket}{\bra
x\ket}\big)\,,&\textrm{if\enskip  $\v B\v\geq 2$}\,.
\end{cases}
\end{equation}
Let us set
\begin{equation}\label{eqIV.4.12}\left\{
\begin{array}{ll}
\Phi_1(a,b)&=\frac{-2\theta\,y}{1+4\theta^2}+\xi(\theta;\alpha)-\frac
1{1+4\theta^2} \,\Big[\alpha_\xi+2\theta(x(\theta;\alpha)-\alpha_x)
+2\theta\,z(-\theta;\rho_{y,a}(\theta;\alpha))\\
&\hbox to 1,3cm{} +\frac{\partial z}{\partial
\xi}\,(-\theta;\rho_{y,a}(\theta;\alpha))\cdot
b\Big]+\zeta(-\theta;\rho_{y,a}(\theta;\alpha))+H_3\,b\cdot b\\
\Phi_2(a,b)&=\frac{-y}{1+4\theta^2}-\frac
1{1+4\theta^2}\,\big[x(\theta;\alpha)-\alpha_x-2\theta\,\alpha_\xi+z
(-\theta;\rho_{y,a}(\theta;\alpha))\big]\\
&\hbox to 1,3cm{} -\frac{\partial \zeta}{\partial
\xi}(-\theta;\rho_{y,a}(\theta;\alpha))\cdot b+\frac{2\theta}{1+4\theta^2}\,
\frac{\partial z}{\partial \xi}\,(-\theta;\rho_{y,a}(\theta;\alpha))\cdot
b+H_4\,b\cdot b\,.
\end{array}\right.
\end{equation}
Then the system (\ref{eqIV.4.3}) to be solved is equivalent in $E$ to the
equation $(\Phi_1(a,b),\Phi_2(a,b))=(a,b)$. Let us set
$\Phi(a,b)=(\Phi_1(a,b),\Phi_2(a,b))$. We shall use the fixed point theorem in
$E$ (see (\ref{eqIV.4.4})).

(i) $\Phi(E)\subset E$.

Let us recall that $(y+x(\theta;\alpha))\cdot \alpha_\xi\leq c_0\,\bra
y+x(\theta;\alpha)\ket\,\v\alpha_\xi\v$.

{\bf Case 1 } Assume that
\begin{equation}\label{eqIV.4.13}
x(\theta;\alpha)\cdot \alpha_\xi\geq 2\,c_0\,\bra
x(\theta;\alpha)\ket\,\v\alpha_\xi\v\,.
\end{equation}
It follows that
\begin{equation}\label{eqIV.4.14}
\v y\v\geq \frac{c_0}2\,\bra x(\theta;\alpha)\ket\geq \frac{c_0}2\,.
\end{equation}
Indeed one can write
\begin{equation*}
\begin{aligned}
2\,c_0\,\bra x(\theta;\alpha)\ket\,\v\alpha_\xi\v&\leq
(x(\theta;\alpha)+y)\cdot \alpha_\xi-y\cdot \alpha_\xi\\
&\leq c_0\,\bra y+x(\theta;\alpha)\ket\v\alpha_\xi\v+\v y\v\cdot
\v\alpha_\xi\v\\
&\leq c_0\,\bra x(\theta;\alpha)\ket\,\v\alpha_\xi\v+c_0\,\v y\v\cdot
\v\alpha_\xi\v+\v y\v\cdot \v\alpha_\xi\v\,.
\end{aligned}
\end{equation*}
Therefore $c_0\,\bra x(\theta;\alpha)\ket\leq 2\,\v y\v$. Here we have used
the inequality $\bra a+b\ket \leq \bra a\ket+\v b\v$.

Let $\theta^*\geq 0$ be such that $x(\theta^*;\alpha)\cdot \alpha_\xi=0$ (this
is possible since $\alpha_x\cdot \alpha_\xi\leq 0$ and $x(\theta;\alpha)\cdot
\alpha_\xi\geq \alpha_x\cdot \alpha_\xi+\theta\,\v\alpha_\xi\v^2\rightarrow
+\infty $ if $\theta\rightarrow +\infty )$. Then the point
$(x(\theta^*;\alpha),\alpha_\xi)$ is outgoing for $\theta\geq 0$ and
$\theta\leq 0$.

We can write by Proposition \ref{pIII.3.1},
\begin{equation*}
\begin{aligned}
x(\theta;\alpha)&=x(\theta-\theta^*;x(\theta^*;\alpha),\,\xi(\theta^*;\alpha))\\
&=x(\theta^*;\alpha)+2
(\theta-\theta^*)\,\xi(\theta-\theta^*;x(\theta^*;\alpha),\xi(\theta^*;\alpha))
+z(\theta-\theta^*,\cdots)\\
&=x(\theta^*;\alpha)+2(\theta-\theta^*)\,\xi(\theta;\alpha)+z(\theta-\theta^*;x
(\theta^*;\alpha),\,\xi(\theta^*;\alpha))\,.
\end{aligned}
\end{equation*}
It follows that
$$
x(\theta;\alpha)\cdot \alpha_\xi=2(\theta-\theta^*)\,\v\alpha_\xi\v^2+\CO
(\varepsilon)\,\v\theta-\theta^*\v+z(\theta-\theta^*;x(\theta^*,\alpha),\,
\xi(\theta^*;\alpha))\cdot \alpha_\xi\,.
$$
Since $\v z(\theta-\theta^*,\cdots )\v\leq
C\,\varepsilon\,\v\theta-\theta^*\v$ we deduce the estimate,
$$
2\,\v\theta-\theta^*\v\,\v\alpha_\xi\v^2\leq \v x(\theta,\alpha)\v\cdot
\v\alpha_\xi\v+C'\,\varepsilon\,\v\theta-\theta^*\v\,.
$$
Therefore if $\varepsilon$ is small enough we obtain
\begin{equation}\label{eqIV.4.15}
\v\theta-\theta^*\v\leq 5\,\v x(\theta,\alpha)\v\,.
\end{equation}
Now let us introduce
\begin{equation}\label{eqIV.4.16}
u(\theta)=x(\theta;\alpha)-\alpha_x-2\theta\,\alpha_\xi\,.
\end{equation}
We claim that
\begin{equation}\label{eqIV.4.17}
\v u(\theta)\v\leq C\,\varepsilon\,\bra\theta-\theta^*\ket\,.
\end{equation}
Indeed we have $u(0)=0$ and for all $\theta$ in $\R$,
$$
\dot u(\theta)=\dot
x(\theta;\alpha)-2\,\alpha_\xi=2(\xi(\theta;\alpha)-\alpha_\xi)+2\varepsilon\,b(x(\theta;\alpha))\cdot
\xi(\theta;\alpha)\,.
$$
It follows from Proposition \ref{pIII.4.1} that
$$
\v\dot u(\theta)\v\leq C\,\varepsilon\,.
$$
Now since $x(\theta^*;\alpha)\cdot \alpha_\xi=0$ it follows from Proposition
\ref{pIII.5.2} that
$$
x(\theta^*;\alpha)=\alpha_x+2\theta^*\,\alpha_\xi-z(-\theta^*;x(\theta^*;\alpha),\,\xi(\theta^*;\alpha))\,.
$$
This implies that $\v u(\theta^*)\v\leq C_1\,\varepsilon$. Now we write
$$
\v u(\theta)-u(\theta^*)\v\leq \Big\v\int^\theta_{\theta^*}\,\v \dot
u(s)\v\, ds\Big\v\leq C_2\,\varepsilon\,\v\theta-\theta^*\v
$$
and
$$
\v u(\theta)\v\leq C_1\,\varepsilon+C_2\,\varepsilon\,\v\theta-\theta^*\v\leq
C_3\,\varepsilon\,\bra\theta-\theta^*\ket\,.
$$
It follows then from (\ref{eqIV.4.14}) to (\ref{eqIV.4.17}) that,
\begin{equation*}
\begin{aligned}
\v x(\theta;\alpha)-\alpha_x-2\theta\,\alpha_\xi\v&\leq
C_3\,\varepsilon\,\bra\theta-\theta^*\ket\leq 5\,C_3\,\varepsilon\,\bra
x(\theta;\alpha)\ket\\
&\leq \frac{10\,C_3}{c_0}\,\varepsilon\,\v y\v\,.
\end{aligned}
\end{equation*}
Using (\ref{eqIV.4.12}) we see that
\begin{equation*}
\begin{split}
\Big\v\Phi_2(a,b)&+\frac y{1+4\theta^2}\Big\v\leq \underbrace{\frac{\v
x(\theta;\alpha)-\alpha_x-2\theta\,\alpha_\xi\v}{1+4\theta^2}}_{(1)}
+\underbrace{\frac{\v z(-\theta;\rho_{y,a})\v}{1+4\theta^2}}_{(2)}\\
&+\underbrace{\Big\v\frac{\partial \zeta}{\partial
\xi}\,(-\theta,\rho_{y,a})\Big\v}_{(3)}\cdot \v b\v +\frac 1{\bra
2\theta\ket}\,\underbrace{\Big\v\frac{\partial z}{\partial
\xi}\,(-\theta;\rho_{y,a})\Big\v}_{(4)}\cdot \v b\v+\underbrace{\Vert
H_4\Vert \cdot \v b\v^2}_{(5)}\,.
\end{split}
\end{equation*}
We have $(1)\leq C_4\,\varepsilon\,\frac{\v y\v}{\bra\theta\ket^2}$, $(2)\leq
\frac{C\,\varepsilon}{\bra\theta\ket^2}\leq \frac{C'\,\varepsilon\,\v
y\v}{\bra\theta\ket^2}$ since by (\ref{eqIV.4.14}) $\v y\v\geq \frac{c_0}2$.
Moreover
$$
(3)\leq C\,\varepsilon\,\v b\v\leq \frac{C'\,\varepsilon\,\v
y\v}{\bra\theta\ket^2}\,,\quad (4)\leq \frac{C\,\varepsilon}{\bra
2\theta\ket}\,\v b\v\leq \frac{C'\,\varepsilon\,\v y\v}{\bra\theta\ket^2}
$$
and, by (\ref{eqIV.4.11}), $(5)\leq C\,\frac{\v
y\v}{\bra\theta\ket^2}\,\frac{\v y\v}{\bra\theta\ket^2}$. Since $\v y\v\leq
\delta\,\bra\theta\ket$ it follows that $(5)\leq C\,\delta\,\frac{\v
y\v}{\bra\theta\ket^2}$.

Summing up we obtain
\begin{equation}\label{eqIV.4.18}
\Big\v\Phi_2(a,b)+\frac y{1+4\theta^2}\Big\v\leq
C\,(\varepsilon+\delta)\,\frac{\v y\v}{1+4\theta^2}
\end{equation}
so we take $\varepsilon,\delta$ so small that $C\,(\varepsilon+\delta)\leq
\sqrt\delta$.

Let us look to the term
$$
(II)=\Big\v\Phi_1(a,b)+\frac{2\theta\,y}{1+4\theta^2}\Big\v\,.
$$
We have
\begin{equation*}
\begin{split}
\Phi_1(a,b)&+\frac{2\theta\,y}{1+4\theta^2}=\underbrace{\xi(\theta;
\alpha)-\alpha_\xi}_{(1)}
-\underbrace{\frac{2\theta}{1+4\theta^2}\,
[x(\theta;\alpha)-\alpha_x-2\theta\,\alpha_\xi]}_{(2)}\\
&-\underbrace{\frac{2\theta}{1+4\theta^2}\,z(-\theta;\rho_{y,a})}_{(3)}
-\underbrace{\frac 1{1+4\theta^2}\,\frac{\partial z}{\partial
\xi}\,(-\theta;\rho_{y,a})\cdot
b}_{(4)}+\underbrace{\zeta(-\theta;\rho_{y,a})}_{(5)}+\underbrace{H_3\,b\cdot
b}_{(6)}\,.
\end{split}
\end{equation*}
We have $\v(1)\v\leq C\,\varepsilon$, $\v(2)\v\leq
\frac{C\,\varepsilon\,\theta\bra\theta-\theta^*\ket}{1+4\theta^2}\leq
C'\,\varepsilon$ (see (\ref{eqIV.4.17})), $\v(3)\v\leq
\frac{C\,\varepsilon}{\bra\theta\ket}$ (by Proposition \ref{pIII.3.2}),
$\v(4)\v\leq \frac{C\,\varepsilon}{\bra\theta\ket^2}\,\v b\v\leq
C''\,\varepsilon$, $\v(5)\v\leq C\,\varepsilon$, $\v(6)\v\leq C\,\delta^2$. It
follows that if $\varepsilon$ and $\delta$ are small enough we have,
\begin{equation}\label{eqIV.4.19}
(II)\leq C\,(\varepsilon+\delta^2)\leq C\,(\varepsilon+\delta^2)\,\frac
2{c_0}\,\inf (1,\,\v y\v)\leq \sqrt\delta\,\inf (1,\,\v y\v)
\end{equation}
since $\v y\v\geq \frac{c_0}{2}$. It follows from (\ref{eqIV.4.18}) and
(\ref{eqIV.4.19}) that $\Phi$ maps $E$ into $E$.

We show now that one can find a constant $k\in]0,1[$ such that
\begin{equation}\label{eqIV.4.20}
\v\Phi(a,b)-\Phi(a',b')\v\leq k\,\v (a,b)-(a',b')\v\,,\quad \forall
\,(a,b)\,,\enskip (a',b')\in E\,.
\end{equation}
We have
\begin{equation*}
\begin{split}
\v\Phi_1(a,b)-\Phi_1(a',b')\v\leq \underbrace{\frac{1}{1+4\theta^2}\,
\v z(-\theta;\rho_{y,a}(\theta;\alpha))-z(-\theta;\rho_{y,a'}
(\theta,\alpha))\v}_{(1)}\\
+\underbrace{\frac{1}{1+4\theta^2} \Big\v\,\frac{\partial z}{\partial
\xi}\,(-\theta;\rho_{y,a}(\theta;\alpha))\cdot b-\frac{\partial z}{\partial
\xi}\,(-\theta;\rho_{y,a'}(\theta;\alpha))\cdot b'}_{(2)}\Big\v\\
+\underbrace{\v\zeta(-\theta;\rho_{y,a}(\theta;\alpha))-\zeta(-\theta;
\rho_{y,a'}(\theta;\alpha))}_{(3)}\v +\v\underbrace{H_3\,b\cdot
b-H'_{3}\,b'\cdot b'}_{(4)}\v\,.
\end{split}
\end{equation*}
Since the point $(y+x(\theta;\alpha),\alpha_\xi)$ is outgoing we can use
Proposition \ref{pIII.3.2}, (\ref{eqIV.4.11}) and the fact that $\v b\v\leq
\frac{2\,\v y\v}{\bra\theta\ket^2}\leq 2\delta$, to write
\begin{equation*}
\begin{array}{ll}
(1)\leq \frac{C\,\theta}{1+4\theta^2}\,\varepsilon\,\v a-a'\v\,,\quad &(2)\leq
\frac{C\,\theta\,\varepsilon}{1+4\theta^2}\,(\v b-b'\v+\v a-a'\v)\\
(3)\leq C\,\varepsilon\,\v a-a'\v\, ,\quad &(4)\leq C\,\delta\,(\v a-a'\v+\v
b-b'\v)\,. 
\end{array}
\end{equation*}
It follows that
\begin{equation}\label{eqIV.4.21}
\v\Phi_1(a,b)-\Phi_1(a',b')\v\leq C\,(\varepsilon+\delta)(\v a-a'\v+\v
b-b'\v)\,.
\end{equation}
Now we have
\begin{equation*}
\begin{split}
\v\Phi_2(a,b)&-\Phi_2(a',b')\v\leq \Big\v\frac{\partial \zeta}{\partial
\xi}\,(-\theta;\rho_{y,a}(\theta;\alpha))\cdot b-\frac{\partial
\zeta}{\partial \xi}\,(-\theta;\rho_{y,a'}(\theta;\alpha))\cdot b'\Big\v\\
&+\frac{2\theta}{1+4\theta^2}\, \Big\v\frac{\partial
z}{\partial\xi}(-\theta;\rho_{y,a}(\theta;\alpha))\cdot b-\frac{\partial
z}{\partial \xi}\,(-\theta;\rho_{y,a'}(\theta;\alpha))\cdot b'\Big\v\\
&+\v H_4\,b\cdot b-H'_4\,b'\cdot b'\v+\frac
1{1+4\theta^2}\, \v
z((-\theta;\rho_{y,a}(\theta;\alpha))-z(-\theta;\rho_{y,a'}(\theta;\alpha))\v\,.
\end{split}
\end{equation*}
The same estimates as those used in the first case show that
\begin{equation}\label{eqIV.4.22}
\v\Phi_2(a,b)-\Phi_2(a',b')\v\leq C\,(\varepsilon+\delta)(\v a-a'\v+\v
b-b'\v)\,.
\end{equation}
Thus (\ref{eqIV.4.20}) is proved if $(\varepsilon+\delta)$ is small enough.
Then the fixed point theorem shows that the system (\ref{eqIV.4.3}) has a
unique solution in $E$.

{\bf Case 2 } Assume that
$$
x(\theta;\alpha)\cdot \alpha_\xi\leq 2\,c_0\,\bra
x(\theta;\alpha)\ket\,\v\alpha_\xi\v\,.
$$
It follows that we can apply Proposition \ref{pIII.5.2} with $(y,\eta)=\alpha$
which allows us to write
\begin{equation*}\left\{
\begin{array}{l}
x(\theta;\alpha)=\alpha_x+2\theta\,\alpha_\xi-z(-\theta;\rho(\theta;\alpha))\\
\xi(\theta;\alpha)=\alpha_\xi-\zeta(-\theta;\rho(\theta;\alpha))
\end{array}\right.
\end{equation*}
where $\rho(\theta;\alpha)=(x(\theta;\alpha),\,\xi(\theta;\alpha))$.

Using (\ref{eqIV.4.12}) we obtain the following expressions of $\Phi_1,\Phi_2$.
\begin{equation*}
\begin{aligned}
\Phi_1(a,b)&=\frac{-2\theta\,y}{1+4\theta^2}-\underbrace{\frac{2\theta}{1+4\theta^2}\,\big[
z((-\theta;\rho_{y,a}(\theta;\alpha))-z(-\theta;\rho(\theta;\alpha))
}_{(1)}]\\
&\quad-\underbrace{\frac 1{1+4\theta^2}\,\frac{\partial z}{\partial
\xi}\,(-\theta;\rho_{y,a}(\theta;\alpha))\cdot
b}_{(2)}
+\underbrace{\zeta(-\theta;\rho_{y,a}(\theta;\alpha))-\zeta(-\theta;\rho
(\theta;\alpha))}_{(3)}+\underbrace{H_3\,b\cdot b}_{(4)}\,.\\
\Phi_2(a,b)&=-\frac y{1+4\theta^2}-\underbrace{\frac
1{1+4\theta^2}\,\big[z(-\theta;\rho_{y,a}(\theta;\alpha))-z(-\theta;
\rho(\theta;\alpha))}_{(5)} \big]\\
&\quad-\underbrace{\frac{\partial \zeta}{\partial
\xi}\,(-\theta;\rho_{y,a}(\theta;\alpha))\cdot b}_{(6)}
+\underbrace{\frac{2\theta}{1+4\theta^2}\,\frac{\partial z}{\partial
\xi}\,(-\theta;\rho_{y,a}(\theta;\alpha))\cdot b}_{(7)}
+\underbrace{H_4\,b\cdot b}_{(8)}\,.
\end{aligned}
\end{equation*}
Let us show that
\begin{equation}\label{eqIV.4.23}\left\{
\begin{array}{l}
\big\v\Phi_1(a,b)+\frac{2\theta\,y}{1+4\theta^2}\big\v\leq \sqrt\delta\,\inf
(1,\v y\v)\\
\big\v\Phi_2(a,b)+\frac y{1+4\theta^2}\big\v\leq \sqrt\delta\,\frac{\v
y\v}{\bra\theta\ket^2}\,.
\end{array}\right.
\end{equation}
If $\v y\v\geq c_0$ then $\inf (1,\v y\v)\geq c_0$. It follows that
\begin{equation*}
\begin{aligned}
\v(1)\v&\leq \frac{2\theta}{1+4\theta^2}\,C\,\varepsilon\leq
\frac{C\,\varepsilon}{c_0}\,\inf (1,\v y\v)\,,\\
\v(2)\v&\leq C\,\varepsilon\leq 
\frac{C\,\varepsilon}{c_0}\,\inf (1,\v y\v)\,,\\
\v(3)\v&\leq C\,\varepsilon\,\v b\v\leq \frac{2C\,\varepsilon\,\v
y\v}{\bra\theta\ket^2}\leq 2C\,\varepsilon\,\delta\leq
\frac{2C\,\varepsilon\,\delta}{c_0}\,\inf (1,\v y\v)\,,\\
\v(4)\v&\leq \frac{C\,\delta}{c_0}\,\inf (1,\v y\v)\,,\\
\v(5)\v&\leq\frac{C\,\varepsilon}{\bra\theta\ket^2}\leq
\frac{C\,\varepsilon}{c_0}\,\frac{\v y\v}{\bra\theta\ket^2}\,,\\
\v(6)\v&\leq C\,\varepsilon\,\v b\v\leq C'\,\varepsilon\,\frac{\v
y\v}{\bra\theta\ket^2}\,,\\
\v(7)\v&\leq C\,\varepsilon\,\frac{\v y\v}{\bra\theta\ket^2}\,,\\
\v(8)\v&\leq C\,\delta\,\frac {\v y\v}{\bra\theta\ket^2}\,.
\end{aligned}
\end{equation*}
These estimates imply (\ref{eqIV.4.23}).

Assume now that $\v y\v\leq c_0$. It follows that for every $t$ in $[0,1]$
the point $(t\,y+x(\theta;\alpha),\xi(\theta;\alpha)-t\,a)$ is outgoing for
$\theta\leq 0$ ({\it i.e.} belongs to $\CS_-$). Indeed we have

$$
\xi(\theta;\alpha)-t\,a=\alpha_\xi+\CO (\varepsilon+\delta)
$$
$$
\v x(\theta;\alpha)\v\leq \v y\v+\v y+x(\theta;\alpha)\v\leq 1+\v
y+x(\theta;\alpha)\v
$$
so
\begin{equation*}
\begin{aligned}
(t\,y+x(\theta;\alpha))\cdot \alpha_\xi&\leq \v y\v\cdot
\v\alpha_\xi\v+2\,c_0\,\bra x(\theta;\alpha)\ket\,\v\alpha_\xi\v\\
&\leq 3\,c_0\,\bra x(\theta;\alpha)\ket\,\v\alpha_\xi\v\\
&\leq 6\,c_0\,\bra y+x(\theta;\alpha)\ket\,\v\alpha_\xi\v\,.
\end{aligned}
\end{equation*}
Then we have the following estimates.
\begin{equation*}
\begin{aligned}
\v(1)\v&\leq \int^1_0 \Big[\v y\v\,\Big\v\frac{\partial z}{\partial
y}\Big\v+\v a\v\,\Big\v\frac{\partial z}{\partial
\xi}\Big\v\Big](-\theta,t\,y+x(\theta;\alpha),\xi(\theta;\alpha)-ta)\,dt\\
\v(1)\v&\leq C\,\varepsilon\,\v y\v\quad \textrm{since}\quad \v a\v\leq
C'\,\v\ y\v\,.
\end{aligned}
\end{equation*}

By the same way we have $\v(3)\v\leq C\,\varepsilon\,\v y\v$.

 Moreover $\v
(2)\v\leq C\,\varepsilon\,\v y\v$, $\v(4)\v\leq C\,\delta\,\v y\v$ since $\v
b\v\leq \frac{2\,\v y\v}{\bra\theta\ket^2}\leq 2\delta$.

On the other hand we have
\begin{equation*}
\begin{array}{ll}
\v(5)\v\leq \frac{C\,\varepsilon\,\v y\v}{\bra\theta\ket^2}\,,\quad &\v
(6)\v\leq C\,\varepsilon\,\v b\v\leq \frac{C'\varepsilon\,\v
y\v}{\bra\theta\ket^2}\\
\v(7\v)\leq C\,\varepsilon\,\v b\v\leq C'\,\varepsilon\,\frac{\v
y\v}{\bra\theta\ket^2}\,,\quad &\v(8)\v\leq C\,\delta\,\frac{\v
y\v}{\bra\theta\ket^2}\,.
\end{array}
\end{equation*}
These estimates imply (\ref{eqIV.4.23}) since in this case $\v y\v=\inf (1,\v
y\v)$. Summing up we have proved that $\Phi$ maps $E$ into itself.

We show now the estimate (\ref{eqIV.4.22}). But it is easy to see that the
proof given in Case~1 works also in Case~2. Using again the fixed point
theorem we see that the system (\ref{eqIV.4.3}) has a unique solution in $E$.
This proves the points (i) and (ii) of Theorem \ref{tIV.4.2}.

To prove (iii) we use an induction on $\v A\v$ starting with $\v A\v=1$. Let us
set for fixed $(\theta,\alpha)$
\begin{equation}\label{eqIV.4.24}
\gamma_{y,\tilde a}=\Big(y+x(\theta,\alpha),\xi(\theta,\alpha)-\tilde
a(\theta,y,\alpha)+\frac{2\theta\,y}{1+4\theta^2}\Big)\,.
\end{equation}
Using (\ref{eqIV.4.12}) we see that $(\tilde a,\tilde b)$ satisfy the system
\begin{equation}\label{eqIV.4.25}
\begin{split}
\tilde
a&=\xi(\theta,\alpha)-\frac{\alpha_\xi+2\theta(x(\theta,\alpha)-\alpha_x)}{1+4\theta^2}
-\underbrace{\frac{2\theta}{1+4\theta^2}\,z(-\theta,\gamma_{y,\tilde
a})}_{(1)} -\underbrace{\frac 1{1+4\theta^2}\,\frac{\partial z}{\partial \xi}\,
(-\theta,\gamma_{y,\tilde a})
\frac 1{\bra\theta\ket}\,\tilde b}_{(2)}\\
&\quad +\underbrace{\frac 1{1+4\theta^2}\,\frac{\partial
z}{\partial \xi}\,
(-\theta,\gamma_{y,\tilde a})\, \frac{y}{1+4\theta^2}}_{(3)}
+\underbrace{\zeta(-\theta,\gamma_{y,\tilde a})}_{(4)}\\
&\quad +\underbrace{\frac 1{\bra\theta\ket^2}\,\tilde H_3\,\tilde b\cdot \tilde
b}_{(5)}
+\underbrace{\frac 2{\bra\theta\ket}\,\tilde H_3,\tilde b\,\frac
y{1+4\theta^2}}_{(6)}+\underbrace{\tilde H_3\,\frac y{1+4\theta^2}\,\frac
y{1+4\theta^2}}_{(7)}\,.
\end{split}
\end{equation}
\begin{equation}\label{eqIV.4.26}
\begin{split}
\tilde
b=&-\frac{\bra\theta\ket(x(\theta,\alpha)-\alpha_x-2\theta\,\alpha_\xi)}{1+4\theta^2}
-\underbrace{\frac{\bra\theta\ket}{1+4\theta^2}\,z(-\theta,\gamma_{y,\tilde a})
}_{(8)}\\
&-\underbrace{\frac{\partial \zeta}{\partial \xi}(-\theta,\gamma_{y,\tilde
a})\,\tilde b }_{(9)}+\underbrace{\frac{\partial \zeta}{\partial
\xi}\,(-\theta,\gamma_{y,\tilde
a})\,\frac{y\bra\theta\ket}{1+4\theta^2}}_{(10)}+\underbrace{\frac
1{\bra\theta\ket}\,
\tilde H_4\,\tilde b\cdot \tilde b}_{(11)}+\underbrace{2\,\tilde
H_4\,\tilde b\,\frac y{1+4\theta^2}}_{(12)}\\
&+\underbrace{\bra\theta\ket\,\tilde H_4\,\frac y{1+4\theta^2}\cdot
\frac{y}{1+4\theta^2}}_{(13)}+\underbrace{\frac{2\theta}{1+4\theta^2}\,\frac{\partial
z}{\partial
\xi}(-\theta,\gamma_{y,\tilde a})\,\tilde b}_{(14)}
-\underbrace{\frac{2\theta}{1+4\theta^2}\,\frac{\partial z}{\partial \xi}
(-\theta,\gamma_{y,\tilde
a})\,\frac{y\bra\theta\ket}{1+4\theta^2}}_{(15)}
\end{split}
\end{equation}
where for $j=3,4$ $\tilde H_j=H_j\big(\theta,y,\alpha,\tilde 
a-\frac{2\theta\,y}{1+4\theta^2}$, $\bra\theta\ket\,\tilde b-\frac
y{1+4\theta^2}\big)$ and $H_j$ satisfies (\ref{eqIV.4.11}). We claim that, for
$j=3,4$,
\begin{equation}\label{eqIV.4.27}
\v\partial _y[\tilde H_j]\v\leq C \Big(\v\partial _y\,\tilde
a\v+\bra\theta\ket\,\v\partial _y\,\tilde b\v+\frac
1{\bra\theta\ket}+\frac\varepsilon{\bra x\ket^{2+\sigma_0}}\Big )\,.
\end{equation}
Indeed, skipping the index $j$ for convenience, we have
$$
\frac{\partial }{\partial y_k}\,[\tilde H]=\frac{\partial H}{\partial 
y_k}+\frac{\partial H}{\partial a}\,\frac{\partial \tilde a}{\partial
y_k}-\frac{2\theta}{1+4\theta^2}\,\frac{\partial H}{\partial
a_k}+\bra\theta\ket\,\frac{\partial H}{\partial \tilde b}\,\frac{\partial
\tilde b}{\partial y_k}-\frac 1{1+4\theta^2}\,\frac{\partial H}{\partial
b_k}\,.
$$
Now we use (\ref{eqIV.4.11}). The first term in the right hand side is bounded
by $\frac {C}{\bra\theta\ket}+\frac{C\,\varepsilon}{\bra x\ket^{\sigma_0+2}}$,
the second by $C\,\v\nabla \tilde a\v$, the third by $\frac
C{\bra\theta\ket}$, the fourth by $C\,\bra\theta\ket \v\nabla _y\,\tilde
b\v$ and the last one by $\bra\theta\ket^{-2}$.

For $\ell \in\N$ let us introduce the following space
\begin{equation}\label{eqIV.4.28}
\CF_\ell  =\Big\{F\in C^\infty (\R^n\times\R^n) : \v\partial
^A_x\,\partial ^B_\xi\,F(x,\xi)\v\leq \frac{C_{AB}\,\varepsilon}{\bra x\ket^{\v
A\v+\ell +\sigma_0}}\,,\enskip \forall \,(x,\xi)\in\R^n\times\R^n\Big\}\,.
\end{equation}
For example $\zeta,\,\frac{\partial \zeta}{\partial \xi}\in \CF_1$,
$z,\,\frac{\partial z}{\partial \xi}\in\CF_0$ according to Proposition
\ref{pIII.3.2}. 

Let us set now,
\begin{equation}\label{eqIV.4.29}\left\{
\begin{array}{l}
g(y)=y+x(\theta;\alpha)\\
h(y)=\xi(\theta;\alpha)-\tilde
a(\theta,y,\alpha)+\frac{2\theta\,y}{1+4\theta^2}\,.
\end{array}\right.
\end{equation}
Then for $F\in\CF_\ell $ and $k=1,\ldots ,n$ we have,
\begin{equation}\label{eqIV.4.30}
\Big\v\frac{\partial }{\partial y_k}\,[F(g(y),h(y))]\Big\v\leq C\,\varepsilon
\Big (\frac 1{\bra x\ket^{\ell +1+\sigma_0}}+\frac 1{\bra x\ket^{\ell
+\sigma_0}\bra\theta\ket}+\v\nabla _y\,\tilde a\v\Big)
\end{equation}
where $x=y+x(\theta;\alpha)$.

Let us prove (iii) for $\v A\v=1$. We differentiate the equations
(\ref{eqIV.4.25}), (\ref{eqIV.4.26}) with respect to $y_k$ and we use
(\ref{eqIV.4.27}), (\ref{eqIV.4.30}) and the fact that $\v\tilde b\v\leq
\frac{3\,\v y\v}{\bra\theta\ket}\leq 3\delta\leq 1$. We have, with the
notations in (\ref{eqIV.4.25}), (\ref{eqIV.4.26}),
\begin{equation*}
\begin{array}{l}
\v\partial _{y_k}\,(1)\v+\v\partial _{y_k}\,(8)\v\leq C\,\varepsilon
\,\v\nabla _y\,\tilde a\v+C\,\varepsilon\, \frac 1{\bra
x\ket^{\sigma_0}}\,\frac 1{\bra\theta\ket}\Big (\frac 1{\bra x\ket}+\frac
1{\bra\theta\ket}\Big)\\
\v\partial _{y_k}\,(2)\v\leq C\,\varepsilon (\v\nabla _y\,\,\tilde
a\v+\v\nabla _y\,\tilde b\v)+\frac {C\,\varepsilon}{\bra
x\ket^{\sigma_0}\bra\theta\ket^3}\,\Big (\frac 1{\bra x\ket}+\frac
1{\bra\theta\ket}\Big)\\
\v\partial _{y_k}\,(3)\v\leq C\,\varepsilon \v\nabla _y\,\,\tilde
a\v+\frac {C\,\varepsilon}{\bra
x\ket^{\sigma_0}\bra\theta\ket^3}\,\Big (\frac 1{\bra x\ket}+\frac
1{\bra\theta\ket}\Big)\\
\v\partial _{y_k}\,(4)\v\leq C\,\varepsilon \, \v\nabla _y\,\,\tilde
a\v+\frac {C\,\varepsilon}{\bra
x\ket^{\sigma_0+1}}\,\Big (\frac 1{\bra x\ket}+\frac
1{\bra\theta\ket}\Big)\\
\v\partial _{y_k}\,(5)\v+\v\partial _{y_k}\,(6)\v+\v\partial _{y_k}\,(7)\v
\leq C\,\delta (\v\nabla \,\,\tilde a\v+\v\nabla \,\tilde
b\v)+\frac{C\,\varepsilon}{\bra
x\ket^{2+\sigma_0}\bra\theta\ket^2}+\frac{C\,\delta}{\bra\theta\ket^3}\\
\v\partial _{y_k}\,(9)\v\leq C\,\varepsilon (\v\nabla\,\tilde
a\v+\v\nabla \,\tilde b\v)+\frac {C\,\varepsilon}{\bra
x\ket^{\sigma_0+1}}\,\Big (\frac 1{\bra x\ket}+\frac
1{\bra\theta\ket}\Big)\\
\v\partial _{y_k}\,(10)\v\leq C\,\varepsilon \v\nabla\,\tilde
a\v+\frac {C\,\varepsilon}{\bra
x\ket^{\sigma_0}}\,\Big (\frac 1{\bra x\ket}+\frac
1{\bra\theta\ket}\Big)^2\\
\v\partial _{y_k}\,(11)\v+\v\partial _{y_k}\,(12)\v+\v\partial _{y_k}\,(13)\v
\leq C\,\delta\,(\v\nabla \,\tilde a\v+\v\nabla \,\tilde
b\v)+\frac{C\,\varepsilon}{\bra
x\ket^{2+\sigma_0}\bra\theta\ket}+\frac{C\,\delta}{\bra\theta\ket^2}\\
\v\partial _{y_k}\,(14)\v+\v\partial _{y_k}\,(15)\v\leq
C\,(\varepsilon+\delta)(\v\nabla _y\,\tilde a\v+\v\nabla\, \tilde 
b\v)+\frac{C\,\varepsilon}{\bra x\ket^{\sigma_0}}\Big (\frac 1{\bra x\ket}+\frac
1{\bra\theta\ket}\Big)^2\,.
\end{array}
\end{equation*}
It follows from (\ref{eqIV.4.25}), (\ref{eqIV.4.26}), that
$$
\v\nabla _y\,\tilde a\v+\v\nabla _y\,\tilde b\v\leq
C\,(\varepsilon+\delta)(\v\nabla \,\tilde a\v+\v\nabla \,\tilde
b\v)+\frac{C\,\varepsilon}{\bra x\ket^{\sigma_0}}\Big (\frac 1{\bra x\ket}+\frac
1{\bra\theta\ket}\Big)^2+\frac{C\,\delta}{\bra\theta\ket^2}
$$
which is (iii) for $\v A\v=1$.

Let us assume that (iii) is true for $1\leq \v A\v\leq k$ and let $\v
A\v=k+1\geq 2$.

We claim that for $F\in\CF_\ell $ and $\v B\v\leq k$ we have
\begin{equation}\label{eqIV.4.31}
\big\v\partial ^B_y[F(g(y),h(y))]\big\v\leq \frac{C_B\varepsilon}{\bra
x\ket^{\sigma_0}}\,\Big (\frac 1{\bra x\ket}+\frac
1{\bra\theta\ket}\Big)^{\v B\v+\ell }\,.
\end{equation}
Indeed the term we want to estimate is a finite sum of terms of the following
form (see Section \ref{ssVIII.1}).
\begin{equation}\label{eqIV.4.32}
R_{\beta\gamma\,s}=(\partial ^\beta_x\,\partial ^\gamma_\xi\,F)(g(y),h(y))
\prod^s_{j=1} \big(\partial ^{\ell _j}_y\,g\big)^{k'_j} \,\big(\partial ^{\ell
_j}_y\,h\big)^{k_j}
\end{equation}
where $1\leq s\leq \v B\v$, $1\leq \v\beta\v+\v\gamma\v\leq \v B\v$,
$\som^s_{j=1} k'_j=\beta$, $\som^s_{j=1} k_j=\gamma$, $\som^s_{j=1} (\v
k_j\v+\v k'_j\v)\,\ell _j=B$, $\ell _j\not =0$, $(k'_j,k_j)\not =0$,
$j=1,\ldots ,s$.

Let us write $\{1,\ldots ,s\}=I_1\cup I_2$ where
$$
I_1=\{j : \v\ell _j\v=1\}\,,\quad \v I_2\v=\{j : \v\ell _j\v\geq 2\}\,.
$$
For $j\in I_1$ we have $\partial ^{\ell _j}_y\,g_k=\CO(1)$, $\partial ^{\ell
_j}_y\,h_k=-\partial ^{\ell _j}_y\,\tilde a+\CO \big (\frac
1{\bra\theta\ket}\big)$. For $j\in I_2$ we have $\partial ^{\ell
_j}_y\,g_k\equiv 0$. Therefore the only terms which are present are those for
which $k'_j=0$. It follows that $\som^s_{j=1} k'_j=\som_{j\in I_1}
k'_j=\beta$. Moreover $\partial ^{\ell _j}_yh=-\partial ^{\ell _j}_y\,\tilde
a$.
It follows from these facts and the definition of $\CF_\ell $ that
$$
\v R_{\beta\gamma\,s}\v\leq \frac{C_{\beta\gamma}\varepsilon}{\bra
x\ket^{\v\beta\v+\ell +\sigma_0}}\left[\frac{(C\,\varepsilon)^{\som^s_{j=1}
\v k_j\v}} {\bra
x\ket^{{\sigma_0}\som^s_{j=1}\v k_j\v}}
\Big (\frac 1{\bra x\ket}+\frac 1{\bra \theta\ket}\Big)^{\som^s_{j=1}\v
k_j\v(\v\ell _j\v+1)}+\frac{(C\,\delta)^{\som^s_{j=1}\v k_j\v}}{\bra
\theta\ket^{\som^s_{j=1}\v k_j\v(\v\ell _j\v+1)}}\right]\,.
$$
Now $\v\beta\v=\som^s_{j=1} \v k'_j\v=\som^s_{j=1} \v k'_j\v\,\v\ell _j\v$
since $\v\ell _j\v=1$ in $I_1$ and $k'_j=0$ in $I_2$. It follows that
$$
\v R_{\beta\gamma\,s}\v\leq \frac{C\,\varepsilon}{\bra x\ket^{\sigma_0}}\Big
(\frac 1{\bra x\ket}+\frac 1{\bra\theta\ket}\Big)^{\som^s_{j=1} \v
k'_j\v\,\v\ell '_j\v+\ell }\Big
(\frac 1{\bra x\ket}+\frac 1{\bra\theta\ket}\Big)^{\som^s_{j=1} \v
k_j\v\,\v\ell _j\v}\,.
$$
The result follows then, since $\som^s_{j=1} (\v k_j\v+\v k'_j\v)\v\ell _j\v=\v
B\v$.

On the other hand, for $F\in\CF_\ell $ and $\v A\v=k+1\geq 2$, we have
\begin{equation}\label{eqIV.4.33}
\big\v\partial ^A_y(F(g(y),h(y))\big\v\leq C_0\varepsilon\,\v\partial
^A_y\,\tilde a\v+\frac{C_A\varepsilon}{\bra x\ket^{\sigma_0}}\Big
(\frac 1{\bra x\ket}+\frac 1{\bra\theta\ket}\Big)^{\v A\v+\ell }\,.
\end{equation}
Indeed
$$
\partial ^A_y(F(g(y),h(y))-\som^n_{k=1} \Big (\frac{\partial F}{\partial
x_k}\,(g(y),h(y))\,\partial ^A_y\,g_k+\frac{\partial F}{\partial
\xi_k}\,(g(y),h(y))\,\partial ^A_y\,h_k)
$$
is a finite sum of terms $R_{\beta\gamma\,s}$ given by (\ref{eqIV.4.32}) where
$1\leq s\leq \v A\v$, $2\leq \v\beta\v+\v\gamma\v\leq \v A\v$,
$\som^s_{j=1}\,k'_j=\beta$, $\som^s_{j=1}\,k_j=\gamma$, $\som^s_{j=1} (\v
k_j\v+\v k'_j\v)\,\ell _j=A$, $\v\ell _j\v\geq 1$, $\v k_j\v+\v
k'_j\v\geq 1$, $j=1,\ldots ,s$.
Since $\v\beta\v+\v\gamma\v\geq 2$ we have $\v\ell _j\v\leq \v A\v-1$ so the
term $R_{\beta\gamma\,s}$ is bounded, using the induction, by
$$
\frac{C_A\varepsilon}{\bra x\ket^{\sigma_0}}\Big
(\frac 1{\bra x\ket}+\frac 1{\bra\theta\ket}\Big)^{\v A\v+\ell \,.}
$$
On the other hand, since $\v A\v\geq 2$ we have $\partial ^A_y\,g_k=0$ and
$\partial ^A_y\,h_k=-\partial ^A_y\,\tilde a_k$. By (\ref{eqIV.4.28}) we have
$\big\v\frac{\partial F}{\partial \xi_k}\,(g(y),h(y))\big\v\leq
C\,\varepsilon$ so (\ref{eqIV.4.33}) is proved. Note that $C$ is independent
of $A$.

Let us now prove the last step of the induction. We apply $\partial ^A_y$ to
both members of (\ref{eqIV.4.25}), (\ref{eqIV.4.26}). Then we estimate each
term of the right hand side using (\ref{eqIV.4.31}), (\ref{eqIV.4.33}). We
obtain
\begin{equation*}
\begin{array}{l}
\v\partial ^A_y\,(1)\v+\v\partial ^A_\gamma\,(8)\v\leq
C_0\varepsilon\,\v\partial ^A_y\,\tilde a\v+\frac{C_A\varepsilon}{\bra
x\ket^{\sigma_0}}\,\frac 1{\bra\theta\ket}\,\Big
(\frac 1{\bra x\ket}+\frac 1{\bra\theta\ket}\Big)^{\v A\v}\,.\\
\v\partial ^A_y\,(4)\v\leq C_0\varepsilon\,\v\partial ^A_y\,\tilde a\v
+\frac{C_A\varepsilon}{\bra x\ket^{\sigma_0}}\,\Big
(\frac 1{\bra x\ket}+\frac 1{\bra\theta\ket}\Big)^{\v A\v+1}\,.\\
\v\partial ^\alpha_y\,(2)\v+\v\partial ^A_y\,(3)\v\leq C_0\varepsilon\, \big
(\v\partial ^A_y\,\tilde a\v+\v\partial ^A_y\,\tilde
b\v\big)+\frac{C_A\varepsilon}{\bra x\ket^{\sigma_0}}\,\frac 1{\bra
\theta\ket^3}\,\Big
(\frac 1{\bra x\ket}+\frac 1{\bra\theta\ket}\Big)^{\v A\v}\,.\\
\v\partial ^A_y\,(9)\v+\v\partial ^A_y\,(10)\v\leq C_0\varepsilon (\big
(\v\partial ^A_y\,\tilde a\v+\v\partial ^A_y\,\tilde b\v\big
)+\frac{C_A\varepsilon}{\bra x\ket^{\sigma_0}}\,\Big
(\frac 1{\bra x\ket}+\frac 1{\bra\theta\ket}\Big)^{\v A\v+1}\,.
\end{array}
\end{equation*}
Finally
\begin{equation*}
\begin{split}
\v\partial ^A_y\,(5)\v+\v\partial ^A_y\,(6)\v&+\v\partial
^A_y\,(7)\v+\v\partial ^A_y\,(11)\v+\v\partial ^A_y\,(12)\v+\v\partial
^A_y\,(13)\v\leq C_0\delta\big (\v\partial ^A_y,\tilde a\v+\v\partial
^A_y\,\tilde b\v\big )\\
&+\frac{C_A\varepsilon}{\bra x\ket^{\sigma_0}}\,\Big
(\frac 1{\bra x\ket}+\frac 1{\bra\theta\ket}\Big)^{\v
A\v+1}+\frac{C_A\delta}{\bra\theta\ket^{\v A\v+1}\,.}
\end{split}
\end{equation*}
If $\varepsilon+\delta$ is small enough (compared with a finite number of fixed
constants) we can absorb the term $C_0(\varepsilon+\delta)(\v\partial
^A_y\,\tilde a\v+\v\partial ^A_y\,\tilde b\v)$ by the left hand side and we
obtain the estimate given in (iii).

Let us now prove (iv). First of all since the point $(y+x(\theta,\alpha)$,
$\eta+\xi(\theta,\alpha))$ belongs to $\CS_-$ and $-\theta\leq 0$ we deduce
from Proposition \ref{pIII.3.2} that $\v\partial
^\gamma_\xi\,\xi(-\theta,y+x(\theta,\alpha),\eta+\xi(\theta,\alpha))\v+\v\partial
^\gamma_\xi\,x(-\theta,y+x(\theta,\alpha),\eta+\xi(\theta,\alpha)\v\leq
C_\gamma\,\bra\theta\ket$. It follows then, from (\ref{eqIV.4.2}) and Lemma
\ref{lIV.2.1} that
$$
\Big\v\frac{\partial q_k}{\partial
\xi_j}\,(\eta,a,b,g_j)\Big\v+\Big\v\frac{\partial q_k}{\partial b_\ell
}\,(\eta,a,b,g_j)\Big\v\leq C\,\bra\theta\ket\,.
$$
Therefore we will have,
$$
\v q_k(\eta,a,b,g_j)-q_k(-a,a,0,g_j)\v\leq C\,\bra\theta\ket\,(\v\eta+a\v+\v
b\v)\,.
$$
Now if $\v\eta\v\leq \sqrt\delta$ we will have $\v\eta+a\v+\v b\v\leq
5\sqrt\delta$. On the other hand (\ref{eqIV.2.7}) shows that, when
$\v\eta\v\leq \sqrt\delta<\frac 12\,c_0$,
$$
q_k(-a,a,0,g_j)=\frac{\partial g_j}{\partial \xi_k}\,(-a)=\Big (\frac{\partial
\xi_j}{\partial \xi_k}-i\,\frac{\partial x_j}{\partial x_k}\Big
)(-\theta,y+x(\theta,\alpha),-a+\xi(\theta,\alpha))\,.
$$
By Corollary \ref{cIII.3.3} we have
$$
q_k(-a,a,0,g_j)=(1+2i\theta)\,\delta_{jk}+\CO(\varepsilon\,\bra\theta\ket)\,.
$$
Finally we obtain,
$$
\v q_k(\eta,a,b,g_j)-(1+2i\theta)\,\delta_{jk}\v\leq
C\,(\varepsilon+\sqrt\delta)\,\bra\theta\ket\,,
$$
which is precisely the claim of point (iv).

The last point (v) can be easily deduced from Proposition \ref{pIII.3.2} and
Lemma \ref{lIV.2.1}. This ends the proof of Theorem \ref{tIV.4.2}. \cqfd

We continue the proof of Theorem \ref{tIV.1.2} in the case (\ref{eqIV.4.1}).
We follow basically the same method as in Section \ref{ssIV.3} with small
changes. For the convenience of the reader we give some details.

 Let us set,
\begin{equation}\label{eqIV.4.34}\left\{
\begin{array}{l}
\CO=\big\{(\theta,(y,\eta))\in\R_+\times\R^{2n}\ : \v
y\v<\delta\bra\theta\ket,(y+x(\theta,\alpha))\cdot 
\alpha_\xi\leq c_1\,\bra
y+x(\theta\,\alpha)\ket\v\alpha_\xi\v;\v\eta\v\leq \sqrt\delta\big\}\\
\tilde \Lambda=\big\{\alpha\in T^*\R^n : \frac 12\leq \v\alpha_\xi\v\leq
2,\,\alpha_x\cdot \alpha_\xi\leq -c_0\,\bra\alpha_x\ket\v\alpha_\xi\v\}~.
\end{array}\right.
\end{equation}
We consider families $(f(\cdot ,\alpha))_{\alpha\in\tilde\Lambda}$ of
function on $\CO$. 
\begin{definition}\sl\label{dIV.4.3}
We say that $(f(\cdot ,\alpha))_{\alpha\in\tilde \Lambda}$ belongs to $\CH$ if
\begin{itemize}
\item[(i)] for all $\alpha\in\tilde \Lambda$, $(\theta,y,\eta)\mapsto
f(\theta,y,\eta,\alpha)$ is $C^\infty $ on $\CO$.
\item[(ii)] For every $A,B$ in $\N^n$ there exists $C_{AB}>0$  such that
$$
\sup_{(\theta,y,\alpha)\in\CO\times\tilde \Lambda}\, \v\partial ^A_y\,\partial
^B_\eta\,f(\theta,y,\eta,\alpha)\v\leq C_{AB}\,.
$$
\end{itemize}
\end{definition}
\begin{remark}\sl\label{rIV.4.4}

\noindent
1) $\CH$ is closed under multiplication and derivation with respect to
$(y,\eta)$.

2) If we set, with the notation (\ref{eqIV.2.10}),
\begin{equation*}
\begin{aligned}
f(\theta,y,\eta,\alpha)&=\frac
1{\bra\theta\ket}\,v_j(\theta,y+x(\theta,\alpha),\eta+\xi(\theta,\alpha))\\
&=\frac
1{\bra\theta\ket}\,\big[\xi_j(-\theta,y+x(\theta,\alpha),\eta+\xi(\theta,\alpha))
-\alpha^j_\xi-i(x(-\theta,y+x(\theta,\alpha),
\eta+\xi(\theta,\alpha))-\alpha^j_x)\big]
\end{aligned}
\end{equation*}
then $(f(\cdot ,\alpha))_{\alpha\in\tilde \Lambda}\in\CH$. This is a
consequence of Proposition \ref{pIII.3.2}.
\end{remark}
\begin{definition} (\bf (Lagrangian ideals) \sl\label{dIV.4.5}
The Lagrangian ideal $\CJ$ is defined as the set of families $F=(F(\cdot
,\alpha))_{\alpha\in\tilde \Lambda}$ which can be written as
$$
F(\theta,y,\eta,\alpha)=\sum^n_{j=1} f_j(\theta,y,\eta,\alpha)\, \frac
1{\bra\theta\ket}\,v_j(\theta,y+x(\theta,\alpha),\eta+\xi(\theta,\alpha))
$$
for all $(\theta,y,\eta)$ in $\CO$ and $\alpha$ in $\tilde \Lambda$, where
$(f(\cdot ,\alpha))_{\alpha\in\tilde \Lambda}\in\CH$.
\end{definition}
\begin{example}\sl\label{eIV.4.6}
Let us set
$$
F(\theta,y,\eta,\alpha)=\eta_k-\psi_k(\theta,y,\alpha)
$$
with
\begin{equation} \label{eqIV.4.35}
\psi_k(\theta,y,\alpha)=-(a_k+i\,b_k)(\theta,y,\alpha)\,,
\end{equation}
where $a_k$, $b_k$ are those given in Theorem \ref{tIV.4.2}.

Then $(F((\cdot ,\alpha))_{\alpha\in\Lambda}\in\CJ$.

Indeed if $\v\eta\v\leq \sqrt\delta<\frac 12\,c_0$ then $\chi_0(\eta)=1$ so it
follows from (\ref{eqIV.4.2}) and Theorem \ref{tIV.4.2} that we have
$$
v_j(\theta,y+x(\theta,\alpha),\xi(\theta,\alpha)+\eta)=\sum^n_{k=1}
q_k(\eta,a,b,g_j)(\eta_k-\psi_k(\theta,y,\alpha))\,.
$$
Since $q_k(\eta,a,b,g_j)=(1+2i\theta)\,\delta_{jk}+\CO
((\varepsilon+\sqrt\delta)\bra\theta\ket)$ (by (v)) it follows that the matrix
$(q_k(\eta,a,b,g_j))^{-1}=(d_{jk}(\theta,y,\eta,\alpha))_{1\leq j,k\leq n}$
exists. Moreover $(\bra\theta\ket\,d_{jk}(\cdot ,\alpha))_{\alpha\in\tilde
\Lambda} \in \CH$. Now we have
$$
\eta_k-\psi_k(\theta ,y,\alpha)=\sum^n_{j=1}
\bra\theta\ket\,d_{jk}(\theta,y,\alpha)\cdot \frac
1{\bra\theta\ket}\,v_j(\theta,y+x(\theta,\alpha),\eta+\xi(\theta,\alpha))\,.
$$
This proves our claim.
\end{example}
\begin{lemma}\sl \label{lIV.4.7}
For $F$ and $G$ in $\CJ$ let us define
$$
\{F,G\}(\theta,y,\eta,\alpha)=\sum^n_{j=1} \Big (\frac{\partial F}{\partial
\eta_j}\,\frac{\partial G}{\partial y_j}-\frac{\partial F}{\partial
y_j}\,\frac{\partial G}{\partial \eta_j}\Big )(\theta,y,\eta,\alpha)\,.
$$
Then $\{F,G\}\in\CJ$.
\end{lemma}

{\bf Proof } Since
$v_j(\theta,y+x(\theta,\alpha),\eta+\xi(\theta,\alpha)=u_j\circ\chi_{-\theta}(y,\eta)$
where $u_j(x,\xi,\alpha)=\xi_j-\alpha^j_\xi-i(x_j-\alpha^j_x)$ and
$\chi_{-\theta}(y,\eta)=(x(-\theta,y,\eta),\xi(-\theta,y,\theta))$ is the
symplectic map defined by the flow we have
$$
\{v_j,v_k\}(\theta,y,\eta,\alpha)=\{u_j,u_k\}(\chi_{-\theta}(y,\eta))=0
$$
 because $\{u_j,u_k\}\equiv 0$.

Let $F=\som \,f_j\,\frac 1{\bra\theta\ket}\,v_j$, $G=\som\,g_k\,\frac
1{\bra\theta\ket}\,v_k$ be two elements of $\CJ$ with $f_j\in\CH$, $g_k\in\CH$.
Then a straightforward computation and the Remark \ref{rIV.4.4} give the
conclusion (see the proof of Lemma \ref{lIV.3.10}). 
\cqfd

\begin{lemma}\sl\label{lIV.4.8}
Let $R=(R(\cdot ,\alpha))_{\alpha\in\tilde \Lambda}\in \CJ$ and assume that
$R(\cdot ,\alpha)$ does not depend on $\eta$. Then for every $N\in\N$ one can
find $C_N>0$ such that for every $(\theta,y)$ in $\tilde \Omega_\delta$ and
$\alpha$ in $\tilde \Lambda$ we have
$$
\v R(\theta,y,\alpha)\v\leq C_N\,\v \Im \psi(\theta,y,\alpha)\v^N\,.
$$
\end{lemma}
{\bf Proof } We are going to show by induction on $N\geq 1$ that we can write
\begin{equation}\label{eqIV.4.36}
R(\theta,y,\alpha)=\som_{0<\v \gamma\v<N}
h_\gamma(\theta,y,\alpha)(\eta-\psi)^\gamma+\som_{\v\gamma\v=N}
g_\gamma(\theta,y,\alpha,\eta)(\eta-\psi)^\gamma
\end{equation}
where $(h_\gamma(\cdot ,\alpha))_{\alpha\in\tilde \Lambda}$ and
$(g_\gamma(\cdot ,\alpha))_{\alpha\in\tilde \Lambda})$ belong to $\CH$.

For $N=1$ the first sum in the right hand side of (\ref{eqIV.4.36}) is empty
and by assumption we have
$$
R(\theta,y,\alpha)=\som^n_{j=1} f_j(\theta,y,\eta,\alpha)\,\frac
1{\bra\theta\ket}\,v_j(\theta,y+x(\theta,\alpha),\eta+\xi(\theta,\alpha))\,.
$$
Using Theorem \ref{tIV.4.2} we obtain since $\chi_0(\eta)=1$ when
$\v\eta\v\leq \sqrt\delta$,
$$
R(\theta,y,\alpha)=\som^n_{k=1} \Big (\som^n_{j=1} \frac
1{\bra\theta\ket}\,f_j(\theta,y,\eta,\alpha)\,q_k(\eta,a,b,g_j)\Big)
(\eta_k-\psi_k(\theta,y,\alpha))\,.
$$
Since $f_j$ and $\frac 1{\bra\theta\ket}\,q_k$ belong to $\CH$ this shows that
(\ref{eqIV.4.36}) is true when $N=1$. Assume now it is true up to the order
$N$. We can apply Lemma \ref{lIV.2.1} to the function
$$
\tilde
g_\gamma(\theta,y,\eta,\alpha)=\chi_0(\eta)\,g_\gamma(\theta,y,\eta,\alpha)\,,
\quad \v\gamma\v=N
$$
with $z_j=-\psi_j(\theta,y,\alpha)$. It follows that
\begin{equation}\label{eqIV.4.37}
\tilde g_\gamma(\theta,y,\eta,\alpha)=\som^n_{k=1} q_k(\eta,a,b,\tilde
g_\gamma)(\eta_k-\psi_k(\theta,y,\alpha))+r(a,b,\tilde g_\gamma)\,.
\end{equation}
For the $q'_ks$ and $r$ we have the estimates (\ref{eqIV.2.2}). Let us set
\begin{equation}\label{eqIV.4.38}\left\{
\begin{array}{l}
h_\gamma(\theta,y,\alpha)=r(a(\theta,y,\alpha),b(\theta,y,\alpha),\tilde
g_\gamma(h,y,\cdot ,\alpha))\\
g_\gamma(\theta,y,\eta,\alpha)=q_k(\eta,a(\theta,y,\alpha),b(\theta,y,\alpha),\tilde
g_\gamma(\theta,y,\cdot ,\alpha))\,.
\end{array}\right.
\end{equation}
It follows from (\ref{eqIV.2.2}) and Theorem \ref{tIV.4.2} that
$(h_\gamma(\cdot ,\alpha))_\alpha$ and $(g_{\gamma\,k}(\cdot ,\alpha))_\alpha$
belong to $\CH$. Using (\ref{eqIV.4.36}) at the level $N$ and
(\ref{eqIV.4.37}), (\ref{eqIV.4.38}) we deduce that (\ref{eqIV.4.36}) holds at
the level $N+1$.

Now let us take in (\ref{eqIV.4.36}) $\eta=(\Re \psi+s\, \Im
\psi)(\theta,y,\alpha)$ when $s\in[0,1]$. Then the same argument as in the end
of the proof of Lemma \ref{lIV.3.11} gives the result. \cqfd

\begin{corollary}\sl\label{cIV.4.9}
For every $N\in\N$ there exists a constant $C_N>0$ such that
$$
\Big\v\Big(\frac{\partial \psi_j}{\partial y_k}-\frac{\partial
\psi_k}{\partial y_j}\Big)(\theta,y,\alpha)\Big\v\leq C_N\,\v \Im
\psi(\theta,y,\alpha)\v^N
$$
for every $(\theta,y)$ in $\tilde \Omega_\delta$ and $\alpha$ in $\tilde
\Lambda$.
\end{corollary}

{\bf Proof }  Identic to the proof of Corollary \ref{cIV.3.12}. \cqfd

Now we go back to the original coordinates
$$
x=y+x(\theta,\alpha)\,,\quad \xi=\eta+\xi(\theta,\alpha)
$$
and we set for $k=1,\ldots ,n$,
\begin{equation}\label{eqIV.4.39}
\Phi_k(\theta,x,\alpha)=\psi_k(\theta,x-x(\theta,\alpha),\alpha)
=\xi_k(\theta,\alpha)-(a_k+i\,b_k)(\theta,x-x(\theta,\alpha),\alpha)
\end{equation}
where $a_k$, $b_k$ have been described in Theorem \ref{tIV.4.2}.

Then we can state the following result.

\begin{theorem}\sl\label{tIV.4.10}
We can write for $(\theta,x)\in\Omega_\delta$ and
$\v\xi-\xi(\theta,\alpha)\v\leq \sqrt \delta$,
\begin{itemize}
\item[(i)] $\ds{\xi_k-\Phi_k(\theta,x,\alpha)=\som^n_{j=1}
e_{jk}(\theta,x,\xi,\alpha)\,v_j(\theta,x,\xi)}$

where $e_{jk}$ are smooth functions which satisfy
\item[(ii)] $\ds{\v\partial ^\ell _\theta\,\partial
^A_x\,e_{jk}(\theta,x,\xi,\alpha)\v\leq \frac{C_A}{\bra\theta\ket}}$ for all
$\alpha\in\N$ and $\ell =0,1$.

Moreover we have for $(\theta,x)$ in $\Omega_\delta$ and
$\v\xi-\xi(\theta,\alpha)\v<\sqrt\delta$,
\item[(iii)] $\v\Phi_k(\theta,x,\alpha)-\alpha_\xi\v\leq
C_0\,(\varepsilon+\sqrt\delta)$.
\item[(iv)] $\ds{\Big\v\Im
\Phi_k(\theta,x,\alpha)-\frac{x_k-x_k(\theta,\alpha)}{1+4\theta^2}\Big\v\leq
\sqrt\delta \,\frac{\v x-x(\theta,\alpha)\v}{\bra\theta\ket^2}}$.
\item[(v)] $\v\Phi(\theta,x,\alpha)\v\leq C_0$,

$\ds{\v\partial ^A_x\,\Phi(\theta,x,\alpha)\v\leq C_A\Big (\frac
1{\bra\theta\ket^{\v A\v+1}}+\frac 1{\bra x\ket^{\v A\v+1+\sigma_0}}\Big )}$ if
$A\in\N^n$,\enskip $\v A\v\geq 1$.
\item[(vi)] $\Phi_k(\theta,x(\theta,\alpha),\alpha)=\xi_k(\theta,\alpha)$.
\item[(vii)] $\ds{\Big\v\Big(\frac{\partial \Phi_k}{\partial
x_j}-\frac{\partial \Phi_j}{\partial x_k}\Big )(\theta,x,\alpha)\Big\v\leq
C_N\Big (\frac 1{\bra x\ket^{\frac 3 2}}+\frac{1}{\bra\theta\ket^{\frac 3
2}}\Big)\,\frac{\v x-x(\theta,\alpha)\v^N}{\bra\theta\ket^{2N}}}$, \enskip
$N\in\N$

where the constants $C_A, C_0, C_N$ are independent of
$(\theta,x,\xi,\alpha)$. 
\end{itemize}
\end{theorem}

{\bf Proof } (i) follows immediately from the computations made in Example
\ref{eIV.4.6} as well as (ii). The point (iii) is obvious since
$\xi(\theta,\alpha)=\alpha_\xi+O(\varepsilon)$ and $\v a_k\v+\v b_k\v=\CO
(\sqrt\delta)$ by (\ref{eqIV.4.4}). Then (iv) follows from Theorem
\ref{tIV.4.2} (ii) as well as (v). The point (vi) is obvious since $a_k=b_k=0$
when $y=0$. To prove (vii) we use Corollary \ref{cIV.4.9}, (\ref{eqIV.4.35}),
(\ref{eqIV.4.39}) and (iv) of Theorem \ref{tIV.4.10}. We obtain
$$
\Big\v\Big (\frac{\partial \Phi_k}{\partial x_j}-\frac{\partial
\Phi_j}{\partial x_k}\Big )(\theta,x,\alpha)\Big\v\leq C_N\, \frac{\v
x-x(\theta,\alpha)\v^N}{\bra\theta\ket^{2N}}\,.
$$
Now by (v) of the present theorem we have,
$$
\v\partial _j\,\Phi(\theta,x,\alpha)\v\leq C\Big (\frac
1{\bra\theta\ket^2}+\frac 1{\bra x\ket^2}\Big )\,.
$$
Writing $\v a\v=\v a\v^{\frac 3 4}\,\v a\v^{\frac 14}$ and using the above
estimates we obtain (vii). \cqfd

\begin{definition}\sl\label{dIV.4.11}
Let $(F(\cdot ,\alpha))_{\alpha\in\tilde \Lambda}$ be a family of $C^\infty $
functions for $(\theta,x)$ in $\Omega_\delta$ and
$\v\xi-\xi(\theta,\alpha)\v<\sqrt\delta$. We shall say that $F\in
\CJ_{(x,\xi)}$ if we can write
$$
F(\theta,x,\xi,\alpha)=\som^n_{j=1} f_j(\theta,x,\xi,\alpha)\,\frac
1{\bra\theta\ket}\,v_j(\theta,x,\xi,\alpha)
$$
where
$$
\v\partial ^A_x\,\partial ^B_\xi\,f_j(\theta,x,\xi,\alpha)\v\leq C_{AB}
$$
for all $(\theta,x)$ in $\Omega_\delta$, 
$\v\xi-\xi(\theta,\alpha)\v<\sqrt\delta$, $\alpha\in\tilde \Lambda$ where
$C_{AB}$ is independent of $(\theta,x,\xi,\alpha)$.
\end{definition}

Then exactly as in Lemma \ref{lIV.4.7} $\CJ_{(x,\xi)}$ is closed under the 
Poisson bracket in $(x,\xi)$ and we have the analogue of Lemma \ref{lIV.4.8}.
In fact
$\CJ_{(x,\xi)}$ is just the image of $\CJ$ under the diffeomorphism
$x=y+x(\theta,\alpha),\xi=\eta+\xi(\theta,\alpha)$. Then we have
\begin{theorem}\sl\label{tIV.4.12}
With $\Phi$ defined in (\ref{eqIV.4.39}) we have for $k=1,\ldots ,n$,
$$
\Big (-\frac{\partial p}{\partial
x_k}\,(x,\Phi(\theta,x,\alpha))-\frac{\partial
\Phi_k}{\partial \theta}\,(\theta,x,\alpha)-\frac{\partial \Phi_k}{\partial
x}\,(\theta,x,\alpha)\,\frac{\partial p}{\partial
\xi}\,(x,\Phi(\theta,x,\alpha))\Big)_\alpha\in\CJ_{(x,\xi)}\,.
$$
\end{theorem}

{\bf Proof } We follow word by word the proof of Lemma \ref{lIV.3.16} and
Corollary \ref{cIV.3.17}. Let us just sketch the proof.

We begin by the proof of the theorem when $\Phi(\theta,x,\alpha)$ is replaced
by $\xi$. By Theorem \ref{tIV.4.10} (i) we have
$$
\xi_k-\Phi_k(\theta,x,\alpha)=\som\,
e_{jk}\,(\theta,x,\xi,\alpha)\,v_j(\theta,x,\xi)\,.
$$
Then we set $x(-\theta,x,\xi)=X$, $\xi(-\theta,x,\xi)=\Xi$ that is
$x(\theta,X,\Xi)=x$, $\xi(\theta,X,\Xi)=\xi$. The above identity reads
$$
\xi_k(\theta,X,\Xi)-\Phi_k(\theta,x(\theta,X,\Xi),\alpha)=\som^n_{k=1}
e_{jk}(\theta,x(\theta,X,\Xi),\xi(\theta,X,\Xi),\alpha)\,u_j(X,\Xi)\,.
$$
We differentiate this equality with respect to $\theta$ using the equations of
the flow given by (\ref{eqIII.1.2}). Then we use Theorem \ref{tIV.4.10} (ii)
and we come back to the original coordinates $(x,\xi)$. Finally we write
$\xi=\xi-\Phi(\theta,x,\alpha)+\Phi(\theta,x,\alpha)$ as in the proof of
Corollary \ref{cIV.3.17}. Details are left to the reader. \cqfd
\begin{corollary}\sl\label{cIV.4.13}
For every $N\in\N$ one can find $C_N>0$ such that
$$
\Big\v-\frac{\partial p}{\partial
x_k}\,(x,\Phi(\theta,x,\alpha))-\frac{\partial \Phi_k}{\partial
\theta}\,(\theta,x,\alpha)-\frac{\partial \Phi_k}{\partial
x}\,(\theta,x,\alpha)\,\frac{\partial p}{\partial
\xi}\,(x,\Phi(\theta,x,\alpha))\Big\v
\leq C_N \Big (\frac{\v x-x(\theta,\alpha)\v}{\bra\theta\ket^2}\Big)^N\,.
$$
\end{corollary}

{\bf Proof } Use Theorem \ref{tIV.4.12} and Lemma \ref{lIV.4.8} in the
coordinates $(x,\xi)$. \cqfd

We are ready now to define the phase $\varphi $, as we did in Proposition
\ref{pIV.3.19} for the outgoing points, but we find here a slight problem.
Indeed if we look to formula (\ref{eqIV.3.48})  we see that $\varphi $ is
defined by mean of $\Phi(\theta,s\,x+(1-s)\,x(\theta,\alpha))$, $s\in[0,1]$.
In the  present case when $\theta\geq 0$, $\Phi(\theta,z,\alpha)$ is defined
for $z\cdot \alpha_\xi\leq c_0\,\bra z\ket\,\v \alpha_\xi\v$ and $\v
z-x(\theta,\alpha)\v\leq \delta\,\bra\theta\ket$ and it is easily seen that
the point $z=s\,x+(1-s)\,x(\theta,\alpha)$ does not satisfy these conditions.
Therefore we have to modify the expression of $\varphi $ in (\ref{eqIV.3.48}) 
to take care of this problem. We split the discussion into several cases
giving in each of them a different expression of $\varphi $. Our purpose is to
prove the following result.

Let us set,
\begin{equation}\label{eqIV.4.40}
\CO_\delta=\Big\{(\theta,x)\in\R_+\times\R^n : x\cdot \alpha_\xi\leq \frac
{c_0}{10}\, \bra x\ket\,\v\alpha_\xi\v\,,\enskip \v x-x(\theta,\alpha)\v\leq
\frac\delta {10}\,\bra\theta\ket\Big\}\,.
\end{equation}
\begin{proposition}\sl\label{pIV.4.14}
There exists a smooth function $\varphi =\varphi (\theta,x,\alpha)$ defined on
$\CO_\delta$ such that,
\begin{itemize}
\item[(i)]$\ds{\varphi (0,x,\alpha)=(x-\alpha_x)\,\alpha_\xi+\frac i2\,
(x-\alpha_x)^2+\frac 1{2i}\,\alpha^2_\xi+\CO(\v x-\alpha_x\v^N)}$,\enskip
$\forall \,N\in\N$.

For every $N\in\N$ there exists $C_N>0$ such that
\item[(ii)] $\ds{\Big\v\frac{\partial \varphi }{\partial
x}\,(\theta,x,\alpha)-\Phi(\theta,x,\alpha)\Big\v\leq C_N \Big (\frac{\v
x-x(\theta,\alpha)\v}{\bra\theta\ket}\Big)^N }$
\item[(iii)] $\ds{\Big\v\frac{\partial \varphi }{\partial
\theta}\,(\theta,x,\alpha)+p\Big (x,\,\frac{\partial \varphi }{\partial x}\,
(\theta,x,\alpha)\Big)\Big\v\leq C_N \Big (\frac{\v
x-x(\theta,\alpha)\v}{\bra\theta\ket}\Big)^N }$
uniformly in $(\theta,x,\alpha)$.

Moreover
\item[(iv)] $\ds{\Big\v\frac{\partial \varphi }{\partial
x}\,(\theta,x,\alpha)-\alpha_\xi\Big\v\leq C\,(\varepsilon+\sqrt\delta)}$.
\item[(v)] $\ds{\v\partial ^A_x\,\varphi (\theta,x,\alpha)\v\leq C_A}$,
$\forall \,A\in\N^n$, \enskip $\v A\v\geq 1$.
\item[(vi)] $\ds{\Big\v\Im \varphi (\theta,x,\alpha)-\frac 12\,\frac{\v
x-x(\theta,\alpha)\v^2}{1+4\theta^2}+\frac 12\,\alpha^2_\xi \Big\v\leq
C\,(\varepsilon+\sqrt\delta)\,\frac{\v
x-x(\theta,\alpha)\v^2}{\bra\theta\ket^2}}$.
\end{itemize}
\end{proposition}

We split the proof  into several cases which are summarized in the following
figures.

{\bf Case 1 } (See Figure 1). Let $(\theta,x)\in\CO_\delta$ be such that
$$
x(\theta,\alpha)\cdot \alpha_\xi\leq \frac{c_0}{3}\,\bra
x(\theta,\alpha)\ket\,\v\alpha_\xi\v\enskip \textrm{and}\enskip \v
x-x(\theta,\alpha)\v\leq \bra x(\theta,\alpha)\ket\,.
$$

\begin{lemma}\sl \label{lIV.4.15}
In the case 1 we have $s\,x+(1-s)\,x(\theta,\alpha)\in\Omega_\delta$, for $s\in
[0,1]$, that is
$$
(s\,x+(1-s)\,x(\theta,\alpha))\cdot \alpha_\xi\leq c_0\,\bra
s\,x+(1-s)\,x(\theta,\alpha)\ket\,\v\alpha_\xi\v
$$
and
$$
\v s\,x+(1-s)\,x(\theta,\alpha)-x(\theta,\alpha)\v\leq
\delta\,\bra\theta\ket\,.
$$
\end{lemma}

{\bf Proof } We use the following elementary lemma.
\begin{lemma}\sl\label{lIV.4.16}
Let $a,b\in\R^n$ be such that $\v a-b\v\leq \bra a\ket$. Then for all $s$ in
$[0,1]$
$$
(1-s) \v a\v +s\,\v b\v\leq \sqrt 2\, \bra (1-s)\,a+s\,b\ket\,.
$$ 
\end{lemma}

{\bf Proof } Since $\v a-b\v^2\leq \v a\v^2+1$ we have $2a\cdot b\geq -1$, it
follows that
\begin{equation*}
\begin{split}
\v (1-s)\,a+s\,b\v^2=(1-s)^2\,\v a\v^2+2s(1-s)\,a\cdot b+s^2\,\v b\v^2\geq
(1-s)^2\, \v a\v^2\\
+s^2\, \v b\v^2-s(1-s)\geq \frac 12\,((1-s)\, \v
a\v +s\, \v b\v))^2-\frac 12\,.
\end{split}
\end{equation*}
Therefore $\bra (1-s)\,a+s\,b\ket^2\geq \frac 12 ((1-s)\v a\v+s\v b\v)^2$.
\cqfd

  

Let us now apply Lemma \ref{lIV.4.16} to $a=x(\theta,\alpha)$, $b=x$. Using
our hypotheses we obtain
\begin{equation*}
\begin{split}
(s\,x+(1-s)\,x(\theta,\alpha))\cdot \alpha_\xi&\leq \frac{c_0}3\,(s\,\bra x\ket
+(1-s)\bra x(\theta,\alpha)\ket)\,\v\alpha_\xi\v\\
&\leq \frac{c_0}3\,(s+s\,\v x\v+(1-s)+(1-s)\v x(\theta,\alpha)\v)\,\v
\alpha_\xi\v\\
&\leq \frac{c_0}3\,(1+\sqrt 2\,\bra s\,x+(1-s)\,
x(\theta,\alpha)\ket)\cdot \v \alpha_\xi\v\\
&\leq c_0\,\bra s\,x+(1-s)\,x(\theta,\alpha)\ket\,\v\alpha_\xi\v\,,
\end{split}
\end{equation*}
On the other hand $\v s\,x+(1-s)\,x(\theta,\alpha)-x(\theta,\alpha)\v=s\,\v
x-x(\theta,\alpha)\v\leq \delta\,\bra \theta\ket$. \cqfd

In the set defined in case 1 we can therefore define $\varphi $ by the same
formula as in Proposition \ref{pIV.3.19}. We set
\begin{equation}\label{eqIV.4.41}
\varphi (\theta,x,\alpha)=\int^1_0 (x-x(\theta,\alpha))\cdot
\Phi(\theta,s\,x+(1-s)\,x(\theta,\alpha),\alpha)\,ds
+\theta\,p(\alpha)+\frac 1{2i}\,\alpha^2_\xi\,.
\end{equation}
The proof of the  points (i) to (vi) is
exactly the same as the corresponding points in Proposition \ref{pIV.3.19}
using Theorem \ref{tIV.4.10}. \cqfd

{\bf Case 2 } Let $(\theta,x)\in\CO_\delta$ be such that
$$
x(\theta,\alpha)\cdot \alpha_\xi\leq \frac{c_0}3\,\bra
x(\theta,\alpha)\ket\,\v\alpha_\xi\v\enskip \textrm{and}\enskip \v
x-x(\theta,\alpha)\v\geq \frac 12\,\v x(\theta,\alpha)\v\,.
$$
(See Figure 1).

In this set we have
\begin{equation}\label{eqIV.4.42}\left\{
\begin{array}{l}
\v x(\theta,\alpha)\v\leq 2\,\v x-x(\theta,\alpha)\v \leq
\delta\,\bra\theta\ket\,,\\
\v x\v\leq 3\,\v x-x(\theta,\alpha)\v\leq \delta\,\bra\theta\ket\,.
\end{array}\right.
\end{equation}
Moreover for $y\in[0,x]\cup [0,x(\theta,\alpha)]$ the point $(\theta,y)$ 
belongs to the set $\Omega_\delta$ on which $\Phi$ is defined. Indeed if
$s\in[0,1]$ we have $s\,x\cdot \alpha_\xi\leq s\,\frac{c_0}2\,\bra
x\ket\,\v\alpha_\xi\v$ and $\v s\,x-x(\theta,\alpha)\v\leq s\,\v
x-x(\theta,\alpha)\v+(1-s)\v x(\theta,\alpha)\v\leq 2\,\v
x-x(\theta,\alpha)\v\leq \delta\, \bra\theta\ket$ by (\ref{eqIV.4.42}). On
the other hand, $s\,x(\theta,\alpha)\cdot \alpha_\xi\leq s\,\frac{c_0}3\,\bra
x(\theta,\alpha)\ket\,\v\alpha_\xi\v\leq c_0\,\bra
x(\theta,\alpha)\ket\,\v\alpha_\xi\v$ and $\v
s\,x(\theta,\alpha)-x(\theta,\alpha)\v=(1-s)\,\v x(\theta,\alpha)\v\leq
\delta\,\bra\theta\ket$. Therefore we can define the phase $\varphi $ by the
following formula.
\begin{equation}\label{eqIV.4.43}
\varphi (\theta,x,\alpha)=\int^1_0 x\cdot
\Phi(\theta,s\,x,\alpha)\,ds-\int^1_0 x(\theta,\alpha)\cdot
\Phi(\theta,s\,x(\theta,\alpha),\alpha)\,ds+\theta\,p(\alpha)+\frac
1{2i}\,\alpha^2_\xi\,.
\end{equation}
Let us show that $\varphi $ satisfies the conditions of Proposition
\ref{pIV.4.14}. It follows from Theorem \ref{tIV.4.2} and (\ref{eqIV.4.39})
that
$$
\Phi(0,z,\alpha)=\alpha_\xi+i\,(z-\alpha_x)+\CO(\v z-\alpha_x\v^N)\,.
$$
Therefore
\begin{equation*}
\begin{aligned}
\varphi (0,x,\alpha)&=\int^1_0 \big (x\cdot
\alpha_\xi+i\,x(s\,x-\alpha_x)-\alpha_x\,\alpha_\xi-i\,\alpha_x(s\,\alpha_x-\alpha_x)\big
\v x\v\,\CO(\v s\,x-\alpha_x\v^N)\\
&\qquad +\v
\alpha_x\v\,\CO(\v s\,x-\alpha_x\v^N))\,ds+\frac 1{2i}\,\alpha^2_\xi\\
&=(x-\alpha_x)\cdot \alpha_\xi+\frac i 2\, (x-\alpha_x)^2+\frac
1{2i}\,\alpha^2_\xi+\CO(\v x-\alpha_x\v^N)
\end{aligned}
\end{equation*}
because $\v x\v\leq \v x-\alpha_x\v+\v\alpha_x\v$, $\v s\,x-\alpha_x\v\leq \v
x-\alpha_x\v$ and $\v \alpha_x\v\leq \v x-\alpha_x\v$. Thus (i) is proved. Now 
$$
\frac{\partial \varphi }{\partial x_j}\,(\theta,x,\alpha)=\int^1_0
\Phi_j(\theta,s\,x,\alpha)\,ds+\som^n_{k=1} \int^1_0 s\,x_k\,\frac{\partial
\Phi_k}{\partial x_j}\,(\theta,s\,x,\alpha)\,ds\,.
$$
Using Theorem \ref{tIV.4.10}, (vii) we obtain
\begin{equation*}
\begin{split}
\Big\v\frac{\partial \varphi }{\partial x_j}\,(\theta,x,\alpha)
-\int^1_0 \Phi_j(\theta,s\,x,\alpha)\,ds
-\int^1_0
s\,\frac d{ds}\,(\Phi_j(\theta,s\,x,\alpha))\,ds\Big\v\\
\leq C_N \int^1_0 s\,\v x\v \Big (\frac 1{\bra s\,x\ket^{\frac 3 2}}+\frac
1{\bra\theta\ket^{\frac 32}}\Big )\,\frac{\v
s\,x-x(\theta,\alpha)\v^N}{\bra\theta\ket^{2N}}\,ds\,.
\end{split}
\end{equation*}
Now $s\,\v x\v \leq \bra s\,x\ket^{\frac 32}$ and by (\ref{eqIV.4.42}), $s\,\v
x\v\leq \delta\,\bra\theta\ket$, $\v s\,x-x(\theta,\alpha)\v\leq s\,\v
x-x(\theta,\alpha)\v+(1-s)\v x(\theta,\alpha)\v\leq 2\,\v
x-x(\theta,\alpha)\v$. Therefore the right hand side of the above inequality
is bounded by
$C_N\,\frac{\v x-x(\theta,\alpha)\v^N}{\bra\theta\ket^{2N}}$. Integrating by
parts in the second integral of the left hand side we obtain (ii).

Now we have
\begin{equation*}
\begin{split}
\frac{\partial \varphi }{\partial \theta}&=\underbrace{\int^1_0
x\,\frac{\partial
\Phi}{\partial \theta}\,(\theta,s\,x,\alpha)\,ds}_{(1)}-\underbrace{\int^1_0
\dot x(\theta,\alpha)\cdot
\Phi(\theta,s\,x(\theta,\alpha),\alpha)\,ds}_{(2)}\\
&-\underbrace{\int^1_0 x(\theta,\alpha)\cdot \frac{\partial
\Phi}{\partial
\theta}\,(\theta,s\,x(\theta,\alpha),\alpha)\,ds}_{(3)}-\som^n_{k,\ell
=1}\underbrace{ \int^1_0 x_k(\theta,\alpha)\,\frac{\partial \Phi_k}{\partial
x_\ell }\,(\theta,s\,x(\theta,\alpha),\alpha)
\,s\,\dot x_\ell (\theta,\alpha)\,ds}_{(4)}+p(\alpha)\,.
\end{split}
\end{equation*}
We use Corollary \ref{cIV.4.13} to write
$$
(1)=\som^n_{k=1} \int^1_0 x_k\,\frac{\partial p}{\partial
x_k}\,(s\,x,\Phi(\theta,s\,x,\alpha))\,ds-\som^n_{k,\ell =1} \int^1_0
x_k\,\frac{\partial \Phi_k}{\partial x_\ell }\,(\theta,s\,x,\alpha)
\cdot \frac{\partial p}{\partial \xi_\ell
}\,(s\,x,\Phi(\theta,s\,x,\alpha))\,ds+R_0
$$
with
$$
\v R_0\v\leq \int^1_0 \v x\v\,\frac{\v
s\,x-x(\theta,\alpha)\v^N}{\bra\theta\ket^{2N}}\,.
$$
By (\ref{eqIV.4.42}) we have $\v x\v\leq \delta\,\bra\theta\ket$\,; since
$\bra\theta\ket^{2N-1}\geq \bra\theta\ket^N$ if $N\geq 1$, and $\v
s\,x-x(\theta,\alpha)\v\leq s\,\v x-x(\theta,\alpha)\v+(1-s)\v
x(\theta,\alpha)\v\leq 2\,\v x-x(\theta,\alpha)\v$, we obtain $\v R_0\v\leq
C_N\,\frac{\v x-x(\theta,\alpha)\v^N}{\bra\theta\ket^N}$ if $N\geq 1$. Now $\v
x-x(\theta,\alpha)\v\leq \delta\,\bra\theta\ket$ so the same estimate is valid
for $N=0$. Finally
\begin{equation}\label{eqIV.4.44}
\v R_0\v\leq C_N\,\frac {\v x-x(\theta,\alpha)\v^N}{\bra\theta\ket^N}\,,\quad
\forall \,N\geq 0\,.
\end{equation}
Using Theorem \ref{tIV.4.10} (vii) we obtain
$$
(1)=-\int^1_0 \frac d{ds}\,(p(s\,x,\Phi(\theta,s\,x,\alpha)))\,ds+R_1
$$
where $R$ satisfies (\ref{eqIV.4.44}). Therefore
\begin{equation}\label{eqIV.4.45}
\big\v(1)+p(\Phi(\theta,x,\alpha)-p(0,\Phi(\theta,0,\alpha))\big\v\leq
C_N\, \frac{\v x-x(\theta,\alpha)\v^N}{\bra\theta\ket^N}\,.
\end{equation}
Let us look the term (4)\,; we use Theorem \ref{tIV.4.10} (vii) again and we
obtain
$$
(4)=\som^n_{\ell =1} \int^1_0 s\,\dot x_\ell (\theta,\alpha)\,\frac
d{ds}\,(\Phi_\ell (\theta,s\,x(\theta,\alpha),\alpha))\,ds+R_2\,,
$$
with
$$
\v R_2\v\leq \int^1_0 s\,\v x(\theta,\alpha)\v\Big (\frac 1{\bra
s\,x(\theta,\alpha)\ket^{\frac 32}}+\frac 1{\bra\theta\ket^{\frac 32}}\Big
)\Big (\frac{\v x(\theta,\alpha)\v^N}{\bra\theta\ket^{2N}}\Big
)(s-1)^N\,ds\,.
$$
It follows from (\ref{eqIV.4.42}) that $R_2$ satisfies (\ref{eqIV.4.44}).
Therefore integrating by parts in the above integral we obtain
$$
(4)=-\int^1_0 \dot x(\theta,\alpha)\cdot
\Phi(\theta,s\,x(\theta,\alpha),\alpha))\,ds+\dot
x(\theta,\alpha)\,\Phi(\theta,x(\theta,\alpha),\alpha)+\CO\Big (\frac{\v
x-x(\theta,\alpha)\v^N}{\bra\theta\ket^N}\Big )\,.
$$
Using Theorem \ref{tIV.4.10}, (vi) and the Euler identity we can write
$$
\dot
x(\theta,\alpha)\,\Phi(\theta,x(\theta,\alpha),\alpha)=\xi(\theta,\alpha)\cdot
\frac{\partial p}{\partial
\xi}\,(x(\theta,\alpha),\xi(\theta,\alpha))=2p(\alpha)\,.
$$
It follows that
\begin{equation}\label{eqIV.4.46}
\v(2)+(4)-2p(\alpha)\v\leq C_N\,\frac{\v
x-x(\theta,\alpha)\v^N}{\bra\theta\ket^N}\,.
\end{equation}
Now, by Corollary \ref{cIV.4.13} we have
\begin{equation*}
\begin{split}
(3)=&-\int^1_0 x(\theta,\alpha)\cdot \frac{\partial p}{\partial
x}\,(s\,x(\theta,\alpha),\Phi(\theta,s\,x(\theta,\alpha),\alpha))\,ds\\
&-\som^n_{\ell
,k=1} \int^1_0 x_k(\theta,\alpha)\cdot 
\frac{\partial \Phi_k}{\partial x_\ell
}\,(\theta,s\,x(\theta,\alpha),\alpha)\,\frac{\partial p}{\partial \xi_\ell
}\,(s\,x(\theta,\alpha),\Phi(\theta,s\,x(\theta,\alpha),\alpha)\,ds+R_3
\end{split}
\end{equation*}
where
$$
\v R_3\v\leq C_N\,\v x(\theta,\alpha)\v \int^1_0 \frac{\v
x(\theta,\alpha)\v^N\,\v s-1\v^N}{\bra\theta\ket^{2N}}\,ds\,.
$$
If $N\geq 1$ we have $\frac{\v x(\theta,\alpha)\v}{\bra\theta\ket^{2N}}\leq
\frac 1{\bra\theta\ket^N}$ so $R_3$ satisfies (\ref{eqIV.4.44}) using
(\ref{eqIV.4.42}).

Using again Theorem \ref{tIV.4.10} (vii) we obtain
$$
(3)=-\int^1_0 \frac d{ds}\,
\big[p(s\,x(\theta,\alpha),\Phi(\theta,s\,x(\theta,\alpha),\alpha)\big]\,ds+\CO
\Big (\frac{\v x-x(\theta,\alpha)\v^N}{\bra\theta\ket^N}\Big )
$$
so
$$
\big\v(3)+p(x(\theta,\alpha),\Phi(\theta,x(\theta,\alpha),\alpha)-p(0,\Phi(\theta,0,\alpha)\big\v\leq
C_N\,\frac{\v x-x(\theta,\alpha)\v^N}{\bra\theta\ket^N}\,.
$$
Finally we obtain
\begin{equation}\label{eqIV.4.47}
\big\v(3)+p(\alpha)-p(0,\Phi(\theta,0,\alpha))\big\v\leq C_N\,\frac{\v
x-x(\theta,\alpha)\v^N}{\bra\theta\ket^N}\,.
\end{equation}
Since $\frac{\partial \varphi }{\partial
\theta}\,(\theta,x,\alpha)=(1)-(2)-(3)-(4)+p(\alpha)$ we deduce from
(\ref{eqIV.4.45}) to (\ref{eqIV.4.47}) that
$$
\Big\v\frac{\partial \varphi }{\partial
\theta}\,(\theta,x,\alpha)+p(x,\Phi(\theta,x,\alpha))\Big\v\leq C_N\,
\frac{\v x-x(\theta,\alpha)\v^N}{\bra\theta\ket^N}\,.
$$
Using the point (ii) already proved in Proposition \ref{pIV.4.14} we obtain
the point (iii).

The point (iv) follows easily from (ii) since by (\ref{eqIV.4.39}) and Theorem
\ref{tIV.4.10} we have $\Phi(\theta,x,\alpha)=\alpha_\xi+\CO (\varepsilon+\sqrt
\delta)$.

Let us prove (v). To bound $\partial ^A_x\,\varphi $ when $\v A\v\geq 1$ we
have to bound, according to (\ref{eqIV.4.43}) the quantities $(1)=\int^1_0
s^{\v A'\v}(\partial ^{A'}_x\Phi)(\theta,s\,x,\alpha)\,ds$, $\v A'\v\leq \v
A\v-1$ and $(2)=\int^1_0 \v x\v\,\v\partial
^A_x\,\Phi(\theta,s\,x,\alpha)\v\,ds$. Using Theorem \ref{tIV.4.10},
(v), we see easily that (1) is uniformly bounded and
$$
\v(2)\v\leq C_A\,\v x\v\,\frac 1{\bra\theta\ket^{\v A\v}}+C_A\,\v x\v \int^1_0
\frac{ds}{\bra s\,x\ket^{\v A\v+1+\sigma_0}}\,.
$$
By (\ref{eqIV.4.42}) we have $\v x\v\leq 2\delta\,\bra\theta\ket\leq
2\delta\,\bra\theta\ket^{\v A\v}$ since $\v A\v\geq 1$ and setting $t=\v
x\v\,s$ in the integral above we see that (2) is uniformly bounded in
$(\theta,x,\alpha)$. This shows (v). Finally by (\ref{eqIV.4.43}),
$$
\Im \varphi (\theta,x,\alpha)=\int^1_0 x\cdot \Im
\Phi(\theta,s\,x,\alpha)\,ds-\int^1_0 x(\theta,\alpha)\cdot \Im
\Phi(\theta,s\,x(\theta,\alpha),\alpha)\,ds-\frac 12\,\alpha^2_\xi\,.
$$
Using Theorem \ref{tIV.4.10}, (iv) and (\ref{eqIV.4.42}) we obtain (vi).
This completes the proof of Proposition \ref{pIV.4.14} in the case 2. \cqfd

{\bf Case 3 } We consider here the case where
\begin{equation}\label{eqIV.4.48}
(\theta,x)\in \CO_\delta\enskip \textrm{and}\enskip x(\theta,\alpha)\cdot
\alpha_\xi>\frac{c_0}3\,\bra x(\theta,\alpha)\ket\,\v \alpha_\xi\v\,.
\end{equation}
(See Figure 2).

Let us recall that we are dealing in this Section IV.4 with the case where
$\alpha_x\cdot \alpha_\xi\leq -c_0\,\bra\alpha_x\ket\,\v\alpha_\xi\v$, (see
(\ref{eqIV.4.1})). 

(1) The continuous function $t\mapsto x(t,\alpha)\cdot \alpha_\xi$ is then
strictly negative for $t=0$ and strictly positive for $t=\theta$. It follows
that
\begin{equation}\label{eqIV.4.49}
\textrm{there exists}\enskip \theta^*\in]0,\theta[\enskip \textrm{depending
only on}\enskip \alpha\enskip \textrm{such that}\enskip
x(\theta^*,\alpha)\cdot \alpha_\xi=0\,.
\end{equation}
Then we have the following Lemma.
\begin{lemma}\sl\label{lIV.4.17}
\noindent 
\begin{itemize}
\item[(i)] $\ds{\frac 32\, \v\theta-\theta^*\v\,\v\alpha_\xi\v\leq \v
x(\theta,\alpha)-x(\theta^*,\alpha)\v\leq 3\,\v
\theta-\theta^*\v\,\v\alpha_\xi\v}$,

\item[(ii)] $\ds{\v\theta-\theta^*\v\geq \frac{c_0}{50}}$,

\item[(iii)] $\v x-x(\theta,\alpha)\v\leq \v x-x(\theta^*,\alpha)\v+\v
x(\theta^*,\alpha)-x(\theta,\alpha)\v\leq 5\,\v x-x(\theta,\alpha)\v$,
\item[(iv)] $K_1\,\bra\theta\ket\leq
\bra\theta^*\ket\leq K_2\,\bra\theta\ket$.
\end{itemize}
\end{lemma}

{\bf Proof } It follows from (\ref{eqIV.4.48}) and Definition \ref{dIII.2.2}
that the point $\rho^*=(x(\theta^*,\alpha),\xi(\theta^*,\alpha))$ belongs to
$\CS_+\cap\CS_-$. By the group property and Proposition \ref{pIII.3.1} we have
$$
x(\theta,\alpha)=x(\theta-\theta^*,\rho^*)=x(\theta^*,\alpha)+2(\theta-\theta^*)\,
\alpha_\xi+\CO (\varepsilon\,\v\theta-\theta^*\v)+\CO (\varepsilon)\,.
$$
It follows that
\begin{equation}\label{eqIV.4.50}
x(\theta,\alpha)-x(\theta^*,\alpha)=2(\theta-\theta^*)\,\alpha_\xi+\CO
(\varepsilon\,\v\theta-\theta^*\v)+\CO (\varepsilon)\,.
\end{equation}
Now we deduce from (\ref{eqIV.4.48}) and (\ref{eqIV.4.49}) that
$$
(x(\theta,\alpha)-x(\theta^*,\alpha))\cdot \alpha_\xi>\frac{c_0}3\,\bra
x(\theta,\alpha)\ket\,\v\alpha_\xi\v\geq \frac{c_0}3\,\v\alpha_\xi\v\geq
\frac{c_0}6
$$
since $\frac 12 \leq \v\alpha_\xi\v\leq 2$. So by (\ref{eqIV.4.50}),
$$
2(\theta-\theta^*)\,\v\alpha_\xi\v^2\geq
\frac{c_0}6-C_1\,\varepsilon\,\v\theta-\theta^*\v-C_2\,\varepsilon\,.
$$
Taking $\varepsilon$ small compared to $c_0$ and $C_1$ we obtain (ii). Then
(i) follows easily from (\ref{eqIV.4.50}) if $\varepsilon\ll c_0$. Now the
first inequality in (iii) being obvious, let us prove the second one. We write
\begin{equation}\label{eqIV.4.51}\left\{
\begin{array}{l}
\v x-x(\theta,\alpha)\v^2=(1)+(2)\\
(1)=\v x-x(\theta^*,\alpha)\v^2+\v x(\theta^*,\alpha)-x(\theta,\alpha)\v^2\\
(2)=2(x-x(\theta^*,\alpha))(x(\theta^*,\alpha)-x(\theta,\alpha))\,.
\end{array}\right.
\end{equation}
If we use (\ref{eqIV.4.50}), (i) and (ii) we obtain,
$$
(2)=-4(\theta-\theta^*)(x-x(\theta^*,\alpha))\cdot \alpha_\xi+\CO (\varepsilon
(1))\,,
$$
so by (\ref{eqIV.4.49}),
$$
(2)=-4(\theta-\theta^*)\,x\cdot \alpha_\xi+\CO(\varepsilon\,(1))\,.
$$
Now since $(\theta,x)$ belongs to $\CO_\delta$ (see (\ref{eqIV.4.40})) we have
$$
x\cdot \alpha_\xi\leq \frac{c_0}{10}\,\bra x\ket\,\v\alpha_\xi\v\leq
\frac{c_0}{10}\,\bra x(\theta,\alpha)\ket\,\v\alpha_\xi\v+\frac{c_0}{10}\,\v
x-x(\theta,\alpha)\v\,\v\alpha_\xi\v\,.
$$
It follows from (\ref{eqIV.4.48}) that
$$
x\cdot \alpha_\xi\leq \frac 3{10}\,x(\theta,\alpha)\cdot
\alpha_\xi+\frac{c_0}{10}\,\v x-x(\theta,\alpha)\v\,\v\alpha_\xi\v\,,
$$
and we deduce from (\ref{eqIV.4.49}) that
\begin{equation*}
\begin{aligned}
x\cdot \alpha_\xi&\leq  \frac 3{10}\,
(x(\theta,\alpha)-x(\theta^*,\alpha))\cdot
\alpha_\xi+\frac{c_0}{10}\,\v x-x(\theta,\alpha)\v\,\v\alpha_\xi\v\,,\\
x\cdot \alpha_\xi&\leq \frac 3{10}\, \v x(\theta,\alpha)-x(\theta^*,\alpha)\v\,
\v\alpha_\xi\v+\frac{c_0}{10}\,\v x-x(\theta,\alpha)\v\,\v\alpha_\xi\v\,,\\
x\cdot \alpha_\xi&\leq \Big ( \frac 3{10}+\frac {c_0}{10}\Big ) \v
x(\theta,\alpha)-x(\theta^*,\alpha)\v\,\v\alpha_\xi\v+\frac{c_0}{10}\,\v
x-x(\theta^*,\alpha)\v\,\v\alpha_\xi\v\,,
\end{aligned}
\end{equation*}
$$
(2)\geq -4\Big ( \frac 3{10}+\frac{c_0}{10}\Big )\v
x(\theta,\alpha)-x(\theta^*,\alpha)\v\,
\v\theta-\theta^*\v\,\v\alpha_\xi\v-\frac{2c_0}5\,\v
x-x(\theta^*,\alpha)\v\,\v\theta-\theta^*\v\, \v\alpha_\xi\v+\CO
(\varepsilon\,(1))\,.
$$
Using the first inequality in (i) we obtain
$$
(2)\geq -\Big (\frac 4{5}+\frac{4c_0}{15}\Big ) \v
x(\theta,\alpha)-x(\theta^*,\alpha)\v^2-\frac{4c_0}{15}\, \v
x-x(\theta^*,\alpha)\v\,\v x(\theta,\alpha)-x(\theta^*,\alpha)\v+\CO
(\varepsilon\,(1))\,.
$$
Finally
$$
(2)\geq -\Big (\frac 4{5}+\frac {12\,c_0}{15}+C\,\varepsilon\Big )\,(1)\,.
$$
If $c_0$ and $\varepsilon$ are small enough we find $(1)+(2)\geq \frac
16\,(1)$ so (iii) is proved using (\ref{eqIV.4.51}). Finally (iv) follows from
(i) and (iii) taking $\delta$ small enough.
\cqfd

Then we split the case 3 in two subcases.

{\bf Case 3.1 } $(\theta,x)\in\CO_\delta,\,x(\theta,\alpha)\cdot
\alpha_\xi>\frac{c_0}3\,\bra x(\theta,\alpha)\ket\,\v\alpha_\xi\v$ and $\v
x-x(\theta^*,\alpha)\v\leq \langle x(\theta^*,\alpha)\rangle$. 

It follows then that
\begin{equation}\label{eqIV.4.52}
s\,x+(1-s)\,x(\theta^*,\alpha)\in\Omega_\delta\enskip \textrm{for all}\enskip
s\in[0,1]\,.
\end{equation}
Indeed, using Lemma \ref{lIV.4.16} with $a=x(\theta^*,\alpha)$, $b=x$ we
obtain $s\,\v x\v\leq \sqrt 2\,\langle s\,x+(1-s)\,x(\theta^*,\alpha)\rangle$
so if
$(\theta,x)\in\CO_\delta$ we get
$$
(s\,x+(1-s)\,x(\theta^*,\alpha))\cdot \alpha_\xi=s\,x\cdot \alpha_\xi\leq
\frac{c_0}{10}\, s\,\bra x\ket\,\v\alpha_\xi\v\leq c_0\,\bra
s\,x+(1-s)\,x(\theta^*,\alpha)\ket\,\v\alpha_\xi\v\,.
$$
Moreover by Lemma \ref{lIV.4.17},
\begin{equation*}
\begin{aligned}
\v s\,x+(1-s)\,x(\theta^*,\alpha)-x(\theta,\alpha)\v&\leq s\,\v
x-x(\theta,\alpha)\v+(1-s)\,\v x(\theta^*,\alpha)-x(\theta,\alpha)\v\\
&\leq s\,\v x-x(\theta,\alpha)\v+(1-s)\,5\,\v x-x(\theta,\alpha)\v\\
&\leq 5\,\v x-x((\theta,\alpha)\v\leq \delta\,\bra\theta\ket
\end{aligned}
\end{equation*}
since in $\CO_\delta$, $\v x-x(\theta,\alpha)\v\leq \frac \delta {10}
\,\bra\theta\ket$.

Therefore we can define $\varphi $ on this part of $\CO_\delta$ by the
following formula.
\begin{equation}\label{eqIV.4.53}
\begin{split}
\varphi (\theta,x,\alpha)&=\int^1_0 (x-x(\theta^*,\alpha))\cdot
\Phi(\theta,s\,x+(1-s)\,x(\theta^*,\alpha),\alpha)\,ds\\
&\quad -\int^\theta_{\theta^*}
p(x(\theta^*,\alpha),\,\Phi(s,x(\theta^*,\alpha),\alpha)\,ds+\theta^*\,p(\alpha)+\frac
1{2i}\,\v\alpha_\xi\v^2\,.
\end{split}
\end{equation}
Our goal now is to show that $\varphi $ satisfies the claims (i) to (vi) of
Proposition \ref{pIV.4.14}.

The point $\theta=0$ does not belong, by (\ref{eqIV.4.48}), to this part of
$\CO_\delta$. Thus the claim (i) is empty. Let us check (ii). We have
\begin{equation*}
\begin{split}
\frac{\partial \varphi }{\partial x_k}\,(\theta,x,\alpha)&=\int^1_0
\Phi_k(\theta,s\,x+(1-s)\,x(\theta^*,\alpha),\alpha)\,ds\\
&\quad +\som^n_{\ell =1}
\int^1_0 s(x_\ell -x_\ell (\theta^*,\alpha))
\frac{\partial \Phi_\ell }{\partial
x_k}\,(\theta,s\,x+(1-s)\,x(\theta^*,\alpha),\alpha)\,ds\,.
\end{split}
\end{equation*}
Using Theorem \ref{tIV.4.10}, (vii) we see easily that
\begin{equation}\label{eqIV.4.54}
\begin{split}
\frac{\partial \varphi }{\partial x_k}\,(\theta,x,\alpha)=\int^1_0
\Phi_k(\theta,s\,x+(1-s)\,x(\theta^*,\alpha),\alpha)\,ds\\
+\int^1_0 s\,\frac
d{ds}\,[\Phi_k(\theta,s\,x+(1-s)\,x(\theta^*,\alpha),\alpha)]\,ds+R
\end{split}
\end{equation}
with
$$
\v R\v\leq C_N\,\v x-x(\theta^*,\alpha)\v \int^1_0 \frac{\v
s\,x+(1-s)\,x(\theta^*,\alpha)-x(\theta,\alpha)\v^N}{\bra\theta\ket^{2N}}\,ds\,.
$$
It follows from Lemma \ref{lIV.4.17}, (iii) that
\begin{equation}\label{eqIV.4.55}
\v R\v\leq C_N\, \frac{\v x-x(\theta,\alpha)\v^N}{\bra\theta\ket^N}\,,\quad
N\geq 0\,.
\end{equation}
Integrating by parts in the second integral of the right hand side of
(\ref{eqIV.4.54}) we obtain the claim (ii). Let us prove (iii). We have
$$
\frac{\partial \varphi }{\partial \theta}\,(\theta,x,\alpha)=\int^1_0
(x-x(\theta^*,\alpha))\cdot \frac{\partial \Phi}{\partial
\theta}\,(\theta,X_s,\alpha)\,ds-p(x(\theta^*,\alpha),\Phi(\theta,x(\theta^*,\alpha),\alpha)
$$
where $X_s=s\,x+(1-s)\,x(\theta^*,\alpha)$.

Using Corollary \ref{cIV.4.13} we obtain
\begin{equation*}
\begin{split}
\frac{\partial \varphi }{\partial \theta}\,(\theta,x,\alpha)=-\int^1_0
(x-x(\theta^*,\alpha))\cdot \frac{\partial p}{\partial
x}\,(X_s,\Phi(\theta,X_s,\alpha))\,ds-\som^n_{k,\ell =1} \int^1_0
(x_k-x_k(\theta^*,\alpha))\\
\frac{\partial \Phi_k}{\partial x_\ell
}\,(\theta,X_s,\alpha)\,\frac{\partial p}{\partial \xi_\ell
}\,(X_s,\Phi(\theta,X_s,\alpha))\,ds-p(x(\theta^*,\alpha),\Phi(\theta,x(\theta^*,
\alpha),\alpha))+R
\end{split}
\end{equation*}
where $R$ satisfies (\ref{eqIV.4.55}).
We use again Theorem \ref{tIV.4.10}, (vii), and we obtain
$$
\frac{\partial \varphi }{\partial \theta}\,(\theta,x,\alpha)=-\int^1_0 \frac
d{ds}\,
[p(X_s,\Phi(\theta,X_s,\alpha))]\,ds-p(x(\theta^*,\alpha),\Phi(\theta,x(\theta^*,\alpha),\alpha)+R'
$$
where $R'$ satisfies also (\ref{eqIV.4.55}). This implies immediately (iii).

The points (iv), (v) follow easily from Theorem \ref{tIV.4.10}.  Let us check
(vi). According to (\ref{eqIV.4.53}) we can write
\begin{equation}\label{eqIV.4.56}\left\{
\begin{array}{l}
\varphi (\theta,x,\alpha)=A+B\enskip \textrm{with,}\\
A=\int^1_0 (x-x(\theta^*,\alpha))\cdot
\Phi(\theta,s\,x+(1-s)\,x(\theta^*,\alpha),\alpha)\,ds\,.
\end{array}\right.
\end{equation}
Using Theorem \ref{tIV.4.10}, (iv) and Lemma \ref{lIV.4.17} (iii) we see that
\begin{equation}\label{eqIV.4.57}\left\{
\begin{array}{l}
\Im A=\frac 12\,\frac 1{1+4\theta^2}\,
[(x-x(\theta^*,\alpha))^2+2(x-x(\theta^*,\alpha))
(x(\theta^*,\alpha)-x(\theta,\alpha))]+R\\
\v R\v\leq C\,\sqrt\delta\, \frac{\v
x-x(\theta,\alpha)\v^2}{\bra\theta\ket^2}\,.
\end{array}\right.
\end{equation}
To check the term $B$ let us set $x(t)=x(t,\alpha)$ and
$$
F(t)=-\int^\theta_t\, p(x(t),\Phi(s,x(t),\alpha))\,ds+t\,p(\alpha)\,.
$$
Then
\begin{equation*}
\begin{split}
F'(t)&=\underbrace{p(x(t),\Phi(t,x(t),\alpha))}_{(1)}-\underbrace{\int^\theta_t
\som^n_{k=1}
\frac{\partial p}{\partial
x_j}(x(t),\Phi(s,x(t),\alpha))\,\dot x_k(t)\,dt}_{(2)}\\
&\quad \quad -\underbrace{\int^\theta_t \som^n_{k,\ell =1} \frac{\partial
p}{\partial
\xi_\ell }\,(x(t),\Phi(s,x(t),\alpha))\,\frac{\partial
\Phi_\ell }{\partial x_k}\,(s,x(t),\alpha)\,\dot x_k(t)\,ds}_{(3)}+p(\alpha)\,.
\end{split}
\end{equation*}
By Theorem \ref{tIV.4.10} (vi) we have,
$$
(1)=p(x(t),\,\xi(t))=p(\alpha)\,.
$$
By the point (vii) of the same theorem we have,
$$
(3)=\som^n_{k,\ell =1} \int^\theta_t \frac{\partial p}{\partial \xi_\ell
}\,(x(t),\Phi(s,x(t),\alpha))\,\frac{\partial \Phi_k}{\partial x_\ell
}\,(s,x(t),\alpha)\,\dot x_k(t)\,ds+R_0
$$
with
$$
\v R_0\v\leq C_N \int^\theta_t \frac{\v x(t)-x(s)\v^N}{\bra s\ket^{2N}}\,.
$$
Since $\v x(t)-x(s)\v\leq \int^s_t \v\dot x(\sigma)\v\,d\sigma\leq C(s-t)$ we
obtain
\begin{equation}\label{eqIV.4.58}
\v R_0\v\leq C'_N \int^\theta_t \frac{(s-t)^N}{\bra s\ket^{2N}}\,ds\,,\quad
\theta^*\leq t\leq \theta\,.
\end{equation}
Using Corollary \ref{cIV.4.13} we obtain
$$
(3)=-\som^n_{k,\ell =1} \int^\theta_t \frac{\partial p}{\partial
x_k}\,(x(t),\Phi(s,x(t),\alpha),\alpha)\,\dot x(t)\,ds-\som^n_{k=1}
\int^\theta_t \frac{\partial \Phi_k}{\partial s}\,(s,x(t),\alpha)\,\dot
x_k(t)\,ds+R_1$$
where $R_1$ satisfies (\ref{eqIV.4.58}).

It follows that
$$
(3)=-(2)-\som^n_{k=1} \dot
x_k(t)(\Phi_k(\theta,x(t),\alpha)-\Phi_k(t,x(t),\alpha))+R_1\,.
$$
Now
$$
\som^n_{k=1} \dot x_k(t)\,\Phi_k(t,x(t),\alpha)=\som^n_{k=1}
\xi_k(t)\,\frac{\partial p}{\partial \xi_k}\,(x(t),\xi(t))=2p(\alpha)\,.
$$
Therefore we obtain,
$$
F'(t)=(1)-(2)-(3)+p(\alpha)=\som^n_{k=1} \dot
x_k(t)\,\Phi_k(\theta,x(t),\alpha)+R_1\,.
$$
Now by Theorem \ref{tIV.4.10} (iv),
$$
\Im
\Phi_k(\theta,x(t),\alpha)=\frac{x_k(t)-x_k(\theta)}{1+4\theta^2}+\CO\,(\sqrt
\delta)\,\frac{\v x_k(t)-x_k(\theta)\v}{\bra\theta\ket^2}\,.
$$
Since $\dot x(t)$ is uniformly bounded we deduce that
\begin{equation*}
\begin{aligned}
\Im F'(t)&=\frac 12\,\frac d{dt}\,\frac{\v
x(t)-x(\theta)\v^2}{1+4\theta^2}+G(t)\quad \textrm{with}\\
\v G(t)\v&\leq C\,\sqrt\delta\,\frac{\theta-t}{\bra\theta\ket^2}+C_N
\int^\theta_t \,\frac{(s-t)^N}{\bra s\ket^{2N}}\,ds\,.
\end{aligned}
\end{equation*}
Integrating between $\theta^*$ and $\theta$ we obtain
\begin{equation}\label{eqIV.4.59}
\Big\v\Im F(\theta^*)-\frac 12\,\frac{\v
x(\theta^*)-x(\theta)\v^2}{1+4\theta^2}\Big\v\leq
C'\,\sqrt\delta\,\frac{(\theta-\theta^*)^2}{\bra\theta\ket^2}+C_N
\int^\theta_{\theta^*} \int^\theta_t\,\frac{(s-t)^N}{\bra
s\ket^{2N}}\,ds\,dt\,.
\end{equation}
Let us call $I$ (resp. $II$) the first (resp. the second) term in the right
hand side of (\ref{eqIV.4.59}). By Lemma \ref{lIV.4.17} we have
\begin{equation}\label{eqIV.4.60}
\v I\v\leq C\,\sqrt\delta\,\frac{\v
x-x(\theta,\alpha)\v^2}{\bra\theta\ket^2}\,.
\end{equation}
Now
$$
\v II\v\leq C'_N \int^\theta_{\theta^*} \Big (\int^s_{\theta^*}
(s-t)^N\,dt\Big )\,\frac{ds}{(1+s)^{2N}}\leq C''_N
\int^\theta_{\theta^*}\,\frac{(s-\theta^*)^{N+1}}{(1+s)^{2N}}\,ds\,.
$$
Now it follows from Lemma \ref{lIV.4.17} and (\ref{eqIV.4.40}) that
$\theta-\theta^*\leq 2\delta\,\bra\theta\ket\leq 2\delta\,(1+\theta)$ which
means that $(1-2\delta)\,\theta\leq \theta^*+2\delta$. Since $2\delta\leq
\frac 12$ we have $\theta\leq 2\theta^*+1$. It is then easy to see that the
function $s\mapsto \frac{(s-\theta^*)^{N+1}}{(1+s)^{2N}}$
 is increasing on $(\theta^*,\theta)$. Therefore
$$
\v II\v\leq C_N\, \frac{(\theta-\theta^*)^{N+2}}{\bra\theta\ket^{2N}}\leq
C_N(\theta-\theta^*)^2\, \frac{(2\delta)^N}{\bra\theta\ket^N}\,.
$$
Taking $N=2$ and using Lemma \ref{lIV.4.17} (i) and (iii) we obtain
\begin{equation}\label{eqIV.4.61}
\v II\v\leq C\,\delta^2\,\frac{\v x-x(\theta,\alpha)\v^2}{\bra\theta\ket^2}\,.
\end{equation}
It follows from (\ref{eqIV.4.59}) to (\ref{eqIV.4.61}) and from
(\ref{eqIV.4.53}), (\ref{eqIV.4.56}), (\ref{eqIV.4.57}) that
$$
\Im \varphi (\theta,x,\alpha)=\frac 12\, \frac{\v
x-x(\theta,\alpha)\v^2}{1+4\theta^2}-\frac 12\, \v\alpha_\xi\v^2+\CO \Big
(\sqrt\delta\,\frac{\v x-x(\theta,\alpha)\v^2}{\bra\theta\ket^2}\Big)
$$
which is precisely the claim of point (vii) of Proposition \ref{pIV.4.14}.

{\bf Case 3.2 } $(\theta,x)\in\CO_\delta$, $x(\theta,\alpha)\cdot
\alpha_\xi>\frac{c_0}3\,\bra x(\theta,\alpha)\ket\v\alpha_\xi\v$ and $\v
x-x(\theta^*,\alpha)\v\geq \frac 12\v x(\theta^*,\alpha)\v$.

According to Lemma \ref{lIV.4.17} (iii) we have
\begin{equation}\label{eqIV.4.62}\left\{
\begin{array}{l}
\v x(\theta^*,\alpha)\v\leq 10\,\v x-x(\theta,\alpha)\v\leq 2
\delta\,\bra\theta\ket\\
\v x\v\leq \v x-x(\theta^*,\alpha)\v+\v x(\theta^*,\alpha)\v\leq 15\,\v
x-x(\theta,\alpha)\v\leq \frac 32\,\delta\, \bra\theta\ket\,.
\end{array}\right.
\end{equation}
On the other hand if $y$ belongs to the union of the two segments $[0,x]$ and
$[0,x(\theta^*,\alpha)]$ then $(y,\theta)$ belongs to $\Omega_\delta$, the set
(defined in (\ref{eqIV.1.4})) on which $\Phi$ is defined. Indeed, by
(\ref{eqIV.4.40}), if $s\in (0,1)$ then $s\,x\cdot \alpha_\xi\leq s\cdot
\frac{c_0}{10}\,\bra x\ket\v\alpha_\xi\v\leq c_0\,\bra
s\,x\ket\v\alpha_\xi\v$. Moreover $\v s\,x-x(\theta,\alpha)\v\leq \v
s\,x-x(\theta^*,\alpha)\v+\v x(\theta^*,\alpha)-x(\theta,\alpha)\v\leq s\v
x-x(\theta^*,\alpha)\v+(1-s)\v x(\theta^*,\alpha)\v+\v
x(\theta^*,\alpha)-x(\theta,\alpha)\v$. Since we are in case 3.2 we have by
Lemma \ref{lIV.4.17}, $\v s\,x-x(\theta,\alpha)\v
\leq 10\,\v
x-x(\theta,\alpha)\v\leq \delta\,\bra\theta\ket$. On the other hand, if $s\in
(0,1)$ we have, by (\ref{eqIV.4.49}), $s\,x(\theta^*,\alpha)\cdot
\alpha_\xi=0$. Moreover $\v s\,x(\theta^*,\alpha)-x(\theta,\alpha)\v\leq \v
x(\theta^*,\alpha)-x(\theta,\alpha)\v+(1-s)\v x(\theta^*,\alpha)\v\leq \v
x(\theta^*,\alpha)-x(\theta,\alpha)\v+2(1-s)\v x-x(\theta^*,\alpha)\v\leq
10\,\v x-x(\theta,\alpha)\v\leq \delta\,\bra\theta\ket$, by Lemma
\ref{lIV.4.17}, (iii).

Therefore in the present case we can set
\begin{equation}\label{eqIV.4.63}
\begin{split}
\varphi (\theta,x,\alpha)&=\int^1_0 x\cdot
\Phi(\theta,s\,x,\alpha)\,ds-\int^1_0 x(\theta^*,\alpha)\cdot
\Phi(\theta,s\,x(\theta^*,\alpha),\alpha)\,ds\\
&\quad -\int^\theta_{\theta^*}
p(x(\theta^*,\alpha),\Phi(s,x(\theta^*,\alpha),\alpha))\,ds+\theta^*\,
p(\alpha) +\frac 1{2i}\,\alpha^2_\xi\,.
\end{split}
\end{equation} 
Our goal is to show that $\varphi $ satisfies all the requirements of
Proposition \ref{pIV.4.14}.

The point (i) is empty since $\theta=0$ does not belong to this part of
$\CO_\delta$. Let us check (ii). We have
$$
\frac{\partial \varphi }{\partial x_j}\,(\theta,x,\alpha)=\int^1_0
\Phi_j(\theta,x,\alpha)\,ds+\int^1_0 \som^n_{k=1} s\,x_k\,\frac{\partial
\Phi_k}{\partial x_j}\,(\theta,s\,x,\alpha)\,ds\,.
$$
By Theorem \ref{tIV.4.10} and (\ref{eqIV.4.62}) we have
\begin{equation}\label{eqIV.4.64}
\Big\v\frac{\partial \Phi_k}{\partial
x_j}\,(\theta,s\,x,\alpha)-\frac{\partial \Phi_j}{\partial
x_k}\,(\theta,s\,x,\alpha)\Big\v\leq C_N\,\frac{\v
s\,x-x(\theta,\alpha)\v^N}{\bra\theta\ket^{2N}}\leq C'_N\,\frac{\v
x-x(\theta,\alpha)\v^N}{\bra\theta\ket^{2N}}\,.
\end{equation}
It follows that
$$
\frac{\partial \varphi }{\partial x_j}\,(\theta,x,\alpha)=\int^1_0
\Phi_j(\theta,s\,x,\alpha)\,ds+\int^1_0 s\,\frac
d{ds}\,(\Phi_j(\theta,s\,x,\alpha))\,ds+\CO \Big (\frac{\v
x-x(\theta,\alpha)\v^{N+1}}{\bra\theta\ket^{2N}}\Big)\,.
$$
Integrating by parts and using the bound $\v x-x(\theta,\alpha)\v\leq
\delta\,\bra\theta\ket$ we obtain
$$
\Big\v\frac{\partial \varphi }{\partial
x_j}\,(\theta,x,\alpha)-\Phi_j(\theta,x,\alpha)\Big\v\leq C_N\,\frac{\v
x-x(\theta,\alpha)\v^N}{\bra\theta\ket^N}\,,\quad \forall \,N\in\N\,.
$$
Thus (ii) is proved. Let us prove (iii). We have
\begin{equation}\label{eqIV.4.65}
\begin{split}
\frac{\partial \varphi }{\partial
\theta}\,(\theta,x,\alpha)=\underbrace{\int^1_0 x\,\frac{\partial
\Phi}{\partial
\theta}\,(\theta,s\,x,\alpha)\,ds}_{(1)}-\underbrace{\int^1_0
x(\theta^*,\alpha)\,\frac{\partial \Phi}{\partial \theta}\,
(\theta,s\,x(\theta^*,\alpha),\alpha)\,ds}_{(2)}\\
-\underbrace{p(x(\theta^*,\alpha),\Phi\,(\theta,x(\theta^*,\alpha),\alpha)}_{(3)}\,.
\end{split}
\end{equation}
Using Corollary \ref{cIV.4.13} we can write
\begin{equation*}
\begin{split}
(1)=&-\som^n_{k=1} \int^1_0 x_k\,\frac{\partial p}{\partial
x_k}\,(s\,x,\Phi(\theta,s\,x,\alpha))\,ds\\
&+\som^n_{k,\ell =1} \int^1_0
x_k\,(\theta,s\,x,\alpha)\cdot \frac{\partial p}{\partial x_\ell
}\,(s\,x,\Phi(\theta,s\,x,\alpha))\,\frac{\partial \Phi_k}{\partial
x_k}\,(\theta,s\,x,\alpha)\,ds\,.
\end{split}
\end{equation*}
Using again (\ref{eqIV.4.64}) we obtain,
$$
(1)=-\int^1_0 \frac{d}{ds} \big
(p(s\,x,\Phi(\theta,s\,x,\alpha))\big)\,ds+\CO\Big (\frac{\v
x-x(\theta,\alpha)\v^{N+1}}{\bra\theta\ket^{2N}}\Big)\,.
$$
Finally, since $\v x-x(\theta,\alpha)\v\leq \delta\,\bra\theta\ket$, we have
\begin{equation}\label{eqIV.4.66}
(1)=p(0,\Phi(\theta,0,\alpha))-p(x,\Phi(\theta,x,\alpha))+\CO\Big (\frac{\v
x-x(\theta,\alpha)\v^N}{\bra\theta\ket^N}\Big)\,.
\end{equation}
By exactly the same computation (using (\ref{eqIV.4.62})) we obtain
\begin{equation}\label{eqIV.4.67}
(2)=p(0,\Phi(\theta,0,\alpha))-p(x(\theta^*,\alpha),\Phi(\theta,x(\theta^*,\alpha),\alpha)
+\CO \Big (\frac{\v
x-x(\theta,\alpha)\v^N}{\bra\theta\ket^N}\Big)\,.
\end{equation}
So using (\ref{eqIV.4.65}) to (\ref{eqIV.4.67}) we derive the point (iii). The
last non trivial point to be proved is the point (vi).

Using the expression of $\Im \Phi$ given by Theorem \ref{tIV.4.10},
(\ref{eqIV.4.62}) and (\ref{eqIV.4.59}) to (\ref{eqIV.4.61}) which are valid
also in the Case~3.2,  we can write
\begin{equation*}
\begin{split}
\Im \varphi (\theta,x,\alpha)=\frac 12\, \frac 1{1+4\theta^2} \big (\v
x-x(\theta,\alpha)\v^2-\v x(\theta,\alpha)\v^2-\v
x(\theta^*,\alpha)-x(\theta,\alpha)\v^2\\
+\v x(\theta,\alpha)\v^2+\v
x(\theta^*,\alpha)-x(\theta,\alpha)\v^2\big)+\CO\Big ( \sqrt\delta\,\frac{\v
x-x(\theta,\alpha)\v^2}{\bra\theta\ket^2}\Big)
\end{split}
\end{equation*}
which is exactly what is needed.

To finish the proof of Proposition \ref{pIV.4.14} we must show that the phases
$\varphi $ which have been constructed by the formulas (\ref{eqIV.4.41}),
(\ref{eqIV.4.43}), (\ref{eqIV.4.53}), (\ref{eqIV.4.63}) in different regions
can be matched in only one phase. We begin by a Lemma.
\begin{lemma}\sl \label{lIV.4.18}
Let $I=\big[\frac{c_0}{10}\,, \frac{c_0}{2}\big]$ and let us consider the
function on $[0,+\infty [\times I$,
$$
g(s,c)=x(s,\alpha)\cdot \alpha_\xi-c\,\bra x(s,\alpha)\ket\,\v\alpha_\xi\v\,.
$$
\begin{itemize}
\item[(i)]  For all $c$ in $I$ the function $[0,+\infty [\times \R$, $s\mapsto
g(s,c)$ is strictly increasing.

\item[(ii)] For all $c$ in $I$ there exists a unique $\theta(c)>0$ such that
$$
g(\theta(c),c)=0\,.
$$ 
\item[(iii)]  The function $I\rightarrow [0,+\infty [$, $c\mapsto \theta(c)$
is strictly increasing.

Moreover we have the following estimates
\item[(iv)] $\ds{\frac 32\,\v\theta(c)-\theta^*\v\,\v\alpha_\xi\v\leq \v
x(\theta(c),\alpha)-x(\theta^*,\alpha)\v\leq
3\,\v\theta(c)-\theta^*\v\,\v\alpha_\xi\v}$.
\item[(v)] $\ds{\v\theta(c)-\theta^*\v\geq \frac{c_0}{120}}$.
\item[(vi)] For all $x$ in $\CO_\delta$ and $c$ in
$\big[\frac{c_0}3\,\frac{c_0}2\big]$,
$$
\v x-x(\theta(c),\alpha)\v\leq \v x-x(\theta^*,\alpha)\v+\v
x(\theta^*,\alpha)-x(\theta(c),\alpha)\v\leq 4\,\v x-x(\theta(c),\alpha)\v\,.
$$
\item[(vii)] If $c\in\big[\frac{c_0}3\,\frac{c_0}2\big]$ we have
$$
\frac c{10}\, \bra x(\theta^*,\alpha)\ket\leq \v\theta(c)-\theta^*\v\leq
4c\,\bra x(\theta^*,\alpha)\ket\,.
$$
\end{itemize}
\end{lemma}

\noindent {\bf Proof } 

(i) We have $\frac{\partial g}{\partial s}\,(s,c)=\dot x(s,\alpha)\cdot
\alpha_\xi-c\,\frac{x(s,\alpha)\cdot \dot x(s,\alpha)}{\bra
x(s,\alpha)\ket}\,\v\alpha_\xi\v=2\,\v\alpha_\xi\v^2+\CO (\varepsilon+c_0)$. 
Thus $\frac{\partial g}{\partial s}\,(s,c)\geq \frac 1{10}$ if $\varepsilon$
and $c_0$ are small enough.

(ii) It follows from above that $g(s,c)\geq \frac 1{10}\,s+g(0,c)$ so
$g(s,c)\rightarrow +\infty $ if $s\rightarrow +\infty $. Moreover
$g(0,c)=\alpha_x\cdot \alpha_\xi-c\,\bra\alpha_x\ket\,\v\alpha_\xi\v\leq
\alpha_x\cdot \alpha_\xi\leq -c_0\,\bra \alpha_x\ket\,\v\alpha_\xi\v<0$.
Therefore there exists a unique $\theta(c)$ such that $g(\theta(c),c)=0$ and
$c\mapsto \theta(c)$ is $C^\infty $. Differentiating this equality with
respect to $c$ we obtain
$$
\theta'(c)\,\frac{\partial g}{\partial s}\,(\theta(c),c)+\frac{\partial
g}{\partial c}\,(\theta(c),c)=0\,.
$$
By the above computation of $\frac{\partial g}{\partial s}$ we can write
$$
\theta'(c)=\frac{\bra
x(\theta(c),\alpha)\ket\,\v\alpha_\xi\v}{2\,\v\alpha_\xi\v^2+\CO
(\varepsilon+c_0)}
$$
which proves (iii). Now we have
$$
x(\theta(c),\alpha)-x(\theta^*,\alpha)=\int^{\theta(c)}_{\theta^*}
x(s,\alpha)\,ds=2\alpha_\xi(\theta(c)-\theta^*)+\CO
(\varepsilon\v\theta(c)-\theta^*\v)
$$
from which (iv) follows easily. Let us prove (v). By definition of $\theta(c)$
and $\theta^*$ we can write
\begin{equation}\label{eqIV.4.68}
(x(\theta(c),\alpha)-x(\theta^*,\alpha))\cdot \alpha_\xi=c\,\bra
x(\theta(c),\alpha)\ket\,\v\alpha_\xi\v\geq \frac{c_0}{20}
\end{equation}
so by (iv),
$$
\frac{c_0}{20}\leq 3\,\v\theta(c)-\theta^*\v\cdot \v\alpha_\xi\v\leq
6\,\v\theta(c)-\theta^*\v\,.
$$
The first inequality in (vi) beeing trivial let us prove the second one. We
write
$$
\v x-x(\theta,\alpha)\v^2=\underbrace{\v x-x(\theta^*,\alpha)\v^2+\v
x(\theta^*,\alpha)-x(\theta(c),\alpha)\v^2}_{(1)}+\underbrace{2(x-x(\theta^*,\alpha))\cdot
(x(\theta^*,\alpha)-x(\theta(c),\alpha))}_{(2)}\,.
$$
We have
$$
(2)=2(x-x(\theta^*,\alpha)) [2(\theta^*-\theta(c))\,\alpha_\xi+\CO
(\varepsilon\,\v\theta^*-\theta(c)\v]\,.
$$
It follows from (iv) that
$$
(2)=-4(\theta(c)-\theta^*)\,x\cdot \alpha_\xi+\CO (\varepsilon(1))\,,\quad
\textrm{where}\enskip \theta(c)-\theta^*\geq 0\,.
$$
Now in $\CO_\delta$ we have $x\cdot \alpha_\xi\leq \frac{c_0}{10}\,\bra
x\ket\,\v \alpha_\xi\v$. It follows that 
$$
x\cdot \alpha_\xi\leq
\frac{c_0}{10}\,\bra
x(\theta(c),\alpha)\ket\,\v\alpha_\xi\v+\frac{c_0}{10}\,\v
x-x(\theta(c),\alpha)\v\,\v\alpha_\xi\v\,.
$$
Using (\ref{eqIV.4.68}) we obtain
$$
x\cdot \alpha_\xi\leq
\frac{c_0}{10c}\,(
x(\theta(c),\alpha)-x(\theta^*,\alpha))\cdot \alpha_\xi+\frac{c_0}{10}\,\v
x-x(\theta(c),\alpha)\v\,\v\alpha_\xi\v\,,
$$
so we will have
\begin{equation*}
\begin{aligned}
(2)&\geq
-\frac{2c_0}{5c}\,\v
x(\theta(c),\alpha)-x(\theta^*,\alpha)\v\,\v\alpha_\xi\v\,\v\theta(c)-\theta^*\v
-\frac{2c_0}{5}\,\v
x-x(\theta(c),\alpha)\v\,\v\alpha_\xi\v\,\v\theta(c)-\theta^*\v-\CO
(\varepsilon(1))\\
(2)&\geq
-\frac{4c_0}{15c}\,\v
x(\theta(c),\alpha)-x(\theta^*,\alpha)\v^2
-\frac{8c_0}{30}\,\v
x-x(\theta(c),\alpha)\v\,\v x(\theta(c),\alpha)-x(\theta^*,\alpha)\v-\CO
(\varepsilon(1))\\
(2)&\geq \Big (-\frac{4c_0}{15c}-\frac{4c_0}{30}\Big)\,\v
x(\theta(c),\alpha)-x(\theta^*,\alpha)\v^2-\frac{2}{15}\,c_0\,\v
x-x(\theta(c),\alpha)\v^2-\CO (\varepsilon(1))\,.
\end{aligned}
\end{equation*}
If $c\geq \frac{c_0}3$ then $\frac{4c_0}{15c}\leq \frac 45$ so we obtain
$$
(1)+(2)\geq \frac1{10} \big (\v x(\theta(c),\alpha)-x(\theta^*,\alpha)\v^2+\v
x-x(\theta,\alpha)\v^2\big )
$$
if $c_0$ is small enough. This implies (vi). Let us prove (vii). We have
$$
x(\theta(c),\alpha)=x(\theta^*,\alpha)+2(\theta(c)-\theta^*)\cdot
\alpha_\xi+\CO (\varepsilon)\,.
$$
It follows that
$$
c\,\bra x(\theta(c),\alpha)\ket\,\v\alpha_\xi\v=x(\theta(c),\alpha)\cdot
\alpha_\xi=2(\theta(c)-\theta^*)\v\alpha_\xi\v^2+\CO(\varepsilon)
$$
so
$$
\frac 12\,\v\theta(c)-\theta^*\v\leq c\,\bra x(\theta(c),\alpha)\ket\leq
5\,\v\theta(c)-\theta^*\v\,.
$$
Moreover we have
\begin{equation*}
\begin{aligned}
\bra x(\theta(c),\alpha)\ket&\leq \bra x(\theta^*,\alpha)\ket+\v
x(\theta^*,\alpha)-x(\theta(c),\alpha)\v\\
&\leq \bra x(\theta^*,\alpha)\ket+3\,\v\theta(c)-\theta^*\v\,\v\alpha_\xi\v\\
&\leq \bra x(\theta^*,\alpha)\ket+6c\,\bra
x(\theta(c),\alpha)\ket\,\v\alpha_\xi\v\\
&\leq \bra x(\theta^*,\alpha)\ket+12c\,\bra x(\theta(c),\alpha)\ket
\end{aligned}
\end{equation*}
which implies that $\bra x(\theta(c),\alpha)\ket\leq 2\,\bra
x(\theta^*,\alpha)\ket$ if $c_0$ is small enough. By the same way $\bra
x(\theta^*,\alpha)\ket\leq 2\,\bra x(\theta(c),\alpha)\ket$. Thus we obtain
(vii). \cqfd

Now let us set
\begin{equation}\label{eqIV.4.69}\left\{
\begin{array}{l}
\widetilde {\CO}_\delta(\theta)=\Big\{x\in\R^n : x\cdot \alpha_\xi\leq
\frac{c_0}{10}\,\bra x\ket\,\v\alpha_\xi\v\,,\enskip \v
x-x(\theta,\alpha)\v\leq \frac{\delta}{40}\,\bra\theta\ket\Big\}\\
\overline \theta=\theta\Big (\frac{c_0}2\Big )\,,\enskip \theta'=\theta\Big
(\frac{c_0}{10}\Big )\,.
\end{array}\right.
\end{equation}
In the beginning of this Section we have constructed the different $\varphi $
assuming $x\cdot \alpha_\xi\leq \frac{c_0}{10}\,\bra x\ket\,\v\alpha_\xi\v$,
$\v x-x(\theta,\alpha)\v\leq \frac\delta{5}\, \bra\theta\ket$.

In the proof of Proposition \ref{pIV.4.14} we have constructed
\begin{equation*}
\begin{aligned}
\varphi _1\enskip &\textrm{when}\enskip \theta\in[0,\overline \theta]\enskip
\textrm{and}\enskip \v x-x(\theta,\alpha)\v\leq \bra
x(\theta,\alpha)\ket\,,\hbox to 1,2cm{}  \textrm{(case 1)}\\
\varphi _2\enskip &\textrm{when}\enskip \theta\in[0,\overline \theta]\enskip
\textrm{and}\enskip \v x-x(\theta,\alpha)\v\geq \frac 12\,\v
x(\theta,\alpha)\v\, ,\hbox to 1,2cm{}  \textrm{(case 2)}\\
\varphi _4\enskip &\textrm{when}\enskip \theta\in[\theta'+\infty [\enskip
\textrm{and}\enskip \v x-x(\theta^*,\alpha)\v\leq \bra
x(\theta^*,\alpha)\ket\,,\hbox to 0,4cm{}  \textrm{(case 3.1)}\\
\varphi _5\enskip &\textrm{when}\enskip \theta\in[\theta'+\infty [\enskip
\textrm{and}\enskip \v x-x(\theta^*,\alpha)\v\geq \frac 12\,\v
x(\theta^*,\alpha)\v\,,\enskip \textrm{(case 3.2)}\,.
\end{aligned}
\end{equation*}
We are going first to match $\varphi _1$ and $\varphi _2$, $\varphi _4$ and
$\varphi _5$. The matched phase will be defined on a smallest set than
$\CO_\delta$ defined in (\ref{eqIV.4.40}) namely for $(\theta,x)$ where
$x\in\widetilde{\CO}_\delta(\theta)$ (see (\ref{eqIV.4.69})). We show first
that the point $(\theta,0)$ belongs to the sets where $\varphi _1$ and
$\varphi _2$ are defined. According to (\ref{eqIV.4.40}) and what we recalled
above it will be the case if $\v x(\theta,\alpha)\v\leq \frac\delta
5\,\bra\theta\ket$. We may assume that the domain where $\varphi _2$ is
defined contains points $(\theta,x)$ where
$x\in\widetilde{\CO}_\delta(\theta)$ otherwise we don't match $\varphi _1$ and
$\varphi _2$ and we take only $\varphi _1$. So let $(\theta,x)$ be such $\v
x-x(\theta,\alpha)\v\leq \frac\delta{40}\,\bra\theta\ket$ and $\v
x(\theta,\alpha)\v\leq 2\,\v x-x(\theta,\alpha)\v$. Then $\v
x(\theta,\alpha)\v\leq \frac\delta{40}\,\bra\theta\ket$ which implies our
claim.

Now it follows from (\ref{eqIV.4.41}) and (\ref{eqIV.4.43}) that
\begin{equation}\label{eqIV.4.70}
\varphi _1(\theta,0,\alpha)=\varphi _2(\theta,0,\alpha)\,.
\end{equation}
By the same way we may assume that the domain where $\varphi _5$ is defined
contains points $(\theta,x)$ where $x\in\widetilde{\CO}_\delta(\theta)$. So
let $x$ be such that $\v x-x(\theta,\alpha)\v\leq
\frac\delta{40}\,\bra\theta\ket$ and $\v x-x(\theta^*,\alpha)\v\geq \frac
12\,\v x(\theta^*,\alpha)\v$. Then we write, using Lemma \ref{lIV.4.18},
\begin{equation*}
\begin{aligned}
\v x(\theta,\alpha)\v&\leq \v x(\theta^*,\alpha)\v+\v
x(\theta^*,\alpha)-x(\theta,\alpha)\v\,,\\
\v x(\theta,\alpha)\v&\leq2 (\v x-x(\theta^*,\alpha)\v+\v
x(\theta^*,\alpha)-x(\theta,\alpha)\v)\leq 8\,\v x-x(\theta,\alpha)\v\leq
\frac \delta 5\,\bra\theta\ket\,.
\end{aligned}
\end{equation*}
So the point $(\theta,0)$ belongs also to the sets where $\varphi _4$ and
$\varphi _5$ are defined and by (\ref{eqIV.4.53}), (\ref{eqIV.4.63}) we have
\begin{equation}\label{eqIV.4.71}
\varphi _4(\theta,0,	\alpha)=\varphi _5(\theta,0,\alpha)\,.
\end{equation}
Let us match $\varphi _1$ and $\varphi _2$. Let $x\in\widetilde{\CO}_\delta$,
$\theta\in[0,\overline \theta]$ be such that
\begin{equation}\label{eqIV.4.72}
\frac 12\,\v x(\theta,\alpha)\v\leq \v x-x(\theta,\alpha)\v\leq \bra
x(\theta,\alpha)\ket\,.
\end{equation}
We are going to show then that
\begin{equation}\label{eqIV.4.73}
\forall\, N\in\N\enskip \exists \,  C_N>0 : \v\varphi
_1(\theta,x,\alpha)-\varphi _2(\theta,x,\alpha)\v\leq C_N\,\frac{\v
x-x(\theta,\alpha)\v^{N+1}}{\bra\theta\ket^N}\,.
\end{equation}
Indeed let $\gamma(\sigma,x)$ be a regular path such that
\begin{equation}\label{eqIV.4.74}
\gamma(0,x)=0\,,\quad \gamma(1,x)=x
\end{equation}
and there exists $K\geq 0$ such that for all $\sigma$ in $[0,1]$,
\begin{equation}\label{eqIV.4.75}
\Big\v\frac{\partial \gamma}{\partial \sigma}\,(\sigma,x)\Big\v\leq K\,\v
x-x(\theta,\alpha)\v
\end{equation}
\begin{equation}\label{eqIV.4.76}
\gamma(\sigma,x)\cdot \alpha_\xi\leq
\frac{c_0}{10}\,\bra\gamma(\sigma,x)\ket\,\v\alpha_\xi\v
\end{equation}
\begin{equation}\label{eqIV.4.77}
\v\gamma(\sigma,x)-x(\theta,\alpha)\v\leq \frac \delta 5\,\bra\theta\ket
\end{equation}
\begin{equation}\label{eqIV.4.78}\left\{
\begin{array}{l}
\textrm{if}\enskip \v x-x(\theta,\alpha)\v\geq \v x(\theta,\alpha)\v\enskip
\textrm{then,}\\
\v x(\theta,\alpha)\v\leq \v\gamma(\sigma,x)-x(\theta,\alpha)\v\leq \v
x-x(\theta,\alpha)\v 
\end{array}\right.
\end{equation}
\begin{equation}\label{eqIV.4.79}\left\{
\begin{array}{l}
\textrm{if}\enskip \v x-x(\theta,\alpha)\v\leq \v x(\theta,\alpha)\v\enskip
\textrm{then,}\\
\v x-x(\theta,\alpha)\v\leq \v\gamma(\sigma,x)-x(\theta,\alpha)\v\leq \v
x(\theta,\alpha)\v\,.
\end{array}\right.
\end{equation}
The construction of this path will be made at the end of this Section.

It follows from (\ref{eqIV.4.78}) or (\ref{eqIV.4.79}) that
$$
\frac{1}{2}\,\v x(\theta,\alpha)\v\leq
\v\gamma(\sigma,x)-x(\theta,\alpha)\v\leq \bra x(\theta,\alpha)\ket\,.
$$
We write for $j=1,2$,
$$
\varphi _j(\theta,x,\alpha)=\varphi _j(\theta,0,\alpha)+\int^1_0
\frac{\partial \gamma}{\partial \sigma}\,(\sigma,x)\,\frac{\partial \varphi
_j}{\partial x}\,(\theta,\gamma(\sigma,x),\alpha)\,d\sigma\,.
$$

Using Proposition \ref{pIV.4.14} (ii) and (\ref{eqIV.4.75}) we obtain
\begin{equation*}\left\{
\begin{array}{l}
\varphi _j(\theta,x,\alpha)=\varphi _j(\theta,0,\alpha)+\int^1_0
\frac{\partial \gamma}{\partial
\sigma}\,(\sigma,x)\,\Phi(\theta,\gamma(\sigma,x),\alpha)\,d\sigma+R_j\,,\\
\v R_j\v\leq C_N \int^1_0 \v
x-x(\theta,\alpha)\v\,\frac{\v\gamma(\sigma,x)-x(\theta,\alpha)\v^N}{\bra\theta\ket^N}\,
d\sigma\,.
\end{array}\right.
\end{equation*}
Then by (\ref{eqIV.4.72}), (\ref{eqIV.4.78}) or (\ref{eqIV.4.79}) and
(\ref{eqIV.4.70}) we obtain (\ref{eqIV.4.73}).

Now let $\chi_0\in C^\infty _0(\R)$ be such that $0\leq \chi_0\leq 1$ and
$\chi_0(\sigma)=1$ if $\v\sigma\v\leq \frac 12$, $\chi_0(\sigma)=0$ if
$\v\sigma\v\geq 1$. Let us set
$$
\chi_1(\theta,x)=\chi_0\Big (\frac{x-x(\theta,\alpha)}{\bra
x(\theta,\alpha)\ket}\Big )\,.
$$
Now for $x$ in $\widetilde{\CO}_\delta(\theta)$ we set
\begin{equation}\label{eqIV.4.80}
\varphi
_3(\theta,x,\alpha)=\chi_1(\theta,x)\,\varphi _1(\theta,x,\alpha)+(1-\chi_1
(\theta,x))\,\varphi_2(\theta,x,\alpha)\,.
\end{equation}
On the support of $\chi_1$ we have $\v x-x(\theta,\alpha)\v\leq \bra
x(\theta,\alpha)\ket$ thus $\varphi _1$ is well defined. On the support of
$1-\chi_1(\theta,x)$ we have $\v x-x(\theta,\alpha)\v\geq \frac 12\, \bra
x(\theta,\alpha)\ket\geq \frac 12\,\v x(\theta,\alpha)\v$ so $\varphi _2$ is
well defined. Therefore $\varphi _3$ is well defined when
$x\in\widetilde{\CO}_\delta(\theta)$. We show now that $\varphi _3$ satisfies
all the conditions in Proposition \ref{pIV.4.14}. We have
\begin{equation}\label{eqIV.4.81}
\frac{\partial \varphi _3}{\partial
x}\,(\theta,x,\alpha)=\Big[\chi_1\,\frac{\partial \varphi _1}{\partial
x}+(1-\chi_1)\,\frac{\partial \varphi _2}{\partial x}+\frac{\partial
\chi_1}{\partial x}\,(\varphi _1-\varphi _2)\Big](\theta,x,\alpha)\,.
\end{equation}
On the support of $\frac{\partial \chi_1}{\partial x}$ we have $\v(\varphi
_1-\varphi _2)(\theta,x,\alpha)\v\leq C_N\,\frac{\v
x-x(\theta,\alpha)\v^{N+1}}{\bra\theta\ket^N}$ by (\ref{eqIV.4.73}). Moreover
we have by (\ref{eqIV.4.72})
$$
\Big\v\frac{\partial \chi_1}{\partial x}\,(\theta,x)\Big\v \leq \frac C{\bra
x(\theta,\alpha)\ket}\leq \frac C{\v x-x(\theta,\alpha)\v}
$$
and
\begin{equation*}
\begin{aligned}
\Big\v\Big[\chi_1\,\frac{\partial \varphi _1}{\partial
x}&+(1-\chi_1)\,\frac{\partial \varphi _2}{\partial x}-\Phi\Big
](\theta,x,\alpha)\Big\v\leq \v\chi_1\v\Big\v\frac{\partial \varphi
_1}{\partial x}\,(\theta,x,\alpha)-\Phi(\theta,x,\alpha)\Big\v\\
&+(1-\chi_1)\Big\v\frac{\partial \varphi _2}{\partial
x}\,(\theta,x,\alpha)-\Phi(\theta,x,\alpha)\Big\v\,.
\end{aligned}
\end{equation*}
It follows that the claim (ii) in Proposition \ref{pIV.4.14} holds for
$\varphi _3$. The point (iv) follows from (\ref{eqIV.4.81}) and
(\ref{eqIV.4.73}) for $N=1$ which gives $\v (\varphi _1-\varphi
_2)(\theta,x,\alpha)\v\leq C\,\delta$. The points (v) and (vi) are
straightforward. Let us show (iii). We have
$$
\frac{\partial \chi_1}{\partial \theta}\,(\theta,x)=\Big[-\frac{\dot
x(\theta,\alpha)}{\bra x(\theta,\alpha)\ket}-\frac{x(\theta,\alpha)\,\dot
x(\theta,\alpha)}{\bra x(\theta,\alpha)\ket^3}\,(x-x(\theta,\alpha))\Big]\,
\frac{\partial \chi_0}{\partial \sigma}\,(\cdots )\,.
$$
Since $\dot x(\theta,\alpha)$ is bounded we deduce from (\ref{eqIV.4.73}) that
$$
\Big\v\frac{\partial \chi_1}{\partial \theta}\,(\theta,x)\Big\v\leq \frac C{\v
x-x(\theta,\alpha)\v}\,.
$$
It follows then that
$$
\bigg\v\bigg(\frac{\partial \varphi _3}{\partial \theta}-\Big
(\chi_1\,\frac{\partial \varphi _1}{\partial
\theta}+(1-\chi_1)\,\frac{\partial \varphi _2}{\partial \theta}\Big
)\bigg)(\theta,x,\alpha)\bigg\v\leq C_N\,\frac{\v
x-x(\theta,\alpha)\v^N}{\bra\theta\ket^N}\,.
$$
Then (iii) follows easily.

Let us now take $\theta\in[\theta',+\infty [$ and
$x\in\widetilde{\CO}_\delta(\theta)$. Assume that


\begin{equation}\label{eqIV.4.82}
\frac{1}{2}\,\v x(\theta^*,\alpha)\v\leq \v x-x(\theta^*,\alpha)\v\leq \bra
x(\theta^*,\alpha)\ket\,.
\end{equation} 
We shall show that
\begin{equation}\label{eqIV.4.83}
\v(\varphi _4-\varphi _5)(\theta,x,\alpha)\v\leq C_N\,\frac{\v
x-x(\theta,\alpha)\v^N}{\bra\theta\ket^N}\,\v x-x(\theta^*,\alpha)\v\,.
\end{equation}
Let us take a path $\gamma$ satisfying (\ref{eqIV.4.74}), (\ref{eqIV.4.76}),
(\ref{eqIV.4.77}) and
\begin{equation}\label{eqIV.4.84}\left\{
\begin{array}{l}
\textrm{if}\enskip \v x-x(\theta^*,\alpha)\v\geq \v
x(\theta^*,\alpha)\v\enskip \textrm{then}\\
\v x(\theta^*,\alpha)\v\leq \v\gamma(\sigma,x)-x(\theta^*,\alpha)\v\leq \v
x-x(\theta^*,\alpha)\v\,,
\end{array}\right.
\end{equation}
\begin{equation}\label{eqIV.4.85}\left\{
\begin{array}{l}
\textrm{if}\enskip \v x-x(\theta^*,\alpha)\v\leq \v
x(\theta^*,\alpha)\v\enskip \textrm{then}\\
\v x-x (\theta^*,\alpha)\v\leq \v\gamma(\sigma,x)-x(\theta^*,\alpha)\v\leq \v
x(\theta^*,\alpha)\v\,,
\end{array}\right.
\end{equation}
\begin{equation}\label{eqIV.4.86}
\Big\v\frac{\partial \gamma}{\partial \sigma}\,(\sigma,x)\Big\v\leq K\,\v
x-x(\theta^*,\alpha)\v\,.
\end{equation}
Let us remark that in the two cases (\ref{eqIV.4.84}) or (\ref{eqIV.4.85}) we
have
$$
\v\gamma(\sigma,x)-x(\theta^*,\alpha)\v\leq 2\,\v x-x(\theta^*,\alpha)\v\leq
8\,\v x-x(\theta,\alpha)\v
$$
by Lemma \ref{lIV.4.18}.

Then using the same method as before we obtain easily (\ref{eqIV.4.83}). To
match $\varphi _4$ and $\varphi _5$ we set
$$
\chi_1(x)=\chi_0\Big (\frac{x-x(\theta^*,\alpha)}{\bra
x(\theta^*,\alpha)\ket}\Big )
$$
and we deduce from (\ref{eqIV.4.83}) that $\big\v\frac{\partial
\chi_1}{\partial x}\,(x)\big\v\leq \frac C{\v x-x(\theta^*,\alpha)\v}$. Then
we set
\begin{equation}\label{eqIV.4.87}
\varphi _6(\theta,x,\alpha)=[\chi_1\,\varphi _4+(1-\chi_1)\,\varphi
_5](\theta,x,\alpha)\,.
\end{equation}
It is then easy to see that $\varphi _6$ satisfies all the requirements of
Proposition \ref{pIV.4.14}.

Our last step is to match $\varphi _3$ and $\varphi _6$. With the notation
$\theta(c)$ introduced in Lemma \ref{lIV.4.18} let us set
$$
\theta_1=\theta\Big (\frac{11\,c_0}{30}\Big)\,,\quad \theta_2=\theta\Big
(\frac{12\,c_0}{30}\Big)\,,\quad \theta_3=\theta_0\Big
(\frac{14\,c_0}{30}\Big )\,.
$$
We have therefore according to (\ref{eqIV.4.69})
$$
\theta'<\theta_1<\theta_2<\theta_3<\overline \theta\,.
$$
Using (\ref{eqIV.4.87}), the fact that $\chi_1(x(\theta^*,\alpha))=1$ and
(\ref{eqIV.4.53}) we get
\begin{equation}\label{eqIV.4.88}
\varphi _6(\theta_2,x(\theta^*,\alpha),\alpha)=-\int^{\theta_2}_{\theta^*}
p(x(\theta^*,\alpha),\Phi(s,x(\theta^*,\alpha),\alpha)\,ds+\theta^*\,p(\alpha)+\frac
1{2i} \,\v\alpha_\xi\v^2\,.
\end{equation}
On the other hand we have
\begin{equation}\label{eqIV.4.89}
p(x(\theta^*,\alpha),\Phi(s,x(\theta^*,\alpha),\alpha))=-\frac{\partial
\varphi _3}{\partial s}\,(s,x(\theta^*,\alpha),\alpha)+A+B
\end{equation}
where
\begin{equation*}
\begin{aligned}
A&=p\big (x(\theta^*,\alpha),\Phi(s,x(\theta^*,\alpha),\alpha)\big )-p\Big
(x(\theta^*,\alpha),\,\frac{\partial \varphi _3}{\partial
x}\,(s,x(\theta^*,\alpha),\alpha)\Big)\\
B&=p\Big (x(\theta^*,\alpha),\,\frac{\partial \varphi _3}{\partial
x}\,(s,x(\theta^*,\alpha),\alpha)\Big )+\frac{\partial \varphi _3}{\partial s}
\,(s,x(\theta^*,\alpha),\alpha)\,.\end{aligned}
\end{equation*}
By the estimates proved in Proposition \ref{pIV.4.14} we have
\begin{equation}\label{eqIV.4.90}
\v A\v+\v B\v\leq C_N\,\frac{\v x(\theta^*,\alpha)-x(s,\alpha)\v^N}{\bra
s\ket^N}\leq C_N\,\frac{\v
x(\theta^*,\alpha)-x(\theta_2,\alpha)\v^N}{\bra\theta^*\ket^N}\,.
\end{equation}
Here we used the fact that for $s\in [\theta',\overline \theta]$ we have $\bra
s\ket\sim \bra\theta^*\ket$ (see Lemma \ref{lIV.4.17} (ii)) and
$\v x(\theta^*,\alpha)-x(s,\alpha)\v\sim\v
x(\theta^*,\alpha)-x(\theta_2,\alpha)\v$ by Lemma \ref{lIV.4.18}.

It follows from (\ref{eqIV.4.88}) and (\ref{eqIV.4.89}) that
$$
\varphi _6(\theta_2,x(\theta^*,\alpha),\alpha)=\int^{\theta_2}_{\theta^*}
\frac{\partial \varphi _3}{\partial
s}\,(s,x(\theta^*,\alpha),\alpha)\,ds+\theta^*\,p(\alpha)+\frac
1{2i}\,\v\alpha_\xi\v^2
+\int^{\theta_2}_{\theta^*} (\v A\v+\v B\v)\,ds\,.
$$
Therefore we obtain
\begin{equation}\label{eqIV.4.91}
\varphi _6(\theta_2,x(\theta^*,\alpha),\alpha)=\varphi
_3(\theta_2,x(\theta^*,\alpha),\alpha)-\varphi
_3(\theta^*,x(\theta^*,\alpha),\alpha)+\theta^*\,p(\alpha)+\frac
1{2i}\,\v\alpha_\xi\v^2+R
\end{equation}
where by (\ref{eqIV.4.90})
\begin{equation}\label{eqIV.4.92}
\v R\v\leq C_N\,\frac{\v
x(\theta^*,\alpha)-x(\theta_2,\alpha)\v^{N+1}}{\bra\theta^*\ket^N}\,.
\end{equation}
Now using (\ref{eqIV.4.80}) and (\ref{eqIV.4.53}) we have
$$
\varphi _3(\theta^*,x(\theta^*,\alpha),\alpha)=\varphi
_1(\theta^*,x(\theta^*,\alpha),\alpha)=\theta^*\,p(\alpha)+\frac 1{2i}\,\v
\alpha_\xi\v^2\,,
$$
so we obtain
\begin{equation}\label{eqIV.4.93}
\varphi _6(\theta_2,x(\theta^*,\alpha),\alpha)=\varphi
_3(\theta_2,x(\theta^*,\alpha),\alpha)+R
\end{equation}
where $R$ satisfies (\ref{eqIV.4.92}).

Now let $\theta\in [\theta_1,\theta_3]$. We set
$$
\CO'_\delta(\theta)=\Big\{x\in\R^n : x\cdot \alpha_\xi\leq
\frac{c_0}{20}\,\bra x\ket\,\v\alpha_\xi\v\,, \v x-x(\theta,\alpha)\v\leq
\frac{\delta}{40}\,\bra\theta\ket\Big\}\,.
$$
Let $x\in \CO'_\delta(\theta)$. We can find a path $\gamma$ joining $x$ to
$x(\theta^*,\alpha))$ such that $\gamma\subset \CO'_\delta(\theta)$ and there
exists $K\geq 0$ such that
\begin{equation}\label{eqIV.4.94}
\Big\v\frac{\partial \gamma}{\partial \sigma}\,(\sigma,x)\Big\v\leq K\,\v
x-x(\theta^*,\alpha)\v\,, \v\gamma(\sigma,s)-x(\theta_2,\alpha)\v\leq
K\v x-x(\theta_2,\alpha)\v
\end{equation}
Indeed if $\v x-x(\theta^*,\alpha)\v\leq \v x(\theta^*,\alpha)\v$ we set
$$
\gamma(\sigma,x)=\sigma\,x+(1-\sigma)\,x(\theta^*,\alpha)
$$
and Lemma \ref{lIV.4.16} show that
\begin{equation*}
\begin{aligned}
\gamma(\sigma,x))\cdot \alpha_\xi=\sigma\,x\cdot \alpha_\xi&\leq
\frac{c_0}{20}\,\sigma(1+\v x\v)\,\v\alpha_\xi\v\\
&\leq \frac{c_0}{20}\,(\sigma+\sigma\,\v x\v+(1-\sigma)\v
x(\theta^*,\alpha)\v))\,\v\alpha_\xi\v\\
&\leq \frac{c_0}{20}\,\sqrt 2\,(1+\v\gamma(\sigma,x)\v)\v \alpha_\xi\v\leq
\frac{c_0}{10}\,\bra\gamma(\sigma,x)\ket\,\v\alpha_\xi\v\,.
\end{aligned}
\end{equation*}
If $\v x-x(\theta^*,\alpha)\v>\v x(\theta^*,\alpha)\v$ we take $\gamma$ to be
the union of the two segments joining $x$ and $x(\theta^*,\alpha)$ to $0$ and
we obtain with $y=x$ or $y=x(\theta^*,\alpha)$ for $t\in[0,1]$,
$$
t\,y\cdot \alpha_\xi\leq t\,\frac{c_0}{20}\,\bra y\ket\,\v\alpha_\xi\v\leq
\frac{c_0}{20}\,\bra t\,y\ket\,\v\alpha_\xi\v\,.
$$
Since $0$ belongs to $\CO'_\delta(\theta)$ these two segments are contained in
$\CO'_\delta(\theta)$.

Let us prove the estimate on $\gamma$ given in (\ref{eqIV.4.94}). 

If $\v
x-x(\theta^*,\alpha)\v\leq \v x(\theta^*,\alpha)\v$ we have
\begin{equation*}
\begin{aligned}
\v\sigma\,x+(1-\sigma)\,x(\theta^*,\alpha)-x(\theta_2,\alpha)\v&\leq \sigma\v
x-x(\theta^*,\alpha)\v+(1-\sigma) \v x(\theta^*,\alpha)-x(\theta_2,\alpha)\v\\
&\leq K\,\v x-x(\theta_2,\alpha)\v\,.
\end{aligned}
\end{equation*}
If $\v x-x(\theta^*,\alpha)\v>\v x(\theta^*,\alpha)\v$ we have
\begin{equation*}
\begin{aligned}
\v t\,x-x(\theta_2,\alpha)\v&\leq t\,\v x-x(\theta^*,\alpha)\v+(1-t)\,\v
x(\theta^*,\alpha)\v+\v x(\theta^*,\alpha)-x(\theta_2,\alpha)\v\\
&\leq \v x-x(\theta^*,\alpha)\v+\v x(\theta^*,\alpha)-x(\theta_2,\alpha)\v\\
&\leq K\,\v x-x(\theta_2,\alpha)\v
\end{aligned}
\end{equation*}
again by Lemma \ref{lIV.4.18}. Moreover
\begin{equation*}
\begin{aligned}
\v t\,x(\theta^*,\alpha)-x(\theta_2,\alpha)\v&\leq \v
x(\theta_2,\alpha)-x(\theta^*,\alpha)\v+(1-t)\,\v x(\theta^*,\alpha)\v\\
&\leq K\,\v x-x(\theta_2,\alpha)\v\,.
\end{aligned}
\end{equation*}
Concerning the estimate on $\frac{\partial \gamma}{\partial \sigma}$, if $\v
x-x(\theta^*,\alpha)\v\leq \v x(\theta^*,\alpha)\v$ it is straightforward by
Lemma \ref{lIV.4.18}. If $\v x(\theta^*,\alpha)\v\leq \v
x-x(\theta^*,\alpha)\v$ the same Lemma shows that $\v x(\theta^*,\alpha)\v\leq
K\,\v\theta^*-\theta_2\v\leq K\,\v x(\theta^*,\alpha)-x(\theta_2,\alpha)\v$
and $\v x\v\leq \v x-x(\theta_2,\alpha)\v+\v
x(\theta_2,\alpha)-x(\theta^*,\alpha)\v+\v x(\theta^*,\alpha)\v\leq K'\,\v
x-x(\theta_2,\alpha)\v$. Thus (\ref{eqIV.4.94}) is entirely proved.

Now for $j=3$ or $6$ we can write
\begin{equation*}
\begin{aligned}
\varphi _j(\theta_2,x,\alpha)&=\varphi
_j(\theta_2,x(\theta^*,\alpha),\alpha)+\int^1_0 \frac{\partial \gamma}{\partial
\sigma}\, (\sigma,x)\,\frac{\partial \varphi _j}{\partial
x}\,(\theta_2,\gamma(\sigma,x),\alpha)\,d\sigma\\
&=\varphi _j(\theta_2,x(\theta^*,\alpha),\alpha)
+\int^1_0 \frac{\partial \gamma}{\partial
\sigma}\, (\sigma,x)\,\Phi(\theta_2,\gamma(\sigma,x),\alpha)\,d\sigma+R_j
\end{aligned}
\end{equation*}
where
$$
\v R_j\v\leq C_N\,\frac{\v x-x(\theta_2,\alpha)\v^{N+1}}{\bra\theta_2\ket^N}
$$
by Proposition \ref{pIV.4.14} and (\ref{eqIV.4.94}). By  (\ref{eqIV.4.91}) and
(\ref{eqIV.4.92}) we have
\begin{equation}\label{eqIV.4.95}
\v(\varphi _3-\varphi _6)(\theta_2,x,\alpha)\v\leq C_N\,\frac{\v
x-x(\theta_2,\alpha)\v^{N+1}}{\bra\theta_2\ket^N}\,.
\end{equation}
Now for $\theta\in [\theta_1,\theta_3]$ and $j=3$ or $6$ we can write
$$
\varphi _j(\theta,x,\alpha)=\varphi _j(\theta_2,x,\alpha)+\int^1_0
(\theta-\theta_2)\,\frac{\partial \varphi _j}{\partial
\theta}\,(\sigma\theta+(1-\sigma)\,\theta_2,x,\alpha)\,d\sigma
$$
so by Proposition \ref{pIV.4.14} and Lemma \ref{lIV.4.18} we have
\begin{equation}\label{eqIV.4.96}
\varphi _j(\theta,x,\alpha)=\varphi _j(\theta_2,x,\alpha)-\int^1_0
(\theta-\theta_2)\,p(x,\Phi(\sigma\theta+(1-\sigma)\,\theta_2,x,\alpha))\,d\sigma+A\,,
\end{equation}
where
$$
\v A\v\leq C_N\,\frac{\v
x-x(\sigma\theta+(1-\sigma)\,\theta_2,\alpha)\v^N}{\bra\sigma\theta+(1-\sigma)\,\theta_2\ket^N\,.}
$$
\begin{proposition}\sl\label{pIV.4.19}
\begin{equation}\label{eqIV.4.97}\left\{
\begin{array}{l}
\textrm{For all}\enskip \theta\enskip \textrm{in}\enskip
[\theta_1,\theta_3]\enskip \textrm{and all}\enskip \sigma\enskip
\textrm{in}\enskip [0,1]\enskip \textrm{we have}\\
\enskip (i)\enskip \bra\sigma\theta+(1-\sigma)\,\theta_2\ket\geq
K_1\,\bra\theta\ket\\ 
\enskip (ii)\enskip \v
x-x(\sigma\theta+(1-\sigma)\,\theta_2,\alpha)\v\leq K_2\,\v
x-x(\theta,\alpha)\v\,.
\end{array}\right.
\end{equation}
\end{proposition}

\noindent {\bf Proof}

{\bf Case 1 :} $\theta\in[\theta_2,\theta_3]$. We write with
$\theta^\sigma=\sigma\theta+(1-\sigma)\,\theta_2$,
\begin{equation}\label{eqIV.4.98}\left\{
\begin{array}{l}
\v x-x(\theta,\alpha)\v^2=I+II\quad \textrm{where}\\
\enskip I=\v x-x(\theta^\sigma,\alpha)\v^2+\v
x(\theta^\sigma,\alpha)-x(\theta,\alpha)\v^2\\
II=2(x-x(\theta^\sigma,\alpha))\cdot
(x(\theta^\sigma,\alpha)-x(\theta,\alpha))\,.
\end{array}\right.
\end{equation}
Since $x(\theta,x)-x(\theta^\sigma,\alpha)=\int^\theta_{\theta^\sigma} \dot
x(s,\alpha)\,ds=2(\theta-\theta^\sigma)\,\alpha_\xi+\CO(\varepsilon\,\v\theta-\theta^\sigma\v)$
we have $II=-4(\theta-\theta^\sigma)(x-x(\theta^\sigma,\alpha))\cdot
\alpha_\xi+\CO(\varepsilon\,I)$. Now in $\CO'_\delta$ we have $x\cdot
\alpha_\xi\leq \frac{c_0}{20}\,\bra x\ket\,\v\alpha_\xi\v$~; moreover by Lemma
\ref{lIV.4.18} (i) we have $x(\theta^\sigma,\alpha)\cdot \alpha_\xi\geq
\frac{12\,c_0}{30}\,\bra x(\theta^\sigma,\alpha)\ket\,\v\alpha_\xi\v$ since
$\theta^\sigma\geq \theta_2=\theta\big (\frac{12\,c_0}{30}\big )$. It follows
that $II\geq -\frac{c_0}5\,(\theta-\theta^\sigma)\bra
x\ket\,\v\alpha_\xi\v+\frac{8\,c_0}5\,\bra
x(\theta^\sigma,\alpha)\ket\,\v\alpha_\xi\v-\CO(\varepsilon\,I)$. Therefore we
obtain
$$
II\geq -\frac{c_0}5\,(\theta-\theta^\sigma)\bra
x(\theta^\sigma,\alpha)\ket\,\v\alpha_\xi\v-\frac{c_0}5\,(\theta-\theta^\sigma)\,\v
x-x(\theta^\sigma,\alpha)\v\,\v\alpha_\xi\v+\frac{8\,c_0}5\,\bra
x(\theta^\sigma,\alpha)\ket\,\v\alpha_\xi\v-\CO(\varepsilon\,I)\,.
$$
The second term in the right hand side can be bounded by $\frac{c_0}{10}\,I$.
Using (\ref{eqIV.4.98}) we obtain
$$
\v x-x(\theta,\alpha)\v^2=I+II\geq \Big
(1-\frac{c_0}{10}-\varepsilon\,K)\,I+\frac{7\,c_0}5\,\bra
x(\theta^\sigma,\alpha)\ket\,\v\alpha_\xi\v\,.
$$
Taking $c_0$ and $\varepsilon$ small enough we obtain $I\leq 2\,\v
x-x(\theta,\alpha)\v^2$ which implies since $\v\theta-\theta^\sigma\v\leq
2\,\v x(\theta,\alpha)-x(\theta^\sigma,\alpha)\v\leq 2\sqrt I$, that
$\v\theta-\theta^\sigma\v\leq 2\,\v x-x(\theta,\alpha)\v\leq
2\delta\,\bra\theta\ket$ so $\bra\theta\ket\leq
\bra\theta^\sigma\ket+2\delta\,\bra\theta\ket$ and therefore
$\bra\theta^\sigma\ket\geq \frac 12\,\bra\theta\ket$ since $\delta$ is small.
This proves the claim (i) of (\ref{eqIV.4.97}).

To prove (ii) we just use the fact that
$$
\v x(\theta,\alpha)-x(\theta^\sigma,\alpha)\v\leq 3\,\v
\theta-\theta^\sigma\v\leq 6\,\v x-x(\theta,\alpha)\v\,.
$$

\noindent {\bf Case 2 :} $\theta\in[\theta_1,\theta_2]$. 

The point (i) in
(\ref{eqIV.4.97}) is obvious in this case since
$\sigma\theta+(1-\sigma)\,\theta_2\geq \theta$.

By the same computation as above, since $\theta\geq \theta_1$, we will have
\begin{equation}\label{eqIV.4.99}\left\{
\begin{array}{l}
\v x-x(\theta,\alpha)\v\geq \frac 12\,(\v x-x(\theta_1,\alpha)\v+\v
x(\theta_1,\alpha)-x(\theta,\alpha)\v)\,,\\
\frac 12\,\v \theta-\theta_1\v\leq \v x(\theta,\alpha)-x(\theta_1,\alpha)\v\leq
6\,\v\theta-\theta_1\v\,.
\end{array}\right.
\end{equation}
On the other hand we claim that we have
\begin{equation}\label{eqIV.4.100}
\frac{9\,c_0}{30}\,\bra x(\theta,\alpha)\ket\,\v\alpha_\xi\v\leq 3\,\v
x-x(\theta,\alpha)\v\,.
\end{equation}
Indeed we have
$$
(x-x(\theta,\alpha))\cdot \alpha_\xi=x\cdot \alpha_\xi-x(\theta,\alpha)\cdot
\alpha_\xi\leq \frac{c_0}{20}\,\bra
x\ket\,\v\alpha_\xi\v-\frac{11\,c_0}{30}\,\bra
x(\theta,\alpha)\ket\,\v\alpha_\xi\v
$$
by Lemma \ref{lIV.4.18} (i) since $\theta\geq \theta_1=\theta\big
(\frac{11\,c_0}{30}\big)$. Thus
$$
(x-x(\theta,\alpha))\cdot \alpha_\xi\leq \frac{c_0}{20}\,\bra
x(\theta,\alpha)\ket\,\v\alpha_\xi\v+\frac{c_0}{20}\,\v
x-x(\theta,\alpha)\v\,\v\alpha_\xi\v-\frac{11\,c_0}{30}\,\bra
x(\theta,\alpha)\ket\,\v\alpha_\xi\v\,.
$$
It follows that
$$
\frac{9\,c_0}{30}\,\bra x(\theta,\alpha)\ket\,\v\alpha_\xi\v\leq
\frac{c_0}{20}\,\v
x-x(\theta,\alpha)\v\,\v\alpha_\xi\v-(x-x(\theta,\alpha))\cdot \alpha_\xi
$$
from which (\ref{eqIV.4.100}) follows easily since $\v\alpha_\xi\v\leq 2$ and
$\frac{c_0}{20}\,\v\alpha_\xi\v\leq 1$. Now
$\theta^\sigma=\sigma\theta+(1-\sigma)\,\theta_2$ belongs to
$[\theta_1,\theta_2]$ for $\sigma\in[0,1]$. Since by Lemma \ref{lIV.4.18}
(iii) the function $\theta(c)$ is strictly increasing there exists a unique
$c_\sigma\in\big[\frac{11\,c_0}{30}\,,\frac{12\,c_0}{30}\big]$ such that
$\theta^\sigma=\theta(c_\sigma)$. Now we have $\v
x-x(\theta^\sigma,\alpha)\v\leq \v x-x(\theta_1,\alpha)\v+\v
x(\theta_1,\alpha)-x(\theta^\sigma,\alpha)\v$ which implies
\begin{equation}\label{eqIV.4.101}
\v x-x(\theta^\sigma,\alpha)\v\leq \v
x-x(\theta_1,\alpha)\v+6\,\Big\v\theta\Big
(\frac{11\,c_0}{30}\Big)-\theta(c_\sigma)\Big\v\,.
\end{equation}
We claim that
\begin{equation}\label{eqIV.4.102}
\Big\v\theta\Big (\frac{11\,c_0}{30}\Big )-\theta(c_\sigma)\Big\v\leq
\frac{c_0}{15}\,\sup_{s\in[\theta_1,\theta_2]}\,\bra x(s,\alpha)\ket\,.
\end{equation}
To see this we compute $\theta'(c)$. Recall (see Lemma \ref{lIV.4.18}) that
$g(\theta(c),c)=0$ for $c\in\big[\frac{c_0}{10}\, ,\frac{c_0}2\big]$.It
follows that $\frac{\partial g}{\partial
s}\,(\theta(c),c)\,\theta'(c)+\frac{\partial g}{\partial c}\,(\theta(c),c)=0$.
Now we have $\frac{\partial g}{\partial c}\,(s,c)=-\bra
x(s,\alpha)\ket\,\v\alpha_\xi\v$ and $\frac{\partial g}{\partial
s}\,(s,c)=\dot x(s,\alpha)\cdot \alpha_\xi-c\,\frac{x(s,\alpha)\cdot \dot x
(s,\alpha)}{\bra x(s,\alpha)\ket}\,\v\alpha_\xi\v$, which shows that
$\frac{\partial g}{\partial s}\,(s,c)=2 (\v\alpha_\xi\v^2+\CO
(\varepsilon+c)$. Therefore we have $\v\theta'(c)\v\leq 2\,\bra
x(\theta(c),\alpha\ket$ and we obtain (\ref{eqIV.4.102}). The last step
consists in showing that
\begin{equation}\label{eqIV.4.103}
\sup_{\theta\in[\theta_1,\theta_2]}\,\bra x(s,\alpha)\ket\leq 2\,\bra
x(\theta,\alpha)\ket\,.
\end{equation}
To see this let us set $h(s)=\bra x(s,\alpha)\ket$. Then
$h'(s)=\frac{x(s,\alpha)\cdot \dot x(s,\alpha)}{\bra x(s,\alpha)\ket}$. Thus
$h'(s)=\frac{2x(s,\alpha)\cdot \alpha_\xi}{\bra x(s,\alpha)\ket}+\CO
(\varepsilon)$. Now since $s\in[\theta_1,\theta_2]$ we have
$\frac{11\,c_0}{30}\,\bra x(s,\alpha)\ket\,\v\alpha_\xi\v\leq x(s,\alpha)\cdot
\alpha_\xi\leq \frac{12\,c_0}{30}\,\bra x(s,\alpha)\ket\,\v\alpha_\xi\v$ so
$0<h'(s)\leq 2\,c_0+\CO(\varepsilon)$ and therefore if
$s_1,s_2\in[\theta_1,\theta_2]$, $\v h(s_1)-h(s_2)\v\leq (2\,c_0+\CO
(\varepsilon))\,\v s_1-s_2\v$. Let us take $s_1=\theta(c)$, $s_2=\theta(c')$
with $\frac{11\,c_0}{30}\leq c$, $c'\leq \frac{12\,c_0}{30}$. Then
\begin{equation}\label{eqIV.4.104}
\v\bra x(\theta(c),\alpha)\ket-\bra x(\theta(c'),\alpha)\ket\v\leq (2\,c_0+\CO
(\varepsilon))\v \theta(c)-\theta(c')\v\,.
\end{equation}
On the other hand
$$
x(\theta(c),\alpha)-x(\theta(c'),\alpha)=2(\theta(c)-\theta(c'))\cdot
\alpha_\xi+\CO (\varepsilon\,\v\theta(c)-\theta(c')\v)
$$
which implies that
$$
x(\theta(c),\alpha)\cdot \alpha_\xi-x(\theta(c'),\alpha)\cdot
\alpha_\xi=2(\theta(c)-\theta(c'))\,\v\alpha_\xi\v^2+\CO
(\varepsilon\,\v\theta(c)-\theta(c')\v
$$
and, by definition of $\theta(c)$ (see Lemma \ref{lIV.4.18})
$$
\big[c \bra x(\theta(c),\alpha\ket-c'\,\bra
x(\theta(c'),\alpha)\ket\big]\,\v\alpha_\xi\v=2(\theta(c)-\theta(c'))\,\v\alpha_\xi\v^2+\CO
(\varepsilon\,\v \theta(c)-\theta(c')\v)\,.
$$
Combining with (\ref{eqIV.4.104}) we obtain
$$
\v \bra x(\theta(c),\alpha)-\bra x(\theta(c'),\alpha)\ket\v\leq 2 (2\,c_0+\CO
(\varepsilon))\big[c\,\v\bra x(\theta(c),\alpha)\ket-\bra
x(\theta(c'),\alpha\ket\v
+\v c-c'\v \bra x(\theta(c'),\alpha\ket\big]\,.
$$
Since $c_0$ and $\varepsilon$ are small enough we obtain
$$
\v\bra x(\theta(c),\alpha)\ket-\bra x(\theta(c'),\alpha)\ket\v\leq 2\,\v
c-c'\v\bra x(\theta(c'),\alpha\ket\leq \frac{c_0}{15}\,\bra
x(\theta(c'),\alpha)\ket
$$
which shows that all the $\bra x(\theta(c),\alpha)\ket$ are equivalent in
$[\theta_1,\theta_2]$, more precisely taking $s=\theta(c'),\theta=\theta(c)$
we obtain
$$
\sup_{s\in[\theta_1,\theta_2]}\,\bra x(s,\alpha)\ket\leq 2\,\bra
x(\theta,\alpha)\ket\,,
$$
which is (\ref{eqIV.4.103}).

Finally using (\ref{eqIV.4.101}), (\ref{eqIV.4.99}), 
(\ref{eqIV.4.102}),  (\ref{eqIV.4.103}) and (\ref{eqIV.4.100}) we obtain
$$
\v x-x(\theta^\sigma,\alpha)\v\leq K\,\v x-x(\theta,\alpha)\v
$$  
which is Proposition \ref{pIV.4.19} (iii) in the case 2. \cqfd

Now using Proposition \ref{pIV.4.19}, (\ref{eqIV.4.95}), (\ref{eqIV.4.96}) we
obtain 
\begin{equation}\label{eqIV.4.105}
\v(\varphi _3-\varphi _6)(\theta_2,x,\alpha)\v\leq C_N\,\frac{\v
x-x(\theta,\alpha)\v^{N+1}}{\bra\theta\ket^N}
\end{equation}

\begin{equation}\label{eqIV.4.106}\left\{
\begin{array}{l}
\varphi _j(\theta,x,\alpha)=\varphi _j(\theta_2,x,\alpha)-\int^1_0
p(x,\Phi(\theta^\sigma,x,\alpha)\,d\sigma+A\,,\quad j=3,6\,,\\
\v A\v\leq C_N\,\frac{\v x-x(\theta,\alpha)\v^{N+1}}{\bra\theta\ket^N}\,.
\end{array}\right.
\end{equation}
Let now $\chi_2\in C^\infty (\R)$ be such that $\chi_2(s)=1$ if $s\geq 1$,
$\chi_2(s)=0$ if $s\leq 0$ and set $\chi_3(\theta)=\chi_2\big
(\frac{\theta-\theta_1}{\theta_3-\theta_1}\big)$. Then let us set
\begin{equation}\label{eqIV.4.107}
\varphi (\theta,x,\alpha)=\chi_3(\theta)\,\varphi
_3(\theta,x,\alpha)+(1-\chi_3(\theta))\,\varphi _6(\theta,x,\alpha)\,.
\end{equation}
We have
$$
\frac{\partial \varphi }{\partial
\theta}\,(\theta,x,\alpha)=\chi_3(\theta)\,\frac{\partial \varphi _3}{\partial
\theta}\,(\theta,x,\alpha)+(1-\chi_3(\theta))\,\frac{\partial \varphi
_6}{\partial \theta}\,(\theta,x,\alpha)+\frac{\partial \chi_3}{\partial
\theta}\,(\theta)(\varphi _3-\varphi _6)(\theta,x,\alpha)\,.
$$
Now we deduce from (\ref{eqIV.4.105}) and (\ref{eqIV.4.106}) that on the
support of $\frac{\partial \chi_3}{\partial \theta}$ we have
$$
\v(\varphi _3-\varphi _6)(\theta,x,\alpha)\v\leq C_N\,\frac{\v
x-x(\theta,\alpha)\v^{N+1}}{\bra\theta\ket^N}\,.
$$
By Proposition \ref{pIV.4.14} for $\varphi _3$ and $\varphi _6$ we have for
$j=3$ or $6$,
$$
\Big\v\frac{\partial \varphi _j}{\partial
x}\,(\theta,x,\alpha)-\Phi(\theta,x,\alpha)\Big\v\leq C_N\,\frac{\v
x-x(\theta,\alpha)\v^N}{\bra\theta\ket^N}
$$
therefore this is also true for $\varphi $ and since for $j=3$ or $6$,
$$
\Big\v\frac{\partial \varphi _j}{\partial
\theta}\,(\theta,x,\alpha)-p(x,\Phi(\theta,x,\alpha))\Big\v\leq C_N\,\frac{\v
x-x(\theta,\alpha)\v^N}{\bra\theta\ket^N}\,,
$$
the function $\varphi $ defined in (\ref{eqIV.4.107}) satisfies all the
requirements of Proposition \ref{pIV.4.14}.

The proof of Proposition \ref{pIV.4.14} will be therefore complete when we
will construct the path $\gamma(\sigma,x)$ satisfing (\ref{eqIV.4.74}) to
(\ref{eqIV.4.79}).

\noindent {\bf Construction of $\gamma(\sigma,x)$}

Let us set $a=x(\theta,\alpha)$. We first show that we can join any point $x$
to a point $a-\v x-a\v\,\alpha_\xi$ by path remaining in the set
$$
\Big\{y\in\R^n : y\cdot \alpha_\xi\leq \frac{c_0}{10}\,\bra
y\ket\,\v\alpha_\xi\v\,,\enskip \v y-a\v=\v x-a\v\Big\}\,.
$$
Making rotations we may without loss of generality assume that
$\frac{\alpha_\xi}{\v\alpha_\xi\v}=(-1,0,\ldots, 0)$, $a=(a_1,a_2,0,\ldots
,0)$,
$x=(x_1,x_2,x_3,0,\ldots ,0)$. Therefore it will be suficient to restrict
ourselves to the dimension three. We will construct our path on planes so we
begin by the dimension two. Let us set with $D\in]0,1[$, $k>0$,
\begin{equation*}
\begin{aligned}
\CC&=\big\{y\in\R^2 : \v y-a\v^2=\v x-a\v^2\big\}\,,\\
\CH&=\big\{y\in\R^2 : -y_1=D\,\sqrt{k^2+y^2_2}\big\}\,,\\
\CD&=\big\{y\in\R^2 : -y_1\leq D\,\sqrt{k^2+y^2_2}\big\}\,.
\end{aligned}
\end{equation*}
\begin{lemma}\sl\label{lIV.4.20}
$\CD^c=\R^2\setminus\CD$ is strictly convex.
\end{lemma}

\noindent {\bf Proof } This follows easily from the strict convexity of the
function $g(t)=\sqrt{k^2+t^2}$. \cqfd
\begin{lemma}\sl\label{lIV.4.21}
\begin{itemize}
\item[(i)] Let $b\in\CD$ and $u=(1,y)$ with $\v y\v\leq 1$. 

Then for all $t>0$
we have $b+t\,u\in\CD\setminus\CH$.
\item[(ii)] Let $b\in\overline {\CD^c}$ and $v=(-1,y)$ with $\v y\v\leq 1$.
Then for all $t>0$ we have $b+t\,v\in\CD^c$.
\end{itemize}
\end{lemma}

\noindent {\bf Proof }

\noindent (i) Let $b=(b_1,b_2)$ and $h(t)=b_1+t+D\,\sqrt{k^2+(b_2+t\,y)}$.
Then $h(0)=b_1+D\,\sqrt{k^2+b^2_2}\geq 0$ and
$h'(t)=1+\frac{D(b_2+t\,y)\,y}{\sqrt{k^2+(b_2+t\,y)^2}}$. Since
$D\,\v y\v<1$ we have $\frac{D\,\v b_2+t\,y\v\,\v
y\v}{\sqrt{k^2+(b_2+t\,y)^2}}<1$. It follows that $h'(t)>0$ so $h(t)>h(0)\geq
0$.

The proof of (ii) is the same. \cqfd

Assume that $\CC\cap\CH$ contains at least two different points (otherwise
$\CC\setminus (\CC\cap\CH)$ would be connected). Let us set
\begin{equation}\label{eqIV.4.108}\left\{
\begin{array}{l}
M_\theta=a+\v x-a\v \begin{pmatrix}\cos\theta\\
\sin\theta\end{pmatrix}\,,\enskip \theta\in[0,2\pi [\,,\\
\theta_1=\inf \{\theta\in[0,2\pi [ : M_\theta\in\CC\cap\CH\}\,,\\
\theta_2=\sup \{\theta\in[0,2\pi [ : M_\theta\in\CC\cap\CH\}\,.
\end{array}\right.
\end{equation}
\begin{remark}\sl\label{rIV.4.22}
(i) If $\theta\in\big[0,\frac{\pi}{4}\big] \cup \big[2\pi -\frac\pi 4\,,2\pi
\big]$ we have $\overrightarrow{a\,M_\theta}=\v x-a\v\begin{pmatrix} \cos
\theta\\ \sin\theta\end{pmatrix}$ with $\cos \theta>0$ and
$\big\v\frac{\sin\theta}{\cos\theta}\big\v\leq 1$. Since $a\in\CD$ Lemma
\ref{lIV.4.21} (i) implies that $M_\theta\in\CD\setminus\CH$. It follows that
we have
$$
\frac\pi 4<\theta_1<\theta_2<2\pi -\frac\pi 4\,.
$$

(ii) We cannot have $\theta_1\in\big]\frac\pi 4\,,\frac\pi 2\big]$ and
$\theta_2\in\big[\frac{3\pi }2\,, 2\pi -\frac\pi 4\big[$.
Indeed if this was true then by Lemma \ref{lIV.4.20} the segment
$]M_{\theta_1},M_{\theta_2}[$ would be in $\CD^c$. But $\sin \theta_1>0$ and
$\sin \theta_2<0$ so there exists $t\in]0,1[$ such that $t\,\sin
\theta_1+(1-t)\,\sin \theta_2=0$~; then
$N_t=t\,M_{\theta_1}+(1-t)\,M_{\theta_2}=a+\v x-a\v \begin{pmatrix}
t\,\cos \theta_1+(1-t)\,\cos \theta_2\\
0\end{pmatrix}=a+\alpha\begin{pmatrix}1\\0\end{pmatrix}$ with $\alpha>0$ since
$\cos \theta_1>0$ and $\cos \theta_2>0$. By Lemma \ref{lIV.4.21}
(i) $N_t\in\CD$ since $a\in\CD$ which is in contradiction with
$N_t\in]M_{\theta_1},M_{\theta_2}[\subset\CD^c$.

(iii) If $\theta_1\in\big]\frac\pi 4\,,\frac\pi 2\big]$ then for all $\theta$
in $]\theta_1,\pi ]$ we have $M_\theta\in\CD^c$ which implies
$\theta_2\in\big]\pi ,\,\frac{3\pi }2\big]$ by (ii). Indeed we have
$\overrightarrow{M_{\theta_1}M_\theta}=\v x-a\v\begin{pmatrix}\cos \theta-\cos
\theta_1\\ \sin\theta-\sin \theta_1\end{pmatrix}$~; since for $\theta\in[0,\pi
]$, $\cos \theta$ is decreasing we have $\cos \theta-\cos \theta_1<0$ and
$$
\Big\v\frac{\sin \theta-\sin \theta_1}{\cos \theta-\cos
\theta_1}\Big\v=\Big\v\cotg \Big(\frac{\theta+\theta_1}2\Big)\Big\v\leq 1
$$
since $\frac \pi 4\leq \frac{\theta+\theta_1}2\leq \frac{3\pi }4$. Then Lemma
\ref{lIV.4.21} (ii) implies that $M_\theta\in\CD^c$.
\end{remark}

It follows from Remark \ref{rIV.4.22}, (i), (ii), (iii) that we have else
$\theta_1\in\big]\frac\pi 2\,,\pi \big[$ or $\theta_2\in\big]\pi ,\,\frac{3\pi
}2\big[$. By symmetry it is enough  to consider one case. Therefore we shall
assume in the sequel that $\frac \pi 2<\theta_1<\theta_2<2\pi -\frac \pi 4$,
$\theta_1\in \big]\frac\pi 2\,,\pi \big[$. We claim that
\begin{equation}\label{eqIV.4.109}
M_\theta\in\CD^c\enskip \textrm{for all}\enskip \theta\enskip
\textrm{in}\enskip ]\theta_1,\theta_2[\,.
\end{equation}
We split the proof in two cases.

\noindent {\bf Case 1 :} $\theta_2\leq \frac{3\pi }2$. 

Since $\sin \theta$ is decreasing on $\big[\frac \pi 2\,,\frac{3\pi }2\big]$
we have

 if $\frac \pi 2< \theta_1<\theta_2\leq \frac{3\pi }2$, $\sin
\theta_1>\sin\theta>\sin\theta_2$,

if $\frac \pi 2<\theta_1<\theta<3\pi -\theta_2\leq \frac{3\pi }2$, $\sin
\theta_1>\sin \theta>\sin (3\pi -\theta_2)=\sin \theta_2$.

Let us set
$$
N_t=t\,M_{\theta_1}+(1-t)\,M_{\theta_2}=a+\v x-a\v\begin{pmatrix}t\,\cos
\theta_1+(1-t)\,\cos \theta_2\\t\,\sin \theta_1+(1-t)\,\sin
\theta_2\end{pmatrix}\,.
$$
Now there exists $t\in]0,1[$ such that
$$
t\,\sin \theta_1+(1-t)\,\sin \theta_2=\sin \theta
$$
and since $\theta\in\big[\frac\pi 2\,,\frac{3\pi }2\big]$ we have
$\theta=\Arc\sin (t\,\sin \theta_1+(1-t)\,\sin \theta_2)+\pi$. Then $\cos
\theta=-\cos (\Arc\sin (t\,\sin \theta_1+(1-t)\,\sin
\theta_2))=-\sqrt{1-(t\,\sin \theta_1+(1-t)\,\sin \theta_2)^2}$, $\cos
\theta<-t\sqrt{1-\sin^2\,\theta_1}-(1-t\sqrt{1-\sin^2\,\theta_2}\leq t\,\cos
\theta_1+(1-t)\,\cos \theta_2$. Here we have used the strict convexity of the
function $\sqrt{1-x^2}$. Since $N_t\in\CD^c$ and
$M_\theta=N_t+\alpha\begin{pmatrix}-1\\0\end{pmatrix}$ where $\alpha>0$ (see
(\ref{eqIV.4.108})) we deduce from Lemma \ref{lIV.4.21} (ii) that
$M_\theta\in\CD^c$ which proves (\ref{eqIV.4.109}) in case 1.

\noindent {\bf Case 2 :} $\theta_2>\frac{3\pi }2$ and $\theta<3\pi -\theta_2$.

Since $\theta_2<2\pi $ we have $\theta>\pi $. Now by (\ref{eqIV.4.108}),
$$
\overrightarrow{M_{\theta_2}M_\theta}=\v x-a\v\begin{pmatrix}\cos \theta-\cos
\theta_2\\ \sin \theta-\sin \theta_2\end{pmatrix}\,.
$$
Since $\cos \theta$ is increasing for $\theta\in[\pi ,2\pi ]$ we have $\cos
\theta-\cos \theta_2<0$. Moreover, $\big\v\frac{\sin \theta-\sin
\theta_2}{\cos \theta-\cos
\theta_2}\big\v=\big\v\cotg\,\frac{\theta+\theta_2}2\big\v\leq 1$ since
$\frac{\theta+\theta_2}2\geq \frac{3\pi }2$, $\theta\leq 2\pi -\frac\pi 4$,
$\theta_2\leq 2\pi -\frac\pi 4$ so $\frac{\theta+\theta_2}2\leq 2\pi -\frac\pi
4$. It follows from Lemma \ref{lIV.4.21} (ii) that $M_\theta\in\CD^c$ which
proves (\ref{eqIV.4.109}) in case 2.

We conclude that if $x\in\CC\cup(\CD\setminus\CH)$ then $x=a+\v
x-a\v\begin{pmatrix}\cos \theta\\ \sin \theta\end{pmatrix}$ with
$\theta\notin ]\theta_1,\theta_2[$ and there exists a path joining the point
$x$ to the point
$a+\v x-a\v \begin{pmatrix}1\\ 0\end{pmatrix}$ with lenght less than $2\pi \,\v
x-a\v$.

\noindent {\bf Construction of the path in dimension 3}

We have $\alpha_\xi=(-1,0,0)$, $a=(a_1,a_2,0)$, $x=(x_1,x_2,x_3)$ and $-a\leq
D_0\sqrt{1+\v a\v^2}$, $-x_1\leq D_0\sqrt{1+\v x\v^2}$ with
$D_0=\frac{c_0}{10}$. We first construct a path in the plane $y_3=x_3$. We set
\begin{equation*}
\begin{aligned}
\CD&=\Big\{(y_1,y_2,x_3) : -y_1\leq D_0\sqrt{1+\v x_3\v^2+y^2_1+y^2_2}\Big\}\\
&=\Big\{(y_1,y_2,x_3) : -y_1\leq \frac{D_0}{\sqrt{1-D_0^2}}\,\sqrt{1+\v
x_3\v^2+y^2_2}\Big\}\,.
\end{aligned}
\end{equation*}
Since $c_0$ is small enough we have $\frac{D_0}{\sqrt{1-D^2_0}}<1$.
By the same way we see that the point $a$ is such that $-a_1\leq
\frac{D_0}{\sqrt{1-D^2_0}}\,\sqrt{1+a^2_2}$. Therefore $\tilde
a=(a_1,a_2,x_3)\in\CD$. Since $x\in\CD$, by the construction made in two
dimensions there exists a path lying in the set $\big\{y : y_3=x_3,\,\v
y-\tilde a\v=\v x-\tilde a\v=\sqrt{(x_1-a_1)^2+(x_2-a_2)^2}\big\}$ joining $x$
to $z=\big (a_1+\sqrt{(x_1-a_1)^2+(x_2-a_2)^2},a_2,x_3\big)$ of lenght smaller
than
$2\pi \,\v x-\tilde a\v\leq 2\pi \,\v x-a\v$.

Let us construct now a path in the plane $y_2=a_2$. Let us set
\begin{equation*}
\begin{aligned}
\CD&=\Big\{y=(y_1,a_2,y_3) : -y_1\leq D_0\,\sqrt{1+y^2_1+a^2_2+y^2_3}\Big\}\\
&=\Big\{y=(y_1,a_2,y_3) : -y_1\leq
\frac{D_0}{\sqrt{1-D_0^2}}\,\sqrt{1+a^2_2+y^2_3}\Big\}\,.
\end{aligned}
\end{equation*}
We have $z\in\CD$, $a\in\CD$. There exists a path joining $z$ to $(a_1+\v
x-a\v,\,a_2,0)$ lying in the set $\{y=(y_1,y_2,y_3) : y_2=a_2,\,\v y-a\v=\v
x-a\v\}$ with lenght smaller than $2\pi \,\v x-a\v$.

Now to join $0$ to $x$ we join $0$ to $z_1=(a_1+\v a\v,a_2,0)$ then $x$ to
$z_2=(a_1+\v x-a\v,a_2,0)$ and since the segment $[z_1,z_2]$ is included in
$\CD$ by Lemma \ref{lIV.4.21} (i), the path joins $0$ to $x$ and its lenght is
smaller than $C(\v x-a\v+\v a\v)$. Now by (\ref{eqIV.4.78}) we have $\v
a\v\leq \v x-a\v$ or $\v x-a\v\leq \v a\v\leq 2\,\v x-a\v$, so  the lenght of
the path is smaller than $C\,\v x-a\v=c\,\v x-x(\theta,\alpha)\v$. Moreover
$\v\gamma(\sigma,x)-a\v\leq 2\,\v x-a\v\leq \frac\delta{20}\,\bra\theta\ket$
so (\ref{eqIV.4.77}) is satisfied. Finally (\ref{eqIV.4.78}) and
(\ref{eqIV.4.79}) are obviously satisfied.

This ends the proof of Proposition \ref{pIV.4.14}.


\subsection{The phase for small $\theta$}\label{ssIV.5}

We shall need the following precision on the phase when $\v\theta\v\leq 1$.
\begin{theorem}\sl\label{tIV.5.1}
Let $\varphi $ be the phase given by Theorem \ref{tIV.1.2}. Then one can find
positive constants such that for $\v\theta\v\leq 1$, $\v
x-x(\theta,\alpha)\v\leq \delta\,\bra\theta\ket$ and $\v\alpha_\xi\v\leq 2$ one
can write
$$
\varphi (\theta,x,\alpha)=\frac{(x-\alpha_x)\cdot
\alpha_\xi-\theta\,\v\alpha_\xi\v^2+\frac i 2\,\v x-\alpha_x\v^2}{1+2i\theta}+
\frac 1 {2i}\,\v \alpha_\xi\v^2+R(\theta,x,\alpha)\,;
$$
where
\begin{equation*}
\begin{aligned}
\Big\v\frac{\partial R}{\partial \alpha_x}\Big\v&\leq
C\,(\varepsilon+\delta)(\v x-\alpha_x\v^2+\v\theta\v)\,,\enskip
\Big\v\frac{\partial R}{\partial \alpha_\xi}\Big\v\leq
C\,(\varepsilon+\delta)\,\v\theta\v\,,\\
\Big\v\frac{\partial^2 R}{\partial \alpha^2_x}\Big\v&\leq
C\,(\varepsilon+\delta)(\v x-\alpha_x\v^2+\v\theta\v)\,,\enskip
\Big\v\frac{\partial^2 R}{\partial \alpha_x\,\partial \alpha_\xi}\Big\v\leq
C\,(\varepsilon+\delta)\,\v\theta\v\,,\\
\Big\v\frac{\partial^2 R}{\partial \alpha^2_\xi}\Big\v&\leq
C\,(\varepsilon+\delta)\,\v\theta\v\,,
\end{aligned}
\end{equation*}
and
$$
\v\partial ^{A_1}_{\alpha_x}\,\partial
^{A_2}_{\alpha_\xi}\,R(\theta,x,\alpha)\v\leq 
\begin{cases}
C_{A_1}&\textrm{if \enskip  $A_2=0$}\\
C_{A_1,A_2}\,\v\theta\v &\textrm{if \enskip $\v A_2\v\geq 1$}\,.
\end{cases}
$$
\end{theorem}

{\bf Proof } Let us introduce the following space of functions.
\begin{equation}\label{eqIV.5.1}\left\{
\begin{array}{l}
\CE=\big\{ Z\in C^\infty (\R\times\R^n\times\R^n) : \v\partial ^\ell
_t\,\partial ^{A_1}_x\,\partial ^{A_2}_\xi\,Z(t,x,\xi)\v\leq C_{\ell
,A_1,A_2}\,\varepsilon\,\v t\v^{1-\ell }\,, \textrm{ for all}\\
 A_j\in\N^n\,,\enskip \ell =0,1\,,\enskip \v t\v\leq
1\,,\enskip x\in\R^n\,,\enskip \xi\in\R^n \textrm{ with } \v\xi\v\leq 2
\textrm{ and } Z(0,x,\xi)=0\big\}
\end{array}\right.
\end{equation}
Let us also recall that Proposition \ref{pIII.2.1} gives the following
description of the flow for $\v t\v\leq 1$, $x\in\R^n$, $\xi\in\R^n$ with
$\v\xi\v\leq 3$.
\begin{equation}\label{eqIV.5.2}\left\{
\begin{array}{l}
x(t,x,\xi)=x+2t\,\xi+r(t,x,\xi)\,,\\
\xi(t,x,\xi)=\xi+\zeta(t,x,\xi)\,,\\
z,\zeta\in\CE\,.
\end{array}\right.
\end{equation}
It follows that, with $f=x$ or $\xi$, we have
\begin{equation}\label{eqIV.5.3}
\big\v\partial ^\ell _t\,\partial ^{A_1}_x\,\partial
^{A_2}_\xi\,f(t,x,\xi)\big\v\leq C_{\ell ,A_1,A_2} \enskip \textrm{if}\enskip
\ell +\v A_1\v+\v A_2\v\geq 1\,.
\end{equation}
Let us set now
$$
g_j(\eta)=\chi_0\big(\frac\eta{\mu_0}\big)\big[(\xi_j-i\,x_j)(-\theta,y+x(\theta,\alpha),\eta+\xi(\theta,\alpha))
-(\alpha^j_\xi-i\,\alpha^j_x)\big]
$$
where $\chi_0\in C^\infty _0(\R^n)$, $\chi_0(\eta)=1$ if $\v\eta\v\leq \frac
12$, $\chi_0(\eta)=0$ if $\v\eta\v\geq 1$ and $\v y\v\leq \delta$. Setting
$x=y+x(\theta,\alpha)$, $\xi=\eta+\xi(\theta,\alpha)$ and using
(\ref{eqIV.5.2}) we obtain
\begin{equation}\label{eqIV.5.4}
\begin{split}
g_j(\eta)=\chi_0(\eta)\big[(1+2i\,\theta)\,\eta_j-i\,y_j+(1+2i\,\theta)\,\zeta_j
(\theta,\alpha)-i\,r_j(\theta,\alpha)\\
+(1+2i\theta)\,\zeta_j(-\theta,x,\xi)-i\,z_j(-\theta,x,\xi)\big]\,.
\end{split}
\end{equation}
We claim that we have the following estimates for $\ell =0,1$,
\begin{equation}\label{eqIV.5.5}
\v \partial ^\ell _\theta\, \partial ^\gamma_y\,\partial ^{\mu}_\eta\,\partial
^A_\alpha\,g_j\v\leq 
\begin{cases}
C_{\ell ,\gamma,\mu}&\textrm{if \enskip  $A=0$}\,,\\
C_{\ell ,\gamma,\mu,A}\,\varepsilon\,\v\theta\v^{1-\ell } &\textrm{if \enskip
$\v A\v\geq 1$}\,.
\end{cases}
\end{equation}
These estimates are obvious for the four first terms of $g_j$. So we are left
with the estimate of
$$
(1)=\partial ^\ell _\theta\,\partial ^\gamma_y\,\partial ^\mu_\eta\,\partial
^A_\alpha\big[Z(-\theta,y+x(\theta,\alpha),\eta+\xi(\theta,\alpha))\big]\,,
\quad Z\in \CE\,.
$$
To handle this term we shall make use of the Faa di Bruno formula given in
Appendix VIII.1, with $F=Z$, $Y=(\theta,y,\eta,\alpha)$, $U_1(Y)=-\theta$,
$U_{1+j}(Y)=y_j+x_j(\theta,\alpha)$,
$U_{1+n+j}(Y)=\eta_j+\xi_j(\theta,\alpha)$, $j=1,\ldots ,n$. Since $Z\in\CE$
we find easily, using (\ref{eqIV.5.3}) that $(1)\leq
C\,\varepsilon\,\v\theta\v^{1-\ell }$ which proves our claim.

Another property of $g_j$ which will be used in the sequel is the following.
\begin{equation}\label{eqIV.5.6}
\textrm{For } \theta=0\,,\enskip
g_j(\eta)=\chi_0(\eta)(\eta_j-i\,y_j)\enskip \textrm{is independent of}\enskip
\alpha\,.
\end{equation}
Now according to our procedure we have solved the equations (see
(\ref{eqIV.3.14})),
\begin{equation}\label{eqIV.5.7}
0=r(a,b,g_j)=g_j(-a)-i\,\som^n_{k=1} \frac{\partial g_j}{\partial
\eta_k}\,(-a)\,b_k+\som^n_{p,q=1} H^j_{pq}(\theta,y,\alpha,a,b)\,b_p\,b_q 
\end{equation}
in the set
$$
E=\Big\{(a,b)\in\R^n\times\R^n : \Big\v
a+\frac{2\theta\,y}{1+4\theta^2}\Big\v\leq \sqrt \delta\, \frac{\v
y\v}{\bra\theta\ket}\,,\enskip \Big\v b+\frac y{1+4\,\theta^2}\Big\v\leq
\sqrt\delta\,\frac{\v y\v}{\bra\theta\ket^2}\Big\}\,.
$$
Let us recall that we have the following bounds on $H^j_{pq}$ (see
(\ref{eqIV.3.13}) and (\ref{eqIV.3.16}))
$$
\big\v\partial ^\ell _\theta\,\partial ^\gamma_y\,\partial ^A_\alpha\,\partial
^M_{(a,b)}\,H^j_{p,q}\big\v\leq \som_{\v\mu\v\leq \v M\v+3n+2} 
\int\big\v\partial ^\ell _\theta\,\partial ^\gamma_y\,\partial
^A_\alpha\,\partial ^\mu_\eta\,g_j(\eta)\big\v\,d\eta\,,\quad \ell =0,1\,.
$$
Here we have used the fact that $r(a,b,g_j)$ is linear with respect to $g_j$.
It follows from (\ref{eqIV.5.5}), since $g_j$ has compact support in $\eta$,
that
\begin{equation}\label{eqIV.5.8}
\v \partial ^\ell _\theta\, \partial ^\gamma_y\,\partial ^A_\alpha\,\partial
^{M}_{(a,b)}\,H^j_{pq}\v\leq 
\begin{cases}
C_{\ell ,\gamma,M}&\textrm{if \enskip  $A=0$}\,,\\
C_{\ell ,\gamma,A,M}\,\varepsilon\,\v\theta\v^{1-\ell } &\textrm{if \enskip
$\v A\v\geq 1$}\,.
\end{cases}
\end{equation}
Using (\ref{eqIV.5.4}) we see easily that the equations (\ref{eqIV.5.7}) are
equivalent to the following system
\begin{equation}\label{eqIV.5.9}\left\{
\begin{array}{l}
a_j=\frac{-2\theta\,y_j}{1+4\theta^2}+Z^a_{j,1}(\theta,\alpha)+Z^a_{j,2}(-\theta,y
+x(\theta,\alpha),-a+\xi(\theta,\alpha))\\
\hbox to 1,88cm{}
+Z^a_{j,3}(-\theta,y+x(\theta,\alpha),-a+\xi(\theta,\alpha))\,b+H^a_j
(\theta,y,\alpha,a,b)\,b\cdot b\\
b_j=\frac{-y_j}{1+4\theta^2}+Z^b_{j,1}(\theta,\alpha)+Z^b_{j,2}(-\theta,y
+x(\theta,\alpha),-a+\xi(\theta,\alpha))\\
\hbox to 1,88cm{}
+Z^b_{j,3}(-\theta,y+x(\theta,\alpha),-a+\xi(\theta,\alpha))\,b+H^b_j
(\theta,y,\alpha,a,b)\,b\cdot b
\end{array}\right.
\end{equation}
where the $Z's$ belong to the space $\CE$ defined in (\ref{eqIV.5.1}) and the
$H'_js$ satisfy the estimates (\ref{eqIV.5.8}).

According to (\ref{eqIV.5.6}), (\ref{eqIV.4.2}), Theorem \ref{tIV.4.2} and
Theorem \ref{tIV.3.1} for
$\theta=0$ $a_j$ and $b_j$ do not depend on $\alpha$ and moreover we have,
\begin{equation}\label{eqIV.5.10}\left\{
\begin{array}{l}
a_j(0,y,\alpha)=\tilde a_j(y)=\CO(\v y\v^N)\,,\\
b_j(0,y,\alpha)=\tilde b_j(y)=-y_j+\CO (\v y\v^N)\,,
\end{array}\right.
\end{equation}
for every $N\in\N$ and $\v y\v\leq \delta$.

Let us set
\begin{equation}\label{eqIV.5.11}
G^*_j(\theta,y,\alpha)=H^*_j(\theta,y,\alpha,a(\theta,y,\alpha),b(\theta,y,\alpha))
\cdot b(\theta,y,\alpha)\cdot b(\theta,y,\alpha),*=a\textrm{ or }b\,.
\end{equation}
Then, since the $Z'_js$ vanish for $\theta=0$, (\ref{eqIV.5.10}) implies that
\begin{equation}\label{eqIV.5.12}
G^*_j(0,y,\alpha)=\tilde G^*_j(y)=\CO (\v y\v^N)\,,\quad \forall \,N\in\N\,.
\end{equation}
Therefore we can write
\begin{equation}\label{eqIV.5.13}
G^*_j(\theta,y,\alpha)=G^*_j(y)+\int^\theta_0 \frac d{d\theta}\,
G^*_j(\sigma,y,\alpha)\,d\sigma\,.
\end{equation}
We claim that we have the following estimates on $a_j$, $b_j$. Let us set for
convenience $f_j=a_j$ or $b_j$.
\begin{equation}\label{eqIV.5.14}\left\{
\begin{array}{l}
\Big\v\frac{\partial f_j}{\partial \theta}\,(\theta,y,\alpha)\Big\v\leq
C\,(\varepsilon+\delta)\\
\v \partial ^\ell _\theta\, \partial ^\gamma_y\,\partial ^A_\alpha\,
f_j(\theta,y,\alpha)\v\leq 
\begin{cases}
C_{\ell ,\gamma}&\textrm{if \enskip  $A=0$}\,,\\
C_{\ell ,\gamma,A}\,(\varepsilon+\delta)\,\v\theta\v^{1-\ell } &\textrm{if
\enskip $\v A\v\geq 1,\enskip \ell =1$}\,.
\end{cases}
\end{array}\right.
\end{equation}
To prove the first estimate we differentiate both sides of (\ref{eqIV.5.9})
with respect to $\theta$. Since the terms in $Z$ belong to $\CE$ and using
(\ref{eqIV.5.8}), the fact that $\dot x(\theta,\alpha),\dot\xi(\theta,\alpha)$
are bounded, we obtain
$$
\Big\v\frac{\partial a_j}{\partial \theta}\Big\v+\Big\v\frac{\partial
b_j}{\partial \theta}\Big\v\leq C_1\,\v
y_j\v+C_2\,\varepsilon\,+C_3 
(\varepsilon+\delta)\Big(\Big\v\frac{\partial a}{\partial
\theta}\Big\v+\Big\v\frac{\partial b}{\partial \theta}\Big\v\Big)\,.
$$
Taking $\varepsilon+\delta$ small enough and since $\v y\v\leq \delta$ we
obtain our first claim. To prove the second estimate we use the Faa di Bruno
formula (see Appendix VIII.1) and an induction procedure.

Let us set $Y=(\theta,y,\alpha),\Lambda=(\ell ,\gamma,A)$ and let us apply the
operator $\partial ^\Lambda_Y$ to both sides of (\ref{eqIV.5.9}). We have
\begin{equation*}
\begin{array}{l}
\partial ^\Lambda_Y\,f_0=
\begin{cases}
0(1)&\textrm{if \enskip  $A=0$}\,,\\
0 &\textrm{if
\enskip $\v A\v\geq 1$}\,.
\end{cases}\,,\enskip f_0=\frac{-2\theta\,y_j}{1+4\theta^2}\quad 
\textrm{or}\quad  \frac{-y_j}{1+4\theta^2}\\
\v\partial ^\Lambda_Y\,Z(\theta,\alpha)\leq C\, \varepsilon\,\v
\theta\v^{1-\ell }\,.
\end{array}
\end{equation*}
Assume now that our estimate is true for $\v\Lambda\v\leq k$ and let
$\v\Lambda\v=k+1$. Then
$$
\partial
^\Lambda_Y\big[Z(-\theta,y+x(\theta,\alpha),-a+\xi(\theta,\alpha))
\,b(\theta,y,\alpha)\big]=(1)+(2)+(3)
$$
where
\begin{equation*}
\begin{aligned}
(1)&=Z(-\theta,\cdots )\,\partial ^\Lambda_Y\,b(Y)\\
(2)&=\partial ^\Lambda_Y
\big[Z(-\theta,y+x(\theta,\alpha),-a(Y)+\xi(\theta,\alpha))\big]\,b(Y)\\
(3)&=\som_{\Lambda_1+\Lambda_2=\Lambda\atop \Lambda_j\not =0} \begin{pmatrix}
\Lambda\\\Lambda_1\end{pmatrix}\,\partial ^{\Lambda_1}_ Y [Z(-\theta,\cdots
)]\,\partial ^{\Lambda_2}_Y\, b\,.
\end{aligned}
\end{equation*}
Using the Faa di Bruno formula in the terms (2) and (3) we see that
\begin{equation*}
\v\partial ^\Lambda_Y\,[Z(-\theta,\cdots)]\v\leq C\,\varepsilon\,\v\theta\v
(\v\partial ^\Lambda_Y\,a\v+\v\partial ^\Lambda_Y\,b\v)+
\begin{cases}
0(1)&\textrm{if \enskip  $A=0$}\,,\\
0((\varepsilon+\delta)\v\theta\v^{1-\ell }) &\textrm{if \enskip $\v A\v\not =
0$}\,.
\end{cases}
\end{equation*}
By the same way, using (\ref{eqIV.5.8}) we obtain
\begin{equation*}
\begin{split}
\v\partial ^\Lambda_Y[H(Y,a(Y),\,b(Y))&\cdot b(Y)\cdot b(Y)\v]\leq
C\,(\varepsilon+\delta)\v\theta\v(\v\partial ^\Lambda_Ya(Y)\v\\
&+\v\partial
^\Lambda_Y\,b(Y)\v)+\begin{cases}
0(1)&\textrm{if \enskip  $A=0$}\,,\\
0((\varepsilon+\delta)\v\theta\v^{1-\ell }) &\textrm{if \enskip $\v A\v\geq 
1$}\,.
\end{cases}
\end{split}
\end{equation*}
Taking $\varepsilon+\delta$ small enough we obtain the second estimate of
(\ref{eqIV.5.14}).

Now it follows from (\ref{eqIV.5.9}), (\ref{eqIV.5.11}), (\ref{eqIV.5.12}) and
(\ref{eqIV.5.13}) that we can write, with $Y=(\theta,y,\alpha)$,
\begin{equation}\label{eqIV.5.16}\left\{
\begin{array}{l}
(a_j+i\,b_j)(Y)=-\frac{y_j}{2\theta-i}+U_j(Y)\enskip \textrm{where}\\
U_j(Y)=Z_{j,1}(\theta,\alpha)+Z_{j,2}(-\theta,y+x(\theta,\alpha),-a(Y)+\xi(\theta,\alpha))\\
+Z_{j,3} (-\theta,y+x(\theta,\alpha),-a(Y)+\xi(\theta,\alpha))\cdot b(Y)\\
+\tilde G_j(y)+\int^\theta_0\Big\{\Big[\frac{\partial H}{\partial
\theta}+\frac{\partial H}{\partial a}\cdot \frac{\partial a}{\partial
\theta}+\frac{\partial H}{\partial b}\,\frac{\partial b}{\partial
\theta}\Big]\,b\cdot b+2H\,b\cdot \frac{\partial b}{\partial
\theta}\Big\}(\sigma,y,\alpha)\,d\sigma\\
\textrm{with}\enskip \v\partial ^\gamma_y\,\tilde G_j(y)\v\leq C_N\,\v
y\v^N\enskip \textrm{for all}\enskip N\in\N\enskip \textrm{and}\enskip
\gamma\in\N^n\,.
 \end{array}\right.
\end{equation}
Using (\ref{eqIV.5.1}), (\ref{eqIV.5.3}), (\ref{eqIV.5.8}), (\ref{eqIV.5.14})
we deduce the following estimates
\begin{equation}\label{eqIV.5.16}\left\{
\begin{array}{l}
\v U(Y)\v+\big\v\frac{\partial U}{\partial y}\,(Y)\big\v\leq C_N\,\v
y\v^N+C\,(\varepsilon+\delta)\,\v\theta\v\leq C'\,(\varepsilon+\delta)\\
\v\partial ^\gamma_y\,U(Y)\v\leq C_{N,\gamma}\,\v
y\v^N+C_\gamma\,\v\theta\v\enskip \textrm{if}\enskip \v\gamma\v\geq 2\\
\v\partial ^\gamma_y\,\partial  ^A_\alpha\,U(Y)\v\leq
C_{\gamma,A}\,(\varepsilon+\delta)\,\v\theta\v\enskip \textrm{if}\enskip
\v\gamma\v\geq 0\enskip \textrm{and}\enskip \v A\v\geq 1\,.
\end{array}\right.
\end{equation}
Now recall that,
\begin{equation}\label{eqIV.5.17}\left\{
\begin{array}{l}
\Phi_j(\theta,x,\alpha)=\xi_j(\theta,\alpha)-(a_j+i\,b_j)(\theta,x-
x(\theta,\alpha),\alpha)\\
\varphi (\theta,x,\alpha)=\int^1_0 (x-x(\theta,\alpha))\cdot
\Phi(\theta,s\,x+(1-s)\,x(\theta,\alpha),\alpha)\,ds+\theta\,p(\alpha)+\frac
1{2i}\,\v\alpha_\xi\v^2\,.
\end{array}\right.
\end{equation}
It follows that,
\begin{equation*}
\begin{split}
\varphi (\theta,x,\alpha)&=\underbrace{(x-x(\theta,\alpha))}_{(1)}\cdot
\xi(\theta,\alpha)+\underbrace{\frac 12\, \frac{\v
x-x(\theta,\alpha)\v^2}{2\theta-i}}_{(2)}\\
&-\int^1_0 (x-x(\theta,\alpha))\cdot
U(\theta,s(x-x(\theta,\alpha)),\alpha)\,ds+\theta p(\alpha)+\frac
1{2i}\,\v\alpha_\xi\v^2\,.
\end{split}
\end{equation*}
We have,
\begin{equation}\label{eqIV.5.18}\left\{
\begin{array}{l}
(1)=(x-\alpha_x)\alpha_\xi-2\theta\,\v\alpha_\xi\v^2-r(\theta,\alpha)\cdot
\alpha_\xi+(x-x(\theta,\alpha))\cdot \zeta(\theta,\alpha)\\
(2)=\frac 12\,\frac 1{2\theta-i} \big[(x-\alpha_x)^2-4\theta(x-\alpha_x)\cdot
\alpha_\xi+4\theta^2\,\v\alpha_\xi\v^2\\
\hbox to 5cm{}+2(x-\alpha_x-2\theta\,\alpha_\xi)
\cdot r(\theta,\alpha)+\v r(\theta,\alpha)\v^2\big]\,.
\end{array}\right.
\end{equation}
Let us consider the term in $\varphi $ which does not contain any error term.
It can be written as
$$
(x-\alpha_x)\cdot
\alpha_\xi-2\theta\,\v\alpha_\xi\v^2+\frac{(x-\alpha_x)^2-4\theta(x-\alpha_x)\cdot
\alpha_\xi+4\theta^2\,\v\alpha_\xi\v^2}{2(2\theta-i)}+\theta\,\v\alpha_\xi\v^2+\frac
1{2i}\,\v\alpha_\xi\v^2
$$
which is equal to
$$
\frac{(x-\alpha_x)\cdot \alpha_\xi+\frac
i2\,(x-\alpha_x)^2-\theta\,\v\alpha_\xi\v^2}{1+2i\,\theta}+\frac
1{2i}\,\v\alpha_\xi\v^2\,.
$$
It follows that
$$
\varphi (\theta,x,\alpha)=\frac{(x-\alpha_x)\cdot \alpha_\xi+\frac
i2\,(x-\alpha_x)^2-\theta\,\v\alpha_\xi\v^2}{1+2i\theta}+\frac
1{2i}\,\v\alpha_\xi\v^2+R(\theta,x,\alpha)
$$
where
\begin{equation*}
\begin{split}
R(\theta,x,\alpha)=\underbrace{-r(\theta,\alpha)\cdot
\alpha_\xi}_{(1)}+\underbrace{(x-x(\theta,\alpha))\cdot
\zeta(\theta,\alpha)}_{(2)}\\
\underbrace{\frac{(x-\alpha_x-2\theta\,\alpha_\xi)\cdot r(\theta,\alpha)+\frac
12\,\v
r(\theta,\alpha)\v^2}{2\theta-i}}_{(3)}+\underbrace{\theta(p(\alpha)
-\v\alpha_\xi\v^2)}_{(4)}\\
-\int^1_0\underbrace{(x-x(\theta,\alpha))\cdot
U(\theta,s(x-x(\theta,\alpha)),\alpha)\,ds}_{(5)}\,.
\end{split}
\end{equation*}
We are ready now to show that the remainder $R$ satisfies the estimates given
in Theorem \ref{tIV.5.1}.

First of all, since $r$ and $\zeta$ belong to $\CE$ (see (\ref{eqIV.5.1})),
since the functions $x-x(\theta,\alpha)$ and $x-\alpha_x-2\theta\,\alpha_\xi$
are bounded together with all their derivatives with respect to $\alpha$ and
since $p(\alpha)=\v\alpha_\xi\v^2+\varepsilon \som^n_{j,k=1}
b_{jk}(\alpha_x)\,\alpha^j_\xi\,\alpha^k_\xi$ we have
\begin{equation}\label{eqIV.5.19}
\v\partial ^A_\alpha\,(i)\v\leq C_A\,\varepsilon\,\v\theta\v\enskip
\textrm{for}\enskip i=1,2,3,4\,.
\end{equation}
So we are left with the term (5). Let us note that if we set
$f_0(\theta,x,\alpha)=x-x(\theta,\alpha)$ then,
\begin{equation}\label{eqIV.5.20}\left\{
\begin{array}{l}
\v f_0\v\leq 2\delta\,,\enskip  \v f_0\v\leq \v x-\alpha_x\v+C\,\v\theta\v\,,
\enskip 
\big\v\frac{\partial f_0}{\partial \alpha_x}\big\v\leq C\,,\enskip
\big\v\frac{\partial f_0}{\partial \alpha_\xi}\big\v\leq C\,\v\theta\v\\
\v\partial ^A_\alpha\,f_0\v\leq C\,\varepsilon\,\v\theta\v\enskip
\textrm{if}\enskip \v A\v\geq 2,\enskip \textrm{uniformly in}\enskip
(x,\theta,\alpha)\,.
\end{array}\right.
\end{equation}
With this notation one has
$$
(5)=f_0(\theta,x,\alpha) \int^1_0
U(\theta,s\,f_0(\theta,x,\alpha),\alpha)\,ds\,.
$$
Then, with $i=x$ or $\xi$,
$$
\frac{\partial }{\partial \alpha_i}\,(5)=\frac{\partial f_0}{\partial
\alpha_i}\int^1_0 U(\theta,s\,f_0,\alpha)\,ds+f_0 \int^1_0 \Big (s\,
\frac{\partial f_0}{\partial \alpha_i}\, \frac{\partial U}{\partial
y}+\frac{\partial U}{\partial \alpha_i}\Big )(\theta,s\,f_0,\alpha)\,ds\,.
$$
Now it follows from (\ref{eqIV.5.18}) and (\ref{eqIV.5.20}) that
\begin{equation}\label{eqIV.5.21}\left\{
\begin{array}{l}
\big\v\frac\partial {\partial \alpha_x}\,(5)\big\v\leq
C\,(\varepsilon+\delta)(\v x-\alpha_x\v^2+\v\theta\v)\,,\\
\big\v\frac\partial {\partial \alpha_\xi}\,(5)\big\v\leq
C\,(\varepsilon+\delta)\,\v\theta\v)\,,
\end{array}\right.
\end{equation}
since $\v\frac{\partial f_0}{\partial \alpha_\xi}\,(5)\v\leq C\,\v\theta\v$
and $\v x-\alpha_x\v\leq \delta$.
Now, with $i,j=x$ or $\xi$, we have
\begin{equation*}
\begin{split}
\frac{\partial ^2}{\partial \alpha_i\,\partial \alpha_j}\,(5)=\frac{\partial
^2\,f_0}{\partial \alpha_i\,\partial \alpha_j} \int^1_0
U(\theta,s\,f_0,\alpha)\,ds+\frac{\partial f_0}{\partial \alpha_i} \int^1_0
\Big (s\,\frac{\partial f_0}{\partial \alpha_j}\,\frac{\partial U}{\partial y}+
\frac{\partial U}{\partial \alpha_j}\Big)\,ds\\
+\frac{\partial f_0}{\partial \alpha_j} \int ^1_0 \Big (s\,\frac{\partial f_0}
{\partial \alpha_i}\,\frac{\partial U}{\partial y}+
\frac{\partial U}{\partial \alpha_i}\Big)\,ds+
f_0 \int^1_0\Big
[s\,\frac{\partial^2 f_0}{\partial \alpha_i\,\partial \alpha_j}
\,\frac{\partial U}{\partial y}+ s\frac{\partial f_0}{\partial \alpha_i}
\Big(\frac{\partial f_0}{\partial \alpha_j}\,\frac{\partial ^2U}{\partial
y^2}+\frac{\partial U}{\partial \alpha_j}\Big )\\
+\frac{\partial ^2U}{\partial
\alpha_i\,\partial y}\,\frac{\partial f_0}{\partial \alpha_j}+\frac{\partial
^2U}{\partial \alpha_i\,\partial \alpha_j}\Big]\,ds\,.
\end{split}
\end{equation*}
Using (\ref{eqIV.5.16}) and (\ref{eqIV.5.20}) we check easily that
\begin{equation}\label{eqIV.5.22}\left\{
\begin{array}{l}
\Big\v\frac{\partial ^2}{\partial \alpha^2_x}\,(5)\Big\v\leq
C\,(\varepsilon+\delta)(\v x-\alpha_x\v^2+\v\theta\v)\\
\Big\v\frac{\partial ^2}{\partial \alpha_x\,\partial \alpha_\xi}\,(5)\Big\v
\leq C\,(\varepsilon+\delta)\v\theta\v\\
\Big\v\frac{\partial ^2}{\partial  \alpha^2_\xi}\,(5)\Big\v
\leq C\,(\varepsilon+\delta)\v\theta\v\,.
\end{array}\right.
\end{equation}
Combining (\ref{eqIV.5.19}), (\ref{eqIV.5.21}) and (\ref{eqIV.5.22}) we obtain
the claimed estimates on the two first derivatives of $R$.

Finally using again (\ref{eqIV.5.18}) and (\ref{eqIV.5.20}) we deduce the
following estimates of the higher derivatives
\begin{equation*}
\v\partial ^{A_1}_{\alpha_x}\,\partial ^{A_2}_{\alpha_\xi}\,(5)\v\leq 
\begin{cases}
C_{A_1}&\textrm{if \enskip  $A_2=0$}\,,\\
C_{A_1\,A_2}\,\v\theta\v &\textrm{if \enskip $\v A_2\v\geq 1$}\,.
\end{cases}
\end{equation*}
The gain of $\v\theta\v$ when $\v A_2\v\geq 1$ coming from the fact that a
derivative of $x-x(\theta,\alpha)$ and $U$ with respect to $\alpha_\xi$ makes
appear a $\theta$. Thus the proof of Theorem \ref{tIV.5.1} is complete.
\cqfd

\section{The transport equations}\label{sV}
\subsection{Statement of the result and preliminaries}\label{ssV.1}

Let $P=\som^n_{j,k=1} g^{jk}(x)\,D_j\,D_k$ be a
second order differential operator of the form (\ref{eqII.2.7}) We shall denote
by
${^t}P$ the transposed operator.

In Section \ref{sIV} we have constructed a phase function for $P$. The purpose
now is to construct an amplitude.

Recall that the set $\Omega_\delta$ have been
introduced in Definition \ref{dIV.1.1}.

The main result of this Section is the following.
\begin{theorem}\sl \label{tV.1.1}
For every $\alpha\in T^*\R^n$ with $\frac 12\leq \vert \alpha_\xi\vert
\leq 2$, every $N\in\N$ and every $\lambda\geq 1$ one can find an amplitude
$e_N(\theta,y,\alpha,\lambda)$ which is $C^\infty $ on $\tilde \Omega_\delta$
such that
\begin{itemize}
\item[(i)] $e_N(0,y,\alpha,\lambda)=1$.
\item[(ii)] $\big (i\lambda\,\frac \partial{\partial \theta}-i\lambda\,\frac n
2\,\frac{\theta}{1+\theta^2}-{}^tP\big )\big (e^{i\lambda\varphi
(\theta,x,\alpha)}\,e_N(\theta,x-x(\theta,\alpha),\alpha,\lambda)\big)=R_N(\theta
,x,\alpha,\lambda)\,e^{i\lambda\varphi (\theta,x,\alpha)}$
where
\begin{equation*}
\vert R_N(\theta,x,\alpha,\lambda)\vert \leq C_N\Big
(\lambda^{-N}+\lambda^2\,\frac{\vert x-x(\theta,\alpha)\vert
^N}{\bra\theta\ket^N}\Big )
\end{equation*}
for every $(\theta,x,\alpha)$ in $\Omega_\delta$, $\lambda\geq 1$ and
$C_N$ is independent of $(\theta,x,\alpha,\lambda)$.
\item[(iii)] $\vert \partial
^A_x\,e_N(\theta,x-x(\theta,\alpha),\alpha,\lambda)\vert \leq C_{N,A}$
uniformly with respect to $(\theta,x,\alpha,\lambda)$.
\end{itemize}
\end{theorem}
\begin{corollary} \sl \label{cV.1.2}
For every $\alpha\in T^*\R^n$ with $\frac 12\leq \vert \alpha_\xi\vert
\leq 2$, every $N\in\N$ and every $\lambda\geq 1$ one can find an amplitude
$a_N(\theta,x,\alpha,\lambda)$ which is $C^\infty $ on $\Omega_\delta$ such that
\begin{itemize}
\item[(i)] $a_N(0,x,\alpha,\lambda)=1$,
\item[(ii)] $\big (i\lambda\, \frac \partial {\partial \theta}-{^t}P\big)\big
(e^{i\lambda\varphi (\theta,x,\alpha)}\,a_N(\theta,x,\alpha,\lambda)\big
)=R'_N(\theta,x,\alpha,\lambda)\,e^{i\lambda\varphi (\theta,x,\alpha)}$
where
\begin{equation*}
\vert R'_N(\theta,x,\alpha,\lambda)\vert \leq C'_N\Big
(\lambda^{-N}+\lambda^2\,\frac{\vert x-x(\theta,\alpha)\vert
^N}{\bra\theta\ket^N}\Big )
\end{equation*}
uniformly with respect to $(\theta,x,\alpha,\lambda)$.
\item[(iii)] $\vert \partial^A_x\,a_N(\theta,x,\alpha,\lambda)\vert \leq
C_{N,A}\,\bra \theta\ket^{-\frac n2}$, uniformly with
respect to
$(\theta,x,\alpha,\lambda)$.
\end{itemize}
\end{corollary}

\noindent {\bf Proof } We have just to set
$a_N(\theta,x,\alpha,\lambda)=\bra\theta\ket^{-\frac
n2}\,e_N(\theta,x-x(\theta,\alpha),\alpha,\lambda)$ where $e_N$ has been
defined in Theorem \ref{tV.1.1}.
\cqfd

\noindent {\bf Proof of Theorem \ref{tV.1.1}} We have $^tP=\som^n_{j,k=1}
g^{jk}(x)\,D_j\,D_k+\som^n_{j=1} g_j(x)\,D_j+g_0(x)$, where
$g^{jk}=\delta_{jk}+\varepsilon\,b_{jk}$ and $b_{jk}\in \BB^1_{\sigma_0}$,
$g_j\in \BB^2_{\sigma_0}$, $1\leq j,k\leq n$, $g_0\in\BB^3_{\sigma_0}$ where,
\begin{equation}\label{eqV.1.1}
\BB^\ell _{\sigma_0}=\Big\{g\in C^\infty (\R^n):\vert \partial ^A_x\,g(x)\vert
\leq
\frac{C_A}{\bra x\ket^{\vert A\vert +\ell +\sigma_0}}\,,\enskip \forall
x\in\R^n\,,\enskip \forall A\in\N^n\Big\}\,.
\end{equation}
A straightforward computation shows that
\begin{equation}\label{eqV.1.2}
\begin{split}
\Big (i\lambda\,\frac\partial {\partial \theta}&-i\lambda\,\frac
n2\,\frac\theta{1+\theta^2}-{}^tP\Big )(e^{i\lambda\varphi
}f)=e^{i\lambda\varphi }\bigg[-\lambda^2\Big (\frac{\partial \varphi }{\partial
\theta}+p\Big (x,\,\frac{\partial \varphi }{\partial x}\Big)\Big)\,f+i\lambda
\Big (\frac{\partial f}{\partial \theta}+2 \sum^n_{j,k=1}
g^{jk}\,\frac{\partial \varphi }{\partial x_j}\,\frac{\partial f}{\partial
x_k}\\
&+\sum^n_{j,k=1} g^{jk}\, \frac{\partial ^2\varphi }{\partial x_j\,\partial
x_k}\,f-\frac n2\,\frac \theta{1+\theta^2}\,f+i \sum^n_{j=1} g_j \,\frac
{\partial \varphi }{\partial x_j}\,f\Big )-{^t}Pf\bigg]\,.
\end{split}
\end{equation}
According to Theorem \ref{tIV.1.2} the coefficient of $\lambda^2$ in the right
hand side of (\ref{eqV.1.2}) is bounded by $C_N\big (\frac{\vert
x-x(\theta,\alpha)\vert }{\bra\theta\ket}\big)^N$, for any $N$. Therefore if we
set
\begin{equation}\label{eqV.1.3}\left\{
\begin{array}{l}
I=e^{-i\lambda\varphi }\Big (i\lambda\,\frac\partial {\partial
\theta}-i\lambda\,\frac n2\,\frac \theta{1+\theta^2}-{}^tP\Big
)\big(e^{i\lambda\varphi }f\big)\\
X=\frac\partial {\partial \theta}+2\sum^n_{j,k=1}g^{jk}\,\frac{\partial \varphi
}{\partial x_j}\cdot \frac\partial {\partial x_k}
\end{array}\right.
\end{equation}
we obtain
\begin{equation}\label{eqV.1.4}
\begin{split}
\bigg\vert I-i\lambda\,Xf-i\lambda\Big (\sum^n_{j,k=1} g^{jk}\,\frac{\partial
^2\varphi }{\partial x_j\,\partial x_k}-\frac
n2\,\frac{\theta}{1+\theta^2}+i\sum^n_{j=1} g_j\,\frac{\partial \varphi
}{\partial x_j}\Big )\,f-{}^tPf\bigg\vert \\
\leq C_N\,\lambda^2\,\Big (\frac{\vert x-x(\theta,\alpha)\vert
}{\bra\theta\ket}\Big)^N\,.
\end{split}
\end{equation}
To pursue the proof we consider separately the two cases.

\subsection{The case of outgoing points}\label{ssV.2}

For convenience we shall set
\begin{equation*}
\begin{aligned}
\tilde \CS_+&=\big\{\alpha\in T^*\R^n : \frac 12\leq \v\alpha_\xi\v\leq
2\,, \alpha_x\cdot \alpha_\xi\geq
-c_0\,\bra\alpha_x\ket\,\v\alpha_\xi\v\big\}\\
\tilde \CS_-&=\big\{\alpha\in T^*\R^n : \frac 12\leq \v\alpha_\xi\v\leq
2\,, \alpha_x\cdot \alpha_\xi\leq
c_0\,\bra\alpha_x\ket\,\v\alpha_\xi\v\big\}\,.
\end{aligned}
\end{equation*}

Here we assume $\alpha\in\tilde \CS_\pm$.
Let us set $R_j=\frac{\partial \varphi }{\partial x_j}-\Phi_j$. It follows
from Theorem \ref{tIV.3.13} (v) and (\ref{eqIV.3.48}) that $\v\partial
^A_x\,R_j(\theta,x,\alpha)\v\leq \frac{C_A}{\bra\theta\ket^{\v A\v}}$ for
$A\in\N^n$ and from Proposition \ref{pIV.3.19} (ii) that $\v
R_j(\theta,x,\alpha)\v\leq C_N \big(\frac{\v
x-x(\theta,\alpha)\v}{\bra\theta\ket}\big)^N$ for all $N\in\N$.

First of all, according to (\ref{eqIV.3.43}) and Taylor's formula, we have
\begin{equation}\label{eqV.2.1}\left\{
\begin{array}{l}
\frac{\partial \varphi }{\partial
x_j}\,(\theta,x,\alpha)=\xi_j(\theta,\alpha)+\frac{\frac 12 \sgn \theta
(1-\chi_1(\theta))
(x_j-x_j(\theta,\alpha))-(a_j+i\,b_j)(\theta,x-x(\theta,\alpha),
\alpha)} {\bra\theta\ket}+R_j(\theta,x,\alpha)\,,\\
\vert \partial ^A_x\,R_j(\theta,x,\alpha)\vert \leq C_{A,N} \big (\frac{\vert
x-x(\theta,\alpha)\vert }{\bra\theta\ket}\big)^N\,.
\end{array}\right.
\end{equation}
Using Theorem \ref{tIV.3.1} (iii) we deduce,
\begin{equation}\label{eqV.2.2}\left\{
\begin{array}{l}
\frac{\partial ^2\varphi }{\partial x_j\,\partial x_k}\,(\theta,x,\alpha)=\frac
12\, \frac{\sgn\theta}{\bra\theta\ket}\,(1-\chi_1(\theta))\,\delta_{jk}+d_{jk}
(\theta,x-x(\theta,\alpha),\alpha)+R_{jk}(\theta,x,\alpha)\, ,\\
\vert \partial ^A_x\,d_{jk}(\theta,x-x(\theta,\alpha),\alpha)\vert \leq
\frac{C_A}{\bra\theta\ket^{\vert A\vert +2}}\,,\\
\vert \partial ^A_x\,R_{jk}(\theta,x,\alpha)\vert \leq C_{AN}\big(\frac{\vert
x-x(\theta,\alpha)\vert}{\bra\theta\ket}\big)^N
\end{array}\right.
\end{equation}
uniformly with respect to $(\theta,x,\alpha)$.

Now if $g\in\BB_{\sigma_0}$ (see (\ref{eqV.1.1})) and
$(\theta,x)\in\Omega_\delta$ we can write $\vert \partial ^A_x\,g(x)\vert
=\vert (\partial ^A_xg)(y+x(\theta,\alpha))\vert \leq \frac{C_A}{\bra
y+x(\theta,\alpha)\ket^{\vert A\vert +1+\sigma_0}}$. Since $\vert y\vert \leq
\delta\bra\theta\ket$ we can use Proposition \ref{pIII.3.1} to write
$$
\bra y+x(\theta,\alpha)\ket\geq \bra x(\theta,\alpha)\ket-\vert y\vert \geq
\frac 1{\sqrt 3}\, \bra\theta\ket-\delta\bra\theta\ket\geq \frac 12 \,\bra
\theta\ket\,.
$$
It follows that
\begin{equation}\label{eqV.2.3}
\vert \partial ^A_x\,g(x)\vert \leq \frac{C_A}{\bra \theta\ket^{\vert A\vert
+1+\sigma_0}}\,.
\end{equation}
This can be applied to the functions $g^{jk}-\delta_{jk}$, $g_j$, $g_0$, $1\leq
j,k\leq n$. It follows from (\ref{eqV.2.2}) that
$$
\sum^n_{j,k=1} g^{jk}(x)\,\frac{\partial ^2\varphi }{\partial x_j\,\partial
x_k}\,(\theta,x,\alpha)=\frac n2\,\frac{\sgn \theta}{\bra
\theta\ket}(1-\chi_1(\theta))+d(\theta,x-x(\theta,\alpha),\alpha)+R(\theta,x,\alpha)
$$
where $d$ and $R$ satisfy the same estimates as in (\ref{eqV.2.2}). Now we
have
\begin{equation}\label{eqV.2.4}
\begin{aligned}
\frac 12\,\frac{\sgn \theta}{\bra\theta\ket}(1-\chi_1(\theta))-\frac
12\,\frac{\theta}{1+\theta^2}&=\frac 12\,\frac{\sgn
\theta}{\bra\theta\ket}\,(1-\chi_1(\theta))\Big (1-\frac{\vert \theta\vert
}{\bra\theta\ket}\Big
)-\frac{\theta\,\chi_1(\theta)}{2\bra\theta\ket^2}\\
&=\frac
12\,\frac{\sgn\theta}{\bra\theta\ket^2}\,\frac{1-\chi_1(\theta)}{\bra\theta\ket+\vert
\theta\vert }=\CO\Big (\frac{1}{\bra\theta\ket^3}\Big)\,.
\end{aligned}
\end{equation}
Summing up, we have proved
\begin{equation}\label{eqV.2.5}\left\{
\begin{array}{l}
\som^n_{j,k=1} g^{jk}(x)\,\frac{\partial ^2\varphi }{\partial x_j\,\partial
x_k}\,(\theta,x,\alpha)-\frac n2
\,\frac{\theta}{1+\theta^2}=D_1(\theta,x-x(\theta,\alpha),\alpha)
+R_1(\theta,x,\alpha)\,,\\
\som^n_{j=1} g_{j}(x)\,\frac{\partial\varphi }{\partial
x_j}\,(\theta,x,\alpha)
=D_2(\theta,x-x(\theta,\alpha),\alpha)
+R_2(\theta,x,\alpha)\,,\\
\vert \partial ^A_x\,D_j(\theta,x-x(\theta,\alpha),\alpha)\vert \leq
\frac{C_A}{\bra\theta\ket^{\vert \alpha\vert +1+\sigma_0}}\,,\quad j=1,2\,,\\
\vert R_j(\theta,x,\alpha)\vert \leq C_N \Big (\frac{\vert
x-x(\theta,\alpha)\vert }{\bra\theta\ket}\Big )^N\,, \quad j=1,2\,,\enskip
\forall N\in\N.\,
\end{array}\right.
\end{equation}
We are going now to simplify the vector fields $X$ introduced in
(\ref{eqV.1.3}). Let us set
\begin{equation}\label{eqV.2.6}\left\{
\begin{array}{l}
s=\theta\\
y=x-x(\theta,\alpha)\,.
\end{array}\right.
\end{equation}
Since $\dot x_k(\theta,\alpha)=\frac{\partial p}{\partial
\xi_k}\,(x(\theta,\alpha),\xi(\theta,\alpha))=2\som^n_{j=1}
g^{jk}(x(\theta,\alpha))\xi_k(\theta,\alpha)$, we obtain
$$
X=\frac{\partial }{\partial s}-2 \sum^n_{j,k=1}
\Big\{g^{jk}(x(s,\alpha))\,\xi_j(s,\alpha)-g^{jk}(y+x(s,\alpha))\,\frac{\partial
\varphi }{\partial x_j}\,(s,y+x(s,\alpha),\alpha)\Big\}\,\frac{\partial
}{\partial y_k}\,.
$$
Now using (\ref{eqV.2.1}) and $g^{jk}=\delta_{jk}+\varepsilon\,b_{jk}$,
$b_{jk}\in\BB_{\sigma_0}$ we can write
\begin{equation*}
\begin{split}
X=\frac\partial {\partial s}&+\frac{\sgn s}{\bra s\ket}(1-\chi_1(s))
\sum^n_{j=1} y_j\,\frac\partial {\partial y_j}-2\varepsilon \sum^n_{j,k=1}
\big\{b_{jk}(x(s,\alpha))-b_{jk}(y+x(s,\alpha))\big\}\,\xi_j(s,\alpha)\,\frac\partial 
{\partial y_k}\\
&\quad -2\sum^n_{j=1} \frac{(a_j+i\,b_j)(s,y,\alpha)}{\bra s\ket}\cdot
\frac\partial {\partial y_j}+2\varepsilon \sum^n_{j,k=1}
b_{jk}(y+x(s,\alpha))\,\frac 1{\bra s\ket }\Big (\frac 12\,\sgn s\,y_j\\
&\quad -(a_j+i\,b_j)(s,y,\alpha)\Big
)\,\frac\partial {\partial y_k}+2 \sum^n_{j,k=1}
g^{jk}(y+x(s,\alpha))\,R_j(s,y+x(s,\alpha),\alpha)\,\frac\partial {\partial
y_k}\,.
\end{split}
\end{equation*}
\begin{definition}\sl \label{dV.2.1}
We shall say that a function $f=f(s,y,\alpha)$ on $\tilde
\Omega_\delta\times\tilde \CS_\pm$ belongs to $\CE$ if
\begin{equation}\label{eqV.2.7}\left\{
\begin{array}{l}
f(s,0,\alpha)=0\\
\vert \partial ^A_y\,f(s,y,\alpha)\vert \leq \frac{C_A}{\bra s\ket^{\vert
A\vert +1}}\,,\quad A\in\N^n
\end{array}\right.
\end{equation}
uniformly when $(s,y)\in\tilde \Omega_\delta$ and $\alpha\in\tilde \CS_\pm$.
\end{definition}

According to (\ref{eqV.2.4}), (\ref{eqV.2.3}) and Theorem \ref{tIV.3.1}, (ii),
(iii) we have, 
\begin{equation}\label{eqV.2.8}\left\{
\begin{array}{l}
(1-\chi_1(s))\,\frac{\sgn s}{\bra s\ket}\,y_j-\frac s{1+s^2}\,y_j\in\CE\\
\varepsilon(b_{jk}(x(s,\alpha)-b_{jk}(y+x(s,\alpha))))\,\xi_j(s,\alpha)\in\CE\\
\varepsilon\,b_{jk}(y+x(s,\alpha))\,\frac 1{\bra s\ket} \Big (\frac 12\cdot 
\sgn s\cdot y_j-(a_j+i\,b_j)(s,y,\alpha)\Big )\in\CE\\
\frac 1{\bra s\ket}\,(a_j+i\,b_j)(s,y,\alpha)\in\CE\,.
\end{array}\right.
\end{equation}
Then we have,
\begin{equation}\label{eqV.2.9}\left\{
\begin{array}{l}
X=\frac\partial {\partial s}+\frac s{1+s^2} \som^n_{j=1} y_j\,\frac{\partial
}{\partial y_j}+\som^n_{j=1} E_j(s,y,\alpha)\,\frac\partial {\partial
y_j}+\som^n_{j=1} R'_j(s,y,\alpha)\,\frac\partial {\partial y_j}\\
\textrm{where } E_j\in\CE\textrm{ and } \vert R'_j(s,y,\alpha)\vert \leq C_N
\Big (\frac{\vert y\vert }{\bra s\ket}\Big)^N\,.
\end{array}\right.
\end{equation}
Now we perform another change of variables. We set
\begin{equation}\label{eqV.2.10}\left\{
\begin{array}{l}
\theta=s\\
z=\frac y{\bra s\ket}\,.
\end{array}\right.
\end{equation}
Then we have $\frac\partial {\partial \theta}=\frac\partial{\partial s}+\frac
s{1+s^2} \som^n_{j=1} y_j\,\frac\partial {\partial y_j}$. It follows that
\begin{equation}\label{eqV.2.11}\left\{
\begin{array}{ll}
\textrm{(i)}&\quad X=\frac\partial {\partial \theta}+\som^n_{j=1}
h_j(\theta,z,\alpha)\,\frac\partial {\partial z_j}+\som^n_{j=1}\tilde
R_j(\theta,z,\alpha)\,\frac\partial {\partial z_j}\,,\\
\textrm{(ii)}&\quad h_j(\theta,0,\alpha)=0\,,\\
\textrm{(iii)}&\quad \vert \partial ^A_z\,h_j(\theta,z,\alpha)\vert \leq
\frac{C_A}{\bra\theta\ket^2}\,,\quad A\in\N^n\,,\\
\textrm{(iv)}&\quad \vert \tilde R_j(\theta,z,\alpha)\vert \leq
\frac{C_N}{\bra\theta\ket}\,\vert z\vert ^N\,,\textrm{ for all } N\textrm{ in
}\N\,,\textrm{ uniformly when}\\
&\quad \theta\geq 0\textrm{ (resp. }\theta\leq
0),\enskip \vert z\vert \leq \delta,\enskip \alpha\in\Lambda_{\pm},\enskip
j=1,\ldots ,n\,.
\end{array}\right.
\end{equation}
Moreover, since $\frac\partial {\partial y_j}=\frac
1{\bra\theta\ket}\,\frac\partial {\partial z_j}$ we have by (\ref{eqV.2.3}),
\begin{equation}\label{eqV.2.12}\left\{
\begin{array}{l}
^tP=\som_{\vert \nu\vert \leq 2} k_\nu(\theta,z,\alpha)\,\partial ^\nu_z\,,\\
\vert \partial ^\gamma_z\,k_\nu(\theta,z,\alpha)\vert \leq
\frac{C_\gamma}{\bra\theta\ket^{1+\sigma_0}}\,, \quad \gamma\in\N^n\,.
\end{array}\right.
\end{equation}
Let us set
\begin{equation}\label{eqV.2.13}
X_0=\frac\partial {\partial \theta}+\sum^n_{j=1}
h_j(\theta,z,\alpha)\,\frac\partial {\partial z_j}\,,
\end{equation}
where $h_j$ satisfies (\ref{eqV.2.11}).

It follows from (\ref{eqV.1.4}), (\ref{eqV.2.5}), (\ref{eqV.2.6}),
(\ref{eqV.2.10}) and (\ref{eqV.2.11}) (i) that
\begin{equation}\label{eqV.2.14}\left\{
\begin{array}{l}
\Big\vert I-i\lambda\Big (X_0\,f+d(\theta,z,\alpha)\,f-\frac i\lambda
{}^tPf\Big)\Big\vert \leq C_N\,\lambda^2\,\vert z\vert ^N\,(\vert f\vert +\vert
\nabla_zf\vert )\\
\vert \partial ^\gamma_z\,d(\theta,z,\alpha)\vert \leq
\frac{C_\gamma}{\bra\theta\ket^{1+\sigma_0}}\,,\quad \gamma\in\N^n~.
\end{array}\right.
\end{equation} 
Now let us fix an integer $N_0$ large enough depending only on the dimension
$n$ (and chosen later on). For the coefficients $h_j$, $k_\nu$, $d$ in
(\ref{eqV.2.13}), (\ref{eqV.2.12}) and (\ref{eqV.2.14}) we write
\begin{equation}\label{eqV.2.15}\left\{
\begin{array}{l}
F(\theta,z,\alpha)=F^{N_0}(\theta,z,\alpha)+r^{N_0}(\theta,z,\alpha)\textrm{
where,}\\
F^{N_0}(\theta,z,\alpha)=\som_{\vert \gamma\vert \leq N_0-1}\partial
^\gamma_z\,F(\theta,0,\alpha)\,\frac{z^\gamma}{\gamma!}\,,\\
\vert r^{N_0}(\theta,z,\alpha)\vert \leq C_{N_0}\,\vert z\vert ^{N_0}\,.
\end{array}\right.
\end{equation}
Let us set
\begin{equation}\label{eqV.2.16}\left\{
\begin{array}{l}
L=\frac\partial {\partial \theta}+\som^n_{j=1}
h_j^{N_0}(\theta,z,\alpha)\,\frac\partial {\partial
z_j}+d^{N_0}(\theta,z,\alpha)\,,\\
Q=\som_{\vert \nu\vert \leq 2} k^{N_0}_\nu(\theta,z,\alpha)\,\partial^\nu_z\,.
\end{array}\right.
\end{equation}
Using (\ref{eqV.2.14}) to (\ref{eqV.2.16}) we see that
\begin{equation}\label{eqV.2.17}
\Big\vert I-i\lambda\Big (L\,f-\frac i\lambda\,Q\,f\Big)\Big\vert \leq
C_{N_0}\,\lambda^2\,\vert z\vert ^{N_0}\sum_{\vert \gamma\vert \leq 2} \vert
\partial ^\gamma_z\,f(\theta,z,\alpha)\vert \,.
\end{equation}
Now we have the following result.
\begin{lemma}\sl \label{lV.2.2}
There exist functions $A_\ell =A_\ell (\theta,z,\alpha)$, $\ell =0,\ldots
,N_0+1$ which are $C^\infty $ in $(\theta,z)$ in the set
$\CO=\{(\theta,z):\theta\in\R^{\pm},\,\vert z\vert \leq \delta\}$ such that
\begin{itemize}
\item[(i)] $A_0(0,z,\alpha)=1$, $A_\ell (0,z,\alpha)=0$, $\ell =1,\ldots
,N_0+1$,
\item[(ii)] $\vert \partial ^\gamma_z\,A_\ell (\theta,z,\alpha)\vert \leq
C_{\ell ,\gamma}$, uniformly in $\CO\times\tilde \CS_{\pm}$, $(\ell \in\N$,
$\gamma\in\N^n)$,
\item[(iii)] $L\,A_0=0$, $L\,A_\ell =i\,Q\,A_{\ell -1}$, $\ell =1,\ldots
,N_0+1$.
\end{itemize}
\end{lemma}

Let us assume for a moment this lemma proved. Let us set
\begin{equation}\label{eqV.2.18}
f=f_{N_0}=A_0+\frac 1\lambda\,A_1+\ldots +\frac
1{\lambda^{N_0+1}}\,A_{N_0+1}\,.
\end{equation}
Then Lemma \ref{lV.2.2} shows that
\begin{equation}\label{eqV.2.19}\left\{
\begin{array}{l}
f_{N_0}(0,z,\alpha,\lambda)=1\\
\vert \partial ^\gamma_z\,f_{N_0}(\theta,z,\alpha,\lambda)\vert \leq
C_{\gamma,N_0}\textrm{ if }(\theta,z)\in\CO\,,\enskip
\alpha\in\tilde \CS_{\pm}\,,\enskip \lambda\geq 1\\
\Big\vert L\,f_{N_0}-\frac i\lambda\,Q\,f_{N_0}\Big\vert
\leq \lambda^{-N_0-2}\,\vert Q\,A_{N+1}\vert \leq
C_{N_0}\,\lambda^{-N_0-2}\,.
\end{array}\right.
\end{equation}
It follows from (\ref{eqV.2.17}) and (\ref{eqV.2.19}) that
\begin{equation}\label{eqV.2.20}
\vert I\vert \leq C'_{N_0}\,\lambda^2\,\vert z\vert
^{N_0}+C''_{N_0}\,\lambda^{-N_0-1}\,.
\end{equation}
Coming back to the variables $(\theta,x)$ we set
$$
e_{N_0}(\theta,x-x(\theta,\alpha),\alpha,\lambda)=f_{N_0}\Big
(\theta,\,\frac{x-x(\theta,\alpha)}{\bra\theta\ket}\,,\alpha,\lambda\Big)\,.
$$
Then it follows from (\ref{eqV.2.20}), (\ref{eqV.2.19}) that $e_{N_0}$
satisfies the conditions (i), (ii), (iii) in Theorem \ref{tV.1.1}.

So we are left with the proof of Lemma \ref{lV.2.2}.

\noindent {\bf Proof of Lemma \ref{lV.2.2} } We are going to straighten the
principal part of the operator $L$ given by (\ref{eqV.2.16}). Recall that we
have $L=L_0 +d^{N_0}(\theta,z,\alpha)$ with 
$$
L_0=\frac\partial {\partial \theta}+\sum^n_{j=1}
h^{N_0}_j(\theta,z,\alpha)\,\frac\partial {\partial z_j}\,.
$$  
Moreover according to (\ref{eqV.2.11}) and (\ref{eqV.2.15}) we have
\begin{equation}\label{eqV.2.21}\left\{
\begin{array}{ll}
\textrm{(i)}&\quad h^{N_0}_j(\theta,z,\alpha)=\som^n_{k=1} \frac{\partial
h_j}{\partial z_k}\,(\theta,0,\alpha)\,z_k+g_j(\theta,z,\alpha)\,,\\
\textrm{(ii)}&\quad g_j(\theta,z,\alpha)=\som_{2\leq \vert \gamma\vert \leq
N_0-1}\frac 1{\gamma!}\,\partial ^\gamma_z\,h_j(\theta,0,\alpha)\,z^\gamma\,,\\
\textrm{(iii)}&\quad \som^n_{j=1} \vert \partial
^\gamma_z\,h_j(\theta,0,\alpha)\vert \leq
\frac{C_\gamma}{\bra\theta\ket^2}\,,\quad \forall \gamma\in\N^n\,.
\end{array}\right.
\end{equation}
In that follows all the objects will depend on $\alpha\in\tilde \CS_{\pm}$ but
all the estimates will be uniform with respect to $\alpha$.

Let us set
\begin{equation}\label{eqV.2.22}
H(\theta)=\Big (\frac{\partial h_j}{\partial z_k}\,(\theta,0,\alpha)\Big
)_{1\leq j,k\leq n}\,.
\end{equation}
If $\theta_0\in\R^{\pm}$ we shall denote by $Y(\theta,\theta_0)$ the unique
$n\times n$ matrix solution of the problem
\begin{equation}\label{eqV.2.23}\left\{
\begin{array}{l}
\dot Y(\theta,\theta_0)=H(\theta)\,Y(\theta,\theta_0)\,,\enskip
\theta\in\R^\pm\,,\\
Y(\theta_0,\theta_0)=\Id\,.
\end{array}\right.
\end{equation}
Since by (\ref{eqV.2.21}) (iii) the entries of the matrix $H(\theta)$ are
bounded by $\frac C{\bra\theta\ket^2}$, the Gronwall inequality shows that there
exists $M_0\geq 1$ such that
\begin{equation}\label{eqV.2.24}
\Vert Y(\theta,\theta_0)\Vert \leq M_0\,, \textrm{ for all }
\theta,\theta_0\in\R^\pm\textrm{ and } \alpha\in\tilde \CS_\pm\,.
\end{equation}
Moreover since $Y(\theta,\theta_0)^{-1}=Y(\theta_0,\theta)$ we have also,
\begin{equation}\label{eqV.2.25}
\Vert Y(\theta,\theta_0)^{-1}\Vert \leq M_0\,, \textrm{ for all }
\theta,\theta_0\in\R^\pm\textrm{ and } \alpha\in\tilde \CS_\pm\,.
\end{equation}
Now using (\ref{eqV.2.21}) we see that the problem
\begin{equation}\label{eqV.2.26}\left\{
\begin{array}{l}
\dot z_j(\theta)=h^{N_0}_j(\theta,z(\theta),\alpha)\,,\enskip
\theta\in\R^\pm\,,\enskip 1\leq j\leq n\,,\\
z_j(0)=y_j\,,
\end{array}\right.
\end{equation}
is equivalent, setting $z=(z_j)_{1\leq j\leq n}$, $g=(g_j)_{1\leq j\leq n}$, to
\begin{equation}\label{eqV.2.27}
z(\theta)=Y(\theta,0)\,y+\int^\theta_0 Y(\theta,t)\,g(t,z(t),\alpha)\,dt\,.
\end{equation}
Then we have the following Lemma.
\begin{lemma}\sl \label{lV.2.3}
One can find $\eta>0$ such that for all $y\in\C^n$ such that $\vert y\vert
\leq
\eta$, the problem (\ref{eqV.2.27}) has a unique global solution $z$ such that
$\vert z(\theta)\vert \leq 2M_0\eta$ for all $\theta\in\R^\pm$. This solution
will be denoted by $z(\theta,y)$. Moreover one can find a constant
$C(N_0,M_0)\geq 0$ such that
\begin{itemize}
\item[(i)] $\big\Vert \big (\frac{\partial z_j}{\partial y_k}\big
)(\theta,y)-Y(\theta,0)\big\Vert \leq C(M_0,N_0)\,\eta$,

and for every $\gamma\in\N^n$, one can find a constant $C_\gamma\geq 0$ such
that
\item[(ii)] $\vert \partial ^\gamma_y\,z(\theta,y)\vert \leq C_\gamma $, for
all $\theta\geq 0$ and $\vert y\vert \leq \eta$.
\end{itemize}
\end{lemma}

\noindent {\bf Proof } Let $\eta>0$ (to be chosen). Assume $\vert y\vert \leq
\eta$ and set $A=\{T>0$ such that (\ref{eqV.2.27}) has a solution for
$\theta\in[0,T]$ satisfying $\vert z(\theta)\vert \leq 2M_0\,\eta\}$. Since
(\ref{eqV.2.27}) (which is equivalent to (\ref{eqV.2.26})) has a
continuous solution for small $\theta$ and since $\vert z(0)\vert \vert =\vert
y\vert \leq \eta$ there exists $\varepsilon_0>0$ such that $\varepsilon_0\in
A$. Thus $A$ is non empty and it is obviously an interval. Let $T^*=\sup A$.
If $T^*=+\infty $ we are done so assume $T^*<+\infty $. Let us take
$T\in]0,T^*[$. Then on $[0,T]$ we have
$$
\vert z(\theta)\vert \leq \vert Y(\theta,0)\,y\vert +\int^\theta_0 \Vert
Y(\theta,t)\Vert \,\vert g(t,z(t),\alpha)\vert \,dt\,.
$$
Now by (\ref{eqV.2.21}) (ii) and (iii) we have
$$
\vert g(t,z(t),\alpha)\vert \leq \frac{K_{N_0}}{\bra
t\ket^2}\,2M_0\,\eta\,\vert z(t)\vert 
$$
if $2M_0\,\eta\leq 1$, where $K_{N_0}=\som_{2\leq \vert \gamma\vert \leq
N_0-1}\frac{C_\gamma}{\gamma \kern 1pt_ !}$.

It follows from (\ref{eqV.2.24}) that,
$$
\vert z(\theta)\vert \leq M_0\,\eta+\int^\theta_0
\frac{2M^2_0\,K_{N_0}\,\eta}{\bra t\ket^2}\,\vert z(t)\vert \,dt\,.
$$
Then the Gronwall inequality implies that
$$
\vert z(\theta)\vert \leq M_0\,\eta\,\exp \Big[2M^2_0\,K_{N_0}\,\eta
\int^{+\infty }_0 \frac{dt}{\bra t\ket^2}\Big]\,.
$$
Therefore taking $\eta$ small (compared to $M_0$ and $K_{N_0}$) we can achieve
that $\vert z(\theta)\vert \leq \frac 32\,M_0\,\eta$ for all $\theta\in[0,T]$.
A classical argument shows that $z(T^*)$ can be defined and $\vert z(T^*)\vert
\leq \frac 32\, M_0\,\eta$. Then solving again (\ref{eqV.2.26}) with data
$z(T^*)$ we see easily that this contradicts the definition of $T^*$ as the
supremum of $A$. Therefore $T^*=+\infty $. 

Now differentiating (\ref{eqV.2.27})
 with respect to $y_k$ yields
\begin{equation}\label{eqV.2.28}
\frac{\partial z}{\partial y_k}\,(\theta,y)-Y(\theta,0)\,e_k+\int^\theta_0
Y(\theta,t)\sum^n_{\ell =1} \frac{\partial g}{\partial z_\ell
}\,(t,z(t),\alpha)\,\frac {\partial z_\ell }{\partial y_k}\,(t,y)\,dt\,.
\end{equation}
First of all (\ref{eqV.2.21}) (ii) show that
\begin{equation}\label{eqV.2.29}
\sum^n_{\ell =1} \Big\vert \frac{\partial g}{\partial z_\ell
}\,(t,z(t),\alpha)\Big\vert \leq \sum_{2\leq \vert \gamma\vert \leq N_0-1}
\frac{C'_\gamma}{\bra t\ket^2}\, \vert z(t,y)\vert ^{\vert \gamma\vert -1}\leq
\frac{C_{N_0}\,M_0\,\eta}{\bra t\ket^2}
\end{equation}
if $2M_0\,\eta\leq 1$. It follows that
$$
\Big\vert \frac{\partial z}{\partial y_k}\,(\theta,y)\Big\vert \leq
M_0+\int^\theta_0 \frac{C'_{N_0}\,M^2_0\,\eta}{\bra t\ket^2} \Big\vert
\frac{\partial z}{\partial y_k}\,(t,y)\Big\vert \,dt\,.
$$
The Gronwall inequality shows that one can find $K=K(N_0,M_0)$ such that
$\big\vert \frac{\partial z}{\partial y_k}\,(\theta,y)\vert \leq K$, for all
$\theta\in\R^\pm$ and $\vert y\vert \leq \eta$.

Using again (\ref{eqV.2.28}) and (\ref{eqV.2.29}) we see that
$$
\Big\vert \frac{\partial z}{\partial y_k}\,(\theta,y)-Y(\theta,0)\,e_k\Big\vert
\leq \int^\theta_0 \frac{C_{N_0}\,M^2_0\cdot K\eta}{\bra t\ket^2}\,dt\leq
C(M_0,N_0)\,\eta\,.
$$
Finally the estimate on $\partial ^\gamma_y\,z$, which is true for $\vert
\gamma\vert =0,1$ by the above results, can be easily obtained by induction on
$\vert \gamma\vert $ using (\ref{eqV.2.27}), (\ref{eqV.2.21}) and the Gronwall
Lemma. \cqfd

In the sequel we shall take $\eta$ so small that $C(M_0,N_0)\,\eta\leq \frac
12$. 

Let us now consider the map
\begin{equation}\label{eqV.2.30}\left\{
\begin{array}{ll}
\Phi:\R^\pm \times\{y\in\C^n: \vert y\vert \leq \eta\}&\longrightarrow
\R^\pm\times\C^n\,,\\
\hbox to 3,6cm{}(\theta,y)&\longmapsto (\theta,z(\theta,y))\,.
\end{array}\right.
\end{equation}
We claim that $\Phi$ is injective. Indeed for a fixed $\theta\in\R^\pm$ if we
have $y_j$, $j=1,2$ such that $\vert y_j\vert \leq \eta$ and
$z(\theta,y_1)=z(\theta,y_2)$ then
$$
0=\sum^n_{k=1} \int^1_0 \frac{\partial z}{\partial
y_k}\,(\theta,t\,y_1+(1-t)\,y_2)(y^k_1-y^k_2)\,dt\,.
$$
Since $\vert t\,y_1+(1-t)\,y_2\vert \leq \eta$ when $t\in[0,1]$ we can use the
estimate given in Lemma \ref{lV.2.3} to ensure that
$$
\vert Y(\theta,0)(y_1-y_2)\vert \leq C'(N_0,M_0)\,\eta\,\vert y_1-y_2\vert
\,.
$$
According to (\ref{eqV.2.25}) this implies that $y_1=y_2$ if $\eta$ is small
enough. 

It follows that $\Phi$ is bijective on its range. We show now that
\begin{equation}\label{eqV.2.31}\left\{
\begin{array}{l}
\textrm{If } \delta\textrm{ is small enough we have}\\
\R^\pm\times\{z\in\C^n:\vert z\vert \leq
\delta\}\subset\Phi(\R^\pm\times\{y\in\C^n:\vert y\vert \leq \eta\})\,.
\end{array}\right.
\end{equation}
This equivalent to show that for fixed $\theta\in\R^\pm$,
\begin{equation}\label{eqV.2.32}\left\{
\begin{array}{l}
\textrm{for all }z\in\C^n\textrm{ such that }\vert z\vert \leq \delta
\textrm{ there exists }y\in\C^n\\
\textrm{such that }\vert y\vert \leq \eta\textrm{ and }z(\theta,y)=z\,.
\end{array}\right.
\end{equation}

According to (\ref{eqV.2.27})
 the equation to solve is equivalent to the equation $y=F(y)$ where
\begin{equation}\label{eqV.2.32}
F(y)=Y(\theta,0)^{-1}\,z-Y(\theta,0)^{-1} \int^\theta_0
Y(\theta,t)\,g(t,z(t,y),\alpha)\,dt\,.
\end{equation}
Let $B=\{y\in\C^n:\vert y\vert \leq \eta\}$. We shall show that if $\delta$ and
$\eta$ are small enough compared to $N_0,M_0$ then $F$ maps $B$ into $B$ and
there exists $\varepsilon_0<1$ such that $\vert F(y_1)-F(y_2)\vert \leq
\varepsilon_0\,\vert y_1-y_2\vert $ for all $y_1,y_2$ in $B$. Then
(\ref{eqV.2.32}) will follow from the fixed point Theorem. Since
$Y(\theta,0)^{-1}=Y(0,\theta)$ and $Y(0,\theta)\,Y(\theta,t)=Y(0,t)$ it follows
from (\ref{eqV.2.24}) and (\ref{eqV.2.21}) that if $\vert z\vert \leq \delta$
we have
$$
\vert F(y)\vert \leq M_0\,\delta+\Big\vert \int^\theta_0
\frac{C(M_0)\,\eta^2}{\bra t\ket^2}\,dt\Big \vert 
$$
since $\vert z(t,y)\vert \leq 2M_0\,\eta$ by Lemma \ref{lV.2.3}. Then $\vert
F(y)\vert \leq \eta$ if $\eta$ is small enough in terms of $M_0$ and
$M_0\,\delta\leq \frac 12\,\eta$.

Moreover if $y_1,y_2$ belong to $B$ we have
$$
\vert F(y_1)-F(y_2)\vert \leq \Big\vert \int^\theta_0 \frac{C(M_0,N_0)}{\bra
t\ket^2}\,\eta\,\vert z(t,y_1)-z(t,y_2)\vert \,dt\Big\vert \,.
$$
Since by Lemma \ref{lV.2.3} we have $\vert z(t,y_1)-z(t,y_2)\vert \leq
C'(M_0)\vert y_1-y_2\vert $ we obtain finally
$$
\vert F(y_1)-F(y_2)\vert \leq C'(M_0,N_0)\,\eta\,\vert y_1-y_2\vert \,.
$$
Taking $\eta$ small enough we obtain (\ref{eqV.2.32}).

We can now straighten the vector field $L_0$ which is the principal part of $L$
given in (\ref{eqV.2.16}). Let us make the change of variables,
$(\theta',y)\mapsto (\theta,z(\theta,y))$. Then we have, according to
(\ref{eqV.2.26}) 
$$
\frac{\partial }{\partial \theta'}=\frac\partial {\partial \theta}+\sum^n_{j=1}
\dot z_j(\theta,y)\,\frac\partial {\partial z_j}=\frac\partial {\partial
\theta}+\sum^n_{j=1} h^{N_0}_j(\theta,z(\theta,y))\,\frac\partial {\partial
z_j}=L_0\,.
$$
In the new coordinates $(\theta',y)$ the operator $L$ has therefore the form
$$
L=\frac\partial {\partial
\theta'}+d^{N_0}(\theta',z(\theta',y),\alpha)=\frac\partial {\partial
\theta'}+\tilde d(\theta',y,\alpha)\,.
$$
Now we note that
$$
\frac \partial {\partial \theta'} \Big (e^{\int^{\theta'}_0 \tilde
d(t,y,\alpha)\,dt}\,u(t,y,\alpha)\Big )=e^{\int^{\theta'}_0 \tilde
d(t,y,\alpha)\,dt}\,Lu(t,y,\alpha)\,.
$$
It follows that the problem
$$
L\,\tilde A_0=0\,,\quad \tilde A_0(0,y,\alpha)=1
$$
has the (unique) solution $\tilde A_0(\theta',y,\alpha)=e^{-\int^{\theta'}_0
\tilde d(t,y,\alpha)\,dt}$. By the same way the problems
$$
L\,\tilde A_\ell =i\,\tilde Q\,\tilde A_{\ell -1}\,,\quad \tilde A_\ell
(0,y,\alpha)=0\,,\quad \ell =1,\ldots ,N_0+1\,,
$$
are solved by
$$
\tilde A_\ell (\theta',y,\alpha)=e^{-\int^{\theta'}_0 \tilde
d(t,y,\alpha)\,dt} \int^{\theta'}_0 i(\tilde Q\,\tilde A_{\ell
-1})(t,y,\alpha)\,e^{\int^{t}_0 \tilde
d(s,y,\alpha)\,ds}\,dt\,.
$$
To end the proof of Lemma \ref{lV.2.2} we are left with the uniform estimates
(ii).

First of all, using the estimate in (\ref{eqV.2.14}), (\ref{eqV.2.15}), Lemma
\ref{lV.2.3} (ii) and the Faa di Bruno formula we see that,
\begin{equation}\label{eqV.2.34}
\vert \partial ^\gamma_y(d^{N_0}(\theta',z(\theta',y),\alpha)\vert \leq
\frac{C_\gamma}{\bra\theta'\ket^2}\,.
\end{equation}
Denoting by $\kappa(\theta,z)$ the inverse map of $y\mapsto z(\theta,y)$, that
is
$\kappa(\theta,z(\theta,y))=y$ and using Lemma \ref{lV.2.3} we see that,
\begin{equation}\label{eqV.2.35}
\vert \partial ^\gamma_z\,\kappa(\theta,z)\vert \leq C_\gamma\textrm{ for }
\theta\geq 0\textrm{ and }\vert z\vert \leq \delta\,.
\end{equation}
Then let us set for $\ell =0,\ldots ,N_0+1$
$$
A_\ell (\theta,z,\alpha)=\tilde A_\ell (\theta,\kappa(\theta,z),\alpha)\,.
$$
Using (\ref{eqV.2.34}), (\ref{eqV.2.35}), (\ref{eqV.2.16}), (\ref{eqV.2.15})
and the estimate in (\ref{eqV.2.12}) we see that $(A_\ell )_{\ell =0,\ldots
,N_0+1}$ satisfy all the requirements of Lemma \ref{lV.2.2}. This ends the
proof of Theorem \ref{tV.1.1} in the case of outgoing points.\cqfd

We consider now the case of incoming points.

\subsection{The case of incoming points}\label{ssV.3}

We assume here that $\alpha\in T^*\R^n$ is such that $\frac 12\leq \vert
\alpha_\varepsilon\vert \leq 2$ and
\begin{equation}\label{eqV.3.1}
\alpha_x\cdot \alpha_\xi\leq -c_0\,\bra\alpha_x\ket\,\vert \alpha_\xi\vert \,.
\end{equation}
Since such points belong to $\CS_-$ (see Definition \ref{dIII.2.2}) the case
where
$\theta\leq 0$ is covered by the Section \ref{ssV.2}. We focus now on the case
$\theta\geq 0$. Here the method used in Section \ref{ssV.2} does not work for
many technical reasons. For instance, when $\vert \alpha_x\vert $ is very
large, $\alpha_\xi=-\frac{\alpha_x}{\vert \alpha_x\vert }$ and $\theta=\frac
12\, \vert \alpha_x\vert $ we can see that $\bra x(\theta,\alpha)\ket$ is of
magnitude one. Therefore we are far from the estimate $\bra
x(\theta,\alpha)\ket\geq \frac 1{\sqrt 2}\,\bra\theta\ket$ used for instance to
get (\ref{eqV.2.3}). Here also we shall use the method which consists to
straighten the vector field $X$, defined in (\ref{eqV.1.3}). This is done by a
change of variables in $(\theta,x)$ deduced from the flow of $X$. The problem
here is that $X$ has non real coefficients (because of $\frac{\partial \varphi
}{\partial x_j}$) which are merely $C^\infty $. Therefore we are led to push
the problem in the complex domain by extending all the functions almost
analytically as in  \cite{MS} for instance. So we begin our Section by a
Lemma on almost analytic extensions adapted to our situation. In that follows
we shall consider together two cases.
Case 1 : $\Omega=\R^n_x$, case 2 : $\Omega=\Omega_\delta$ (see
Definition \ref{dIV.1.1}). We shall denote by $X$ the variable in $\Omega$ that
is $X=x$ in the first case, $X=(\theta,x)$ in the second one.
\begin{lemma}\sl \label{lV.3.1}
Let $f$ be a function defined on $\Omega$ which is $C^\infty$ in $X$ and
satisfies for all $X$ in $\Omega$, in case~1 (resp. case~2),
\begin{equation}\label{eqV.3.2}\left\{
\begin{array}{l}
\vert f(X)\vert \leq \frac{M_0}{\bra x\ket^{\sigma_1}}\enskip  \Big
(\textrm{resp. } M_0\Big (\frac 1{\bra x\ket^{\sigma_1}}+\frac 1{\bra
\theta\ket^{\sigma_2}}\Big ) \Big )\\
\som_{\vert \gamma\vert =k} \vert \partial ^\gamma_x\,f(X)\vert \leq
\frac{M_k}{\bra x\ket^{k+\sigma_3}} \enskip \Big (\textrm{resp. } M_k\Big
(\frac 1{\bra x\ket^{k+\sigma_3}}+\frac 1{\bra\theta\ket^{k+\sigma_3}}\Big
)\,,\quad k\geq 1
\end{array}\right.
\end{equation}
where $(M_k)_{k\geq 0}$ is an increasing sequence in $]0,+\infty [$ and $0\leq
\sigma_1\leq \sigma_3$, $0\leq \sigma_2\leq \sigma_3$. Then there exists
$F=F(X,y)$ defined on $\Omega\times\R^n_y$ which is $C^\infty $ in $(X,y)$ and
satisfies for all $(X,y)$ in $\Omega\times\R^n_y$,
\begin{itemize}
\item[(i)] $F(X,0)=f(X)$.
\item[(ii)] $\v F(X,y)\vert \leq \frac{C_0}{\bra x\ket^{\sigma_1}}$ (resp.
$C_0\big (\frac 1{\bra x\ket^{\sigma_1}}+\frac 1{\bra\theta\ket^{\sigma_2}}\big
)$).
\item[(iii)] For every $A,B$ in $\N^n$ with $\vert A\vert +\vert B\vert \geq 1$
there exists $C_{AB}>0$ such that 

$\vert \partial ^A_x\,\partial
^B_y\,F(X,y)\vert \leq \frac{C_{AB}}{\bra x\ket^{\vert A\vert +\vert B\vert
+\sigma_3}}$ (resp. $C_{AB}\big (\frac 1{\bra x\ket^{\vert A\vert +\vert B\vert
+\sigma_3}}+\frac 1{\bra \theta\ket^{\vert A\vert +\vert B\vert
+\sigma_3}}\big )$).
\item[(iv)] For every $N\in\N$ there exists $C_N>0$ such that for $j=1,\ldots
,n$, 

$\vert \overline \partial _j\,F(X,y)\vert \leq C_N\big (\frac{\vert y\vert
}{\bra x\ket}\big )^N\cdot \frac 1{\bra x\ket^{1+\sigma_3}}$ (resp. $C_N\,\vert
y\vert ^N\big[\big (\frac 1{\bra x\ket}+\frac 1{\bra\theta\ket}\big )^N\big
(\frac 1{\bra x\ket^{1+\sigma_3}}+\frac 1{\bra\theta\ket^{1+\sigma_3}}\big)
\big]$) 

where $\overline \partial _j=\frac 12 \big (\frac\partial {\partial
x_j}+i\,\frac\partial {\partial y_j}\big )$.
\end{itemize}
\end{lemma}

\noindent {\bf Proof } See Section \ref{ssVIII.4}. in the Appendix.

 Now recall
that for $\alpha\in T^*\R^n$ such that $\alpha_x\cdot \alpha_\xi\leq
-c_0\,\bra\alpha_x\ket \v\alpha_\varepsilon\vert $ and $\frac 12\leq \vert
\alpha_\varepsilon\vert \leq 2$ we have constructed in Theorem \ref{tIV.4.10} a
function $\Phi=\Phi(\theta,x,\alpha)$ uniformly bounded on the set
$\Omega_\delta$. By Lemma \ref{lV.3.1} we can extend $\Phi$ almost analytically
as a function, which we denote by $\Phi(\theta,z,\alpha)$, on the set
\begin{equation}\label{eqV.3.3}
\Omega^{\C}_\delta=\big\{(\theta,z)\in\R\times\C^n : \vert
z-x(\theta,\alpha)\vert \leq \delta\,\bra\theta\ket\,, \Re z\cdot
\alpha_\xi\leq c_0\, \bra \Re z\ket \vert \alpha_\xi\vert \,, \vert \Im z\vert
\leq \delta\big\}
\end{equation}
and $\Phi(\theta,z,\alpha)$ is still uniformly bounded on this set.

Again by Lemma \ref{lV.3.1} one can extend almost analytically the coefficients
of our symbol $p$, keeping the bounds of its coefficients. In that follows for
$z\in\C^n$ we shall denote by $X(t,\theta,z)$ the solution, whenever it exists,
of the following problem.
\begin{equation}\label{eqV.3.4}\left\{
\begin{array}{l}
\dot X(t,\theta,z)=\frac{\partial p}{\partial
\xi}\,(X(t,\theta,z),\Phi(t,X(t,\theta,z),\alpha))\,,\\
X(\theta,\theta,z)=z\,.
\end{array}\right.
\end{equation}
Our aim is to prove the following result.
\begin{theorem}\sl \label{tV.3.2}
One can find positive constants $c_1$, $\delta_1$, $K$, $\tilde K$, with
$c_1\ll c_0$, $\delta_1\ll \delta$, such that for all $x\in\R^n$ such that
$$
\vert x-x(\theta,\alpha)\vert \leq \delta_1\,\bra\theta\ket\,,\quad x\cdot
\alpha_\xi\leq c_1\,\bra x\ket\,\vert \alpha_\xi\vert \,,
$$
the solution of (\ref{eqV.3.4}) exists on $[0,\theta]$ and satisfies the
estimates,
\begin{equation}\label{eqV.3.5}\left\{
\begin{array}{ll}
(i)&\quad \vert X(t,\theta,x)-x(t,\alpha)\vert \leq K \,\vert
x-x(\theta,\alpha)\vert \,\frac{\bra t\ket}{\bra\theta\ket}\,,\\
 (ii)&\quad \vert \Im X(t,\theta,x)\vert \leq 
K\,\frac{\vert x-x(\theta,\alpha)\vert }{\bra \theta\ket}\, \\
(iii)&\quad  \bra x\ket+\bra\theta-t\ket\leq K\,\bra\Re X(t,\theta,x)\ket\,, \\
(iv)&\quad \Re X(t,\theta,x)\cdot \alpha_\xi\leq \frac 1{\tilde K}\,\bra\Re
X(t,\theta,x)\ket \vert \alpha_\xi\vert \,, 
\end{array}\right.
\end{equation}
uniformly for $t\in[0,\theta]$.
\end{theorem}

Let us remark that the estimates (\ref{eqV.3.5}) ensure in particular that if
$\delta_1$ is small enough  we have $(t,X(t,\theta,\alpha))\in\Omega^\C_\delta$.
With $0<c_1\ll c_2\ll c_0$ to be chosen, we divide the proof in three cases.

\noindent Case 1 : $x\cdot \alpha_\xi\leq c_2\,\bra x\ket \vert
\alpha_\xi\vert$, $x(\theta,\alpha)\cdot \alpha_\xi\leq
c_2\,\bra x(\theta,\alpha)\ket \vert \alpha_\xi\vert $, $\vert
x-x(\theta,\alpha)\vert \leq \vert x(\theta,\alpha)\vert $.

\noindent Case 2 : $x\cdot \alpha_\xi\leq c_2\,\bra x\ket \vert
\alpha_\xi\vert $, $x(\theta,\alpha)\cdot \alpha_\xi\leq
c_2\,\bra x(\theta,\alpha)\ket \vert \alpha_\xi\vert $, $\vert
x-x(\theta,\alpha)\vert >\vert x(\theta,\alpha)\vert $.

\noindent Case 3 : $x\cdot \alpha_\xi\leq c_1\,\bra x\ket\vert
\alpha_\xi\vert $, $x(\theta,\alpha)\cdot \alpha_\xi>c_2\, \bra
x(\theta,\alpha)\ket \vert \alpha_\xi\vert $.

Here is the geometrical interpretation of case 1 and 2. We denote by $[a,b]$
the segment joining two points $a,b\in\R^n$.
\begin{lemma}\sl \label{lV.3.3}
Let $c_2>0$ and assume that $x\in\R^n$ is such that $x\cdot
\alpha_\xi\leq c_2\,\bra x\ket\vert \alpha_\xi\vert $, $\vert
x-x(\theta,\alpha)\vert \leq \delta\,\bra\theta\ket$ and that
$x(\theta,\alpha)\cdot \alpha_\xi\leq c_2\,\bra x(\theta,\alpha)\ket
\vert \alpha_\xi\vert $. Then we have :
\begin{itemize}
\item[(i)] either $\vert x-x(\theta,\alpha)\vert \leq \vert
x(\theta,\alpha)\vert $ and then,
$$
\forall \,y\in[x,x(\theta,\alpha)]\,,\quad y\cdot \alpha_\varepsilon\leq
2\,c_2\,\bra y\ket \vert \alpha_\xi\vert \,,
$$
\item[(ii)] or $\vert x-x(\theta,\alpha)\vert > \vert x(\theta,\alpha)\vert
$ and then, $[0,x(\theta,\alpha)]\cup[0,x]\subset\{y\in\R^n : \vert
y-x(\theta,\alpha)\vert \leq \delta_0\,\bra\theta\ket$ and $y\cdot
\alpha_\xi\leq c_2\,\bra y\ket \vert \alpha_\xi\vert \}$.
Moreover
$$
\vert x\vert +\vert x(\theta,\alpha)\vert \leq 3\,\vert x-x(\theta,\alpha)\vert
\leq 3\,(\vert x\vert +\vert x(\theta,\alpha)\vert )\,.
$$
\end{itemize}
\end{lemma}

\noindent {\bf Proof } In the first case applying Lemma \ref{lIV.4.16} we
obtain for $t\in[0,1]$
\begin{equation*}
\begin{aligned}
(tx+(1-t)\,x(\theta,\alpha))\cdot \alpha_\xi&\leq c_2 (t\bra
x\ket+(1-t)\bra x(\theta,\alpha)\ket)\vert \alpha_\xi\vert \\
&\leq c_2 (1+t\vert x\vert +(1-t)\vert x(\theta,\alpha)\vert )\vert
\alpha_\xi\vert \\
&\leq c_2 (1+\sqrt 2 \vert t x+(1-t)\,x(\theta,\alpha)\vert )\vert
\alpha_\xi\vert \\
&\leq 2\,c_2\bra tx+(1-t)\,x(\theta,\alpha)\ket \vert \alpha_\xi\vert
\,.
\end{aligned}
\end{equation*}
Assume now that $\vert x-x(\theta,\alpha)\vert >\vert x(\theta,\alpha)\vert
$. Then $\vert x(\theta,\alpha)\vert \leq \delta_0 \bra \theta\ket$. Therefore
$[0,x]\cup[0,x(\theta,\alpha)]\subset
B(x(\theta,\alpha),\delta_0\bra\theta\ket)$. Now if $t\in[0,1]$ and $Z=x$ or
$x(\theta,\alpha)$ we have $t\,Z\cdot \alpha_\xi\leq t\,c_2\bra Z\ket
\vert \alpha_\xi\v \leq c_2\,\bra t\,Z\ket \vert
\alpha_\xi\vert $. Moreover $3\,\vert x-x(\theta,\alpha)\vert =\vert
x-x(\theta,\alpha)\vert +2\,\vert x-x(\theta,\alpha)\vert \geq \vert x\vert
-\vert x(\theta,\alpha)\vert +2\,\vert x(\theta,\alpha)\vert =\vert x\vert
+\vert x(\theta,\alpha)\vert $. 

\cqfd

\noindent {\bf 1) Proof of Theorem \ref{tV.3.2} in case 1 and 2 } 

Let us take
$c_2$, $\delta_2$ such that $0<c_2\ll c_0$, $0<\delta_2\ll\delta$. Let $A$ be
the set of $T\in[0,\theta]$ such that for every $z\in\C^n$ such that $\vert
z-x(\theta,\alpha)\vert \leq \delta_2\,\bra\theta\ket$, $\Re z\cdot
\alpha_\xi\leq c_2\,\bra\Re z\ket \vert \alpha_\xi\vert $,
$\vert \Im z\vert \leq \delta_2$, $x(\theta,\alpha)\cdot \alpha_\xi\leq
c_2 \bra x(\theta,\alpha)\ket \vert \alpha_\xi\vert $ the problem
(\ref{eqV.3.5}) has a unique solution on $[T,\theta]$ which satisfies for
$t\in[T,\theta]$, in case 1 :
\begin{equation}\label{eqV.3.6}
\vert X(t,\theta,z)-x(t,\alpha)\vert \leq M_1\,\vert z-x(\theta,\alpha)\vert
\,\frac{\bra t\ket}{\bra\theta\ket}
\end{equation}
\begin{equation}\label{eqV.3.7}
\bra \theta-t\ket+\bra\Re (sz+(1-s)\,x(\theta,\alpha))\ket\leq M_1\bra s\,\Re
X(t,\theta,z)+(1-s)\,x(t,\alpha)\ket\,,\quad s\in[0,1]
\end{equation}
\begin{equation}\label{eqV.3.8}
(s\,\Re X(t,\theta,z)+(1-s)\,x(\theta,\alpha))\cdot \alpha_\xi\leq
M_2\,\bra s\,\Re X(t,\theta,z)+(1-s)\,x(\theta,\alpha)\ket\,\vert
\alpha_\xi\vert 
\end{equation}
\begin{equation}\label{eqV.3.9}
\vert \Im X(t,\theta,z)\vert \leq M_3 \Big (\frac{\vert z-x(\theta,\alpha)\vert
}{\bra\theta\ket}+\vert \Im z\vert \Big )
\end{equation}
in case 2 :
\begin{equation}\label{eqV.3.10}
\vert X(t,\theta,z)-X(t,\theta,0)\vert \leq M_1\,\vert z\vert
\,\frac{\bra t\ket}{\bra\theta\ket}
\end{equation}
\begin{equation}\label{eqV.3.11}
\bra\theta-t\ket+\bra\Re (sz)\ket\leq M_1 \bra s\Re X(t,\theta,z)+(1-s)\,\Re
X(t,\theta,0)\ket\,, \quad s\in[0,1]
\end{equation}
\begin{equation}\label{eqV.3.12}
\big (s\,\Re X(t,\theta,z)+(1-s)\,\Re X(t,\theta,0)\big )\cdot
\alpha_\xi\leq M_2\,\bra s\,\Re X(t,\theta,z)+(1-s)\,\Re
X(t,\theta,0)\ket
\end{equation}
and (\ref{eqV.3.9}).

Our aim is to show that if $M_1, M_2, M_3$ are correctly chosen then
$A=[0,\theta]$.

Let us show that the set $A$ is not empty. Indeed if $t=\theta$ the estimates
(\ref{eqV.3.6}) to (\ref{eqV.3.12}) are satisfied with strict inequalities if
$M_1>1$, $M_2>2\,C_2$, $M_3>1$ (using Lemma \ref{lV.3.3}). It follows that
they still hold for $T=\theta-\varepsilon$, if $\varepsilon$ is small enough.

On the other hand $A$ is an interval. Let $T_*=\inf A$. If $T_*=0$ then the
theorem \ref{tV.3.2} is proved. Assume then that $T_*>0$ and let $T\geq T_*$.
Then on $[T,\theta]$, (\ref{eqV.3.6}) to (\ref{eqV.3.12}) hold.
\begin{remark}\sl \label{rV.3.4}
If the case 2 is not empty then the point $z_0=0$ satisfies all the
requirements of case\,1. Indeed if there exists $z_1$ such that $\vert
z_1-x(\theta,\alpha)\vert >\vert x(\theta,\alpha)\vert $ then $\vert
x(\theta,\alpha)\vert <\delta_2\,\bra\theta\ket$ so $\vert
0-x(\theta,\alpha)\vert <\delta_2\,\bra\theta\ket$ and the other requirements
are trivial. Therefore if the case\,2 is not empty then $X(t,\theta,0)$ is well
defined on $[T,\theta]$ and satisfies (\ref{eqV.3.6}) to (\ref{eqV.3.9}).
\end{remark}

Let us show that we have $(t,X(t,\theta,z))\in\Omega^\C_\delta$ (see
(\ref{eqV.3.3})). This is the case if $M_1\,\delta_2\leq \delta$, $M_2\leq
c_0$,
 $2\,M_3\,\delta_2\leq \delta$. Indeed the only non trivial point is to prove
that $\vert X(t,\theta,z)-x(t,\alpha)\vert \leq \delta\,\bra t\ket$ in
case~2. We have
$$
\vert X(t,\theta,z)-x(t,\alpha)\vert \leq \vert
X(t,\theta,z)-X(t,\theta,0)\vert +\vert X(t,\theta,0)-x(t,\alpha)\vert
=(1)+(2)\,.
$$
It follows from (\ref{eqV.3.10}) that $(1)\leq M_1\,\vert z\vert \,\frac{\bra
t\ket}{\bra\theta\ket}$ and from (\ref{eqV.3.6}) with $z=0$, that $(2)\leq
M_1\,\vert x(\theta,\alpha)\vert \,\frac{\bra t\ket}{\bra\theta\ket}$. Now by
Lemma \ref{lV.3.3} (ii) we have $\vert \Re z\vert +\vert x(\theta,\alpha)\vert
<3\,\vert \Re z-x(\theta,\alpha)\vert $; since $\vert \Im z\vert \leq \delta_2$
we will have $(1)+(2)\leq M_1(\delta_2+3\,\v\Re
z-x(\theta,\alpha)\v)\,\frac{\bra t\ket}{\bra\theta\ket}$. Since $\v\Re
z-x(\theta,\alpha)\v\leq \delta_2\,\bra\theta\ket$ we obtain finally
$(1)+(2)\leq 4\,M_1\,\delta_2\,\bra t\ket\leq \delta\,\bra t\ket$.

In the sequel we shall denote by $C$ or $O(1)$ the constants which may depend
on bounds of $p,\Phi$ but are independent of $M_1,M_2, M_3$. Moreover for the
sake of symplicity we shall write
\begin{equation}\label{eqV.3.13}\left\{
\begin{array}{l}
X(t)=X(t,\theta,z)\\
\tilde X(t)=x(t,\alpha)\textrm{ in case 1}, X(t,\theta,0)\textrm{ in case 2}\,.
\end{array}\right.
\end{equation}
In particular $\tilde X(\theta)=x(\theta,\alpha)$ in case~1 and $\tilde
X(\theta)=0$ in case~2. Our goal is to show that the estimates (\ref{eqV.3.6})
to (\ref{eqV.3.12}) hold on $[T,\theta]$ with better constants than
$M_1,M_2,M_3$.

\noindent {\bf a) Improvement of (\ref{eqV.3.7}) and (\ref{eqV.3.11})}

By Theorem \ref{tIV.4.10} (iii) we have
$\Phi(\theta,x,\alpha)-\alpha_\xi=\CO(\varepsilon+\delta)$ if
$(\theta,x)\in\Omega_\delta$ and by Lemma \ref{lV.3.1} this estimate still hold
on $\Omega^\C_\delta$~; it follows that
$\Phi(t,X(t),\alpha)-\alpha_\xi=\CO(\varepsilon+\delta)$. On the other
hand $\frac{\partial p}{\partial \xi}\,(x,\xi)-2\xi=\CO(\varepsilon)\v\xi\v$
which also extends for $z\in\C^n$, $\v\Im z\v\leq \delta_2$. It follows then
from (\ref{eqV.3.4}) that $\dot
X(t)=2\alpha_\varepsilon+\CO(\varepsilon+\delta)$. Therefore
\begin{equation}\label{eqV.3.14}\left\{
\begin{array}{l}X(t)=z-2(\theta-t)\,\alpha_\xi+\CO(\varepsilon+\sqrt\delta)(\theta-t)\\
\tilde X(t)=\tilde
X(\theta)-2(\theta-t)\,\alpha_\xi+\CO(\varepsilon+\sqrt\delta)(\theta-t)\,.
\end{array}\right.
\end{equation}
Now for $s\in[0,1]$
\begin{equation*}
\begin{split}
(1)&=\v\Re (sX(t)+(1-s)\,\tilde X(t))\v^2=\v s\,\Re z+(1-s)\,\Re \tilde
X(\theta)\v^2+4 (\theta-t)^2\,\v\alpha_\xi\v^2\\
&\quad -4(\theta-t)(s\,\Re z+(1-s)\,\Re \tilde X(\theta))\cdot
\alpha_\xi+\CO((\varepsilon+\delta)\big[(\theta-t)^2+\v s\,\Re
z+(1-s)\,\Re \tilde X(\theta)\v^2\big]\,.
\end{split}
\end{equation*}
It follows from the conditions on $z$ and the definition of $\tilde X(\theta)$
that $(s\,\Re z+(1-s)\,\Re \tilde X(\theta))\cdot \alpha_\varepsilon\leq
2c_2\bra s\,\Re z+(1-s)\,\Re \tilde X(\theta)\ket\,\v\alpha_\xi\v$ so
$$
(1)\geq \frac 12 \v s\Re z+(1-s)\,\Re \tilde
X(\theta)\v^2+3(\theta-t)^2\,\v\alpha_\xi\v^2-8c_2 \bra s \Re
z+(1-s)\,\Re \tilde X(\theta)\ket (\theta-t)\v\alpha_\xi\v
$$
if $\varepsilon+\delta$ is small enough. It follows that
$$
(1)\geq \Big (\frac 12 -16\,c_2\Big )\v s\Re z+(1-s)\,\Re \tilde
X(\theta)\v^2+(3-16\,c_2)(\theta-t)^2\,\v\alpha_\xi\v^2-16\,c_2\,.
$$
If $c_2$ has been chosen small enough we obtain in particular
\begin{equation}\label{eqV.3.15}\left\{
\begin{array}{ll}
\textrm{(i)}& \v\Re (s\,X(t)+(1-s)\,\tilde X(t))\v^2\geq \frac 14\, \v\Re
(s\,z+(1-s)\,\tilde X(\theta)\v^2+2(\theta-t)^2\,\v\alpha_\xi\v^2-\frac
12\\
\textrm{(ii)}&  \bra\Re (s\,X(t)+(1-s)\,\tilde X(t))\ket^2\geq \frac 1{10}\,
[\bra\theta-t\ket^2+\bra\Re (s\,z+(1-s)\,\tilde X(\theta))\ket^2]\,.
\end{array}\right.
\end{equation}
This improves (\ref{eqV.3.7}) and (\ref{eqV.3.11}) if $M_1>4$.

\noindent {\bf b) Improvement of (\ref{eqV.3.8}), (\ref{eqV.3.12})}

It follows from (\ref{eqV.3.14}) that
$$
(2)=\Re (s\,X(t)+(1-s)\,\tilde X(t))\cdot \alpha_\xi=\Re
(s\,z+(1-s)\,\tilde X(\theta))\cdot
\alpha_\xi-2(\theta-t)\v\alpha_\xi\v^2+\CO(\varepsilon+\delta)(\theta-t)~.
$$
Applying Lemma \ref{lV.3.3} we obtain if $\varepsilon+\delta$ is small,
$$
(2)\leq 2 \,c_2  \,\bra\Re (s\,z+(1-s)\,\tilde X(\theta))\ket
-(\theta-t)\v\alpha_\xi\v^2\,.
$$
Using (\ref{eqV.3.15}) (i) we obtain, $(2)\leq 4\,c_2 \bra\Re
(s\,X(t)+(1-s)\,\tilde X(t)\ket$. Taking $16\,c_2\leq M_2$ we deduce finally
that
\begin{equation}\label{eqV.3.16}
\Re (s\,X(t)+(1-s)\,\tilde X(t))\cdot \alpha_\xi\leq \frac 12\,
M_2\bra\Re (s\,X(t)+(1-s)\,\tilde X(t))\ket\,\v\alpha_\xi\v\,.
\end{equation}
This improves (\ref{eqV.3.8}) and (\ref{eqV.3.12}).

\noindent {\bf c) Improvement of (\ref{eqV.3.6}) and (\ref{eqV.3.10})}

We have
\begin{equation*}\left\{
\begin{array}{l}
\dot X(t)=\frac{\partial p}{\partial \xi}\,(X(t),\,\Phi(t,X(t),\alpha))\,,\\
\dot{\tilde X}(t)=\frac{\partial p}{\partial \xi}\, (\tilde
X(t),\,\Phi(t,\tilde X(t),\alpha))\,,
\end{array}\right.
\end{equation*}
the second equation being true in the case~1 according to the identity
$\Phi(t,x(t,\alpha),\alpha)=\xi(t,\alpha)$. Let us set 
\begin{equation}\label{eqV.3.17}
Z(t)=X(t)-\tilde X(t)\,.
\end{equation}
 Then
$$
\dot Z(t)=2\big[\Phi(t,X(t),\alpha)-\Phi (t,\tilde X(t),\alpha)\big
]+\frac{\partial q}{\partial \xi} \,(X(t),\Phi (t,X(t),\alpha))-\frac{\partial
q}{\partial \xi}\,(\tilde X(t),\Phi (t,\tilde X(t),\alpha))\,,
$$
since $p=\v\xi\v^2+q$.

Now we use (\ref{eqIV.4.39}) and Theorem \ref{tIV.4.2} (i). It follows after
extending almost analytically $\tilde a,\tilde b$  and the coefficients
of $q$ by Lemma
\ref{lV.3.1},
\begin{equation}\label{eqV.3.18}
\Phi (t,z,\alpha)=\xi(t,\alpha)+\frac{z-x(t,\alpha)}{2t-i}-\Big (\tilde a+\frac
i{\bra t\ket}\,\tilde b\Big )(t,z,\alpha)~.
\end{equation}
It follows then that
\begin{equation}\label{eqV.3.19}\left\{
\begin{split}
\dot Z(t)&=\frac{2Z(t)}{2t-i}-\big(\tilde a(t,X(t),\alpha)-\tilde a(t,\tilde
X(t),\alpha)\big )-\frac i{\bra t\ket}\big(\tilde b(t,X(t),\alpha)-\tilde
b(t,\tilde X(t),\alpha)\big )\\
&\enskip +\frac{\partial q}{\partial \xi}\,\big (X(t),\Phi (t,X(t),\alpha)\big
)-\frac{\partial q}{\partial \xi}\,(\tilde X(t),\Phi (t,X(t),\alpha))+
\frac{\partial q}{\partial \xi}\,(\tilde X(t),\Phi (t,X(t),\alpha))\\
&\enskip  -\frac{\partial q}{\partial \xi}\,\big (\tilde X(t),\Phi (t,\tilde
X(t),\alpha)\big )\,.
\end{split}\right.
\end{equation}
We have the following lemma.
\begin{lemma}\sl \label{lV.3.5}
One can find a positive constant $C$ such that
$$
\Big\v\dot Z(t)-\frac{2Z(t)}{2t-i}\Big\v\leq
C(\varepsilon+\delta)\,\v Z(t)\v\Big (\frac 1 {\bra\theta-t\ket^2}+\frac 1 {\bra
t\ket^2}\Big )\,.
$$
\end{lemma}

\noindent {\bf Proof} 

\noindent (i) Estimation of $(1)=\tilde a(t,X(t),\alpha)-\tilde a(t,\tilde
X(t),\alpha)$. We have
\begin{equation*}
\begin{split}
(1)&=\int^1_0 \frac{\partial \tilde a}{\partial z}\,(t,s\, X(t)+(1-s)\,\tilde
X(t),\alpha)(X(t)-\tilde X(t))ds\\
&\quad +\int^1_0 \frac{\partial \tilde a}{\partial \overline z}
\,(t,s\,X(t)+(1-s)\,\tilde X(t),\alpha)\,\overline {(X(t)-\tilde X(t))}\,ds\,.
\end{split}
\end{equation*}
Using the estimates on $\tilde a$ given in Theorem \ref{tIV.4.2} and Lemma
\ref{lV.3.1} with $\sigma_3=1$ we find
$$
\Big(\Big\v\frac{\partial \tilde a}{\partial z}\Big\v+\Big\v\frac{\partial
\tilde a}{\partial \overline z}\Big\v\Big)(t,\cdots )\leq
C(\varepsilon+\delta)\Big (\frac 1{\bra\Re (s\,X(t)+(1-s)\,\tilde
X(t))\ket^2}+\frac 1{\bra t\ket^2}\Big )\, .
$$
Using (\ref{eqV.3.7}) and (\ref{eqV.3.12}) we deduce that
$$
\Big\v\frac{\partial \tilde a}{\partial z}\Big\v+\Big\v\frac{\partial
\tilde a}{\partial \overline z}\Big\v(t,\cdots )\leq
C (M_1) (\varepsilon+\delta)\Big (\frac 1{\bra\theta-t\ket^2}+\frac 1{\bra
t\ket^2}\Big )\,.
$$
It follows that,
\begin{equation}\label{eqV.3.20}
\v(1)\v\leq C (M_1)(\varepsilon+\delta)\,\v Z(t)\v\Big (\frac
1{\bra\theta-t\ket^2}+\frac 1{\bra t\ket^2}\Big )\,.
\end{equation}
Here $C(M_1)$ is a constant depending only on $M_1$.

\noindent (ii) Setting $(2)=\frac 1{\bra t\ket}\,(\tilde
b(t,X(t),\alpha)-\tilde b(t,\tilde X(t),\alpha))$ we have exactly by the same
way
\begin{equation}\label{eqV.3.21}
\v(2)\v\leq C (M_1)(\varepsilon+\delta)\, \frac{\v Z(t)\v}{\bra t\ket} \Big
(\frac 1{\bra\theta-t\ket^2}+\frac 1{\bra t\ket^2}\Big )\,.
\end{equation}

\noindent (iii) Estimation of $(3)=\frac{\partial q}{\partial \xi}\,(X(t),\Phi
(t,X(t),\alpha))-\frac{\partial q}{\partial \xi}\,(\tilde X(t),\Phi
(t,X(t),\alpha))$. We have
\begin{equation*}
\begin{split}
(3)=\int^1_0 \frac{\partial ^2q}{\partial z\,\partial \xi}\,
(s\,X(t)&+(1-s)\,\tilde X(t),\Phi (t,X(t),\alpha))(X(t)-\tilde
X(t))\,ds\\
&+\textrm{ analogue term with } \frac{\partial ^2q}{\partial \overline
z\,\partial \xi}\,.
\end{split}
\end{equation*}
Now the coefficients of $q$ say $b_{jk}$, extended by Lemma \ref{lV.3.1} satisfy
$$
\Big\v\frac{\partial b_{jk}}{\partial z}\,(z)\Big\v+\Big\v\frac{\partial
b_{jk}}{\partial \overline z}\,(z)\Big\v\leq \frac{C\,\varepsilon}{\bra
z\ket^{2+\sigma_0}}\,.
$$
Using again (\ref{eqV.3.7}) and (\ref{eqV.3.13}) we obtain
\begin{equation}\label{eqV.3.22}
\v(3)\v\leq \frac{C\,\varepsilon \v Z(t)\v}{\bra\theta-t\ket^{2+\sigma_0}}\,.
\end{equation}

\noindent (iv) Estimation of $(4)=\frac{\partial q}{\partial \xi}\,(\tilde
X(t),\Phi (t,X(t),\alpha))-\frac{\partial q}{\partial \xi}\,(\tilde X(t),\Phi
(t,\tilde X(t),\alpha))$. We have, by (\ref{eqV.3.18}),
$$
\v\Phi (t,X(t),\alpha)-\Phi (t,\tilde X(t),\alpha)\v\leq C\,\frac{\v
Z(t)\v}{\bra t\ket}+\v(1)\v+\v(2)\v\,.
$$
On the other hand (\ref{eqV.3.7}), (\ref{eqV.3.11}) with $s=0$ imply that
$M_1\,\bra\Re \tilde X(t)\ket\geq \bra\theta-t\ket$. Therefore using the decay
of the coefficients $b_{jk}$ of $q$ and the estimates (\ref{eqV.3.20}),
(\ref{eqV.3.21}) we obtain
$$
\v(4)\v\leq \frac{C\,\varepsilon}{\bra\theta-t\ket^{1+\sigma_0}}\,\v Z(t)\v
\Big (\frac 1{\bra t\ket}+\frac 1 {\bra\theta-t\ket^2}\Big )\,.
$$
It follows then that,
\begin{equation}\label{eqV.3.23}
\v(4)\v\leq C\,\varepsilon\,\v Z(t)\v
\Big (\frac 1{\bra t\ket^2}+\frac 1 {\bra\theta-t\ket^2}\Big )\,.
\end{equation}
Gathering the estimates (\ref{eqV.3.20}) to (\ref{eqV.3.23}) we obtain the
claim of the Lemma. \cqfd

Next we state the following Lemma.
\begin{lemma}\sl \label{lV.3.6}
Let $0\leq T<\theta$. Let $Y(t)=(Y_1(t),\ldots ,Y_n(t))\in\C^n$ be such that
$Y\in C^1([T,\theta])$ and satisfies on $[T,\theta]$ the inequality
$$
\Big\v\dot Y(t)-\frac 2{2t-i}\,Y(t)\Big\v\leq \v h(t)\v\,\v Y(t)\v+\v g(t)\v\,,
$$
for some continuous fonctions $h,g$. Then for all $t$ in $[T,\theta]$ we have
$$
\v Y(t)\v\leq \Big (\frac{\bra 2t\ket}{\bra 2\theta\ket}\,\v
Y(\theta)\v+\bra 2t\ket \int^\theta_T \frac{\v g(s)\v}{\bra 2s\ket}\,ds\Big)
\,\exp
\Big(\int^\theta_T
\v h(s)\v\,ds\Big )\,.
$$
\end{lemma}

\noindent {\bf Proof } Let us set $W(t)=\frac{Y(t)}{2t-i}$. Then $\v
W(t)\v=\frac{\v Y(t)\v}{\bra 2t\ket}$,
$$
\dot W(t)=\frac{\dot Y(t)}{2t-i}-\frac{2Y(t)}{(2t-i)^2}=\frac 1{2t-i}\Big (\dot
Y(t)-\frac{2Y(t)}{2t-i}\Big )\,.
$$
It follows that $\v\dot W(t)\v\leq \frac 1{\bra 2t\ket}\,(\v h(t)\v\,\v
Y(t)\v+\v g(t)\v)\leq \v h(t)\v\,\v W(t)\v+\frac{\v g(t)\v}{\bra 2t\ket}$. 
Then, for $t\geq T$ and $\sigma\in[t,\theta]$,
$$
\v W(\sigma)\v\leq \v W(\theta)\v+\int^\theta_\sigma \v h(s)\v\,\v
W(s)\v\,ds+\int^\theta_t \frac{\v g(s)\v}{\bra 2s\ket}\, ds\,.
$$
By the Gronwall Lemma we obtain
$$
\v W(t)\v\leq \Big (\v W(\theta)\v+ \int^\theta_T \frac {\v g(s)\v}{\bra
2s\ket}
\,ds\Big )\,\exp \Big (\int^\theta_T \v h(s)\v\,ds\Big )\,.
$$
Coming back to $Y(t)$ we obtain the claim of the Lemma. \cqfd
\begin{corollary}\sl \label{cV.3.7}
Let $Z(t)$ be defined by (\ref{eqV.3.17}). Then if $\varepsilon+\delta$ is
small compared to $M_1$ we have
$$
\v Z(t)\v\leq 2\,\frac{\bra 2t\ket}{\bra 2\theta\ket}\,\v z-\tilde
X(\theta)\v\,.
$$
\end{corollary}

\noindent {\bf Proof } We apply Lemma \ref{lV.3.6} and Lemma \ref{lV.3.5} with
$g(t)=0$ and $h(t)=C(\varepsilon+\delta)\big (\frac 1{\bra\theta-t\ket^2}+\frac
1{\bra t\ket^2}\big )$. Then $\int^\theta_T \v h(s)\v\,ds\leq C
(\varepsilon+\delta) \int_\R \frac{d\sigma}{\bra\sigma\ket^2}$. It follows that
$$
\v Z(t)\v\leq e^{C'(\varepsilon+\delta)}\cdot \frac{\bra
2t\ket}{\bra2\theta\ket}\, \v Z(\theta)\v\, .
$$
Since $Z(\theta)=X(\theta)-\tilde X(\theta)=z-\tilde X(\theta)$ our lemma
follows. \cqfd

We can now show the improvement of (\ref{eqV.3.6}) and (\ref{eqV.3.10}). In the
case~1 we have $\tilde X(\theta)=x(\theta,\alpha)$ and in the case~2, $\tilde
X(\theta)=0$. Therefore in case~1 we find by Corollary \ref{cV.3.7}, 
$$\v
X(t,\theta,z)-x(t,\alpha)\v\leq 4\,\frac{\bra t\ket}{\bra\theta\ket}\,\v
z-x(\theta,\alpha)\v\,,
$$
 and in case\,2,
$$
\v X(t,\theta,z)-X(t,\theta,0)\v\leq 4\frac{\bra t\ket}{\bra\theta\ket}\,\v
z\v\,.
$$
Taking $M_1>4$ this shows that (\ref{eqV.3.6}) and (\ref{eqV.3.10}) have been
improved.

\noindent {\bf d) Improvement of (\ref{eqV.3.9})}

Let us set
\begin{equation}\label{eqV.3.24}\left\{
\begin{array}{ll}
U(t)&=\Im (X(t)-\tilde X(t))\,,\\
(1)&=\frac{4t}{1+4t^2}\,U(t)\,,\\
(2)&=\frac{2 \Re (X(t)-\tilde X(t))}{1+4 t^2}\,,\\
(3)&=-\Im [\tilde a(t,X(t),\alpha)-\tilde a(t,\tilde X(t),\alpha)]\,,\\
(4)&=-\frac{1}{\bra t\ket}\,\Re [\tilde b(t,X(t),\alpha)-\tilde b(t,\tilde
X(t),\alpha)]\,,\\
(5)&=\Im \Big[\frac{\partial q}{\partial \xi}\,(\Re X(t),\Phi
(t,X(t),\alpha))-\frac{\partial q}{\partial \xi}\,(\Re X(t),\Phi (t,\tilde
X(t),\alpha))\Big]\,,\\
(6)&=\Im \Big[\frac{\partial q}{\partial \xi}\,(\Re X(t),\Phi
(t,\tilde X(t),\alpha))-\frac{\partial q}{\partial \xi}\,(\Re \tilde X(t),\Phi
(t,\tilde X(t),\alpha))\Big]\,,\\
(7)&=\Im \int^1_0 \frac{\partial ^2q}{\partial \xi\,\partial z}\, (\Re
X(t)+is\,\Im X(t),\Phi (t,X(t),\alpha))\,ds\, (i\,\Im X(t))\,,\\
(8)&=\Im \int^1_0 \frac{\partial ^2q}{\partial \xi\,\partial\overline  z}\, (\Re
X(t)+is\,\Im X(t),\Phi (t,X(t),\alpha))\,ds\, (-i\,\Im X(t))\,,\\
(9)&=-\Im \int^1_0 \frac{\partial ^2q}{\partial \xi\,\partial z}\, (\Re
\tilde X(t)+is\,\Im \tilde X(t),\Phi (t,\tilde X(t),\alpha))\,ds\, (i\,\Im
\tilde X(t))\,,\\
(10)&=-\Im \int^1_0 \frac{\partial ^2q}{\partial \xi\,\partial \overline z}\,
(\Re \tilde X(t)+is\,\Im \tilde X(t),\Phi (t,\tilde X(t),\alpha))\,ds\, (-i\,\Im
\tilde X(t))\,.
\end{array}\right.
\end{equation}
Then it follows from (\ref{eqV.3.17}) and (\ref{eqV.3.19}) that,
\begin{equation}\label{eqV.3.25}
\dot U(t)=\frac{4t}{1+4t^2}\,U(t)+\sum^{10}_{i=2}\,(i)\,.
\end{equation}
\begin{lemma}\sl \label{lV.3.8}
With the above notations, if $\varepsilon+\delta$ is small enough we have
\begin{equation*}
\begin{split}
\Big\v \dot U(t)-\frac{4t}{1+4t^2}\,U(t)\Big\v&\leq \frac{3\,M_1\,\v
z-x(\theta,\alpha)\v}{\bra\theta\ket \bra t\ket}+\frac{\v
z-x(\theta,\alpha)\v}{\bra\theta\ket} \,\Big (\frac
1{\bra\theta-t\ket^{1+\sigma_0}}+\frac 1{\bra t\ket^2}\Big )\\
&\qquad +\Big (\frac 1{\bra\theta-t\ket^2}+\frac 1{\bra t\ket^2}\Big )\,\v
U(t)\v\,.
\end{split}
\end{equation*}
\end{lemma}

\noindent {\bf Proof } We use (\ref{eqV.3.25}) and (\ref{eqV.3.24}). We
estimate the terms $(i)$ for $i=2,\ldots ,10$.

\noindent (i) Estimation of (2). It follows from (\ref{eqV.3.6}) and
(\ref{eqV.3.10}), since $\tilde X(\theta)=x(\theta,\alpha)$ in case~1 and
$\tilde X(\theta)=0$ in case~2, that
\begin{equation*}
\frac{\v X(t)-\tilde X(t)\v}{1+4t^2}\leq M_1\, \frac{\bra t\ket}{\bra\theta\ket
(1+4t^2)} \,\v z-\tilde X(\theta)\v\leq \frac{M_1}{\bra\theta\ket \bra
t\ket}\quad 
\begin{cases}
\v z-x(\theta,\alpha)\v\,,&\textrm{case 1}\\
\v z\v\,, &\textrm{case 2}
\end{cases}\,.
\end{equation*}
But in case~2 according to Lemma \ref{lV.3.3} (ii) we have $\v z\v\leq 3\,\v
z-x(\theta,\alpha)\v$. It follows that in both cases we have
\begin{equation}\label{eqV.3.26}
\v(2)\v\leq \frac{3\,M_1\,\v z-x(\theta,\alpha)\v}{\bra\theta\ket \bra t\ket}\,.
\end{equation}

\noindent (ii) Estimation of (3) and (4). We note that $\tilde a(t,\Re
X(t),\alpha)$ and $\tilde a(t,\Re \tilde X(t),\alpha)$ are real. It follows that
\begin{equation*}
\begin{split}
(3)=&-\Im \int^1_0 \frac{\partial \tilde a}{\partial z}\,(t,\Re X(t)+s\,i\,\Im
X(t),\alpha)\,ds (i\,\Im X(t))\\
&-\Im \int^1_0 \frac{\partial \tilde a}{\partial
\overline z}\quad  \textrm{ (idem )} ds
(-i\,\Im X(t))+\Im \int^1_0 \frac{\partial \tilde a}{\partial z}\,(t,\Re
\tilde X(t)+s\,i\,\Im \tilde X(t),\alpha)\,ds\,(i\,\Im \tilde X(t))\\
&+\Im \int^1_0\frac{\partial \tilde a}{\partial \overline z}\quad \textrm{
(idem) }\,ds\, (-i\,\Im \tilde X(t))~.
\end{split}
\end{equation*}
Now, according to Theorem \ref{tIV.4.2} and Lemma \ref{lV.3.1} we have
$$
\Big\v\frac{\partial \tilde a}{\partial
z}\,(t,w,\alpha)\Big\v+\Big\v\frac{\partial \tilde a}{\partial \overline
z}\,(t,w,\alpha)\Big\v\leq C\,(\varepsilon+\delta)\Big (\frac 1{\bra\Re
w\ket^2}+\frac 1{\bra t\ket^2}\Big )\,.
$$
We use this estimate with $w=\Re X(t)+i\,s\,\Im X(t)$ and $w=\Re \tilde
X(t)+i\,s\,\Im \tilde X(t)$. By (\ref{eqV.3.7}) (with $s=1$) and
(\ref{eqV.3.11}) (with $s=0$) we have $\bra\Re w\ket\geq \frac
1{M_1}\,\bra\theta-t\ket$.

Moreover in case 1, $\Im \tilde X(t)=\Im x(t,\alpha)=0$ and in case~2,
\begin{equation}\label{eqV.3.27}
\v\Im \tilde X(t)\v\leq M_3\,\frac{\v z\v}{\bra\theta\ket}\leq \frac{3M_5\,\v
z-x(\theta,\alpha)\v}{\bra\theta\ket}\,.
\end{equation}
Summing up we obtain
\begin{equation}\label{eqV.3.28}
\v(3)\v\leq C\,(\varepsilon+\delta)\Big (\frac 1{\bra
t\ket^2}+\frac{M^{2+\sigma_0}_1}{\bra\theta-t\ket^2}\Big )\Big (\v
U(t)\v+M_5\,\frac{\v z-x(\theta,\alpha)\v}{\bra\theta\ket}\Big )\,,
\end{equation}
since $\v\Im X(t)\v\leq \v U(t)\v+\v\Im \tilde X(t)\v$.

For the term (4), due to the factor $\frac 1{\bra t\ket}$ we have a better
estimate. Indeed by (\ref{eqV.3.21}), (\ref{eqV.3.17}), (\ref{eqV.3.6}) and
(\ref{eqV.3.10}) we have
\begin{equation}\label{eqV.3.29}
\v(4)\v\leq C\,(M_1)\,(\varepsilon+\delta)\,\frac{\v
z-x(\theta,\alpha)\v}{\bra\theta\ket}\Big (\frac
1{\bra\theta-t\ket^2}+\frac1{\bra t\ket^2}\Big)\,.
\end{equation}

\noindent (iii) Estimation of (5). We note here that $\frac{\partial
q}{\partial \xi}\,(x,\xi)$ is linear in $\xi$ and real if
$(x,\xi)\in\R^n\times\R^n$. It follows that
$$
(5)=\frac{\partial q}{\partial \xi}\,\big(\Re X(t),\,\Im
(\Phi(t,X(t),\alpha)-\Phi(t,\tilde X(t),\alpha))\big )\,.
$$
Using (\ref{eqV.3.18}) we obtain
\begin{equation*}
\begin{split}
\Im (\Phi (t,X(t),\alpha)&-\Phi (t,\tilde X(t),\alpha))=\frac{\Re (X(t)-\tilde
X(t))}{1+4t^2}+\frac{2t}{1+4t^2}\,\Im (X(t)-\tilde X(t))\\
& -\Im \big[\tilde a(t,X(t),\alpha)-\tilde a(t,\tilde
X(t),\alpha)\big]-\frac1{\bra t\ket}\,\Re \big[\tilde b(t,X(t),\alpha)-\tilde
b(t,\tilde X(t),\alpha)\big]\,.
\end{split}
\end{equation*}
By (\ref{eqV.3.6}) and (\ref{eqV.3.10}) we have $\v X(t)-\tilde X(t)\v\leq
3M_1\,\v z-x(\theta,\alpha)\v\,\frac{\bra t\ket}{\bra\theta\ket}$. Moreover we
can use (\ref{eqV.3.28}) and (\ref{eqV.3.29}). Finally we use the fact that the
coefficients of $\frac{\partial q}{\partial \xi}\,(\Re X(t),\cdots )$ are
bounded by $\frac {C\,\varepsilon}{\bra\Re X(t)\ket^{1+\sigma_0}}$ which by
(\ref{eqV.3.7}) (with $s=1$) and (\ref{eqV.3.11}) (with $s=1$) can be estimated
by $\frac{C\,\varepsilon}{\bra\theta-t\ket^{1+\sigma_0}}$. Gathering these
informations we see that
\begin{multline}\label{eqV.3.30}
\v(5)\v\leq \frac{C\,\varepsilon}{\bra\theta-t\ket^{1+\sigma_0}} \Big
(\frac1{\bra t\ket^2}+\frac1{\bra\theta-t\ket^2}\Big )\,\v
U(t)\v+\frac{C(M_1\,\varepsilon \v
z-x(\theta,\alpha)\v}{\bra\theta-t\ket^{1+\sigma_0}\bra\theta\ket}\\
+C\,\varepsilon\,C(M_1,M_3)(\varepsilon+\delta)\Big (\frac1{\bra
t\ket^2}+\frac1{\bra\theta-t\ket^2}\Big )\,\frac{\v
z-x(\theta,\alpha)\v}{\bra\theta\ket}\,.
\end{multline}

\noindent (iv) Estimation of (6). Since again $\frac{\partial q}{\partial
\xi}\,(x,\xi)$ is linear in $\xi$ and real when $(x,\xi)\in\R^n\times\R^n$ we
can write
\begin{equation*}
\v(6)\v\leq \Big\v\frac{\partial q}{\partial \xi}\,(\Re X(t),\Im \Phi(t,\tilde
X(t),\alpha))\Big\v+\Big\v\frac{\partial q}{\partial \xi}\,(\Re\tilde X(t),\Im
\Phi (t,\tilde X(t),\alpha))\Big\v\,.
\end{equation*}
Since the coefficients of $\frac{\partial q}{\partial \xi}$ are bounded by
$\frac{C\,\varepsilon}{\bra\theta-t\ket^{1+\sigma_0}}$ we obtain
$$
\v(6)\v\leq \frac{C\,\varepsilon}{\bra\theta-t\ket^{1+\sigma_0}}\, \v\Im
\Phi(t,\tilde X(t),\alpha)\v\,.
$$
In the case 1, $\tilde X(t)=x(t,\alpha)$ which implies that $\Im \Phi(t,\tilde
X(t),\alpha)=0$. In the case~2, $\v\Im \Phi(t,\tilde X(t),\alpha)\v\leq
M_3\,\frac{\v x(\theta,\alpha)\v}{\bra\theta\ket}\leq 3\,M_3\,\frac{\v
z-x(\theta,\alpha)\v}{\bra\theta\ket}$. Therefore
\begin{equation}\label{eqV.3.31}
\v(6)\v\leq \frac{C\,M_3\,\varepsilon\,\v
z-x(\theta,\alpha)\v}{\bra\theta-t\ket^{1+\sigma_0}\bra\theta\ket}~.
\end{equation}

\noindent (v) Estimation of (7), (8), (9), (10). Using (\ref{eqV.3.7}) and
(\ref{eqV.3.11}) and the estimates on the coefficients of $q$ we find that
$$
\v(7)+(8)+(8)+(10)\v\leq
\frac{C\,\varepsilon}{\bra\theta-t\ket^{2+\sigma_0}}\,(\v U(t)\v+\v\Im \tilde
X(t)\v)\,.
$$
Using (\ref{eqV.3.28}) we obtain finally
\begin{equation}\label{eqV.3.32}
\v(7)+(8)+(9)+(10)\v\leq
\frac{C\,\varepsilon}{\bra\theta-t\ket^{2+\sigma_0}}\,(\v U(t)\v+\frac{3M_3\,\v
z-x(\theta,\alpha)\v}{\bra\theta\ket}\Big )\,.
\end{equation}
Gathering the estimates given by (\ref{eqV.3.26}) to (\ref{eqV.3.32}) and
taking
$\varepsilon+\delta$ small compared to $M_1$, $M_3$ we obtain the conclusion of
Lemma \ref{lV.3.8}. \cqfd
\begin{lemma}\sl \label{lV.3.9}
Let $Y(t)=(Y,(t),\ldots ,Y_n(t))$ be a $C^1$ function from $[T,\theta]$ to
$\R^n$ which satisfies
$$
\Big\v\dot Y(t)-\frac{4t}{1+4t^2}\,Y(t)\Big\v\leq \v h(t)\v\,\v Y(t)\v+\v
g(t)\v+\frac K{\bra 2t\ket}
$$
for some continuous functions $h,g$ and $K\geq 0$. Then
$$
\v Y(t)\v\leq \Big (\frac{\bra 2t\ket}{\bra 2\theta\ket}\,\v
Y(\theta)\v+\int^\theta_T \v g(s)\v\,ds+K\Big )\,\exp \Big (\int^\theta_T
\v h(s)\v\,ds\Big )\,.
$$
\end{lemma}

\noindent {\bf Proof } Let us set $Z(t)=\frac{Y(t)}{\bra 2t\ket}$. Then $\dot
Z(t)=\frac{\dot Y(t)}{\bra 2t\ket}-\frac{4t}{\bra 2t\ket^3}\,Y(t)$. It follows
that
$$
\vert \dot Z(t)\vert \leq \vert h(t)\vert \,\vert Z(t)\vert +\frac{\vert
g(t)\vert }{\bra 2t\ket}+\frac K{1+4t^2}\,.
$$
Therefore for $\sigma\in [t,\theta]$, $t\geq T$ we have
$$
\v Z(\sigma)\v\leq \v Z(\theta)\v +\int^\theta_\sigma \v h(s)\v\,\v
Z(s)\v\,ds+\int^\theta_t \frac{\v g(s)\v}{\bra 2s\ket}\,ds+K\,\int^\theta_t
\frac{ds}{1+4s^2}\,.
$$
Now we have,
$$
\int^{+\infty }_t \frac{ds}{1+4s^2}\leq \frac 1{\bra 2t\ket}\,,\quad 
\int^\theta_t
\frac{\v g(s)\v}{\bra 2s\ket}\,ds\leq \frac 1{\bra 2t\ket} \int^\theta_t \v
g(s)\v\,ds\,.
$$
Using Gronwall's Lemma we obtain
$$
\v Z(\sigma)\v\leq \Big (\v Z(\theta)\v+\frac 1{\bra 2t\ket} \int^\theta_t \v
g(s)\v\,ds+\frac K{\bra 2t\ket}\Big )\,\exp \Big (\int^\theta_t \v
h(s)\v\,ds\Big )\,.
$$
Taking $t=T$ and $\sigma=t$ we obtain, since $Y(t)={\bra 2t\ket}\v Z(t)\v$, the
claim of the Lemma.\cqfd
\begin{corollary}\sl \label {cV.3.10}
With $U(t)=\Im (X(t)-\tilde X(t))$ introduced in (\ref{eqV.3.24}) we have
$$
\v U(t)\v \leq C\Big (\v U(\theta)\v +(6\,M_1+C)\,\frac{\v
z-x(\theta,\alpha)\v}{\bra \theta\ket}\Big )\,.
$$
\end{corollary}

\noindent {\bf Proof } This follows from Lemmas \ref{lV.3.8} and
\ref{lV.3.9}. 
\cqfd

We can now finish the proof of the improvement of (\ref{eqV.3.9}).

Indeed we have $U(\theta)=\Im (X(\theta,\theta,z)-\tilde
X(\theta,\theta,0))=\Im z$. Therefore Corollary \ref{cV.3.10} and Remark
\ref{rV.3.4} show that if $C<M_1$ and $(6\,M_1+C)\cdot C<M_3$ then
(\ref{eqV.3.9}) is improved.

\noindent {\bf End of the proof of Theorem \ref{tV.3.2} in the cases 1 and 2}

The estimates (\ref{eqV.3.5}) to (\ref{eqV.3.12}) improved are true for
$t\in[T,\theta]$ for all $T>T_*$. By continuity they continue to hold on
$[T_*,\theta]$. Now we consider problem (\ref{eqV.3.5}) with data at $t=T_*$
equal to $X(T_*,\theta,\alpha)$. For this problem the estimates (\ref{eqV.3.6})
to (\ref{eqV.3.12}) hold on $[T_*-\varepsilon_0,T_*]$ which contradicts the
fact that $T_*=\inf A$. Therefore $A=[0,\theta]$ which implies Theorem
\ref{tV.3.2} in this case.

\noindent {\bf 2) Proof of Theorem \ref{tV.3.2} in case 3}

Here we shall take $x\in\R^n$ such that $x\cdot \alpha_\xi\leq
c_1\,\bra x\ket \v\alpha_\xi\v$ and $\v x-x(\theta,\alpha)\v\leq
\delta_1\,\bra\theta\ket$, with $0<c_1\ll c_2$, $0<\delta_1\ll \delta_2$.

Let us recall (see (\ref{eqIV.4.49})) that there exists a unique
$\theta^*\in[0,\theta]$ such that $x(\theta^*,\alpha)\cdot
\alpha_\xi=0$. We shall make use of Lemma \ref{lIV.4.17}. To prove the
claim of Theorem \ref{tV.3.2}  we shall use the same method as in the
cases~1 and 2.

We introduce first the set $A$ of $T\geq \theta^*$ such that the problem
(\ref{eqV.3.4}) has a solution on $[T,\theta]$ which satisfies
\begin{equation}\label{eqV.3.33}
\v X(t,\theta,x)-x\v\leq M_4\,\v t-\theta\v
\end{equation}
\begin{equation}\label{eqV.3.34}
\Re X(t,\theta,x)\cdot \alpha_\xi\leq M_5\,\bra\Re
X(t,\theta,x)\ket\,\v\alpha_\varepsilon\v
\end{equation}
\begin{equation}\label{eqV.3.35}
\bra x\ket+\bra t-\theta\ket\leq M_4\,\bra\Re X(t,\theta,x)\ket
\end{equation}
\begin{equation}\label{eqV.3.36}
\v\Im X(t,\theta,x)\v\leq M_6\,\frac{\v x-x(\theta,\alpha)\v}{\bra\theta\ket}\,.
\end{equation}
If $M_4$ is large enough, $M_5>c_1$, $M_6>0$ one can find $\varepsilon_0>0$
such that $\theta-\varepsilon_0\in A$. Let $T_*=\inf A$. We want to prove that
$T_*=\theta^*$ if $M_4,M_5,M_6$ are correctly chosen; let us assume
$T_*>\theta$ and let $t\geq T_*$. Then on $[T,\theta]$ we have a solution
$X(t,\theta,z)$ which satisfies (\ref{eqV.3.33}) to (\ref{eqV.3.36}). Let us
show that this implies that $(t,X(t,\theta,x))\in\Omega^{\C}_\delta$ for
$t\in[T,\theta]$ (see (\ref{eqV.3.3})) if $\delta_1$ is small enough.

From (\ref{eqV.3.36}) we have $\v\Im X(t,\theta,z)\v\leq M_6\,\delta_1\leq
\delta$ if $\delta_1$ is small enough. Moreover by (\ref{eqV.3.34}) we have
$$
\Re X(t,\theta,x)\cdot \alpha_\xi\leq M_5\,\bra\Re
X(t,\theta,x)\ket\,\v\alpha_\varepsilon\v\leq c_0 \bra\Re
X(t,\theta,x)\ket\,\v\alpha_\varepsilon\v
$$
if $M_5\leq c_0$. Finally,
 $$
(1)=\v X(t,\theta,x)-x(t,\alpha)\v\leq \v X(t,\theta,x)-x\v+\v
x-x(\theta^*,\alpha)\v+\v x(\theta^*,\alpha)-x(t,\alpha)\v\,.
$$
From (\ref{eqV.3.33}) we have $\v X(t,\theta,x)-x\v\leq M_4\,\v t-\theta\v\leq
M_4(\theta-\theta^*)$ since $t\geq \theta^*$. Now we use Lemma \ref{lIV.4.17} 
to write
$$
\v X(t,\theta,x)-x\v\leq 10\,M_4\,\v x-x(\theta,\alpha)\v\leq
10\,M_4\,\delta_1\bra\theta\ket\leq \frac{10\,M_4}{K_1}\,\delta_1
\bra\theta^*\ket\leq \frac{10\,M_4}{K_1}\,\delta_1 \bra t\ket\,.
$$
Again by Lemma \ref{lIV.4.17},
$$
\v x-x(\theta^*,\alpha)\v\leq 6 \v x-x(\theta,\alpha)\v\leq 6
\delta_1\bra\theta\ket\leq \frac{6\delta_1}{K_1}\,\bra\theta^*\ket\leq
\frac{6\delta_1}{K_1}\,\bra t\ket\,.
$$
Finally $\v x(t,\alpha)-x(\theta^*,\alpha)\v\leq \int^t_{\theta^*} \v\dot
x(s,\alpha)\v\,ds\leq 5(t-\theta^*)$ if $\varepsilon$ is small enough. It
follows from Lemma \ref{lIV.4.17} that
$$
\v x(t,\alpha)-x(\theta^*,\alpha)\v\leq 5 (\theta-\theta^*)\leq 50\, \v
x-x(\theta,\alpha)\v\leq 50\,\delta_1 \bra\theta\ket\leq
\frac{50\,\delta_1}{K_1} \,\bra t\ket\,.
$$
Summing up we find that if $\delta_1$ is small enough,
\begin{equation}\label{eqV.3.37}
(1)\leq \max \Big (\frac{10 M_4}{K_1}\,,\,\frac{56}{K_1}\Big )\,\delta _1\bra
t\ket\leq \delta \bra t\ket\,.
\end{equation}
As in the cases 1 and 2 our goal is to prove that one can improve the estimates
(\ref{eqV.3.33}) to (\ref{eqV.3.36}).

(i) {\bf Improvement of (\ref{eqV.3.33})}

We have by (\ref{eqV.3.5}), $\dot X(t,\theta,x)=2\alpha_\xi+\CO
(\varepsilon+\delta)$. Therefore
\begin{equation}\label{eqV.3.38}
X(t,\theta,x)=x-2(\theta-t)\,\alpha_\xi+\CO
((\varepsilon+\delta)(\theta-t))\,.
\end{equation}
It follows that $\v X(t,\theta,x)-x\v\leq 5 (\theta-t)$ if $\varepsilon+\delta$
is small enough. We shall take $M_4$ so that $5\leq \frac 12\,M_4$ and then,
(\ref{eqV.3.33}) will be improved.

(ii) {\bf Improvement of (\ref{eqV.3.35})}

We deduce from (\ref{eqV.3.38}) that
$$
1+\v\Re X(t,\theta,z)\v^2=1+\v
x\v^2+4(\theta-t)^2\,\v\alpha_\xi\v^2+\CO ((\varepsilon+\delta)(\v
x\v^2+(\theta-t)^2))-2(\theta-t)\,x\cdot \alpha_\xi\,.
$$
Since $x\cdot \alpha_\varepsilon\leq c_1\bra x\ket\v\alpha_\varepsilon\v$,
taking $\varepsilon+\delta$ small enough we obtain, $1+\v\Re
X(t,\theta,z)\v^2\geq 1+\frac 12 \,\v x\v^2+\frac 12
\,(\theta-t)^2-4c_1(\theta-t)\bra x\ket$, so
\begin{equation}\label{eqV.3.39}
1+\v\Re X(t,\theta,z)\v^2\geq \frac 14 (\bra x\ket^2+(\theta-t)^2)\,,
\end{equation}
if $c_1\leq \frac 1{10}$. In particular $3\bra \Re X(t,\theta,x)\ket\geq
\bra\theta-t\ket+\bra x\ket$, so $\bra\theta-t\ket+\bra x\ket\leq \frac 12
\,M_4\,\bra\Re X(t,\theta,x)\ket$ if $M_4\geq 6$.

(iii) {\bf Improvement of (\ref{eqV.3.34})}

From (\ref{eqV.3.38}) we have
$$
\Re X(t,\theta,x)\cdot \alpha_\xi=x\cdot
\alpha_\xi-2(\theta-t)\,\v
\alpha_\xi\v^2+\CO((\varepsilon+\delta)(\theta-t))\,.
$$
It follows that
$$
\Re X(t,\theta,x)\cdot \alpha_\xi\leq c_1 \bra
x\ket\,\v\alpha_\xi\v-\frac{(\theta-t)}4\leq 10\,c_1 \,\bra\Re
X(t,\theta,x)\ket\,\v\alpha_\xi\v\,,
$$
by (\ref{eqV.3.39}). We shall take $10\,c_1\leq \frac 12 \, M_5$ and
(\ref{eqV.3.34}) will be improved.

(iv) {\bf Improvement of (\ref{eqV.3.36})}

Let us set $X(t,\theta,x)=X(t)=Y_1(t)+i\,Y_2(t)$ where $Y_1,Y_2$ are real.
\begin{lemma}\sl \label{lV.3.11}
There exists positive constants $C,K$ independent of $\varepsilon,\delta$ and
$T$ such that for all $t\in[T,\theta]$ we have
$$
\v\dot Y_2(t)\v\leq C\,M_6 \Big
(\frac{\delta_1}{\bra\theta\ket}+(\varepsilon+\delta)\,\frac{\v
x-x(\theta,\alpha)\v}{{\bra\theta\ket}}\,g(t)\Big
)+\frac{\delta_1}{\bra\theta\ket}
$$
where $g$ is a continuous positive function satisfying $\int^\theta_T
g(s)\,ds\leq K$.
\end{lemma}

{\bf Proof } From (\ref{eqV.3.36}) and (\ref{eqV.3.37}) we get
\begin{equation}\label{eqV.3.40}
\v X(t)-x(t,\alpha)\v\leq C\,(M_4)\,\delta_1\,\bra t\ket\,,\quad \v
Y_2(t)\v\leq M_6\,\delta_1\,.
\end{equation}
Now (\ref{eqV.3.18}) shows that
$$
\Im \Phi(t,X(t),\alpha)=\Im \Big[\frac{X(t)-x(t)}{2t-i}-\tilde
a(t,X(t,X(t),\alpha)-\frac{i}{\bra t\ket}\,\tilde b(t,X(t),\alpha)\Big]
$$
where $x(t)=x(t,\alpha)$. First of all we have,
$$
\Im\,
\frac{X(t)-x(t)}{2t-i}=\frac{2t\,Y_2(t)}{1+4t^2}+\frac{Y_1(t)-x(t)}{1+4t^2}\,.
$$
Using (\ref{eqV.3.40}) we deduce, since $\bra\theta\ket\sim\bra t\ket$,
\begin{equation}\label{eqV.3.41}
\Big\v\Im\,\frac{X(t)-x(t)}{2t-i}\Big\v\leq
C\,M_6\,\frac{\delta_1}{\bra\theta\ket}+C\,\frac{\delta_1}{\bra\theta\ket}\,.
\end{equation}
On the other hand we can write with $f=\tilde a$ or $\tilde b$,
\begin{equation*}
\begin{split}
f(t,X(t),\alpha)&=f(t,Y_1(t),\alpha)+\int^1_0 \frac{\partial f}{\partial
z}\,(t,Y_1(t)+is\,Y_2(t),\alpha)\,ds(i\,Y_2(t))\\
&\quad + \int^1_0 \frac{\partial f}{\partial \overline z}
\,(t,Y_1(t)+is\,Y_2(t),\alpha)\,ds\,(-i\,Y_2(t))\,.
\end{split}
\end{equation*}
Since $\tilde a(t,Y_1(t),\alpha)$ is real, using the estimates on the
derivatives of $\tilde a$ and $\tilde b$ given by Theorem \ref{tIV.4.2} and
their extensions to the complex domain proved in Lemma \ref{lV.3.1} we obtain
\begin{equation}\label{eqV.3.42}
\v\Im \tilde a(t,X(t),\alpha)\v\leq C\,M_6\,(\varepsilon+\delta)\Big (\frac
1{\bra\theta-t\ket^2}+\frac 1{\bra t\ket^2}\Big )\,\frac{\v
x-x(\theta,\alpha)\v}{\bra\theta\ket}\,.
\end{equation}
Here we have used the estimate in (ii) and (\ref{eqV.3.36}).

Moreover we have by Theorem \ref{tIV.4.2} (ii) and (\ref{eqV.3.37}),
$$
\v\tilde b(t,Y_1(t),\alpha)\v\leq \sqrt\delta\,\frac{\v Y_1(t)-x(t)\v}{\bra
t\ket}\leq \sqrt\delta\,C(M_4)\,\delta_1\leq \delta_1
$$
if $\delta$ is small enough. Therefore
\begin{equation}\label{eqV.3.43}
\Big\v \Im\, \frac{i\,\tilde b(t,X(t),\alpha)}{\bra t\ket}\Big\v\leq
C\,\frac{\delta_1}{\bra\theta\ket}+C\,M_6 \Big (\frac
1{\bra\theta-t\ket^2}+\frac 1{\bra t\ket^2}\Big )\,\frac{\v
x-x(\theta,\alpha)\v}{\bra\theta\ket}\,.
\end{equation}
We deduce from (\ref{eqV.3.41}) to (\ref{eqV.3.43}) that
\begin{equation}\label{eqV.3.44}
\v\Im \Phi(t,X(t),\alpha)\v\leq
\frac{C\,\delta_1}{\bra\theta\ket}+C\,M_6\,\frac{\delta_1}{\bra\theta\ket}
+M_6\,\frac{\v
x-x(\theta,\alpha)\v}{\bra\theta\ket}\,g(t)
\end{equation}
where $g(t)=C\big (\frac 1{\bra\theta-t\ket^2}+\frac 1{\bra t\ket^2}\big)$.

It follows from (\ref{eqV.3.4}) that
\begin{equation*}
\begin{split}
\v \dot Y_2(t)\v&\leq 2\,\v\Im \Phi(t,X(t),\alpha)\v+\Big\v\frac{\partial
q}{\partial \xi}\,(Y_1(t),\Im \Phi(t,X(t),\alpha))\Big\v\\
&\qquad +\Big\v\int^1_0 \frac{\partial ^2q}{\partial \xi\,\partial
z}\,(Y_1(t)+is\,Y_2(t),\Phi(t,X(t),\alpha))\,ds\Big\v\,\v
Y_2(t)\v+\Big\v\int^1_0 \frac{\partial ^2q}{\partial \xi\partial \overline
z}\,\textrm{ (idem) }ds\Big\v\,\v Y_2(t)\v\\
& \leq C\,\v\Im \Phi(t,X(t),\alpha\v+\frac{C
\varepsilon}{\bra\theta-t\ket^{2+\sigma_0}}\,\v Y_2(t)\v\,.
\end{split}
\end{equation*}
This estimate together with  (\ref{eqV.3.36}), (\ref{eqV.3.44}) prove the Lemma.
\cqfd

We can now improve (\ref{eqV.3.36}). Indeed, by Lemma \ref{lV.3.11}, we have,
since $X(\theta,\theta,x)=x$ is real
$$
\v Y_2(t)\v\leq \int^\theta_t \v\dot Y_2(s)\v\,ds\leq C\,M_6
\Big[\frac{\delta_1(\theta-t)}{\bra\theta\ket}
+(\varepsilon+\delta)\,K\,\frac{\v
x-x(\theta,\alpha)\v}{\bra\theta\ket}\Big]+\delta_1\,\frac{\theta-t}{\bra\theta\ket}\,.
$$
Moreover by Lemma \ref{lIV.4.17} we have
$$
\theta-t\leq \theta-\theta^*\leq C\,\v x-x(\theta,\alpha)\v\,.
$$
Taking $\delta_1, \delta,\varepsilon$ small enough we obtain $\v Y_2(t)\v\leq
\frac 12\,M_6\,\frac{\v x-x(\theta,\alpha)\v}{\bra\theta\ket}$ which improves
(\ref{eqV.3.36}).

The improvements (i) to (iv) show that the set $A$ where (\ref{eqV.3.33}) to
(\ref{eqV.3.36}) are true is equal to $[\theta^*,\theta]$.

We can now give the proof of Theorem \ref{tV.3.2} in the case~3. Indeed
(\ref{eqV.3.34}) to (\ref{eqV.3.36}) imply the estimates (ii) to (iv) in this
Theorem. To prove (i) we just remark that
\begin{equation*}
\begin{aligned}
\v X(t,\theta,x)-x(t,\alpha)\v&\leq \v X(t,\theta,x)-x\v+\v
x-x(\theta,\alpha)\v\\
&\leq (10M_4+1)\v x-x(\theta,\alpha)\v\leq C\,\v
x-x(\theta,\alpha)\v\,\frac{\bra t\ket}{\bra\theta\ket}
\end{aligned}
\end{equation*}
since $\bra t\ket\sim\bra\theta\ket$ when $t\in[\theta^*,\theta]$.
Therefore we are done for $t\in[\theta^*,\theta]$. For $t\in[0,\theta^*]$ we
first remark that
$$
X(t,\theta,x)=X(t,\theta^*,X(\theta^*,\theta,x))\,.
$$
We would like to apply the cases 1 and 2 already done, with $\theta=\theta^*$
and $z=X(\theta^*,\theta,x)$. So we have to prove that
\begin{itemize}
\item[(i)] $x(\theta^*,\alpha)\cdot \alpha_\xi\leq c_2\,\bra
x(\theta^*,\alpha)\ket\,\v\alpha_\xi\v$,
\item[(ii)] $\v z-x(\theta^*,\alpha)\v\leq \delta_2\,\bra\theta^*\ket$,
\item[(iii)] $\Re z\cdot \alpha_\xi\leq c_2\,\bra\Re
z\ket\,\v\alpha_\xi\v$,
\item[(iv)] $\v\Im z\v\leq \delta_2$.
\end{itemize}
First of all (i) is trivial since $x(\theta^*,\alpha)\cdot
\alpha_\varepsilon=0$. Now we have,
$$
\v X(\theta^*,\theta,x)-x(\theta^*,\alpha)\v\leq \v X(\theta^*,\theta,x)-x\v+\v
x-x(\theta^*,\alpha)\v=(1)+(2)\,.
$$
By (\ref{eqV.3.33}) and Lemma \ref{lIV.4.17} we have if $\delta_1\ll\delta_2$
\begin{equation*}
\begin{aligned}
(1)&\leq M_4(\theta-\theta^*)\leq 10\,M_4\,\v x-x(\theta,\alpha)\v\leq
10\,M_4\,\delta_1 \bra\theta\ket\leq C\,\delta_1 \bra\theta^*\ket\leq
\frac{\delta_2}2\,\bra\theta^*\ket\\
(2)&\leq 5\,\v x-x(\theta,\alpha)\v\leq C'\,\delta_1 \bra\theta^*\ket\leq
\frac{\delta_2}2\,\bra\theta^*\ket
\end{aligned}
\end{equation*}
since $\bra\theta\ket\sim\bra\theta^*\ket$. Thus (ii) is satisfied. Now (iii)
is also satisfied if $M_5\leq c_2$. This is possible since  the only constraint
on $M_5$ (see (iii) improvement of (\ref{eqV.3.34})) was $M_5\geq 20\,c_1$.
Finally by (\ref{eqV.3.36}), $\v\Im X(\theta^*,\theta,x)\v\leq
\delta_1\,M_6\leq \delta_2$ if $\delta_1\ll\delta_2$. Therefore
$X(t,\theta^*,X(\theta^*,\theta,x))$ satisfies the estimates (\ref{eqV.3.5}) to
(\ref{eqV.3.9}) in case~1 and (\ref{eqV.3.9}) to (\ref{eqV.3.12}) in case~2.
Therefore we have the following estimate,
\begin{equation*}
\begin{aligned}
(1)&=\v X(t,\theta^*,X(\theta^*,\theta,x))-x(t,\alpha)\v\leq 3\,M_1\,\v
X(\theta^*,\theta,x)-x(\theta^*,\alpha)\v\,\frac{\bra t\ket}{\bra\theta^*\ket}\\
(1)&\leq 3\,M_1\,\frac{\bra t\ket}{\bra\theta^*\ket}\,(\v
X(\theta^*,\theta,x)-x\v+\v x-x(\theta^*,\alpha)\v)\\
(1)&\leq 3\,M_1\,\frac{\bra t\ket}{\bra\theta^*\ket}\, (M_4(\theta-\theta^*)+\v
x-x(\theta^*,\alpha)\v)\\
(1)&\leq C\,M_1 (1+M_4)\v x-x(\theta,\alpha)\v\,\frac{\bra
t\ket}{\bra\theta\ket}\,.
\end{aligned}
\end{equation*}
Here we have used (\ref{eqV.3.33}), Lemma \ref{lIV.4.17} and
$\bra\theta^*\ket\sim\bra\theta\ket$. Therefore we obtain (i) of
(\ref{eqV.3.5}) if $K\geq C\,M_1(1+M_4)$. Now (\ref{eqV.3.9}) implies that
$$
(2)=\v\Im X(t,\theta^*,X(\theta^*,\theta,x))\v\leq M_3 \Big (\frac{\v
X(\theta^*,\theta,x)-x(\theta^*,\alpha)\v}{\bra\theta^*\ket}+\v\Im
X(\theta^*,\theta,x)\v\Big)\,.
$$
Using (\ref{eqV.3.36}) and the same argument as in the term (1) we obtain,
$(2)\leq C(M_1,M_4,M_6)\,\frac{\v x-x(\theta,\alpha)\v}{\bra\theta\ket}$. Thus
(ii) in (\ref{eqV.3.5}) is satisfied if $K\geq C(M_1,M_4,M_6)$.

From (\ref{eqV.3.15}) (ii) with $s=1$ we have
$$
(3)=\bra\Re X(t,\theta^*,X(\theta^*,\theta,x))\ket\geq \frac 15\,
[\bra\theta^*-t\ket+\bra\Re X(\theta^*,\theta,x)\ket]\,.
$$
So using (\ref{eqV.3.35}) we obtain
$$
(3)\geq \frac 15 \Big[\bra\theta^*-t\ket+\frac
1{M_4}\,\bra\theta^*-\theta\ket\Big]\geq C(M_4)\bra\theta-t\ket\,,
$$
and (iii) satisfied $K\cdot C(M_4)\geq 1$.

Finally let us set $(4)=\Re X(t,\theta^*,X(\theta^*,\theta,x))\cdot
\alpha_\xi$. Using (\ref{eqV.3.8}) and (\ref{eqV.3.12}) with $s=1$ we
can write,
$$
(4)\leq M_2\,\bra\Re X(t,\theta^*,X(\theta^*,\theta,x))\ket\,.
$$
This shows that (\ref{eqV.3.5}) (iv) holds if $\tilde K\geq \frac 1{M_2}$ and
completes the proof of Theorem \ref{tV.3.2}. \cqfd

Having proved in Theorem \ref{tV.3.2} the existence of the solution
$X(t,\theta,x)$ of (\ref{eqV.3.4}) we want to give estimates on its derivatives
with respect to $(\theta,x)$.
\begin{proposition}\sl \label{pV.3.12}
The solution given by Theorem \ref{tV.3.2} is $C^\infty $ with respect to
$y=(\theta,x)$ and satisfies the following estimates,
\begin{equation}\label{eqV.3.45}
\v\partial ^A_y\,X(t,\theta,x)\v\leq 
\begin{cases}
C\,\frac{\bra t\ket}{\bra\theta\ket}&\textrm{if \enskip $\v A\v=1$}\,, \\
C_A\,\frac{\bra t\ket}{\bra\theta\ket} \big (\frac 1{\bra
x\ket^{\v A\v+\sigma_0}}+\frac 1{\bra\theta\ket^{\v
A\v-1}}\big)&\textrm{if\enskip $\v A\v\geq 2$}\,.
\end{cases}\,.
\end{equation}
uniformly in  $(t,\theta,x)\in [0,\theta]\times\Omega_\delta$.
\end{proposition}
To prove this result we need a Lemma.
\begin{lemma}\sl\label{lV.3.13}
For $j=1,\ldots ,n$ let us set $L_j(t,z)=\frac{\partial p}{\partial
\xi_j}\,(z,\Phi(t,z,\alpha))$. Then for any integer $N>0$ one can find
$C_N>0$ such that for $j=1,\ldots ,n$ and all $(t,\theta,x)$ in
$[0,\theta]\times\Omega_\delta$ we have
\begin{equation}\label{eqV.3.46}
\Big\v\frac{\partial L_j}{\partial
z_k}\,(t,X(t,\theta,x))-\frac{2\delta_{jk}}{2t-i}\Big\v\leq C_1 \Big
(\frac 1{(\bra x\ket+\bra\theta-t\ket)^{2+\sigma_0}}+\frac 1{\bra t\ket^2}\Big)
\end{equation}
\begin{equation}\label{eqV.3.47}
\Big\v\frac{\partial L_j}{\partial
\overline z_k}\,(t,X(t,\theta,x))\Big\v\leq C_N \Big
(\frac 1{(\bra x\ket+\bra \theta-t\ket)^{2+\sigma_0}}+\frac 1{\bra
t\ket^2}\Big)\Big (\frac{\v
x-x(\theta,\alpha)\v}{\bra\theta\ket}\Big)^N\,.
\end{equation}
For any $\mu,\nu\in\N^n$,  such that $\v\mu\v+\v\nu\v=N\geq 2$, $j=1,\ldots
,n$,
\begin{equation}\label{eqV.3.48}
\Big\v\frac{\partial ^{\mu+\nu}\,L_j}{\partial z^\mu\,\partial \overline
z^{\nu}}\,(t,X(t,\theta,x))\Big\v\leq C_{\mu\nu}\Big (\frac 1{(\bra
x\ket+\bra\theta-t\ket)^{\v\mu\v+\v\nu\v+1+\sigma_0}}+\frac 1 {\bra
t\ket^{\v\mu\v+\v\nu\v+1}}\Big)\,.
\end{equation}
\end{lemma}

{\bf Proof of Lemma \ref{lV.3.13} } We have
$$
L_j(t,z)=2\,\Phi_j(t,z,\alpha)+2\varepsilon\,\sum^n_{\ell =1} b_{j\ell
}(z)\,\Phi_\ell (t,z,\alpha)\,.
$$
Using (\ref{eqV.3.18}) we obtain
$$
\frac{\partial \Phi_j}{\partial
z_k}\,(t,z,\alpha)=\frac{\delta_{jk}}{2t-i}-\Big (\frac{\partial \tilde
a_j}{\partial z_k}+\frac i{\bra t\ket}\,\frac{\partial \tilde b_j}{\partial
z_k}\Big )\,(t,z,\alpha)\,.
$$
Then (\ref{eqV.3.46}) follows easily from the estimates on $\tilde a,\tilde b$
given in Theorem \ref{tIV.4.2}, the estimates on the coefficients $b_{j\ell }$
and from the inequality (iii) in Theorem \ref{tV.3.2}.

The estimate (\ref{eqV.3.47}) follows from the same arguments and Lemma
\ref{lV.3.1} (iv), Theorem \ref{tV.3.2} (ii), (iii). The same method can also
be used to prove (\ref{eqV.3.48}). \cqfd

\noindent {\bf Proof of Proposition \ref{pV.3.12} } Let us set for $k=1,\ldots
,n$, $q\geq 1$,
$$
Y^q_k(t)=\partial ^A_x\, X_k(t,\theta,x)
$$
where  $\v A\v=q$.

We begin by the case $q=1$. Differentiating one time (\ref{eqV.3.4}) with
respect to $y$ we obtain
\begin{equation}\label{eqV.3.49}
\dot Y^1_k(t)=\sum^n_{j=1} \Big [\frac{\partial L_k}{\partial
z_j}\,(t,X(t,\theta,x))\,Y^1_j(t)+\frac{\partial L_k}{\partial \overline
z_j}\,(t,X(t,\theta,x))\,\overline {Y^1_j(t)}\Big]\,.
\end{equation}
Using (\ref{eqV.3.46}) and (\ref{eqV.3.47}) we see that
$Y^1(t)=(Y^1_1(t),\ldots ,Y^1_n(t))$ satisfies the hypotheses of Lemma
\ref{lV.3.6} with $g\equiv 0$ and $h(t)=\frac 1{(\bra
x\ket+\bra\theta-t\ket)^{2+\sigma_0}}+\frac 1{\bra t\ket^2}$. Since
$\frac{\partial X_j}{\partial x_k}\,(\theta,\theta,x)=\delta_{jk}$ and $\frac
{\partial X_j}{\partial \theta}\,(\theta,\theta,x)$ is bounded we obtain
(\ref{eqV.3.45}) when $\v A\v=1$. Let us consider the case $\v A\v=2$.
Differentiating (\ref{eqV.3.49}) with respect to $y$ we see that $Y^2_k(t)$
satisfies the equation 
$$
\dot Y^2_k(t)=\sum^n_{j=1} \Big[\frac{\partial L_k}{\partial
z_j}\,(t,X(t,\theta,x))\,Y^2_j(t)+\frac{\partial L_k}{\partial \overline
z_j}\,(t,X(t,\theta,x))\,\overline {Y^2_j(t)}\Big]+Z_k(t,\theta,x)
$$
where, by (\ref{eqV.3.48}) and (\ref{eqV.3.45}) for $\v A\v=1$, $Z_k(t)$ is
estimated as follows
$$
\v Z_k(t,\theta,x)\v\leq C\,\frac{\bra t\ket^2}{\bra\theta\ket^2} \Big (\frac
1{(\bra x\ket+\bra\theta-t\ket)^{3+\sigma_0}}+\frac 1 {\bra t\ket^3}\Big )\,.
$$
We want to use Lemma \ref{lV.3.6} (with $T=0$) so we are led to estimate the
quantity $(1)=\int^\theta_0 \frac{\v Z_k(s,\theta,x)\v}{\bra 2s\ket}\,ds$.
Using the above estimation we see that
$$
(1)\leq C'\,\frac 1{\bra\theta\ket^2} \int^\theta_0 \Big
(\frac{\bra\sigma\ket}{(\bra x\ket+\bra\theta-\sigma\ket)^{3+\sigma_0}}+\frac
1{\bra\sigma\ket^2}\Big )\,d\sigma\,.
$$
By a straightforward computation we see that we have
\begin{equation}\label{eqV.3.50}
\int^\theta_0 \frac{\bra\sigma\ket^\ell }{(\bra
x\ket+\bra\theta-t\ket)^{k+\sigma_0}}\,d\sigma\leq
C\,\frac{\bra\theta\ket^\ell }{\bra x\ket^{k-1+\sigma_0}}\,,\quad k,\ell \geq
1\,.
\end{equation}
It follows that
$$
(1)\leq \frac C{\bra\theta\ket} \Big (\frac 1{\bra x\ket^{2+\sigma_0}}+\frac
1{\bra\theta\ket}\Big)\,.
$$
Using Lemma \ref{lV.3.6} and the fact that  $\partial
^A_x\,X(\theta,\theta,x)=0$, since $\v A\v\geq 2$, we obtain (\ref{eqV.3.45})
when $\v A\v=2$.

Now we proceed by induction on $q\geq 2$. Let $\v A\v=q+1$ and let us
differentiate the equation $\dot X_k(t,\theta,x)=L_k(t,X(t,\theta,x))$ $\v
A\v$ times with respect to $x$. Using the Faa di Bruno formula and the
notation $Y_k^{q+1}=\partial ^A_x\,X_k$ we obtain the equation
$$
\dot Y^{q+1}_k(t)=\sum^n_{j=1} \Big[\frac{\partial L_k}{\partial
z_j}\,(t,X(t,\theta,x))\,Y^{q+1}_j(t)+\frac{\partial L_k}{\partial \overline
z_j}\,(t,X(t,\theta,x))\Big]\,\overline {Y^{q+1}_j(t)}+Z_k(t)
$$
where $Z_k(t)$ is a finite linear combination of terms of the form
$$
(2)=\big (\partial ^\beta_{(z,\overline z)}\,L_j\big )(t,X(t,\theta,x))
\prod^s_{\ell =1} \big (\partial ^{L_\ell }_x\,X(t,\theta,x)\big )^{K_\ell }
$$
where $2\leq \v \beta\v\leq q+1$, $1\leq s\leq q+1$, $\v K_\ell \v\geq 1$, $\v
L_\ell \v\geq 1$, $\som^s_{\ell =1} K_\ell=\beta$ , $\som^s_{\ell =1} \v
K_\ell \v\,L_\ell =A$.
It follows that $\v L_\ell \v\leq \v A\v-1=q$.

Since by (\ref{eqV.3.45}) we have different estimates for $\v L_\ell\v =1$ and
$\v L_\ell \v\geq 2$ we must separate these two cases. So let us write
$\{1,\ldots ,s\}=I_1\cup I_2$, $I_1=\{\ell :\v L_\ell \v=1\}$, $\v
I_2\v=\{\ell  : \v L_\ell \v\geq 2\}$.

Now let us use (\ref{eqV.3.48}) and the induction. We obtain
$$
\v(2)\v\leq C\Big (\frac 1{(\bra
x\ket+\bra\theta-t\ket)^{\v\beta\v+1+\sigma_0}}+\frac 1{\bra
t\ket^{\v\beta\v+1}}\Big ) \prod_{\ell \in I_1} \Big (\frac{\bra t\ket}{\bra
\theta\ket}\Big )^{\v K_\ell \v} \prod_{\ell \in I_2}
 \bigg[\frac{\bra t\ket}{\bra\theta\ket}\Big (\frac 1{\bra x\ket^{\v L_\ell \v+\sigma_0}}+\frac 1{\bra\theta\ket^{\v
L_\ell \v-1}}\Big )\bigg]^{\v K_\ell \v}\,.
$$
Since $\som^s_{i=1} \v K_i\v=\v\beta\v$ we have  $\prodi_{\ell \in I_1} \big
(\frac{\bra t\ket}{\bra\theta\ket}\big )^{\v K_\ell \v}\prodi_{\ell \in I_2}
\big (\frac{\bra t\ket}{\bra\theta\ket}\big )^{\v K_\ell \v}=\frac{\bra
t\ket^{\v\beta\v}}{\bra\theta\ket^{\v\beta\v}}$. It follows from
(\ref{eqV.3.50}) that
$$
\int^\theta_0 \frac{\v (2)\v}{\bra\sigma\ket}\,d\sigma\leq \frac
C{\bra\theta\ket}\Big (\frac 1{\bra x\ket^{\v\beta\v+\sigma_0}}+\frac
1{\bra\theta\ket^{\v \beta\v-1}}\Big )\prod_{\ell \in I_2}\Big (\frac 1{\bra
x\ket^{\v L_\ell\v +\sigma_0}}+\frac 1{\bra \theta\ket^{\v L_\ell \v-1}}\Big
)^{\v K_\ell \v}\,.
$$
Now we have
\begin{equation*}
\begin{split}
\prod_{\ell \in I_2}\Big (\frac 1{\bra x\ket^{\v L_\ell \v+\sigma_0}}&+\frac
1{\bra\theta\ket^{\v L_\ell \v-1}}\Big )^{\v K_\ell \v}\leq \prod_{\ell \in
I_2} \Big (\frac 1{\bra x\ket^{\v L_\ell \v-1}}+\frac
1{\bra\theta\ket^{\v L_\ell \v-1}}\Big )^{\v K_\ell \v}\\
&\leq C\prod_{\ell \in I_2} \Big (\frac 1{\bra x\ket^{\v K_\ell \v(\v L_\ell
\v-1)}}+\frac 1{\bra\theta\ket^{\v K_\ell \v(\v L_\ell \v-1)}}\Big )\\
&\leq C' \left[\Big (\frac 1{\bra x\ket}\Big)^{\som_{\ell \in I_2}\v K_\ell
\v(\v L_\ell \v-1)}+\Big (\frac 1{\bra\theta\ket}\Big )^{\som_{\ell \in
I_2}\v K_\ell\v(\v L_\ell \v-1)}\right]\\
&\leq C' \left[\Big (\frac 1{\bra x\ket}\Big )^{\v A\v-\v \beta\v}+\Big (\frac
1{\bra\theta\ket}\Big )^{\v A\v-\v\beta\v}\right]\,.
\end{split}
\end{equation*}
Indeed
$$\v A\v-\v\beta\v=\som^s_{\ell =1} \v K_\ell \v\,\v L_\ell
\v-\som^s_{\ell =1} \v K_\ell \v=\som_{\ell \in I_1} \v K_\ell \v+\som_{\ell
\in I_2} \v K_\ell \v\,\v L_\ell \v-\som_{\ell \in I_1} \v K_\ell
\v-\som_{\ell \in I_2}\v K_\ell \v=\som_{\ell \in I_2} \v K_\ell \v (\v L_\ell
\v-1)\,.
$$
It follows that
\begin{equation*}
\begin{split}
\int^\theta_0 \frac{\v(2)\v}{\bra\sigma\ket}\,d\sigma&\leq \frac
C{\bra\theta\ket} \Big (\frac 1{\bra x\ket^{\v\beta\v+\sigma_0}}+\frac
1{\bra\theta\ket^{\v\beta\v-1}}\Big )\Big (\frac 1{\bra x\ket^{\v
A\v-\v\beta\v}}+\frac 1{\bra \theta\ket^{\v A\v-\v\beta\v}}\Big )\\
&\leq  \frac
{C'} {\bra\theta\ket} \Big (\frac 1{\bra x\ket^{\v A\v+\sigma_0}}+\frac
1{\bra\theta\ket^{\v A\v-1}}\Big )\,.
\end{split}
\end{equation*}
Then using Lemma \ref{lV.3.6} and the fact that $\partial
^A_x\,X(\theta,\theta,x)=0$ since $\v A\v\geq 2$ we obtain (\ref{eqV.3.45})
for $\v A\v=q+1$. \cqfd


We need another lemma. \hfill\break
Let us recall that we have set $L_j(\theta,x)=\frac{\partial p}{\partial
\xi_j}\,(x,\Phi(\theta,x,\alpha))$.
\begin{lemma}\sl \label{lV.3.14}
Let $u_j(t)=\frac{\partial X_j}{\partial\theta}\,(t,\theta,x)+\som^n_{k=1}
\frac{\partial p}{\partial \xi_k}\,(x,\Phi(\theta,x,\alpha))\,\frac{\partial
X_j}{\partial x_k}\,(t,\theta,x)$. Then for every integer $N>0$ one can find a
constant $C_N>0$ such that for all $t\in[0,\theta]$ and all $(\theta,x)$ in
$\Omega_\delta$ we have
$$
\v u_j(t)\v\leq C_N \Big (\frac{\v
x-x(\theta,\alpha)\v}{\bra\theta\ket}\Big)^N\,.
$$
 \end{lemma}

{\bf Proof } First of all we claim that $u_j(\theta)=0$. Indeed since $\dot
X_j(\theta,\theta,x)+\frac{\partial X_j}{\partial \theta}\,(\theta,\theta,x)=0$
we have $\frac{\partial X_j}{\partial
\theta}\,(\theta,\theta,x)=-L_j(\theta,x)$. Then our claim follows from the
fact that $\frac{\partial X_j}{\partial x_k}\,(\theta,\theta,x)=\delta_{jk}$.
Now
$$
\dot u_j(t)=\frac{\partial \dot X_j}{\partial
\theta}\,(t,\theta,x)+\som^n_{k=1} L_k(\theta,x)\,\frac{\partial \dot
X_j}{\partial x_k}\,(t,\theta,x)\,.
$$
Using (\ref{eqV.3.4}) we obtain
\begin{equation*}
\begin{aligned}
\dot u_j(t)&=\sum^n_{\mu=1} \Big[\frac{\partial L_j}{\partial
z_\mu}\,(t,X(t,\theta,x))\,\frac{\partial X_\mu}
{\partial \theta}\,(t,\theta,x)+\frac{\partial L_j}{\partial \overline
z_\mu}\,(t,X(t,\theta,x))\,\overline{\frac{\partial X_\mu}{\partial
\theta}\,(t,\theta,x)}\Big]\\
&\quad  +\sum^n_{k=1} L_k(\theta,x) \sum^n_{\mu=1} \Big[\frac{\partial
L_j}{\partial z_\mu}\,(t,X(t,\theta,x))\,\frac{\partial X_\mu}{\partial
x_k}\,(t,\theta,x)+\frac{\partial L_j}{\partial \overline
z_\mu}\,(t,X(t,\theta,x))\,\overline {\frac{\partial X_\mu}{\partial
x_k}\,(t,\theta,x)}\Big ]\\
\dot u_j(t)&=\sum^n_{\mu=1} \frac{\partial L_j}{\partial
z_\mu}\,(t,X(t,\theta,x))\,u_\mu(t)+\sum^n_{\mu=1} \frac{\partial L_j}{\partial
\overline z_\mu}\,(t,X(t,\theta,x))\Big[\overline {\frac{\partial
X_\mu}{\partial \theta}\,(t,\theta,x)}\\
&\hbox to 9,8cm{}+\sum^n_{k=1} L_k(\theta,x)\,\overline
{\frac{\partial X_\mu}{\partial x_k}\,(t,\theta,x)}\Big]\,.
\end{aligned}
\end{equation*}
It follows then from (\ref{eqV.3.45}) that with $u(t)=(u_1(t),\ldots ,u_n(t))$,
$$
\Big\v\dot u_j(t)-\frac{2\,u_j(t)}{2t-i}\Big\v\leq \v h(t)\v\,\v u(t)\v+\v
g(t)\v\,,
$$
where
\begin{equation*}\left\{
\begin{array}{l}
h(t)=\som^n_{\mu,j=1}\Big\v\frac{\partial L_j}{\partial
z_\mu}\,(t,X(t,\theta,x))-\frac{2\,\delta_{jk}}{2t-i}\Big\v\,,\\
g(t)=C \som^n_{\mu,j=1} \Big\v\frac{\partial L_j}{\partial \overline
z_\mu}\,(t,X(t,\theta,x))\Big\v\,.
\end{array}\right.
\end{equation*}
Now using (\ref{eqV.3.46}) and (\ref{eqV.3.47}) we have $\int^\theta_0
h(t)\,dt\leq C$ and $\int^\theta_0 g(t)\,dt\leq C\,C_N\,\frac{\v
x-x(\theta,\alpha)\v^N}{\bra\theta\ket^N}$, so Lemma \ref{lV.3.14} follows from
Lemma \ref{lV.3.6} since $u_j(\theta)=0$. \cqfd

To solve the transport equations we need to introduce some notations. First of
all we shall set,
\begin{equation}\label{eqV.3.51}
D=\Big\{(\theta,z)\in\R\times\C^n : \v z-x(\theta,\alpha)\v\leq
\frac{\delta_1}{K}\,\bra\theta\ket\,,\enskip \v\Im z\v\leq
\frac{\delta_1}K\,,\enskip 
\Re z\cdot \alpha_\varepsilon\leq c_1\,\bra\Re z\ket\,\v\alpha_\varepsilon\v
\Big\}
\end{equation}
where $\delta_1,c_1,K$ have been introduce in the statement of Theorem
\ref{tV.3.2}.

Let $a\in C^\infty (D)$. We introduce some possible estimates.
\begin{equation}\label{eqV.3.52}\left\{
\begin{array}{l}
\forall \mu,\nu\in\N^n\,,\enskip \exists  
C_{\mu,\nu,\ell }\geq 0\enskip \textrm{ such that} \\
\textrm{ for all }(\theta,z)\in D\,,\enskip \v\partial ^\mu_z\,\partial
^\nu_{\overline z}\,a(\theta,z)\v\leq C_{\mu,\nu}
\end{array}\right.
\end{equation}
\begin{equation}\label{eqV.3.53}\left\{
\begin{array}{l}
\forall  N\in\N\,,\enskip \exists   C_N\geq 0 : \forall j=1,\ldots ,n\,,\enskip
\forall  (\theta,z)\in D\\
\Big\v\frac{\partial a}{\partial \overline z_j}\,(\theta,z)\Big\v\leq
C_N\,\v\Im z\v^N
\end{array}\right.
\end{equation}
\begin{equation}\label{eqV.3.54}\left\{
\begin{array}{l}
\exists  \sigma_0>0 : \forall \mu,\nu\in\N^n\,,\enskip \exists  
C_{\mu,\nu}\geq 0 : \forall  (\theta,z)\in D\\
\big\v\partial ^\mu_z\,\partial ^\nu_{\overline z}\,
a(\theta,z)\big\v\leq \frac{C_{\mu,\nu}}{\bra\Re
z\ket^{1+\sigma_0}}\,.
\end{array}\right.
\end{equation}
We first state the following result.
\begin{proposition}\sl \label{pV.3.15}
Let $u_0=u_0(z)$ be a $C^\infty $ function in a neighborhood of $D_0=\{z\in\C^n
: \v z-\alpha_x\v\leq \delta_1\}$ such that for any $N\in\N$ one can find
$C_N\geq 0$ such that for every $j=1,\ldots ,n$ and $z\in D_0$,
$$
\Big\v\frac{\partial u_0}{\partial \overline z_j}\,(z)\Big\v\leq C_N\,\v\Im
z\v^N\,.
$$
For $(\theta,x)\in\R\times\R^n$, $(\theta,x)\in D$ we set
$u(\theta,x)=u_0(X(0,\theta,x))$. Then for any $N\geq 0$ we can find $C'_N\geq
0$ such that
\begin{equation}\label{eqV.3.55}
\Big\v\frac{\partial u}{\partial \theta}\,(\theta,x)+\sum^n_{k=1}
\frac{\partial p}{\partial \xi_k}\,(x,\Phi(\theta,x,\alpha))\,\frac{\partial
u}{\partial x_k}\,(\theta,x)\Big\v\leq C'_N \Big (\frac{\v
x-x(\theta,\alpha)\v}{\bra\theta\ket}\Big)^N
\end{equation}
\begin{equation}\label{eqV.3.56}
u(0,x)=u_0(x)
\end{equation}
\begin{equation}\label{eqV.3.57}\left\{
\begin{array}{l}
\textrm{For any } \gamma\in\N^n\enskip \textrm{ one can find }\enskip 
C_{\gamma}\geq 0\enskip \textrm{ such that}\\
 \v\partial ^\gamma_x\,u(\theta,x)\v\leq C_{\gamma}\enskip \textrm{ for every
}\enskip  (\theta,x)\enskip \textrm{ in } D\cap \R\times\R^n\,.
\end{array}\right.
\end{equation}
\end{proposition}

{\bf Proof } First of all by (\ref{eqV.3.5}) (i) we have for $(\theta,x)\in
D\cap\R\times\R^n$,
$$
\v X(0,\theta,x)-\alpha_x\v\leq K\,\frac{\v
x-x(\theta,\alpha)\v}{\bra\theta\ket}\leq K\,\frac {\delta_1}K=\delta_1\,.
$$
Therefore $u(\theta,x)=u_0(X(0,\theta,x))$ is well defined and satisfies
(\ref{eqV.3.57}) by Proposition \ref{pV.3.12}, the fact that $u_0$ is $C^\infty
$ in a neighborhood of $D_0$ and the Faa di Bruno formula (Section
\ref{sVIII}). Now since $X(0,0,x)=x$, (\ref{eqV.3.56}) is obvious. Let us check
(\ref{eqV.3.55}). We set
$$
(1)=\frac{\partial u}{\partial \theta}\,(\theta,x)+\sum^n_{k=1} \frac{\partial
p}{\partial \xi_k}\,(x,\Phi (\theta,x,\alpha))\,\frac{\partial u}{\partial
x_k}\,(\theta,x)\,.
$$
Then
\begin{equation*}
\begin{split}
(1)&=\sum^n_{j=1} \frac{\partial u_0}{\partial z_j}\,(X(0,\theta,x))
\Big[\frac{\partial X_j}{\partial \theta}\,(0,\theta,x)+\sum^n_{k=1}
\frac{\partial p}{\partial \xi_k}\,(x,\Phi (\theta,x,\alpha))\,\frac{\partial
X_j}{\partial x_k}\,(0,\theta,x)\Big] \\
&\quad +\sum^n_{j=1} \frac{\partial u_0}{\partial \overline
z_j}\,(X(0,\theta,x))\Big[\overline {\frac{\partial X_j}{\partial
\theta}\,(0,\theta,x)}+\sum^n_{k=1} \frac{\partial p}{\partial \xi_k}\,(x,\Phi
(\theta,x,\alpha))\,\overline {\frac{\partial X_j}{\partial x_k}\,(0,\theta,x)}
\Big]
\end{split}
\end{equation*}
and we write $(1)=(A)+(B)$.

By Lemma \ref{lV.3.14} with $t=0$, the term $(A)$ satisfies (\ref{eqV.3.55}).
By the hypothesis made on $u_0$ and (\ref{eqV.3.5}) (ii) we have
$$
\Big\v\frac{\partial u_0}{\partial \overline z_j}\,(X(0,\theta,x))\Big\v\leq
C_N\,\v\Im X(0,\theta,x)\v^N\leq C'_N \Big (\frac{\v
x-x(\theta,\alpha)\v}{\bra\theta\ket}\Big)^N~.
$$
Using (\ref{eqV.3.45}) and the fact that $\Phi$ and the coefficients of $p$ are
bounded we deduce that the term $(B)$ satisfies also (\ref{eqV.3.55}).\cqfd
\begin{proposition}\sl\label{pV.3.16}
Let $a\in C^\infty $ on $D$ which satisfies (\ref{eqV.3.53}), (\ref{eqV.3.54}).
Let us set
$$
A(s,\theta,x)=\int^s_\theta a(\sigma,X(\sigma,\theta,x))\,d\sigma
$$
for $s\in[0,\theta]$ and $(\theta,x)\in D\cap\R\times\R^n$. Then
\begin{equation}\label{eqV.3.58}\left\{
\begin{array}{l}
\textrm{for any }N\in\N\textrm{ one can find } C_N\geq 0\textrm{ such that}\\
\Big\v\frac{\partial A}{\partial \theta}\,(s,\theta,x)+\som^n_{k=1}
\frac{\partial p}{\partial \xi_k}\,(x,\Phi (\theta,x,\alpha))\,\frac{\partial
A}{\partial x_k}\,(s,\theta,x)+a(\theta,x)\Big\v\leq C_N \Big (\frac{\v
x-x(\theta,\alpha)\v}{\bra\theta\ket}\Big )^N\\
\textrm{for all } s\in[0,\theta] \textrm{ and } (\theta,x)\in
D\cap\R\times\R^n\,.
\end{array}\right.
\end{equation}
\begin{equation}\label{eqV.3.59}
A(\theta,\theta,x)=0\,.
\end{equation}
\begin{equation}\label{eqV.3.60}\left\{
\begin{array}{l}
\textrm{For every }\enskip  \gamma\in\N^n\,,\enskip \textrm{there exists
}\enskip  C_{\gamma}\geq 0\enskip \textrm{ such that}\\
\v\partial ^\gamma_x\,A(s,\theta,x)\v\leq C_{\gamma}\enskip  \textrm{ on
}\enskip  [0,\theta]\times D\cap\R\times\R^n\,.
\end{array}\right.
\end{equation}
\end{proposition}

{\bf Proof } The claim (\ref{eqV.3.59}) is trivial, (\ref{eqV.3.60}) follows
from Proposition \ref{pV.3.12}, (\ref{eqV.3.54}) and (\ref{eqV.3.5}) (iii). Let
us show (\ref{eqV.3.58}). We set
$$
(1)=\frac{\partial A}{\partial \theta}\,(s,\theta,x)+\som^n_{k=1}
\frac{\partial p}{\partial \xi_k}\,(x,\Phi (\theta,x,\alpha))\,\frac{\partial
A}{\partial x_k}\,(s,\theta,x)+a(\theta,x)\,.
$$
Then
\begin{equation*}
\begin{split}
(1)&=\int^s_\theta \som^n_{j=1} \frac{\partial a}{\partial
z_j}\,(\sigma,X(\sigma,\theta,x))\left[\frac{\partial X_j}{\partial
\theta}\,(\sigma,\theta,x)+\som^n_{k=1}\frac{\partial p}{\partial
\xi_k}\,(x,\Phi (\theta,x,\alpha))\,\frac{\partial X_j}{\partial
x_k}\,(\sigma,\theta,x)\right]\,d\sigma\\
&\enskip +\int^s_\theta \som^n_{j=1} \frac{\partial a}{\partial
\overline z_j}\,(\sigma,X(\sigma,\theta,x))\left[\overline {\frac{\partial
X_j}{\partial
\theta}\,(\sigma,\theta,x)}+\som^n_{k=1}\frac{\partial p}{\partial
\xi_k}\,(x,\Phi (\theta,x,\alpha))\,\overline {\frac{\partial X_j}{\partial
x_k}\,(\sigma,\theta,x)}\right]\,d\sigma=(A)+(B)\,.
\end{split}
\end{equation*}
By (\ref{eqV.3.54}) and (\ref{eqV.3.5})(iii) we have $\big\v\frac{\partial
a}{\partial z_j}\,(\sigma,X(\sigma,\theta,x))\big\v\leq \frac
C{\bra\theta-\sigma\ket^{1+\sigma_0}}$. Therefore using Lemma \ref{lV.3.14} we
obtain (\ref{eqV.3.58}) for the term $(A)$.

Now it follows from (\ref{eqV.3.53}) and (\ref{eqV.3.54}) by interpolation that
$$
\Big\v\frac{\partial a}{\partial \overline z_j}\,(\theta,z)\Big\v\leq
C''_{N,\sigma_0}\,\frac{\v\Im z\v^N}{\bra\Re z\ket^{1+\frac{\sigma_0}2}}\,.
$$
Thus using (\ref{eqV.3.45}), (\ref{eqV.3.5}) (ii), (iii), we obtain
(\ref{eqV.3.58}) for the term $(B)$ since $\Phi $ and the coefficients of $p$
are uniformly bounded. \cqfd

This is the last result before the final one solving the transport equations.
\begin{proposition}\sl \label{pV.3.17}
Let $b$ be $C^\infty $ on $D$ satisfying (\ref{eqV.3.52}) and (\ref{eqV.3.53}).
Let us set $B(s,\theta,x)=b(s,X(s,\theta,x))$, $s\in[0,\theta]$, $(\theta,x)\in
D\cap \R\times\R^n$. Then,
\begin{equation}\label{eqV.3.61}\left\{
\begin{array}{l}
\textrm{for every } N\geq 0\textrm{ there exists } C_N\geq 0\textrm{ such
that}\\
\Big\v\frac{\partial B}{\partial \theta}\,(s,\theta,x)+\som^n_{k=1}
\frac{\partial p}{\partial \xi_k}\,(x,\Phi(\theta,x,\alpha))\,\frac{\partial
B}{\partial x_k}\,(s,\theta,x)\Big\v\leq C_N \Big (\frac{\v
x-x(\theta,\alpha)\v}{\bra\theta\ket}\Big )^N\,,
\end{array}\right.
\end{equation}
\begin{equation}\label{eqV.3.62}
B(\theta,\theta,x)=b(\theta,x)\,,
\end{equation}
\begin{equation}\label{eqV.3.63}\left\{
\begin{array}{l}
\textrm{for every } \gamma\in\N^n\,,\enskip \ell ,m\in\N\textrm{ there exists }
C_{\gamma}\geq 0\textrm{ such that}\\
\v\partial ^\gamma_x
\,B(s,\theta,x)\v\leq C_{\gamma}\,,\enskip \forall
(s,\theta,x)\in[0,\theta]\times D\cap \R\times\R^n\,.
\end{array}\right.
\end{equation}
\end{proposition}

{\bf Proof } The claim (\ref{eqV.3.62}) is obvious and (\ref{eqV.3.63}) follows
from Proposition \ref{pV.3.12} and (\ref{eqV.3.52}). Let us show
(\ref{eqV.3.61}). The left hand side of (\ref{eqV.3.61}) can be written,
\begin{equation*}
\begin{split}
(1)&= \som^n_{j=1} \frac{\partial b}{\partial
z_j}\,(s,X(s,\theta,x))\left[\frac{\partial X_j}{\partial
\theta}\,(s,\theta,x)+\som^n_{k=1}\frac{\partial p}{\partial
\xi_k}\,(x,\Phi (\theta,x,\alpha))\,\frac{\partial X_j}{\partial
x_k}\,(s,\theta,x)\right]\\
&\enskip +\som^n_{j=1} \frac{\partial b}{\partial
\overline z_j}\,(s,X(s,\theta,x))\left[\overline {\frac{\partial
X_j}{\partial
\theta}\,(s,\theta,x)}+\som^n_{k=1}\frac{\partial p}{\partial
\xi_k}\,(x,\Phi (\theta,x,\alpha))\,\overline {\frac{\partial X_j}{\partial
x_k}\,(s,\theta,x)}\right]=(A)+(B)\,.
\end{split}
\end{equation*}
The estimation of $(A)$ follows from Lemma \ref{lV.3.14} and (\ref{eqV.3.52}).
Now from (\ref{eqV.3.53}), (\ref{eqV.3.5}) (ii) and Proposition \ref{pV.3.12}
we deduce the estimation of $(B)$ since $\Phi $ and the coefficients of $p$ are
uniformly bounded. \cqfd
\begin{theorem}\sl \label{tV.3.18}
Let $a=a(\theta,z)$ be a $C^\infty $ function on $D$ satisfying
(\ref{eqV.3.53}), (\ref{eqV.3.54}). Let $b=b(\theta,z)$ be a $C^\infty $
function on $D$ satisfying (\ref{eqV.3.52}), (\ref{eqV.3.53}). Let $u_0=u_0(z)$
be a $C^\infty $ function on $D_0$ satisfying the hypothesis of Proposition
\ref{pV.3.15}. With the notations of Propositions \ref{pV.3.15},
\ref{eqV.3.16} and \ref{eqV.3.17} we set
$$
v(\theta,x)=\int^\theta_0
e^{A(s,\theta,x)}\,B(s,\theta,x)\,ds+e^{A(0,\theta,x)}\,u(\theta,x)\,.
$$
Then
\begin{equation}\label{eqV.3.64}\left\{
\begin{array}{l}
\textrm{for every }N\geq 1\textrm{ there exists } C_N\geq 0\textrm{ such that}\\
\Big\v\frac{\partial v}{\partial \theta}\,(\theta,x)+\som^n_{j=1}
\frac{\partial p}{\partial \xi_j}\,(x,\Phi(\theta,x,\alpha))\,\frac{\partial
v}{\partial x_k}(\theta,x)+a(\theta,x)\,v(\theta,x)-b(\theta,x)\Big\v\\
\hbox to 5cm{}\leq C_N\,\frac{\v x-x(\theta,\alpha)\v^N}{\bra\theta\ket^{N-1}}\\
\textrm{for all } (\theta,x)\in D\cap\R\times\R^n\,,
\end{array}\right.
\end{equation}
\begin{equation}\label{eqV.3.65}
v(0,x)=u_0(x)\,,
\end{equation}
\begin{equation}\label{eqV.3.66}\left\{
\begin{array}{l}
\textrm{for all }\gamma\in\N^n\,,\enskip \textrm{there exists }
C_{\gamma}\geq 0 \textrm{ such that}\\
\v\partial ^\nu_x\,v(\theta,x)\v\leq C_{\gamma}\bra\theta\ket\textrm{ for all
} (\theta,x)\in D\cap\R\times\R^n\,.
\end{array}\right.
\end{equation}
\end{theorem}

{\bf Proof } (\ref{eqV.3.65}) is obvious, (\ref{eqV.3.66}) follows from
(\ref{eqV.3.57}), (\ref{eqV.3.60}) and (\ref{eqV.3.63}). Let us show
(\ref{eqV.3.64}). We set
$$
\CL=\frac\partial {\partial \theta}+\sum^n_{k=1} \frac{\partial p}{\partial
\xi_k}\,(x,\Phi(\theta,x,\alpha))\,\frac\partial {\partial x_k}\,.
$$
Then
\begin{equation*}
\begin{split}
\CL v(\theta,x)+a(\theta,x)\,v(\theta,x)&-b(\theta,x)=b(\theta,x)+\int^\theta_0
e^{A(s,\theta,x)}\,[\CL B(s,\theta,x)+\CL A(s,\theta,x)\,B(s,\theta,x)]\,ds\\
&\quad
+e^{A(0,\theta,x)}\,(u(\theta,x)\,\CL
A(0,\theta,x) + \CL u(\theta,x))+a(\theta,x)\,
e^{A(0,\theta,x)}\,u(\theta,x)\\  &\quad +a(\theta,x) \int^\theta_0
e^{A(s,\theta,x)}\,B(s,\theta,x)\,ds -b(\theta,x)\,.
\end{split}
\end{equation*}
So
\begin{equation*}
\begin{split}
\CL v(\theta,x)+a(\theta,x)\,v(\theta,x)-b(\theta,x)=\int^\theta_0
e^{A(s,\theta,x)}\,\big[\CL
B(s,\theta,x)+(\CL  A(s,\theta,x)+a(\theta,x))\,B(s,\theta,x)
\big]\,ds\\
+e^{A(0,\theta,x)}\,\big[\CL
A(0,\theta,x)+a(\theta,x)\big]\,u(\theta,x)+e^{A(0,\theta,x)}\,\CL
u(\theta,x)\,.
\end{split}
\end{equation*}
The Propositions \ref{pV.3.15}, \ref{pV.3.16} and \ref{pV.3.17} show that
$$
\v \CL B(s,\theta,x)\v+\v \CL A(s,\theta,x)+a(\theta,x)\v+\v \CL
u(\theta,x)\v\leq C_N\,\frac{\v x-x(\theta,\alpha)\v^N}{\bra\theta\ket^N}
$$
and $\v A(s,\theta,x)\v \leq C$. Then (\ref{eqV.3.64}) follows. \cqfd

\noindent {\bf Proof of Theorem \ref{tV.1.1} (continued). Case of incoming
points } Let us set
\begin{equation}\label{eqV.3.67}\left\{
\begin{array}{l}
\CL =\frac\partial {\partial \theta}+\som^n_{j=1} \frac{\partial p}{\partial
\xi_j}\,(x,\Phi(\theta,x,\alpha))\,\frac\partial {\partial x_k}\\
a(\theta,x)=\som^n_{j,k=1} g^{jk}(x)\,\frac{\partial \Phi_j}{\partial
x_k}\,(\theta,x,\alpha)-\frac n 2\,\frac \theta{1+\theta^2}+i\,\som^n_{j=1}
g_j(x)\,\Phi_j(\theta,x,\alpha)\,.
\end{array}\right.
\end{equation}
By Proposition \ref{pIV.4.14} (ii) we have
$$
\Big\v\frac{\partial \varphi }{\partial
x}\,(\theta,x,\alpha)-\Phi(\theta,x,\alpha)\Big\v\leq C_N \left (\frac{\v
x-x(\theta,\alpha)\v}{\bra\theta\ket}\right )^N\,.
$$
Therefore using (\ref{eqV.1.3}), (\ref{eqV.1.4}) we see that to prove Theorem
\ref{tV.1.1} it will be sufficient to construct a smooth symbol
$f=f(\theta,x,\lambda)$ (with all derivatives in $x$ bounded) such that
$$
\Big\v\CL f+af+\frac 1 \lambda\,{}^tPf\Big\v\leq C_N\,\lambda^{-N}\,.
$$
We shall take $f$ on the form
$$
f(\theta,x,\lambda)=\sum^n_{k=0} \lambda^{-k}\,f_k(\theta,x,\alpha)\,,
$$
where the $f'_ks$ are the solutions of the problems,
\begin{equation}\label{eqV.3.68}\left\{
\begin{array}{l}
\CL f_0+a f_0=0\,,\quad f_0\v_{\theta=0}=1\,,\\
\CL f_k+a f_k=-{}^tP\,f_{k-1}\,, f_k\v_{\theta=0}=0\,,\quad k\geq 1\,.
\end{array}\right.
\end{equation}
By Theorem \ref{tV.3.18} we have,
\begin{equation}\label{eqV.3.69}\left\{
\begin{array}{l}
f_0(\theta,x)=e^{A(0,\theta,x)}\,,\\
f_k(\theta,x)=\int^\theta_0 e^{A(s,\theta,x)}\,B_k(s,\theta,x)\,ds\enskip
\textrm{where}\,,\\
A(s,\theta,x)=-\int^\theta_s a(\sigma,X(\sigma,\theta,x))\,d\sigma\,,\\
B_k(s,\theta,x)=-{}^tP\,f_{k-1} (s,X(s,\theta,x))\,,\quad k\geq 1\,.
\end{array}\right.
\end{equation}
Our aim is to prove, by induction on $k\geq 0$ that,
\begin{equation}\label{eqV.3.70}
\v\partial ^\gamma_x\,f_k(\theta,x)\v\leq C_{k,\gamma} \left (\frac 1{\bra
x\ket}+\frac 1{\bra\theta\ket}\right)^{\v\gamma\v}\,.
\end{equation}
We claim that for all $\ell \in\N^n$,
\begin{equation}\label{eqV.3.71}
\v\partial ^\ell _x\,A(s,\theta,x)\v\leq C_\ell \left (\frac 1{\bra
x\ket}+\frac 1{\bra\theta\ket}\right)^{\v\ell \v}\,,
\end{equation}
uniformly with respect to $s\in [0,\theta]$ and $(\theta,x)\in\Omega_\delta$.

Indeed using (\ref{eqV.3.18}) and the estimates on $\tilde a,\tilde b$ given
in Theorem \ref{tIV.4.2} we see easily that
\begin{equation}\label{eqV.3.72}
\v\partial ^\beta_{(z,\overline z)}\,a(\sigma,z)\v\leq C_\beta \left (\frac
1{\bra \Re z\ket^{\v\beta\v+2+\sigma_0}}+\frac
1{\bra\sigma\ket^{\v\beta\v+2}}\right)\,.
\end{equation}
Moreover by Theorem \ref{tV.3.2} (iii) we have
\begin{equation}\label{eqV.3.73}
\bra\Re X(\sigma,\theta,x)\ket\geq \frac 1 K\, (\bra
x\ket+\bra\theta-\sigma\ket)\,.
\end{equation}
Using the Faa di Bruno formula we see that $\partial ^\ell _x\,A(s,\theta,x)$
is bounded by a finite sum of terms of the following form.
$$
(1)=\int^\theta_0 \Big\v\partial ^\beta_{(z,\overline
z)}\,a(\sigma,X(\sigma,\theta,x)) \prod^s_{j=1} (\partial ^{\ell
_j}_x\,X(\sigma,\theta,x))^{k_j}\Big\v\,d\sigma
$$
where $1\leq \v\beta\v\leq \v\ell \v$, $1\leq s\leq \v\beta\v$, $\som^s_{j=1}
k_j=\beta$, $\som^s_{j=1} \v k_j\v\,\ell _j=\ell $.

Setting $I_1=\{j\in\{1,2,\ldots ,s\} : \v\ell _j\v=1\}$,
$I_2=\{j\in\{1,2,\ldots ,s\} : \v\ell _j\v\geq 2\}$ and using (\ref{eqV.3.72}),
(\ref{eqV.3.73})  and (\ref{eqV.3.45}) we can write
$$
(1)\leq C \int^\theta_0 \left[\frac 1{(\bra
x\ket+\bra\theta-\sigma\ket)^{\v\beta\v+2+\sigma_0}}+\frac
1{\bra\sigma\ket^{\v\beta\v+2}}\right]
\,\frac{\bra\sigma\ket^{\v\beta\v}}{\bra\theta\ket^{\v\beta\v}}\,d\sigma
\prod_{j\in I_2} \left(\frac 1{\bra x\ket^{\v\ell _j\v-1}}+\frac
1{\bra\theta\ket^{\v\ell _j\v-1}}\right)^{\v k_j\v}
$$
since $\som^s_{j=1} k_j=\beta$.

Using (\ref{eqV.3.50}) and the fact that $\som_{j\in I_2} \v k_j\v\,\v\ell
_j\v=\v\ell \v-\v\beta\v$ we obtain
$$
(1)\leq C\left(\frac 1{\bra x\ket}+\frac
1{\bra\theta\ket}\right)^{\v\beta\v}\,\left (\frac 1{\bra x\ket}+\frac
1{\bra\theta\ket}\right)^{\v\ell \v-\v\beta\v}
$$
which proves (\ref{eqV.3.71}).

As a consequence of (\ref{eqV.3.71}) we claim that
\begin{equation}\label{eqV.3.74}
\big\v\partial ^\gamma_x(e^{A(s,\theta,x)})\big\v\leq C_\gamma \left(\frac
1{\bra x\ket}+\frac 1{\bra\theta\ket}\right)^{\v\gamma\v}\,.
\end{equation}
Indeed by the Faa di Bruno formula the left hand side of (\ref{eqV.3.74}) can
be bounded by a finite sum of terms of the form
$$
(2)=\Big\v e^{A(s,\theta,x)} \prod^s_{j=1} \big (\partial ^{\ell
_j}_x\,A(s,\theta,x)\big)^{k_j}\Big\v
$$
where $1\leq s\leq \v\gamma\v$, $1\leq \som\,\v k_j\v\leq \ell $,
$\som^s_{j=1}\,\v k_j\v\,\ell _j=\v\gamma\v$.

Then (\ref{eqV.3.74}) follows easily from (\ref{eqV.3.71}).

Now (\ref{eqV.3.70}) for $k=0$ follows from (\ref{eqV.3.74}) (take $s=0$). On
the other hand by the Faa di Bruno formula $\partial
^\gamma_x\,f_k(\theta,x,\alpha)$ can be bounded by a finite sum of terms of
the form
$$
(3)=\int^\theta_0 \big\v \partial ^{\gamma_1}_x
(e^{A(s,\theta,x)})\big\v\,\partial ^{\gamma_2}_x
\big[^tP\,f_{k-1}(s,X(s,\theta,x))\big]\,ds\,.
$$
Setting for convenience $b_k={}^tP\,f_{k-1}$ it follows from the induction and
Lemma \ref{lV.3.1} that
\begin{equation}\label{eqV.3.75}
\big\v\partial ^\beta_{(z,\overline z)}\,b_k(s,z)\big\v\leq C_\beta
\left(\frac 1{\bra\Re z\ket}+\frac 1{\bra s\ket}\right)^{\v\beta\v+2}\,.
\end{equation}
Then the Faa di Bruno formula shows that the term $\partial
^{\gamma_2}_x[^tP\,f_{k-1}(s,X(s,\theta,x))]$ can be estimated by a finite sum
of terms of the form
$$
(4)=\Big\v\big (\partial ^\beta_{(z,\overline
z)}\,b_k\big)(s,X(s,\theta,x))\,\prod ^s_{\ell =1} (\partial ^{\ell
_j}_x\,X)^{k_j}\Big\v
$$
where $1\leq \v\beta\v\leq \v\gamma _2\v$, $1\leq s\leq \v\gamma_2\v$,
$\som^s_{j=1} k_j=\beta$, $\som^s_{j=1} \v k_j\v\,\ell _j=\gamma _2$. Then
using (\ref{eqV.3.75}), (\ref{eqV.3.73}), (\ref{eqV.3.45}) we see that
$$
(4)\leq C\left(\frac 1{(\bra x\ket+\bra\theta-s\ket)^{\v\beta\v+2}}+\frac
{1}{\bra s\ket^{\v\beta\v+2}}\right)\,\frac{\bra
s\ket^{\v\beta\v}}{\bra\theta\ket^{\v\beta\v}} \prod_{j\in I_2} \left(\frac
1{\bra x\ket}+\frac 1{\bra\theta\ket}\right)^{\v k_j\v(\v\ell _j\v-1)}\,.
$$
Now by (\ref{eqV.3.71}) we can write
$$
(3)\leq C \int^\theta_0 \left[\frac 1{(\bra
x\ket+\bra\theta-s\ket)^{\v\beta\v+2}}+\frac 1{\bra
s\ket^{\v\beta\v+2}}\right]\,\frac{\bra
s\ket^{\v\beta\v}}{\bra\theta\ket^{\v\beta\v}}\,ds \left(\frac 1{\bra
x\ket}+\frac 1{\bra\theta\ket}\right)^{\v\gamma_1\v+\v\gamma_2\v-\v\beta\v}
$$
since $\som_{j\in I_2} \v k_j\v(\v\ell _j\v-1)=\v\gamma_2\v-\v\beta\v$.

It follows from (\ref{eqV.3.50}) that
$$
(3)\leq C\left(\frac 1{\bra x\ket^{\v\beta\v}}+\frac
1{\bra\theta\ket^{\v\beta\v}}\right) \left(\frac 1{\bra x\ket}+\frac
1{\bra\theta\ket}\right)^{\v\gamma_1\v+\v\gamma_2\v-\v\beta\v}\leq C'
\left(\frac 1{\bra x\ket}+\frac
1{\bra\theta\ket}\right)^{\v\gamma\v}
$$
since $\v\gamma_1\v+\v\gamma_2\v=\v\gamma\v$. This proves (\ref{eqV.3.70}) for
all $k$.


\subsection{The amplitude for short time}\label{ssV.4}

We shall need the following precision on the amplitude when $\v\theta\v\leq 1$.
\begin{proposition}\sl \label{pV.4.1}
Let $a_N$ be the amplitude defined in Corollary \ref{cV.1.2}. Then for every
$\gamma\in \N^{2n}$ one can find a constant $C_\gamma\geq 0$ such that
$$
\v\partial ^\gamma_\alpha [a_N(\theta,x,\alpha,\lambda)]\v\leq C_\gamma
$$
for all $\v\theta\v\leq 1$, $\v x-x(\theta,\alpha)\v\leq \delta
\bra\theta\ket$, $\lambda\geq 1$ and $\alpha\in T^*\R^n$ such that $\frac
12\leq \v\alpha_\varepsilon\v\leq 2$.
\end{proposition}

{\bf Proof } When $\theta\leq 1$ we can use the method of  Section
\ref{ssV.2} no matter is $\alpha$, provided that $\frac12\leq \v
\alpha_\xi\v\leq 2$. Let us recall how the amplitude $a_N$ is produced.
We have
\begin{equation}\label{eqV.4.1}\left\{
\begin{array}{l}
a_N(\theta,x,\alpha,\lambda)=\bra\theta\ket^{-\frac
n2}\,e_N(\theta,x-x(\theta,\alpha),\alpha,\lambda)\,,\\
e_N(\theta,x-x(\theta,\alpha),\alpha,\lambda)=f_N\Big
(\theta,\,\frac{x-x(\theta,\alpha)}{\bra\theta\ket}\,,\alpha,\lambda\Big)\,,\\
f_N(\theta,z,\alpha)=\som^{N+1}_{\ell =0} \lambda^{-\ell }\,A_\ell
(\theta,z,\alpha)\,,\\
L\,A_0=0\,,\enskip L\,A_\ell =i\,Q\,A_{\ell -1}\,,\enskip \ell =1,\ldots
,N+1\,,\\
A_0(0,z,\alpha)=1\,,\enskip A_\ell (0,z,\alpha)=0\,,\\
\v\partial ^\beta_z\,A_\ell (\theta,z,\alpha)\v\leq C_\beta\,,\enskip \forall
\beta\in\N^n\,,\enskip \v\theta\v\leq 1\,,\enskip \v z\v\leq \delta\,,\enskip
\frac12\leq \v\alpha_\xi\v\leq 2\,.
\end{array}\right.
\end{equation}
Now, according to Proposition \ref{pIII.2.1} for every $\gamma\in\N^{2n}$ such
that $\v\gamma\v\geq 1$ one can find $C'_\gamma\geq 0$ such that
\begin{equation}\label{eqV.4.2}
\v\partial ^\gamma_\alpha\,x(\theta,\alpha)\v+\v\partial
^\gamma_\alpha\,\xi(\theta,\alpha)\v\leq C'_\gamma
\end{equation}
if $\v\theta\v\leq 1$, $\alpha\in T^*\R^n$, $\frac12\leq
\v\alpha_\xi\v\leq 2$.

Assume that we show that for all $\beta\in\N^n$, $\gamma\in\N^{2n}$ one can
find $C_{\beta\gamma}\geq 0$ such that
\begin{equation}\label{eqV.4.3}
\v\partial ^\beta_z\, \partial ^\gamma_\alpha\,A_\ell
(\theta,z,\alpha)\v\leq C_{\beta\gamma}
\end{equation}
if $\v\theta\v\leq 1$ and $\frac12\leq \v\alpha_\xi\v\leq 2$.
It will follow from (\ref{eqV.4.1}) to (\ref{eqV.4.3}) and the Faa di Bruno
formula that
\begin{equation}\label{eqV.4.4}
\v\partial ^\gamma_\alpha\,a_N(\theta,x,\alpha,\lambda)\v\leq C_\gamma
\end{equation}
if $\v\theta\v\leq 1$, $\v x-x(\theta,\alpha)\v\leq \delta\,\bra\theta\ket$,
$\frac12\leq \v\alpha_\xi\v\leq 2$, which is the claim of Proposition
\ref{pV.4.1}.

So we are left with the proof of (\ref{eqV.4.3}). By (\ref{eqV.2.8}),
(\ref{eqV.2.9}) and (\ref{eqV.4.2}) for all $\mu\in\N^n$, $\gamma\in\N^{2n}$
there exists $C_{\mu\gamma}\geq 0$ such that
$$
\v\partial ^{\mu}_y\,\partial ^\gamma_\alpha\,E_j(s,y,\alpha)\v\leq
C_{\mu\gamma}\,,
$$
for all $\v s\v\leq 1$, $\v y\v\leq \delta\,\bra s\ket$ and $\frac12\leq
\v\alpha_\xi\v\leq 2$.\hfill\break
It follows from (\ref{eqV.2.10}), (\ref{eqV.2.11}) that for all $\beta\in\N^n$,
$\gamma\in\N^{2n}$ there exists $C_{\beta\gamma}\geq 0$ such that
\begin{equation}\label{eqV.4.5}
\v\partial ^\beta_z\,\partial ^\gamma_\alpha\,h_j(\theta,z,\alpha)\v\leq
C_{\beta\gamma}\,,
\end{equation}
for all $\v\theta\v\leq 1$, $\v z\v\leq \delta$, $\frac12\leq
\v\alpha_\xi\v\leq 2$.\hfill\break
And we see easily from (\ref{eqV.2.5}), (\ref{eqV.4.2}), (\ref{eqV.2.15}),
(\ref{eqV.2.16}) that $h^{N_0}_j$, $d^{N_0}$, $k^{N_0}_\nu$ satisfy also the
bound (\ref{eqV.4.5}).
By induction on the size of derivation, using the Faa di Bruno formula and the
Gronwall Lemma we see easily that the solution $z=z(\theta,y,\alpha)$ of
(\ref{eqV.2.26}) satisfies the bound
\begin{equation}\label{eqV.4.6}
\v\partial ^\mu_y\,\partial ^\gamma_\alpha \,z(\theta,y,\alpha)\v\leq
C_{\mu\gamma}
\end{equation}
uniformly for $\v\theta\v\leq 1$, $\v y\v\leq \eta$, $\frac12\leq
\v\alpha_\xi\v\leq 2$.\hfill\break
Moreover by Lemma \ref{lV.2.3}, if we denote by $\kappa(\theta,z)$ the inverse
map of $z(\theta,y)$ we have also by (\ref{eqV.4.6}),
\begin{equation}\label{eqV.4.7}
\v\partial ^\beta_z\,\partial ^\gamma_\alpha\,\kappa(\theta,z,\alpha)\v\leq
C_{\beta\gamma}
\end{equation}
uniformly for $\v\theta\v\leq 1$, $\v z\v\leq \delta$, $\frac12\leq
\v\alpha_\xi\v\leq 2$.
Finally we have set for $\ell =0,1,\ldots ,N_0+1$,
\begin{equation}\label{eqV.4.8}
A_\ell (\theta,z,\alpha)=\tilde A_\ell (\theta,\kappa(\theta,z,\alpha),\alpha)
\end{equation}
where
\begin{equation*}
\begin{aligned}
\tilde A_0(\theta,y,\alpha)&=\exp \Big[-\int^\theta_0 d^{N_0}
(t,z(t,y,\alpha),\alpha)\,dt\Big]\\
\tilde A_\ell (\theta,y,\alpha)&=\exp \Big[-\int^\theta_0 d^{N_0}
(t,z(t,y,\alpha),\alpha)\,dt\Big] \int^\theta_0 i(\tilde Q\,A_{\ell
-1})(t,y,\alpha)\,dt
\end{aligned}
\end{equation*}
so using (\ref{eqV.4.6}), (\ref{eqV.4.7}), (\ref{eqV.4.8}) and the estimates
(\ref{eqV.4.5}) for $d^{N_0}$, $k^{N_0}_\nu$ we see easily that (\ref{eqV.4.3})
holds, which completes the proof of Proposition \ref{pV.4.1}. \cqfd

\section{Microlocal localizations and the use of the FBI transform}\label{sVI}

In this Section using the phase and the amplitude constructed in Sections
\ref{sIV} and \ref{sV} we shall define general FBI transforms which will lead
to a parametrix for the Schr\"odinger equation. These constructions will be
microlocal so we will need several microlocal localizations.

\subsection{Preliminaries}\label{ssVI.1}

\subsubsection{The semi-classical calculus}\label{sssVI.1.1}

We shall work with semi-classical pseudo-differential operators (p.d.o) and we
shall use the Weyl calculus described by H\"ormander. We refer to \cite{H}
for notations and details.

Let $p\in S^m_{1,0}(\R^n)$ (the usual class of symbol of order $m$) and let us
set $a(x,\xi)=p\big (x,\,\frac \xi\lambda\big )$, $\lambda\geq 1$. It is easy
to see that $a\in S(M,g)$ where,
$$
g=dx^2+\frac{d\xi^2}{\lambda^2+\v\xi\v^2}\,,\quad
M=\lambda^{-m}(\lambda^2+\v\xi\v^2)^{\frac m2}\,.
$$
The p.d.o associated to the symbol $a$ is denoted by $p\big (x,\frac
D\lambda\big )$. Then we have the following symbolic calculus.

i) Let $p\in S^m_{1,0}$, $q\in S^{m'}_{1,0}$. Then one can find $\ell
_\lambda\in S^{m+m'}_{1,0}$ such that
$$
p\Big (x,\,\frac D\lambda\Big )\circ q\Big (x,\,\frac D\lambda\Big )=\ell
_\lambda\Big (x,\,\frac D\lambda\Big)\,.
$$
The semi norms of $\ell _\lambda$ are uniformly bounded when $\lambda\geq 1$
and for any $N\in\N^*$ we have
$$
\ell _\lambda(x,\xi)=\sum_{\v\alpha\v\leq N-1} \frac 1{\alpha!}\,\frac
1{i^{\v\alpha\v}}\,\frac 1{\lambda^{\v\alpha\v}}\,\partial
^\alpha_\xi\,p(x,\xi)\,\partial ^\alpha_x\,q(x,\xi)+\frac
1{\lambda^N}\,r_N(\lambda,x,\xi)
$$
where $r_N\in S^{m+m'-N}_{1,0}$ uniformly for $\lambda\geq 1$.

ii) Let $p\in S^0_{1,0}$. Then there exists $C>0$ such that
$$
\Big\Vert p\Big(x,\,\frac D\lambda\Big)\,u\Big\Vert _{L^2}\leq C\,\Vert
u\Vert _{L^2}
$$
for every $u\in L^2(\R^n)$ and $\lambda\geq 1$. As a consequence, for all
$s\in\R$ and all $p\in S^m_{1,0}$ one can find a constant $C>0$ such that for
every $u\in\CS(\R^n)$,
$$
\Big\Vert \Big(I-\frac\Delta{\lambda^2}\Big)^s\,p\Big
(x,\,\frac{D}{\lambda}\Big)\,u\Big\Vert _{L^2}\leq C\,
\Big\Vert \Big(I-\frac\Delta{\lambda^2}\Big)^{s+m}\,u\Big\Vert _{L^2}
$$
for all $\lambda\geq 1$, where $\Delta=\som^n_{j=1} \partial ^2_j$.

\subsubsection{The FBI transform}\label{sssVI.1.2}

We recall here the definition of the classical FBI transform as described in 
Sj\"ostrand  \cite{Sj}. We set for $\alpha=(\alpha_x,\alpha_\xi)\in
T^*\R^n$, $\lambda\geq 1$, $u\in L^2(\R^n)$ and $c_n=2^{-\frac n2}\pi
^{-\frac{3n}4}$,
\begin{equation}\label{eqVI.1.1}
Tu(\alpha,\lambda)=c_n\,\lambda^{\frac{3n}4} \int_{\R^n}
e^{i\lambda(y-\alpha_x)\cdot \alpha_\varepsilon-\frac \lambda2 \v
y-\alpha_x\v^2+\frac \lambda2 \v\alpha_\xi\v^2}\,u(y)\,dy\,.
\end{equation}
Then $T$ maps continuously the space $L^2(\R^n)$ to the space of functions
$v=v(\alpha) $ such that $e^{-\frac \lambda2 \v\alpha_\xi\v^2}\,v\in
L^2(\R^{2n})$.

The adjoint of $T$ is then given by the formula
\begin{equation}\label{eqVI.1.2}
T^*v(x,\lambda)=c_n\,\lambda^{\frac {3n}4} \int e^{-i\lambda(x-\alpha_x)\cdot
\alpha_\xi-\frac\lambda2 \v x-\alpha_x\v^2-\frac\lambda2
\v\alpha_\xi\v^2}\,v(x)\,d\alpha\,.
\end{equation}
Then we have
\begin{equation}\label{eqVI.1.3}
T^*T\textrm{ is the identity operator of } L^2(\R^n)\,.
\end{equation}
We shall need also the expressions of $T$ and $T^*$ by means of the Fourier
transform. We have,
\begin{equation}\label{eqVI.1.4}\left\{
\begin{array}{l}
Tu(\alpha,\lambda)=\Big(\frac\lambda\pi \Big )^{\frac n4} \int e^{i\sigma\cdot
\alpha_x-\frac\lambda2
\v\alpha_\xi+\frac\sigma\lambda\v^2+\frac\lambda2\v
\alpha_\xi\v^2}\,\hat u(\sigma)\,d\sigma\\
\widehat{T^*v}(\xi,\lambda)=c'_n\,\lambda^{\frac n4} \int e^{-i\xi\cdot
\alpha_x-\frac\lambda2
\v\alpha_\xi+\frac\xi\lambda\v^2-\frac\lambda
2\v\alpha_\xi\v^2}\,v(\alpha)\,d\alpha\,.
\end{array}\right.
\end{equation}
Let us consider now a self adjoint operator,
\begin{equation}\label{eqVI.1.5}
P=\sum^n_{j,k=1} D_j (g^{jk}\,D_k)+\sum^n_{j=1} (D_j\,b_j+b_j\,D_j)+b_0
\end{equation}
with $g^{jk}=\delta_{jk}+\varepsilon\,b_{jk}$, where $\varepsilon$ is a small
positive constant, $\delta_{jk}$ is the Kronecker symbol and,
$$
\sum_{\v\alpha\v=k} \v\partial ^\alpha_x\,b_{jk}(x)\v\leq \frac{A_k}{\bra
x\ket^{k+1+\sigma_0}}\,,\quad k=0,1,\ldots ,x\in\R^n\,,\quad \sigma_0>0\,.
$$
Then by interpolation and duality we can prove the following estimates. For all
$s\in\R$ there exists $C\geq 1$ such that for all $u\in C^\infty _0(\R^n)$
\begin{equation}\label{eqVI.1.6}
\frac 1C \Vert u\Vert _{H^s}\leq \Vert (I+P)^{\frac s 2}\,u\Vert _{L^2}\leq
C\,\Vert u\Vert _{H^s}\,.
\end{equation}
For all $s\geq 0$ there exists $C\geq 1$ such that for all $\lambda\geq 1$ and
$u\in C^\infty _0(\R^n)$,
\begin{align}
\frac{1}{C}\,\lambda^{-s}\,\Vert u\Vert _{H^s}&\leq \Big\Vert \Big (I+P\Big
(x,\,\frac D\lambda\Big)\Big)^{\frac s 2}\,u\Big\Vert _{L^2}\leq C\,\Vert
u\Vert _{H^s}\\
\frac{1}{C}\,\Vert u\Vert _{H^{-s}}&\leq \Big\Vert \Big (I+P\Big
(x,\,\frac D\lambda\Big)\Big)^{-\frac {s}{2}}\,u\Big\Vert _{L^2}\leq
C\,\lambda^s\Vert u\Vert _{H^{-s}}\,.
\end{align}

\subsection{The microlocalization procedure}\label{ssVI.2}

We begin by introducing several cut-off functions. Generally speaking we shall
denote by $\chi$ (resp. $\psi$) cut-off functions in the space (resp.
frequency) variables.

Let $\xi_0\in\R^n$, $\v\xi_0\v=1$ be fixed.

Let $\chi_0\in C^\infty (\R)$ be such that,
\begin{equation}\label{eqVI.2.1}
\chi_0(s)=1 \textrm { if } s\leq \frac34\,,\quad \chi_0(s)=0\textrm{ if } s\geq
1\,,\quad 0\leq \chi_0\leq 1\,.
\end{equation}
With $\delta_1>0$ to be chosen later on we set
\begin{equation}\label{eqVI.2.2}\left\{
\begin{array}{ll}
\chi^+_1(x)=\chi_0\Big (-\frac{x\cdot \xi_0}{\delta_1}\Big)\,,\enskip
\chi^+_2(x)=\chi_0\Big (-\frac{x\cdot \xi_0}{2\delta_1}\Big )\,,\enskip
&\chi^+_3(x)=\chi_0\Big (-\frac{x\cdot \xi_0}{3\delta_1}\Big)\\
\chi^-_1(x)=\chi_0\Big (\frac{x\cdot \xi_0}{\delta_1}\Big)\,,\enskip
\chi^+_2(x)=\chi_0\Big (\frac{x\cdot \xi_0}{2\delta_1}\Big )\,,\enskip
&\chi^+_3(x)=\chi_0\Big (\frac{x\cdot \xi_0}{3\delta_1}\Big)\,.
\end{array}\right.
\end{equation}
These cut-off functions will correspond to outgoing and incoming points. Now
for convenience we shall set,
\begin{equation}\label{eqVI.2.3}
a=\frac 6{10}\,,\quad b=\frac{19}{10}\,,
\end{equation}
and with $\delta_2\leq \frac 1{100}$ (chosen later on) we introduce the
following cut-off functions.

Let $\psi_0\in C^\infty (\R^n)$ be such that $0\leq \psi_0\leq 1$ and
\begin{equation}\label{eqVI.2.4}\left\{
\begin{array}{l}
\psi_0(\xi)=1\textrm { if }\enskip 
\Big\v\frac\xi{\v\xi\v}-\xi_0\Big\v\leq
\delta_2\textrm{ and } \v\xi\v\geq 2\delta_2\\
\supp \psi_0\subset \Big\{\xi : \Big\v\frac\xi{\v\xi\v}-\xi_0\Big\v\leq
2\delta_2\textrm{ and } \v\xi\v\geq \delta_2\Big\}
\end{array}\right.
\end{equation}
Let $\psi_1\in C^\infty _0(\R^n)$ be such that, $0\leq \psi_1\leq 1$ and
\begin{equation}\label{eqVI.2.5}\left\{
\begin{array}{l}
\psi_1(\xi)=1\textrm{ if } a-\delta_2\leq \v\xi\v\leq b+\delta_2\,,\\
\supp \psi_1\subset \big\{\xi : a-2\delta_2\leq \v\xi\v\leq
b+2\delta_2\big\}\,.
\end{array}\right.
\end{equation}
We shall set
\begin{equation}\label{eqVI.2.6}
\psi_2(\xi)=\psi_0(\xi)\,\psi_1(\xi)\,.
\end{equation}
Now we introduce for $t\in\R$ the operators,
\begin{equation}\label{eqVI.2.7}\left\{
\begin{array}{l}
U_+(t)=\chi^+_1(x)\,\psi_2\Big(\frac D\lambda\Big)\,e^{-itP}\,,\\
U_-(t)=\chi^-_1(x)\,\psi_2\Big(\frac D\lambda\Big)\,e^{-itP}\,.
\end{array}\right.
\end{equation}
It follows that if we denote by $U^*$ the adjoint of $U$ then we have
\begin{equation}\label{eqVI.2.8}\left\{
\begin{array}{l}
U_{\pm}(t_1)\,U_{\pm}(t_2)^*=K_\pm(t_1-t_2)\textrm{ where,}\\
K_+(t)=\chi^+_1(x)\,\psi_2\Big(\frac D\lambda\Big)\,e^{-itP}\,\psi_2\Big(\frac
D\lambda\Big)\,\chi^+_1(x)\,,\\
K_-(t)=\chi^-_1(x)\,\psi_2\Big(\frac D\lambda\Big)\,e^{-itP}\,\psi_2\Big(\frac
D\lambda\Big)\,\chi^-_1(x)\,.
\end{array}\right.
\end{equation}
Let now $\psi_3\in C^\infty _0(\R^n)$ be such that $0\leq \psi_3\leq 1$ and
\begin{equation}\label{eqVI.2.9}\left\{
\begin{array}{ll}
\psi_3(\xi)=1\textrm{ if }\enskip 
\Big\v\frac\xi{\v\xi\v}+\xi_0\Big\v\leq
3\,\delta_2&\textrm{ and } a-3\delta_2\leq \v\xi\v\leq b+3\delta_2\\
\supp \psi_3 \subset \Big\{\xi : \Big\v\frac\xi
{\v\xi\v}+\xi_0\Big\v\leq 4\delta_2&\textrm{ and }
a-4\delta_2\leq \v\xi\v\leq b+4\delta_2\Big\}\,.
\end{array}\right.
\end{equation}
The first localization result requires to introduce the following Definition.
\begin{definition}\sl \label{dVI.2.1}
We shall call $\CR$ the set of families of operators $R=(R_\lambda(t))$
depending on $\lambda\in[1,+\infty [$ and $t\in[-T,T]$ such that for every $N$
in $\N$ one can find a constant $C_N\geq 0$ such that for every $u\in\CS(\R^n)$,
\begin{equation}\label{eqVI.2.10}
\Vert R_\lambda(t)\,u\Vert _{H^{2N}(\R^n)}\leq C_N\,\Vert u\Vert
_{H^{-2N}(\R^n)}
\end{equation}
uniformly with respect to $(\lambda,t)\in[1,+\infty [\times [-T,T]$.
\end{definition}
Then we can state the following result.
\begin{theorem}\sl \label{tVI.2.2}
Let $T>0$. For every $t\in[-T,T]$ and $\lambda\geq 1$ one can write
\begin{equation}\label{eqVI.2.11}
K_+(t)=\chi^+_1\,\psi_2\Big(\frac
D\lambda\Big)\,\chi^+_2\,T^*_{\alpha\rightarrow
x}\,\chi^+_3(\alpha_x)\,\psi_3(\alpha_\xi)\,T_{y\rightarrow
\alpha}\Big[e^{-itP}\,\chi^+_2\,\psi_2\Big(\frac D\lambda\Big)\,\chi^+_1
\Big]+R^+_\lambda(t)
\end{equation}
where $(R^+_\lambda(t))\in\CR$.\hfill\break
The same formula is true with the sign $-$ instead of $+$.
\end{theorem}

The proof of this Theorem requires several steps.
\begin{lemma}\sl \label{lVI.2.3}
There exist a constant $C>0$ such that
$$
\Big\Vert \psi_2\Big(\frac
D\lambda\Big)\,T^*\big[(1-\psi_3(\alpha_\xi))\,v\big]\Big\Vert
_{L^2(\R^n)}\leq C\,e^{-\frac \lambda 8\,\delta^2_2} \big\Vert e^{-\frac\lambda
2 \v\alpha_\xi\v^2}\,v\big\Vert _{L^2(\R^{2n})}
$$
for every $v$ such that $e^{-\frac{\lambda }{2}\v \alpha_\xi\v^2}v\in
L^2(\R^{2n})$.
\end{lemma}

{\bf Proof } We claim that on the support of
$\psi_2(\xi)(1-\psi_3(\alpha_\xi))$ we have
\begin{equation}\label{eqVI.2.12}
\v\xi+\alpha_\xi\v \geq \frac12\,\delta_2\,\v\alpha_\xi\v.
\end{equation}
Indeed, according to (\ref{eqVI.2.4}) to (\ref{eqVI.2.6}) and (\ref{eqVI.2.9})
we have on this support $\big\v\frac\xi{\v\xi\v}-\xi_0\v\leq
2\delta_2$, $a-2\delta_2\leq \v\xi\v\leq b+2\delta_2$ and either
$\v\alpha_\xi\v\leq a-3\delta_2$ or $\v\alpha_\xi\v\geq
b+3\delta_2$ or
$\big\v\frac{\alpha_\xi}{\v\alpha_\xi\v}+\xi_0\big\v\geq
3\delta_2$. In the first case we have $\v\xi+\alpha_\xi\v\geq 
\v\xi\v-\v\alpha_\xi\v\geq \delta_2$. In the second case we have
$\v\xi+\alpha_\xi\v\geq \v\alpha_\xi\v-\v\xi\v\geq \delta_2$ and in
the third one we have
$$
\Big\v\frac{\alpha_\xi}{\v\alpha_\xi\v}+\frac\xi{\v\xi\v}\Big\v
\geq
\Big\v\frac{\alpha_\xi}{\v\alpha_\xi\v}+\xi_0\Big\v-
\Big\v\frac\xi{\v\xi\v}-\xi_0\Big\v\geq \delta_2\,.
$$
Therefore
$\big\v\alpha_\xi+\frac{\v\alpha_\xi\v}{\v\xi\v}\,\xi\big\v\geq
\delta_2\,\v\alpha_\xi\v$. It follows that
$$
\v\alpha_\xi+\xi\v\geq
\Big\v\alpha_\xi+\frac{\v\alpha_\xi\v}{\v\xi\v}\,\xi\Big\v
-\Big\v\xi-\frac{\v\alpha_\xi\v}{\v\xi\v}\,\xi\Big\v\geq
\delta_2\,\v\alpha_\xi\v-\big\v\v\xi\v-\v\alpha_\xi\v\big\v\geq
\delta_2\,\v\alpha_\xi\v-\v\alpha_\xi+\xi\v
$$
so $\v\alpha_\xi+\xi\v\geq \frac12\,\delta_2\,\v\alpha_\xi\v$ and our claim is
proved.

Now using (\ref{eqVI.1.4}) we can write
\begin{equation*}
\begin{split}
(1)&=:\CF\Big(\psi_2\Big(\frac
D\lambda\Big)\,T^*\big[(1-\psi_3(\alpha_\xi))\big]\,v\Big)(\xi)\\
&=c_n\,\lambda^{
\frac n4}\,\psi_2 \Big (\frac\xi\lambda\Big) \int_{\R^{2n}}
e^{-i\xi\cdot \alpha_x-\frac\lambda 2
\v\alpha_\xi+\frac\xi\lambda\v^2-\frac \lambda 2
\v\alpha_\xi\v^2} (1-\psi_3(\alpha_\xi))\,v(\alpha)\,d\alpha\,.
\end{split}
\end{equation*}
Let us set
$$
w(\xi,\alpha_\xi)=\int_{\R^{n}} e^{-i\xi\cdot
\alpha_x-\frac\lambda 2
\v\alpha_\xi\v^2}\,v(\alpha)\,d\alpha_x=\CF_{\alpha_x}\big (e^{-\frac
\lambda 2 \v\alpha_\xi\v^2}\,v\big)(\xi,\alpha_\xi)\,.
$$
Using Cauchy-Schwartz inequality we obtain
$$
\v(1)\v^2\leq C\,\lambda^{\frac n 2} \left (\int_{\R^n}
e^{-\lambda\v\alpha_\xi+\frac\xi\lambda\v^2}\,\psi^2_2
\Big(\frac\xi\lambda\Big)\big
(1-\psi_3(\alpha_\xi)\big)^2\,d\,\alpha_\xi\right )\left (\int \v
w(\xi,\alpha_\xi)\v^2\,d\alpha_\xi\right)
$$
so by (\ref{eqVI.2.12})
$$
\v(1)\v^2\leq C\,\lambda^{\frac\lambda n} \,e^{-\frac\lambda 8
\delta^2_2}\left (
\int e^{-\frac\lambda 2
\v\alpha_\xi+\frac\xi\lambda\v^2}d\,\alpha_\xi\right)\left
(\int \v w(\xi,\alpha_\xi)\v^2\,d\alpha_\xi\right)\,.
$$
It follows that
$$
\v(1)\v^2\leq C'\,e^{-\frac\lambda 8 \delta^2_2} \int \v
w(\xi,\alpha_\xi)\v^2\,d\alpha_\xi\,.
$$
Integrating with respect to $\xi$ and using Parseval identity we obtain
$$
\Big\Vert \psi_2\Big(\frac
D\lambda\Big)\,T^*\big[(1-\psi_3(\alpha_\xi))\,v\big]\Big\Vert
_{L^2}\leq C\,e^{-\frac\lambda 8 \delta^2_2} \int_{\R^{2n}} e^{-\frac\lambda 2
\v\alpha_\xi\v^2}\,\v v(\alpha)\v^2\,d\alpha\,.
$$
\cqfd

\begin{corollary}\sl \label{cVI.2.4}
We have for $t\in [-T,T]$ and $\lambda\geq 1$,
$$
K_+(t)=\chi^+_1\,\psi_2\Big(\frac
D\lambda\Big)\,T^*\,\psi_3(\alpha_\xi)\,T\,e^{-itP}\,\psi_2\Big(\frac
D\lambda\Big)\,\chi^+_1+R^+_\lambda(t)
$$
where $(R^+_\lambda(t))\in\CR$. The same is true for the sign $-$.
\end{corollary}

{\bf Proof } Using (\ref{eqVI.1.3}) we write
$\Id=T^*T=T^*\,\psi_3(\alpha_\xi)\,T+T^*(1-\psi_3(\alpha_\xi))\,T$.
So we have to prove that the family of operators
$$
R^+_\lambda(t)=\chi^+_1\,\psi_2\Big(\frac
D\lambda\Big)\,T^*(1-\psi_3(\alpha_\xi))\,T\,e^{-itP}\,\psi_2\Big(\frac
D\lambda\Big)\,\chi^+_1
$$
belong to $\CR$.

Let $\tilde \psi\in C^\infty _0(\R^n\setminus 0)$ be such that $\tilde
\psi(\xi)\,\psi_2(\xi)=\psi_2(\xi)$. Then writing
$\psi_2\big(\frac D\lambda\big)=\tilde \psi\big(\frac
D\lambda\big)\,\psi_2\big(\frac D\lambda\big)$, using Lemma \ref{lVI.2.3} and
the fact that $\big\Vert \tilde \psi\big(\frac D\lambda\big)\,v\big\Vert
_{H^{2N}}\leq C\,\lambda^{2N}\,\Vert v\Vert _{L^2}$ we obtain
$$
\Vert R^+_\lambda(t)\,v\Vert _{H^{2N}}\leq C\,\lambda^{2N}\,e^{-\frac 18
\lambda\,\delta^2_2}\,\Big\Vert e^{-\frac\lambda 2
\v\alpha_\xi\v^2}\,T\,e^{-itP}\,\psi_2\Big(\frac
D\lambda\Big)\,\chi^+_1\,u\Big\Vert _{L^2}\,.
$$
Since $T$ is continuous from $L^2$ to the space of $v$ such that
$e^{-\frac\lambda 2 \v\alpha_\xi\v^2}\,v\in L^2$ and using the
conservation of $L^2$ norm we obtain
$$
\Vert R^+_\lambda(t)\,v\Vert _{H^{2N}}\leq C\,\lambda^{2N}\,e^{-\frac 18
\lambda\,\delta^2_2}\,\Big\Vert \psi_2\Big(\frac
D\lambda\Big)\,\chi^+_1\,u\Big\Vert _{L^2}\leq C'\,\lambda^{4n}\,
e^{-\frac 18\lambda\,\delta^2_2}\,\Vert u\Vert _{H^{-2N}}
$$
so our claim is proved. \cqfd
\begin{lemma}\sl \label{lVI.2.5}
Let $\chi^+_2$ be defined in (\ref{eqVI.2.2}). Then  we have for $t\in[-T,T]$
and $\lambda\geq 1$, 
$$
K_+(t)=\chi^+_1\,\psi_2 \Big(\frac
D\lambda\Big)\,\chi^+_2\,T^*\,\psi_3(\alpha_\xi)\,T\,e^{-itP}\,\chi^+_2\,\psi_2\Big(\frac
D\lambda\Big)\,\chi^+_1+R^+_\lambda(t)
$$
where $(R^+_\lambda(t))\in\CR$.
\end{lemma}

{\bf Proof } Since by (\ref{eqVI.2.2}) the support of $\chi^+_1$ and
$1-\chi^+_2$ are disjoint, the symbolic calculus shows that the operators
$\chi^+_1\,\psi_2\big(\frac D\lambda\big)(1-\chi^+_2)$ and
$(1-\chi^+_2)\,\psi_2\big(\frac D\lambda\big)\,\chi^+_1$ belong to $\frac
1{\lambda^M}\,S^{-M}_{1,0}$ for any $M\in\N$. It follows from (VI.1.7)
that if $M\geq 2N$,
\begin{equation*}
\begin{split}
\Big\Vert\chi^+_1\,\psi_2\Big(\frac
D\lambda\Big)\,&(1-\chi^+_2)\,T^*\,\psi_3(\alpha_\xi)\,T\,e^{-itP}\,\psi_2\Big(\frac
D\lambda\Big)\,\chi^+_1\,u\Big\Vert _{H^{2N}}\\
&\qquad \leq C\,\lambda^{2N} \Big\Vert \Big
(I-\frac\Delta{\lambda^2}\Big)^N\,\chi^+_1\,\psi_2\Big(\frac
D\lambda\Big)(1-\chi^+_2)\,T^*\cdots \Big\Vert _{L^2}\\
&\qquad \leq C\,\lambda^{2N-M} \Big\Vert
T^*\,\psi_3(\alpha_\xi)\,T\,e^{-itP}\,\psi_2\Big(\frac
D\lambda\Big)\,\chi^+_1\,u\Big\Vert _{L^2}\\
&\qquad \leq C\,\lambda^{2N-M} \Big\Vert e^{-\frac \lambda 2
\v\alpha_\xi\v^2}\,\psi_3(\alpha_\xi)\,T\,e^{-itP}\,\psi_2\Big(\frac
D\lambda\Big)\,\chi^+_1\,u\Big\Vert _{L^2(\R^{2n})}\\
&\qquad \leq C\,\lambda^{2N-M} \,\Big\Vert e^{-itP}\,\psi_2\Big(\frac
D\lambda\Big)\,\chi^+_1\,u\Big\Vert _{L^2}\leq C\,\lambda^{2N-M}\,\Big\Vert
\psi_2\Big(\frac D\lambda\Big)\,\chi^+_1\,u\Big\Vert _{L^2}\\
&\qquad \leq C\,\lambda^{4N-M}\,\Vert u\Vert _{H^{-2N}}\,. 
\end{split}
\end{equation*}
Taking  $M\geq 4N$ we conclude that the remainder under consideration belong to
$\CR$. By the same way
\begin{equation*}
\begin{split}
\Big\Vert\chi^+_1\,\psi_2\Big(\frac
D\lambda\Big)\,&\chi^+_2\,T^*\,\psi_3(\alpha_\xi)\,T\,e^{-itP}\,(1-\chi^+_2)\,\psi_2
\Big(\frac D\lambda\Big)\,\chi^+_1\,u\Big\Vert _{H^{2N}}\\
&\leq C\,\lambda^{2N} \Big\Vert \Big (I-\frac\Delta{\lambda^2}\Big)^N\,\psi_2
\Big(\frac D\lambda\Big)\,\chi^+_2\,T^* \cdots u\Big\Vert _{L^2}\\
&\leq C\,\lambda^{2N} \Big\Vert (1-\chi^+_2)\,\psi_2\Big(\frac
D\lambda\Big)\,\chi^+_1\Big(I-\frac\Delta{\lambda^2}\Big)^N\, \Big
(I-\frac\Delta{\lambda^2}\Big)^{-N}\,u\Big\Vert _{L^2}\\
&\leq C'\,\lambda^{2N-M}\,\Big\Vert \Big
(I-\frac\Delta{\lambda^2}\Big)^{-N}\,u\Big\Vert _{L^2}\leq
C''\,\lambda^{4N-M}\,\Vert u\Vert _{H^{-2N}}\,.
\end{split}
\end{equation*}
The proof is complete. \cqfd
\begin{lemma}\sl \label{lVI.2.6}
Let $\psi_a\in C^\infty _0(\R)$ and $\chi_a(x)$, $\chi_b(\alpha_x)$ be $C^\infty
$ functions such that one can find $\mu>0$ such that $\v x-\alpha_x\v\geq \mu$
if $(x,\alpha_x)$ belongs to $\supp  [\chi_a(x)(1-\chi_b(\alpha_x))]$. Then one
can find $\varepsilon>0$, $C>0$ such that
$$
\big\Vert
\chi_a\,T^*[\psi_a(\alpha_\xi)(1-\chi_b(\alpha_x))\,v]\big\Vert
_{L^2(\R^n)}\leq C\,e^{-\varepsilon\lambda}\, \big\Vert e^{-\frac\lambda 2
\v\alpha_\xi\v^2}\,v\big\Vert _{L^2(\R^{2n}_\alpha)}
$$
for all $v$ such that the right hand side norm is finite.
\end{lemma}

{\bf Proof } It follows from (\ref{eqVI.1.2}) that
$$
\chi_a(x)\,T^*(\psi_a(\alpha_\xi)(1-\chi_b(\alpha_x))\,v)(x)=\int
K(x,\alpha)\,e^{-\frac\lambda 2 \v\alpha_\xi\v^2}\,v(\alpha)\,d\alpha
$$
where
$$
K(x,\alpha)=c_n\,\lambda^{\frac{3n}4}\, e^{-i\lambda(x-\alpha_x)\cdot
\alpha_\xi-\frac\lambda 2 \v
x-\alpha_x\v^2}\,\psi_a(\alpha_\xi)\,\chi_a(x)(1-\chi_b(\alpha_x))\,.
$$
Using our assumption we can write
$$
\v K(x,\alpha)\v\leq C\,\lambda^{\frac{3n}4}\,e^{-\frac\lambda 4
\mu^2}\,e^{-\frac\lambda 4 \v x-\alpha_x\v^2}\,\psi_a(\alpha_\xi)\,.
$$
Therefore one can find $\varepsilon>0$ such that,
$$
\sup_x \int \v K(x,\alpha)\v\,d\alpha\leq C\,e^{-\varepsilon\lambda}\,,\quad
\sup_\alpha \int \v K(x,\alpha)\v\,dx\leq C\,e^{-\varepsilon\lambda}\,,
$$
so the Lemma follows from the well known Schur Lemma. \cqfd
\begin{corollary}\sl \label{cVI.2.7}
We have
$$
K_+(t)=\chi^+_1\,\psi\Big(\frac D
\lambda\Big)\,\chi^+_2\,T^*\,\psi_3(\alpha_\xi)\,\chi^+_3(\alpha_x)\,
T\,e^{-itP}\,\chi^+_2\,\psi_2\Big(\frac D
\lambda\Big)\,\chi^+_1+R^+_\lambda(t)
$$
where $R\in\CR$ and the same is true for the sign $-$.
\end{corollary}

{\bf Proof } We have to show that the operator
$$
R^+_\lambda(t)=\chi^+_1\,\psi_2\Big(\frac D
\lambda\Big)\,\chi^+_2\,T^*\,\psi_3(\alpha_\xi)\,(1-\chi^+_3(\alpha_x))\,
T\,e^{-itP}\,\chi^+_2\,\psi_2\Big(\frac D
\lambda\Big)\,\chi^+_1
$$
belongs to $\CR$.

We apply Lemma \ref{lVI.2.6} with $\chi_a=\chi^+_2$, $\psi_a=\psi_3$,
$\chi_b=\chi^+_3$. Then according to (\ref{eqVI.2.2}) we have on the support of
$\chi^+_2(x)(1-\chi^+_3(\alpha_x))$, $-x\cdot \xi_0\leq 2\delta_1$,
$-\alpha_x\cdot \xi_0\geq \frac 94 \,\delta_1$. Therefore we have,
$\frac94\,\delta_1\leq -\alpha_x\cdot \xi_0\leq (x-\alpha_x)\cdot
\xi_0-x\cdot\xi_0\leq \v x-\alpha_x\v\cdot
\v\xi_0\v+2\delta_1$, so $\v x-\alpha_x\v\geq \mu>0$. Then using Lemma
\ref{lVI.2.6} we can write
\begin{equation*}
\begin{aligned}
\Vert R_\lambda(t)\,u\Vert _{H^{2N}}&\leq
C\,\lambda^{2N}\,e^{-\varepsilon\lambda} \Big\Vert
e^{-\frac\lambda2\v\alpha_\xi\v^2}\,T\,e^{-itP}\,\chi^+_2\,\psi_2\Big(\frac
D\lambda\Big)\,\chi^+_1\,u\Big\Vert _{L^2}\\
&\leq C'\,\lambda^{2N}\,e^{-\varepsilon\lambda} \Big\Vert \psi_2\Big(\frac
D\lambda\Big)\,\chi^+_1\,u\Big\Vert _{L^2}\leq
C''\,\lambda^{4N}\,e^{-\varepsilon\lambda}\,\Vert u\Vert _{H^{-2N}}\,.
\end{aligned}
\end{equation*}
\cqfd

{\bf Proof of Theorem \ref{tVI.2.2} } It follows immediately from Corollary
\ref{cVI.2.7}.

\subsection{The one sided parametrix}\label{ssVI.3}

The purpose of this Section is to use the results of Sections \ref{sIV},
\ref{sV} and \ref{ssVI.2} to show that the operators $K_+(t)$ and $K_-(t)$
introduced in (\ref{eqVI.2.5}) can be written as Fourier integral operators
with complex phase functions. We shall take the expression of $K_\pm(t)$ given
by Theorem \ref{tVI.2.2} and we begin by considering the expression
$$
T\Big[e^{-itP}\,\chi^+_2\,\psi_2\Big(\frac
D\lambda\Big)\,\chi^+_1\,u\Big](\alpha,\lambda)\quad \textrm{(see
(\ref{eqVI.1.1}))}\,.
$$
Let us first introduce some other cut-off functions.
\begin{equation}\label{eqVI.3.1}\left\{
\begin{array}{l}
\textrm{If }\enskip  \v\alpha_x\cdot \alpha_\xi\v\leq 
\frac{c_0}{2}\,\bra\alpha_x\ket\,\v\alpha_\xi\v \textrm{ we set }
\chi^{\pm}_4(y)\equiv 1\,.\\
\textrm{If } \enskip \v\alpha_x\cdot \alpha_\xi\v>
\frac{c_0}{2}\,\bra\alpha_x\ket\,\v\alpha_\xi\v \textrm{ we set,}\\
\quad \textrm{(i)}\hbox to 0,5cm{}\textrm{in the $+$ case,}\enskip
\chi^+_4(y)=\chi_0\Big (-\frac{y\cdot \xi_0}{5\delta_1}\Big)\,,\\
\quad \textrm{(ii)}\quad \textrm{in the $-$ case,}\enskip \chi^-_4(y)=\chi_0\Big
(\frac{y\cdot \xi_0}{5\delta_1}\Big)\,,\\
\textrm{where } \chi_0\textrm{ has been defined (in (\ref{eqVI.2.1}))}\,. 
\end{array}\right.
\end{equation}
In all cases let $\chi_5\in C^\infty _0(\R^n)$ be such $0\leq \chi_5\leq 1$ and
\begin{equation}\label{eqVI.3.2}
\chi_5(y)=1\enskip \textrm{ if }\enskip \v y\v\leq \frac\delta 2\,,\quad \supp
\chi_5 \subset \{y : \v y\v\leq \delta\}
\end{equation}
where $\delta$ is the small constant introduced in Theorem \ref{tIV.1.2}. For
the convenience of the reader let us recall the main properties of the phase
 and the symbol constructed in Corollary 
\ref{cV.1.2}, and Theorem \ref{tIV.1.2}. First of all the phase $\varphi $ is
defined on the set $\Omega_\delta$ where,

(i) if $\v\alpha_x\cdot
\alpha_\xi\v\leq c_0\,\bra\alpha_x\ket\v\alpha_\xi\v$  then
$$
\Omega_\delta=\big\{(\theta,y)\in\R\times\R^n :\, 
 \enskip \v y-x(\theta,\alpha)\v<\delta\bra\theta\ket\big\}\,.
$$

(ii) if $\alpha_x\cdot \alpha_\xi>c_0 \bra\alpha_x\ket
\v\alpha_\xi\v$ then,
\begin{equation*}
\begin{split}
\Omega_\delta=\{(\theta,y)\in (0,+\infty )&\times\R^n : \v
y-x(\theta,\alpha)\v<\delta\bra\theta\ket\}\cup\{(\theta,y)\in
(-\infty ,0)\times\R^n :\\
&\v y-x(\theta,\alpha)\v<\delta\bra \theta\ket\enskip \textrm{ and }\enskip
y\cdot \alpha_\xi\geq -c_1 \,\bra y\ket \v\alpha_\xi\v\}\,.
\end{split}
\end{equation*}

(iii) if $\alpha_x\cdot \alpha_\xi< -c_0 \bra\alpha_x\ket \v
\alpha_\xi\v$ then,
\begin{equation*}
\begin{split}
\Omega_\delta=\{(\theta,y)\in (-\infty,0 )&\times\R^n : \v
y-x(\theta,\alpha)\v<\delta \bra\theta\ket\}\cup\{(\theta,y)\in
(0,+\infty )\times\R^n :\\
&\v y-x(\theta,\alpha)\v<\delta\bra \theta\ket\enskip \textrm{ and }\enskip
y\cdot \alpha_\xi\leq c_1 \,\bra y\ket \v\alpha_\xi\v\}\,.
\end{split}
\end{equation*}
Moreover on this domain
\begin{equation}\label{eqVI.3.3}
\Im \varphi (\theta,y,\alpha)\geq \frac 14 \,\frac{\v
y-x(\theta,\alpha)\v^2}{1+4\theta^2}-\frac 12\,\v\alpha_\xi\v^2\,.
\end{equation}
Now if we set, with the notations of Theorem \ref{tV.1.1}, for $N\in\N,$
\begin{equation}\label{eqVI.3.4}
a(\theta,y,\alpha,\lambda)=(1+\theta^2)^{-\frac n4}\,
e_N(\theta,y-x(\theta,\alpha),\alpha,\lambda)\,,
\end{equation}
then $a$ is defined on $\Omega_\delta$ for $\lambda\geq 1$ and satisfies,
\begin{equation}\label{eqVI.3.5}\left\{
\begin{array}{l}
a(0,y,\alpha,\lambda)=1\,,\\
\v a(\theta,y,\alpha,\lambda)\v\leq c\,(1+\theta^2)^{-\frac n 4}\,,\\
\Big (i\lambda\,\frac\partial {\partial \theta}+{}^tP\Big )\big
(e^{i\lambda\varphi (\theta,y,\alpha)}\,a(\theta,y,\alpha,\lambda)\big
)=b_N(\theta,y,\alpha,\lambda)\, e^{i\lambda\varphi (\theta,y,\alpha)}\,,\\
\textrm{with } \v b_N(\theta,y,\alpha)\v\leq c_N (1+\theta^2)^{-\frac n 4} \Big
(\lambda^{-N}+\lambda^2 \Big (\frac{\v
y-x(\theta,\alpha)\v}{\bra\theta\ket}\Big)^N\Big )\,.
\end{array}\right.
\end{equation}
Let us introduce now the following set.
\begin{equation}\label{eqVI.3.6}
\begin{split}
W^{\pm}=\Big\{\alpha\in T^*\R^n : \frac 12 \leq \v\alpha_\xi\v\leq
2\,,\enskip \v\alpha_x\cdot \alpha_\xi\v\leq
c_0\bra\alpha_x\ket\v\alpha_\xi\v\Big\}\cup \big\{\alpha\in T^*\R^n :\\
\v\alpha_x\cdot
\alpha_\xi\v>c_0\,\bra\alpha_x\ket\v\alpha_\xi\v\,,\enskip
(\alpha_x,\alpha_\xi)\in\supp (\chi^\pm_3(\alpha_x)\cdot
\psi_3(\alpha_\xi))\big\}
\end{split}
\end{equation}
where $\chi^\pm_3$ and $\psi_3$ have been defined in (\ref{eqVI.2.2}) and
(\ref{eqVI.2.9}). 

Then we can state the main result of this Section.
\begin{theorem}\sl \label{tVI.3.1}
We have for $t\in[-T,T]$, and $\alpha\in W^+$,
\begin{equation*}
\begin{split}
T\big[e^{-itP}\,\chi^+_2v\big](\alpha,\lambda)=\lambda^{\frac {3n}4} \int
e^{i\lambda\varphi (-\lambda t,y,\alpha)}\,a(-\lambda
t,y,\alpha,\lambda)\,\chi^+_4(y)\\
\chi_5\Big (\frac{y-x(-\lambda t,\alpha)}{\bra\lambda t\ket}\Big) [\chi^+_2
v](y)\,dy+J^+_\lambda(t)\,v(\alpha)
\end{split}
\end{equation*}
where the operator $J^+_\lambda$ is such that, for any $M\in\N$ one can find a
constant $C_M>0$ such that for all $\lambda\geq 1$ and $t\in [-T,T]$
$$
\big\Vert e^{-\frac\lambda 2
\v\alpha_\xi\v^2}\,J^+_\lambda(t)\,v\big\Vert _{L^2(W^+)}\leq
\frac{C_M}{\lambda^M} \Vert v\Vert _{L^2(\R^n)}\,,
$$
and the same is true with the minus sign.
\end{theorem}

{\bf Proof } Let us introduce the following family of operators. We set for
$\alpha\in W^+$
\begin{equation}\label{eqVI.3.7}
Sv(\theta,t,\alpha,\lambda)=\lambda^{\frac {3n}4} \int_{\R^n}
e^{i\lambda\varphi
(\theta,y,\alpha)}\,a(\theta,y,\alpha,\lambda)\,\chi^+_4(y)\,\chi_5\Big
(\frac{y-x(\theta,\alpha)}{\bra\theta\ket}\Big)
\big[e^{-itP}\,\chi^+_2 v\big]\,(y)\,dy\,.
\end{equation}
We must verify that the right hand side is indeed well defined.

On the support of $\chi_5$ we have $\v
y-x(\theta,\alpha)\v<\delta\bra\theta\ket$ (which is one of the conditions for
$(\theta,y)$ to be in $\Omega_\delta$). If $\alpha\in W^+$ then either
$\v\alpha_x\cdot \alpha_\varepsilon\v\leq c_0
\bra\alpha_x\ket\v\alpha_\varepsilon\v$ or $\v\alpha_x\cdot
\alpha_\varepsilon\v>c_0\bra\alpha_x\ket\v\alpha_\varepsilon\v$ and
$(\alpha_x,\alpha_\varepsilon)\in\supp (\chi^+_2\cdot \psi_3)$. In the first
case by (\ref{eqVI.3.1}) $\chi^+_4(y)=1$ but $(\theta,y)\in\Omega_\delta$ for
$\theta\in\R$. In the second case since $\alpha_x\in\supp \chi^+_3$ and
$\alpha_\xi\in\supp \psi_3$ we have, by (\ref{eqVI.2.1}),
(\ref{eqVI.2.2}) and (\ref{eqVI.2.9}), $\alpha_x\cdot \xi_0\geq -3\delta_1$ and
$\big\v\frac{\alpha_\xi}{\v\alpha_
\xi\v}+\xi_0\big\v\leq 4\delta_2$.
It follows that 
\begin{equation}\label{eqVI.3.8}
\alpha_x\cdot \frac{\alpha_\xi}{\v\alpha_\xi\v}=-\alpha_x\cdot
\xi_0+\alpha_x\cdot
\Big(\frac{\alpha_\xi}{\v\alpha_\xi\v}+\xi_0\Big)\leq
7\,\delta_1\,\bra\alpha_x\ket\,.
\end{equation}
On the other hand we have $\v\alpha_x\cdot
\alpha_\xi\v>c_0\,\bra\alpha_x\ket\v\alpha_\xi\v$. According to
(\ref{eqVI.3.8}) we cannot have $\alpha_x\cdot
\alpha_\xi>c_0\,\bra\alpha_x\ket\v\alpha_\xi\v$ if
$\delta_1\ll c_0$~; thus we must have $\alpha_x\cdot \alpha_\xi<-c_0
\bra\alpha_x\ket\v\alpha_\xi\v$. So we are in the case (iii) for the
definition of $\Omega_\delta$. Since in the integral, on the support of
$\chi^+_4(y)$ we have by (\ref{eqVI.3.1}), $y\cdot \xi_0\geq
-5\,\delta_1\geq -5\,\delta_1\,\bra y\ket$, we deduce that
$$
y\cdot \frac{\alpha_\xi}{\v\alpha_\xi\v}=-y\cdot
\xi_0+y\cdot
\Big(\frac{\alpha_\xi}{\v\alpha_\xi\v}+\xi_0\Big)\leq
9\,\delta_1 \bra y\ket<c_0 \bra y\ket\,.
$$
Therefore $(\theta,y)\in\Omega_\delta$ and the right hand side of
(\ref{eqVI.3.7}) is well defined. Of course the same argument is valid in the
case of the minus sign.

Now let us set for $s\in [0,t]$ if $t>0$ (resp. $s\in[t,0]$ if $t<0$)
\begin{equation}\label{eqVI.3.9}
g(s)=Sv(\lambda (s-t),s,\alpha,\lambda)\,.
\end{equation}
 Since $g'(s)=(\lambda\partial _\theta+\partial _t)(Sv)(\lambda
(s-t),s,\alpha,\lambda)$ we obtain
\begin{equation}\label{eqVI.3.10}\left\{
\begin{array}{l}
Sv(0,t,\alpha,\lambda)=Sv(-\lambda t,0,\alpha,\lambda)+\int^t_0(\lambda\partial
_\theta+\partial _t)(Sv)(\lambda(s-t),s,\alpha,\lambda)\,ds\enskip \textrm{ if
}\enskip  t>0\\ Sv(0,t,\alpha,\lambda)=Sv(-\lambda
t,0,\alpha,\lambda)-\int^t_0(\lambda\partial _\theta+\partial
_t)(Sv)(\lambda(s-t),s,\alpha,\lambda)\,ds\enskip \textrm{ if }\enskip  t<0\,.
\end{array}\right.
\end{equation}
Let us set
\begin{equation}\label{eqVI.3.11}
U(t,y,\lambda)=\big[e^{-itP}\,\chi^+_2v\big](y)\,.
\end{equation}
Then according to (\ref{eqVI.3.7}) we can write
\begin{equation}\label{eqVI.3.12}
(\lambda\partial _\theta+\partial
_t)(Sv)(\theta,t,\alpha,\lambda)=\lambda^{\frac{3n}4}\,(A_1+A_2+A_3)
\end{equation}
where
\begin{equation}\label{eqVI.3.13}\left\{
\begin{array}{l}
A_1=\int \lambda\partial _\theta\big (e^{i\lambda\varphi
(\theta,y,\alpha)}\,a(\theta,y,\alpha,\lambda)\big )\,\chi^+_4(y)\,\chi_5 \Big
(\frac{y-x(\theta,\alpha)}{\bra\theta\ket}\Big )\,U(t,y,\lambda)\,dy\\
A_2=\int e^{i\lambda\varphi
(\theta,y,\alpha)}\,a(\theta,y,\alpha,\lambda)\,\chi^+_4(y)\,\chi_5 \Big
(\frac{y-x(\theta,\alpha)}{\bra\theta\ket}\Big )(-iP)\,U(t,y,\lambda)\,dy\\
A_3=\som^n_{j=1} \int e^{i\lambda\varphi
(\theta,y,\alpha)}\,a(\theta,y,\alpha,\lambda)\,\chi^+_4(y)\,(\partial
_j\,\chi_5)
\Big (\frac{y-x(\theta,\alpha)}{\bra\theta\ket}\Big )\,\partial _\theta\Big
(\frac{y_j-x_j(\theta,\alpha)}{\bra\theta\ket}\Big)\,U(t,y,\lambda)\,dy\,.
\end{array}\right.
\end{equation}
Integrating by parts in the integral giving $A_2$ we obtain
\begin{equation*}
\begin{split}
A_2&=\int (-i\,^tP) \big[e^{i\lambda\varphi
(\theta,y,\alpha)}\,a(\theta,y,\alpha,\lambda)\big]\,\chi^+_4 (y)\,\chi_5\Big
(\frac{y-x(\theta,\alpha)}{\bra\theta\ket}\Big )\,U(t,y,\lambda)\,dy\\
&+\enskip \sum_{1\leq \v\beta\v\leq 2\atop \beta=\beta_1+\beta_2} \int
e^{i\lambda\varphi (\theta,y,\alpha)}b_\beta(\theta,y,\alpha,\lambda)(\partial
^{\beta_1}_y\,\chi^+_4)(y)\bra\theta\ket^{-\v\beta_2\v}\,(\partial
^{\beta_2}_y\,\chi_5)\Big(\frac{y-x(\theta,\alpha)}{\bra\theta\ket}\Big)\,
U(t,y,\lambda)\,dy
\end{split}
\end{equation*}
where $b_\beta$ is a symbol satisfying the same estimates as
$a(\theta,y,\alpha,\lambda)$ (see (\ref{eqVI.3.5})).

Combining the term $A_1$ with the first part of $A_2$ and using
(\ref{eqVI.3.5}) we can write
\begin{equation}\label{eqVI.3.14}
(\lambda\,\partial _\theta+\partial
_t)(Sv)(\theta,t,\alpha,\lambda)=\lambda^{\frac{3n}4} (B_1+B_2+B_3)
\end{equation}
where
\begin{equation}\label{eqVI.3.15}\left\{
\begin{array}{l}
B_1(\theta,t,\alpha,\lambda)=\int e^{i\lambda\varphi
(\theta,y,\alpha)}\,b_N(\theta,y,\alpha,\lambda)\,\chi^+_4(y)\,\chi_5\Big(
\frac{y-x(\theta,\alpha)}{\bra\theta\ket}\Big)\,U(t,y,\lambda)\,dy\\
B_2(\theta,t,\alpha,\lambda)=\som_{\v\beta_1+\beta_2\v\leq 2\atop
\beta_2\not =0}\int e^{i\lambda\varphi
(\theta,y,\alpha)}\, C_{\beta_1\beta_2}(\theta,y,\alpha,\lambda)(\partial
^{\beta_1}_y\,\chi^+_4)(y)\bra\theta\ket^{-\v\beta_2\v}\\
\hbox to 6cm{}(\partial
^{\beta_2}_y\,\chi_5)\Big(
\frac{y-x(\theta,\alpha)}{\bra\theta\ket}\Big)\,U(t,y,\lambda)\,dy\\
B_3(\theta,t,\alpha,\lambda)=\som_{1\leq \v\gamma\v\leq 2}\int
e^{i\lambda\varphi (\theta,y,\alpha)}\,d_\gamma
(\theta,y,\alpha,\lambda)(\partial
^{\gamma}_y\,\chi^+_4)(y)\,\\
\hbox to 6cm{}\chi_5\Big(
\frac{y-x(\theta,\alpha)}{\bra\theta\ket}\Big)\,U(t,y,\lambda)\,dy
\end{array}\right.
\end{equation}
where $C_{\beta_1,\beta_2}$ and $d_\gamma$ are bounded symbol. Here we used the
estimate $\big\v\partial
_\theta\big(\frac{y_j-x_j(\theta,\alpha)}{\bra\theta\ket}\big)\big\v\leq \frac
C{\bra\theta\ket}$. Now we state a Lemma.
\begin{lemma}\sl \label{lVI.3.2}
One can find constants $C>0$, $\eta>0$ and for every $M\geq 0$ a constant
$C_M>0$ such that for all $\theta,t,\lambda$ such that $\lambda\geq 1$,
$\v\theta\v\leq\v \lambda t\v$, $\v t\v\leq T$ we have
\begin{equation*}
\begin{aligned}
\big\Vert e^{-\frac\lambda 2
\v\alpha_\xi\v^2}\,B_1(\theta,t,\cdot ,\lambda)\big\Vert
_{L^2(W^+)}&\leq C_M\,\lambda^{-M}\,\Vert v\Vert _{L^2(\R^n)}\\
\big\Vert e^{-\frac\lambda 2
\v\alpha_\xi\v^2}\,B_2(\theta,t,\cdot ,\lambda)\big\Vert
_{L^2(W^+)}&\leq C\,e^{-\eta\lambda}\,\Vert v\Vert _{L^2(\R^n)}
\end{aligned}
\end{equation*}
for all $v$ in $L^2(\R^n)$.
\end{lemma}

{\bf Proof } Let us consider the term $B_1$. Using (\ref{eqVI.3.5}) and
(\ref{eqVI.3.3}) we obtain
\begin{equation*}
\begin{split}
&\Big\v e^{i\lambda\varphi (\theta,y,\alpha)-\frac\lambda 2
\v\alpha_\xi\v^2}\,b_N(\theta,y,\alpha,\lambda)\,\chi_5
\Big(\frac{y-x(\theta,\alpha)}{\bra\theta\ket}\Big)\,\chi^+_4(y)\,\boldsymbol{1}
_{W^+}(\alpha)\Big\v\\
&\qquad \leq C_N\, e^{-\frac\lambda {16} \frac{\v
y-x(\theta,\alpha)\v^2}{\bra\theta\ket^2}}\,\bra\theta\ket^{-\frac n2} \Big
(\lambda^{-N}+\lambda^2\Big (\frac{\v
y-x(\theta,\alpha)\v}{\bra\theta\ket}\Big)^N\Big
)\,\boldsymbol{1}_{W^+(\alpha)}=:K(\theta,y,\alpha,\lambda)\,.   
\end{split}
\end{equation*}
On the other hand we have
$$
e^{-\frac\lambda 2 \v\alpha_\xi\v^2}\,B_1=\int
K(\theta,y,\alpha,\lambda)\,U(t,y,\lambda)\,dy\,.
$$
We want to apply Schur Lemma to this integral operator. First of all making the
change of variables $y-x(\theta,\alpha)=\frac{\bra\theta\ket\,z}{\sqrt\lambda}$
we can write
$$
\int K(\theta,y,\alpha,\lambda)\,dy\leq C_N\,\bra\theta\ket^{\frac n2} \int
e^{-\frac {\v z\v^2}{16}}\,\lambda^{-\frac n2}
(\lambda^{-N}+\lambda^2\,\lambda^{-\frac N 2}\,\v z\v^N)\,dz\,.
$$
Since $\bra\theta\ket\leq 1+\lambda T$ we obtain finally for all $M\geq 0$ and
$\lambda\geq 1$, $\int K(\theta,y,\alpha,\lambda)\,dy\leq C_M\,\lambda^{-M}$.

To estimate $(1)=\int K(\theta,y,\alpha,\lambda)\,d\alpha$ we make the change
of variable $\beta=(x(\theta,\alpha),\xi(\theta,\alpha))$ Since this
transformation is symplectic we have $d\alpha=d\beta$. Therefore we obtain,
$$
(1)\leq C_N \int e^{-\frac\lambda{16}\,\frac{\v
y-\beta_x\v^2}{\bra\theta\ket^2}}\,\bra\theta\ket^{-\frac n2}\, \boldsymbol{1}
_{\frac 12 \leq \v\xi(-\theta,\beta)\v\leq 2}\,(\beta)
\left(\lambda^{-N}+\lambda^2\Big (\frac{\v
y-\beta_x\v}{\bra\theta\ket}\Big)^N\right)\,d\beta\,.
$$
Since $\xi(-\theta,\beta)=\beta_\xi+\CO(\varepsilon)$ we have $\frac
13\leq \v\beta_\xi\v\leq 3$. Setting as before
$y-\beta_x=\frac{\bra\theta\ket}{\sqrt\lambda}$ we deduce easily that $(1)\leq
C_M\,\lambda^{-M}$ for every $M\geq 0$. It follows from these estimates, the
Schur Lemma and (\ref{eqVI.3.11}) that,
$$
\big\Vert e^{-\frac\lambda 2\v\alpha_\xi\v^2}\,B_1\big\Vert
_{L^2(W^+)}\leq C_M\,\lambda^{-M} \big\Vert
e^{-itP}\,\psi^2\Big(\frac{P}{\lambda^2}\Big)\,\chi^+_2
v\big\Vert _{L^2(\R^n)}\leq C'_M\,\lambda^{-M}\,\Vert v\Vert _{L^2(\R^n)}\,. 
$$
To deal with the term $B_2$ we use exactly the same computations and the fact
that on the support of $(\partial
^{\beta_2}_y\,\chi_5)\big(\frac{y-x(\theta,\alpha)}{\bra\theta\ket}\big)$, for
$\v\beta_2\v\geq 1$ we have $\frac\delta 2\leq \frac{\v
y-x(\theta,\alpha)\v}{\bra\theta\ket}\leq \delta$.\cqfd

{\bf Proof of Theorem \ref{tVI.3.1} in the following cases}
\begin{equation}\label{eqVI.3.16}\left\{
\begin{array}{ll}
\textrm{(i)}\quad &\v\alpha_x\cdot \alpha_\xi\v\leq
c_0\,\bra\alpha_x\ket\,\v\alpha_\xi\v\enskip \textrm{and} \enskip t\in
[-T,T]\,,\\
\textrm{(ii)}\quad &\v\alpha_x\cdot \alpha_\xi\v>
c_0\,\bra\alpha_x\ket\,\v\alpha_\xi\v\enskip \textrm{and} \enskip t\in
[0,T]\,.
\end{array}\right.
\end{equation}
In both these cases we are going to show, with the notations of
(\ref{eqVI.3.15}), that $B_3\equiv 0$.

In the case (i) this is obvious since by (\ref{eqVI.3.1}) we took
$\chi^+_4(y)\equiv 1$. Now if $t>0$, by (\ref{eqVI.3.9}) we have
$\theta=\lambda(s-t)\leq 0$.
In the case (ii) we have seen that in the case $+$ we must have $\alpha_x\cdot
\alpha_\varepsilon\leq -c_0\,\bra\alpha_x\ket\,\v\alpha_\varepsilon\v$. Now on
the support of $\chi_5\big (\frac{y-y(\theta,\alpha)}{\bra\theta\ket}\big)$ we
have $\v y-x(\theta,\alpha)\v\leq \delta\bra\theta\ket$ and it follows from
Proposition \ref{pIII.4.1} that
$x(\theta,\alpha)=\alpha_x+2\theta\,\alpha_\xi+\CO
(\varepsilon\bra\theta\ket)$. Then we write
$$y\cdot\xi_0=(y-x(\theta,\alpha))\cdot \xi_0+(\alpha_x+
2\theta\,\alpha_\xi)\cdot
\xi_0+\CO(\varepsilon\bra\theta\ket)\,.$$
 On the support of
$\psi_3(\alpha_\xi)$ we have
$\big\v\frac{\alpha_\xi}{\v\alpha_\xi\v}+\xi_0\big\v\leq
4\,\delta_1$ so
$$
y\cdot \xi_0=(\alpha_x+2\theta\,\alpha_\xi)\cdot
\Big(\xi_0+\frac{\alpha_\xi}{\v\alpha_\xi\v}\Big
)-(\alpha_x+2\theta\,\alpha_\xi)\cdot
\frac{\alpha_\xi}{\v\alpha_\xi\v}+\CO
((\varepsilon+\delta)\bra\theta\ket)
$$
$$
y\cdot \xi_0\geq
c_0\,\bra\alpha_x\ket+2\,\v\theta\v\,\v\alpha_\xi\v-4\,\delta_1\,\bra
\alpha_x\ket-C
(\delta+\varepsilon+\delta_1)\bra\theta\ket\,.
$$
Since $\delta,\varepsilon,\delta_1$ are small compared to $c_0$ we deduce that
$y\cdot \xi_0\geq \frac{c_0}2\,\bra\theta\ket$ in the integral defining
$B_3$ in (\ref{eqVI.3.15}). Since the support of $\partial
^\gamma_y\,\chi^+_4(y)$, for $\gamma\not =0$, is contained in $\frac 34\leq
\frac{-y\cdot \xi_0}{5\,\delta_1}\leq 1$ we deduce that $B_3\equiv 0$
in this case, (see (\ref{eqVI.3.1}) and (\ref{eqVI.2.2})). It follows from
(\ref{eqVI.3.10}) and (\ref{eqVI.3.14}) that
\begin{equation}\label{eqVI.3.17}
Sv(0,t,\alpha,\lambda)=Sv(-\lambda
t,0,\alpha,\lambda)+\int^t_0 (B_1+B_2)(s)\,ds\,.
\end{equation}
It follows from Lemma \ref{lVI.3.2} that,
\begin{equation}\label{eqVI.3.18}
\Big\Vert e^{-\frac\lambda 2 \v\alpha_\xi\v^2} \int^t_0
(B_1+B_2)(s)\,ds\Big\Vert_{L^2(W^+)}\leq \frac C {\lambda^M}\,\Vert v\Vert
_{L^2(\R^n)}\,.
\end{equation}
Now we have
\begin{equation*}\left\{
\begin{array}{l}
\varphi (0,y,\alpha)=\varphi _0(y,\alpha)+g(y-\alpha_x)\,,\textrm{ where}\\
\varphi _0(y,\alpha)=(y-\alpha_x)\cdot \alpha_\xi+\frac i 2 \,\v
y-\alpha_x\v^2+\frac 1{2i}\,\v \alpha_\xi\v^2\,,\\
\v g(x)\v\leq C_N\,\v x\v^N\textrm{ for every } N\in\N\,. 
\end{array}\right.
\end{equation*}
Let $\chi\in C^\infty _0(\R^n)$ be such that $\chi(x)=1$ if $\v x\v\leq 1$,
$\supp \chi\subset \{x : \v x\v\leq 2\}$ and let us fix $N\geq 3$. We can write
\begin{equation}\label{eqVI.3.19}
Sv(0,t,\alpha,\lambda)=A_1-A_2-A_3+A_4
\end{equation}
where
\begin{equation}\label{eqVI.3.20}\left\{
\begin{array}{l}
A_1=c_n\,\lambda^{\frac{3n}4} \int e^{i\lambda\varphi
_0(y,\alpha)}\,\chi^+_4(y)\,\chi_5(y-\alpha_x)\,U(t,y,\lambda)\,dy\\
A_2=c_n\,\lambda^{\frac{3n}4} \int e^{i\lambda\varphi
_0(y,\alpha)}\, (1-\chi(C_N\lambda\v
y-\alpha_x\v^N))\,\chi^+_4(y)\,\chi_5(y-\alpha_x)\,U(t,y,\lambda)\,dy\\
A_3=c_n\,\lambda^{\frac{3n}4} \int e^{i\lambda\varphi
_0(y,\alpha)}\, (1-e^{i\lambda g(y-\alpha_x)})\,\chi(C_N\,\lambda\v
y-\alpha_x\v^N)\,\chi^+_4(y)\\
\hbox to 6,5cm{}\cdot \chi_5(y-\alpha_x)\,U(t,y,\lambda)\,dy\\
A_4=c_n\,\lambda^{\frac{3n}4} \int e^{i\lambda(\varphi
_0(y,\alpha)+g(y-\alpha_x))}\,(1-\chi(C_N\,\lambda \v
y-\alpha_x\v^N))\,\chi^+_4(y)\\
\hbox to 6,5cm{}\cdot \chi_5(y-\alpha_x)\,U(t,y,\lambda)\,dy\,.
\end{array}\right.
\end{equation}
We claim that we have for $j=2,3,4$,
\begin{equation}\label{eqVI.3.21}
\big\Vert e^{-\frac\lambda 2 \v\alpha_\xi\v^2}\,A_j\big\Vert
_{L^2(W^+)}\leq \frac{C_{M_N}}{\lambda^{M_N}} \,\Vert v\Vert _{L^2(\R^n)}\,,
M_N\rightarrow +\infty \textrm{ if } N\rightarrow +\infty \,.
\end{equation}

(i) Term $A_2$

On the support $1-\chi(C_N\lambda \v y-\alpha_x\v^N)$ we have $\v
y-\alpha_x\v\geq C'_N\, \lambda^{-\frac 1 N}$. So,
$$
\Big\v e^{-\frac\lambda 2 \v\alpha_\xi\v^2} \, e^{i\varphi
_0(y,\alpha)}\,(1-\chi(C_N\lambda\v y-\alpha_x\v^N))\Big\v\leq
C\,e^{-\frac\lambda 4 \v y-\alpha_x\v^2}\, e^{-C''_N\,\lambda^{1-\frac 2N}}\,.
$$
Using the Schur Lemma and the inequality $\Vert U(t,\cdot ,\alpha)\Vert
_{L^2}\leq C\,\Vert v\Vert _{L^2}$ (see (\ref{eqVI.3.11})) we obtain
(\ref{eqVI.3.21}) for $A_2$.

(ii) Term $A_3$

On the support of $\chi(C_N\lambda\v y-\alpha_x\v^N)$ we have $\lambda\v
g(y-\alpha_x)\v\leq C_N \lambda\v y-\alpha_x\v^N\leq 2$. Therefore we have $\v
1-e^{i\lambda g(y-\alpha_x)}\v\leq C\,\lambda\,\v g(y-\alpha_x)\v\leq
C'_N\,\lambda \v y-\alpha_x\v^N$. It follows that
\begin{equation*}
\begin{split}
\Big\v e^{-\frac\lambda 2 \v\alpha_\xi\v^2}\, e^{i\lambda\varphi
_0(y,\alpha)}\,(1-e^{i\lambda g(y-\alpha_x})\,\chi(C_N\,\lambda \v
y-\alpha_x\v^N)\Big\v&\leq C'_N~e^{-\frac\lambda 2 \v y-\alpha_x\v^2}\cdot
\lambda \v y-\alpha_x\v^N\\
&\leq \frac{C'_N}{\lambda^{\frac N2 -1}}\,(\lambda^{\frac 12} \v
y-\alpha_x\v)^N\, e^{-\frac \lambda 4 \v y-\alpha_x\v^2}\, e^{-\frac\lambda 4
\v y-\alpha_x\v^2}\\
&\leq \frac{C'_N}{\lambda^{\frac N2 -1}}\, e^{-\frac\lambda 4
\v y-\alpha_x\v^2}\,.
\end{split}
\end{equation*}
The Schur Lemma shows again that $A_3$ satisfies (\ref{eqVI.3.21}).

(iii) Term $A_4$

On the support of $\chi_5(y-\alpha_x)$ we have, according to (\ref{eqVI.3.2}),
$\v y-\alpha_x\v \leq \delta$. It follows that $\v g(\alpha_x-y)\v\leq
\delta\,C_3
\v y-\alpha_x\v^2$ so if $\delta$ is small enough,
$$
-\frac 12 \,\v\alpha_\xi\v^2-\Im \varphi _0 (y,\alpha)-\Im
g(y-\alpha_x)\leq -\frac 14 \,\v y-\alpha_x\v^2
$$
and, as for $A_2$, the Schur Lemma implies that $A_4$ satisfies
(\ref{eqVI.3.21}).

Using (\ref{eqVI.3.19}) and (\ref{eqVI.3.21}) we see that
\begin{equation}\label{eqVI.3.22}
Sv(0,t,\alpha,\lambda)=c_n \,\lambda^{\frac {3n}4} \int e^{i\lambda\varphi
_0(y,\alpha)}\,\chi^+_4(y)\,\chi_5(y-\alpha_x)\,U(t,y,\lambda)\,dy+J^+_\lambda(t)\,v(\alpha)\,.
\end{equation}
\begin{equation}\label{eqVI.3.23}
\big\Vert e^{-\frac\lambda 2
\v\alpha_\xi\v^2}\,J^+_\lambda(t)\,v\big\Vert \leq
C_M\,\lambda^{-M}\,\Vert v\Vert _{L^2(\R^n)}\,.
\end{equation}
Now, on the support of $\chi_5 (y-\alpha_x)-1$ we have $\v y-\alpha_x\v\geq
\frac\delta2$ so modulo a term which satisfies (\ref{eqVI.3.23}) we can remove
$\chi_5(y-\alpha_x)$ in the right hand side of (\ref{eqVI.3.22}). Let us remove
$\chi^+_4$. When $\alpha\in W^+$ we have $\v\alpha_x\cdot
\alpha_\varepsilon\v\leq c_0 \bra\alpha_x\ket\v\alpha_\xi\v$ and
$\chi^+_4(y)\equiv 1$ (see (\ref{eqVI.3.1})) (so there is nothing to remove) or
$\alpha_x\in\supp \chi^+_3$ (see (\ref{eqVI.3.6}))  that is $-\alpha_x\cdot
\xi_0\leq 3\delta_1$. In the later case on the support of
$1-\chi^+_4(y)$ we have $\frac{-y\cdot \xi_0}{5\delta_1}\geq \frac 34$
(see (\ref{eqVI.2.1}), (\ref{eqVI.3.1})) so $\v y-\alpha_x\v\geq \alpha_x\cdot
\xi_0-y\cdot \xi_0\geq \frac 34\,\delta_1$. The corresponding
term, again by the Schur Lemma, satisfies (\ref{eqVI.3.23}).

Using (\ref{eqVI.1.1}) we see therefore that
\begin{equation}\label{eqVI.3.24}
Sv(0,t,\alpha,\lambda)=T\,[e^{-itP}\,\chi^+_2
v](\alpha,\lambda)+J^+_\lambda(t)\,v(\alpha)\,,
\end{equation}
where $J^+_\lambda(t)$ satisfies (\ref{eqVI.3.23}).

Gathering the informations given by (\ref{eqVI.3.17}), (\ref{eqVI.3.18}),
(\ref{eqVI.3.24}) and (\ref{eqVI.3.23}) we obtain the claim of Theorem
\ref{tVI.3.1} in the case (\ref{eqVI.3.16}). 

{\bf Proof of  Theorem \ref{tVI.3.1} in the following case}
\begin{equation}\label{eqVI.3.25}
\v\alpha_x\cdot \alpha_\xi\v>c_0 \bra\alpha_x\ket \v\alpha_\xi\v
\enskip \textrm  { and }\enskip  t\in[-T,0]\,.
\end{equation}
According to (\ref{eqVI.3.14}), (\ref{eqVI.3.15}) and Lemma \ref{lVI.3.2}, we
must prove that for all $N\in\N$ one can find $C_N>0$ such that for $\lambda\geq
1$,
\begin{equation}\label{eqVI.3.26}
\big\Vert e^{-\frac \lambda 2
\v\alpha_\xi\v^2}\,B_3(\theta,t,\cdot ,\lambda)\big\Vert _{L^2(W^+)}\leq
C_N\,\lambda^{-N}\,\Vert v\Vert _{L^2(\R^n)}\,.
\end{equation}
Here $\theta=\lambda(s-t)>0$ since $s\in [t,0]$.

Let us introduce a new cut-off function. Let $\psi_4\in C^\infty _0(\R^n)$ be
such that $0\leq \psi_4\leq 1$ and,
\begin{equation}\label{eqVI.3.27}\left\{
\begin{array}{l}
\psi_4(\xi)=1\enskip \textrm { if }\enskip
\Big\v\frac\xi{\v\xi\v}+\xi_0\Big\v\leq
5\delta_2\,,\enskip a-5\delta_2\leq \v\xi\v\leq b+5\delta_2\\
\supp \psi_4 \subset \Big\{\xi :
\Big\v\frac\xi{\v\xi\v}+\xi_0\Big\v\leq 6\delta_2\,,\enskip
a-6\delta_2\leq \v\xi\v\leq b+6\delta_2\Big\}\,.
\end{array}\right.
\end{equation}
We state a Lemma.
\begin{lemma}\sl\label{lVI.3.3}
Let us set
\begin{equation}\label{eqVI.3.28}
\tilde W^+=\left\{\alpha : \frac 12 \leq \v\alpha_\xi\v\leq 2\,,\enskip
\v\alpha_x\cdot \alpha_\xi\v>c_0
\bra\alpha_x\ket\v\alpha_\xi\v\,,\enskip
(\alpha_x,\alpha_\xi)\in
\supp (\chi^+_3(\alpha_x)\,\psi_3(\alpha_\xi))\right\}\,.
\end{equation}
Let $k(\theta,y,\alpha,\lambda)$ be a symbol and let us set
$$
F(\theta,\alpha,\lambda)=\boldsymbol{1}_{\tilde W^+}(\alpha)\,e^{-\frac\lambda
2 \v\alpha_\xi\v^2} \int e^{i\lambda\varphi
(\theta,y,\alpha)}\,k(\theta,y,\alpha,\lambda)\,\chi_5\Big(\frac{y-x(\theta,\alpha)}
{\bra\theta\ket}\Big)
 \partial ^\gamma_y\,\chi^+_4(y)\Big[I-\psi_4\Big(\frac
D\lambda\Big)\Big]\,v(y)\,dy\,.
$$
Then for every $N\in\N$ one can find $C_N>0$ such that for $\lambda\geq 1$ and
$\v\theta\v\leq \lambda T$ we have
$$
\Vert F(\theta,\cdot ,\lambda)\Vert _{L^2}\leq C_N\,\lambda^{-N}\,\Vert v\Vert
_{L^2(\R^n)}\,.
$$ 
\end{lemma}

{\bf Proof } By (VI.2.12) we have $\v\xi+\alpha_\xi\v\geq
\mu>0$ on the support of $\psi_3(\alpha_\xi)(1-\psi_4(\xi))$. Now
recall that Theorem \ref{tIV.1.2} shows that the phase $\varphi $
 satisfies
\begin{equation}\label{eqVI.3.29}\left\{
\begin{array}{l}
\Big\v\frac{\partial \varphi }{\partial
y}\,(\theta,y,\alpha)-\alpha_\xi\Big\v\leq C\,(\varepsilon+\sqrt\delta
)\\
\v \partial ^\beta_y\,\varphi (\theta,y,\alpha)\v\leq C_\beta\enskip \textrm{
if }\v\beta\v\geq 1
\end{array}\right.
\end{equation}
on the support of $\chi_5\big (\frac{y-x(\theta,\alpha)}{\bra\theta\ket}\big
)\,\partial ^\gamma_y\,\chi^+_4(y)\,\boldsymbol{1}_{\tilde W^+}(\alpha)$.

Let $g\in C^\infty _0(\R^n)$ be such that $g(\xi)=1$ if $\v\xi\v\leq 1$. Then
$$
\Big(I-\psi_4\Big(\frac D\lambda\Big)\Big)\,v\,(y)=\lim_{\varepsilon\rightarrow
0}\Big (\frac\lambda{2\pi }\Big )^n \iint e^{i\lambda (y-z)\cdot
\xi}\,(1-\psi_4(\xi))\,g(\varepsilon\,\xi)\,v(z)\,dz\,d\xi\,.
$$
It follows that
\begin{equation}\label{eqVI.3.30}\left\{
\begin{array}{l}
{\displaystyle F(\theta,\alpha,\lambda)=\lim_{\varepsilon\rightarrow 0} \int
K_\varepsilon(\alpha,z)\,v(z)\,dz}\enskip \textrm{with}\\
 K_\varepsilon(\alpha,z)=\Big(\frac\lambda{2\pi
}\Big)^^n
\iint  e^{-\frac\lambda 2 \v\alpha_\xi\v^2}\,\boldsymbol{1}_{\tilde
W^+}(\alpha)\,e^{i\lambda[\varphi (\theta,y,\alpha)+(y-z)\cdot
\xi]}\,k(\theta,y,\alpha,\lambda)\\
\quad \quad \quad
\quad \quad \quad
\quad \chi_5\Big(\frac{y-x(\theta,\alpha)}{\bra\theta\ket}\Big)\,\partial
^\gamma_y\,\chi^+_4(y)(1-\psi_4(\xi))\,g(\varepsilon\xi)\,d\xi\,dy\,.
\end{array}\right.
\end{equation}
Let us consider the vector field
$$
X=\frac 1{1+\lambda\v y-z\v^2} \left (1+\frac 1 i \sum^n_{j=1} 
(y_j-z_j)\,\frac\partial {\partial \xi_j}\right)\,.
$$
Then it is easy to see that
\begin{equation*}\left\{
\begin{array}{l}
X\,e^{i\lambda(y-z)\cdot \xi}=e^{i\lambda(y-z)\cdot \xi}\\
(^tX)^N=\som_{\v A\v\leq N} \frac{C_A(y-z)^A\,\partial
^\alpha_\xi}{(1+\lambda \v y-z\v^2)^N}\,.
\end{array}\right.
\end{equation*}
Then we can write
\begin{multline}
\int e^{i(y-z)\cdot \xi}\,(1-\psi_4(\xi))\,g(\varepsilon\xi)\,d\xi=\int
e^{i\lambda(y-z)\cdot \xi}\,(^tX)^N
[(1-\psi_4(\xi))\,g(\varepsilon \,\xi)]\,d\xi\\
=\sum_{A=A_1+A_2\atop \v A\v\leq N} \int e^{i\lambda(y-z)\cdot
\xi}\,C_{A_1A_2}\,\frac{(y-z)^A}{(1+\lambda\v y-z\v^2)^N}\,\partial
^{A_1}_\xi (1-\psi_4(\xi))\,\varepsilon^{\v A_2\v}\,(\partial
^{A_2}_\xi\,g)(\varepsilon\,\xi)\,d\xi\,.
\end{multline}
Now on the support of
$\psi_3(\alpha_\xi)(1-\psi_4(\xi))$ we have
$\v\xi+\alpha_\xi\v\geq \mu>0$, it follows from \kern 1pt
(\ref{eqVI.3.29}) that
$$
\Big\v\frac{\partial \varphi }{\partial y}\,(\theta,y,\alpha)+\xi\Big\v\geq
\v\alpha_\xi+\xi\v-\Big\v\frac{\partial \varphi }{\partial
y}\,(\theta,y,\alpha)-\alpha_\xi\Big\v\geq
\mu-C(\varepsilon+\sqrt\delta)\geq \frac \mu 2
$$
if $\varepsilon$ and $\delta$ are small enough.

If $\v\xi\v\geq 2\,\sup \big\v\frac{\partial\varphi }{\partial
y}\,(\theta,y,\alpha)\big\v$ we have $\big\v\frac{\partial \varphi }{\partial
y}\,(\theta,y,\alpha)+\xi\big\v\geq \frac{\v\xi\v}2$. Therefore in all
cases we have, with $\eta_0>0$,
\begin{equation}\label{eqVI.3.32}
\Big\v\frac{\partial \varphi }{\partial y}\,(\theta,y,\alpha)+\xi\Big\v\geq
\eta_0\,\bra\xi\ket\,.
\end{equation}
Let us set then
$$
Y=\frac1{i\lambda\v\xi+\frac{\partial \varphi }{\partial
y}\,(\theta,y,\alpha)\v^2}
\sum^n_{j=1} \Big(\frac{\partial \varphi }{\partial
y_j}\,(\theta,y,\alpha)+\xi_j\Big)\,\frac\partial {\partial y_j}\,.
$$
and $T=\frac{\partial \varphi }{\partial y}\,(\theta,y,\alpha)+\xi$.

Then
\begin{equation}\label{eqVI.3.33}\left\{
\begin{array}{l}
Y\,e^{i\lambda(\varphi (\theta,y,\alpha)+(y-z)\cdot
\xi)}=e^{i\lambda(\varphi (\theta,y,\alpha)+(y-z)\cdot \xi)}\\
(^tY)^N=\frac 1{(i\,\lambda)^N}\bigg\{\som_{\v\nu\v=N}\Big (\frac T{\v
T\v^2}\Big)^\nu\,\partial ^\nu_y+\som_{\v\nu\v\leq
N-1}\frac{P_{3N-2\v\nu\v-1}(\theta,y,\alpha,T,\overline T)}{\v
T\v^{4N-2\v\nu\v}}\,\partial ^\nu_y\bigg\}
\end{array}\right.
\end{equation}
where $P_k(\theta,y,\alpha,T,\overline T)$ is a polynomial in $T,\overline T$
of order $\leq k$ with $C^\infty-$bounded coefficients.

It follows from (\ref{eqVI.3.32}) that on the support of
$\psi_3(\alpha_\xi)(1-\psi_4(\xi))$ we have,
\begin{equation}\label{eqVI.3.34}
\som_{\v\nu\v=N}\bigg\v\Big(\frac T{\v
T\v^2}\Big)^\nu\bigg\v+\som_{\v\nu\v\leq N-1} \frac{\v
P_{3N-2\v\nu\v-1}(\cdots )\v}{\v T\v^{4N-2\v\nu\v}}\leq
\frac{C_N}{\bra\xi\ket^N}\,.
\end{equation}
On the other hand we check by induction that
\begin{equation}\label{eqVI.3.35}
\partial ^\nu_y \Big[\frac{(y-z)^A}{(1+\lambda\v y-z\v^2)^N}\Big]=
\som_{{j\leq \v\nu\v\atop\v\beta\v\leq \v A\v+j}\atop 2j+\v A\v\leq
\v\beta\v+\v\nu\v} b_{N,j,A,\beta}\, \frac{(y-z)^\beta\,\lambda^j}{(1+\lambda
\v y-z\v^2)^{N+j}}\,.
\end{equation}
Now if we insert (VI.3.31) into (\ref{eqVI.3.30}) and if we make
integration by parts with respect to $Y$ using (\ref{eqVI.3.33}) we see using
(\ref{eqVI.3.34}) and (\ref{eqVI.3.3}) that $K_\varepsilon(\alpha,z)$ is
bounded by a finite sum of integrals of the following type
\begin{equation*}
\begin{split}
\int \frac{\lambda^n \v
y-z\v^{\v\beta\v}\,\lambda^j}{\lambda^N\bra\xi\ket^{N}(1+\lambda\v
y-z\v^2)^{N+j}}\,\boldsymbol{1}_{\tilde W^+}(\alpha) \v\partial
^{A_1}_\xi(1-\psi_4(\xi))\v\,\varepsilon^{\v A_2\v}\,(\partial
^{A_2}_\xi g)(\varepsilon\xi)\,\v\partial ^{\nu_1}_yk\v\\
\Big\v(\partial
^{\nu_2}_y\,\chi_5)\Big(\frac{y-x(\theta,\alpha)}{\bra\theta\ket}\Big)\Big\v
\bra\theta\ket^{-\v\nu_2\v}\v \partial ^{\gamma+\nu_3}_y\,\chi^+_4(y)\v\, 
e^{-\frac\lambda{16}\frac{\v
y-x(\theta,\alpha)\v^2}{\bra\theta\ket^2}}\,dy\,d\xi
\end{split}
\end{equation*}
where $\v\beta\v\leq \v A\v+j$, $\nu_1+\nu_2+\nu_3=\nu$, $\v\nu\v\leq N$,
$A=A_1+A_2$, $\v A\v\leq N$, $j\leq \v\nu\v$, $2j+\v A\v\leq \v\beta\v+\v\nu\v$.

{\bf Claim } We have
\begin{equation*}\left\{
\begin{array}{l}
(1)={\displaystyle \sup_\alpha }\int \v K_\varepsilon(\alpha,z)\v\,dz\leq C_N\,
\frac{\bra\theta\ket^n}{\lambda^{N/2}}\\
(2)={\displaystyle \sup_z} \int \v K_\varepsilon(\alpha,z)\v\,d\alpha\leq C'_N\,
\frac{\bra\theta\ket^n}{\lambda^{N/2}}\,.
\end{array}\right.
\end{equation*}
Let us first   remark that
\begin{equation}\label{eqVI.3.36}
\int \frac{\lambda^j\,\v x\v^{\v\beta\v}}{(1+\lambda \v
x\v^2)^{N+j}}\,dx=C\, \lambda^{j-\frac{\v\beta\v}2-\frac n2}\leq
C\,\lambda^{\frac N2-\frac n2}\,.
\end{equation}
Indeed $-\v\beta\v\leq \v\nu\v-2j-\v A\v$ so $j-\frac{\v\beta\v}2\leq
\frac{\v\nu\v}2-\frac{\v A\v}2\leq \frac N2$.

To estimate $(1)$ we use the above estimate of $K_\varepsilon(\alpha,z)$ which
we integrate with respect to $z$. For the integral in $z$ we use
(\ref{eqVI.3.36}). The integral in $\xi$ is estimated thanks to the term $\frac
1{\bra\xi\ket^{N}}$ where $\v\nu\v\leq N$, finally the integral in $y$
is bounded by $\int e^{-\frac\lambda{16} \frac{\v
y-x(\theta,\alpha)\v^2}{\bra\theta\ket^2}}\,dy\leq
C\,\frac{\bra\theta\ket^n}{\lambda^{\frac n2}}$. Therefore we obtain
$$
(1)\leq C\, \frac{\lambda^n}{\lambda^N}\, \lambda^{\frac{N}{2}-\frac {n}{2}}\,
\frac{\bra\theta\ket^n}{\lambda^{\frac
n2}}=C\,\frac{\bra\theta\ket^n}{\lambda^{\frac N2}}\,.
$$
To estimate the term $(2)$ we use the change of variables $\tilde
\alpha=(x(\theta,\alpha),\xi(\theta,\alpha))$ as in the proof of Lemma
\ref{lVI.3.2} and (\ref{eqVI.3.36}). This  gives, 
$(2)\leq C\,\frac{\bra\theta\ket^n}{\lambda^{\frac N2}}$. Since
$\v\theta\v\leq
\lambda\,T$ we obtain
$$
(1)+(2)\leq C\,\lambda^{-\frac N 2+n}\,.
$$
We can therefore use the Schur Lemma and (\ref{eqVI.3.30}) to achieve the proof
of Lemma \ref{lVI.3.3}. \cqfd
\begin{corollary}\sl \label{cVI.3.4}
With the notations of (\ref{eqVI.3.15}) and (\ref{eqVI.3.11}) we have
$$
B_3=\sum_{1\leq \v\gamma\v\leq 2} \int e^{i\lambda\varphi
(\theta,y,\alpha)}\,d_\gamma(\theta,y,\alpha,\lambda)(\partial
^\gamma_y\,\chi^+_4)(y)\,\chi_5 \Big
(\frac{y-x(\theta,\alpha)}{\bra\theta\ket}\Big )\,\psi_4\Big (\frac
D\lambda\Big )\,U(t,y,\lambda)\,dy+J^+_\lambda(t)\,v(\alpha)
$$
\begin{equation}\label{eqVI.3.37}
\Vert J^+_\lambda(t)\,v\Vert _{L^2(W^+)}\leq C_N\,\lambda^{-N}\,\Vert v\Vert
_{L^2(\R^n)}\,.
\end{equation}
\end{corollary}
Now we state the following result.
\begin{lemma}\sl \label{lVI.3.5}
Let $b=b(\theta,y,\alpha,\lambda)$ be a bounded symbol. Let us set
$$
G(\theta,\alpha)=\int e^{i\lambda\varphi (\theta,y,\alpha)}\,\chi_5 \Big
(\frac{y-x(\theta,\alpha)}{\bra\theta\ket}\Big)\,b(\theta,y,\alpha,\lambda)\,v(y)\,dy\,.
$$
Then one can find $C>0$ such that for all $\v\theta\v\leq \lambda\,T$ and $v\in
L^2(\R^n)$,
$$
\big\Vert e^{-\frac\lambda 2 \v\alpha_\xi\v^2}\,G(\theta,\cdot
)\big\Vert _{L^2(W^+)}\leq C\, \frac{\bra\theta\ket^{n}}{\lambda^{n/2}}\,
\Vert v\Vert _{L^2(\R^n)}\,.
$$
\end{lemma}

{\bf Proof } Let us write
$$
G(\theta,\alpha)=\int K(\theta,\alpha,y,\lambda)\,v(y)\,dy\,.
$$
Then using the estimate (\ref{eqVI.3.3}) we see that
$$
e^{-\frac \lambda2 \v\alpha_\xi\v^2}\,\boldsymbol{1}_{W^+}(\alpha)\, \v
K(\theta,\alpha,y,\lambda)\v\leq C\, e^{-\frac\lambda{16} \frac{\v
y-x(\theta,\alpha)\v^2}{\bra\theta\ket^2}}\,\boldsymbol{1}_{\frac 12 \leq
\v\alpha_\xi\v\leq 2}\,.
$$
From this estimate we can use the Schur Lemma (making the change of variables
$\tilde \alpha=(x(\theta,\alpha),\xi(\theta,\alpha))$) to conclude. \cqfd

Let now $\psi_5\in C^\infty _0(\R^n)$ be such that $0\leq \psi_5\leq 1$ and
\begin{equation}\label{eqVI.3.38}\left\{
\begin{array}{ll}
\psi_5(\xi)=1\enskip \textrm { if }\enskip
\Big\v\frac\xi{\v\xi\v}+\xi_0\Big\v\leq
7\delta_2\,,&  a-7\delta_2\leq \v\xi\v\leq b+7\delta_2\\
\supp \psi_5 \subset
\Big\{\Big\v\frac\xi{\v\xi\v}+\xi_0\Big\v\leq
8\delta_2\Big\}\,,& a-8\delta_2\leq \v\xi\v\leq b+8\delta_2
\end{array}\right.
\end{equation}
The analogue of Lemma \ref{lVI.2.3} proves that one can find $C>0$,
$\varepsilon_0>0$ such that
\begin{equation}\label{eqVI.3.39}
\Big\Vert \psi_4\Big(\frac D\lambda\Big)\,T^*_{\beta\rightarrow
y}[(1-\psi_5(\beta_\xi))v]\Big\Vert _{L^2(\R^n_y)}\leq
C\,e^{-\varepsilon_0\lambda}\, \Vert e^{-\frac\lambda 2
\v\beta_\xi\v^2}\,v\Vert _{L^2(\R^n_\beta)}\,.
\end{equation}
\begin{corollary}\sl \label{cVI.3.6}
We have
\begin{equation*}
\begin{split}
B_3=\sum_{1\leq \v\gamma\v\leq 2} \int e^{i\lambda\varphi
(\theta,y,\alpha)}\,d_\gamma(\theta,y,\alpha,\lambda)(\partial
^\gamma_y\,\chi^+_4)(y)\,\chi_5\Big
(\frac{y-x(\theta,\alpha)}{\bra\theta\ket}\Big )\, &T^*_{\beta\rightarrow y}
[\psi_5(\beta_\xi)\,T_{z\rightarrow
\beta}\,U(t,z,\lambda)]\,dy\\
&+J^+_\lambda(t)\,v(\alpha)
\end{split}
\end{equation*}
where $J^+_\lambda(t)$ satisfies (\ref{eqVI.3.41}).
\end{corollary}

{\bf Proof } We use Corollary \ref{cVI.3.4}, (\ref{eqVI.3.39}) and Lemma
\ref{lVI.3.3} to remove $\psi_4\big (\frac D\lambda\big )$. \cqfd

Let now $\chi^+_6=\chi^+_6(\beta_x)\in C^\infty _0(\R^n)$ be such that $0\leq
\chi^+_6\leq 1$ and
\begin{equation}\label{eqVI.3.40}\left\{
\begin{array}{l}
\chi^+_6(\beta_x)=1\enskip \textrm { if }\enskip \frac 7 2 \,\delta_1\leq
-\beta_x\cdot \xi_0\leq 6\,\delta_1\\
\supp \chi^+_6\subset \big\{\beta_x : \frac{17}5\,\delta_1\leq -\beta_x\cdot
\xi_0\leq 7\,\delta_1\big\}\,.
\end{array}\right.
\end{equation}
Then we can apply Lemma \ref{lVI.2.6} with $\chi_a=\partial
^\gamma_y\,\chi^+_4$, $\v\gamma\v\geq 1$, $\chi_b=\chi^+_6$. Indeed on the
support of $\partial ^\gamma_y\,\chi^+_4(y)(1-\chi^+_6(\beta_x))$ we have
$\frac{15}4\,\delta_1\leq -y\cdot \xi_0\leq 5\,\delta_1$ and
$-\beta_x\cdot \xi_0\leq \frac 72\,\delta_1$ or $-\beta_x\cdot
\xi_0\geq 6\,\delta_1$. In the first case we write
$$
\v y-\beta_x\v\geq \beta_x\cdot \xi_0-y\cdot \xi_0\geq
\frac{15}4\,\delta_1-\frac {7 }{2}\,\delta_1=\frac 14\,\delta_1\,,
$$
 and in the second case we have,
$$
\v y-\beta_x\v\geq y\cdot\xi_0-\beta_x\cdot \xi_0\geq
6\,\delta_1-5\,\delta_1=\delta_1\,.
$$
Therefore we obtain
\begin{equation}\label{eqVI.3.41}
\Vert \partial ^\gamma_y\,\chi^+_4(y)\,T^*_{\beta\rightarrow
y}\,\psi_5(\beta_\xi)(1-\chi_6(\beta_x))W\Vert _{L^2}\leq
C\,e^{-\varepsilon\lambda}\,\Vert e^{-\frac\lambda 2
\v\beta_\xi\v^2}\,W\Vert _{L^2}\,.
\end{equation}
Using Corollary \ref{cVI.3.6} we deduce the following Lemma.
\begin{lemma}\sl \label{lVI.3.7}
We have
\begin{equation*}
\begin{split}
B_3=\sum_{1\leq \v\gamma\v\leq 2} \int e^{i\lambda\varphi
(\theta,y,\alpha)}\,d_\gamma(\theta,y,\alpha,\lambda)(\partial
^\gamma_y\,\chi^+_4)(y)\,\chi_5\Big(\frac{y-x(\theta,\alpha)}{\bra\theta\ket}\Big
)\,T^*_{\beta\rightarrow y}
[\psi_5(\beta_\xi)\,\chi^+_6(\beta_x)\\T_{z\rightarrow
\beta}\,U(t,\cdot ,\lambda)](y)\,dy
+J^+_\lambda(t)\,v(\alpha)\,,
\end{split}
\end{equation*}
where $J^+_\lambda(t)$ satisfies (\ref{eqVI.3.37}).
\end{lemma}
Now on the support of $\psi_5(\beta_\xi)\,\chi^+_6(\beta_x)$ we have by
(\ref{eqVI.3.38}), (\ref{eqVI.3.40}),
\begin{equation*}
\begin{aligned}
\v\beta_x\cdot \beta_\xi\v&\leq \Big(\v\beta_x\cdot
\xi_0\v+\v\beta_x\v\Big\v\frac{\beta_\xi}{\v\beta_\xi\v}+
\xi_0\Big\v\Big)\,\v\beta_\xi\v\\
&\leq (16\,\delta_1+8\,\delta_2\,\v\beta_x\v)\,\v\beta_\xi\v\leq
c_0\,\bra\beta_x\ket \,\v\beta_\xi\v
\end{aligned}
\end{equation*}
if $16\,\delta_1+8\,\delta_2\leq c_0$.  Therefore we are in the case (i) of
(\ref{eqVI.3.17}) and since Theorem \ref{tVI.3.1} is already proved in this
case for $t\in[-T,T]$ we can write,
\begin{multline}
T_{z\rightarrow \beta}\,U(t,\cdot ,\lambda)=\lambda^{\frac{3n}4}\int
e^{i\lambda\varphi (-\lambda t,z,\beta)}\,a(-\lambda
t,z,\beta,\lambda)\,\chi^+_4(z)\,\chi_5\Big(\frac{z-x(-\lambda
t,z)}{\bra\lambda t\ket}\Big)\\
(\chi^+_2v)(z)\,dz+J^+_\lambda(t)\,v(\beta)\,,
\end{multline}
where $J^+_\lambda(t)$ satisfies
\begin{equation}\label{eqVI.3.43}\left\{
\begin{array}{l}
\textrm{for every}\enskip N\in\N\enskip \textrm{ one can find}\enskip
C_N>0\enskip \textrm{such that}\\
\big\Vert e^{-\frac\lambda 2
\v\beta_\xi\v^2}\,\psi_5(\beta_\xi)\,\chi^+_6(\beta_x)\,J^+_\lambda(t)
v\big\Vert _{L^2}\leq C_N\,\lambda^{-N}\,\Vert v\Vert _{L^2(\R^n)}\,.
\end{array}\right.
\end{equation}\cqfd

From this we can deduce the following result.
\begin{corollary}\sl \label{cVI.3.8}
We have 
$$
T_{z\rightarrow \beta}\,U(t,\cdot ,\lambda)=\lambda^{\frac{3n}4} \int
e^{i\varphi (-\lambda t,z,\beta)}\,a(-\lambda
t,z,\beta,\lambda)\,\chi_5\Big(\frac{z-x(-\lambda t,\beta)}{\bra\lambda
t\ket}\Big )\,\chi^+_3(z)\,(\chi^+_2v)(z)\,dz+J^+_\lambda(t)\,v(\beta)
$$
where $J^+_\lambda(t)$ satisfies (\ref{eqVI.3.43}).
\end{corollary}

{\bf Proof } We have just to show that we can replace $\chi^+_4$
 by $\chi^+_3$ in (VI.3.42). But this is obvious since (see
(\ref{eqVI.2.1}), (\ref{eqVI.3.1})) we have
$\chi^+_4\,\chi^+_2=\chi^+_4(1-\chi^+_3)\,\chi^+_2+\chi^+_4\,\chi^+_3\,\chi^+_2$
and
$(1-\chi^+_3)\,\chi^+_2\equiv 0$, $\chi^+_4\,\chi^+_3=\chi^+_3$. \cqfd

We are ready now to prove (\ref{eqVI.3.26}).
\begin{lemma}\sl \label{lVI.3.9}
For all $N\in\N$ one can find $C_N>0$ such that
$$
\big\Vert e^{-\frac\lambda 2 \v\alpha_\xi\v^2}\,B_3(\theta,t,\cdot
,\lambda)\big\Vert _{L^2(W^+)}\leq C_N\,\lambda^{-N}\, \Vert v\Vert _{L^2(\R^n)}
$$
for all $\lambda\geq 1$, $\theta=\lambda (s-t)\in[0,\lambda T]$ and all $v\in
L^2(\R^n)$.
\end{lemma}

{\bf Proof } We use first Corollary \ref{cVI.3.8} and Lemma \ref{lVI.3.7}. On
the support of
$\psi_5(\beta_\xi)\,\chi^+_6(\beta_x)\,\chi_5\big (\frac{z-x(-\lambda
t,\beta)}{\bra\lambda t\ket}\big)$ we have by (\ref{eqVI.3.38}),
(\ref{eqVI.3.40}), (\ref{eqVI.3.2}), since $x(-\lambda
t,\beta)=\beta_x-2\lambda t\,\beta_\xi+\CO(\varepsilon\bra t\ket )$, 
\begin{equation*}
\begin{aligned}
z\cdot \xi_0&\leq (z-x(-\lambda t,\beta))\cdot
\xi_0+(\beta_x-2\lambda t\,\beta_\xi)\cdot
\xi_0+C\,\varepsilon\,\bra\lambda t\ket \\
&\leq \beta_x\cdot \xi_0-2\lambda t\,\beta_\xi\cdot
\Big(\xi_0+\frac{\beta_\xi}{\v\beta_\xi\v}\Big)+2\lambda
t\,\v\beta_\xi\v+C\,(\varepsilon+\delta)\bra\lambda t\ket\\
&\leq \frac 7 2\,\delta_1+2\lambda t\,\v\beta_\xi\v+C
(\varepsilon+\delta+\delta_2)\bra\lambda t\ket\,.
\end{aligned}
\end{equation*}
Since $\v\beta_\xi\v\geq a-\delta_2$ we obtain
$$
z\cdot \xi_0\leq -\frac{17}5\,\delta_1-2(a-\delta_2)\,\lambda\,\v t\v+C
(\varepsilon+\delta+\delta_2)\bra\lambda t\ket\,.
$$
Taking $\varepsilon,\delta,\delta_2$ small with respect to $\delta_1$ and $a$
we obtain $z\cdot \xi_0\leq -\frac{10}3\,\delta_1$. Now on the support of
$\chi^+_3(z)$ we have by (\ref{eqVI.2.1}), $z\cdot \xi_0\geq
-3\,\delta_1$.

It follows from Corollary \ref{cVI.3.8} that $T_{z\rightarrow \beta}\,U(t,\cdot
,\lambda)=R^+v$ where $R^+$ satisfies (\ref{eqVI.3.37}). Now we use Lemma
\ref{lVI.3.5} and we obtain since $\v\theta\v\leq \lambda T$,
\begin{equation*}
\begin{split}
\big\Vert e^{-\frac\lambda 2\v \alpha_\xi\v^2}\,B_3(\theta,\cdot
,\lambda)\big\Vert _{L^2(W^+)}&\leq C\,\frac{\bra\theta\ket^{\frac n
2}}{\lambda^{\frac n 2}}\,\big\Vert T_{\beta\rightarrow y}
[\psi_5(\beta_\xi)\,\chi^+_6(\beta_x)\,T_{z\rightarrow
\beta}\,U(t,\cdot ,\lambda)]\big\Vert _{L^2(\R^n)}\\
&\quad \leq \big\Vert e^{-\frac\lambda 2\v
\beta_\xi\v^2}\,\psi_5(\beta_\xi)\,\chi^+_6(\beta_x)\,T_{z\rightarrow
\beta}\,U(t,\cdot )\big\Vert _{L^2}\,.
\end{split}
\end{equation*}
Since $T_{z\rightarrow \beta}\,U(t,\cdot )=R^+v$, where $R^+$ satisfies
(\ref{eqVI.3.37}), we obtain the conclusion of Lemma \ref{lVI.3.9}. 

\cqfd

To complete the proof of Theorem \ref{tVI.3.1} in the case (\ref{eqVI.3.25}) we
use (\ref{eqVI.3.11}), (\ref{eqVI.3.14}), Lemmas \ref{lVI.3.2},
\ref{lVI.3.9} and the same argument as in the end of the proof of the case
(\ref{eqVI.3.17}) to remove the cut-off functions $\chi^+_4(y)$ and
$\chi_5(y-\alpha_x)$. \cqfd

\subsection{Conclusion of Section VI}\label{ssVI.4}

Here we state a result which combines the conclusions of Theorems
\ref{tVI.2.2}, \ref{tVI.3.1} and (\ref{eqVI.3.5}).
\begin{theorem}\sl\label{tVI.4.1}
Let $K_\pm$ the operators defined in (\ref{eqVI.2.8}). Then we can write
$$
K_+(t)\,u(x)=I+II+III
$$
where
\begin{equation*}
\begin{aligned}
I&=\lambda^{\frac{3n}2}\,\chi^+_1(x)\,\psi_2\Big(\frac
{D_x}\lambda\Big)\,\chi^+_2(x)\Bigg[\iint e^{i\lambda F(-\lambda
t,x,y,\alpha)}\,a(-\lambda
t,y,\alpha)\,\chi^+_2(y)\,\chi^+_3(\alpha_x)\,\psi_3(\alpha_\xi)\\
&\hbox to 6cm{}\chi_5\Big(\frac{y-x(-\lambda t,\alpha)}{\bra\lambda t\ket}\Big)
\Big (\psi_2\Big(\frac D\lambda\Big)\,\chi^+_1u\Big)(y)\,dy\,d\alpha\Bigg]\\
II&=\chi^+_1(x)\,\psi_2\Big(\frac
{D_x}\lambda\Big)\,\chi^+_2(x)\,T^*_{\alpha\rightarrow x}
\bigg[\chi^+_3(\alpha_x)(\alpha_x)\,\psi_3(\alpha_\xi)\,J^+_\lambda(t)\Big(\chi^+_2\,\psi_2\Big(\frac
{D_x}\lambda\Big)\,\chi^+_1u\Big)\bigg]\\
III&=R^+_\lambda(t)\,u
\end{aligned}
\end{equation*}
where
\begin{equation}\label{eqVI.4.1}\left\{
\begin{array}{l}
F(-\lambda t,x,y,\alpha)=\varphi (-\lambda t,y,\alpha)-(x-\alpha_x)\cdot
\alpha_\xi+\frac i2 \,\v x-\alpha_x\v^2+\frac i2
\,\v\alpha_\xi\v^2\,,\\
\v a(-\lambda t,y,\alpha,\lambda)\v\leq C\,\bra\lambda t\ket^{-\frac n2}\,,\\
\Vert \chi^+_3(\alpha_x)\,\psi_3(\alpha_\xi)\,e^{-\frac\lambda 2
\v\alpha_\xi\v^2}\,J^+_\lambda(t)\,v\Vert _{L^2}\leq
C\,\lambda^{-N}\,\Vert v\Vert _{L^2}\,,\enskip \forall N\in\N\,,\\
\Vert R^+_\lambda(t)\,u\Vert _{H^{2N}}\leq C_N\, \Vert u\Vert
_{H^{-2N}}\,,\enskip \forall N\in\N\,,
\end{array}\right.
\end{equation}
and $\chi^+_i,\psi_j$ have been defined in (\ref{eqVI.2.1}) to (\ref{eqVI.2.6}),
(\ref{eqVI.2.9}), (\ref{eqVI.3.1}) and (\ref{eqVI.3.2}).
Moreover the same result holds with the minus sign.
\end{theorem}

\section{\ The dispersion estimate
and the end of the proof of
Theorem~I.1}\label{sVII}

\subsection{\ The dispersion estimate for the operators $K_\pm(t)$}\label{ssVII.1}

Let us recall that $K_\pm(t)$ have been introduced in (\ref{eqVI.2.8}). The
purpose of this paragraph is to prove the following result.
\begin{theorem}\sl \label{tVII.1.1}
Let $T>0$. Then there exists a constant $C\geq 0$ such that
$$
\Vert K_\pm (t)\,u\Vert _{L^\infty }\leq \frac{C}{\v t\v^{\frac n 2}}\, \Vert
u\Vert _{L^1}
$$
for all $0<\v t\v\leq T$ and all $u\in L^1(\R^n)$.
\end{theorem}

{\bf Proof } We shall use Theorem \ref{tVI.4.1} and its notation and we shall
consider only $K_+((t)$. Then we can write
\begin{equation}\label{eqVII.1.1}
\Vert K_+(t)\,u\Vert _{L^\infty }\leq \Vert I\Vert _{L^\infty }+\Vert II\Vert
_{L^\infty }+\Vert III\Vert _{L^\infty }\,.
\end{equation}
Let $N_0\in\N$ be such that $2N_0>\frac n2$. By the Sobolev embedding and
(\ref{eqVI.4.1}) we have
\begin{equation}\label{eqVII.1.2}
\Vert III\Vert _{L^\infty }\leq C\,\Vert R^+_\lambda(t)u\Vert _{H^{2N_0}}\leq
C'_{N_0}\,\Vert u\Vert _{H^{-2N_0}}\leq C''\,\Vert u\Vert _{L^1}\leq
\frac{C(T)}{\v t\v^{\frac n2}}\,\Vert u\Vert _{L^1}\,. 
\end{equation}
Let us consider the term $II$. We have
\begin{equation*}
\begin{aligned}
\Vert II\Vert _{L^\infty }&\leq C\,\Big\Vert \chi^+_1\,\psi_2\Big (\frac
D\lambda\Big)\,\chi^+_2\,T^*\bigg (\chi^+_3\,\psi_3\,J^+_\lambda\Big
(\chi^+_2\,\psi_2\Big (\frac D\lambda\Big )\,\chi^+_1u\Big)\bigg)\Big\Vert
_{H^{2N_0}}\\
&\leq C'\,\lambda^{2N_0}\,\Vert T^*(\cdots )\Vert _{L^2(\R^n)}\\
&\leq C''\,\lambda^{2N_0}\,\Big\Vert e^{-\frac \lambda 2
\v\alpha_\xi\v^2}\,\chi^+_3(\alpha_x)\,\psi_3(\alpha_\xi)\,J^+_\lambda(t)\bigg
(\chi^+_2\,\psi_2\Big(\frac D\lambda\Big)\,\chi^+_1u\bigg)\Big\Vert
_{L^2(\R^{2n}_\alpha)}\\
&\leq C_N\,\lambda^{2N_0-N}\,\Big\Vert \chi^+_2\,\psi_2\Big(\frac
D\lambda\Big)\,\chi^+_1u\Big\Vert _{L^2(\R^n)}\\
&\leq C'_N\,\lambda^{2N_0-N}\,\Big\Vert \psi_2\Big(\frac
D\lambda\Big)(I-\Delta)^{-N_0}\,\chi^+_1u\Big\Vert _{L^2(\R^n)}\\
&\leq C''_N\,\lambda^{4N_0-N}\,\Vert u\Vert _{H^{-2N_0}}\,.
\end{aligned}
\end{equation*}
Taking $N=4N_0$ we obtain finally,
\begin{equation}\label{eqVII.1.3}
\Vert II\Vert _{L^\infty }\leq C\,\Vert u\Vert _{L^1}\leq \frac{C(T)}{\v
t\v^{\frac n 2}}\,\Vert u\Vert _{L^1}\,.
\end{equation}
So we are left with the estimation of $\Vert I\Vert _{L^\infty }$. Let us set
\begin{equation}\label{eqVII.1.4}
\begin{split}
k_+(t,x,y,\lambda)=\lambda^{\frac{3n}2} \int e^{i\lambda F(-\lambda
t,x,y,\alpha)}\,&a(-\lambda
t,y,\alpha,\lambda)\,\chi^+_2(y)\,\chi^+_3(\alpha_x)\,\psi_3(\alpha_\xi)\\
&\quad \quad \cdot \chi_5\Big (\frac{y-x(-\lambda t,\alpha)}{\bra\lambda
t\ket}\Big )\,d\alpha
\end{split}
\end{equation}
and
\begin{equation}\label{eqVII.1.5}
\tilde K_+(t)\,v(x)=\int k_+(t,x,y,\lambda)\Big[\psi_2\Big(\frac
D\lambda\Big)\,\chi^+_1v\Big](y)\,dy\,.
\end{equation}
Then
\begin{equation}\label{eqVII.1.6}
I=\chi^+_1\,\psi_2\Big(\frac D\lambda\Big)\,\chi^+_2\,\tilde K_+(t)\,v\,.
\end{equation}
Since the operator $\chi^+_1\,\psi_2\big(\frac D\lambda\big)\,\chi^+_2$ is
bounded from $L^\infty $ to $L^\infty $ with bound independent of $\lambda$ we
have,
\begin{equation}\label{eqVII.1.7}
\Vert I\Vert _{L^\infty }\leq C\,\Vert \tilde K_+(t)\,v\Vert _{L^\infty }\,.
\end{equation}
Assume that the kernel $k_+$ has the following bound,
\begin{equation}\label{eqVII.1.8}
\vert k_+(t,x,y,\lambda)\vert \leq \frac C{\v t\v^{\frac n 2}}\,,
\end{equation}
with $C$ independent of $\lambda$. It will follow from (\ref{eqVII.1.7}),
(\ref{eqVII.1.5}) and (\ref{eqVII.1.8})
 that
$$
\Vert I\Vert _{L^\infty }\leq \frac C{\v t\v^{\frac n2}}\Big\Vert \psi_2\Big
(\frac D\lambda\Big)\,\chi^+_1v\Big\Vert _{L^1}\,.
$$
Since the operator $\psi_2\big (\frac D\lambda\big)\,\chi^+_1$ is uniformly
bounded on $L^1$ we will have
\begin{equation}\label{eqVII.1.9}
\Vert I\Vert _{L^\infty }\leq \frac C{\v t\v^{\frac n2}}\,\Vert v\Vert _{L^1}\,.
\end{equation}
Then Theorem \ref{tVII.1.1} follows from (\ref{eqVII.1.1}), (\ref{eqVII.1.2}),
(\ref{eqVII.1.3}) and (\ref{eqVII.1.9}).

{\bf Proof of (\ref{eqVII.1.8}) }  We divide the proof in three cases : $\lambda
t\geq 1$, $\lambda t\leq -1$, $\v\lambda t\v\leq 1$. 

Let us remark first that in the integral in the right hand side of
(\ref{eqVI.1.4}), on the support of $\chi^+_3(\alpha_x)\cdot
\psi_3(\alpha_\xi)$ we have $-\alpha_x\cdot \xi_0\leq 3\,\delta_1$ and
$\big\vert \frac{\alpha_\xi}{\v\alpha_\xi\v}+\xi_0\big\vert \leq
4\,\delta_2$. Therefore
$$
\alpha_x\cdot \frac{\alpha_\xi}{\v\alpha_\xi\v}=\alpha_x\cdot
\Big(\frac{\alpha_\xi}{\v\alpha_\xi\v}+\xi_0\Big)-\alpha_x\cdot \xi_0\leq
(4\,\delta_2+3\,\delta_1)\bra\alpha_x\ket
$$
so
\begin{equation}\label{eqVII.1.10}
\alpha_x\cdot \alpha_\xi\leq c_0\,\bra\alpha_x\ket\,\v\alpha_\xi\v\,,
\end{equation}
if $\delta_1$ and $\delta_2$ are small compared to $c_0$.

{\bf Case A : } proof of (\ref{eqVII.1.8}) when $\lambda\, t\geq 1$.

In this case $\theta=-\lambda t<0$ and it follows from (\ref{eqVII.1.10}) and
Definition \ref{dIII.2.2} that all the points $\alpha$ in the integral giving
$k_+$ are outgoing for $\theta<0$. It follows from Corollary \ref{cIII.3.3}
that
\begin{equation}\label{eqVII.1.11}
\frac{\partial x_j}{\partial
\alpha^k_\xi}\,(\theta,\alpha)=2\theta\,\delta_{jk}+\CO(\varepsilon\,\bra\theta\ket)\,,
\quad 1\leq j,k\leq n\,.
\end{equation}
Now using Theorem \ref{tVI.4.1} and (\ref{eqVII.1.4}) we obtain
\begin{equation}\label{eqVII.1.12}
\v k^+(t,x,y,\lambda)\v\leq c\,\lambda^{\frac{3n}2} \int
e^{-\frac\lambda{16}\frac{\v y-x(-\lambda t,\alpha)\v^2}{\bra\lambda
t\ket^2}-\frac\lambda 2 \v x-\alpha_x\v^2}\, \psi_3(\alpha_\xi)\,\bra\lambda
t\ket^{-\frac n 2}\,d\alpha\,.
\end{equation}
By (\ref{eqVII.1.11}) we can make the change of variables
$$
\tilde \alpha_x=\alpha_x\,,\quad \tilde \alpha_\xi=x(-\lambda t,\alpha)
$$
and $\big\v\det\,\frac{\partial \tilde \alpha}{\partial \alpha}\big\vert\geq
(\lambda t)^n$ if $\varepsilon$ is small enough (since $\bra\lambda t\ket\leq
\sqrt 2\,\v\lambda t\v$).
It follows from (\ref{eqVII.1.12}) that
$$
\v k^+(t,x,y,\lambda)\v\leq \lambda^{\frac{3n}2}\,\bra\lambda t\ket^{-\frac n
2}\,\v\lambda t\v^{-n} \iint e^{-\frac \lambda{16}\frac{\v y-\tilde
\alpha_\xi\v^2}{\bra\lambda t\ket^2}}\, e^{-\frac\lambda 2 \v x-\tilde
\alpha_x\v^2}\,d\tilde \alpha\,.
$$
Setting $\tilde \alpha_\xi-y=\frac{4\bra\lambda t\ket}{\sqrt\lambda}\,z_1$,
$\tilde \alpha_x-x=\frac{\sqrt2}{\sqrt\lambda}\,z_2$, $Z=(z_1,z_2)$ we obtain
$$
\v k^+(t,x,y,\lambda)\v\leq C\,\lambda^{\frac{3n}2}\,\bra\lambda t\ket^{-\frac n
2}\,(\lambda t)^{-n}\, \bra\lambda t\ket^n\,\lambda^{-n }\int_{\R^{2n}} e^{-\v
Z\v^2}\,dZ
$$
so
$$
\v k^+(t,x,y,\lambda)\v\leq \frac C{t^{\frac n2}}\,,
$$
since $\bra\lambda t\ket\leq \sqrt 2\,\lambda t$. This proves (\ref{eqVII.1.8})
in this case.

{\bf Case B : } proof of (\ref{eqVII.1.4}) when $\lambda t\leq -1$.

In the right hand side of (\ref{eqVII.1.4}) we integrate on the support of
$\chi^+_3(\alpha_x)\cdot \psi_3(\alpha_\xi)$ on which we have
(\ref{eqVII.1.10}). We divide this support in two subsets $U_1$ and $U_2$ where
\begin{equation*}
\begin{aligned}
U_1&=\Big\{\alpha=(\alpha_x,\alpha_\xi)\in\supp
(\chi^+_3(\alpha_x)\,\psi_3(\alpha_\xi)) : -c_0\,
\bra\alpha_x\ket\,\v\alpha_\xi\v\leq \alpha_x\cdot \alpha_\xi\leq
c_0\bra\alpha_x\ket\,\v\alpha_\xi\v\Big\}\\
U_2&=\Big\{\alpha=(\alpha_x,\alpha_\xi)\in\supp
(\chi^+_3(\alpha_x)\,\psi_3(\alpha_\xi)) : \alpha_x\cdot \alpha_\xi\leq -c_0\,
\bra\alpha_x\ket\,\v\alpha_\xi\v\Big\}
\end{aligned}
\end{equation*}
According to Definition \ref{dIII.2.2} we have $U_1\subset \CS_+\cap\CS_-$
(which means that the points in $U_1$ are outgoing both for $\theta\geq 0$ and
$\theta\leq 0$).

According to Corollary \ref{cIII.3.3} we have (\ref{eqVII.1.11}) for
$\theta\in\R$ so in particular for $\theta=-\lambda t\geq 1$. Therefore the
same arguments as those used in case~1 work. It follows that the part of the
integral giving $k_+$ which concerns $U_1$ is bounded by $C\,\v t\v^{-\frac
n2}$. We consider now the integral on $U_2$. Here $\theta=-\lambda t\geq 1$ and
the points in $U_2$ are incoming for $\theta\geq 0$. The needed estimate on
$k_+$ will follow from the following result.
\begin{proposition}\sl \label{pVII.1.2}
One can find a function $g=g(\theta,y,\alpha_x)$ such that for all $\theta\geq
1$, all $\alpha\in U_2$ and all $y\in\supp \big (\chi^+_2(y)\,\chi_5\big
(\frac{y-x(\theta,\alpha)}{\bra\theta\ket}\big)\big)$ we have,
$$
\v\alpha_\xi-g(\theta,y,\alpha_x)\v\leq \frac C{\theta}\,\v
y-x(\theta,\alpha)\v\,.
$$
\end{proposition}
For the proof we need the following Lemma.
\begin{lemma}\sl\label{lVII.1.3}
Let $\alpha\in U_2$ and $\theta\geq 1$. Then for all $Y\in\R^n$ such that,
\begin{equation*}\left\{
\begin{array}{ll}
\textrm{(i)}\quad &Y\cdot \alpha_\xi\leq 20\,c_0\,\bra
Y\ket\,\v\alpha_\xi\v\,,\\
\textrm{(ii)}\quad &\Big\v\frac{Y-\alpha_x}{2\theta}-\alpha_\xi\Big\v\leq c_0\,,
\end{array}\right.
\end{equation*}
there exists a unique $\beta_\xi(\theta,Y,\alpha_x)\in\R^n$ such that
\begin{equation*}\left\{
\begin{array}{l}
\v\beta_\xi(\theta,Y,\alpha_x)-\alpha_\xi\v\leq 2\,c_0\,,\\
x(-\theta,Y,\beta_\xi(\theta,z,\alpha_x))=\alpha_x\,.
\end{array}\right.
\end{equation*}
\end{lemma}

{\bf Proof } Let $E=\{\beta_\xi\in\R^n : \v\beta_\xi-\alpha_\xi\v\leq
2\,c_0\}$. Then we have
$$
Y\cdot \beta_\xi=Y\cdot \alpha_\xi+Y\cdot (\beta_\xi-\alpha_\xi)\leq
20\,c_0\,\bra Y\ket\,\v\alpha_\xi\v+2\,c_0\,\bra Y\ket<\frac 14\,\bra
Y\ket\,\v\alpha_\xi\v\,.
$$
It follows that the point $(Y,\beta_\xi)$ belongs to $\CS_-$ (see Definition
\ref{dIII.2.2}) . Since $-\theta<0$ it follows from Proposition \ref{pIII.3.1}
that the equation $x(-\theta,Y,\beta_\xi)=\alpha_x$ is equivalent to,
$$
Y-2\theta\xi(-\theta,Y,\beta_\xi)+z(-\theta,Y,\beta_\xi)=\alpha_x\,,
$$
that is to the equation,
$Y-2\theta\,\beta_\xi-2\theta\,\zeta(-\theta,Y,\beta_\xi)+z(-\theta,Y,\beta_\xi)=\alpha_x$.
This equation can be written,
$$
\beta_\xi=\frac{Y-\alpha_x}{2\theta}-\zeta(-\theta,Y,\beta_\xi)+\frac
1{2\theta}\,z(-\theta,Y,\beta_\xi)=:F(\beta_\xi)\,.
$$
We show now that we can apply the fixed point theorem to $F$ in the set $E$.
First of all if $\beta_\xi\in E$ we have,
$$
\v F(\beta_\xi)-\alpha_\xi\v\leq
\Big\v\frac{Y-\alpha_x}{2\theta}-\alpha_\xi\Big\v+\v\zeta(-\theta,Y,\beta_\xi)\v+\frac
1{2\theta}\,\v z(-\theta,Y,\beta_\xi)\v\,.
$$
Using our assumption and Proposition \ref{pIII.3.2} we obtain,
$$
\v F(\beta_\xi)-\alpha_\xi\v\leq c_0+2\varepsilon\leq 2\,c_0\enskip
\textrm{if}\enskip 2\varepsilon\leq c_0\,.
$$
Now if $\beta_\xi\in E$, $\beta'_\xi\in E$ we have again by Proposition
\ref{pIII.3.2},
$$
\v F(\beta_\xi)-F(\beta'_\xi)\v\leq C\,\varepsilon\,\v\beta_\xi-\beta'_\xi\v\,.
$$
Taking $\varepsilon$ so small that $C\,\varepsilon<1$ we obtain the conclusion
of the Lemma. \cqfd
\begin{remark}\sl \label{rVII.1.4}
Let $(\theta,\alpha,y)$ as in Proposition \ref{pVII.1.2}. Then they satisfy
the conditions of Lemma \ref{lVII.1.3}. Indeed we have $y\in\supp \chi^+_4$
that is $-y\cdot \xi_0\leq 5\,\delta_1$ and
$\big\v\frac{\alpha_\xi}{\v\alpha_\xi\v}+\xi_0\v\leq 4\,\delta_2$. It follows
that 
$$
y\cdot \frac{\alpha_\xi}{\v\alpha_\xi\v}\leq y\cdot
\Big(\frac{\alpha_\xi}{\v\alpha_\xi\v}+\xi_0\Big)-y\cdot \xi_0\leq
4\,\delta_2\,\v y\v+5\,\delta_1< 10\,c_0\,\bra y\ket
$$
since $\delta_1$ and $\delta_2$ are small compared to $c_0$.

Moreover we have $\v y-x(\theta,\alpha)\v\leq \delta\,\bra\theta\ket\leq \sqrt
2\,\delta\theta$. By Proposition \ref{pIII.4.1} we have
$x(\theta,\alpha)=\alpha_x+2\theta\,\alpha_\xi+\CO(\varepsilon\,\bra\theta\ket)$
so
$$
\v y-\alpha_x-2\theta\,\alpha_\xi+\CO(\varepsilon\bra\theta\ket)\v\leq \sqrt
2\,\delta\theta\,.
$$
It follows that
$$
\Big\v\frac{y-\alpha_x}{2\theta}-\alpha_\xi\Big\v\leq C
\,(\delta+\varepsilon)\leq c_0\,,
$$
if $\varepsilon+\delta$ is small enough. \cqfd
\end{remark}

Now for fixed $\theta\geq 1$ and $\alpha\in U_2$ let us set
\begin{equation}\label{eqVII.1.13}
A=\left\{Y\in\R^n : Y\cdot \alpha_\xi<20\,c_0\,\bra
Y\ket\,\v\alpha_\xi\v\,,\quad \Big\v\frac{Y-\alpha_x}{2\theta}-\alpha_\xi
\Big\v<c_0\right\}\,,
\end{equation}
and for $Y\in A$ let us set
\begin{equation}\label{eqVII.1.14}
H(Y)=\xi(-\theta,Y,\beta_\xi(\theta,Y,\alpha_x))\,.
\end{equation}
\begin{lemma}\sl \label{lVII.1.5}
There exists a constant $C>0$ independent of $\theta,\alpha$ such that
$$
\Big\Vert \frac{\partial H}{\partial  Y}\,(Y)\Big\Vert \leq \frac C\theta
$$
for every $Y$ in $A$.
\end{lemma}

{\bf Proof } By Lemma \ref{lVII.1.3} we have for $j=1,\ldots ,n$,
$$
x_j(-\theta,Y,\beta_\xi(\theta,Y,\alpha_x))=\alpha^j_x\,.
$$
Let us differentiate this equality with respect to $Y_k$. We obtain
$$
\frac{\partial x_j}{\partial
x_k}\,(-\theta,Y,\beta_\xi(\theta,Y,\alpha_x))+\sum_\ell  \frac {\partial
x_j}{\partial \xi_\ell
}\,(-\theta,Y,\beta_\xi(\theta,Y,\alpha_x))\,\frac{\partial \beta^\ell
_\xi}{\partial Y_k}\,(\theta,Y,\alpha_x)=0\,.
$$
Since the point $(Y,\beta_\xi)$ belongs to $\CS_-$ we have by Corollary
\ref{cIII.3.3}
$$
\frac{\partial x_j}{\partial x_k}\,(-\theta,Y,\beta_\xi(\cdots
))=\delta_{jk}+\CO (\varepsilon\theta)\,,\quad \frac{\partial x_j}{\partial
\xi_\ell }\,(-\theta,Y,\beta_\xi(\cdots ))=-2\theta\,\delta_{jk}+\CO
(\varepsilon\theta)\,.
$$
It follows that 
\begin{equation}\label{eqVII.1.15}
\Big\Vert \frac{\partial \beta_\xi}{\partial Y}\,(\theta,Y,\alpha_x)\Big\Vert
\leq C \Big (\frac 1\theta+\varepsilon\Big)\,.
\end{equation}
Now thanks again to the fact that $(Y,\beta_\xi)\in\CO_-$ we deduce from
Proposition \ref{pIII.3.1} that
$$
\alpha^j_x=x_j(-\theta,Y,\beta_\xi(\theta,Y,\alpha_x))=Y_j-2\theta\,H_j(Y)+z_j(-\theta,Y,\beta_\xi
(\cdots))\,.
$$
Differentiating with respect to $Y_k$ yields
$$
\delta_{jk}-2\theta\,\frac{\partial H_j}{\partial Y_k}\,(Y)+\frac{\partial
z_j}{\partial x_k}\,(-\theta,Y,\beta_\xi(\cdots ))+\sum^n_{\ell =1}
\frac{\partial z_j}{\partial \xi_\ell }\,(-\theta,Y,\beta_\xi(\cdots))\,
\frac{\partial \beta_\xi}{\partial Y_k}\,(\theta,Y,\alpha_x)=0\,.
$$
Using Proposition \ref{pIII.3.2} and (\ref{eqVII.1.15}) we obtain the claim of
the Lemma. \cqfd

{\bf Proof of Proposition \ref{pVII.1.2} } We shall prove that the function $g$
defined by,
\begin{equation}\label{eqVII.1.16}
g(\theta,y,\alpha_x)=H(y)=\xi(-\theta,y,\beta_\xi(\theta,y,\alpha_x))\,,
\end{equation}
satisfies the claim of the Proposition.

To do so we must consider several cases.

{\bf Case 1 } Here we assume that
\begin{equation}\label{eqVII.1.17}\left\{
\begin{array}{l}
x(\theta,\alpha)\cdot \alpha_\xi\leq 10\,c_0\,\bra
x(\theta,\alpha)\ket\,\v\alpha_\xi\v\,,\\
\v y-x(\theta,\alpha)\v<\v x(\theta,\alpha)\v\,.
\end{array}\right.
\end{equation}
Let us show that this implies that
\begin{equation}\label{eqVII.1.18}
t\,y+(1-t)\,x(\theta,\alpha)\in A\enskip \textrm{for all}\enskip t\in[0,1]
\end{equation}
where $A$ has been defined in (\ref{eqVII.1.13}).

First of all using the Remark \ref{rVII.1.4} and (\ref{eqVII.1.17}) we can
write,
\begin{equation*}
\begin{aligned}
(ty+(1-t)\,x(\theta,\alpha))\cdot \alpha_\xi&<10\,c_0(t\,\bra
y\ket+(1-t)\bra x(\theta,\alpha)\ket)\,\v\alpha_\xi\v\\
&\leq 10\,c_0 [t(1+\v y\v)+(1-t)(1+\v x(\theta,\alpha)\v)]\,\v\alpha_\xi\v\,.
\end{aligned}
\end{equation*}
Now we use Lemma \ref{lIV.4.16} and we obtain
\begin{equation*}
\begin{aligned}
(ty+(1-t)\,x(\theta,\alpha))\cdot \alpha_\xi&< 10\,c_0(1+\sqrt 2\cdot \v
ty+(1-t)\,x(\theta,\alpha)\v) \v\alpha_\xi\v\\ 
&\leq 20\,c_0 \bra ty+(1-t)\,x(\theta,\alpha)\ket\,\v\alpha_\xi\v\,.
\end{aligned}
\end{equation*}
On the other hand we have
\begin{equation*}
\begin{aligned}
\v ty+(1-t)\,x(\theta,\alpha)-\alpha_x-2\theta\,\alpha_\xi\v&\leq t\,\v
y-x(\theta,\alpha)\v+\v x(\theta,\alpha)-\alpha_x-20\,\alpha_\xi\v\\
&\leq 2t\,\delta\theta+C\,\varepsilon\theta< 2c_0\,\theta\,,
\end{aligned}
\end{equation*}
since $y\in\supp \chi_5\big (\frac{y-x(\theta,\alpha)}{\bra\theta\ket}\big)$
and $\varepsilon,\delta$ are small compared to $c_0$.

In particular (\ref{eqVII.1.18}) for $t=0$ and $t=1$ shows that we can apply
Lemma \ref{lVII.1.3} to $Y=x(\theta,\alpha)$, $Y=y$. Therefore we have
$$
x(-\theta,x(\theta,\alpha),\beta_\xi(\theta,x(\theta,\alpha),\alpha_x))=\alpha_x\,.
$$
But we have also
$$
x(-\theta,x(\theta,\alpha),\xi(\theta,\alpha))=\alpha_x
$$
and since $\v\xi(\theta,\alpha)-\alpha_\xi\v\leq C\,\varepsilon\leq 2\,c_0$ it
follows from the uniqueness of $\beta_\xi$ in the set $E$ (see the proof of
Lemma \ref{lVII.1.3}) that
$\beta_\xi(\theta,x(\theta,\alpha),\alpha_x)=\xi(\theta,\alpha)$. Therefore we
have by (\ref{eqVII.1.16}),
$$
\alpha_\xi=\xi(-\theta,x(\theta,\alpha),\beta_\xi(\theta,x(\theta,\alpha),
\alpha_x))=H(x(\theta,\alpha))\,.
$$
Finally we write
\begin{equation*}
\v\alpha^j_\xi-g_j(\theta,y,\alpha_x)\v =\v H_j(x(\theta,\alpha))-H_j(y)\v\leq
\v y-x(\theta,\alpha)\v \int^1_0 \sum^n_{k=1} \Big\v\frac{\partial
H_j}{\partial Y_k}\,(ty+(1-t)\,x(\theta,\alpha))\Big\v\,dt
\end{equation*}
so using (\ref{eqVII.1.18}) and Lemma \ref{lVII.1.5} we obtain
$$
\v\alpha_\xi-g(\theta,y,\alpha_x)\v\leq \frac C\theta\,\v
y-x(\theta,\alpha)\v\,.
$$

{\bf Case 2 } Here we assume that
\begin{equation}\label{eqVII.1.19}\left\{
\begin{array}{l}
x(\theta,\alpha)\cdot \alpha_\xi\leq 10\,c_0\,\bra
x(\theta,\alpha)\ket\,\v\alpha_\xi\v\,,\\
\v y-x(\theta,\alpha)\v>\v x(\theta,\alpha)\v\,.
\end{array}\right.
\end{equation}
It follows that
\begin{equation}\label{eqVII.1.20}
\v y\v\leq 2\,\v y-x(\theta,\alpha)\v\,.
\end{equation}
We claim that in this case we have
\begin{equation}\label{eqVII.1.21}
ty\in A\,,\quad tx(\theta,\alpha)\in A\enskip \textrm{for}\enskip t\in[0,1]\,.
\end{equation}
Indeed we have by Remark \ref{rVII.1.4}
$$
ty\cdot \alpha_\xi\leq t\cdot 10\,c_0\,\bra y\ket\,\v\alpha_\xi\v\leq
10\,c_0\,\bra ty\ket\,\v\alpha_\xi\v
$$
\begin{equation*}
\begin{aligned}
\v ty-\alpha_x-2\theta\,\alpha_\xi\v&\leq t\, \v y\v+\v x(\theta,\alpha)\v+\v
x(\theta,\alpha)-\alpha_x-2\theta\,\alpha_\xi\v\\
&\leq 3\,\v y-x(\theta,\alpha)\v+C\,\varepsilon\,\bra\theta\ket\leq
2(C\,\varepsilon+3\delta)\,\theta\leq 2\,c_0\theta\,.
\end{aligned}
\end{equation*}
$$
tx(\theta,\alpha)\cdot \alpha_\xi\leq t\cdot 10\,c_0\,\bra
x(\theta,\alpha)\ket\,\v\alpha_\xi\v\leq 10\,c_0\,\bra
tx(\theta,\alpha)\ket\,\v\alpha_\xi\v
$$
$$
\v tx(\theta,\alpha)-\alpha_x-2\theta\,\alpha_\xi\v\leq (1-t)\v
x(\theta,\alpha)\v+\v x(\theta,\alpha)-\alpha_x-2\theta\,\alpha_\xi\v\leq
2\,c_0\,\theta\,. 
$$
As before we have
$$
(1)=\v\alpha_\xi-g(\theta,y,\alpha_x)\v=\v H(x(\theta,\alpha))-H(y)\v
$$
and we write
\begin{equation*}
\begin{aligned}
(1)&\leq \v H(x(\theta,\alpha))-H(0)\v+\v H(0)-H(y)\v\\
&\leq \v x(\theta,\alpha)\v \int^1_0 \Big\Vert \frac{\partial H}{\partial
Y}\,(tx(\theta,\alpha))\Big\Vert \,dt+\v y\v \int^1_0 \Big\Vert \frac{\partial
H}{\partial Y}\,(t,y)\Big\Vert \,dt\,.
\end{aligned}
\end{equation*}
Then we use (\ref{eqVII.1.21}), Lemma \ref{lVII.1.5} and (\ref{eqVII.1.19}),
(\ref{eqVII.1.20}) to conclude that 
$$
(1)\leq \frac {C}{\theta}\,\v y-x(\theta,\alpha)\v\,.
$$
The last case is the following.

{\bf Case 3 } We assume that
\begin{equation}\label{eqVII.1.22}
x(\theta,\alpha)\cdot \alpha_\xi>10\,c_0\,\bra
x(\theta,\alpha)\ket\,\v\alpha_\xi\v\,.
\end{equation}
Recall that $\alpha\in U_2$ that is in particular $\alpha_x\cdot \alpha_\xi\leq
-c_0\,\bra\alpha_x\ket\,\v\alpha_\xi\v$. By (\ref{eqIV.4.49}) there exists
$\theta^*\in]0,\theta[$ depending only on $\alpha$ such that
$x(\theta^*,\alpha)\cdot \alpha_\xi=0$. Then by Lemma \ref{lIV.4.17} we have
\begin{equation}\label{eqVII.1.23}\left\{
\begin{array}{l}
\frac32\, \v\theta-\theta^*\v\,\v\alpha_\xi\v\leq \v
x(\theta,\alpha)-x(\theta^*,\alpha)\v\leq
3\,\v\theta-\theta^*\v\,\v\alpha_\xi\v\\
\v\theta-\theta^*\v\geq \frac{c_0}{40}\\
\v y-x(\theta,\alpha)\v\leq \v y-x(\theta^*,\alpha)\v+\v
x(\theta^*,\alpha)-x(\theta,\alpha)\v\leq 5\,\v y-x(\theta,\alpha)\v\,.
\end{array}\right.
\end{equation}
Here again we consider two subcases.

{\bf Case 3.1 } $x(\theta,\alpha)\cdot \alpha_\xi>10\,c_0\,\bra
x(\theta,\alpha)\ket\,\v\alpha_\xi\v$, $\v y-x(\theta^*,\alpha)\v< \v
x(\theta^*,\alpha)\v$. 

In this case, (\ref{eqVII.1.23}) ensures that
\begin{equation}\label{eqVII.1.24}
ty+(1-t)\,x(\theta^*,\alpha)\in A\,.
\end{equation}
This follows from Lemma \ref{lIV.4.16}, since $t\,\v y\v\leq \sqrt2\,\v
ty+(1-t)\,x(\theta^*,\alpha)\v$ and from the following estimates,
\begin{equation*}
\begin{split}
\v ty+(1-t)\,&x(\theta^*,\alpha)-\alpha_x-2\theta\,\alpha_\xi\v\leq t\,\v
y-x(\theta^*,\alpha)\v+\v x(\theta^*,\alpha)-\alpha_x-2\theta\,\alpha_\xi\v\\
&\qquad \qquad \leq 5t\,\v
y-x(\theta,\alpha)\v+2\,\v\theta-\theta^*\v\,\v\alpha_\xi\v+\CO(\varepsilon\,\theta^*)\leq
C(\delta+\varepsilon)\,\theta\leq 2\,c_0\theta\,.
\end{split}
\end{equation*}
Then we write
$$
\v\alpha_\xi-g(y,\theta,\alpha_x)\v=\v H(y)-\alpha_\xi\v\leq
\underbrace{\v H(y)-H(x(\theta^*,\alpha))\v}_{(1)}+\underbrace{\v
H(x(\theta^*,\alpha))-\alpha_\xi\v}_{(2)}\,.
$$
By (\ref{eqVII.1.24}), Lemma \ref{lVII.1.5} and (\ref{eqVII.1.23}) we have
\begin{equation}\label{eqVII.1.25}
(1)\leq \frac C\theta\,\v y-x(\theta,\alpha)\v\,.
\end{equation}
Now we have
$$
x\Big(-\theta,x(\theta^*,\alpha),\frac{\theta^*}\theta\,\xi(\theta^*,\alpha)\Big)
=x(-\theta^*,x(\theta^*,\alpha),\xi(\theta^*,\alpha))=\alpha_x\,.
$$
Therefore
$$
\beta_\xi(\theta,x(\theta^*,\alpha),\alpha_x)=\frac{\theta^*}\theta\,\xi(\theta^*,\alpha)\,.
$$
It follows that
$$
H(x(\theta^*,\alpha))=\xi\Big
(-\theta,x(\theta^*,\alpha),\frac{\theta^*}\theta\,\xi(\theta^*,\alpha)\Big)=\frac
{\theta^*}\theta\,\xi(-\theta^*,x(\theta^*,\alpha),\xi(\theta^*,\alpha))
=\frac{\theta^*}\theta\,\alpha_\xi\,.
$$
Therefore
\begin{equation}\label{eqVII.1.26}
(2)=\v\alpha_\xi\v\,\frac{\v\theta-\theta^*\v}\theta\leq \frac C\theta\,
\v y-x(\theta,\alpha)\v
\end{equation}
by (\ref{eqVII.1.23}). The estimates obtained on (1) and (2) show that
$$
\v\alpha_\xi-g(\theta,y,\alpha_x)\v\leq \frac C\theta\,\v y-x(\theta,\alpha)\v
$$
which is the claim of Proposition \ref{pVII.1.2}.

The last step concerns the following case.

{\bf Case 3.2 } $x(\theta,\alpha)\cdot \alpha_\xi>10\,c_0\,\bra
x(\theta,\alpha)\ket\,\v\alpha_\xi\v$, $\v x(\theta^*,\alpha)\v\leq \v
y-x(\theta^*,\alpha)\v$.

It follows then from (\ref{eqVII.1.23}) that we have
\begin{equation}\label{eqVII.1.27}
\v y\v +\v x(\theta^*,\alpha)\v\leq C\,\v y-x(\theta,\alpha)\v\,.
\end{equation}
Moreover we have
\begin{equation}\label{eqVII.1.28}
ty\in A\,,\quad tx(\theta^*,\alpha)\in A\enskip \textrm{for}\enskip
t\in[0,1]\,.
\end{equation}
Indeed we can write,
$$
ty\cdot \alpha_\xi\leq t\cdot 10\,c_0\,\bra y\ket\,\v\alpha_\xi\v\leq
10\,c_0\,\bra ty\ket\,\v\alpha_\xi\v
$$
\begin{equation*}
\begin{aligned}
\v ty-\alpha_x-2\theta\,\alpha_\xi\v&\leq t\, \v y\v+\v x(\theta^*,\alpha)\v+\v
x(\theta^*,\alpha)-\alpha_x-2\theta\,\alpha_\xi\v\\
&\leq C(\v
y-x(\theta,\alpha)\v+\v\theta-\theta^*\v+C\,
\varepsilon\,\theta^*)\leq 2\,c_0\theta\,
\end{aligned}
\end{equation*}
$$
tx(\theta^*,\alpha)\cdot \alpha_\xi=0
$$
\begin{equation*}
\begin{aligned}
\v tx(\theta^*,\alpha)-\alpha_x-2\theta\,\alpha_\xi\v&\leq (1-t)\v
x(\theta^*,\alpha)\v+\v x(\theta^*,\alpha)-\alpha_x-2\theta\,\alpha_\xi\v\\
&\leq C (1-t) \v y-x(\theta,\alpha)\v +C \v
\theta-\theta^*\v+C\,\varepsilon\,\theta^*\\
&\leq 2\,c_0\,\theta
\end{aligned}
\end{equation*}
by (\ref{eqVII.1.23}). Then we can write
\begin{equation*}
\begin{aligned}
\v\alpha_\xi-g(\theta,y,\alpha_x)\v&\leq \v H(y)-H(x(\theta^*,\alpha)\v+\v
H(x(\theta^*,\alpha)-\alpha_\xi\v\\
&\leq \v H(y)-H(0)\v+\v H(0)-H(x(\theta^*,\alpha))\v+\v
H(x(\theta^*,\alpha)-\alpha_\xi\v\,.
\end{aligned}
\end{equation*}
Using (\ref{eqVII.1.28}), Lemma \ref{lVII.1.5}, (\ref{eqVII.1.26}) we obtain
$$
\v\alpha_\xi-g(\theta,y,\alpha_x)\v\leq \frac C \theta \v\,
y-x(\theta,\alpha)\v\,.
$$
This completes the proof of Proposition \ref{pVII.1.2}. \cqfd

{\bf End of the proof of (\ref{eqVII.1.8}) in case B, $(\lambda t\leq -1)$}

For the part of the integral in (\ref{eqVII.1.4}) (giving $k_+$) where
$\alpha\in U_2$ we use Proposition \ref{pVII.1.2}. Let us call it (1). As in
(\ref{eqVII.1.12}) we have
\begin{equation*}
\begin{aligned}
\v(1)\v&\leq C\,\lambda^{\frac{3n}2} \int e^{-\frac{1}{16}\frac{\v
y-x(-\lambda\, t,\alpha)\v^2}{\bra\lambda t\ket^2}-\frac{\lambda}{2}\,\v
x-\alpha_x\v^2}\,\bra\lambda t\ket^{-\frac {n}{2}}\,d\alpha\\
\v(1)\v&\leq C\,\lambda^{\frac{3n}2} \bra\lambda t\ket^{-\frac n 2} \iint 
e^{-\varepsilon_0\,\frac{\v\lambda\, t\v^2\lambda}{\bra\lambda
t\ket^2}\,\v\alpha_\xi-g(-\lambda t,y,\alpha_x)\v^2}\,e^{-\frac\lambda 2 \,\v
x-\alpha_x\v^2}\,d\alpha_x\,d\alpha_\xi\\
\v(1)\v&\leq C'\, \frac{\lambda^{\frac{3n}2}\bra\lambda\, t\ket^{-\frac
n2}\bra\lambda \,t\ket^n}{\lambda^n\,\v\lambda\, t\v^n}\leq \frac{C''}{\v
t\v^{\frac n 2}}\,.
\end{aligned}
\end{equation*}

{\bf Case C : } proof of (\ref{eqVII.1.8}) when $\v\lambda t\v\leq 1$.

In this case the proof made above does not give the needed result since
$\bra\lambda t\ket \approx 1$. We will use instead a stationnary phase method.
Let us set $\theta=-\lambda t$ and let us recall (see (\ref{eqVII.1.4})) that
we have to bound the following kernel
\begin{equation}\label{eqVII.1.29}
\begin{split}
k_+(t,x,y,\lambda)=\lambda^{\frac{3n}2} \int e^{i\lambda
F(\theta,x,y,\alpha)}\,a(\theta,y,\alpha)\,\chi^+_2(y)\,&\chi^+_3(\alpha_x)\,\psi_3(\alpha_\xi)\\
&\chi_5\Big (\frac{y-x(\theta,\alpha)}{\bra\theta\ket}\Big)\,d\alpha\,.
\end{split}
\end{equation}
Let us recall also that, according to the Theorems \ref{tVI.4.1} and
\ref{tIV.1.2} we have,
\begin{equation}\label{eqVII.1.30}
\Im F(\theta,x,y,\alpha)\geq \frac 12\,\v x-\alpha_x\v^2\,.
\end{equation}
Let $\chi_6\in C^\infty _0\R^n)$ be such that,
\begin{equation}\label{eqVII.1.31}\left\{
\begin{array}{ll}
\chi_6(x)=1\enskip \textrm{ if } \enskip \v x\v\leq \frac 12\,,\\
\chi_6(x)=0\enskip \textrm{ if } \enskip \v x\v\geq1\,,
\end{array}\right.
\end{equation}
We write in the integral in the right hand side of (\ref{eqVII.1.29})
$1=\chi_6(x-\alpha_x)+1-\chi_6(x-\alpha_x)$.

The part of the integral containing $1-\chi_6(x-\alpha_x)$ can be bounded,
using (\ref{eqVII.1.30}), by the quantity
$$
C\,\lambda^{\frac{3n}2}\, e^{-\frac {\lambda}{16}} \iint e^{-\frac \lambda
4\,\v x-\alpha_x\v^2}\,\v\psi_3(\alpha_\xi)\v\,d\alpha
$$
which is $0(1)$ uniformly in $(t,x,y,\lambda)$. Setting
$a_1(\theta,y,\alpha)=a(\theta,y,\alpha)\,\chi^+_2(y)\,\chi^+_3(\alpha_x)$ we
see therefore that we are left with the bound of the following kernel.
\begin{equation}\label{eqVII.1.32}
\tilde k_+(t,x,y,\lambda)=\lambda^{\frac{3n}2} \int e^{i\lambda
F(\theta,x,y,\alpha)}\,a_1(\theta,y,\alpha)\,\chi_6(x-\alpha_x)\,\psi_3(\alpha_\xi)
\chi_5\Big (\frac{y-x(\theta,\alpha)}{\bra\theta\ket}\Big)\,d\alpha\,.
\end{equation}
Now, according to Theorem \ref{tVI.4.1} we have,
$$
F(\theta,x,y,\alpha)=\varphi (\theta,y,\alpha)-(x-\alpha_x)\cdot
\alpha_\xi+\frac i 2\,\v x-\alpha_x\v^2-\frac 1{2i}\,\v\alpha_\xi\v^2\,.
$$
Using Theorem \ref{tIV.5.1} we obtain
\begin{equation}\label{eqVII.1.33}
\begin{split}
F(\theta,x,y,\alpha)&=\frac{(y-x)\cdot \alpha_\xi-2i\theta(x-\alpha_x)\cdot
\alpha_\xi-\theta\,\v\alpha_\xi\v^2-\theta \v x-\alpha_x\v^2+\frac i 2\,\v
x-\alpha_x\v^2+\frac i 2\, \v y-\alpha_x\v^2}{1+2i\theta}\\
&+R(\theta,y,\alpha)
\end{split}
\end{equation}
where
\begin{equation}\label{eqVII.1.34}\left\{
\begin{array}{ll}
\Big\v\frac{\partial R}{\partial \alpha_x}\Big\v+\Big\v\frac{\partial
^2R}{\partial \alpha^2_x}\Big\v&\leq C (\varepsilon+\delta)(\v
y-\alpha_x\v+\v\theta\v)\,,\\
\Big\v\frac{\partial R}{\partial \alpha_\xi}\Big\v+\Big\v\frac{\partial
^2R}{\partial \alpha^2_\xi}\Big\v+\Big\v\frac{\partial
^2R}{\partial \alpha_x\partial \alpha_\xi}\Big\v&\leq C
(\varepsilon+\delta)\,\v\theta\v\,,\\
\v\partial ^{A_1}_{\alpha_x}\,\partial ^{A_2}_{\alpha_\xi}
\,R\v&\leq 
\begin{cases}
C_{A_1}& \textrm{ if  \enskip $A_2=0$}\,,\\
C_{A_1,A_2}\v\theta\v& \textrm{ if \enskip  $\v A_2\v\geq 1$}\,.
\end{cases}
\end{array}\right.
\end{equation}
It follows that we have,
\begin{equation}\label{eqVII.1.35}\left\{
\begin{array}{ll}
\frac{\partial F}{\partial \alpha_\xi}&=\frac 1{1+2i\theta}
[(y-x)-2i\theta(x-\alpha_x)-2\theta\,\alpha_\xi]+\frac{\partial R}{\partial
\alpha_\xi}\\
\frac{\partial F}{\partial \alpha_x}&=\frac 1{1+2i\theta}
[2i\theta\,\alpha_\xi+2\theta(x-\alpha_x)-2i(x-\alpha_x)-i(y-x)]+\frac{\partial
R}{\partial\alpha_x}\\
\frac{\partial ^2F}{\partial \alpha^j_\xi\,\partial
\alpha^k_\xi}&=\frac{-2\theta\,\delta_{jk}}{1+2i\theta}
+\frac{\partial ^2R}{\partial \alpha^j_\xi\,\partial \alpha^k_\xi}\\
\frac{\partial ^2F}{\partial \alpha^j_\xi\,\partial \alpha^k_x}&=
\frac{2i\theta\,\delta_{jk}}{1+2i\theta}+
\frac{\partial ^2R}{\partial \alpha^j_\xi\,\partial \alpha^k_x}\\
\frac{\partial ^2F}{\partial \alpha^j_x\,\partial
\alpha^k_x}&=\frac{(2i-2\theta)\delta_{jk}}{1+2i\theta}+\frac{\partial
^2R}{\partial \alpha^j_x\,\partial \alpha^k_x}\\
\frac{\partial ^{\v A\v}F}{\partial \alpha^A}&=\frac{\partial^ {\v
A\v}R}{\partial \alpha^A}\enskip \textrm{ if }\enskip \v A\v\geq 3\,.
\end{array}\right.
\end{equation}
Let us recall that (\ref{eqVI.2.9}) shows that
\begin{equation}\label{eqVII.1.36}
\supp \psi_3(\alpha_\xi)\subset\big\{\alpha_\xi : a-4\delta_2\leq
\v\alpha_\xi\v\leq b+4\delta_2\big\}\,.
\end{equation}
We shall divide the proof into three cases
\begin{equation}\label{eqVII.1.37}\left\{
\begin{array}{l}
\textrm{case 1 : }\enskip \frac{\v x-y\v}{2\v\theta\v}\leq a-5\delta_2\\
\textrm{case 2 : }\enskip \frac{\v x-y\v}{2\v\theta\v}\geq b+5\delta_2\\
\textrm{case 3 : }\enskip a-5\delta_2\leq  \frac{\v x-y\v}{2\v\theta\v}\leq
b+5\delta_2
\end{array}\right.
\end{equation}

{\bf Case 1 } We have the following result.
\begin{lemma}\sl \label{lVII.1.6}
When $x-\alpha_x\in \supp \chi_6$, $y-x(\theta,\alpha)\in \supp \chi_5$,
$\alpha_\xi\in\supp \psi_3$ we have
$$
Q=:\Big\v\frac{\partial F}{\partial \alpha_x}\Big\v^2+\frac
1{\v\theta\v}\,\Big\v\frac{\partial F}{\partial \alpha_\xi}\Big\v^2\geq C (\v
x-\alpha_x\v^2+\v\theta\v)\,,
$$
$$
\big\v\partial ^{A_1}_{\alpha_x}\,\partial ^{A_2}_{\alpha_\xi}\,F\big\v\leq 
\begin{cases}
C_{A_1}& \textrm{ if  \enskip $\v A_1\v\geq 1,\enskip A_2=0$}\,,\\
C_{A_1,A_2}\v\theta\v& \textrm{ if \enskip  $\v A_2\v\geq 1$}\,.
\end{cases}
$$
uniformly in $(\theta,x,y,\alpha)$.
\end{lemma}

{\bf Proof } Let us set $X=x-\alpha_x$, $Y=\alpha_\xi-\frac{y-x}{2\theta}$.
Using (\ref{eqVII.1.35}) we see that $\frac{\partial F}{\partial
\alpha_\xi}=\frac{-2\theta}{1+2i\theta}\,(Y+i\,X)+\frac{\partial R}{\partial
\alpha_\xi}$. It follows from (\ref{eqVII.1.34}) that
$\frac{2\v\theta\v}{\v1+2i\theta\v}\,\v Y+i\,X\v\leq \big\v\frac{\partial
F}{\partial \alpha_\xi}\big \v+C (\varepsilon+\delta)\v\theta\v$. Therefore
\begin{equation}\label{eqVII.1.38}
\frac{4\v\theta\v}{1+4\theta^2}\,(\v X\v^2+\v Y\v^2)\leq \frac
2{\v\theta\v}\,\Big\v\frac{\partial F}{\partial
\alpha_\xi}\Big\v^2+C_1(\varepsilon+\delta)^2\,\v\theta\v\,.
\end{equation}
By the same way we have
$$
\frac{\partial F}{\partial \alpha_x}=\frac 1{1+2i\theta}
\,(2i\theta Y-2iX+2\theta X)+\frac{\partial R}{\partial \alpha_x}\,.
$$
Since $\v y-\alpha_x\v\leq \v y-x\v+\v x-\alpha_x\v\leq \v
x-\alpha_x\v+C\,\v\theta\v$, we obtain
$$
\frac 4{1+4\theta^2}\, (\v\theta Y-X\v^2+\theta^2\,\v X\v^2)\leq
2\Big\v\frac{\partial F}{\partial \alpha_x}\Big\v^2+C_2
(\varepsilon+\delta)^2(\v\theta\v+\v X\v^2)\,.
$$
It follows that
\begin{equation}\label{eqVII.1.39}
2\Big(\Big\v\frac{\partial F}{\partial \alpha_x}\Big\v^2+\frac 1{\v\theta\v }
\Big(\Big\v\frac{\partial F}{\partial \alpha_\xi}\Big\v^2\Big)+C_3
(\varepsilon+\delta)^2 (\v\theta\v+\v X\v^2)\geq \frac 4{1+4\theta^2} (\v\theta
Y-X\v^2+\v\theta\v\,\v Y\v^2)\,.
\end{equation}
Since $2\,\v\theta\v\v X\v\v Y\v\leq \frac 23\,\v X\v^2+\frac 32 \,\theta^2\,\v
Y\v^2$ we can write
$$
\v\theta Y-X\v^2+\v\theta\v\v Y\v^2\geq \frac 13 \,\v X\v^2-\frac
12\,\theta^2\,\v Y\v^2+\v \theta\v\,\v Y\v^2\geq \frac 13\,\v X\v^2+\frac
12\,\v \theta\v\v Y\v^2
$$
since $\v\theta\v\leq 1$. We deduce from (\ref{eqVII.1.39}) that
$$
Q\geq \frac 2{15}\,(\v X\v^2+\v\theta\v\,\v
Y\v^2)-C_4(\varepsilon+\delta)^2(\v\theta\v+\v X\v^2)\,.
$$
Now in case~1 we have, according to (\ref{eqVII.1.37})
$$
\v Y\v=\Big\v\alpha_\xi-\frac{y-x}{2\theta}\Big\v\geq \v\alpha_\xi\v-\frac{\v
y-x\v} {2\v\theta\v}\geq \delta_2\,.
$$
Taking $\varepsilon$ and $\delta$ small compared to $\delta_2$ we deduce that
$Q\geq c (\v X\v^2+\v\theta\v)$ which is the first claim of the Lemma. The
second claim follows easily from the fact that $\v Y\v$, $\v X\v$ are uniformly
bounded and from (\ref{eqVII.1.35}), (\ref{eqVII.1.34}). \cqfd

{\bf Case 2 } We have the following result.
\begin{lemma}\sl \label{lVII.1.7}
When $x-\alpha_x\in\supp \chi_6$, $y-x(\theta,\alpha)\in \supp \chi_5$,
$\alpha_\xi\in\supp \psi_3$ we have
$$
Q=:\Big\v\frac{\partial F}{\partial \alpha_x}\Big\v^2+\frac 1{\v
y-x\v}\,\Big\v\frac{\partial F}{\partial \alpha_\xi}\Big\v^2\geq c (\v
x-\alpha_x\v^2+\v\theta\v)\,,
$$
$$
\v\partial ^A_{\alpha_x}F\v\leq C_A\enskip \textrm{ if }\enskip \v A\v\geq 1\,,
$$
$$
\v\partial ^{A_1}_{\alpha_x}\,\partial ^{A_2}_{\alpha_\xi}F\v\leq
C_{A_1\,A_2}\,\v y-x\v\enskip \textrm{ if }\enskip \v A_2\v\geq 1\,.
$$
\end{lemma}

{\bf Proof } Here we set $X=x-\alpha_x$, $Y=\frac{y-x}{\v
y-x\v}-2\theta\,\frac{\alpha_\xi}{\v y-x\v}$. Then we can write,
\begin{equation*}
\begin{aligned}
\frac{\partial F}{\partial \alpha_\xi}&=\frac 1{1+2i\theta}\, (\v
y-x\v\,Y-2i\theta X)+\frac{\partial R}{\partial \alpha_\xi}\,,\\
\frac{\partial F}{\partial \alpha_x}&=\frac 1{1+2i\theta}\, (2\theta X-i(2X+\v
y-x\v\,Y))+\frac{\partial R}{\partial \alpha_x}\,.
\end{aligned}
\end{equation*}
It follows that
$$
\frac 1{1+4\theta^2}\,(\v y-x\v^2\,\v Y\v^2+4\theta^2\,\v X\v^2)\leq 2\Big
(\Big\v\frac{\partial F}{\partial
\alpha_\xi}\Big\v^2+C_1(\varepsilon+\delta)^2\,\v\theta\v^2\Big)\,.
$$
Therefore we have the estimate
$$
\v y-x\v\,\v Y\v^2\leq 10\,\frac 1{\v y-x\v}\,\Big\v\frac{\partial F}{\partial
\alpha_\xi}\Big\v^2+C_2 (\varepsilon+\delta)^2\,\v\theta\v\,.
$$
On the other hand since $\v y-\alpha_x\v\leq \v y-x\v+\v x-\alpha_x\v$
$$
\v 2X+\v y-x\v\,Y\v^2\leq 10\,\Big\v\frac{\partial F}{\partial
\alpha_x}\Big\v^2+C_3 (\varepsilon+\delta)^2(\v\theta\v^2+\v y-x\v^2+\v
x-\alpha_x\v^2)\,.
$$
Summing up we see that
$$
10\,Q+C_4 (\varepsilon+\delta)^2(\v\theta\v^2+\v y-x\v^2+\v X\v^2)\geq \v 2X+\v
y-x\v\,Y\v^2+\v y-x\v\,\v Y\v^2\,.
$$
Writing $y-x=y-x(\theta,\alpha)+x(\theta,\alpha)-\alpha_x+\alpha_x-x$, we see
that one can find a constant $K_0$ such that $\v y-x\v\leq K_0$. Let $\eta_0>0$
be such that $\frac{\eta_0}{4-\eta_0}\,K_0\leq \frac 12$. Then we have
$$
4Ê\,\v y-x\v\,\v X\v\cdot \v Y\v\leq (4-\eta_0)\,\v X\v^2+\frac 4{4-\eta_0}\,\v
y-x\v^2\,\v Y\v^2\,.
$$
It follows that
$$
\v 2X+\v y-x\v\,Y\v^2\geq \eta_0\,\v X\v^2-\frac{\eta_0}{4-\eta_0}\,\v
y-x\v^2\,\v Y\v^2\,.
$$
This implies that
\begin{equation*}
\begin{aligned}
10\,Q+C_4(\varepsilon+\delta)^2(\v\theta\v^2+\v y-x\v^2+\v X\v^2)&\geq
\eta_0\,\v X\v^2+\v y-x\v\,\v Y\v^2\Big (1-\frac{\eta_0}{4-\eta_0}\,K_0\Big)\\
&\geq \min \Big(\eta_0,\,\frac12\Big)(\v X\v^2+\v y-x\v\,\v Y\v^2)\,.
\end{aligned}
\end{equation*}
On the other hand we have
$$
\v Y\v\geq 1-2\,\v\theta\v\, \frac{\v\alpha_\xi\v}{\v y-x\v}\geq
1-\frac{b+4\delta_2}{b+5\delta_2}=\frac{\delta_2}{b+5\delta_2}
$$
because $\v y-x\v\geq 2(b+5\delta_2)\v\theta\v$, since we are in case 2.
Therefore
$$
10\,Q+C_5(\varepsilon+\delta)(\v\theta\v+\v y-x\v+\v X\v^2)\geq C_6(\v X\v^2+\v
y-x\v)\,.
$$
Taking $\varepsilon+\delta$ small with respect to $C_6$ we obtain the first
part of the Lemma.

The bounds on the derivatives of $F$ can be easily obtained since
$\v\theta\v\leq (b+5\,\delta_2)^{-1}\,\v y-x\v$.

{\bf Case 3 } Recall that we have $a-5\,\delta_2\leq \frac{\v
x-y\v}{2\v\theta\v}\leq b+5\,\delta_2$. Then we have the following result.
\begin{lemma}\sl \label{lVII.1.8}
One can find $\rho_0>0$ such that the equation $x(\theta,x,\alpha_\xi)=y$ has a
unique solution in the set $E=\big\{\alpha_\xi :
\big\v\alpha_\xi-\frac{y-x}{2\theta}\big\v\leq \rho_0\big\}$.
\end{lemma}

{\bf Proof } First of all one can find $\rho_0$ such that if $\alpha_\xi\in E$
then $\frac 12\leq \v\alpha_\xi\v\leq 2$. It follows from Proposition
\ref{pIII.2.1} that
$x(\theta,x,\alpha_\xi)=x+2\theta\,\alpha_\xi+r(\theta,x,\alpha_\xi)$ where $\v
r\v\leq C\,\varepsilon\,\v\theta\v$. The equation to be solved is in $E$
equivalent to the equation $\alpha_\xi=\frac{y-x}{2\theta}+\frac
1{2\theta}\,r(\theta,x,\alpha_\xi)$. If we set
$\Phi(\alpha_\xi)=\frac{y-x}{2\theta}+\frac 1{2\theta}\,u(\theta,x,\alpha_\xi)$
then $\Phi$ maps $E$ in itself if $C\,\varepsilon<\rho_0$. Moreover, again by
Proposition \ref{pIII.2.1} we have $\v\Phi(\alpha_\xi)-\Phi(\alpha'_\xi)\v\leq
C'\,\varepsilon\,\v\alpha_\xi-\alpha'_\xi\v$. Taking $\varepsilon$ small
enough the Lemma follows from the fixed point theorem. \cqfd

We shall set
\begin{equation}\label{eqVII.1.40}
\alpha_c=(x,\alpha^c_\xi)
\end{equation}
where $\alpha^c_\xi$ is the unique solution of $x(\theta,x,\alpha_\xi)=y$ given
by Lemma \ref{lVII.1.8}. Then we have the following result.
\begin{lemma}\sl\label{lVII.1.9}
We have $\frac{\partial F}{\partial
\alpha_x}\,(\theta,x,y,\alpha^c)=\frac{\partial F}{\partial
\alpha_\xi}\,(\theta,x,y,\alpha^c)=0\,.$
\end{lemma}

{\bf Proof } Let us recall that
\begin{equation}\label{eqVII.1.41}\left\{
\begin{array}{l}
\frac{\partial F}{\partial\alpha_x}=\alpha_\xi-i(x-\alpha_x)+\frac{\partial
\varphi }{\partial\alpha_x}\\
\frac{\partial F}{\partial\alpha_\xi}=-(x-\alpha_x)+i\alpha_\xi+\frac{\partial
\varphi }{\partial\alpha_\xi}\,.
\end{array}\right.
\end{equation}
On the other hand we have (see (\ref{eqIV.5.17})) 
\begin{equation}\label{eqVII.1.42}
\varphi (\theta,x(\theta,\alpha),\alpha)=\theta\,p(\alpha)+\frac
1{2i}\,\v\alpha_\xi\v^2\,.
\end{equation}
Therefore
$$
\frac{\partial \varphi }{\partial
\alpha^j_x}\,(\theta,x(\theta,\alpha),\alpha)=\theta\,\frac{\partial
p}{\partial x_j}\,(\alpha)-\sum^n_{k=1}\frac{\partial \varphi }{\partial
x_k}\,(\theta,x(\theta,\alpha),\alpha)\,\frac{\partial x_k}{\partial
x_j}\,(\theta,\alpha)\,.
$$
Moreover we have
$$
\frac{\partial \varphi }{\partial
x_k}\,(\theta,x(\theta,\alpha),\alpha)=\Phi(\theta,x(\theta,\alpha),\alpha)=\xi(\theta,\alpha)
$$
so
\begin{equation}\label{eqVII.1.43}
\frac{\partial \varphi }{\partial
\alpha^j_x}\,(\theta,x(\theta,\alpha),\alpha)=\theta\,\frac{\partial
p}{\partial x_j}\,(\alpha)-\sum^n_{k=1} \xi_k(\theta,\alpha)\,\frac{\partial
x_j}{\partial x_j}\,(\theta,\alpha)\,.
\end{equation}
Now by the definition of the flow and the Euler relation we have
\begin{equation*}\left\{
\begin{array}{l}
\som^n_{k=1} \dot x_k(\theta,\alpha)\,\xi_k(\theta,\alpha)=\som^n_{k=1}
\xi_k(\theta,\alpha)\,\frac{\partial p}{\partial
\xi_k}\,(x(\theta,\alpha),\xi(\theta,\alpha))=2p(\alpha)\\
p(x(\theta,\alpha),\xi(\theta,\alpha))=p(\alpha)\,.
\end{array}\right.
\end{equation*}
Differentiating these two relations with respect to $\alpha^j_x$ we obtain
\begin{equation}\label{eqVII.1.44}
\begin{split}
&\frac{\partial \dot x}{\partial x_j}\,(\theta,\alpha)\cdot
\xi(\theta,\alpha)+\dot x(\theta,\alpha)\cdot \frac{\partial \xi}{\partial
x_j}\,(\theta,\alpha)=2\,\frac{\partial p}{\partial x_j}\,(\alpha)\\
&\frac{\partial p}{\partial x}\,(x(\theta,\alpha),\xi(\theta,\alpha))\cdot
\frac{\partial x}{\partial x_j}\,(\theta,\alpha)+ \frac{\partial
p}{\partial \xi}\,(x(\theta,\alpha),\xi(\theta,\alpha))\,\frac{\partial
\xi}{\partial x_j}\,(\theta,\alpha)=\frac{\partial p}{\partial x_j}\,(\alpha)\,.
\end{split}
\end{equation}
Using the equations of the flow the last equality can be written as
$$
-\dot\xi(\theta,\alpha)\cdot \frac{\partial x}{\partial
x_j}\,(\theta,\alpha)+\dot x(\theta,\alpha)\frac{\partial \xi}{\partial
x_j}(\theta,\alpha)=\frac{\partial p}{\partial x_j}(\alpha)\,.
$$
Combining with (\ref{eqVII.1.44}) we obtain
$$
\frac{\partial \dot x}{\partial x_j}\,(\theta,\alpha)\cdot
\xi(\theta,\alpha)+\dot\xi(\theta,\alpha)\cdot \frac{\partial x}{\partial
x_j}\,(\theta,\alpha)=\frac{\partial p}{\partial x_j}\,(\alpha)
$$
which can be written as
$$
\frac d{d\theta} \Big[\frac{\partial x}{\partial x_j}\,(\theta,\alpha)\cdot
\xi(\theta,\alpha)\Big]=\frac{\partial p}{\partial x_j}\,(\alpha)\,.
$$
Integrating both side we obtain
$$
\frac{\partial x}{\partial x_j}(\theta,\alpha)\cdot
\xi(\theta,\alpha)=\theta\,\frac{\partial p}{\partial
x_j}\,(\alpha)+\alpha^j_\xi\,.
$$
Using (\ref{eqVII.1.43}) we deduce
$$
\frac{\partial \varphi }{\partial
\alpha^j_x}\,(\theta,x(\theta,\alpha),\alpha)=-\alpha^j_\xi\,.
$$
Since $\alpha^c=(x,\alpha^c_\xi)$ where $x(\theta,x,\alpha^c_\xi)=y$ we deduce
from (\ref{eqVII.1.41}) that
$$
\frac{\partial F}{\partial
\alpha^j_x}\,(\theta,x,y,\alpha^c)=\alpha^j_\xi+\frac{\partial \varphi
}{\partial
\alpha^j_x}\,(\theta,x(\theta,\alpha^c),\alpha^c)=\alpha^j_\xi-\alpha^j_\xi=0\,.
$$
Now differentiating (\ref{eqVII.1.42}) with respect to $\alpha^j_\xi$ yields
$$
\frac{\partial \varphi }{\partial
\alpha^j_\xi}\,(\theta,x(\theta,\alpha),\alpha)+\frac{\partial \varphi
}{\partial x}\,(\theta,x(\theta,\alpha),\alpha)\cdot \frac{\partial x}{\partial
\alpha^j_\xi}\,(\theta,\alpha)=\theta\,\frac{\partial p}{\partial
\xi_j}\,(\alpha)+\frac 1 i\,\alpha^j_\xi\,.
$$
As above we see easily that,
$$
\frac d{d\theta} \Big[\frac{\partial x}{\partial
\alpha^j_\xi}\,(\theta,\alpha)\cdot \xi(\theta,\alpha)\Big]=\frac{\partial
p}{\partial \xi_j}\,(\alpha)\,,
$$
from which we deduce that $\frac{\partial x}{\partial
\alpha^j_\xi}\,(\theta,\alpha)\cdot \xi(\theta,\alpha)=\theta\,\frac{\partial
p}{\partial \xi_j}\,(\alpha)$. Since $\frac{\partial \varphi }{\partial
x}\,(\theta,x(\theta,\alpha),\xi(\theta,\alpha))=\xi(\theta,\alpha)$, we obtain
$$
\frac{\partial \varphi }{\partial
\alpha^j_\xi}\,(\theta,x(\theta,\alpha),\alpha)=\frac 1 i\,\alpha^j_\xi
$$
which implies that $\frac{\partial F}{\partial
\alpha^j_\xi}\,(\theta,x,y,\alpha^c)=0$. \cqfd
\begin{lemma}\sl\label{lVII.1.10}
Let us set
\begin{equation*}\left\{
\begin{array}{l}
Q=\Big[\Big\v\frac{\partial F}{\partial \alpha_x}\Big\v^2+\frac
1{\v\theta\v}\,\Big\v\frac{\partial F}{\partial
\alpha_\xi}\Big\v^2\Big](\theta,x,y,\alpha)\\
X=\alpha_x-x\\
Y=\alpha_\xi-\alpha^c_\xi
\end{array}\right.
\end{equation*}
where $\alpha^c_\xi$ has been introduced in (\ref{eqVII.1.40}). Then we have
\begin{equation*}\left\{
\begin{array}{l}
Q\geq C (\v X\v^2+\v\theta\v\,\v Y\v^2)\,,\\
\v\partial ^A_{\alpha_x}F\v\leq C_A\enskip \textrm{ if}\enskip \v A\v\geq 1\,,\\
\v\partial ^{A_1}_{\alpha_x}\,\partial ^{A_2}_{\alpha_\xi}F\v\leq
C_{A_1A_2}\v\theta\v\enskip \textrm{if}\enskip \v A_2\v\geq 1\,,
\end{array}\right.
\end{equation*}
uniformly in $(\theta,x,y,\alpha)$ when $\v\theta\v\leq 1$, $x-\alpha_x\in\supp
\chi_6$, $y-x(\theta,\alpha)\in\supp \chi_6$, $\alpha_\xi\in\supp \psi_3$.
\end{lemma}

{\bf Proof } For $t\in (0,1)$ let us set
$m_t=(\theta,x,y,t\alpha+(1-t)\,\alpha^c)$. Then using Lemma (\ref{lVII.1.9})
we can write for $j=1,\ldots ,n$,
$$
\frac{\partial F}{\partial \alpha^j_x}\,(\theta,x,y,\alpha)=\sum^n_{k=1}\left
(\int^1_0 \left[\frac{\partial ^2F}{\partial \alpha^j_x\partial
\alpha^k_x}\,(m_t)(\alpha^k_x-x_k)+\frac{\partial ^2F}{\partial
\alpha^j_x\partial
\alpha^k_\xi}\,(m_t)(\alpha^k_\xi-(\alpha^c_\xi)^k)^k)\right]\,dt\right)
$$
$$
\frac{\partial F}{\partial \alpha^j_\xi}\,(\theta,x,y,\alpha)=\sum^n_{k=1}\left
(\int^1_0 \left[\frac{\partial ^2F}{\partial \alpha^j_\xi\partial
\alpha^k_x}\,(m_t)(\alpha^k_x-x_k)+\frac{\partial ^2F}{\partial
\alpha^j_\xi\partial
\alpha^k_\xi}\,(m_t)(\alpha^k_\xi-(\alpha^c_\xi)^k)\right]\,dt\right)\,.
$$
It follows from (\ref{eqVII.1.35}) and (\ref{eqVII.1.34}) that
$$
\frac{\partial F}{\partial
\alpha^j_x}=\frac{2i-2\theta}{1+2i\theta}\,X_j+\frac{2i\theta}{1+2i\theta}\,Y_j+
\CO[(\varepsilon+\delta)(\v X\v+\v\theta\v\v Y\v]
$$
$$
\frac{\partial F}{\partial
\alpha^j_\xi}=\frac{2i\theta}{1+2i\theta}\,X_j-\frac{2\theta}{1+2i\theta}\,Y_j+
\CO[(\varepsilon+\delta)\v\theta\v(\v X\v+\v Y\v)]\,.
$$
Therefore we have
$$
\theta^2\,\v X\v^2+\v X-\theta Y\v^2+\v\theta\v\,\v X\v^2+\v\theta\v\,\v
Y\v^2\leq \frac 52\,Q+C(\varepsilon+\delta)^2(\v X\v^2+\v\theta\v^2\,\v
Y\v^2)\,.
$$
Taking $\varepsilon+\delta$ small enough we obtain $Q\geq C(\v
X\v^2+\v\theta\v\,\v Y\v^2)$ which the first claim of the Lemma. The other
claims follows from (\ref{eqVII.1.35}) and the fact that $\v x-y\v\leq 2(
b+5\,\delta_2)\v\theta\v$ since we are in the case~3. \cqfd

We shall set in that follows,
\begin{equation}\label{eqVII.1.45}\left\{
\begin{array}{l}
L=\frac 1i\,\som^n_{j=1} \Big (\frac{\partial \overline F}{\partial
\alpha^j_x}\,(\theta,x,y,\alpha)\,\frac{\partial }{\partial \alpha^j_x}+\frac
1D\,\frac{\partial \overline F}{\partial
\alpha^j_\xi}\,(\theta,x,y,\alpha)\,\frac\partial {\partial \alpha^j_\xi}\Big)\\
Q=\Big(\Big\v\frac{\partial F}{\partial \alpha_x}\Big\v^2+\frac
1D\,\Big\v\frac{\partial F}{\partial
\alpha_\xi}\Big\v^2\Big)(\theta,x,y,\alpha)\\
\textrm{where}\enskip D=\v\theta\v\enskip \textrm{in the cases 1 and 3,}\enskip
D=\v x-y\v\enskip \textrm{in the case 2}\,.
\end{array}\right.
\end{equation}
Let us note that, according to Lemmas \ref{lVII.1.6}, \ref{lVII.1.7} and
\ref{lVII.1.10} $L$ is a vector field whose coefficients are uniformly bounded
together with their derivatives with respect to $\alpha$. Moreover we have
\begin{equation}\label{eqVII.1.46}
L\,e^{i\lambda F}=\lambda\,Q\,e^{i\lambda F}\,.
\end{equation}
Our first goal is to prove the following result.
\begin{proposition}\sl\label{pVII.1.11}
For any $N$ in $\N^n$ we can write
$$
L^N\,e^{i\lambda F(\theta,x,y,\alpha)}=\left (\lambda^N\,Q^N+\sum^{N-1}_{k=1}
h_{k,N}(\theta,x,y,\alpha)\lambda^k\right)\,e^{i\lambda F(\theta,x,y,\alpha)}
$$
where $h_{k,N}$ are smooth functions satisfying
$$
\v L^j\,h_{k,N}(\theta,x,y,\alpha)\v\leq C_{j,N}\,Q^k\,,\quad 0\leq k\leq
N-1\,,\enskip j\in\N\,,
$$
uniformly when $0<\v\theta\v\leq 1$, $\v x-\alpha_x\v\leq 1$, $\v
y-x(\theta,\alpha)\v\leq 2\delta$, $\v\alpha_\xi\v\leq 2$.
\end{proposition}

{\bf Proof }

{\bf Step 1 } Let $\CH$ be the set of $C^\infty $ functions
$c=c(\theta,x,y,\alpha)$ such that for any $\gamma\in\N^{2n}$, $\partial
^\gamma_\alpha c$ is uniformly bounded when $0<\v\theta\v\leq 1$,
 $\v x-\alpha_x\v\leq 1$, $\v y-x(\theta,\alpha)\v\leq 2\delta$,
$\v\alpha_\xi\v\leq 2$. For instance $\frac{\partial F}{\partial \alpha^j_x}$
and $\frac 1D\,\frac{\partial F}{\partial \alpha^j_\xi}$ belong to $\CH$. 

Let us set $T=\frac{\partial F}{\partial \alpha_x}$, $\overline
T=\frac{\partial \overline F}{\partial \alpha_x}$, $S=\frac 1{\sqrt
D}\,\frac{\partial F}{\partial \alpha_\xi}$, $\overline S=\frac 1{\sqrt
D}\,\frac{\partial \overline F}{\partial \alpha_\xi}$.

Let $\CP$ be the set of homogeneous polynomials of order $2$ in $T,\overline T$,
$S,\overline S$ with coefficients in $\CH$. For instance we have $Q\in\CP$.
We claim that $L$ sends $\CP$ into $\CP$. First of all if $P\in\CP$ then
$\frac{\partial \overline F}{\partial \alpha^j_x}\,\frac{\partial P}{\partial
\alpha^j_x}\in\CP$ since $\frac{\partial \overline F}{\partial
\alpha^j_x}\in\CH$, $\frac{\partial ^2F}{\partial \alpha^2}\in\CH$ and
$\frac{\partial \overline F}{\partial \alpha^j_x}=\overline T_j$. On the other
hand $\frac 1 D\,\frac{\partial \overline F}{\partial
\alpha^j_\xi}\,\frac{\partial P}{\partial \alpha^j_\xi}\in
\CP$ since $\frac 1 D\,\frac{\partial \overline F}{\partial
\alpha^j_\xi}\in\CH$,
$\frac 1{\sqrt D}\,\frac{\partial \overline F}{\partial \alpha^j_\xi}=\overline
T_j$ and $\frac 1D\,\frac{\partial ^2F}{\partial \alpha_\xi\partial
\alpha}\in\CH$. It follows that $L$ maps $\CP$ into $\CP$.

 Now if $P\in\CP$ then $\v P\v\leq
C\,Q$. It follows that for all $j\in\N$,
\begin{equation}\label{eqVII.1.47}
\v L^jQ\v\leq C_jQ
\end{equation}
uniformly in $(\theta,x,y,\alpha)$.

{\bf Step 2 } We claim that
\begin{equation}\label{eqVII.1.48}\left\{
\begin{array}{l}
\textrm{for all}\enskip N\in\N\,,\enskip j\in\N\enskip \textrm{we have}\enskip
\v L^j\,Q^N\v\leq C_{jN}\,Q^N\enskip \textrm{uniformly when}\\
0<\v\theta\v<1\,,\enskip \v x-\alpha_x\v\leq 1\,,
\enskip \v y-x(\theta,\alpha)\v\leq 2\delta\,,\enskip \v\alpha_\xi\v\leq 2\,.
\end{array}\right.
\end{equation}
Indeed by the Faa di Bruno formula (since $L$ is a homogeneous vector field)
$L^j\,Q^N$ is a finite linear combination of terms of the form
$$
Q^{N-M} \prod^s_{i=1} (L^{\ell _i}Q)^{k_i}
$$
where $1\leq M\leq N$, $1\leq s\leq M$, $\som^s_{i=1} k_i=M$, $\som^s_{i=1}
k_i\,\ell _i=j$. Each of such terms is bounded, according to
(\ref{eqVII.1.47}) by $C\,Q^{N-M}\,Q^{\som^s_i k_i}=C\,Q^{N-M}\,Q^M=C\,Q^N$
which proves our claim.

{\bf Step 3 } Proof of Proposition \ref{pVII.1.11}

We use an induction on $N$. For $N=1$ the result follows from
(\ref{eqVII.1.46}). Let us assume it is true up to the order $N$. Then
\begin{equation*}
\begin{split}
L^{N+1}(e^{i\lambda F})=e^{i\lambda F}\left
(\lambda^N\,N\,Q^{N-1}\,LQ+\sum^{N-1}_{k=0}
(L\,h_{k,N})\,\lambda^k+\lambda^{N+1}\,Q^{N+1}
+\sum^{N-1}_{k=0} Q\,h_{k,N}\,\lambda^{k+1}\right )\\
=e^{i\lambda F}\left (\lambda^{N+1}\,Q^{N+1}+\sum^N_{k=0}
h_{k,N+1}\,\lambda^k\right )
\end{split}
\end{equation*}
where
\begin{equation*}\left\{
\begin{array}{l}
h_{N,N+1}=N\,Q^{N-1}\,L\,Q+Q\,h_{N-1,N}\,,\\
h_{k,N+1}=L\,h_{k,N}+Q\,h_{k-1,N}\,,\enskip 1\leq k\leq N-1\,,\\
h_{0,N+1}=L\,h_{0,N}\,.
\end{array}\right.
\end{equation*}
Now we have
$$
\v L^j\,h_{N,N+1}\v\leq N\sum^j_{i=0}
\begin{pmatrix}j\\i\end{pmatrix}\v L^i\,Q^{N-1}\v\v L^{j-i}Q\v+\sum^j_{i=0}
\begin{pmatrix}j\\i\end{pmatrix}\v L^iQ\v\,\v L^{j-i}\,h_{N-1,,N}\v\,.
$$
Using (\ref{eqVII.1.47}), (\ref{eqVII.1.48}) and the induction hypothesis we
deduce that $\v L^j\,h_{N,N+1}\v\leq C_{j,N}\,Q^N$. The estimates for the other
terms are completely analogous. \cqfd

We need another Lemma.
\begin{lemma}\sl\label{lVII.1.12}
\begin{itemize}
\item[(i)] For any $N\in\N^*$ we can write
$$
L^N(e^{i\lambda F})=G_N(\theta,x,y,\alpha,\lambda)\,e^{i\lambda F}\,.
$$
\item[(ii)] There exists a constant $K_N>0$ such that
$$
\v K_N+G_N(\theta,x,y,\alpha,\lambda)\v\geq \frac 12 (\lambda^N\,Q^N+1)\,.
$$
\item[(iii)] For any $j\in\N$, $N\in\N^*$ we can find a constant $C_{j,N}>0$
such that
$$
\v L^j\,G_N(\theta,x,y,\alpha,\lambda)\v\leq C_{j,N}\,\v
G_N(\theta,x,y,\alpha,\lambda)+K_N\v\,.
$$
\item[(iv)] For any $j\in\N$, $N\in\N^*$ we can find a constant $C'_{jN}>0$
such that
$$
\Big\v(^tL)^j\Big(\frac 1{K_N+G_N}\Big)\Big\v\leq \frac{C'}{\v K_N+G_N\v}
$$
where $^tL$ denotes the transposed of $L$.
\end{itemize}
\end{lemma}

{\bf Proof } (i) is a consequence of Proposition \ref{pVII.1.11} with
$$
G_N(\theta,x,y,\alpha,\lambda)=\lambda^N\,Q^N+\sum^{N-1}_{k=0}
h_{k,N}(\theta,x,y,\alpha)\,\lambda^k
$$
and $\v h_{k,N}\v\leq C_N\,Q^k$, $0\leq k\leq N-1$, $C_N\geq 1$.

(ii) If $\lambda\,Q\leq 2N\,C_N$ then $\som^{N-1}_{k=0} \v
h_{k,N}\v\,\lambda^k\leq C_N \som^{N-1}_{k=0} (\lambda\,Q)^k\leq C_N
\som^{N-1}_{k=0} (2N\,C_N)^k$. If $\lambda\,Q\geq 2N\,C_N$ then
$(\lambda\,Q)^N=(\lambda\,Q)^{N-k}(\lambda\,Q)^k\geq 2N\,C_N(\lambda\,Q)^k$ if
$0\leq k\leq N-1$. Then
$$
\sum^{N-1}_{k=0} \v h_{k,N}\v\,\lambda^k\leq \frac{C_N}{2N\,C_N} \left
(\sum^{N-1}_{k=0}\,1\right )(\lambda\,Q)^N=\frac 12 (\lambda\,Q)^N\,.
$$
Therefore
$$
\sum^{N-1}_{k=0} \v h_{k,N}\v\,\lambda^k\leq \frac 12 (\lambda\,Q)^N+K'_N\,.
$$
This implies
$$
\Big\v K'_N+\frac 12+G_N\Big\v\geq K'_N+\frac 12 +\lambda^N\,Q^N-\frac
12\,\lambda^N\,Q^N-K'_N\geq \frac 12 (\lambda^N\,Q^N+1)\,.
$$

(iii) We have, by (\ref{eqVII.1.48}) and Proposition \ref{pVII.1.11}
\begin{equation*}
\begin{aligned}
\v L^j\,G_N\v&=\Big\v\lambda^N\,L^j\,Q^N+\sum^{N-1}_{k=0}
L^j\,h_{k,N}\,\lambda^k\Big\v\leq
C_{j,N}\,\lambda^N\,Q^N+\sum^{N-1}_{k=0}C'_{jN}\,\lambda^k\,Q^k\\
&\leq C''_{jN}(1+\lambda^N\,Q^N)\leq \tilde C_{j,N}\,\v K_N+G_N\v\,.
\end{aligned}
\end{equation*}

(iv) Let us recall that $L=\frac{\partial \overline F}{\partial \alpha_x}\cdot
\frac{\partial }{\partial \alpha_x}+\frac 1D\,\frac{\partial \overline
F}{\partial \alpha_\xi}\cdot \frac\partial {\partial \alpha_\xi}$.
According to the Lemmas \ref{lVII.1.6}, \ref{lVII.1.7} and \ref{lVII.1.10}, we
see that $^tL=-L+c(\theta,x,y,\alpha)$ where
$$
\v\partial ^\gamma_\alpha\,c(\theta,x,y,\alpha)\v\leq C_\gamma
$$
uniformly in $(\theta,x,y,\alpha)$. We claim that for $j\in\N$
\begin{equation}\label{eqVII.1.49}\left\{
\begin{array}{l}
(^tL)^j=\som^j_{k=0} c_{j,k}(\theta,x,y,\alpha)\,L^k\,,\enskip \textrm{with}\\
\v\partial ^\gamma_\alpha\,c_{jk}(\theta,x,y,\alpha)\v\leq c_{j,k,\gamma}\,.
\end{array}\right.
\end{equation}
We saw above that this is true when $j=1$. Let us assume that
(\ref{eqVII.1.49}) is true up to the order $j$. Then
$$
(^tL)^{j+1}=(-L+c) \sum^j_{k=0} c_{jk}(\theta,x,y,\alpha)\,L^k=\sum^{j+1}_{k=0}
c_{j+1,k}(\theta,x,y,\alpha) \,L^k
$$
where
\begin{equation*}\left\{
\begin{array}{l}
c_{j+1,0}=c\,c_{j,0}-L\,c_{j,0}\\
c_{j+1,k}=c\,c_{jk}-L\,c_{jk}-c_{j,k-1}\,,\enskip 1\leq k\leq j\\
c_{j+1,j+1}=-c_{j,k}\,.
\end{array}\right.
\end{equation*}
Then the estimate on $\v\partial ^\gamma_\alpha\,c_{j+1,k}\v$ follows from the
induction and the fact that the coefficients of $L$ have all their derivatives
with respect to $\alpha$ bounded uniformly.
It follows that
\begin{equation}\label{eqVII.1.50}
\Big\v(^tL)^j\Big (\frac 1{K_N+G_N}\Big )\Big\v\leq C_{j,N} \sum^j_{k=0}
\Big\v L^k\Big (\frac 1{K_N+G_N}\Big)\Big\v\,.
\end{equation}
Now by the Faa di Bruno formula, $L^k\big (\frac 1{K_N+G_N}\big )$ is a finite
linear combination of terms of the form
$$
\frac 1{(K_N+G_N)^{1+\beta}} \prod^s_{i=1} (L^{\ell _i}\,G_N)^{k_i}
$$
where $1\leq \beta\leq k$, $1\leq s\leq k$, $\som^s_{i=1} k_i=\beta$,
$\som^s_{i=1}k_i\,\ell _i=k$. It follows from (iii) that
$$
\Big\v L^k\Big(\frac 1{K_N+G_N}\Big)\Big\v\leq C_{j,N}\,\frac 1{\v
K_N+G_N\v^{1+\beta}}\,\v K_N+G_N\v^{\som^s_{i=1}k_i}\leq \frac{C_{j,N}}{\v
K_N+G_N\v}\,.
$$
Then (iv) follows from (\ref{eqVII.1.50}).

We can now state the estimate on the kernel for $\v\theta\v\leq 1$.
\begin{lemma}\sl\label{lVII.1.13}
Let $\tilde k_+=\tilde k_+(t,x,y,\lambda)$ be the kernel defined by
(\ref{eqVII.1.32}). Then one can find a positive constant $C$ such that
$$
\v\tilde k_+(t,x,y,\lambda)\v\leq \frac C{\v t\v^{n/2}}
$$
for all $(t,x,y,\lambda)$ such that $\lambda\geq 1$, $\v\lambda t\v\leq 1$,
$x\in\R^n$, $y\in\R^n$.
\end{lemma}

{\bf Proof } The kernel $\tilde k_+$ can be written as
$$
\tilde k_+(t,x,y,\lambda)=\int e^{i\lambda F(-\lambda t,x,y,\alpha)}\,\tilde
a(\lambda t,x,y,\alpha)\,d\alpha
$$
where
$$
\tilde a(\lambda t,x,y,\alpha)=a_1(-\lambda
t,y,\alpha)\,\chi_6(x-\alpha_x)\,\psi_3(\alpha_\xi)\,\chi_5\Big
(\frac{y-x(-\lambda t,\alpha)}{\bra\lambda t\ket}\Big )\,.
$$
Let us note that all the derivatives of $\tilde a$ with respect to $\alpha$ are
uniformly bounded.

We use Lemma \ref{lVII.1.12} to write
\begin{equation*}
\begin{aligned}
\tilde k_+(t,x,y,\lambda)&=\int (K_N+L^N)\,e^{i\lambda F(-\lambda
t,x,y,\alpha)}\,\frac{\tilde a(\lambda t,x,y,\alpha)}{K_N+G_N(-\lambda
t,x,y,\alpha,\lambda)}\,d\alpha\\
&=\int e^{i\lambda F(\cdots)}(K_N+(^tL)^N)\left[\frac{\tilde
a(\cdots)}{K_+G_N(\cdots)}\right]\,d\alpha\,.
\end{aligned}
\end{equation*}
Now using (ii), (iv) of Lemma \ref{lVII.1.12}, (\ref{eqVII.1.49}) and the fact
that $\Im F\geq 0$ we obtain
\begin{equation}\label{eqVII.1.51}
\v\tilde k_+(t,x,y,\lambda)\v\leq C_N\,\lambda^{\frac{3n}2} \int \frac{\v\tilde
\psi(\alpha_\xi)\v}{1+\lambda^N\,Q^N}\,d\alpha=C_N\,I
\end{equation}
where $Q$ has been introduced in (\ref{eqVII.1.45}) and $\tilde
\psi(\alpha_\xi)$ is a smooth function with compact support. Now according to
the Lemmas \ref{lVII.1.6}, \ref{lVII.1.7}  and \ref{lVII.1.10} we have
\begin{equation*}\left\{
\begin{array}{ll}
Q\geq C(\v x-\alpha_x\v^2+\v\lambda t\v)\enskip &\textrm{(case 1 and 2)}\\
Q\geq C(\v x-\alpha_x\v^2+\v\lambda t\v\,\v\alpha_\xi-\alpha^c_\xi\v^2)\enskip
&\textrm{(case 3)}\,.
\end{array}\right.
\end{equation*}
Let us fix the integer $N$ such that $N>n$. Since $\tilde \psi$ has
compact support we have in case 1 and 2
$$
I\leq \lambda^{\frac{3n}2} \int_{\R^n} \,\frac{dX}{1+\lambda^N\,\v
X\v^{2N}+\lambda^N\,\v\lambda t\v^N}\,.
$$
Let us set $X=\frac{(1+\lambda^N\,\v\lambda t\v^N)^{\frac
1{2N}}}{\sqrt\lambda}\,z$. Then
$$
I\leq \lambda^{\frac{3n}2}\,\frac{(1+\lambda^N\,\v\lambda t\v^N)^{\frac
n{2N}}}{\lambda^{\frac n2}}\,\frac 1{1+\lambda^N\,\v\lambda t\v^N} \int_{\R^n}
\frac{dz}{1+\v z\v^{2N}}\,.
$$
If $\lambda\,\v\lambda t\v\leq 1$, that is $\lambda^2\leq \frac 1{\v t\v}$, we
have,
$$
I\leq C_N\,\lambda^n\leq \frac{C_N}{\v t\v^{\frac n2}}\,.
$$
If $\lambda\,\v\lambda t\v\geq 1$, that is $\frac 1{\lambda^2}\leq \v t\v$, we
can write for $N>n$,
$$
I\leq C_N\,\lambda^n\,\frac{(\lambda\,\v\lambda t\v)^{\frac
n2}}{(\lambda\,\v\lambda t\v)^N}\leq \frac{C_N}{\lambda^{2(N-n)}}\,\v
t\v^{\frac n2-N}\leq C_N\,\v t\v^{N-n+\frac n2-N}\leq \frac{C_N}{\v t\v^{\frac
n2}}\,.
$$
It follows from (\ref{eqVII.1.51}) that $\v\tilde k_+(t,x,y,\lambda)\v\leq
\frac C{\v t\v^{\frac n2}}$.

In case 3 we have
$$
I\leq C\,\lambda^{\frac{3n}2}\int_{\R^{2n}}\,\frac{d\alpha}{1+(\lambda\,\v
x-\alpha_x\v^2+\lambda\,\v\lambda t\v\,\v\alpha_\xi-\alpha^c_\xi\v^2)^N}\,.
$$
Setting $X=\sqrt\lambda\,(x-\alpha_x)$, $Y=\sqrt\lambda\,\sqrt{\v\lambda
t\v}\,(\alpha_\xi-\alpha^c_\xi)$ we obtain
$$
I\leq C\,\lambda^{\frac{3n}2}\,\lambda^{-\frac n2}\,\lambda^{-\frac
n2}\,\v\lambda t\v^{-\frac n2} \int_{\R^{2n}}\,\frac{dX\,dY}{1+(\v X\v^2+\v
Y\v^2)^N}\,.
$$
It follows that
$$
\v\tilde k(t,x,y,\lambda)\v\leq \frac C{\v t\v^{\frac n2}}\,.
$$
This completes the proof of Lemma \ref{lVII.1.13} thus the proof of Theorem
\ref{tVII.1.1}. \cqfd

\subsection{\ End of the proof of Theorem \ref{tII.2.1}}\label{ssVII.2}

Let us recall that we have set in (\ref{eqVI.2.7}),
\begin{equation}\label{eqVII.2.1}
U_\pm(t)=\chi^+_1\,\psi_2 \Big(\frac D\lambda\Big)\,e^{-itP}\,.
\end{equation}
It follows from Theorem \ref{tVII.1.1} that there exists $C\geq 0$ such that
$$
\Vert U_\pm(t_1)\,U_\pm(t_2)^*f\Vert _{L^\infty (\R^n)}\leq \frac C{\v
t_1-t_2\v^{\frac n2}}\,\Vert f\Vert _{L^1(\R^n)}
$$
for all $t_1,t_2$ in $[-T,T]$ with $t_1\not =t_2$ and all $f\in L^1(\R^n)$. 
Moreover there is a conservation of the $L^2$ norm for $e^{-itP}$ which implies,
$$
\Vert U_\pm(t)\,u_0\Vert _{L^2(\R^n)}\leq C\,\Vert u_0\Vert _{L^2(\R^n)}
$$
for all $\v t\v\leq T$ and $u_0\in L^2(\R^n)$.

We can therefore apply Lemma \ref{lII.1.3} with $X=\R^n$, $H=L^2(\R^n)$ and we
obtain
\begin{equation}\label{eqVII.2.2}
\Vert U_\pm(\cdot )\,u_0\Vert _{L^q([-T,T],L^r(\R^n))}\leq C\,\Vert u_0\Vert
_{L^2(\R^n)}
\end{equation}
for all $u_0\in L^2(\R^n)$, where $q\geq 2$ and $\frac 2q=\frac n2-\frac nr$.
Now it follows from (\ref{eqVI.2.2}) that we have
$$
\chi^+_1(x)+\chi^-_1(x)\geq 1\enskip \textrm{for all}\enskip x\in\R^n\,.
$$
Then (\ref{eqVII.2.1}) and (\ref{eqVII.2.2}) show that
\begin{equation}\label{eqVII.2.3}
\Big\Vert \psi_2 \Big (\frac D\lambda\Big)\,e^{-itP}\,u_0\Big\Vert
_{L^q(I,L^r(\R^n))}\leq C\,\Vert u_0\Vert _{L^2(\R^n)}
\end{equation}
for all $u_0\in L^2(\R^n)$, where $I=[-T,T]$.

Now let us recall (see (\ref{eqVI.2.6})) that
$\psi_2(\xi)=\psi_0(\xi)\,\psi_1(\xi)$ where
\begin{equation*}\left\{
\begin{array}{l}
\psi_0(\xi)=1\quad  \textrm{if}\quad  \Big\v\frac\xi{\v\xi\v}-\xi_0\Big\v\leq
\delta_2\,,\enskip \v\xi\v\geq 2\,\delta_2\,,\enskip \v\xi_0\v=1\,,\\
\supp \psi_0\subset \Big\{\xi :\frac\xi{\v\xi\v}-\xi_0\Big\v\leq
2\,\delta_2\Big\}\quad 
\textrm{and}\quad  \v\xi\v\geq \delta_2
\end{array}\right.
\end{equation*}
\begin{equation*}\left\{
\begin{array}{l}
\psi_1(\xi)=1\quad  \textrm{if}\enskip a-\delta_2\leq \v\xi\v\leq
b+\delta_2\,,\quad  a=\frac 6{10}\,,\enskip b=\frac{19}{10}\,,\\
\supp \psi_1\subset \{\xi : a-2\,\delta_2\leq \v\xi\v\leq b+2\,\delta_2\}\,.
\end{array}\right.
\end{equation*}
By a finite partition of unity we deduce easily from (\ref{eqVII.2.3}) that
\begin{equation}\label{eqVII.2.4}
\Big\Vert \psi_1\Big(\frac D\lambda\Big)\,e^{-itP}\,u_0\Big\Vert
_{L^q(I,L^r(\R^n))}\leq C\,\Vert u_0\Vert _{L^2(\R^n)}\,.
\end{equation}
Now let us recall that $p(x,\xi)=\v\xi\v^2+\varepsilon\,\som^n_{j,k=1}
b_{jk}(x)\,\xi_j\,\xi_k$. Therefore if $\varepsilon$ is small enough we have
\begin{equation}\label{eqVII.2.5}
\frac 9{10}\,\v\xi\v^2\leq p(x,\xi)\leq \frac{11}{10}\,\v\xi\v^2\,.
\end{equation} 
Let $\varphi _0(t)\in C^\infty _0(\R)$ be such that
$$
\varphi _0(t)=1\enskip \textrm{if}\enskip \v t\v\leq \frac 9 5\,\enskip
\varphi (t)=0\enskip \textrm{if}\enskip \v t\v\geq 3
$$
and set
\begin{equation}\label{eqVII.2.6}
\varphi (t)=\varphi _0(t)-\varphi _0(4t)\,.
\end{equation}
It follows that $\supp \varphi \subset\{t\in\R : \frac 9{20}\leq \v t\v\leq
3\}$.

Let $\tilde \varphi \in C^\infty (\R)$ be such that
\begin{equation}\label{eqVII.2.7}\left\{
\begin{array}{l}
\supp \tilde \varphi \subset\Big\{t\in\R : \frac 4{10}\leq \v t\v\leq
\frac{31}{10}\Big\}\\
\tilde \varphi =1\enskip \textrm{on a neighborhood of}\enskip \supp \varphi
\,,\enskip \textrm{so}\\
\varphi (t)\,\tilde \varphi (t)=\varphi (t)\,.
\end{array}\right.
\end{equation}
We claim that for every $(x,\xi)\in T^*\R^n$,
\begin{equation}\label{eqVII.2.8}
(1-\psi_1(\xi))\,\tilde \varphi (p(x,\xi))\equiv 0\,.
\end{equation}
Indeed on the support of $1-\psi_1(\xi)$ we have
$$
\textrm{(i)}\qquad \v\xi\v\leq
\frac{6}{10}-\delta_2\qquad\qquad\textrm{or}\qquad\qquad(ii)\qquad\v\xi\v\geq
\frac{19}{10}+\delta_2\,.
$$
In the case (i), (\ref{eqVII.2.5}) shows that
$$
0\leq p(x,\xi)\leq \frac{11}{10}\,\v\xi\v^2\leq \frac{11}{10}\Big (\frac
6{10}\Big)^2<\frac 4{10}\,.
$$
In the case (ii) we have
$$
p(x,\xi)\geq \frac 9{10} \Big (\frac{19}{10}\Big)^2>\frac{31}{10}\,.
$$
Thus $\tilde \varphi (p(x,\xi))=0$ by (\ref{eqVII.2.7}).

Now, with $\varphi $ introduced in (\ref{eqVII.2.6}) we claim that
\begin{equation}\label{eqVII.2.9}
(1)=\Big\Vert \varphi \Big(\frac P{\lambda^2}\Big)\,e^{-itP}\,u_0\Big\Vert
_{L^q(I,L^r(\R^n))}\leq C\,\Big\Vert \varphi \Big(\frac P{\lambda^2}\Big
)\,u_0\Big\Vert _{L^2(\R^n)}
\end{equation}
for all $u_0\in L^2(\R^n)$.

Indeed we can write
$$
(1)\leq \Big\Vert \psi_1\underbrace{\Big(\frac D\lambda\Big)\,\varphi
\Big (\frac P{\lambda^2}\Big )}_{(2)}\,e^{-itP}\,u_0\Big\Vert _{L^q(I,L^r)}+
\Big\Vert \Big( I-\psi_1\underbrace{\Big(\frac D\lambda\Big)\Big)\,\varphi
\Big (\frac P{\lambda^2}\Big )}_{(3)}\,e^{-itP}\,u_0\Big\Vert _{L^q(I,L^r)}\,.
$$
Since $\varphi \big (\frac P{\lambda^2}\big)$ commutes with $e^{-itP}$, we
deduce from (\ref{eqVII.2.4}) that 
\begin{equation}\label{eqVII.2.10}
(2)\leq \Big\Vert \varphi \Big (\frac P{\lambda^2}\Big)\,u_0\Big\Vert
_{L^2(\R^n)}\,.
\end{equation}
Using (\ref{eqVII.2.7}) we can write
$$
(3)=\Big\Vert \Big(I-\psi_1\Big(\frac D\lambda\Big)\Big)\,\tilde \varphi
\Big(\frac P{\lambda^2}\Big)\,e^{-itP}\,\varphi \Big(\frac
P{\lambda^2}\Big)\,u_0\Big\Vert _{L^q(I,L^r)}\,.
$$
Now, since $r\geq 2$, there exists, by the Sobolev embedding, $s\geq 0$ such
that $H^s(\R^n)\hookrightarrow L^r(\R^n)$. We fix such an $s$. Then
$$
(3)\leq \Big\Vert \Big(I-\psi_1\Big(\frac D\lambda\Big)\Big)\,\tilde \varphi
\Big(\frac P{\lambda^2}\Big)\,e^{-itP}\,\varphi \Big(\frac
P{\lambda^2}\Big)\,u_0\Big\Vert _{L^\infty (I,H^s(\R^n)}\,.
$$
Now we use (\ref{eqVII.2.8}) and Proposition A.1 in  \cite{BGT}. It
follows that
$$
(3)\leq C' \,\Big\Vert e^{-itP}\,\varphi \Big(\frac
P{\lambda^2}\Big)\,u_0\Big\Vert _{L^\infty (I,L^2(\R^n))}\leq 
C'' \,\Big\Vert \varphi \Big(\frac
P{\lambda^2}\Big)\,u_0\Big\Vert _{L^2(\R^n)}\,,
$$
which, together with (\ref{eqVII.2.10}) proves our claim (\ref{eqVII.2.9}).
By (\ref{eqVII.2.6}) we have
$$
\varphi _0(t)+\sum^{+\infty }_{k=1} \varphi (2^{-2k}\,t)=1\,.
$$
Then using Corollary 2.3 in  \cite{BGT} we can write
$$
\Vert e^{-itP}\,u_0\Vert _{L^q(I,L^r(\R^n))}\leq C\,\Vert u_0\Vert _{L^2(\R^n)}+
\bigg\Vert \Big (\sum^{+\infty }_{k=1} \big\Vert e^{-itP}\,\varphi
(2^{-2k}P)\,u_0\big\Vert ^2_{L^r(\R^n)}\Big)^{\frac 12}\bigg\Vert _{L^q(I)}\,.
$$
By the Minkowski inequality and (\ref{eqVII.2.9}) we obtain
$$
\Vert e^{-itP}\,u_0\Vert _{L^q(I,L^r(\R^n))}\leq C\,\Vert u_0\Vert _{L^2(\R^n)}+
 \left (\sum^{+\infty }_{k=1} \big\Vert \varphi (2^{-2k}P)\,u_0\big\Vert
^2_{L^2(\R^n)}\right)^{\frac 12}\,.
$$
This implies that
$$
\big\Vert e^{-itP}\,u_0\big\Vert _{L^q(I,L^r(\R^n))}\leq C'\,\Vert u_0\Vert
_{L^2}\,,
$$
which is the estimate claimed in Theorem \ref{tII.2.1}.

\section{\quad Appendix}\label{sVIII}

\subsection{\quad The Faa di Bruno Formula}\label{ssVIII.1}

We shall make a repeat use of the following formula which can be found in even
a more precise form in the paper of Constantin and Savits  \cite{CS}.

Let $m\in\N^*$, $F\in C^m(\R^N,\C)$, $U_k\in C^m(\R^p,\R)$, $k=1,\ldots ,N$.
Then for $\v\Lambda\v\leq m$ we have
$$
\partial ^\Lambda_Y[F(U(Y))]=\sum^N_{k=1}\,\frac{\partial F}{\partial
U_k}\,(U(Y))\,\partial^\Lambda _Y\,U_k(Y)+(1)
$$
where $(1)$ is a finite linear combination of terms of the form
$$
(\partial ^\beta_UF)(U(Y)) \prod^s_{j=1} \big (\partial
^{L_j}_Y\,U(Y)\big)^{K_j}
$$
where $2\leq \v\beta\v\leq \v\Lambda\v$, $1\leq s\leq \v\Lambda\v$, $\v
K_j\v\geq 1$, $\v L_j\v\geq 1$ and
$$
\sum^s_{j=1} K_i=\beta\,,\quad \sum^s_{j=1} \v K_j\v\,L_j=\Lambda\,.
$$
The precise coefficients of this sum can be found in \cite{CS}.

\subsection{\quad Proof of Proposition \ref{pIII.2.1}}\label{ssVIII.2}

The system satisfied by $r$ and $\zeta$ is the following.
\begin{equation*}\left\{
\begin{array}{ll}
\dot r_j(t)=2\zeta_j(t)+2\varepsilon \som^n_{k=1}
b_{jk}(x(t))\,\xi_k(t)\,,&  r_j(0)=0\\
\dot \zeta_j(t)=-\varepsilon \som^n_{p,q=1} \frac{\partial b_{pq}}{\partial x_j}
(x(t))\,\xi_p(t)\,\xi_q(t)\,,&  \zeta_j(0)=0\,.
\end{array}\right.
\end{equation*}
We prove our claim by induction on $\v A\v+\v B\v=k$. Moreover it is clear that
(ii) implies (i) since $r(0)=\zeta(0)=0$. Setting $\Phi(t)=\v\dot
r(t)\v+\v\dot\zeta(t)\v$ the above equations show that $\Phi(t)\leq
C_0\varepsilon+C_1\int^t_0 \Phi(\sigma)\,d\sigma$. Thus by Gronwall inequality,
$\Phi(t)\leq C_2(T)\,\varepsilon$. This shows that (i) and (ii) are true for
$k=0$. Let us assume they are true up to the order $k-1$. Let us set
$X=(x,\xi)$ and if $\alpha=(A,B)\in\N^n\times\N^n$, $\partial
^\alpha_X=\partial ^A_x\,\partial ^B_\xi$. It follows from the induction that
\begin{equation*}\left\{
\begin{array}{l}
\v\xi(t,x,\xi)\v\leq C(T)\\
\v\partial ^\alpha_X\,x(t,x,\xi)\v+\v\partial ^\alpha_X\,\xi(t,x,\xi)\v\leq
C_\alpha(T)
\end{array}\right.
\end{equation*}
if $1\leq \v\alpha\v\leq k-1$ and $k\geq 2$.

It follows that if $1\leq \v\alpha\v\leq k-1$ we have $\v\partial
^\alpha_X[F(x(t))]\v\leq C_\alpha(T)$ if $F=b_{jk}$ or $\frac{\partial
b_{pq}}{\partial x_j}$.

Let us now take $\v\alpha\v=k\geq 1$ and let us set $\Phi(t)=\v\partial
^\alpha_X\,\dot r(t)\v+\v\partial ^\alpha_X\,\dot\zeta(t)\v$. Then using the
Leibniz and Faa di Bruno formulas (see Section \ref{ssVIII.1}) and
differentiating the above differential system we find
$$
\Phi(t)\leq C_\alpha\,\varepsilon+C'_\alpha \int^t_0 \Phi(\sigma)\,d\sigma
$$
after using the induction. We conclude again by the Gronwall inequality.
\cqfd

\subsection{\quad Proof of Proposition \ref{pIII.3.2}}\label{ssVIII.3}

The system satisfied by $z$ and $\zeta$ is the following.
\begin{equation}\label{eqVIII.3.1}\left\{
\begin{array}{l}
\dot z_j(t)=2\varepsilon\som^n_{k=1} b_{jk}(x(t))\,\xi_k(t)-2\varepsilon \,t
\som^n_{p,q=1} \frac{\partial b_{pq}}{\partial
x_j}\,(x(t))\,\xi_p(t)\,\xi_q(t)\\
\dot\zeta_j(t)=-\varepsilon \som^n_{p,q=1} \frac{\partial b_{pq}}{\partial
x_j}\,(x(t))\,\xi_p(t)\,\xi_q(t)\\
z_j(0)=\zeta_j(0)=0\,.
\end{array}\right.
\end{equation}
By Proposition \ref{pIII.3.1} we have $\bra x(t,x,\xi)\ket\geq C (1+\v x\v+t)$.
This implies that for $\ell \geq 1$,
\begin{equation}\label{eqVIII.3.2}
\int^t_0 \frac{ds}{\bra x(s)\ket^{\ell +\sigma_0}}\leq \frac{C_\ell }{\bra
x\ket^{\ell -1+\sigma_0}}\,.
\end{equation}
We proceed by induction on $k$. Let us begin by the case $k=0$. We deduce from  
(\ref{eqVIII.3.1}) that
\begin{equation*}
\begin{aligned}
\v\dot z_j(t)\v&\leq C(A_0,A_1)\,\varepsilon \Big (\frac 1{\bra
x(t)\ket^{1+\sigma_0}}+\frac t{\bra x(t)\ket^{2+\sigma_0}}\Big)\leq
\frac{C'(A_0,A_1)\,\varepsilon}{\bra x(t)\ket^{1+\sigma_0}}\\
\v\dot \zeta_j(t)\v&\leq \frac{C(A_1)\,\varepsilon}{\bra x(t)\ket^{2+\sigma_0}}
\end{aligned}
\end{equation*}
since $\v\xi(t,x,\xi)\v\leq 2\,\v\xi\v\leq 4$.
Then the estimates in Proposition \ref{pIII.3.2} when $k=0$ follow from
(\ref{eqVIII.3.2}). Assume now that these estimates are true when $\v A\v+\v
B\v\leq k-1$ and let us deduce several facts.

For $\ell \in\N$ let us introduce the space,
$$
\BB^\ell _{\sigma_0}=\Big\{F\in C^\infty (\R^n) : \v\partial ^\beta F(y)\v\leq
\frac{C_\beta}{\bra y\ket^{\ell +\v\beta\v+\sigma_0}}\,,\enskip \textrm{for all}
\enskip y\in\R^n\Big\}\,.
$$
Let us set $X=(x,\xi)$ and if $\alpha=(A,B)\in\N^n\times\N^n$, $\partial
^\alpha_X=\partial ^A_x\,\partial ^B_\xi$. Then a staightforward
computation shows that for $\alpha\in\N^n\times\N^n$,
\begin{equation}\label{eqVIII.3.3}
\v\partial ^\alpha_X(\xi_j)\v\leq \frac 1{\bra x\ket^{\v A\v}}\,,\quad
\v\partial ^\alpha_X(x_j)\v\leq \frac 1{\bra x\ket^{\v A\v-1}}\,.
\end{equation}
It follows that for $\v\alpha\v\leq k-1$,
\begin{equation}\label{eqVIII.3.4}
\v\partial ^\alpha_X\,\xi_j(t,x,\xi)\v\leq \frac{2M_{k-1}}{\bra x\ket^{\v
A\v}}\,.
\end{equation}
Indeed
$$
\v\partial ^\alpha_X(\xi_j(t,x,\xi)\v\leq \v\partial
^\alpha_X(\xi_j)\v+\v\partial ^\alpha_X\,\zeta_j(t,x,\xi)\v\leq \frac 1{\bra
x\ket^{\v A\v}}+\frac{\varepsilon\,M_{k-1}}{\bra x\ket^{\v A\v+1+\sigma_0}}\,.
$$
Therefore if $\v\alpha\v\leq k-1$
\begin{equation}\label{eqVIII.3.5}
\v\partial ^\alpha_X(\xi_p(t)\,\xi_q(t))\v\leq \frac{C_k\,M^2_{k-1}}{\bra
x\ket^{\v A\v}}\,.
\end{equation}
Now if $\v\alpha\v=k$ we have from (\ref{eqVIII.3.3})
\begin{equation}\label{eqVIII.3.6}
\partial ^\alpha_X\,\xi_j(t,x,\xi)=\partial
^\alpha_X\,\zeta_j(t,x,\xi)+R_{j,\alpha}\,,\quad \v R_{j,\alpha}\v\leq \frac 1
{\bra x\ket^{\v A\v}}
\end{equation}
\begin{equation}\label{eqVIII.3.7}\left\{
\begin{array}{l}
\partial ^\alpha_X(\xi_p(t)\cdot \xi_q(t))=\xi_p(t)\,\partial
^\alpha_X\,\xi_q(t)+\xi_q(t)\,\partial ^\alpha_X\,\xi_p(t)+R_{p,q,\alpha}\\
\v R_{p,q,\alpha}\v\leq \frac{C_k\,M^2_{k-1}}{\bra x\ket^{\v A\v}}\,.
\end{array}\right.
\end{equation}
Now we claim that if $F\in\BB^{\ell +1}_{\sigma_0}$ and $\v\alpha\v\leq k-1$ we
have
\begin{equation}\label{eqVIII.3.8}
\v\partial ^\alpha_X[F(x(t))]\v\leq \frac{C_k\,M^{\v\alpha\v}_{k-1}}{\bra
x(t)\ket^{\ell +1+\sigma_0}{\bra x\ket^{\v A\v}}}\,.
\end{equation}
This estimate is easy if $k=1$ and if $k\geq 2$ we use the Faa di Bruno formula
(see Section \ref{ssVIII.1}). It follows that $\partial ^\alpha_X[F(x(t))]$ is
a finite sum of terms of the following form
$$
(\partial ^\beta F)(x(t)) \prod^s_{j=1} (\partial ^{\ell
_j}_X\,x(t))^{k_j}\,,\quad \ell _j\in\N^n\times\N^n\,,\enskip k_j\in\N^n
$$
where $1\leq \v\beta\v\leq \v\alpha\v$, $1\leq s\leq \v\alpha\v$, $\v k_j\v\geq
1$ and
$$
\sum^s_{j=1} k_j=\beta\,,\quad \sum^n_{j=1} \v k_j\v\,\ell _j=\alpha\,.
$$
If we write $\ell _j=(a_j,b_j)$, $\alpha=(A,B)$ we have in particular
\begin{equation}\label{eqVIII.3.9}
\sum^s_{j=1} \v k_j\v\,a_j=A\,.
\end{equation}
Now we write $\{1,2,\ldots ,s\}=I_1\cup I_2\cup I_3$ where
$$
I_1=\{j : \v\ell _j\v\geq 2\}\,,\quad I_2=\{j : \v\ell _j\v=1,\enskip  \ell
_j=(a_j,0)\}\,,\quad I_3=\{j : \v\ell _j\v=1,\enskip \ell _j=(0,b_j)\}
$$
and we denote by $\Sigma_i$ the sum $\som_{j\in I_i}$, $i=1,2,3$.
When $j\in I_1$ we have
$$
\partial ^{\ell _j}_X(x_j(t))=2t\,\partial ^\alpha_X(\zeta_j(t))+\partial
^\alpha_X(z_j(t))\,.
$$
Since $\v\ell _j\v\leq \v\alpha\v\leq k-1$ it follows from the induction that
$$
\v\partial ^\alpha_X(x(t))\v\leq \frac{2t\,\varepsilon\,M_{k-1}}{\bra
x\ket^{\v a_j\v+\sigma_0+1}}+\frac{\varepsilon\,M_{k-1}}{\bra
x\ket^{\v a_j\v+\sigma_0}}=\frac{\varepsilon\,M_{k-1}}{\bra
x\ket^{\v a_j\v+\sigma_0}}\,\Big (1+\frac{2t}{\bra x\ket}\Big)\,.
$$
It follows that
\begin{equation}\label{eqVIII.3.10}
\Big\v\prod_{j\in I_1} \big(\partial ^{\ell _j}_X(x(t))\big)^{k_j}\Big\v\leq
\frac{C_k(\varepsilon\,M_{k-1})^{\Sigma_1\v k_j\v}}{\bra x\ket^{\Sigma_1\v
k_j\v\v a_j\v}}\Big (1+\frac{t^{\Sigma_1\v k_j\v}}{\bra x\ket^{\Sigma_1\v
k_j\v}}\Big)\,.
\end{equation}
Now when $j\in I_2$ we have $\ell _j=(a_j,0)$, $\v a_j\v=1$. Then
$$
\big\v\partial ^{a_j}_x(x_p(t))\big\v\leq 1+\frac{2\varepsilon\,t\,M_1}{\bra
x\ket^{2+\sigma_0}}+\frac{\varepsilon\,M_1}{\bra x\ket^{1+\sigma_0}}\leq
2M_1\Big (1+\frac t{\bra x\ket^{2+\sigma_0}}\Big)\,.
$$
Therefore we have
\begin{equation}\label{eqVIII.3.11}
\Big\v\prod_{j\in I_2} \big(\partial ^{\ell _1}_X(x(t))\big)^{k_j}\Big\v\leq
C(2M_1)^{\Sigma_2 \v k_j\v}\Big (1+\frac{t^{\Sigma_2\v k_j\v}}{\bra
x\ket^{\Sigma_2\v k_j\v(2+\sigma_0)}}\Big)\,.
\end{equation}
Finally for $j\in I_3$ we have $\ell _j=(0,b_j)$, $\v b_j\v=1$. Then
$$
\v\partial ^{b_j}_\xi(x_p(t))\v\leq 2M_1\,\bra t\ket\,.
$$
It follows that
\begin{equation}\label{eqVIII.3.12}
\Big\v\prod_{j\in I_3} \big(\partial ^{\ell _1}_X(x(t))\big)^{k_j}\Big\v\leq
(2M_1)^{\Sigma_3\v k_j\v}\,\bra t\ket^{\Sigma_3\v k_j\v},.
\end{equation}
Using (\ref{eqVIII.3.10}), (\ref{eqVIII.3.11}) and (\ref{eqVIII.3.12}) we obtain
\begin{equation*}
\begin{split}
\hbox to -3cm{}(1)=\Big\v(\partial ^\beta F)(x(t)) \prod^s_{j=1} \big (\partial
^{\ell _j}_X(x(t))\big)^{k_j}\Big\v\leq \frac{C_\beta\,M^{k-1}_{k-1}}{\bra
x(t)\ket^{\v\beta\v+\ell +1+\sigma_0}}\cdot \frac 1{\bra x\ket^{\Sigma_1\v
k_j\v\v a_j\v}}\\
\Big (1+\frac{t^{\Sigma_1\v k_j\v}}{\bra x\ket^{\Sigma_1\v
k_j\v}}\Big)
\Big (1+\frac{t^{\Sigma_2\v k_j\v}}{\bra x\ket^{\Sigma_2\v k_j\v}}\Big)\,\bra
t\ket^{\Sigma_3\v k_j\v}\,.
\end{split}
\end{equation*}
Now we have
$$
\Sigma_1\,\v k_j\v+\Sigma_2\,\v k_j\v+\Sigma_3\,\v k_j\v=\v\beta\v\,,\enskip
\bra x(t)\ket\geq C(t)\enskip \textrm{and}\enskip \bra x(t)\ket\geq C\,\bra
x\ket\,.
$$
It follows that
\begin{equation}\label{eqVIII.3.13}
(1)\leq \frac{C_k\,M^{k-1}_{k-1}}{\bra x(t)\ket^{\ell +1+\sigma_0}}\enskip
\frac 1{\bra x\ket^{\Sigma_1\,\v k_j\v\v a_j\v+\Sigma_2\,\v k_j\v}}\,.
\end{equation}
On the other hand,
$$
\v A\v=\sum^s_{j=1} \v k_j\v\,\v a_j\v=\Sigma_1\,\v k_j\v\,\v
a_j\v+\Sigma_2\,\v k_j\v\enskip \textrm{since}\enskip \v a_j\v=1\enskip
\textrm{for}\enskip t\in I_2
$$
and $a_j=0$ if $j\in I_3$. Therefore (\ref{eqVIII.3.13}) implies
(\ref{eqVIII.3.8}) and our claim is proved.

Moreover if $\v\alpha\v=k$ we can write
\begin{equation}\label{eqVIII.3.14}\left\{
\begin{array}{l}
\partial ^\alpha_X[F(x(t))]=\som^n_{\ell =1} \frac{\partial F}{\partial y_\ell
}\,(x(t))\,\partial ^\alpha_X\,x_\ell (t)+R\enskip \textrm{where}\\
\v R\v\leq \frac{C_k\,M^{k-1}_{k-1}}{\bra x(t)\ket^{\ell +1+\sigma_0}}\cdot
\frac 1 {\bra x\ket^{\v A\v}}\,.
\end{array}\right.
\end{equation}
Indeed $R$ is a finite sum of terms of the form $(\partial^\beta
F)(x(t)){ \prodi^{s}_{j=1}}(\partial ^{\ell _j}_X\,x)^{k_j}$ where
$2\leq
\v\beta\v\leq \v\alpha\v$. It follows then that $\v\ell _j\v\leq
\v\alpha\v-1=k-1$ and the above computations are valid.

Let us now prove Proposition \ref{pIII.3.2} for $\v\alpha\v=k$. Let us set
$$
Z(t)=\partial ^\alpha_X\,z(t)\,,\quad \Xi(t)=\bra x\ket\,\partial
^\alpha_X\,\zeta(t)\,.
$$
We can write
\begin{equation*}
\begin{split}
\dot Z_j(t)&=2\varepsilon
\sum^n_{j,k=1}\bigg\{\underbrace{b_{jk}(x(t))\,\partial
^\alpha_X\,\xi_k(t)}_{(1)}+\underbrace{\partial
^\alpha_X[b_{jk}(x(t))]\,\xi_k(t)}_{(2)}\\ 
&\quad +\sum_{\alpha=\alpha_1+\alpha_2\atop
\alpha_j\not =0}
\begin{pmatrix}\alpha\\ \alpha_1\end{pmatrix}\,\underbrace{\partial
^{\alpha_1}_X[b_{jk}(x(t))]\,\partial
^{\alpha_2}_X\,\xi_k(t)}_{(3)}\bigg\}-2\varepsilon\,t\sum^n_{p,q=1}\bigg\{\underbrace
{\partial ^\alpha_X\Big[\frac{\partial b_{pq}}{\partial
x_j}\,(x(t))\Big]\,\xi_p(t)\,\xi_q(t)}_{(4)}\\
&\quad +\underbrace{\frac{\partial b_{pq}}{\partial x_j}\,(x(t))\,\partial
^\alpha_X(\xi_p(t)\,\xi_q(t))}_{(5)}+\sum_{\alpha=\alpha_1+\alpha_2\atop
\alpha_j\not =0}\begin{pmatrix} \alpha\\
\alpha_1\end{pmatrix}\,\underbrace{\partial ^{\alpha_1}_X\Big[\frac{\partial
b_{pq}}{\partial x_j}\,(x(t))\Big]\,\partial
^{\alpha_2}_X(\xi_p(t)\,\xi_q(t))}_{(6)}\bigg\}\\
\dot\Xi_j(t)&=-\varepsilon\,\bra x\ket \sum^n_{p,q=1}\bigg\{\underbrace{\partial
^\alpha_X\Big[\frac{\partial b_{pq}}{\partial
x_j}\,(x(t))\Big]\,\xi_p(t)\,\xi_q(t) }_{(7)}+\underbrace{\frac{\partial
b_{pq}}{\partial x_j}\,(x(t))\,\partial ^\alpha_X(\xi_p(t)\,\xi_q(t))}_{(8)}\\
&\quad +\sum_{\alpha=\alpha_1+\alpha_2\atop \alpha_j\not
=0}\begin{pmatrix}\alpha\\
\alpha_1\end{pmatrix}\,\underbrace{\partial ^{\alpha_1}_X\Big[\frac{\partial
b_{pq}}{\partial x_j}\,(x(t))\Big]\,\partial
^{\alpha_2}_X(\xi_p(t)\,\xi_q(t))}_{(9)}\bigg\} \,.
\end{split}
\end{equation*}
We shall use the fact that $b_{jk}\in\BB^1_{\sigma_0}$ and $\frac{\partial
b_{pq}}{\partial x_j}\in\BB^2_{\sigma_0}$. We deduce from (\ref{eqVIII.3.3})
that
\begin{equation}\label{eqVIII.3.15}
\v(1)\v\leq \frac{C\,\varepsilon}{\bra x(t)\ket^{1+\sigma_0}\bra x\ket^{\v
A\v}}+\frac{C\,\varepsilon}{\bra x(t)\ket^{1+\sigma_0}}\enskip
\frac{\v\Xi(t)\v}{\bra x\ket}\,.
\end{equation}
To estimate the term (2) we use (\ref{eqVIII.3.14}) with $F=b_{jk}$ and the
equality
$$
\partial ^\alpha_X\,x_j(t)=\partial
^\alpha_X(x_j+2t\,\xi_j)+2t\,\frac{\Xi_j(t)}{\bra x\ket}+Z_j(t)\,.
$$
We obtain
\begin{equation}\label{eqVIII.3.16}
\v(2)\v\leq \frac{C\,\varepsilon}{\bra x(t)\ket^{2+\sigma_0}}\,\v Z(t)\v+
\frac{C\,\varepsilon}{\bra x(t)\ket^{1+\sigma_0}}\enskip \frac{\v\Xi(t)\v}{\bra
x\ket}+\frac{C\,\varepsilon}{\bra x(t)\ket^{1+\sigma_0}\bra x\ket^{\v A\v}}\,.
\end{equation}
To estimate (3) we use (\ref{eqVIII.3.4}) and (\ref{eqVIII.3.8})~; we obtain,
\begin{equation}\label{eqVIII.3.17}
\v(3)\v\leq \frac{C_k\,\varepsilon\,M^{k-1}_{k-1}}{\bra
x(t)\ket^{1+\sigma_0}\bra x\ket^{\v A\v}}\,.
\end{equation}
By the same way we have,
\begin{equation}\label{eqVIII.3.18}
\v(6)\v\leq \frac{C_k\,\varepsilon\,M^{k-1}_{k-1}}{\bra
x(t)\ket^{1+\sigma_0}\bra x\ket^{\v A\v}}\,.
\end{equation}
Then $(4)$ has the same estimate as $(2)$ and $(5)$ as $(1)$ since $\bra
x(t)\ket\geq C\,\bra t\ket$.

To take care of $(7)$ we use (\ref{eqVIII.3.14}) with $F=\frac{\partial
b_{pq}}{\partial x_j}\in\BB^2_{\sigma_0}$. We obtain
\begin{equation}\label{eqVIII.3.19}
\v(7)\v\leq \frac{C\,\varepsilon}
{\bra x(t)\ket^{2+\sigma_0}}\,\v Z(t)\v+\frac{C\,\varepsilon}
{\bra x(t)\ket^{2+\sigma_0}}\,\v\Xi(t)\v+\frac{C\,\varepsilon}
{\bra x(t)\ket^{1+\sigma_0}\bra x\ket^{\v A\v}}
\end{equation}
since $\bra x(t)\ket\geq C\,\bra x\ket$ and $\bra x(t)\ket\geq C\,\bra t\ket$.
Finally
\begin{equation}\label{eqVIII.3.20}\left\{
\begin{array}{l}
\v(8)\v\leq \frac{C\,\varepsilon}
{\bra x(t)\ket^{1+\sigma_0}}\,\v\Xi(t)\v+\frac{C\,\varepsilon}
{\bra x(t)\ket^{1+\sigma_0}\bra x\ket^{\v A\v}}\,,\\
\v(9)\v\leq \frac{C\,\varepsilon\,M^{k-1}_{k-1}}{\bra x(t)\ket^{1+\sigma_0}\bra
x\ket^{\v A\v}}\,.
\end{array}\right.
\end{equation} 
Gathering the estimates obtained in (\ref{eqVIII.3.15}) to (\ref{eqVIII.3.20})
we obtain
$$
\v\dot Z(t)\v+\v\dot\Xi(t)\v\leq \frac{C\,\varepsilon}{\bra
x(t)\ket^{1+\sigma_0}}\,(\v
Z(t)\v+\v\Xi(t)\v)+\frac{C_k\,\varepsilon\,M^{k-1}_{k-1}}{\bra
x(t)\ket^{1+\sigma_0}\bra x\ket^{\v A\v}}\,.
$$
It follows from Gromwall's Lemma, (\ref{eqVIII.3.2}) and the estimate $\bra
x(t)\ket\geq C\,\bra t\ket$ that
$$
\v Z(t)\v+\v\Xi(t)\v\leq \frac{C(M_{k-1})\,\varepsilon}{\bra x\ket^{\v
A\v+\sigma_0}}
$$
which, according to the definition of $Z$ and $\Xi$ proves Proposition
\ref{pIII.3.2} when $\v\alpha\v=k$. \cqfd

\subsection{\quad Proof of Lemma \ref{lV.3.1}}\label{ssVIII.4}

The proof is the same for the two cases so we shall consider the more general
case where $f=f(x,\theta)$.

Let $\chi\in C^\infty _0(\R^n)$, $\chi(\xi)=1$ if $\v\xi\v\leq \frac 12$,
$\chi(\xi)=0$ if $\v\xi\v\geq 1$. We set
\begin{equation}\label{eqVIII.4.1}
\sum_{\v\beta\v\leq \ell }\, \sup_{\R^n}\, \v\partial ^\beta \chi\v=D_\ell
\,,\quad \ell \in\N\,.
\end{equation}
We want to show that one can find an increasing sequence $(L_k)_{k\geq 1}$ in
$]1,+\infty [$ such that if we set for $(\theta,x,y)$ in $\Omega\times\R^n_y$
\begin{equation}\label{eqVIII.4.2}\left\{
\begin{array}{l}
F_\gamma(\theta,x,y)=\partial
^\gamma_x\,f(\theta,x)\,\frac{(iy)^\gamma}{\gamma!}\, \chi\Big
(L_{\v\gamma\v}\,y\Big (\frac 1{\bra x\ket}+\frac
1{\bra\theta\ket}\Big)\Big)\,,\quad \v\gamma\v\geq 1\,,\\
F(\theta,x,y)=f(\theta,x)+\som_{\gamma\not =0} F_\gamma(x,\theta,y)
\end{array}\right.
\end{equation}
then $F$ is well defined in $\Omega\times\R^n_y$ and satisfies all the
requirements of Lemma \ref{lV.3.1}.

First of all in the expression of $F_\gamma$, on the support of $\chi$ we have
$\v y\v\leq \frac 1{L_{\v\gamma\v}}\,\frac{\bra\theta\ket \bra
x\ket}{\bra\theta\ket+\bra x\ket}$. It follows that
\begin{equation}\label{eqVIII.4.3}
\frac{\v y\v}{\bra x\ket}\leq \frac 1{L_{\v\gamma\v}}\,,\quad \frac{\v
y\v}{\bra\theta\ket}\leq \frac 1{L_{\v\gamma\v}}\,.
\end{equation}
Using (\ref{eqV.3.2}) we can write,
$$
\v F_\gamma(\theta,x,y)\v\leq M_{\v\gamma\v}\Big(\frac 1{\bra
x\ket^{\v\gamma\v+\sigma_3}}+\frac 1{\bra\theta\ket^{\v\gamma\v+\sigma_3}}\Big
)\,\frac{\v y\v^{\v\gamma\v}}{\gamma!}\, \v\chi(L_{\v\gamma\v}\,y(\cdots ))\v
$$
so
$$
\v F_\gamma(\theta,x,y)\v\leq
\frac{D_0\,M_{\v\gamma\v}}{\gamma!\,L^{\v\gamma\v}_{\v\gamma\v}}\Big (\frac
1{\bra x\ket^{\sigma_3}}+\frac 1{\bra\theta\ket^{\sigma_3}}\Big )\,.
$$
Taking
\begin{equation}\label{eqVIII.4.4}
L^{\v\gamma\v}_{\v\gamma\v}\geq D_0\,M_{\v\gamma\v}
\end{equation}
we deduce that $F$ defined in (\ref{eqVIII.4.2}) is well defined and satisfies
\begin{equation}\label{eqVIII.4.5}
\v F(\theta,x,y)\v\leq C_0 \Big (\frac 1{\bra x\ket^{\sigma_1}}+\frac
1{\bra\theta\ket^{\sigma_2}}\Big)
\end{equation}
since $\sigma_3\geq \sigma_1$ and $\sigma_3\geq \sigma_2$.

Therefore (ii) in Lemma \ref{lV.3.1} is satisfied and (i) follows immediately
from (\ref{eqVIII.4.2}).
We shall strenghten the condition (\ref{eqVIII.4.4}) on $L_{\v\gamma\v}$ to
obtain a $C^\infty $ function $F$. First all there exists absolute constants
$C_{i,\ell }$, $(i\in\N$, $\ell \in\N$) such that
\begin{equation}\label{eqVIII.4.6}
\sum_{\v\gamma\v=i}\enskip  \bigg\v\partial ^\gamma_x\Big[\Big(\frac 1{\bra
x\ket}+\frac 1{\bra\theta\ket}\Big)^\ell \Big]\bigg\v\leq C_{i,\ell }\Big(\frac
1{\bra x\ket}+\frac 1{\bra\theta\ket}\Big)^{\ell +i}\,.
\end{equation}
Let $\beta\in\N^n$. For any $\mu\in\N^n$ one can find an absolute constant
$K_{\v\mu\v}$ independent of $(L_k)$) such that for all $(\theta,x,y)$ in
$\R\times\R^n_x\times\R^n_y$ we have
\begin{equation}\label{eqVIII.4.7}
\bigg\v\partial ^\mu_x\Big[(\partial
^\beta_\xi\chi)\Big(L_{\v\gamma\v}y\Big(\frac 1{\bra x\ket}+\frac
1{\bra\theta\ket}\Big)\Big)\Big]\bigg\v\leq
K_{\v\mu\v}\,D_{\v\beta\v+\v\mu\v}\Big (\frac 1{\bra x\ket}+\frac
1{\bra\theta\ket}\Big)^{\v\mu\v}\,.
\end{equation}
Indeed let us set $h(\theta,x,y)=L_{\v\gamma\v}\,y\big (\frac 1{\bra
x\ket}+\frac 1{\bra\theta\ket}\big)$. By the Faa di Bruno formula, $\partial
^\mu_x[(\partial ^\beta\chi)(h(\theta,x,y))]$ is a finite linear combination
with absolute coefficients of terms of the form $(\partial
^{(\beta+\nu)}_\xi\chi)(h(\theta,x,y))\break \prodi^s_{j=1} (\partial
^{\ell _j}_x\,h(\theta,x,y))^{k_j}$ where $1\leq \v\nu\v\leq \v\mu\v$, $1\leq
s\leq
\v\mu\v$, $\v k_j\v\geq 1$, $\v\ell _j\v\geq 1$ and
$$
\sum^s_{j=1}\,k_j=\nu\,,\quad \sum^s_{j=1}\, \v k_j\v\,\ell _j=\mu\,.
$$
Since $\v\nu\v+\v\beta\v\leq \v\mu\v+\v\beta\v$ we have
\begin{equation}\label{eqVIII.4.8}
\v\partial _\xi^{(\beta+\nu)}\chi(h(\theta,x,y))\v\leq D_{\v\mu\v+\v\beta\v}\,.
\end{equation}
On the other hand it follows from (\ref{eqVIII.4.6}) that
$$
(1)=:\Big\v\prod^s_{j=1} (\partial ^{\ell _j}_x\,h(\theta,x,y))^{k_j}\Big\v\leq
\prod^s_{j=1}\, C^{\v k_j\v}_{\v\ell _j\v,1}(\v y\v\,L_{\v\gamma\v})^{\som^s_1
\v k_j\v}\Big (\frac 1{\bra x\ket}+\frac 1{\bra\theta\ket}\Big)^{\som^s_1 \v
k_j\v(1+\v\ell _j\v)}\,.
$$
On the support of $\chi$ we use the estimates (\ref{eqVIII.4.3}). Moreover we
have $\som^s_1\,\v k_j\v=\v\nu\v$, $\som^s_1\,\v k_j\v\,\v\ell _j\v=\v\mu\v$. It
follows then that,
\begin{equation}\label{eqVIII.4.9}
(1)\leq C'_{\v\mu\v}\Big (\frac 1{\bra x\ket}+\frac1
{\bra\theta\ket}\Big)^{\v\mu\v}\,.
\end{equation}
Then (\ref{eqVIII.4.7}) follows from (\ref{eqVIII.4.8}) and (\ref{eqVIII.4.9}).

Now with $F_\gamma$ defined in (\ref{eqVIII.4.2}) we can write
\begin{equation*}
\begin{split}
\partial ^B_y\,F_\gamma(\theta,x,y)&=i^{\v\gamma\v}\,\frac{\partial
^\gamma_x\,f(\theta,x)}{\gamma!} \sum_{B_1\leq B\atop B_1\leq \gamma}
\begin{pmatrix}B\\ B_1\end{pmatrix} \frac{Á\,y^{\gamma-B_1}}{(\gamma-B_1)!}
\Big[L_{\v\gamma\v}\Big(\frac 1{\bra x\ket}+\frac
1{\bra\theta\ket}\Big)\Big]^{\v B\v-\v B_1\v}\\
&\qquad\qquad(\partial ^{B-B_1}_\xi\chi)\Big(L_{\v\gamma\v}\,y\Big(\frac 1{\bra
x\ket}+\frac 1{\bra\theta\ket}\Big)\Big)\,.
\end{split}
\end{equation*}
Then
$$
\partial ^A_x\,\partial ^B_y\,F_\gamma(\theta,x,y)=i^{\v\gamma\v} 
\sum_{B_1\leq B\atop B_1\leq \gamma}\enskip \sum_{A_1\leq A}\enskip
\sum_{A_2\leq A-A_1}\begin{pmatrix}B\\ B_1\end{pmatrix}
\begin{pmatrix}A\\ A_1\end{pmatrix}\begin{pmatrix}A-A_1\\ A_2\end{pmatrix}
\frac 1{(\gamma-B_1)!}\,y^{\gamma-B_1}
$$
$$
L^{\v B\v-\v B_1\v}_{\v\gamma\v}\,\partial
^{\gamma+A_1}_x\,f(\theta,x)\,\partial ^{A_2}_x\Big[\Big (\frac 1{\bra
x\ket}+\frac 1{\bra\theta\ket}\Big)^{\v B\v-\v B_1\v}\Big]\,\partial
^{A-A_1-A_2}_x \Big[(\partial ^{B-B_1}_\xi\,\chi)\Big (y\,L_{\v\gamma\v}\Big (\frac 1{\bra
x\ket}+\frac 1{\bra\theta\ket}\Big)\Big)\Big]\,.
$$
Now we use (\ref{eqVIII.4.3}), (\ref{eqVIII.4.6}) and (\ref{eqVIII.4.7}). Since
$\v\gamma\v+\v A_1\v\leq \v\gamma\v+\v A\v$, $\v A-A_1-A_2\v\leq \v A\v$, $\v
B-B_1\v\leq \v B\v$ we obtain
$$
\v\partial ^A_x\,\partial ^B_y\,F_\gamma(\theta,x,y)\v\leq 
\sum_{B_1\leq B\atop B_1\leq \gamma}\enskip \sum_{A_1\leq A\atop A_2\leq
A-A_1}\enskip
\begin{pmatrix}B\\ B_1\end{pmatrix}
\begin{pmatrix}A\\ A_1\end{pmatrix}\begin{pmatrix}A-A_1\\ A_2\end{pmatrix}
\,2^{\v\gamma\v-\v B_1\v}\,L^{\v B_1\v-\v\gamma\v}_{\v\gamma\v}
\Big (\frac 1{\bra x\ket}+\frac 1{\bra\theta\ket}\Big)^{\v B_1\v-\v\gamma\v}
$$
\begin{equation*}
\begin{split}
L^{\v B\v-\v B_1\v}_{\v\gamma\v}\,M_{\v\gamma\v+\v A\v}&\Big
(\frac 1{\bra x\ket^{\v\gamma\v+\v A_1\v+\sigma_3}}+\frac 1{\bra
\theta\ket^{\v\gamma\v+\v A_1\v+\sigma_3}}\Big)\,C_{\v A_2\v,\v B-B_1\v}\\
&\qquad\qquad\qquad\Big (\frac 1{\bra x\ket}+\frac 1{\bra
\theta\ket}\Big)^{\v A_2\v+\v B\v-\v B_1\v}\,K_{\v A\v}\,D_{\v A\v+\v B\v}
\Big (\frac 1{\bra x\ket}+\frac 1{\bra \theta\ket}\Big)^{\v A\v-\v A_1\v-\v
A_2\v}
\end{split}
\end{equation*}
where $\v y\v$ satisfies (\ref{eqVIII.4.3}).

It follows that we can find a constant $\tilde C_{\v A\v,\v B\v}\leq 1$,
depending only on $\v A\v$, $\v B\v$ and the dimension such that
\begin{equation}\label{eqVIII.4.10}
\v\partial ^A_x\,\partial ^B_y\,F_\gamma(\theta,x,y)\v\leq \frac{\tilde C_{\v
A\v,\v B\v}\,D_{\v A\v+\v B\v}\,M_{\v\gamma\v+\v A\v}}{L^{\v\gamma\v-\v
B\v}_{\v\gamma\v}}
\Big(\frac 1{\bra x\ket^{\v A\v+\v B\v+\sigma_3}}+\frac 1{\bra\theta\ket^{\v
A\v+\v B\v+\sigma_3}}\Big)\,.
\end{equation}
Let us set
\begin{equation}\label{eqVIII.4.11}
\tilde C_k=\max_{\v A\v+\v B\v\leq k}\,\tilde C_{\v A\v,\v B\v}\,.
\end{equation}
We shall take the increasing sequence $(L_k)_{k\geq 0}$ such that
\begin{equation}\label{eqVIII.4.12}
L_k\geq \max \big (1,2^k\,\tilde C_k\,D_k\,M_{2k}\big)\,.
\end{equation}
Then we write, according to (\ref{eqVIII.4.2})
\begin{equation}\label{eqVIII.4.13}
F(\theta,x,y)=f(\theta,x)+\sum_{\v\gamma\v\leq \v A\v+\v B\v}
F_\gamma(\theta,x,y)+\sum_{\v\gamma\v>\v A\v+\v B\v} F_\gamma(\theta,x,y)\,.
\end{equation}
The first two terms in the right hand side of (\ref{eqVIII.4.13}) define a
$C^\infty $ function. For the third one we deduce from (\ref{eqVIII.4.10}),
(\ref{eqVIII.4.11}) and (\ref{eqVIII.4.12}) that
$$
\v\partial ^A_x\,\partial ^B_y\,f_\gamma(\theta,x,y)\v\leq \frac{\tilde
C_{\v\gamma\v}\,D_{\v\gamma\v}\,M_{2\v\gamma\v}}{L_{\v\gamma\v}}\leq \frac
1{2^{\v\gamma\v}}\,.
$$
This shows that the third term define also a $C^\infty $ function. Thus $F$ is
$C^\infty $ in $(x,y)$. Let us prove (iii). According to (\ref{eqVIII.4.2}) if
$\v A\v+\v B\v\geq 1$ we have
$$
(1)=:\partial ^A_x\,\partial ^B_y\,F(\theta,x,y)=\partial ^A_x\,\partial
^B_y\,f(\theta,x)+\sum_{\gamma\not =0} \partial ^A_x\,\partial
^B_y\,F_\gamma(\theta,x,y)\,.
$$
If $\v B\v\geq 1$ we use (\ref{eqVIII.4.10}), (\ref{eqVIII.4.13}) and
(\ref{eqVIII.4.12}). We get
$$
\v(1)\v\leq C_{AB}\Big (\frac 1{\bra x\ket^{\v A\v+\v B\v+\sigma_3}}+\frac
1{\bra\theta\ket^{\v A\v+\v B\v+\sigma_3}}\Big)\,.
$$
If $B=0$ we use furthermore (\ref{eqV.3.2}). We obtain
$$
\v(1)\v\leq M_{\v A\v}\Big (\frac 1{\bra x\ket^{\v A\v+\sigma_3}}+
\frac 1{\bra \theta\ket^{\v A\v+\sigma_3}}\Big )+C'_{AB}\Big (\frac 1{\bra
x\ket^{\v A\v+\sigma_3}}+\frac 1{\bra \theta\ket^{\v A\v+\sigma_3}}\Big )\,.
$$
This shows that (iii) holds.

Finally let us prove (iv). Let us set again
$$
h_{\v\gamma\v}(\theta,x,y)=L_{\v\gamma\v}\,y\Big (\frac 1{\bra x\ket}+\frac
1{\bra\theta\ket}\Big)\,.
$$
According to (\ref{eqVIII.4.2}) we have for $j\in \{1,2,\ldots ,n\}$,
\begin{equation*}
\begin{aligned}
\partial _{x_j}F&=\partial _{x_j}f(\theta,x)+\sum_{\gamma\not =0} \Big
(\partial _{x_j}\,\partial
^\gamma_x\,f(\theta,x)\,\chi(h_{\v\gamma\v}(\theta,x,y))\\
&\qquad\qquad\qquad -\sum^n_{\ell =1}
L_{\v\gamma\v}\,\frac{x_j\,y_\ell }{\bra x\ket^3}\,\partial
^\gamma_x\,f(\theta,x)\cdot \frac{\partial \chi}{\partial \xi_\ell
}\,(h_{\v\gamma\v}(\theta,x,y))\Big)\,\frac{(iy)^\gamma}{\gamma!}\\
\partial _{y_j}F&=\sum_{\gamma\not =0} \partial
^\gamma_x\,f(\theta,x)\,\frac{(iy)^\gamma}{\gamma!}\,L_{\v\gamma\v}\Big ( \frac
1{\bra x\ket}+\frac 1{\bra\theta\ket}\Big)\,\frac{\partial \chi}{\partial
\xi_j}\,(h_{\v\gamma\v}(\theta,x,y))\\
&\qquad\qquad\qquad +i\,\sum_{\gamma_j\geq 1} \partial
^\gamma_x\,f(\theta,x)\,\frac
{(iy)^{\gamma-e_j}}{(\gamma-e_j)!}\,\chi(h_{\v\gamma\v}(\theta,x,y))\,,
\end{aligned}
\end{equation*}
where $e_j=(0,0,\ldots ,1,\ldots ,0)$. Setting $\gamma'=\gamma-e_j$ the sum
above can be written
$$
i\,\partial _{x_j}\,f(\theta,x)\,\chi\Big(L_1\,y\Big(\frac 1{\bra x\ket}+\frac
1{\bra\theta\ket}\Big)\Big)+i\,\sum_{\gamma\not =0} \partial _{x_j}\,\partial
^\gamma_x\,f(\theta,x)\,\frac{(iy)^\gamma}{\gamma!}\,\chi
(h_{\v\gamma\v+1}(\theta,x,y))\,.
$$
It follows that
\begin{equation}\label{eqVIII.4.14}
(\partial _{x_j}F+i\,\partial _{y_j}F)(\theta,x,y)=(1)+(2)+(3)
\end{equation}
where
\begin{equation}\label{eqVIII.4.15}\left\{
\begin{array}{l}
(1)=\Big (1-\chi\Big (L_1\,y\Big(\frac 1{\bra x\ket}+\frac
1{\bra\theta\ket}\Big)\Big)\,\partial _{x_j}\,f(\theta,x)\\
(2)=\som_{\gamma\not =0} \partial _{x_j}\,\partial
^\gamma_x\,f(\theta,x)\,\frac{(iy)^\gamma}{\gamma!}
[\chi(h_{\v\gamma\v}(\theta,x,y))-\chi(h_{\v\gamma\v+1}(\theta,x,y))]\\
(3)=\som_{\gamma\not =0} \partial
_x^\gamma\,f(\theta,x)\,\frac{(iy)^\gamma}{\gamma!}\,L_{\v \gamma\v}\Big
(-\som^n_{\ell =1}
\frac{x_j\,y_\ell }{\bra x\ket^3}\,\frac{\partial \chi}{\partial \xi_\ell
}\,(h_{\v\gamma\v}(\theta,x,y))\\
\hbox to 5cm{}+i\Big (\frac 1{\bra x\ket}+\frac
1{\bra\theta\ket}\Big)\,\frac{\partial
\chi}{\partial \xi_j}\,(h_{\v\gamma\v}(\theta,x,y))\Big)\,.
\end{array}\right.
\end{equation}

Let us set for convenience
\begin{equation}\label{eqVIII.4.16}
R=\v y\v \Big(\frac 1{\bra x\ket}+\frac 1{\bra\theta\ket}\Big)\,.
\end{equation}
On the support of $1-\chi\big(L_1\,y\big(\frac 1{\bra x\ket}+\frac
1{\bra\theta\ket}\big)\big)$ in the term (1) above, we have $L_1R\geq \frac 12$.
Therefore
\begin{equation}\label{eqVIII.4.17}
\frac{(1)}{R^N}\leq (2L_1)^N\,M_1\Big (\frac 1{\bra x\ket^{1+\sigma_3}}+\frac
1{\bra\theta\ket^{1+\sigma_3}}\Big)\,.
\end{equation}
On the support of
$\chi(h_{\v\gamma\v}(\theta,x,y))-\chi(h_{\v\gamma\v+1}(\theta,x,y))$ we have
$\frac 1{2L_{\v\gamma\v+1}}\leq R\leq \frac 1{L_{\v\gamma\v}}$. Now we write
with $N\geq 2$,
$$
(2)=\sum_{1\leq \v\gamma\v\leq N-1} G_\gamma(\theta,x,y)+\sum_{\v\gamma\v\geq
N} G_\gamma(\theta,x,y)\,.
$$
When $\v\gamma\v\leq N-1$ we have $L_{\v\gamma\v+1}\leq L_N$ so $R\geq \frac
1{2L_N}$. It follows that
$$
\sum_{1\leq \v\gamma\v\leq N-1} \v G_\gamma(\theta,x,y)\v\cdot \frac 1{R^N}\leq
(2L_N)^N \sum_{1\leq \v\gamma\v\leq N-1} M_{\v\gamma\v+1}\cdot \frac
1{L^{\v\gamma\v}_{\v\gamma\v}}\Big(\frac 1{\bra x\ket^{1+\sigma_3}}+\frac
1{\bra\theta\ket^{1+\sigma_3}}\Big)\,.
$$
Therefore we obtain
\begin{equation}\label{eqVIII.4.18}
\sum_{1\leq \v\gamma\v\leq N-1} \v G_\gamma(\theta,x,y)\v\leq C_N\,R^N 
\Big(\frac 1{\bra
x\ket^{1+\sigma_3}}+\frac 1{\bra\theta\ket^{1+\sigma_3}}\Big)\,.
\end{equation}
On the other hand we have
\begin{equation*}
\begin{aligned}
\sum_{\v\gamma\v\geq N} \v G_\gamma(\theta,x,y)\v&\leq \Big(\frac 1{\bra
x\ket^{1+\sigma_3}}+\frac 1{\bra\theta\ket^{1+\sigma_3}}\Big)\,R^N
\sum_{\v\gamma\v\geq N}
M_{\v\gamma\v}\,R^{\v\gamma\v-N}[\chi(h_{\v\gamma\v})-\chi(h_{\v\gamma\v+1})]\\
\sum_{\v\gamma\v\geq N} \v G_\gamma(\theta,x,y)\v&\leq \Big(\frac 1{\bra
x\ket^{1+\sigma_3}}+\frac 1{\bra\theta\ket^{1+\sigma_3}}\Big)\,R^N \sum_\beta 
M_{\v\beta\v+N}\,R^{\v\beta\v}[\chi(h_{\v\beta\v+N})-\chi(h_{\v\beta\v+N+1})]\,.
\end{aligned}
\end{equation*}
On the support of $\chi(h_{\v\beta\v+N})-\chi(h_{\v\beta\v+N+1})$, we have
$L_{\v\beta\v+N}\,R\leq 1$. It follows that
\begin{equation}\label{eqVIII.4.19}
\sum_{\v\gamma\v\geq N} \v G_\gamma(\theta,x,y)\v\leq C_N\,R^N \Big (\frac
1{\bra x\ket^{1+\sigma_3}}+\frac 1{\bra\theta\ket^{1+\sigma_3}}\Big)\,.
\end{equation}
Combining (\ref{eqVIII.4.18}) and (\ref{eqVIII.4.19}) we obtain
\begin{equation}\label{eqVIII.4.20}
\v(2)\v\leq C'_N\,R^N\Big (\frac
1{\bra x\ket^{1+\sigma_3}}+\frac 1{\bra\theta\ket^{1+\sigma_3}}\Big)\,.
\end{equation}
Finally we consider the term (3) in (\ref{eqVIII.4.15}). We have $\frac{\v
x_j\v}{\bra x\ket}\leq 1$ and on the support of $\frac{\partial \chi}{\partial
\xi_j}\,(h_{\v\gamma\v}(\theta,x,y))$ we have $L_{\v\gamma\v}\,\frac{\v
y\v}{\bra x\ket}\leq 1$ it follows that (3) is bounded by a finite sum of terms
of the form
$$
(3)'=\Big (\frac 1{\bra x\ket^{1+\sigma_3}}+\frac1{\bra\theta\ket^{1+\sigma_3}}
\Big) \sum_{\gamma\not =0} M_{\v\gamma\v}
\,R^{\v\gamma\v}\,L_{\v\gamma\v}\,\Big\v\frac{\partial\chi}{\partial 
\xi_j}\,(h_{\v\gamma\v}(\theta,x,y))\Big\v\,.
$$
As before we write the above sum as $\som_{1\leq \v\gamma\v\leq
N-1}+\som_{\v\gamma\v\geq N}$. If $\v\gamma\v\leq N-1$ then $L_{\v\gamma\v}\leq
L_{N-1}$ so on the support of $\frac{\partial \chi}{\partial \xi_j}$ we have
$R\geq \frac 1{2L_{\v\gamma\v}}\geq \frac 1{2L_{N-1}}$ and $R\leq \frac
1{L_{\v\gamma\v}}$. It follows that
$$
\frac 1{R^N} \Big\v\sum_{1\leq \v\gamma\v\leq N-1}\Big\v\leq (2L_{N-1})^N
\sum_{1\leq \v\gamma\v\leq N-1} M_{\v\gamma\v}\,\frac
1{L^{\v\gamma\v-1}_{\v\gamma\v}}\,D_1=C_N\,.
$$
For the second sum we write
\begin{equation*}
\begin{aligned}
\Big\v\sum_{\v\gamma\v\geq N}\Big\v&\leq R^N \sum_{\v\gamma\v\geq N}
M_{\v\gamma\v}\,R^{\v\gamma\v-N}\,L_{\v\gamma\v}\,\Big\v\frac{\partial
\chi}{\partial \xi_j}\,(h_{\v\gamma\v})\Big\v\\
&\leq R^N \sum_\beta M_{\v\beta\v+N}\,\frac{L_{\v\beta\v+N}}{L^{\v\beta\v}_{\v
\beta\v+N}}\,D_1=C'_N\,R^N\,.
\end{aligned}
\end{equation*}
It follows that
\begin{equation}\label{eqVIII.4.21}
\v(3)\v\leq C_N\,R^N \Big (\frac 1{\bra x\ket^{1+\sigma_3}}+\frac
1{\bra\theta\ket^{1+\sigma_3}}\Big)\,.
\end{equation}
Using (\ref{eqVIII.4.14}) to (\ref{eqVIII.4.21}) we obtain the part (iv) of
Lemma \ref{lV.3.1}. The proof is complete. \cqfd

\end{document}

\end{document}